\baselineskip=13pt plus 2pt 
\magnification=1150
\font\bigtenrm=cmr10 scaled \magstep2
\def\sqr#1#2{{\vcenter{\vbox{\hrule height.#2pt\hbox{\vrule width.#2pt
height#1pt\kern#1pt \vrule width.#2pt}\hrule height.#2pt}}}}
\def\square{\mathchoice\sqr64\sqr64\sqr{2.1}3\sqr{1.5}3}

\noindent{\bigtenrm On tau functions associated with 
linear systems}
\vskip.1in
\noindent {Gordon Blower*}\par

\noindent {\sl Department of Mathematics and Statistics} {\sl Lancaster University, Lancaster LA1 4YF,  UK}\par
\noindent {\sl E-mail address: g.blower@lancaster.ac.uk}\par
\vskip.1in
\noindent {Samantha Newsham}\par
\noindent {\sl Lancaster University, Lancaster LA1 4YF, UK}\par
\noindent {\sl samantha.l.newsham@gmail.com}\par
\vskip.1in
\noindent {26th June 2017}\par
\vskip.1in
\noindent ${}^*$ corresponding author\par
\vskip.1in
\hrule
\vskip.1in 
\noindent{ABSTRACT}\par
\noindent  Let $(-A,B,C)$ be a linear system in continuous time $t>0$ with input and output space ${\bf C}$ and state space $H$. The function $\phi_{(x)}(t)=Ce^{-(t+2x)A}B$ determines a Hankel integral operator $\Gamma_{\phi_{(x)}}$ on $L^2((0, \infty ); {\bf C})$; if   $\Gamma_{\phi_{(x)}}$ is trace class, then the Fredholm determinant $\tau (x)=\det (I+ \Gamma_{\phi_{(x)}})$ defines the tau function of $(-A,B,C)$. Such tau functions arise in Tracy and Widom's theory of matrix models, where they describe the fundamental probability distributions of random matrix theory. Dyson considered such tau functions in the inverse spectral problem for Schr\"odinger's equation $-f''+uf=\lambda f$, and derived the formula for the potential $u(x)=-2{{d^2}\over{dx^2}}\log \tau (x)$ in the self-adjoint scattering case  {\sl Commun. Math. Phys.} {\bf 47} (1976), 171--183. This paper introduces a operator function $R_x$ that satisfies Lyapunov's equation ${{dR_x}\over{dx}}=-AR_x-R_xA$ and $\tau (x)=\det (I+R_x)$, without assumptions of self-adjointness. When $-A$ is sectorial, and $B,C$ are Hilbert--Schmidt, there exists a non-commutative differential ring
${\cal A}$ of operators in $H$ and a differential ring homomorphism $\lfloor\,\,\rfloor :{\cal A}\rightarrow {\bf C}[u,u', \dots ]$ such that $u=-4\lfloor A\rfloor$, which provides a substitute for the multiplication rules for Hankel operators considered by P\"oppe, and McKean  {\sl Cent. Eur. J. Math.} {\bf
9} (2011), 205--243. The paper obtains conditions on $(-A,B,C)$ for Schr\"odinger's equation with meromorphic $u$ to be integrable by quadratures. Special results apply to the linear systems associated with scattering $u$, periodic $u$ and elliptic $u$. The paper constructs a family of solutions to the Kadomtsev--Petviashivili differential equations, and proves that certain families of tau functions satisfy Fay's identities.\par
\vskip.05in

\noindent 5, 34B25\par
\par
\vskip.1in
\hrule
\vskip.1in
\noindent {\bf 1. Introduction}\par
\vskip.05in
\indent This paper is concerned with Fredholm determinants which arise in the theory of linear systems and their application to differential equations such as the Kadomtsev--Petviashvili equation. For $\phi \in L^2((0, \infty ); {\bf R})$, the Hankel integral operator corresponding to $\phi$ is 
$\Gamma_\phi$ where 
$$\Gamma_{\phi }f(x)=\int_0^\infty \phi (x+y)f(y)\, dy\qquad (f\in L^2((0, \infty ); {\bf C}).\eqno(1.1)$$
\noindent Using the Laguerre system of orthogonal functions as in [60], one can express $\Gamma_\phi$ as a matrix $[\gamma_{j+k}]_{j,k=1}^\infty$ on $\ell^2$, which has the characteristic shape of a Hankel matrix, and one can establish criteria for the operator to be bounded on $L^2((0,\infty ); {\bf C})$. Megretski, Peller and Treil [52] determined the possible spectrum and spectral multiplicity function that can arise from a bounded and self-adjoint Hankel operator. Thus they characterized the class of bounded self-adjoint Hankel operators up to unitary equivalence. Their method involved introducing suitable linear systems on a state space $H$, and this motivated the approach of our paper.\par
\indent Previously, Dyson [20] considered the inverse spectral problem for Schr\"odinger's equation $-f''+uf=\lambda f$, where $u\in C^2({\bf R}; {\bf R})$ that decays rapidly as $x\rightarrow\pm \infty$. From the asymptotic solutions, he introduced a scattering function $\phi$, and considered the translations $\phi_{(x)}(y)=\phi (y+2x)$. He showed that the potential can be recovered from the scattering data by means of the formula
$$u(x)=-2{{d^2}\over{dx^2}}\log\det (I+\Gamma_{\phi_{(x)}}).\eqno(1.2)$$
\indent These results were developed further by Ercolani, McKean [22] and others [36, 69, 73] to describe the inverse spectral problem for self-adjoint Schr\"odinger operators on ${\bf R}$. Remarkably, some of the methods of inverse scattering theory do not really need self-adjointness. However, a significant obstacle in this approach is that Hankel operators do not have a  natural product structure, so it is unclear as to how one can fully exploit the multiplicative properties of determinants. This paper seeks to address this issue, by realizing Hankel operators from linear systems, and then introducing algebras of operators which reflect the properties of Hankel operators and their Fredholm determinants. As in [52], the Lyapunov differential equation is fundamental to the development of the theory.\par   

\vskip.05in
\noindent {\bf Definition} (i) {\sl (Lyapunov equation).} Let $H$ be a complex Hilbert space, known as the state space,  and ${\cal L}(H)$ the space of
bounded linear operators on $H$ with the usual operator norm. Let $(e^{-tA})_{t\geq 0}$ be a strongly continuous $(C_0)$ 
semigroup of bounded linear operators on $H$ such
that  $\Vert e^{-tA}\Vert_{ {\cal L}(H)}\leq M$ for all $t\geq 0$ and some $M<\infty$.
Let ${\cal D}(A)$ be the domain of the generator $-A$ so that ${\cal
D}(A)$ is itself a Hilbert space for the graph norm $\Vert \xi
\Vert_{{\cal D}(A)}^2=\Vert\xi\Vert^2_H+\Vert A\xi\Vert_H^2$, and let
$A^\dagger$ be the adjoint of $A$. Let $R:(0, \infty )\rightarrow {\cal L}(H)$ be a
differentiable function. The Lyapunov equation is
$$-{{dR_z}\over{dz}}=AR_z+R_zA\qquad  (z>0),\eqno(1.3)$$
\noindent where the right-hand side is to be interpreted as a bounded bilinear form on ${\cal D}(A)\times {\cal D}(A^\dagger )$.\par
\indent (ii) {\sl  (Operator ideals).} Let ${\cal L}^2(H)$ be the space of Hilbert--Schmidt operators on $H$, and ${\cal L}^1(H)$ be the space  of trace class operators on $H$, so ${\cal L}^1(H)=\{ T: T=VW; V,W\in {\cal L}^2(H)\}$and let $\det$ be the Fredholm determinant defined on 
$\{ I+T: T\in {\cal L}^1(H)\}$; see [67].\par
\indent (iii) {\sl (Tau function).} Suppose further that $R_x\in {\cal L}^1(H)$ for all $x>x_0$ for some $x_0\in {\bf R}$. Then the tau
function is $\tau (x)=\det (I+R_x)$ for $x>x_0$.\par
\vskip.05in
\indent The significant applications of this equation arises for linear systems.\par
\vskip.05in
\noindent {\bf Definition} (i) {\sl  (Linear system).}  Let $H_0$
be a complex separable Hilbert space which serves as the input and
output spaces; let $B:H_0\rightarrow H$ 
and $C:H\rightarrow H_0$ be bounded linear operators. The continuous-time linear system $(-A, B,C)$ is 
$$\eqalignno{ {{dX}\over{dt}}&=-AX+BU\cr 
                   Y&=CX,\cr
                  X(0)&=0.&(1.4)\cr}$$
\indent (ii) {\sl(Scattering function).} The scattering function is $\phi (x)=Ce^{-xA}B$, which is a bounded and weakly continuous  function $\phi: (0, \infty )\rightarrow {\cal L}(H_0)$.  In control theory, the transfer function is the Laplace transform of $\phi$.\par
\indent (iii) {\sl (Hankel operator).} Suppose that
 $\phi \in L^2((0, \infty ); {\cal L}(H_0))$. Then  the corresponding Hankel operator is $\Gamma_\phi$ on $L^2((0,
\infty ); H_0)$, where $\Gamma_\phi f(x)=\int_0^\infty \phi (x+y)f(y)\, dy$; see [58, 60] for boundedness criteria.\par
\vskip.05in
\noindent {\bf Definition}  {\sl (Admissible linear system).} Let $(-A,B,C)$ be a linear system as
above; suppose furthermore that the observability operator $\Theta_0:L^2((0,\infty );H_0)\rightarrow H$ is
bounded, where 
$$\Theta_0f=\int_0^\infty e^{-sA^\dagger}C^\dagger f(s)\,
ds;\eqno(1.5)$$
\noindent suppose that the controllability operator  $\Xi_0:L^2((0,\infty );H_0)\rightarrow H$ is also 
bounded, where 
$$\Xi_0f=\int_0^\infty e^{-sA}B f(s)\, ds.\eqno(1.6)$$
\indent (i) Then $(-A,B,C)$ is an admissible linear system.\par
\indent (ii)  Suppose furthermore that $\Theta_0$ and $\Xi_0$ belong
to the ideal ${\cal L}^2$ of Hilbert--Schmidt operators. Then
we say that $(-A,B,C)$ is $(2,2)$-admissible.\par 
\vskip.05in
\indent The scattering map associates to any $(2,2)$ admissible linear system  $(-A,B,C)$ the corresponding scattering function $\phi (x)=Ce^{-xA}B$. In Corollary 2.3, we show that for a suitable addition and multiplication on the linear systems, this map is additive and multiplicative.\par
\indent The inverse scattering problem involves recovering data about $u$ from $\phi$, as in (1.2).  In section 2 of this paper, we analyze the existence and uniqueness problem for the Lyapunov equation, and 
show that for any
$(2,2)$ admissible linear system, the
operator
$$R_x=\int_x^\infty e^{-tA}BCe^{-tA}\, dt\eqno(1.7)$$
\noindent is trace class and gives the unique solution to (1.3)  with the initial condition
$$\Bigl({{dR_x}\over{dx}}\Bigr)_{x=0}=-AR_0-R_0A=-BC.\eqno(1.8)$$
\noindent Also,  $R_x\in {\cal L}^1(H)$ and the Fredholm determinant satisfies
$$\det (I+\lambda R_x)=\det (I+\lambda\Gamma_{\phi_{ (x)}})\qquad  (x>0, \lambda\in {\bf C}).\eqno(1.9)$$
\vskip.05in 
\noindent {\bf Definition} {\sl (Tau function).} Given an $(2,2)$ admissible linear system $(-A,B,C)$, we define 
$$\tau (x)=\det (I+R_x).\eqno(1.10)$$
\vskip.05in
\indent Using this general definition of $\tau$, we can unify several results from the scattering theory of ordinary
differential equations. 
Such tau functions are strongly analogous to the tau functions introduced by Jimbo, Miwa and Ueno [55, 33] to describe the isomonodromy 
of rational differential equations. \par
\indent The  Gelfand--Levitan--Marchenko equation provides the  linkage between $\phi$ and $u$ via $R_x$. Consider
$$T(x,y)+\Phi (x+y)+\mu \int_x^\infty T(x,z)\Phi (z+y)\, dz=0\qquad
(0<x<y)\eqno(1.11)$$
\noindent where $T(x,y)$ and $\Phi (x+y)$ are $m\times m$ 
 matrices with scalar entries. In the context of $(-A, B,C)$ we assume that $\Phi (x)=Ce^{-xA}B$ is known and aim to find $T(x,y)$. In section two, we use $R_x$ to construct solutions to the associated
Gelfand--Levitan equation (1.11), and introduce a potential
$$u(x)=-2{{d^2}\over{dx^2}}\log\det (I+R_x).\eqno(1.12)$$
\noindent Then we obtain a differential equation
linking $\phi (x)$ to $u(x)$. In examples of interest in scattering theory, one can
calculate $\det (I+\lambda R_x)$ more easily than the Hankel
determinant of  $\Gamma_{\phi_{(x)}}$  directly [20, 22, 47, 61], since  $R_x$ has additional properties that originate 
from Lyapunov's equation. In section two, we establish properties of $R_x$ and $\tau$ for $(2,2)$ admissible linear systems, and use $R_x$ to solve the Gelfand--Levitan equation.\par
\vskip.05in

\noindent {\bf Definition}  {\sl (Sectorial operator).} For $0<\theta \leq\pi$, we introduce the sector $S_\theta=\{ z\in {\bf C}\setminus \{0\}: \vert \arg z\vert <\theta\}$. A closed and densely defined linear operator $-A$ is sectorial [21] if there exists $\pi/2<\theta<\pi$ such that $S_\theta$ is contained in the resolvent set of $-A$ and $\vert\lambda \vert\Vert (\lambda I+A)^{-1}\Vert_{{\cal L}(H)}\leq M$ for all $\lambda\in S_\theta$. Let ${\cal D}(A)$ be the domain of $A$ and ${\cal D}(A^\infty )=\cap_{n=0}^\infty {\cal D}(A^n)$.\par
 
\vskip.05in
\noindent {\bf Definition}  {\sl (Deformations).} There are three basic deformations of a $(2,2)$-admissible linear system $\Sigma =(-A,B,C)$ with $H_0={\bf C}$ and $-A$ a sectorial operator, with corresponding effects on tau functions. \par
\indent (i)  Translation takes $\Sigma\mapsto\Sigma (a)$, where $\Sigma (a)= (-A,Be^{-aA}, Ce^{-aA})$ and $\tau (x)\mapsto \tau (x+a)$ 
for $a\in\Omega$. The translation operation is accounted for in the properties of Hankel operators, and is used in section two. \par
\indent (ii) The Miura transform is the involution $(-A,B,C)\mapsto (-A,B,-C)$. In section three, we show that the tau functions $\tau_\infty$ for $(-A,B,C)$ and $\tau_0$ for $(-A,B,-C)$ satisfy
$$\tau_0''(x)\tau_\infty (x)-2\tau_0'(x)\tau_\infty'(x)+\tau_0(x)\tau_\infty''(x)=0.\eqno(1.13)$$
\indent (iii)  The Darboux addition (spectral shift) [47] maps $\Sigma  \mapsto \Sigma_\zeta $ where $$\Sigma_\zeta =(-A, (\zeta I+A)(\zeta I-A)^{-1}B,C)\qquad (-\zeta\in S_\theta )\eqno(1.14)$$
\noindent has tau function $\tau_\zeta$. In section three, we use this operation to construct families of solutions to Schr\"odinger's equation.\par 
$$-\psi''_\zeta (x) +u(x)\psi_\zeta (x) =-\zeta^2\psi_\zeta (x)\qquad (x>0),\eqno(1.15)$$
\noindent where $u$ is typically complex-valued. When interpreting results, it is important to realize that there is no simple connection between the spectrum of $A$ on $H$ and the spectrum of Schr\"odinger's equation on $L^2((0, \infty ); {\bf C})$.\par 
\vskip.05in
 \noindent {\bf Definition} {\sl (State ring).}  Let $\Omega$ be a domain in ${\bf C}$ such that $z,w\in\Omega$ implies $z+w\in\Omega$. Suppose that 
$I+R_x$ is invertible for all $x\in \Omega$ and so $F_x=(I+R_x)^{-1}\in {\cal L}(H)$, and $F_x-I\in {\cal L}^1(H)$.  Suppose momentarily that $A\in {\cal L}(H)$. We introduce  an algebra ${\cal A}_\Sigma$  of holomorphic functions from 
$\Omega$ to
${\cal L}(H)$, which  is generated by $I, A\in {\cal L}(H)$ and $F_x$ $(x\in \Omega )$, and 
is a differential ring for the usual pointwise multiplication and differentiation $d/dx$ over $\Omega$ as in [63].  We define a new multiplication 
$$P\ast Q=P(AF+FA-2FAF)Q\eqno(1.16)$$
 on ${\cal A}_\Sigma$ , and a new differentiation
$\partial :{\cal A}_\Sigma\rightarrow{\cal A}_\Sigma$ by 
$$\partial P=A(I-2F)P+{{dP}\over{dx}}+P(I-2F)A\eqno(1.17)$$
\noindent so that $({\cal 
A}, \ast , \partial )$ is a differential ring.  Let ${\cal M}_{\Omega}$ be the meromorphic
complex functions on $\Omega$, with the usual pointwise operations. Then we define the bracket
$$\lfloor P\rfloor
=Ce^{-xA}F_xPF_xe^{-xA}B,\eqno(1.18)$$
\noindent so that
 $\lfloor \, .\, \rfloor :({\cal A}_\Sigma, \ast , \partial )\rightarrow ({\cal M}_\Omega , \cdot , d/dx)$ is a
differential ring homomorphism, so $\lfloor P\ast Q\rfloor =\lfloor P\rfloor\lfloor Q\rfloor$
 and $(d/dx)\lfloor P\rfloor =\lfloor \partial P\rfloor .$ \par
\indent This algebra $({\cal A}_\Sigma, \ast , \partial )$ is generally non commutative, is realized as an algebra of operators 
on the space space $H$ of $(-A,B,C)$, and provides a substitute for the multiplication structure that is lacking in the theory of Hankel operators. 
Our main result is as follows.\par

\vskip.05in
\noindent {\bf Theorem 1.1.} {\sl Suppose that $\Sigma=(-A,B,C)$ is a
$(2,2)$ admissible linear system with $-A$ sectorial for $S_\theta$ and $H_0={\bf C}$.\par
\indent (i) Then there
exist $x_0>0$ and a solution $R_x$ to (1.3) and (1.8) such that $\tau (x)=\det (I+R_x)$
is holomorphic and $u(x)$ is meromorphic for $x\in \Omega$ where $\Omega =\{x_0+z: z\in S_{\theta -\pi/2}\}$.}\par
\indent {\sl (ii) The tau functions for $\Sigma$ and $\Sigma_\zeta $ give $\psi_\zeta (x)=e^{\zeta x}\tau_\zeta (x)/\tau (x)$  that satisfies (1.15) for $-\zeta\in S_\theta$.} \par
\indent {\sl (iii) There exists a differential algebra $({\cal A}_\Sigma, \ast ,\partial )$ on $\Omega$ which contains
$F_x=(I+R_x)^{-1}$ and there exists a differential ring homomorphism $\lfloor \, .\, \rfloor :{\cal A}_\Sigma \rightarrow {\cal M}_{\Omega}$ such that $u=-4\lfloor A\rfloor$.}\par 
\indent {\sl (iv) If $A$ satisfies $\lfloor p(A)\rfloor=0$ for some non-zero odd complex polynomial $p$, then\par
\noindent  ${\bf C}[u, {\partial u}/\partial x, \dots,  ]$ is a Noetherian differential ring for $ {\partial }/\partial x$ and the standard multiplication.}\par
\indent {\sl (v) If $\lfloor A^{2m-1}\rfloor =0$ for some $m$, then  (1.15) can be integrated by quadratures.}\par
\vskip.05in
\indent In sections 4 and 5 we show that $\lfloor \, \cdot \,\rfloor$ maps $(-4\partial^jA)_{j=0}^\infty$ to $(u^{(j)})_{j=0}^\infty$, and $((-1)^j2A^{2j-1})_{j=1}^\infty$ to $(f_j)_{j=1}^\infty$, where $(f_j)$ satisfies the stationary $KdV$ hierarchy. This provides the crucial link between algebraic properties of ${\cal A}_\Sigma$ for $\Sigma= (-A,B,C)$ and the spectral  Schr\"odinger equation. \par
\indent We recall that a compact Riemann surface ${\cal E}$ is 
hyperelliptic if and only if there exists a meromorphic function on ${\cal E}$ that has precisely two poles. In this case, there is a
two-sheeted cover ${\cal E}\rightarrow {\bf P}^1$ with $2g+2$ branch
points, where $g$ is the genus of ${\cal E}$. 
The elliptic case has $g=1$. In case (v), $u$ is associated a meromorphic function on a hyperelliptic curve, and the potential is said to be finite-gap or algebro-geometric; see [14, 29].\par
\indent Theorem 1.1 applies to scattering potentials such that $u(x)\rightarrow 0$ as $x\rightarrow\infty$. The case of periodic potentials is considered in sections 6 and 7 of this paper.
\vskip.05in
\noindent {\bf Definition} {\sl (Periodic linear system).}  Let $A\in {\cal L}(H)$ be such that $(e^{xA})_{x\in {\bf R}}$ is a periodic group, and let $B\in {\cal L}(H_0, H)$ and $C\in {\cal L}(H,H_0)$ satisfy $AE+EA=BC$ for some $E\in {\cal L}^1(H)$. Then $\Sigma =(-A,B,C;E)$ is a uniformly periodic linear  system, and $\tau (x)=\det (I+e^{xA}Ee^{xA})$ is the corresponding tau function.\par
\vskip.05in 
\indent In section 6 of this paper, we introduce a ring $({\cal A}_\Sigma, \ast , \partial )$ for periodic linear systems. The periodic linear system $\Sigma$ has a tau function $\tau$ and a periodic potential $u$, as in Hill's equation $-f''+uf=\lambda f$. Hence $\Sigma$ is associated with Hill's discriminant $\Delta (\lambda )$ and a spectral curve ${\cal E}$, which is typically a transcendental hyperelliptic curve of infinite genus. The Jacobi variety ${\bf X}$ of ${\cal E}$ is then an infinite dimensional complex torus.\par
\indent We show that if $u$ is an elliptic function, then there exist uniformly periodic linear systems with tau functions $\tau_1$ and $\tau_0$ such that $u(x)=\tau_1(x)/\tau_0(x)$. Also, we show that if Hill's equation for this $u$ has general solution $f(x;\lambda )$ that is meromorphic in $x$ for all but finitely many $\lambda\in {\bf C}$, then there exist uniformly periodic linear systems with tau function $\tau_3(x;\lambda )$ and $\tau_4(x;\lambda )$ such that $f(x;\lambda )=\tau_3(x;\lambda )/\tau_4(x;\lambda )$.\par
\indent Using results of Gesztesy and Weikard [31], we prove a partial converse, namely that if Hill's equation has a general solution of the form  $f(x;\lambda )=\tau_3(x;\lambda )/\tau_4(x;\lambda )$ for all $\lambda\in {\bf C}$, where
 $\tau_3(x;\lambda )$ and $\tau_4(x;\lambda )$ are the tau functions of uniformly periodic linear systems, then $u$ is an algebro geometric potential and is associated with a hyperelliptic spectral curve ${\cal E}$ of finite genus.\par
\indent  The notion of a tau function of a linear system generalizes the classical concept of a theta function for an algebraic curve. We recall that a complete complex algebraic curves  is associated with 
a finite dimensional  Jacobian variety which has dimension determined by the genus of the curve.  Riemann's theta functions may be defined on such a variety, and they satisfy some addition rules known as Fay's identities, which reflect the geometry of the underlying curve [23].  Ercolani and McKean [22] observed that Fay's identities can be deduced from secant identities which relate to the properties of Wronskians of suitably chosen functions. In this spirit, we prove that tau functions of linear systems satisfy some Wronskian identities which are counterparts of Fay's identities. Mumford [56] observed that Fay's identities give rise to nontrivial solutions of certain partial differential equations.\par 
\indent The Kadomtsev--Petviashvili equations describe the waves in a two
dimensional dissipative medium where the scale of the propagation of the
wave along the $y$-axis is much larger than the longitudinal scale along
the $x$-axis. We write
$$KP\qquad {{\partial}\over{\partial x}}\Bigl( {{\partial^3 u}\over{\partial
x^3}}-6u{{\partial u}\over{\partial x}}+4\lambda {{\partial u}\over{\partial x}}
+4\alpha {{\partial u}\over{\partial s}}\Bigr)
 +3\beta^2{{\partial^2u}\over{\partial y^2}}=0,\eqno(1.19)$$
\noindent where the $\alpha, \beta, \lambda\in {\bf C}$ are parameters.  Krichever showed that any tau function that
arises from the theta function on the Jacobi variety of a complete
complex algebraic curve satisfies $KP$. Shiota [66] and Mulase [56] proved the converse, that if $\tau$ is the theta
function of a finite-dimensional Abelian manifold, and $\tau$ gives a
potential that satisfies $KP$, then the Abelian manifold arises as a
Jacobi variety of a complex algebraic curve. \par
\indent In section 9, we show that the $\tau$ functions of general linear systems under a group of deformations  give rise to a set of solutions of the $KP$ equations.
\vskip.05in
\noindent {\bf Theorem 1.2.} {\sl  Let $(-A_1,B_0,C_0)$ and $(-A_2,B_0,C_0)$ be $(2,2)$ admissible linear 
systems with input and output spaces ${\bf C}$, where $A_1,A_2\in {\cal L}(H)$. Let $$C(y;t)=C_0e^{t(A_1^3+\lambda A_1)/\alpha -yA_1^2/\beta},\qquad 
 B(y;t)=e^{t(A_2^3+\lambda A_2)/\alpha +yA_2^2/\beta}B_0;$$
\noindent  then with
 $$R_x(y,t)=\int_x^\infty e^{-A_2s}B(y;t)C(y;t)e^{-A_1 s}ds,\eqno(1.20)$$
\noindent let
$$u(x,y,t)=-2{{\partial^2}\over{\partial x^2}}\log\det (I+R_x(y,t)).\eqno(1.21)$$
\indent Then $u$ satisfies $KP$ in the form (1.19).}\par
\vskip.05in
\indent  The $KP$ differential equation is the first in a sequence of nonlinear partial differential equations known as the $KP$ hierarchy. Shiota [66] proved that these are related by an integral involving a family of tau functions which are subject to a group of deformations involving infinitely many parameters.  We introduce a family of linear systems which are subject to a group of deformations, and also show that such a family gives tau functions which make Shiota's integral vanish; this condition is known to give solutions of the $KP$ hierarchy.
\vskip.05in

\indent 
\vskip.05in

\vskip.05in

\noindent {\bf 2. $\tau$ functions in terms of Lyapunov's equation and the Gelfand--Levitan equation}\par
\vskip.05in
\noindent The following proves uniqueness of solutions of the Lyapunov equation (1.3), in a style suggested by
[60, p 503]. Peller [60] discusses scattering functions that produce bounded self-adjoint Hankel operators $\Gamma_\phi$, and their realization in terms of continuous time linear systems. He observes that in some cases one needs a bounded semigroup with unbounded generator $(-A)$.  We prove the uniqueness results for bounded and strongly continuous semigroups, then specialize to holomorphic semigroups. The main application is to the Gelfand--Levitan equation,  and associated determinants. \par
\vskip.05in
\noindent {\bf Proposition 2.1.} {\sl Let $(e^{-tA})_{t\geq 0}$ be a strongly continuous and weakly 
asymptotically stable semigroup on a complex Hilbert space $H$, so $e^{-tA}f\rightarrow 0$ weakly as $t\rightarrow \infty$ for all
$f\in H$. Then\par
\indent (i) $S_t:R\mapsto e^{-tA}Re^{-tA}$ for $t\geq 0$ defines a strongly continuous semigroup 
on ${\cal L}^1(H)$, which has generator 
$(-L)$, with dense domain of definition ${\cal D}(L)$ such that 
$$L(R)=AR+RA\qquad (R\in {\cal D}).\eqno(2.1)$$
\indent (ii) The linear operator $L:{\cal D}(L)\rightarrow {\cal L}^1(H)$ is injective, and each $R_0\in {\cal D}(L)$ 
with $L(R_0)=X$, there exists a weakly convergent improper integral} 
$$R_0=\int_0^\infty e^{-tA}Xe^{-tA}\, dt.\eqno(2.2)$$
\indent {\sl (iii) Suppose moreover that $\Vert e^{-t_0A}\Vert_{{\cal L}(H)}<1$ for some $t_0>0$. Then $L:{\cal D}(L)\rightarrow {\cal L}^1(H)$
is surjective, the integral (2.2) converges absolutely in ${\cal L}^1(H)$ and $R_0$ gives
the unique solution to $AR_0+R_0A=X.$}\par 
\vskip.05in
\noindent {\bf Proof.} (i) First observe that by the uniform boundedness theorem, there exists $M$ such that
 $\Vert e^{-tA}\Vert_{{\cal L}(H)}\leq M$ for all $t\geq 0$, so $(e^{-tA})_{t\geq 0}$ is uniformly bounded. Also, the adjoint semigroup  $(e^{-tA^\dagger})_{t\geq 0}$ is
also strongly continuous and uniformly bounded, so $A$ and $A^\dagger$ have dense
domains ${\cal D}(A)$ and ${\cal D}(A^\dagger)$ in $H$.\par
\indent Now ${\cal L}^1(H)=H\hat\otimes H$, the projective tensor product, so for all $X\in {\cal L}^1(H)$, there exists a nuclear decomposition 
$X=\sum_{j=1}^\infty B_jC_j$ where $B_j, C_j\in H$ satisfy\par
\noindent  $\Vert X\Vert_{{\cal L}^1(H)}=\sum_{j=1}^\infty
 \Vert B_j\Vert_H\Vert C_j\Vert_H$. Then 
$$S_t(X)-X =\sum_{j=1}^\infty (e^{-tA}B_jC_je^{-tA}-B_jC_je^{-tA})
+\sum_{j=1}^\infty (B_jC_je^{-tA}-B_jC_j)\eqno(2.3)$$
\noindent where $(e^{-tA})$ is bounded,  $\Vert e^{-tA}B_j-B_j\Vert_H\rightarrow 0$ 
and $\Vert e^{-tA^\dagger }C_j-C_j\Vert_H\rightarrow 0$ as $t\rightarrow0+$, 
so $\Vert S_t(X)-X\Vert_{{\cal L}^1(H)}\rightarrow 0$ as $t\rightarrow 0+$, 
so $(S_t)_{t\geq 0}$ is strongly continuous on ${\cal L}^1(H)$. By semigroup theory, 
there exists a dense linear subspace  ${\cal D}(L)$ of ${\cal L}^1(H)$ such that 
$S_t(R)$ is differentiable at $t=0+$ for all $R\in {\cal D}$, and $(d/dt)_{t=0+} S_t(R)=-AR-RA$, so the generator
is $(-L)$, where $L(R)=AR+RA$.
\par
\indent (ii) Certainly ${\cal D}$ contains ${\cal D}(A^\dagger )\hat\otimes {\cal D}(A)$ in
${\cal L}^1(H)=H\hat\otimes H$. Choosing $f\in {\cal D}(A)$ and $g\in
{\cal D}(A^\dagger)$, we find that 
$$\eqalignno{{{d}\over{dt}}\bigl\langle e^{-tA}R_0e^{-tA}f,g\bigr\rangle&=-  
\bigl\langle e^{-tA}(AR_0+R_0A)e^{-tA}f,g\bigr\rangle\cr
&=-\bigl\langle e^{-tA}Xe^{-tA}f,g\bigr\rangle&(2.4)\cr}$$
\noindent a continuous function of $t>0$, so integrating we obtain
$$\bigl\langle R_0f,g\bigr\rangle -\bigl\langle e^{-sA}R_0e^{-sA}f,g\bigr\rangle
=\int_0^s \bigl\langle e^{-tA}Xe^{-tA}f,g\bigr\rangle \, dt.\eqno(2.5)$$
\noindent We extend this identity to all $f,g\in H$ by joint continuity; then we
let $s\rightarrow\infty$ and observe that $R_0:H\rightarrow H$ is trace class and hence is completely continuous, hence $R_0$ maps the weakly null
family $(e^{-sA}f)_{s\rightarrow\infty }$ to the norm convergent family $(R_0e^{-sA}f)_{s\rightarrow\infty}$, so
$\langle e^{-sA}R_0e^{-sA}f,g\rangle\rightarrow 0$ as $s\rightarrow\infty$, hence we have a weakly convergent improper integral
$$\bigl\langle R_0f,g\bigr\rangle=
\lim_{s\rightarrow\infty}\int_0^s \bigl\langle e^{-tA}Xe^{-tA}f,g\bigr\rangle 
\, dt\qquad (f,g\in H).\eqno(2.6)$$

\indent (iii) The function $t\mapsto e^{-tA}Xe^{-tA}$ takes values in the separable space ${\cal L}^1(H)$  and is weakly continuous,
hence strongly measurable, by Pettis's theorem. By considering the 
spectral radius, Engel and Nagel [21] show that there exist $\delta >0$ and $M_\delta >0$ such that
$\Vert e^{-tA}\Vert_{{\cal L}(H)}\leq M_\delta e^{-\delta t}$ for all $t\geq 0$;
hence (2.2) converges as a Bochner--Lebesgue integral with 
$$\eqalignno{\Vert R_x\Vert_{{\cal L}^1(H)}&\leq \int_x^\infty 
M_\delta^2 \Vert X\Vert_{{\cal L}^1(H)}e^{-2\delta t}\, dt\cr
&\leq {{M_\delta^2}\over{2\delta}}\Vert X\Vert_{{\cal L}^1(H)}e^{-2\delta x}.&(2.7)\cr}$$
\noindent Furthermore, $A$ is a closed linear operator and satisfies
$$\eqalignno{A\int_x^s e^{-tA}Xe^{-tA} \, dt +\int_x^s e^{-tA}Xe^{-tA} \,
dtA&=\int_x^s -{{d}\over{dt}} \bigl(e^{-tA}Xe^{-tA}\bigr)\, dt\cr
&=e^{-xA}Xe^{-xA}-e^{-TA}Xe^{-TA}\cr
&\rightarrow e^{-xA}Xe^{-xA}&(2.8)}$$ 
\noindent as $s\rightarrow\infty$ where $\int_x^s e^{-tA}Xe^{-tA}dt\rightarrow R_x$; so $AR_x+R_xA=e^{-xA}Xe^{-xA}$ for
all $x\geq 0$. We deduce that $x\mapsto R_x$ is a differentiable
function from $(0, \infty )$ to ${\cal L}^1(H)$ and that the modified Lyapunov
equation (1.3) holds. \par 
\rightline{$\square$}
\indent The hypotheses (i) and (ii) are symmetrical under the adjoint $(A, R_0)\mapsto (A^\dagger, 
R_0^\dagger )$; however, the hypothesis (iii) is rather stringent, and will be replaced in examples by sharper
conditions.

\vskip.05in
\noindent We introduce Lyapunov's equation, and the existence of
solutions for suitable $(-A,B,C)$. The solution $R_x$ is defined by
a formula suggested by Heinz's theorem [8] and has properties
analogous to the resolvent operator of a semigroup.\par
\vskip.05in
\noindent {\bf Definition} {\sl ($(2,2)$  admissible linear systems).}  (i) Let $H$ be a complex Hilbert space and let $\Sigma=(-A,B,C)$ be a linear system with
state space $H$. Suppose that there is a weakly convergent integral 
$$W_c=\int_0^\infty e^{-tA}BB^\dagger e^{-tA^\dagger}\, dt\eqno(2.9)$$
\noindent which defines a bounded linear operator on $H$; then $W_c$ is the controllability Gramian. Suppose further 
that there exists a weakly convergent integral   
$$W_o=\int_0^\infty e^{-tA^\dagger}C^\dagger Ce^{-tA}\, dt\eqno(2.10)$$
\noindent which defines a bounded linear operator on $H$; then $W_o$ is the observability Gramian.\par  
\indent (ii) Then we define $R_x$ to be the bounded
linear operator on $H$ determined by the weakly convergent integral
$$R_x=\int_x^\infty e^{-tA}BCe^{-tA}\, dt.\eqno(2.11)$$
\indent (iii) Then $\Sigma$ satisfying (i) is said to be balanced if $W_c=W_o$ and 
${\hbox{ker}}(W_c)=0$.\par
\indent (iv) Also, $\Sigma$ satisfying (i) is said to be $(2,2)$ admissible if $W_c$ and
$W_o$ are trace class.\par
\indent (v) We introduce the scattering function $\phi (t)=Ce^{-tA}B$ and the shifted scattering function $\phi_{(x)}(t)=\phi (t+2x)$ for $x,t>0$.\par
\indent (vi) The tau function of $\Sigma$ is $\tau (x)=\det (I+R_x)$.\par
\indent (vii) For $0<\delta <\pi$, we introduce the sector $S_\delta =\{ z\in {\bf
C}\setminus \{0\}: \vert \arg  z\vert <\delta \}$. For $\pi/2<\delta <\pi$, we introduce $X_\delta =\{\zeta \in S_\delta : -\zeta \in S_\delta \}$ which is an open set, symmetrical about $i{\bf R}$ and bounded by lines passing through $0$. \par
\vskip.1in

\noindent {\bf Theorem 2.2.} {\sl Let $(-A,B,C)$ be a linear
system such that $\Vert e^{-t_0A}\Vert_{{\cal L}(H)}<1$
for some $t_0>0$, and that $B$ and $C$ are Hilbert--Schmidt
operators  such that $\Vert B\Vert_{{\cal L}^2(H_0;H)}\Vert C\Vert_{{\cal L}^2(H;H_0)}\leq 1$. Suppose further that $-A$ is sectorial on $S_\theta$ for some $\pi/2 <\theta<\pi$. \par
\indent (i) Then $(-A,B,C)$ is $(2,2)$-admissible, so the 
trace class operators $R_x$ give the solution to Lyapunov's equation (1.3) for $x>0$ that
satisfies the initial condition (1.8), and the solution to (1.3) with (1.8)  is unique.\par
\indent (ii) The function $\tau (x) =\det (I+R_x)$ is 
differentiable for $x\in (0, \infty )$.}\par 
\indent {\sl (iii) Then $R_z$ extends to a holomorphic function which satisfies 
(1.3) on $S_{\theta -\pi/2}$, and $R_z\rightarrow 0$ as $z\rightarrow\infty$ 
in $S_{\theta -\varepsilon-\pi/2}$ for all $0<\varepsilon <\theta -\pi/2$.}\par 

\vskip.05in
\noindent {\bf Proof.} (i) Since $BC\in {\cal L}^1(H)$, the integrand of (2.11) takes values in ${\cal L}^1(H)$, and 
we can apply Proposition 2.1(iii) to $X=BC$.\par
\indent (ii) The Fredholm determinant $R\mapsto \det (I+R)$ is a continuous function on ${\cal L}^1(H)$. Also the
integral $R_x=\int_x^\infty e^{-tA}BCe^{-tA}\, dt$ belongs to ${\cal D}(L)$ and gives a
differentiable function of $x>0$ with values in ${\cal L}^1(H)$.\par
\indent (iii) By classical results of Hille, $(e^{-zA})_{z\in S_{\theta -\pi/2}}$ defines an analytic semigroup 
on $S_{\theta -\pi/2}$, bounded on $S_{\nu}$ for all $0<\nu<\theta -\pi/2$, so we can define $R_z=e^{-zA}R_0e^{-zA}$ and obtain an 
analytic solution to Lyapunov's equation. For all $0<\varepsilon <\theta -\pi /2$, 
there exists $M_{\varepsilon}'$ such that 
$\Vert e^{-zA}\Vert_{{\cal L}(H)}\leq M_\varepsilon'$ for all 
$z\in S_{\delta}$ where  $\delta =\theta-\varepsilon -\pi/2$. Now for $z\in
S_{\delta/2}$, we write $z=x/2+(x/2+iy)$ with $x/2+iy\in S_\delta$ and use the
bound $\Vert e^{-zA}\Vert_{{\cal L}(H)}\leq \Vert e^{-xA/2}\Vert_{{\cal L}(H)}\Vert e^{-(x/2+iy)A}\Vert_{{\cal L}(H)}$ to obtain 
$\Vert e^{-zA}\Vert_{{\cal L}(H)}\leq  M_{\varepsilon}'^2\Vert e^{-t_0A}\Vert_{{\cal L}(H)}^{x/(4t_0)}$, so $\Vert
e^{-zA}\Vert_{{\cal L}(H)}\rightarrow 0$ exponentially fast as $z\rightarrow\infty$ in the sector 
$S_{\delta/2}$. Hence $R_z$ is holomorphic and bounded on 
$S_{(\theta -\varepsilon -\pi/2)}$  and 
by (2.7), $R_z\rightarrow 0$ as $z\rightarrow\infty$ in 
$S_{(\theta -\varepsilon -\pi/2)/2}$.\par
\rightline{$\square$}
\vskip.05in
\noindent {\bf Definition} Given  $M, \varepsilon
>0$ and $0<\delta <\pi /2$, let ${\cal R}(\delta , M, \varepsilon )$ be the space of all the holomorphic
functions $\psi :S_\delta \rightarrow {\bf C}$ such that $e^{\varepsilon \vert
z\vert}\vert \psi (z)\vert \leq M$ for all $z\in S_\delta$. Then let ${\cal R}=\cup_{0<M, 0<\varepsilon, 0<\delta<\pi }{\cal R}(\delta, M, \varepsilon)$. \par
\indent  We note that by Cauchy's estimates,
${\cal R}$ has the following properties:\par
\indent (i) $\phi'\in {\cal R}$ for all $\phi\in {\cal R}$;\par
\indent (ii) for all $\phi\in {\cal R}$, there exists $\psi\in {\cal R}$ such that
$\psi'=\phi$;\par
\indent (iii) if $\phi \in {\cal R}$, then $e^\phi -1\in {\cal R}$;\par
\indent (iv) ${\cal R}+{\bf C}$ is an integral domain under pointwise addition and
multiplication;\par
\indent (v) ${\cal R}+{\bf C}$ is closed under taking of
exponentials.\par
\noindent Hence  we can therefore form the field of fractions of ${\cal R}+{\bf C}$ to obtain a
differential field ${\cal F}$, so that every elements of ${\cal F}$ is meromorphic on some sector $S_\eta$.\par
\vskip.05in
\noindent {\bf Corollary 2.3.} {\sl The set of scattering functions $\phi$ which arise from the $(2,2)$ admissible linear systems as in Theorem 2.2 gives a subring of ${\cal R}$.}\par
\vskip.05in 
\noindent {\bf Proof.}  Let $(-A_j, B_j, C_j)$ be $(2,2)$ admissible linear system with scattering function $\phi_j$ as in Theorem 2.2 for $j=1,2$. Then by Theorem 2.2(iii) the
scattering functions satisfy $\phi_j\in {\cal R}$.\par
\indent The linear system 
$$\Bigl(\left[\matrix{-A_1&0\cr 0&-A_2\cr}\right], \left[\matrix{B_1\cr B_2\cr}\right], \left[\matrix{C_1&C_2}\right]\Bigr).\eqno(2.12)$$
\noindent is $(2,2)$-admissible with scattering function $\phi_1(x)+\phi_2(x)$, as one checks by direct calculation.\par
\indent  The linear system
$$\bigl( -(A_1\otimes I+I\otimes A_2), B_1\otimes B_2, C_1\otimes C_2\bigr).\eqno(2.13)$$
\noindent is also $(2,2)$  admissible with scattering function $\phi_1(x)\phi_2(x)$.  One checks that the semigroups $(e^{-tA_1}\otimes I)_{t>0}$ and likewise $(I\otimes e^{-tA_2})_{t>0}$ are strongly continuous on the tensor product Hilbert space $H\otimes H$, then $(e^{-tA_1}\otimes e^{-tA_2})_{t>0}$ is strongly continuous on $H\otimes H$ by B15 of {Engel]. Also, $B_1\otimes B_2$ and $C_1\otimes C_2$ are Hilbert--Schmidt, as one checks by considering orthonormal bases.\par
\rightline{$\square$}
\vskip.05in
 \noindent {\bf Example.} Let $\Delta= -d^2/dx^2$ be the usual Laplace operator which is essentially self-adjoint and non-negative on $C_c^\infty ({\bf R}; {\bf C})$  in $L^2({\bf R}; {\bf C})$. We introduce $A=\sqrt{I+\Delta}$ which is given by the Fourier multiplier ${\cal F} Af(\xi )=\sqrt {1+\xi^2}{\cal F}f(\xi)$. Then $(e^{-zA})$ and $(e^{-zA^2})$  give bounded holomorphic semigroups on $H$, as in Theorem 2.2, 
on the right half-plane $\{z\in {\bf C}: \Re z\geq 0\}$, which is the closure of $S_{\pi/2}$. On the imaginary axis, we have unitary groups $(e^{itA})$ and $(e^{-itA^2})$. By classical results from wave equations, we can write $e^{itA}+e^{-itA}=2\cos (tA)$ 
where $u(x,t)=\cos (tA)f(x)$ is given by
$$u(x,t)={{1}\over{2}}\bigl( f(x+t)+f(x-t)\bigr) +{{t}\over{2}}\int_{x-t}^{x+t} f(y){{J'_0(\sqrt {t^2-(x-s)^2})}\over {\sqrt{t^2-(x-s)^2}}}ds\qquad (f\in C_c^\infty ({\bf R}; {\bf C})),.\eqno(2.14)$$
\noindent with $J_0$ is Bessel's function of the first kind of order zero, and $u$ satisfies 
$$\eqalignno{{{\partial^2 u}\over{\partial x^2}}-{{\partial^2u}\over{\partial t^2}}&=u(x,t)\cr
u(x,0)&=f(x);\cr
 {{\partial u}\over{\partial t}}(x,0)&=0.&(2.15)\cr}$$ 
\noindent Note that $(\exp (t(iA)^{2j-1}))$ gives a unitary group on $H$ for $j=0,1,2, \dots $.\par  

\vskip.05in
\noindent {\bf Definition} (i) {\sl  (Block Hankel operators).} Say that  $\Gamma\in {\cal L}(H)$ is block Hankel if there exists $1\le m<\infty$ such that $\Gamma$ is unitarily equivalent to the block matrix $[A_{j+k-2}]_{j,k=1}^\infty$ on $\ell^2({\bf C}^m)$ where $A_{j}$ is a $m\times m$ complex matrix for $j=0, 1, \dots $.\par
 \indent (ii) Let $(-A,B,C)$ be a $(2,2)$ admissible linear system with input and output space
 $H_0$, where the dimension of $H_0$ over ${\bf C}$ is $m<\infty$. Then $m$ is the 
number of outputs of the system, and systems with finite $m>1$ are known as MIMO for 
multiple input, multiple output, and give rise to block Hankel 
operators with $\Phi (x)=Ce^{-xA}B$; see [60]. \par
\indent (iii) The
Gelfand--Levitan integral equation for $(-A,B,C)$ as in (ii) is
$$T(x,y)+\Phi (x+y)+\mu \int_x^\infty T(x,z)\Phi (z+y)\, dz=0\qquad
(0<x<y)\eqno(2.16)$$
\noindent where $T(x,y)$ and $\Phi (x+y)$ are $m\times m$ 
 matrices with scalar entries, and $\mu\in {\bf C}$.\par
\vskip.05in
\noindent {\bf Proposition 2.4.} {\sl (i) In the notation of Theorem 2.2, there exists $x_0>0$ such that 
$$T_\mu (x,y)=-Ce^{-xA}(I+\mu R_x)^{-1}e^{-yA}B\eqno(2.17)$$
\noindent satisfies the 
integral equation (2.16) for $x_0<x<y$ and $\vert \mu\vert <1.$\par 
\indent (ii) The determinant satisfies $\det (I+\mu R_x)=\det
(I+\mu\Gamma_{\Phi_{(x)}})$ and} 
$$\mu {\hbox{trace}}\, T_\mu (x,x)={{d}\over{dx}}\log\det (I+\mu
R_x).\eqno(2.18)$$ 
\indent {\sl (iii) Suppose that $t\mapsto U(t)$ is a continuous
function $[0,1]\rightarrow {\cal{L}}(H)$ such that $U(t)A=AU(t)$ and $\Vert U(t)\Vert_{{\cal L}(H)}\leq 1$. Then there is a family of $(2,2)$ admissible linear systems
$$\Sigma (t)=(-A,U(t)B, CU(t) )\qquad (t\in [0,1]);\eqno(2.19)$$
\noindent the corresponding tau
function $\tau (x,t)$ is continuous for $(x,t)\in (0,\infty )\times
[0,1]$.}\par 
\vskip.05in
\noindent {\bf Proof.} (i) We choose $x_0$ so large that $e^{\delta x_0}\geq 
M_\delta /2\delta$, then by (2.7), we have $\vert \mu\vert\Vert R_x\Vert_{{\cal L}(H)}
<1$ for $x>x_0$, so $I+\mu R_x$ is invertible. Substituting into the
integral equation, we obtain  
$$\eqalignno{ Ce^{-(x+y)A}B&-Ce^{-xA}(I+\mu
R_x)^{-1}e^{-yA}B\cr
&-\mu Ce^{-xA}(I+\mu R_x)^{-1}\int_x^\infty
e^{-zA}BCe^{-zA}\, dze^{-yA}B\cr
&=Ce^{-(x+y)A}B-Ce^{-xA}(I+\mu R_x)^{-1}e^{-yA}B-\mu
Ce^{-xA}(I+\mu R_x)^{-1}R_xe^{-yA}B\cr
&=0.&(2.20)\cr}$$ 
\indent (ii) As in (1.5), the operator $\Theta_x:L^2(0, \infty )
\rightarrow H$
 is Hilbert--Schmidt; likewise $\Xi_x :L^2(0, \infty )\rightarrow H$ 
is Hilbert--Schmidt; so $(-A,B,C)$ is $(2,2)$-admissible. Hence 
$\Gamma_{\Phi_{(x)}}=\Theta_x^\dagger\Xi_x$ and 
$R_x=\Xi_x\Theta_x^\dagger$ are trace class, $(I+\mu R_x)$ is a holomorphic
function of $x$ on some sector 
$S_\delta $ as in Theorem 2.2 and 
$$\det (I+\mu R_x)=\det (I+\mu \Xi_x\Theta_x^\dagger)=\det
(I+\mu\Theta_x^\dagger \Xi_x)=\det (I+\mu
\Gamma_{\Phi_{(x)}}).\eqno(2.21)$$
\noindent By the Riesz functional calculus, $(I+\mu R_x)^{-1}$ is meromorphic for $x$ in some $S_\delta$. Correcting a typographic error in [9, p. 324], we rearrange
terms and calculate the derivative 
$$\eqalignno{ \mu T_\mu
(x,x)&=-\mu {\hbox{trace}}\, \Bigl(Ce^{-xA}(I+\mu
R_x)^{-1}e^{-xA}B\Bigr)\cr
&= -\mu {\hbox{trace}}(I+\mu R_x)^{-1}e^{-xA}BCe^{-xA}\cr
&=\mu {\hbox{trace}}\, \Bigl((I+\mu R_x)^{-1}{{dR_x}\over{dx}}\Bigr)\cr
&={{d}\over{dx}}{\hbox{trace}}\log (I+\mu
R_x).&(2.22)\cr}$$
\noindent This identity is proved for $\vert \mu\vert <1$ and
extends by analytic continuation to the maximal domain of $T_\mu (x,x)$.\par
\indent (iii) Since $A$ commutes with $U(t)$, the domain ${\cal D}(A)$ is invariant under $U(t)$, and the multiplications $B\mapsto U(t)B$, $C\mapsto
CU(t)$ and $e^{-xA}\mapsto U(t)e^{-xA}U(t)$
preserve the hypotheses of Theorem 2.2, so $(-A, U(t)B, CU(t))$ is
$(2,2)$ admissible. By commutativity, we have 
$$\tau (x,t)=\det (I+U(t)R_xU(t)),\eqno(2.23)$$
\noindent which depends continuously on $(x,t)$.\par
\rightline{$\square$}\par
\noindent We refer to $\Sigma (t)=(-A,U(t)B, CU(t) )$ as a deformation of $\Sigma$, and analyze particular cases below.\par

\vskip.05in
\noindent {\bf 3. The Baker--Akhiezer function  of an admissible  linear system}\par
\vskip.05in

\noindent In this section, we consider the Darboux addition rule for potentials and  analyze the transformation $(-A,B,C)\mapsto (-A,B,-C)$ and the effect on the ratios and derivatives of $\tau$ functions.\par

\vskip.05in
\noindent {\bf Definition} {\sl (Baker--Akhiezer function).} (i) Let $(-A,B,C)$ be as in Theorem 2.2, and let 
$$\Sigma_\zeta =(-A, (\zeta I+A)(\zeta I-A)^{-1}B, C)\qquad (\zeta\in
{\bf C}\cup \{ \infty \}\setminus {\hbox{Spec}}(A))\eqno(3.1)$$
\noindent so that $\Sigma_\zeta$ defines a $(2,2)$ admissible linear
systems for $\zeta$ in an open subset of ${\bf C}\cup\{ \infty \}$
which includes $\{ \zeta \in {\bf C} : -\zeta \in S_\theta\}$ for some $\pi/2<\theta <\pi$. We identify $\Sigma_\infty$
with $(-A, B,C)$, and $\Sigma_0$ with $(-A, B,-C)$.\par
\indent (ii) Let $\tau_\zeta $ be the tau function of $\Sigma_\zeta$, and let the Baker--Akhiezer function for the family of linear systems be
$$\psi_\zeta (x)={{\tau_\zeta (x)}\over{\tau_\infty
(x)}}\exp\bigl(\zeta x\bigr).\eqno(3.2)$$
\indent (iii) Let $\tau^*_\zeta (x)=\overline{\tau_{\bar\zeta}(\bar x)}$ as in Schwarz's reflection principle, and let
$$\Sigma^*_\zeta =(-A^\dagger, C^\dagger, B^\dagger  (\zeta I+A^\dagger)(\zeta I-A^\dagger)^{-1} )\qquad (\zeta\in
{\bf C}\cup \{ \infty \}\setminus {\hbox{Spec}}(A^\dagger))\eqno(3.3)$$
\noindent so $\Sigma_\zeta\mapsto\Sigma_\zeta^*$ is an involution, and $\Sigma_\zeta^*$ has tau function $\tau^*$.\par  
\vskip.05in
\indent The
following result introduces a family of solutions of
Schr\"odinger equation corresponding to the $\Sigma_\zeta$ with
an addition rule in the style of Darboux.\par  

\vskip.05in
\noindent {\bf Proposition 3.1.} {\sl Let $(-A,B,C)$ be as in Theorem 2.2.\par
\indent (i) Then for $-\zeta\in S_\theta$, the linear system  $\Sigma_\zeta $ is also $(2,2)$ admissible, and 
the Baker--Akhiezer function satisfies}
$$  -{{d^2}\over{dx^2}}\psi_\zeta (x)+u_\infty (x)
\psi_\zeta (x)=-\zeta^2 \psi_\zeta (x).\eqno(3.4)$$
\indent {\sl (ii) There exist $h_j\in C^\infty ((0, \infty ); {\bf C})$ such that there is an asymptotic expansion 
$$\psi_\zeta (x)\asymp e^{\zeta x} \Bigl( 1+{{h_1(x)}\over{\zeta }}+{{h_2(x)}\over{\zeta^2}}+\dots \Bigr)\eqno(3.5)$$
\noindent as $\zeta\rightarrow\pm i\infty $, and the expansion is uniform for $x$ in compact subsets of $(0, \infty )$.}\par
\vskip.05in

\noindent {\bf Proof} (i) For all $\zeta\in {\bf C}\setminus {\hbox{Spec}}(A),$ there exists $x_0(\zeta )$ such that $\Vert (\zeta I+A)(\zeta I-A)^{-1}R_x\Vert_{{\cal L}^1(H)} <1$ for all $x>x_0(\zeta )$, so that $\tau_\zeta (x)$ is continuously differentiable and non-zero as a function of $x\in (x_0(\zeta ), \infty )$. In particular,  suppose that $\Re \zeta <0$, then $-\zeta\in S_\theta$ so $\zeta I-A$ is invertible. Using the $R$ function for $\Sigma_\zeta$,
we write
$$\eqalignno{ {{\tau_\zeta (x)}\over{\tau_\infty
(x)}}&={{\det \bigl(I+(\zeta I+A)(\zeta I-A)^{-1}R_x\bigr)}\over{\det
\bigl(I+R_x\bigr)}}\cr
&={{\det \bigl(I+(\zeta I-A)^{-1}( (\zeta I-A)R_x+AR_x+R_xA)\bigr)}\over{\det
\bigl(I+R_x\bigr)}}\cr
&={{\det \bigl(I+R_x+(\zeta I-A)^{-1}(AR_x+R_xA)\bigr)}
\over{\det \bigl(I+R_x\bigr)}}&(3.6)\cr}$$
\noindent so that when $AR_x+R_xA$ has rank one, the perturbing term 
$(\zeta I-A)^{-1}(AR_x+R_xA)$ has rank one; continuing we find
$$\eqalignno{ {{\tau_\zeta (x)}\over{\tau_\infty
(x)}}&=\det \bigl(I+(\zeta I-A)^{-1}e^{-xA}BCe^{-xA}(I+R_x)^{-1}\bigr)\cr
&=\det \bigl(I+Ce^{-xA}(I+R_x)^{-1}(\zeta I-A)^{-1}e^{-xA}B\bigr)\cr
&=1 +Ce^{-xA}(I+R_x)^{-1}(\zeta I-A)^{-1}e^{-xA}B, &(3.7)\cr}$$
\noindent since $B:{\bf C}\rightarrow H$ and $C:H\rightarrow {\bf C}$ have rank one. Hence
$$\eqalignno{\psi_\zeta (x)&={{\tau_\zeta (x)}\over{\tau_\infty
(x)}}\exp\bigl(\zeta x\bigr)\cr
&=\exp\bigl(\zeta x\bigr)+
Ce^{-xA}(I+R_x)^{-1}(\zeta I-A)^{-1}e^{-xA}B\exp\bigl(\zeta\bigr)\cr
&=\exp\bigl(\zeta x\bigr)
-\int_x^\infty Ce^{-xA}(I+R_x)^{-1}e^{-yA}B\exp\bigl(
\zeta y\bigr)\,dy\cr
&=\exp\bigl(
 \zeta x\bigr)
+\int_x^\infty T(x,y)\exp\bigl(\zeta y\bigr)\, dy.&(3.8)\cr}$$
\noindent Here $T$ satisfies the Gelfand--Levitan equation, and by integrating by parts, we see that 
$${{\partial^2T}\over{\partial x^2}}-{{\partial^2T}\over{\partial y^2}}=u(x)T(x,y)\eqno(3.9)$$
\noindent where $u(x)=-2{{d^2}\over{dx^2}}\log\tau (x)$. Then by integrating by parts, we see that
$\psi_\zeta$ satisfies\par
\noindent Schr\"odinger's equation.\par 
\indent  The solutions of the differential equation depend analytically on $\zeta$ at those points where the potential depends analytically on $\zeta$; note that $\zeta\mapsto \tau_\zeta (x)$ is holomorphic and non zero for $\Vert R_x\Vert <1$ and $-\zeta \in S_\theta $.  Then we continue the 
solutions analytically to all $-\zeta$ in the sector $S_\theta$, on which $\psi_\zeta (x)$ is holomorphic as a function of $\zeta$ for $x>0$. \par
\indent (ii)  Observe that $X_\theta= S_\theta\cap (-S_\theta )$  contains $i{\bf R}\setminus \{ 0\}$. For $\zeta\in S_\theta\cap (-S_\theta )$, by (i) there exist solutions $\psi_\zeta (x)$ and $\psi_{-\zeta }(x)$ to (3.3). In particular, $\psi_{ik}$ and $\psi_{-ik}(x)$ are solutions for $k>0$. 
We integrate by parts repeatedly
$$\eqalignno{e^{-xA}(\zeta I-A)^{-1}&=e^{-xA}\int_0^\infty e^{\zeta s}e^{-sA}\, ds\cr
&={{e^{-xA}}\over{\zeta}}+{{Ae^{-xA}}\over{\zeta^2}}+\dots +{{A^{k-1}e^{-xA}}\over{\zeta^k}}+\int_0^\infty {{A^ke^{-xA}}\over{\zeta^k}} e^{\zeta s} e^{-sA}\, ds,&(3.10)\cr}$$
\noindent where the integral converges by the hypothesis of Theorem 2.2. Also, $(e^{-zA})$ is an analytic semigroup in the sector $S_{\theta -\pi /2}$, so ${\cal D}(A^j)$ is a dense linear subspace of $H$ for all $j=1,2,\dots$ and  $A^je^{-xA}\in {\cal L}(H)$ and by Cauchy's estimates there exists $C>0$ such that $\Vert A^je^{-xA}\Vert_{{\cal L}(H)}\leq Cj!/x^j$ for all $x>0$. So we can generate an asymptotic expansion of (3.7) with terms 
 $$h_j(x)=Ce^{-xA} (I+R_x)^{-1}A^{j-1}e^{-xA}B\eqno(3.11)$$
\noindent which is bounded on compact subsets of $(0,\infty )$.\par
\rightline{$\square$}\par
\vskip.05in

\noindent {\bf Definition} {\sl(Darboux transforms).}  Let $(-A,B,C)$ be an
$(2,2)$ admissible linear system with tau function $\tau_{\infty}(x; \mu )=\det(I+\mu R_x)$. 
Define the Darboux transform of $(-A,B,C)$ to be
$(-A,B,-C)$ with tau function
transform $\tau_{0}(x;\mu )=\det (I-\mu R_x)$. Let
$$v={{1}\over{\mu}}{{d}\over{dx}}\log {{\tau_{\infty}}\over{\tau_{0}}}, 
\quad w={{1}\over{\mu}}{{d}\over{ dx}}\log \bigl(\tau_{0}\tau_{\infty}\bigr),$$
$$ u_\infty =-{{2}\over{\mu^2}}{{d^2}\over{ dx^2}}
 \log\tau_{\infty},\quad u_0=-{{2}\over{\mu^2}}{{d^2}\over{dx^2}}
 \log\tau_{0}.\eqno(3.12)$$
\vskip.05in
\indent  In the following result, we show how products and quotients of $\tau$ 
functions can be linked by the Gelfand--Levitan equation for $2\times 2$ matrices, and satisfy the
identities usually associated with Darboux transforms in the theory of integrable systems. Identities such as (3.18) also appear in Appendix A16 of Mehta [54].\par
\vskip.05in
\noindent {\bf Theorem 3.2.} {\sl Let $(-A,B,C)$ be a $(2,2)$-admissible linear system with 
input and output spaces ${\bf C}$, and let $\phi (x)=Ce^{-xA}B$.\par
\indent (i) Then there exists $\delta>0$ such that for all $\mu\in {\bf C}$ such 
that $\vert \mu\vert <\delta$, the integral equation (2.16) with
$$T(x,y)=\left[\matrix{W(x,y)&V(x,y)\cr V(x,y)&W(x,y)\cr}\right],\eqno(3.13)$$
$$ \Phi (x+y)
=\left[\matrix{0&\phi (x+y)\cr \phi (x+y)&0\cr}\right]\eqno(3.14)$$
\noindent has a solution such that}
$$W(x,x)={{d}\over{dx}}{{1}\over{2\mu}}\log\bigl(\tau_\infty (x;\mu)\tau_0(x;\mu )\bigr) ,\eqno(3.15)$$
$$V(x,x)={{d}\over{dx}}{{1}\over{2\mu}}\log {{\tau_\infty (x;\mu )}
\over{\tau_0(x;\mu )}}.\eqno(3.16)$$
\noindent {\sl and} 
$${{1}\over{2\mu}}{{d}\over{dx}}W(x,x)=-V(x,x)^2;\eqno(3.17)$$
\indent {\sl (ii) also Toda's equation holds in the form} 
$$\tau_{0}''\tau_{\infty}-2\tau_{0}'\tau_{\infty}'+\tau_{0}^{}\tau_{\infty}''=0.\eqno(3.18)$$


\vskip.05in
\noindent {\bf Proof.}  (i) Let 
$$T_{\infty}(x,y)=-Ce^{-xA}(I+\mu R_x)^{-1}e^{-yA}B,\eqno(3.19)$$
$$T_{0}(x,y)=Ce^{-xA}(I-\mu R_x)^{-1}e^{-yA}B\eqno(3.20)$$
 and
$$\Phi (x)=\left[\matrix{0&\phi (x)\cr \phi (x)&0\cr}\right].\eqno(3.21)$$
Now let 
$$T(x,y)={{1}\over{2}}\left[\matrix{ T_{\infty}+T_{0}&T_{\infty}-T_{0}\cr
T_{\infty}-T_{0}& T_{\infty}+T_{0} \cr}\right]\eqno(3.22)$$
so that 
$$T(x,y)=-\left[\matrix{C&0\cr 0&C\cr}\right]\left[\matrix{e^{-xA}&0\cr 0&e^{-xA}
\cr}\right]
\left[\matrix{I&\mu R_x\cr \mu R_x&I\cr}\right]^{-1}\left[\matrix{e^{-yA}&0\cr 0&e^{-yA}\cr}\right]
\left[\matrix{0&B\cr B&0\cr}\right]\eqno(3.23)$$
\noindent hence $T$ satisfies the Gelfand--Levitan equation 
$$T(x,y)+\Phi (x+y)+\mu \int_x^\infty T(x,z)\Phi (z,y)\, dz=0.\eqno(3.24)$$
\indent ii)  As in Proposition 2.4, 
$$T_{\infty}(x,x)={{1}\over{\mu}}{{d}\over{dx}}\log \tau_{\infty}(x),\eqno(3.25)$$
$$T_{0}(x,x)={{1}\over{\mu}}{{d}\over{dx}}\log \tau_{0}(x);\eqno(3.26)$$
\noindent hence (3.18) is equivalent to the condition
$${{d}\over{dx}}T_{0}(x,x)+\mu \bigl( T_{0}(x,x)-T_{\infty}(x,x)\bigr)^2
+{{d}\over{dx}} T_{\infty}(x,x)=0,\eqno(3.27)$$
\noindent which we now verify. The left-hand side equals
$$Ce^{-xA}\bigl( -A(I-\mu R_x)^{-1}-(I-\mu R_x)^{-1}\mu (AR_x+R_xA)(I-\mu R_x)^{-1}-(I-\mu R_x)^{-1}A\bigr)e^{-xA}B$$
$$+Ce^{-xA}\bigl( (I-\mu R_x)^{-1}+(I+\mu R_x)^{-1}\bigr) e^{-xA}\mu BCe^{-xA}\bigl( (I-\mu R_x)^{-1}+(I+\mu R_x)^{-1}\bigr) e^{-xA}B$$
$$+Ce^{-xA}\bigl( A(I+\mu R_x)^{-1}-(I+\mu R_x)^{-1}\mu (AR_x+R_xA)(I+\mu R_x)^{-1}+
(I+\mu R_x)^{-1}A\bigr)e^{-xA}B\eqno(3.28)$$
All of the terms begin with $Ce^{-xA}$ and end with $e^{-xA}B$, and we can replace $e^{-xA}\mu BCe^{-xA}$ by
$\mu (AR_x+R_xA)$ to obtain
$$\eqalignno{{(3.28)}&=Ce^{-xA}\Bigl(
-2(I-\mu R_x)^{-1}A(I-\mu R_x)^{-1}+4(I-\mu^2R_x^2)^{-1}\mu (AR_x+R_xA)(I-\mu^2R_x^2)^{-1}\cr
&\qquad\qquad \quad +2(I+\mu R_x)^{-1}A(I+\mu R_x)^{-1}\Bigr)e^{-xA}B\cr
&=0.&(3.29)\cr}$$
\noindent This proves (3.18), and one can easily check that (3.18) is equivalent to  
$$u_{0}(x)={{1}\over{\mu}}{{dv}\over{dx}}+v(x)^2, \quad v(x)^2=-{{1}\over{\mu}}{{dw}\over{dx}}.\eqno(3.30)$$
\indent The entries of $T$ satisfy the pair of coupled integral equations
$$\eqalignno{ 0&=W(x,y)+\mu \int_x^\infty V(x,s)\phi (s+y)\, ds\cr
0&=V(x,y)+\phi (x+y)+\mu \int_x^\infty W(x,s)\phi (s+y)\, ds;&(3.31)\cr}$$
\noindent so $W$ satisfies
$$0=-W(x,z)+\mu \int_x^\infty \phi (x+y)\phi (y+z)\, dy+\mu^2 \int_x^\infty W(x,s) \int_x^\infty \phi (s+y)\phi (y+z)\, dy ds,\eqno(3.32)$$
\noindent which explains how $\mu^2\Gamma_{\phi}^2$ enters the discussion in several determinant formulas [67].\par

\rightline{$\square$}

\vskip.05in

\noindent {\bf Definition} (i) {\sl  (Darboux Addition).}  For $-\zeta\in S_\theta\cup\{0\}$ we define the Darboux addition rule on $(2,2)$ admissible linear systems by $M_\zeta :(-A,B,C)\mapsto (-A, (\zeta I+A) (\zeta I-A)^{-1}B, C)$ and on potentials by 
$$u_\infty\mapsto u_\zeta =u_\infty-2(\log \psi_\zeta )''.\eqno(3.33)$$
\indent (ii) Let ${\hbox{Wr}}(\varphi, \psi )$ be the Wronskian of $\psi, \varphi\in C^1((0,\infty );{\bf C})$.\par
\vskip.05in
\noindent {\bf Corollary 3.3.} {\sl The set $\{ M_\zeta$,  $(\zeta\in X_\theta ), M_0, M_\infty =I\} $ generates a group such that $M_0^2=I$,  $M_\zeta M_{-\zeta }=I$ and $M_\zeta M_\eta$ corresponds to adding} $-2{{d^2}\over{dx^2}}\log {\hbox{Wr}}(\psi_\zeta, \psi_\eta )$ {\sl to the potential.}\par
\vskip.05in 
\noindent {\bf Proof.}  The definition is consistent with  [22, p
484], and p. 414 of [47]. In particular, $\psi_0(x)=\tau_0(x)/\tau_\infty (x)$, and  
$u_0(x) =u_\infty (x)-2{{d^2}\over{dx^2}}\log\psi_0(x)$, which is consistent with (3.18).\par
For $\zeta_1\neq \zeta_2$, let $\Psi (x)=Wr( \psi_{\zeta_1}, \psi_{\zeta_2})/\psi_{\zeta_2}$, and observe that
$$\Psi''=\bigl(\zeta_2^2+u_\infty -2(\log \psi_{\zeta_1})''\bigr)\Psi .$$
\noindent This gives the basic composition rule for $M_{\zeta_2}M_{\zeta_1}.$ The other statements follow from Proposition 3.1 and Theorem 3.2.\par
\rightline{$\square$}\par
\vskip.05in
\noindent {\bf Definition} (i) Let $\Omega$ be a domain in ${\bf C}$. A divisor $\delta $ on $\Omega$ is a function
$\delta:\Omega\rightarrow {\bf Z}$ such that the restriction of $\delta $ to $K$ has finite support for all
compact subsets $K$ of $\Omega$.\par
\indent (ii)  For a meromorphic function $f$ on $\Omega$, we let $\nu (z)$ be the order
of $z$ as a zero of $f$, and $\nu (p)$ be the order of $p$ as a pole. Then $(f)=\sum_z \nu (z)\delta_z -\sum_p \nu (p)\delta_p$ defines the principal
divisor corresponding to $f$; see [65]. \par
\indent (iii) Let $\log_+x=\max \{ 0, \log x\}$ for $x>0$. For any meromorphic function $f$ let
$$m(r;f)=\int_0^{2\pi} \log_+ \vert f(re^{i\theta })\vert {{d\theta}\over{2\pi}}\qquad (r>0).\eqno(3.34)$$
\vskip.05in
\noindent {\bf Proposition 3.4.} {\sl Suppose that $\phi :(0,\infty )\rightarrow {\bf R}$ arises from a $(2,2)$ admissible linear system. Let $\tau_\infty  (x;\mu )=\det
(I+\mu\Gamma_{\phi_{(x)}}), $ and $\tau_0 (x;\mu )=\tau_\infty (x;-\mu )$ for all $x>0$ and $\mu\in {\bf C}$.\par
\indent (i) Then $\mu\mapsto \tau_\infty  (x;\mu )$ is entire, and $\overline{ \tau_\infty  (x;\bar\mu )}=\tau_\infty  (x;\mu
)$ for all $x>0$.\par
\indent (ii) Let 
$$q(\mu )= {{\tau_\infty  (0;\mu )}\over{\tau_0 (0;\mu )}}.\eqno(3.35)$$
\indent $\bullet$ Then $q$ is meromorphic, $q^*(\mu )=q(\mu )$ and 
$q(\mu )q(-\mu )=1$;\par
\indent $\bullet$ all the zeros and poles of $q$ are simple and lie in ${\bf R}\setminus \{ 0\}$;\par 
\indent $\bullet$ the zeros $(z_j)$ of $q$ satisfy $\sum_{j=1}^\infty 1/\vert z_j\vert <\infty$, and \par
\indent $\bullet$ $m(r;q)/r\rightarrow 0$ as $r\rightarrow\infty$.\par
\indent (iii) Conversely, let $q$ satisfy the conclusions of (ii). Then there exists a balanced linear system $(-A,B,C)$ with input and output space ${\bf C}$, such that
$q(\mu )=\kappa {{\tau_\infty (0; \mu )}/{\tau_0(0;\mu )}}$ for some constant $\kappa$.}\par

\vskip.05in
\noindent {\bf Proof.} (i) This follows from Theorem 2.2 (ii).\par
\indent (ii)  Since $\Gamma_\phi$ is self-adjoint and trace 
class, the spectrum consists of $0$,  together with non-zero real eigenvalues 
$\lambda_j$ with multiplicity $\nu (\{ \lambda_j\})$, and for a self-adjoint 
$\Gamma_\phi$ the algebraic and geometric multiplicity are equal. For a Hankel operator $\Gamma_\phi$, the dimension of 
$\{ \xi: \Gamma \xi =0\}$ is either zero or infinity; so by compressing $\Gamma_\phi$ to the closure of its range, we can assume $0$  is not an eigenvalue. 
Then the function $q(\mu )=\tau_ \infty (\mu) /\tau_0(\mu )$ is meromorphic 
on ${\bf C}$, the identity $q(\mu )q(-\mu )=1$ is trivially true, while $q=q^*$ holds by (i).\par
\indent  The formal difference of the zeros and poles of $q$ on 
$\{\mu\in {\bf C}:\vert\mu\vert<\rho\}$ may be represented by the divisor 
$$\eqalignno{(q)&= \sum_{j:\vert 1/\lambda_j\vert<\rho} 
(\nu (\{\lambda_j\})\delta_{-1/\lambda_j}-\nu (\{\lambda_j\})\delta_{1/\lambda_j})\cr
&=\sum_{j:\lambda_j \rho>1} (\nu
(\lambda_j)-\nu(-\lambda_j))(\delta_{-1/\lambda_j}
-\delta_{1/\lambda_j}),&(3.36)\cr}$$ 
where $\nu (\{ \lambda_j\})-\nu (\{-\lambda_j\})$ belongs to $\{-1, 0, 1\}$ by a 
theorem of Megretskii, Peller and Treil [53].  Consider $\lambda>0$. If $\nu (\lambda )-\nu (-\lambda )=-1$, then  $q$ has a simple zero at $-1/\lambda$ and a simple pole at $1/\lambda$; whereas if $\nu (\lambda )-\nu (-\lambda )=1$, then $q$ has a simple zero at $1/\lambda$ and a simple pole at $-1/\lambda$.\par
\indent Also $\sum_{j=1}^\infty \vert\lambda_j\vert $ converges since $\Gamma_\phi$ is trace class, so $\sum_{j=1}^\infty 1/\vert z_j\vert$ converges. Finally, standard results on convergent infinite products from [31', Theorem 1.9] show that $m(r; q)/r\rightarrow 0$ as $r\rightarrow\infty$.\par
\indent (iii) Suppose conversely that $q$ satisfies the conclusions of (ii). Now $q(0)^2=1$ so  $0$  is neither a zero nor a pole of $q$; also $z$ is a zero of $q$ if and only if $p=-z$ is a pole of $q$, so the divisor of $q$ is 
$$(q)=\sum_z \delta_z-\delta_{-z},\eqno(3.37)$$
\noindent where we have summed over the set of zeros $z$ of $q$. Now $z\in {\bf R}\setminus \{ 0\}$ for all the zeros, and the $z$ have no point of accumulation on ${\bf R}$ since $q$ is meromorphic, so we can list the zeros in a (possibly finite) sequence $(z_j)_{j=1}^\infty$ such that $(\vert z_j\vert)_{j=1}^\infty$ is strictly increasing. By hypothesis $\sum_{j=1}^\infty 1/\vert z_j\vert$ converges.\par
\indent Let $\lambda_{j}=-1/z_j$  and define $\nu :{\bf R}\rightarrow {\bf Z}$ by $\nu (\lambda_j)=1$ for $j=1, 2,\dots $, and $\nu (\lambda )=0$ otherwise. Then there exists a self-adjoint Hankel operator $\Gamma$ on $L^2((0, \infty ); {\bf C})$  with spectral multiplicity function $\nu$ and spectrum $\{ \lambda_j: j=1,2, \dots\}\cup\{ 0\}$ . Since $\sum_{j=1}^\infty \vert\lambda_j\vert$ converges, $\Gamma$ is trace class and has a Fredholm determinant, so  
$${{\det (I+\mu \Gamma )}\over {\det (I-\mu \Gamma )}}=\prod_{j=1}^\infty {{1-\mu/z_j}\over{1+\mu/z_j}}\eqno(3.38)$$
\noindent is meromorphic with only simple zeros at $z_j$ and simple poles at $-z_j$. Hence the product in (3.38) has the same divisor as $f$, and we deduce that there exists an entire function $f$ such that $q(\mu )e^{-f(\mu )}$ equals the product in (3.38). The next step is to eliminate this $f$ by using elementary Nevanlinna theory.\par
\indent Let $n(r)=\sharp \{ z_j: \vert z_j\vert \leq r\}$ be the counting function for zeros, or equivalently poles of $q$, and then let $N(r;q)=\int_0^r n(s)s^{-1}ds.$ By standard results, we have convergent expressions 
$$\sum_{j=1}^\infty {{1}\over {\vert z_j\vert }}=\int_0^\infty {{n(r)dr}\over {r^2}},\eqno(3.39)$$
\noindent from  which we deduce that $ n(r)/r\rightarrow 0$ and $N(r;q)/r\rightarrow 0$ as $r\rightarrow\infty$; so the Nevanlinna characteristic satisfies $N(r;q)+m(r;q)=o(r)$ as $r\rightarrow\infty$. By  [31', Theorem 1.9], we deduce that $f(\mu )$ is a polynomial of degee zero, namely a constant.\par
\indent By Theorem 1.1 of [53], $\Gamma$ may be realized by a balanced
linear system $(-A,B,C)$ with input and output space ${\bf C}$, where $A\in {\cal L}(H)$ bounded, and $\Gamma =\Gamma_\phi$ as in (1.1), where $\phi (t)=Ce^{-tA}B$.\par 
\rightline{$\square$}\par
\vskip.05in
\indent We now show how a Schr\"odinger differential equation of scattering type gives rise to an admissible linear system as in Theorem 3.2;
this justifies the terminology 'scattering function' as applied to $\phi$.  In section 4 of [9], we realized the scattering data from Schr\"odinger's equation from a linear system with unbounded $A$; in the following result, we realize the data with a linear system with $A\in {\cal L}(H)$ .The differential equation
$-f''+(v'+v^2)f=k^2f$ may be written as 
$${{d}\over{dx}}\left[\matrix{f\cr g\cr}\right]=\left[\matrix{v&ik\cr 
ik&-v\cr}\right]\left[\matrix{f\cr g\cr}\right]\eqno(3.40)$$
\noindent  where we suppose that $v\in C^\infty _c({\bf R}; {\bf R})$ for simplicity. For
$\lambda =k^2>0$, there exists a solution $f(x)=\psi (x;k)$ such that
$$\psi (x;k)\asymp \cases{e^{-ikx}+s_{21}(k)e^{ikx},& as $x\rightarrow\infty$;\cr
\bar s_{11}(k)e^{-ikx},& as $x\rightarrow -\infty$.\cr}\eqno(3.41)$$
\noindent There may also be a discrete spectrum $\lambda_j=-\kappa^2_j$ with $\kappa_n\geq \dots \geq \kappa_1>0$, where each $\kappa_j^2$ is associated to a real eigenfunction $\psi (x;-\kappa_j^2)$ called a bound state that is asymptotic to $c(-\kappa_j^2) e^{-\kappa_jx}$ as $x\rightarrow\infty$, where $c(-\kappa_j^2)$ is normalized by taking $\int_{-\infty}^\infty \psi (x;-\kappa_j^2)^2 dx=1$. The discrete spectrum gives rise to a linear system with finite dimensional state space, as we discuss in Proposition 4.5. Note that $-\kappa_j-\kappa_\ell <0$ for all $j, \ell$.\par
\indent So we suppose for the moment that the discrete spectrum is absent, and the inverse spectral problem is to recover $v$ from  $s_{21}(k)$ and
$s_{11}(k)$, up to equivalence. We aim to introduce an admissible linear system from the spectral data and determine the potential.\par 
\indent As in [22], the scattering matrix is 
$$S(k)=\left[\matrix{s_{11}(k)& -\bar s_{21}(k)\cr s_{21}(k)& \bar
s_{11}(k)\cr}\right] \in SU(2)\eqno(3.42)$$
\noindent where  $s_{21}(k)$ is called the reflection coefficient and $s_{11}(k)$ is the
transmission coefficient for $k\in {\bf R}$. We suppose
that $s_{21}(-k)=\overline{s_{21}(k)}$, and that $s_{21}\in L^2({\bf
R}; {\bf C})$ is absolutely continuous with $s_{21}'\in L^2({\bf R};
{\bf C})$. Let
$$a_n ={{(-1)^n}\over{2\pi}}\int_{-\infty}^\infty s_{21}(k) 
{{(1/2+ik)^n}\over{(1/2-ik)^{n+1}}}dk,\eqno(3.43)$$
\vskip.05in

\noindent {\bf Proposition 3.5.} {\sl  Suppose that (3.41) has no bound states, and suppose that $(a_n)_{n=0}^\infty$ satisfy 
$$\lim\sup_{n\rightarrow\infty }\vert a_n\vert^{1/n}<1/3.\eqno(3.44)$$
\indent  Then there exists a $(2,2)$ admissible linear system $(-A,B,C)$
with $H_0={\bf C}$ and bounded $A,B,C$ such that the potential $u_0=v'+v^2$ 
of $-f''+u_0f=\lambda f$ is determined by the $\tau$ function of $(-A,B,-C)$.}\par
\vskip.05in
\noindent {\bf Proof.}  We introduce the scattering function
$$\phi (x)={{1}\over{2\pi}}\int_{-\infty}^\infty s_{21}(k)e^{ixk}\, dk\qquad (x>0),\eqno(3.45)$$
\noindent so that $\phi\in L^2((0, \infty ); {\bf R})$ and $x\phi (x)\in
L^2((0, \infty ); {\bf R})$; hence $\Gamma_\phi$ determines a self-adjoint
and Hilbert--Schmidt operator on $L^2((0, \infty ); {\bf C})$. We now
introduce a linear system which realizes this $\phi$.\par
\indent Let $L_n(x)=e^{x}(d/dx)^n (x^ne^{-x})/n!$ be the Laguerre
polynomial of order zero and degree $n$; then the $n^{th}$ Laguerre
function is 
$$e^{-x/2}L_n(x)={{1}\over{2\pi i}}\int_{C(1,\delta )}e^{-xz/2}
{{(1+z)^n}\over{(1-z)^n}} {{dz}\over{1-z}},\eqno(3.46)$$
\noindent where $C(1,\delta )$ is the circle of centre $1$ and radius
$0<\delta <1$ in ${\bf C}$. This is a special case of a formula of Tricomi,
and follows from Cauchy's integral formula. Then $(e^{-x/2}L_n(x))_{n=0}^\infty$ gives an
orthonormal basis for $L^2(0, \infty )$. By Plancherel's formula we find that
$$a_n=\int_0^\infty \phi (x) e^{-x/2}L_n(x)dx\qquad (n=0, 1, 2, \dots )\eqno(3.47)$$
\noindent are the coefficients of $\phi$ with respect to $(e^{-x/2}L_n(x))_{n=0}^\infty$ and we introduce
$$b(z)= \sum_{n=0}^\infty a_n{{(1+z)^n}\over{(1-z)^n}}\qquad (z\in
C(1,\delta ));\eqno(3.48)$$  
\noindent by the hypothesis, $\lim\sup_{n\rightarrow\infty }\vert a_n\vert^{1/n}=\rho<1/3$ so the series for $b(z)$ converges for all $z\in {\bf C}$ outside of the disc $D((1+\rho^2)/(1-\rho^2); 2\rho /(1-\rho^2)).$ Hence can choose $0<\delta<1$ so $(2+\delta )\lim\sup_{n\rightarrow\infty }\vert a_n\vert^{1/n}/\delta <1$ so that this series converges absolutely and uniformly on $C(1, \delta )$.\par
 \indent We parametrize
$C(1, \delta )$ by $z=1+\delta e^{i\theta }$ for $\theta \in
[0, 2\pi ]$, and introduce the Hilbert space 
$H=L^2( C(1,\delta ); d\theta ; {\bf C})$. Then we introduce bounded linear operators 
$$\eqalignno{ A:H\rightarrow H:&\quad f(z)\mapsto zf(z)/2\cr
               B: {\bf C}\rightarrow H: &\quad \beta\mapsto \beta
b(z)\cr
C: H\rightarrow {\bf C}: &\quad g\mapsto {{1}\over{2\pi
i}}\int_{C(1,\delta )} g(z) {{dz}\over{1-z}}.&(3.49)\cr}$$                      

\noindent Thus $(-A,B,C)$ gives a linear system with state space $H$,
input and output space ${\bf C}$. \par
\indent From the orthogonal series expansion $\phi
(t)=\sum_{n=0}^\infty a_n e^{-t/2} L_n(t)$, we note that 
$$\phi (t)={{1}\over{2\pi i}}\int_{C(1, \delta )} \sum_{n=0}^\infty 
{{a_n(1+z)^ne^{-zt/2}dz}\over{(1-z)^{n+1}}}
=Ce^{-tA}B.\eqno(3.50)$$
\noindent  Also, the corresponding $R$  operator has integral kernel 
$$R_x=\int_x^\infty e^{-tA} BCe^{-tA} dt\leftrightarrow
{{e^{-xz/2-xw/2} b(z)}\over{z+w}}.\eqno(3.51)$$ 
\noindent We observe that $\Re z\geq 1-\delta$ for all $z$ on $C(1, \delta )$; hence $\Vert
e^{-tA}\Vert_{{\cal L}(H)}\leq e^{-t(1-\delta )/2}$; so $(-A,B,C)$ is $(2,2)$ admissible by 
Theorem  2.2(i).\par
\indent For all $\phi \in L^2((0, \infty ); {\bf R})$  we have $(a_n)_{n=0}^\infty \in \ell^2$ and hence $b(iy)/(1-iy)$ gives a function in $L^2({\bf R}; {\bf C})$. Hence by deforming (3.48) to an integral along the imaginary axis $\Re z=0$, we can obtain an $L^2$ Plancherel integral 
$$\phi (t)={{1}\over{2\pi }}\int_{-\infty}^\infty e^{-ity} {{b(iy)dy}\over{1-iy}}.\eqno(3.52)$$ 
\noindent However, the operation of multiplication by $iy$ is clearly unbounded on this deformed contour.\par 
\indent We now recover $u$ from $(-A,B,C)$ via the Gelfand--Levitan equation. The kernel
$$T_{0}(x,y)=Ce^{-xA}(I-R_x)^{-1}e^{-yA}B\qquad (x_0<x<y)\eqno(3.53)$$
\noindent as in (2.17) gives the unique solution of the integral equation
$$-T_{0}(x,y)+\phi (x+y)+\int_x^\infty T_{0}(x,y) \phi (y+z)dz=0\qquad (x_0<x<y)\eqno(3.54)$$
\noindent for some $x_0>0$. Then one checks that 
$${{\partial^2 T_{0}}\over{\partial x^2}}-{{\partial^2 T_{0}}\over{\partial y^2}}-q(x)T_{0}(x,y)=0,\eqno(3.55)$$
\noindent where $q(x)=-2(d/dx)T_{0}(x,x)$. Also, one has $T_{0}(x,x)=(d/dx)\log\det
(I-R_x)$, so $q(x)=-2(d^2/dx^2)\log\tau (x)$. 
Given this partial equation for $T_{0}(x,y)$ one can check that 
$$f(x;k) =e^{ikx} +\int_x^\infty e^{iky} T_{0}(x,y)dy\eqno(3.56)$$
\noindent satisfies $-f''(x;k)+q(x)f(x;k)=k^2f (x;k)$. Hence we can identify 
$q(x)$ with $u_0(x)=v'(x)+v(x)^2$, and thus obtain a solution of the inverse spectral
problem.\par

\rightline{$\square$}\par

\vskip.05in

\noindent {\bf 4 The state ring associated with an admissible linear
system}\par
\vskip.05in
\indent Gelfand and Dikii [26] considered the algebra  
${\cal A}_u={\bf C}[u,u',u'',
\dots ]$ of complex polynomials in a smooth potential $u$ and its derivatives. They
showed that if $u$ satisfies the stationary higher
order $KdV$ equations (5.1), then ${\cal A}_u$ is a Noetherian 
ring [6] and the associated Schr\"odinger equation is integrable by quadratures; see [14, 63]. In this section, we introduce an analogue ${\cal A}_\Sigma$ for an
admissible linear system. In the subsequent section, we link this to the result from [26].\par
\indent We introduce these state rings in this section, 
and develop a calculus for $R_x$ which is the counterpart of 
P\"oppe's functional calculus for Hankel operators from [61, 62]. 
As we see in subsequent sections, our theory of state rings has 
wider scope for generalization.\par
\vskip.05in
\noindent {\bf Definition} {\sl (Differential rings).}  (i) Let ${\cal R}$ be a ring with ideal
${\cal J}$, and let $\partial :{\cal R}\rightarrow {\cal R}$ be a
derivation. Then ${\cal R}_{\cal J} =\{r\in {\cal R}: \partial (r)\in {\cal
J}\}$ gives a subring of ${\cal R}$, the ring of constants relative to
${\cal J}$. When ${\cal R}$ is an algebra over ${\bf C}$ and ${\cal
J}=0$, we call ${\cal R}_0 $ the constants; see [63].\par
 
\indent (ii) Let $H$ and $H_0$ be separable complex Hilbert spaces, let
${\cal L}(H)$ be the ring of bounded linear operators on $H$. Let ${\cal S}$ be a subring of $C^\infty ((0,\infty );{\cal L}(H))$  and ${\cal B}$ be a subring of  $C^\infty ((0,\infty );{\cal L}(H_0))$; that is we suppose that each 
$T\in {\cal S}$ is a differentiable function of $x\in (0,\infty )$ as we
indicate by writing $T_x$; we suppose further that $dT_x/dx\in
{\cal S}$, and that $(d/dx)(ST)=(dS/dx)T+S(dT/dx)$. Then
${\cal S}$ is a differential ring. When $I\in {\cal S}$ and $\theta\in {\bf C}$, we
identify $\theta I$ with $\theta$ to simplify notation.\par
\vskip.05in

\noindent {\bf Definition} {\sl (State ring of a linear system).} Let $(-A,B,C)$ be a linear system such that 
$A\in {\cal L}(H)$. Suppose that:\par
\indent (i)  ${\cal S}$ is a differential subring of
$C^\infty ((0,\infty ); {\cal L}(H))$;\par
\indent (ii) $I, A$ and $BC$ are constant elements of ${\cal S}$;\par
\indent (iii) $e^{-xA}$, $R_x$ and 
$F_x=(I+R_x)^{-1}$ belong to
${\cal S}$.\par
\noindent Then ${\cal S}$ is a state ring for $(-A,B,C)$.\par
\vskip.05in
\noindent {\bf Lemma 4.1.} {\sl Suppose that $(-A,B,C)$ is a linear
system with bounded $A$ and that $R_x$ gives a solution of Lyapunov's
equation (1.3) such that $I+R_x$ is invertible for $x>0$ with inverse
$F_x$. }\par

\indent {\sl Then the free associative algebra ${\cal S}$ generated by $I,R_0, A,
F_0, e^{-xA},R_x$
and $F_x$ is a state ring for $(-A,B,C)$ on $(0, \infty
)$. For all $t>0$, there exists a ring homomorphism $S_t:{\cal S}\rightarrow {\cal S}$ given by $S_t:G(x)\mapsto G(x+t)$ such that $S_t$ commutes with $d/dx$}\par
\vskip.05in
\noindent {\bf Proof.}  We can regard ${\cal S}$ as a subring of $C_b((0, \infty ), {\cal L}(H))),$ so the multiplication is
well defined. Then we note that $BC=AR_0+R_0A$ belongs to ${\cal S}$, as required. We also note that $(d/dx)e^{-xA}=-Ae^{-xA}$ and that 
Lyapunov's equation (1.3) gives
$${{d}\over{dx}}(I+R_x)^{-1}=(I+R_x)^{-1}(AR_x+R_xA)(I+R_x)^{-1},\eqno(4.1)$$
\noindent which implies 
$${{dF_x}\over{dx}}=AF_x+F_xA-2F_xAF_x.\eqno(4.2)$$
\noindent {with the initial condition}
$$AF_0+F_0A-2F_0AF_0=F_0BCF_0.\eqno(4.3)$$
\noindent Hence ${\cal S}$ is a differential ring.\par
\indent We can map $I\mapsto I$, $e^{-xA}\mapsto e^{-(x+t)A}$, $R_0\mapsto e^{-tA}R_0e^{-tA}$, $R_x\mapsto e^{-tA}R_{x}e^{-tA}$ and 
$F_x\mapsto (I+e^{-tA}R_xe^{-tA})^{-1}$, and thus produce a ring 
homomorphism $G(x)\mapsto G(x+t)$ which satisfies 
$(d/dx)S_tG(x)=G'(x+t)=S_t(d/dx)G(x)$.\par
\rightline{$\square$}\par
\vskip.05in
\noindent {\bf Definition} {\sl (Products and brackets).} (i) Given a state ring ${\cal S}$ for $(-A,B,C)$, and let ${\cal B}$ be any 
differential ring of
functions from  $(0,\infty)\rightarrow {\cal L}(H_0)$.
Let
$${\cal A}_\Sigma={\hbox{span}}_{\bf C}\{A^{n_1}, A^{n_1}F_xA^{n_2}\dots F_xA^{n_r}: n_j\in {\bf
N}\}.\eqno(4.4)$$
\indent (ii) On ${\cal S}$ we introduce the associative product $\ast$ by 
$$P\ast Q=P(AF+FA-2FAF)Q,\eqno(4.5)$$
which is distributive over the standard addition, and 
 the derivation ${\partial }: {\cal S}\rightarrow {\cal S}$ by
$$\partial P=A(I-2F)P+{{dP}\over{dx}}+P(I-2F)A,\eqno(4.6)$$
\indent (iii)  Let  $\lfloor\,\cdot\,\rfloor :{\cal S}\rightarrow {\cal B}$ be 
the linear map
$$\lfloor Y\rfloor =Ce^{-xA}F_xYF_xe^{-xA}B \qquad (Y\in {\bf
S}),\eqno(4.7)$$
\noindent so that $x\mapsto \lfloor Y\rfloor$ is a function $(x_0, x_1)\rightarrow {\cal L}(H_0)$.\par
\vskip.05in
\noindent {\bf Proposition 4.2.} {\sl Then  $({\cal A}_\Sigma, \ast , \partial )$  is a
differential ring, and there is a 
homomorphism of differential rings 
$({\cal A}_\Sigma , \ast , \partial )\rightarrow ({\cal B}, \cdot , d/dx)$ 
given by $P\mapsto \lfloor P\rfloor .$}\par

\vskip.05in
\noindent {\bf Proof.} The basic observation is that
$dF/dx=AF+FA-2FAF$, so one can check that
$$\partial (P\ast Q)=(\partial P)\ast Q+P\ast (\partial Q);\eqno(4.8)$$
\noindent hence $({\cal S}, \ast , \partial )$ is a
differential ring. \par
\indent We can multiply elements in ${\cal S}$
by concatenating words and taking linear combinations. Since all
words in ${\cal A}_\Sigma$ begin and end with $A$, we obtain words of the
required form, hence ${\cal A}_\Sigma$ is a subring. To differentiate a 
word in ${\cal A}_\Sigma$ we add words
in which we successively replace each $F_x$ by $AF_x+F_xA-2F_xAF_x$,
giving a linear combination of words of the required form.\par

\indent From the definition of $R_x$, we have
$AR_x+R_xA=e^{-xA}BCe^{-xA}$, and hence
$$F_xe^{-xA}BCe^{-xA}F_x=AF_x+F_xA-2F_xAF_x,\eqno(4.9)$$
\noindent which implies
$$\eqalignno{ \bigl\lfloor P\bigr\rfloor
\bigl\lfloor Q\bigr\rfloor&=Ce^{-xA}F_xPF_xe^{-xA}BCe^{-xA}F_xQF_xe^{-xA}B
\cr
&= Ce^{-xA}F_xP(AF_x+F_xA-2F_xAF_x)QF_xe^{-xA}B\cr&=\bigl\lfloor
 P(AF_x+F_xA-2F_xAF_x)Q\bigr\rfloor\cr
&=\lfloor P\ast
Q\rfloor.&(4.10)\cr}$$
\noindent Moreover, the first and last terms in $\lfloor P\rfloor$ have derivatives
$${{d}\over{dx}}Ce^{-xA}F_x=Ce^{-xA}F_xA(I-2F_x), \qquad
{{d}\over{dx}}F_xe^{-xA}B=(I-2F_x)AF_xe^{-xA}B,\eqno(4.11)$$
\noindent which implies (3.7).\par
\indent 
The bracket operation satisfies
$${{d}\over{dx}}\bigl\lfloor P\bigr\rfloor=\Bigl\lfloor
A(I-2F_x)P+{{dP}\over{dx}}+P(I-2F_x)A\Bigr\rfloor =\lfloor
\partial P\rfloor.\eqno(4.12)$$
\rightline{$\square$}\par
\noindent For $x_0\geq 0$ and $0<\phi<\pi$, let $S_\delta^{x_0}$ be the translated sector $S_\delta^{x_0}=\{ z=x_0+w: w\in {\bf C}\setminus \{ 0\}; \vert\arg w\vert<\delta\}$ and let $H^\infty ( S_\delta^{x_0})$ the the bounded holomorphic complex functions on $S_\delta^{x_0}$. Then let $H^\infty_\infty =\cup_{x_0>0} H^\infty ( S_\delta^{x_0})$ be the algebra of complex functions which are bounded on some translated sector  $S_\delta^{x_0}$, with the usual pointwise multiplication.\par
\vskip.05in

\noindent {\bf Theorem 4.3.} {\sl  Let $(-A,B,C)$ be a $(2,2)$-admissible linear system 
with $H_0={\bf C}$ as in Theorem 2.2, so $(e^{-zA})$ for $z\in S_\phi^0$ is a bounded holomorphic semigroup  on $H$. Let $\Theta_0 =\{ P\in {\cal A}_\Sigma:\lfloor P\rfloor=0\}$.\par
\indent (i) Then $({\cal A}_\Sigma, \ast, \partial )$ is a differential ring with bracket $\lfloor \cdot\rfloor$;\par 
 \indent (ii) there is a homomorphism of differential rings  $\lfloor\,
\cdot \, \rfloor: ({\cal A}_{\Sigma}, \ast , \partial )\rightarrow (H_\infty^\infty, \cdot, d/dz)$;\par 
\indent (iii) $\Theta_0$ is a differential ideal in $({\cal A}_\Sigma, \ast,\partial )$ such that ${\cal A}_\Sigma/\Theta_0$ is a
commutative differential ring, and an integral domain.}\par
\vskip.05in
\noindent {\bf Proof} (i) In this case $A$ is possibly unbounded as an operator, so we use the holomorphic semigroup to ensure that products and brackets are well defined. We observe that ${\cal A}_\Sigma$ has a grading ${\cal A}_\Sigma=\oplus_{n=1}^\infty A_n$, where $A_n$ is the span of the elements that have total degree $n$ when viewed as products of $A$ and $F$. For $X_n\in A_n$ and $Y_m\in A_m$, we have $X_n\ast Y_m\in A_{n+m+2}\oplus A_{n+m+3}$ and $\partial X_n\in A_{n+1}\oplus A_{n+2}$. \par
\indent Also we have $A^ke^{-zA}\in {\cal L}(H)$ for all $z\in S_\phi^0$ and $\Vert A^ke^{-zA}\Vert_{{\cal L}(H)}\rightarrow 0$ as $z\rightarrow\infty $ in $S_\phi^0$; hence $R_zA^k\rightarrow 0$ and $A^kR_z\rightarrow 0$ in ${\cal L}(H)$ as $z\rightarrow \infty$  in $S^0_\phi$. Hence there exists an increasing positive sequence $(x_k)_{k=0}^\infty$ such that  $A^kF_z-A^k\in {\cal L}(H)$ for all $z\in S^{x_k}_\phi$ and   $A^kF_z-A^k\rightarrow 0$  in ${\cal L}(H)$ as $z\rightarrow\infty $ in $S^{x_k}_\phi$. Let $X_n\in A_n$ and consider a typical summand $AF_zA^kF_z\dots A$ in $X_n$; we replace each factor like $A^kF_z$ by the sum of $A^k(F_z-I)$ and $A^k$ where $k\leq n$; then we observe that there in an initial factor $Ce^{-zA}$ and a final factor $e^{-zA}B$ in $\lfloor X_n\rfloor$; hence $\lfloor X_n\rfloor$ determines an element of $H^\infty (S_\phi^{x_n})$. \par
\indent (ii) We can identify $H^\infty_\infty$ with the algebraic direct limit $H^\infty_\infty =\lim_{n\rightarrow\infty} H^\infty ( S_\phi^{x_0+n})$.  By the principle of isolated zeros, the multiplication on $H^\infty_\infty $ is consistently defined, and $H^\infty_\infty $ is an integral domain. Now each $f\in H^\infty_\infty$ gives $f\in H^\infty (S_\phi^{x_0})$ so $f'\in H^\infty (S^{x_0+1}_\phi )$ by Cauchy's estimates, so $f'\in H^\infty-\infty $. From (i) we deduce that $\lfloor\,\cdot\, \rfloor: \oplus_{n=1}^\infty A_n\rightarrow \cup_{n=1}^\infty H^\infty (S_\phi^{x_n})$ is well-defined and the bracket is multiplicative with respect to $\ast$, and behaves naturally with respect to differentiation.\par  
\indent  (iii) We check that $\lfloor\,
.\, \rfloor$ is a trace on $({\cal A}_{\Sigma}, \ast , \partial )$, by computing 
$$\eqalignno{\lfloor P\ast Q\rfloor
&={\hbox{trace}}\bigl(Ce^{-xA}FPFe^{-xA}BCe^{-xA}FQFe^{-xA}B\bigr)\cr
&={\hbox{trace}}\bigl(Ce^{-xA}FQFe^{-xA}BCe^{-xA}FPFe^{-xA}B\bigr)\cr
&=\lfloor Q\ast P\rfloor .&(4.13)\cr}$$
\noindent Hence  $\Theta_0$ contains all the commutators $P\ast Q-Q\ast P$, and $\Theta_0$ is the kernel of the homomorphism $\lfloor\, \cdot\,\rfloor$, hence is an ideal  for $\ast$. Also,  we observe that for all $Q\in \Theta_0$, we have
$\partial Q\in \Theta_0$ since $\lfloor \partial Q\rfloor =(d/dx)\lfloor Q\rfloor =0$. 
Hence $\Theta_0$ is a differential ideal which contains the commutator subspace of 
$({\cal A}_\Sigma, \ast )$, so ${\cal A}_\Sigma/\Theta_0$ is a commutative algebra. Also, $\partial $
determines a unique derivation $\bar \partial $ on ${\cal S}/\Theta_0$ by $\bar\partial
Q=\partial Q+\Theta_0$ for all $Q\in {\cal S}$; hence ${\cal A}_\Sigma /\Theta_0$ is a
differential algebra. \par
\indent We can identify  ${\cal A}_\Sigma /\Theta_0$ with a subalgebra of $H_\infty^\infty$, which is an integral domain.\par
\rightline{$\square$}\par


\noindent {\bf Remarks 4.4} (i) P\"oppe [61, 62] introduced a linear functional 
$\lceil\,.\, \rceil$ on Fredholm kernels $K(x,y)$ on $L^2(0, \infty )$ 
by $\lceil K\rceil=K(0,0).$ In particular, let $K,G,H,L$ be 
integral operators on $L^2(0, \infty )$ that
have smooth kernels of compact support, 
let $\Gamma =\Gamma_{\phi_{(x)}}$ have kernel $\phi (s+t+2x)$, let 
$\Gamma'={{d}\over{dx}}\Gamma$ and
 $G=\Gamma_{\psi_{(x)}}$ be another Hankel
operator; then the trace satisfies
$$\eqalignno{ \lceil\Gamma \rceil&=-
{{d}\over{dx}}{\hbox{trace}}\, \Gamma
&(4.14)\cr
\lceil\Gamma  KG\rceil&=-{{1}\over{2}}{{d}\over{dx}}{\hbox{trace}}\, \Gamma K
G&(4.15)\cr
\lceil(I+\Gamma )^{-1}\Gamma \rceil&=-{\hbox{trace}}\bigl( (I+\Gamma
)^{-1}\Gamma'\bigr),&(4.16)\cr
\lceil K\Gamma \rceil\lceil GL\rceil&=-{{1}\over{2}}
\lceil K (\Gamma'G+\Gamma G')L\rceil,&(4.17)\cr}$$
\noindent where (4.17) is known as the product formula. The easiest way to prove  (4.15)-(4.18)
is to observe that $\Gamma 'G+\Gamma G'$ is
the integral operator with kernel 
$-2\phi_{(x)}(s)\psi_{(x)}(t)$, which has rank one. These ideas were
subsequently revived by McKean [48]. Our formulas (4.5) and (4.6) incorporate a similar idea, and are the basis of the proof of Proposition 4.2. The results we obtain appear to be more general than those of P\"oppe, and extend to periodic linear systems.\par
\indent (ii) Mumford [57] considers the ring $R_1$ of complex functions that are holomorphic on some neighbourhood of zero,  
and the ring $R_2=C_c^\infty ({\bf R}; {\bf C})$ of compactly supported smooth functions; Mulase [56] considers the differential ring $R_3={\bf C}[[x]]$ of
formal complex power series; McKean and van Moerbeke [50] consider the ring $R_4=C^\infty ({\bf T}; {\bf C})$ of smooth periodic function; then 
$R_j$ is a differential rings with respect to $D=d/dx$, and one can form the rings $R_j\{ D\}=\{ \sum_{k=-\infty}^n a_k(x)D^j: n<\infty; a_k\in R_j, -\infty ; k\leq n\}$ of 
 pseudo differential operators. In this paper we use a differential ring $({\cal A}, \ast ,\partial )$ of operators on state space.\par

\vskip.05in
\indent In the literature on inverse scattering, as in [4], the operator $R_x$ appears implicitly in various formulas, especially in the special case in which $A$ has finite rank and may be represented by a matrix. The following result extends a
special case of the Sylvester--Rosenblum theorem [8]. The formula (4.18) resembles the expressions used to obtain soliton solutions of $KdV$, as in [35, (14.12.11)].\par

\vskip.05in
\indent  For the remainder of this section, we let $A$ be a $n\times n$ 
complex matrix
with eigenvalues $\lambda_j$ $(j=1, \dots, m)$ with geometric multiplicity $n_j$ such 
that $\lambda_j+\lambda_k\neq 0$ for all $j,k\in \{1, \dots, m\}$; let  
 ${\bf K}={\bf
C}(e^{-\lambda_1t}, \dots , e^{-\lambda_mt},t)$. Also, let
$B\in {\bf C}^{n\times 1}$ and $C\in {\bf C}^{1\times
n}$.\par

\vskip.05in
\noindent {\bf Proposition 4.5.} {\sl (i) There exists a solution $R_t$ to Lyapunov's equation (1.3) with initial condition $BC$, such that the entries of $R_t$ belong to ${\bf K}$, and $\tau (t)\in {\bf K}$;}\par   
\indent {\sl (ii) $\phi\in {\bf K}$ satisfies a linear 
differential equation with constant coefficients.}\par 
\indent {\sl (iii) Suppose further that all the eigenvalues of $A$ are simple. Then there exists an invertible matrix $S$ such that $S^{-1}B= (b_j)_{j=1}^n\in {\bf C}^{n\times 1}$ and $CS=(c_j)_{j=1}^n\in {\bf C}^{1\times n}$ and the tau function is given by}
$$\eqalignno{\tau (t;\mu )=&1+ \sum_{j=1}^n
{{b_jc_je^{-2\lambda_jt}}\over{2\lambda_j}}&\cr
&+\sum_{(j,k),
(m,p): j\neq m; k\neq
p}(-1)^{j+k+m+p}{{b_jb_mc_kc_pe^{-(\lambda_j+\lambda_k+\lambda_m+\lambda
_p)t}}\over{(\lambda_j+\lambda_m)(\lambda_k+\lambda_p)}}+\dots\cr
&+\prod_{j=1}^n{{b_jc_j}\over{2\lambda_j}}\prod_{1\leq
j<k\leq n}{{(\lambda_j-\lambda_k)^2}\over{ 
(\lambda_j+\lambda_k)^2}}{e^{-2\sum_{j=1}^n\lambda_jt}}.&(4.18)}$$
\vskip.05in

\noindent {\bf Proof.}(i) By the hypothesis, we can introduce a chain of circles  $\Sigma$  that go once
round each $\lambda_j$ in the positive sense and have all the points
$-\lambda_k$ in their exterior. Then by [8], the matrix
$$R_0={{-1}\over{2\pi i}}\int_\Sigma (A+\lambda I)^{-1}BC(A-\lambda
I)^{-1}d\lambda \eqno(4.19)$$
\noindent gives a solution to the equation $-AR_0-R_0A=-BC.$ To see this, one considers $(A+\lambda I)R_0+R_0(A-\lambda I)$ and then uses the calculus of residues. \par
\indent By the Riesz functional calculus, we also have
$$e^{-tA}={{1}\over{2\pi i}}\int_\Sigma \bigl( \lambda I-A\bigr)^{-1}
e^{-t\lambda } d\lambda ;\eqno(4.20)$$
\noindent hence by Cauchy's residue theorem, there exist complex polynomials $p_j$ and $q_j$, and  integers
$m_j\geq 0$ such that 
$$e^{-tA}=\sum_{j=1}^m q_j(t)e^{-t\lambda_j}p_j(A),\eqno(4.21)$$
\noindent where $q_j(t)$ is constant if the corresponding eigenvalue
is simple. We let $R_t=e^{-tA}R_0e^{-tA}$, which gives a solution to Lyapunov's equation with initial condition $-BC$. From (4.21), we see that all the entries of $R_t$ belong to ${\bf K}$. By the Laplace expansion of the determinant, we see that all entries of $\tau (t)=\det (I+R_t)$ also belong to ${\bf K}$. \par
\indent (ii) We take $\phi (t)=Ce^{-tA}B\in {\bf K}$ by (4.21). Also, we introduce the characteristic polynomial of $(-A)$ by $\det (\lambda I+A)=\sum_{j=0}^na_j\lambda^j$. Then by the Cayley--Hamilton theorem.
$$\sum_{j=0}^na_j{{d^j\phi (t)}\over{dt^j}}=0.\eqno(4.22)$$
\indent (iii) We recall Cauchy's determinant formula. For
$x_r$ and $y_s$ complex numbers such that $x_ry_s\neq 1,$ the determinant satisfies
$$\det
\Big[{{1}\over{1-x_jy_k}}\Bigr]_{j,k=1}^n={{\prod_{1\leq
j<k\leq n}(x_j-x_k)\prod_{1\leq m<p\leq n} (y_m-y_p)}\over{\prod_{1\leq
r,s\leq n}(1-x_ry_s)}}.\eqno(4.23)$$
\indent There exists an invertible matrix $S$ such that $SAS^{-1}$ is the $n\times n$ diagonal matrix $D={\hbox{diag}}(\lambda_1, \dots, \lambda_n)$, and we observe that 
$$R_t=\Bigl[
{{b_jc_ke^{-(\lambda_j+\lambda_k)t}}\over{\lambda_j+\lambda_k}}\Bigr]_{j
,k=1}^n\eqno(4.24)$$
\noindent satisfies ${{d}\over {dt}}R_t=-[b_jc_ke^{-(\lambda_j+\lambda_k)t}]_{j,k=1}^n$ and $-DR_t-R_tD=-[b_jc_ke^{-(\lambda_j+\lambda_k)t}]_{j,k=1}^n$; so $R_t$ gives a solution of the Lyapunov equation with generator $-D$ and initial condition given by the rank-one matrix $-S^{-1}BCS=-[b_jc_k]_{j,k=1}^n$. Hence the tau function is given by  $\tau (t;\mu )=\det (I+\mu R_t)$ for this matrix, and there is an
expansion
$$\det\Bigl[\delta_{jk}+
{{\mu b_jc_ke^{-(\lambda_j+\lambda_k)x}}\over{\lambda_j+\lambda_k}}\Bigr]_{j
,k=1}^n=\sum_{\sigma\subseteq \{ 1, \dots ,n\}}\mu^{\sharp \sigma}\det
\Bigl[{{b_jc_ke^{-\lambda_jx-\lambda_kx}}\over{\lambda_j+
\lambda_k}}\Bigr]_{j,k\in \sigma}\eqno(4.25)$$
in which each
subset $\sigma$ of $\{ 1, \dots ,n\}$ of order $\sharp \sigma$, contributes a minor
indexed by $j,k\in \sigma$. Letting $x_r=\lambda_r$ and
$y_r=-1/\lambda_r$ in the Cauchy determinant formula, we obtain the identity
$$\det\Bigl[{{b_jc_ke^{-\lambda_jx-\lambda_kx}}\over{\lambda_j+
\lambda_k}}\Bigr]_{j,k\in \sigma}=\prod_{j\in
\sigma}{{b_jc_je^{-2\lambda_jx}}\over{2\lambda_j}}\prod_{j,k\in \sigma: j\neq
k}{{\lambda_j-\lambda_k}\over{\lambda_j+\lambda_k}}.\eqno(4.26)$$
\rightline{$\square$}\par
\vskip.05in

\vskip.05in

\indent In the next three sections, we give significant examples of differential rings associated with linear systems.\par
\vskip.05in
\noindent {\bf 5. The diagonal Green's function and the stationary  $KdV$ hierarchy}\par
\vskip.05in
\indent In this section, we obtain properties of ${\cal A}_\Sigma$ in terms of $A$. Thus we obtain some sufficient conditions for some differential equations to be integrable.  Throughout this section, we suppose that the hypotheses of Theorem 4.3 are in force, so the any finite set of elements of ${\cal A}_\Sigma$ are holomorphic functions on a some sector $\Omega$ containing $(x_0, \infty )$ for some $x_0\geq 0$. We do not generally require $u$ to be real valued, although in Theorem 5.2(iv) we impose this further condition so that we can compare our results with the classical spectral theory for the Schr\"odinger equation on the real line.\par
\vskip.05in
\noindent {\bf Definition}  {\sl (Stationary $KdV$ hierarchy).}  (i) Let
$f_0=1$ and $f_1=(1/2)u$. Then the 
$KdV$ recursion formula is
$$4{{d}\over{dx}}f_{m+1}(x)=4f_1(x){{d}\over{dx}}f_m(x)+4{{d}\over{dx}}
\bigl( f_1(x)f_m(x)\bigr)-{{d^3}\over{dx^3}}f_m(x).\eqno(5.1)$$ 
\indent (ii) Let $g_k=f'_k$ for $k=1, 2, \dots .$\par 
\indent (iii) If $u$
 satisfies $f_m=0$ for all $m$ greater than or equal to some $m_0$,
then $u$ satisfies the stationary $KdV$ hierarchy and is said to be an
algebro-geometric (finite gap) potential; see [29,14] . \par
\indent (iv) Suppose that $u(x)\rightarrow 0$ as $x\rightarrow\infty$, and likewise for all the partial derivatives ${\partial^\ell u}/\partial x^\ell$; suppose further that $f_j(x)\rightarrow 0$ as $x\rightarrow 0$ as $x\rightarrow\infty$ for all $j=1, 2, \dots$. Then we say that the $f_j$ are homogeneous solutions of the $KdV$ hierarchy, and write $\hat f_j$ for $f_j$ to indicate that the system of differential equations (5.1) has no arbitrary constants of integration.\par 
\vskip.05in

\noindent {\bf Proposition 5.1.}  {\sl Let ${\cal A}_\Sigma$ be as in Theorem 4.3. Then  $f_m=(-1)^m 2\lfloor A^{2m-1}\rfloor$ satisfies
the stationary $KdV$ hierarchy (Novikov's equations), since}
$$\eqalignno{4{{d}\over{d x}}\lfloor
A^{2m+3}\rfloor ={{d^3}\over{d x^3}}\lfloor A^{2m+1}\rfloor &+8\Bigl( 
{{d }\over{d x}}\lfloor A\rfloor\Bigr)\lfloor
A^{2m+1}\rfloor+16  
\lfloor A\rfloor\Bigl(  {{d}\over{d x}}\lfloor
A^{2m+1}\rfloor\Bigr).&(5.2)\cr}$$

\vskip.05in
\noindent {\bf Proof.} (i) We have the basic identities
$$\lfloor A(I-2F)A(I-2F)X\rfloor =\lfloor A^2X\rfloor-2\lfloor A\rfloor
\lfloor X\rfloor ;\eqno(5.3)$$
$$-2A(AF+FA-2FAF)=A(I-2F)A(I-2F)-A^2\eqno(5.4)$$
\noindent and their mirror images. Hence
$${{d}\over{dx}}\lfloor A^{2m+1}\rfloor =\lfloor A(I-2F)A^{2m+1}+
A^{2m+1}(I-2F)A\rfloor ,\eqno(5.5)$$
\noindent so
$$\eqalignno{{{d^2}\over{dx^2}}\lfloor  A^{2m+1}\rfloor
&=\lfloor
A(I-2F)A(I-2F)A^{2m+1}
+2A(I-2F)A^{2m+1}(I-2F)A\cr
&\quad +A^{2m+1}(I-2F)A(I-2F)A\cr
&\quad -2A(AF+AF-2FAF)A^{2m+1}-2A^{2m+1}(AF+FA-2FAF)A\rfloor\cr
&= \lfloor A(I-2F)A(I-2F)A^{2m+1}+2A(I-2F)A^{2m+1}(I-2F)A\cr
&\quad +A^{2m+1}(I-2F)A(I-2F)A \cr
&\quad +A(I-2F)A(I-2F)A^{2m+1}-A^{2m+1}+A^{2m+3}(I-2F)A(I-2F)A-A^{2m+3}\rfloor
\cr
&=2\lfloor A(I-2F)A^{2m+1}(I-2F)A\rfloor -2\lfloor
A^{2m+3}\rfloor\cr
&= 2\lfloor A(I-2F)A(I-2F)A^{2m+1}\rfloor -2\lfloor A^{2m+3}\rfloor\cr
&\quad+\lfloor A(I-2F)A(I-2F)A^{2m+1}\rfloor+2\lfloor
A^{2m+1}(I-2F)A(I-2F)A\rfloor\cr
&= 2\lfloor A(I-2F)A^{2m+1}(I-2F)A\rfloor +2\lfloor A^{2m+3}\rfloor
\cr
&\quad -4\lfloor A^{2m+1}\rfloor \lfloor A\rfloor-4\lfloor A\rfloor\lfloor
A^{2m+1}\rfloor .&(5.6)\cr}$$
Now we differentiate the first summand of the final term
$$\eqalignno{ {{d}\over{dx}}2\lfloor
&A(I-2F)A^{2m+1}(I-2F)A\rfloor\cr
&= 
2\lfloor A(I-2F)A(I-2F)A^{2m+1}(I-2F)A\rfloor+2\lfloor
A(I-2F)A^{2m+1}(I-2F)A(I-2F)A\rfloor\cr
&\quad -4\lfloor A(AF+FA-2FAF)A^{2m+1}(I-2F)A\rfloor-4\lfloor
A(I-2F)A^{2m+1}(AF+FA-2FAF)A\rfloor\cr
&=2\lfloor A(I-2F)A(I-2F)A^{2m+1}(I-2F)A\rfloor+2\lfloor
A(I-2F)A^{2m+1}(I-2F)A(I-2F)A\rfloor\cr
&\quad +2\lfloor A(I-2F)A(I-2F)A^{2m+1}(I-2F)A \rfloor-2\lfloor
A^{2m+3}(I-2F)A\rfloor\cr
&\quad+2\lfloor A(I-2F)A^{2m+1}(I-2F)A(I-2F)A\rfloor-2 \lfloor
A(I-2F)A^{2m+3}\rfloor\cr
&= 4\lfloor A(I-2F)A(I-2F)A^{2m+1}(I-2F)A\rfloor
+4\lfloor A(I-2F)A^{2m+1}(I-2F)A(I-2F)A\rfloor\cr
&\quad-2\lfloor A(I-2F)A^{2m+3}+A^{2m+3}(I-2F)A\rfloor\cr
&=-8\lfloor A\rfloor\lfloor A^{2m+1}(I-2F)A\rfloor+4\lfloor
A^{2m+3}(I-2F)A\rfloor\cr
&\quad -8\lfloor A\rfloor\lfloor A(I-2F)A^{2m+1}\rfloor+4\lfloor
A(I-2F)A^{2m+3}\rfloor-2{{d}\over{dx}}\lfloor A^{2m+3}\rfloor\cr
&= -8\lfloor A\rfloor\lfloor A(I-2F)A^{2m+1}+A^{2m+1}(I-2F)A \rfloor\cr
&\quad +4\lfloor A(I-2F)A^{2m+3}+A^{2m+3}(I-2F)A\rfloor-2{{d}\over{dx}}
\lfloor A^{2m+3}\rfloor\cr
&=-8\lfloor A\rfloor{{d}\over{dx}}\lfloor
A^{2m+1}\rfloor +2{{d}\over{dx}}\lfloor A^{2m+3}\rfloor;&(5.7)\cr}$$
\noindent hence
$${{d^3}\over{dx^3}}\lfloor A^{2m+1}\rfloor=-8\lfloor
A\rfloor{{d}\over{dx}}\lfloor A^{2m+1}\rfloor +4{{d}\over{dx}}\lfloor
A^{2m+3}\rfloor-8{{d}\over{dx}}\Bigl( \lfloor A\rfloor\lfloor
A^{2m+1}\rfloor\Bigr);\eqno(5.8)$$
\noindent which gives the stated result.\par
\rightline{$\square$}\par
\vskip.05in
\noindent {\bf Definition} {\sl (Diagonal Greens function).} Let $(-A,B,C)$ be as in Theorem 2.2. Then the diagonal Green's function is $g_0(x;\zeta )/\sqrt{\zeta}$ where 
$$g_0(x;\zeta )=(1/2)+ \lfloor A(\zeta I-A^2)^{-1}\rfloor.\eqno(5.9)$$
\indent The notation $g_0(x;\zeta )$ is chosen to indicate a generating function and also the diagonal of a Green's function; now we explain the latter connection.  Let ${\bf C}_+=\{\lambda \in {\bf C}: \Im \lambda >0\}$ be the open upper half plane. \par

\vskip.05in
\noindent {\bf Theorem 5.2.} {\sl Let $(-A,B,C)$ be as in Theorem 2.2.\par 
\indent (i)  Then $g_0(x;\zeta )$ is bounded and continuously differentiable in $x$ and has a unique  asymptotic expansion depending on the odd powers of $A$, \par
$$g_0(x; \zeta )\asymp {{1}\over{2}}+{{\lfloor A\rfloor}\over{\zeta}}+ {{\lfloor A^3\rfloor}\over{\zeta^2}}+{{\lfloor A^5\rfloor}\over{\zeta^3}}+\dots \qquad (\zeta\rightarrow-\infty );\eqno(5.10)$$
\indent (ii) $g_0(x;\zeta )$ satisfies Drach's equation}
$${{d^3g_0}\over{dx^3}} = 4(u+\zeta ){{d g_0}\over{dx}}+2{{du}\over{dx}}g_0\qquad (x>x_0; -\zeta >\omega );\eqno(5.11)$$ 
 \indent {\sl (iii) there exists $x_1>0$ such that 
$$\psi_{\pm} (x,\zeta  ) =\sqrt{g_0(x, -\zeta )}\exp\Bigl( \mp \sqrt{-\zeta }\int_{x_1}^x {{dy}\over{2g_0(y; -\zeta )}}\Bigr)\eqno(5.12)$$
\noindent satisfies} 
$$-\psi_{\pm}''(x; \zeta )+u(x)\psi_{\pm} (x, \zeta )=\zeta \psi_\pm (x; \zeta )\qquad (x>x_1, \zeta >\omega ).\eqno(5.13)$$
\indent  {\sl (iv) Suppose that $u$ is a continuous real function that is bounded below,
 and that $\psi_{\pm}$ from (iii) satisfy $\psi_+(x; \zeta )\in L^2((0, \infty ); {\bf C})$ and $\psi_-(x; \zeta )\in L^2((-\infty,0 ); {\bf C})$ for all $\zeta\in {\bf C}_+$. Then 
$L=-{{d^2}\over{dx^2}}+u(x)$ defines an essentially
 self-adjoint operator in $L^2({\bf R}; {\bf C})$, and  the Greens function $G(x,y; \zeta )$ which represents $(\zeta I-L)^{-1}$ has a diagonal that satisfies}
$$G(x,x; \zeta )={{g_0(x;-\zeta )}\over{\sqrt{-\zeta }}}.\eqno(5.14)$$ 

\vskip.05in
\noindent {\bf Proof.}  (i) Let $\pi -\theta <\arg \lambda <\theta$, so $\lambda$ and $-\lambda$ both lie in $S_\theta$, hence  $\zeta =\lambda^2$ satisfies $2\pi -2\theta <\arg \zeta <2\theta $, so $\zeta$ lies close to $(-\infty , 0)$. Then $\zeta I-A^2$ is invertible and $\vert \zeta\vert \Vert (\zeta I-A^2)^{-1}\Vert_{{\cal L}(H)} \leq M$.  The function
$$g_0(x; \zeta )={{1}\over{2}}+ Ce^{-xA}(I+R_x)^{-1}A(\zeta I-A^2)^{-1}(I+R_x)^{-1}e^{-xA}B\qquad (x>0)\eqno(5.15)$$
\noindent is well defined by Theorem 2.2(iii). \par
\indent To obtain the asymptotic expansion,  we note that $e^{-xA}(I+R_x)^{-1}$ and $(I+R_x)e^{-xA}$ involve the factor $e^{-xA}$, where $(e^{-zA})$ is a holomorphic semigroup on $S_{\theta -\pi /2}$. Hence $A^{2j+1}e^{-xA}\in {\cal L}(H)$ and by Cauchy's estimates  there exist $x_0, M_0>0$ such that  $\Vert A^{2j+1}e^{-xA}\Vert_{{\cal L}(H)}\leq M_0(2j+1)!$ for all $x\geq x_0>0$. As in Proposition 3.1, we have an asymptotic expansion of 
$$\eqalignno{e^{-zA}\bigl( (\lambda I -A)^{-1}-(\lambda I+A)^{-1}\bigr)&=-e^{-zA}\int_0^\infty e^{\lambda s} e^{-sA}\, ds-e^{-zA}\int_0^\infty e^{-\lambda s} e^{-sA}\, ds\cr
&=e^{-zA} \Bigl({{A}\over{\lambda^2}}+{{A^3}\over{\lambda^4}}+\dots +{{A^{2j-1}}\over{\lambda^{2j}}}\Bigr)\cr
&\quad +{{e^{-zA}}\over{\lambda^{2j+1}}}\int_0^\infty A^{2j+1}e^{-sA} (e^{s\lambda}-e^{-\lambda s})\, ds,&(5.16)\cr}$$
\noindent  in which all the summands are in ${\cal L}(H)$ due to the factor $e^{-zA}$ for  $z\in S_{\theta -\pi /2}$. Hence 
$$Ce^{-xA}e^{-xA}(I+R_x)^{-1}\int_0^\infty A^{2j+1}e^{-sA} (e^{s\lambda }-e^{-s\lambda })\, ds (I+R_x)^{-1}e^{-xA}B\rightarrow 0\qquad (x>0)\eqno(5.17)$$
\noindent as $\lambda \rightarrow i\infty$, or equivalently $\zeta\rightarrow-\infty$, so
$$g_0(x, \zeta )={{1}\over{2}}+ Ce^{-xA}(I+R_x)^{-1} \Bigl( {{ A}\over{\zeta}}+ {{A^3}\over{\zeta^2}}+\dots +{{A^{2j-1}}\over{\zeta^{j}}}\Bigr) (I+R_x)^{-1}e^{-xA}B+O\Bigl({{1}\over{\zeta^{j+1}}}\Bigr).\eqno(5.18)$$ 
\noindent This gives the asymptotic series; generally, the series is not convergent since the implied constants in the term $O(\zeta^{-(j+1)})$ involve $(2j+1)!$.\par
\indent (ii) From Proposition 5.1 we have
$$\eqalignno{4{{d}\over{dx}}\sum_{m=0}^\infty  {{\lfloor A^{2m+3}\rfloor}\over{\zeta^{m+1}}}&={{d^3}\over{d x^3}}\sum_{m=0}^\infty  {{\lfloor A^{2m+1}\rfloor}\over{\zeta^{m+1}}}\cr
&\quad+8\Bigl( {{d}\over{d x}}\lfloor A\rfloor \Bigr)\sum_{m=0}^\infty  {{\lfloor A^{2m+1}\rfloor}\over{\zeta^{m+1}}}+16\lfloor A\rfloor 
{{d}\over{d x}}\sum_{m=0}^\infty  {{\lfloor A^{2m+1}\rfloor}\over{\zeta^{m+1}}};&(5.19)\cr}$$
\noindent the required result follows on rearranging.\par
\indent Conversely, suppose that $g_0$ as defined in (5.9) has an asymptotic expansion with coefficients in $C^\infty ((0, \infty ); {\bf C})$ as $\zeta\rightarrow-\infty$ and that $g_0(x; \zeta )$ satisfies (5.11). Then the coefficients of $\zeta^{-j}$ satisfy a recurrence relation which is equivalent to the systems of differential equations (5.1).\par
\indent  The asymptotic expansion is unique in the following sense. Suppose momentarily that $t\mapsto \lfloor Ae^{-tA^2}\rfloor$ is bounded and repeatedly differentiable on $(0, \infty )$, with  $M, \omega >0$ such that $\vert \lfloor Ae^{-tA^2}\rfloor \vert \leq Me^{\omega t}$ for $t>0$, and that there is a Maclaurin expansion 
$$\lfloor Ae^{-tA^2}\rfloor=\lfloor A\rfloor -\lfloor A^3\rfloor t+{{\lfloor A^5\rfloor t^2}\over {2!}}-\dots +O(t^k)\eqno(5.20)$$
\noindent on some neighbourhood of $0+$. Then by Watson's Lemma, the integral $\int_0^\infty \lfloor Ae^{-tA^2}\rfloor e^{t\zeta }\, dt$ has an asymptotic expansion as $\zeta\rightarrow-\infty$, where the coefficients give the formula (5.10). \par

\indent  (iii)  Since $(e^{-tA})_{t>0}$ is a contraction semigroup on $H$, we have ${\cal D}(A^2)\subseteq {\cal D}(A)$ and $\Vert Af\Vert_H^2\leq 2\Vert A^2f\Vert_H\Vert f\Vert_H$ for  all $f\in {\cal D}(A^2)$ by the Hardy-Littlewood-Landau inequality, so $\Vert \zeta f+A^2 f\Vert_H\geq \sqrt{\zeta}\Vert Af\Vert_H$ for $\zeta >0$. We deduce that $A^2-2A+\zeta I$ is invertible for $\zeta >9$ and generally for all $\zeta\in {\bf C}$ such that $\Re \zeta $ is sufficiently large. By Proposition 5.1 and the multiplicative property of the bracket, we have
$${{1}\over {2g_0(x; -\zeta )}}=1+\bigl\lfloor 2A(\zeta I+A^2-2A)^{-1}\bigr\rfloor,\eqno(5.21)$$
\noindent and we observe that $g_0(x; -\zeta )\rightarrow 1/2$ as $x\rightarrow\infty$, so there exists $x_1>0$ such that $g_0(x, -\zeta )>0$ for all $x>x_1$ and the differential equation integrates to 
$$g_0{{d^2g_0}\over{dx^2}}-{{1}\over{2}}\Bigl( {{dg_0}\over{dx}}\Bigr)^2 =2(u-\zeta )g_0^2 +{{\zeta}\over{2}}.\eqno(5.22)$$
\noindent So we can define $\psi (x; \zeta  )$ as in (5.13), and then one verifies the differential equation for $\psi (x; \zeta )$ by using (5.22). \par
\indent (iv) By a theorem of Weyl [32, 10.1.4], $L$ is of limit point type at 
$\pm\infty$, and there exist nontrivial solutions 
$\psi_{\pm}(x; \zeta )$ to $-\psi_{\pm}''(x; \zeta )+u(x)\psi_{\pm} (x;\zeta)=\zeta \psi_{\pm} (x;\zeta)$ such that 
$\psi_+(x; \zeta )\in L^2(0, \infty )$ and  $\psi_-(x; \zeta )\in L^2(-\infty ,0)$, and these are unique up to constant multiples. Also the inverse operator $(-\zeta I+L)^{-1}$ may be represented as an integral operator in $L^2( {\bf R}; {\bf C})$ with kernel $G(x,y;\zeta )$, which has diagonal 
$$G(x,x;\zeta )={{\psi_+(x; \zeta )\psi_-(x; \zeta )}\over {{\hbox{Wr}}(\psi_+(\, ; \zeta ), \psi_-(\, ; \zeta ))}}\qquad (\Im\zeta >0),\eqno(5.23)$$
\noindent  Given $\psi_{\mp}$  as in (iii), we can compute  $\psi_+(x; \zeta )\psi_-(x; \zeta )=g_0(x; -\zeta )$ and their Wronskian is
${\hbox{Wr}}(\psi_+, \psi_-)=\sqrt{-\zeta },$ hence the result.\par
\rightline{$\square$}\par

\vskip.05in
\noindent {\bf Remarks 5.3} (i) The importance of the diagonal Greens function is emphasized in [30]. Gesztesy and  Holden [29, Lemma 1.6.1] obtain an asymptotic expansion of the diagonal
$G(x,x;\zeta )$ which is consistent with Theorem 5.2(i). Under conditions discussed in (5.48), we have similar asymptotics as $-\zeta\rightarrow \infty$.\par
\indent (ii) Let $\psi$ and $\varphi$ be solutions of $-f''+uf=\zeta f$ for some $u\in C^2{(\bf R}; {\bf C})$ such that ${\hbox{Wr}}(\psi, \varphi )^2=-\zeta$, and then let $g(x)=\psi (x)\varphi (x)$. By differentiating, one checks that $gg''-(g')^2/2-2(u-\zeta )g^2=C$ for some constant $C$, and one can evaluate $C=-{\hbox{Wr}}(\psi, \varphi )^2/2$. Hence 
$$\rho (x) =\sqrt{g(x)} \exp \Bigl( \pm \sqrt{-\zeta}\int_0^x {{dy}\over {2g(y)}}\Bigr)\eqno(5.24)$$
\noindent gives a solution of $-\rho''+u\rho =\zeta \rho$. \par
\indent (iii) Drach observed that one can start with the differential equation (5.11), and produce the solutions (5.24); see [14].  He showed that Schr\"odinger's equation is integrable by quadratures, if and only if (5.11) can be integrated by quadratures for typical values of $\zeta$, and  Brezhnev translated his results into the  modern theory of finite gap integration [14]. Having established integrability of Schr\"odinger's equation by quadratures, one can introduce the hyperelliptic spectral curve ${\cal E}$ with $g<\infty$ and proceed to express the solution in terms of  the Baker--Akhiezer function. Hence one can integrate the equation and express the solution in terms of the Riemann's theta function on the Jacobian of ${\cal E}$, as in [43]. Our presentation follows Drach's; one advantage is that we can deal with rather degenerate solutions of differential equations, such as occur in the theory of solitons, and can evade technicalities regarding special points on Jacobians of compact Riemann surfaces. \par

\vskip.05in

\noindent {\bf Proposition 5.4.} {\sl (i) In the context of Theorem 4.3, suppose that there exists a non-zero odd complex polynomial $p_0(X)$ such that $\lfloor p_0(A)\rfloor =0.$ 
Then ${\bf C}[u, {{d u}/{d x}}, \dots,  ]$ is a Noetherian differential ring for
$ {{d }/{d x}} $ and the usual multiplication.\par
\indent (ii) In particular, (i) holds when $A^2$ is algebraic.}\par
\vskip.05in
\noindent {\bf Proof.} (i)  It follows from Proposition 5.1 that $\lfloor A^{2j-1}\rfloor =c_j u^{(2j)}+P_j(u, u', \dots, u^{(2j-1)})$ where $P$ is a complex polynomial, and $c_j\neq 0$. Adding multiples of such identities, and using the hypotheses, we deduce that there exists $m$ such that
$${{d^{2m}u}\over{dx^{2m}}} =Q_{2m}\Bigl( u, {{d u}\over {dx}}, \dots,  {{d^{2m-1} u}\over{d x^{2m-1}}}\Bigr),\eqno(5.25)$$
\noindent where $Q_{2m}$ is a complex polynomial which is determined by $p_0$ and the $f_j$. By repeatedly differentiating this identity, and substituting back, one can obtain polynomials $Q_n$ such that 
$u^{(n)}=Q_{n}(u, {{d u}/{dx}}, \dots,  {{d^{n-1} u}/{d x^{n-1}}}),$ for all $n\geq 2m$. Hence  ${\bf C}[u, {{d u}/{d x}}, \dots,  {{d^{2m-1} u}/{d x^{2m-1}}}]$ gives all of 
${\bf C}[u, {{d u}/{d x}}, \dots,  ]$.\par
\indent  Since $u$ is meromorphic on the domain $\Omega$, the algebra ${\bf C}[u, {{d u}/{d x}}, \dots,  {{d^{2m-1} u}/{\partial x^{2m-1}}}]$ is an integral domain, and we have shown it to be closed under differentiation. By mapping $X_j\mapsto {{d^j u}/{d x^j}}$ for $j=0, \dots , 2m-1$,  we obtain a short exact sequence of algebra homomorphisms 
$$0\rightarrow J\rightarrow {\bf C}[X_0, \dots, X_{2m-1}]\rightarrow  {\bf C}\Bigl[ u, {{d u}\over {dx}}, \dots,  {{d^{2m-1} u}\over{d x^{2m-1}}}\Bigr]\rightarrow 0,\eqno(5.26)$$
\noindent where the  ideal $J$ is prime. Hence   ${\bf C}[u, {{d u}/{d x}}, \dots,  {{d^{2m-1} u}/{d x^{2m-1}}}]$ is finitely generated as an algebra. Also, the ${\bf C}[X_0, \dots, X_{2m-1}]/J$ is naturally isomorphic as an algebra to the co ordinate ring ${\bf C}[V]$, where $V$ is the affine variety $\{ z\in {\bf C}^{2m}: f(z)=0, \forall f\in J\}.$\par
\indent (ii) If $A^2$ is algebraic, then there exists a monic complex polynomial such that $p(A^2)=0$, hence $Ap(A^2)=0$ and (i) applies.\par
\rightline{$\square$}
\vskip.05in

\noindent {\bf Proposition 5.5.} {\sl Suppose that $\lfloor A^{2n+1}\rfloor =0$ for some $n\geq 0$.\par
\indent (i)  Then 
$u=-4\lfloor A\rfloor$ is
finite gap;\par
\indent (ii) $-d^2/dx^2+u$ commutes with a differential operator of odd order;\par
 \indent (iii) Schr\"odinger's equation can be integrated by
quadratures over ${\bf C}[u, {{d u}/{d x}}, \dots,  ]$.}\par 
\indent {\sl (iv) The image of  $(A, \partial A, \dots, \partial^{2n}A)\in {\cal A}_\Sigma^{2n+1}$ under $\lfloor\, \cdot\, \rfloor$ equals $(u, u', \dots, u^{(2n)}),$ where 
$(u. u', \dots, u^{(2n)})$ satisfy a system of polynomial equations which determine a complex algebraic variety.}\par 
\vskip.05in
\noindent {\bf Proof.} (i) Note that the hypothesis does not change if we apply the translation $(-A, B,C)\mapsto (-A, e^{-tA}B, C)$ so that $u(x)\mapsto u(x+t)$.Then $u=-4\lfloor A\rfloor$ satisfies
$$\eqalignno{f_1&=(1/2)u+(1/2)c_1=-2\lfloor A\rfloor +(1/2)c_1;\cr
f_2&=-(1/8)u''+(3/8)u^2+(c_1/4)u+(c_2/2)-c_1/8=2\lfloor A^3\rfloor-c_1\lfloor A\rfloor +
 (c_2/2)-c_1/8, \dots .&(5.27)}$$
\noindent  Under the hypotheses of Theorem 4.3, all of the $c_j$ vanish since $u$ and its derivatives converge to zero as $x\rightarrow\infty$, so the $f_j$ are homogeneous solutions of $KdV$ hierarchy.
 By Proposition 5.1, $f_m(x)=0$ for all $m>n$.\par
\indent (ii) With $L=-d^2/dx^2+u$, we have an ordinary differential operator
$$P_{2n+1}=\sum_{j=0}^n \Bigl(f_{n-j}
{{d}\over{dx}}-2^{-1}{{df_{n-j}}\over{dx}}\Bigr) (-L)^j\eqno(5.28)$$
\noindent so that 
$$[L, P_{2n+1}]=2f_{n+1}'.\eqno(5.29)$$
\indent (iii)  By Propositions 5.1 (ii) and 5.4, we have a Noetherian differential ring
 $${\bf C}[u, {{d u}/{d x}}, \dots, ; \zeta; g,g,g'' ].\eqno(5.30)$$
\noindent We aim to prove that we can obtain all elements of this ring by 
quadratures over ${\bf C}[u, {{d u}/{d x}}, \dots, ]$, and thereby solve Schr\"odinger's equation. We observe that $2\zeta^ng_0(x;\zeta )$ is a monic polynomial of degree $n$ in $\zeta$. Changing notation, we introduce
$F(x;\lambda)=\sum_{\ell =0}^n f_{n-\ell}(x)\lambda^\ell$ where the coefficients of $F$, as a monic polynomial of degree $n$  in $\lambda$, are polynomials in $u$ and its spatial derivatives, since $f_j\in {\cal A}$.  Hence
$$p(\lambda)={{1}\over{2}} F(x;\lambda){{\partial^2 F(x;\lambda)}\over{\partial x^2}}-
{{1}\over{4}}\Bigl({{\partial F(x;\lambda) }\over{\partial x}}\Bigr)^2
-(u(x)-\lambda)F(x;\lambda)^2\eqno(5.31)$$
\noindent is a polynomial of degree $2n+1$ in $\lambda$, with coefficients in ${\cal A}$, which is actually independent
of $x$, so $p(\lambda )\in {\bf C}[\lambda ]$. Now we introduce  
$$\psi_{\pm }(x;\lambda )=\sqrt{F(x;\lambda)}\exp\Bigl( \pm
\mu\int^x{{d\xi}\over{F(\xi ;\lambda)}}\Bigr),\eqno(5.32)$$
\noindent which satisfies 
$${{\psi_\lambda''}\over{\psi_\lambda}}={{1}\over{F(x;\lambda )^2}}\bigl( (u-\lambda )F(x;\lambda )^2+p(\lambda )+\mu^2\bigr);\eqno(5.33)$$
\noindent hence $\psi_\lambda$ gives the solution to Schr\"odinger's equation
$\psi''_{\pm }(x;\lambda )=(u-\lambda )\psi''_{\pm }(x;\lambda )$ when $(\lambda, \mu )$ lies on
$${\cal E}=\{ (\lambda,\mu ): \mu^2=-p(\lambda)\}.\eqno(5.34)$$
\indent (iv) This follows from an argument of Mumford [57]. We introduce the polynomials in $\lambda$:
$$U(\lambda ;x)=F; \quad W(\lambda ;x)=(\lambda -u(x))F+{{1}\over{2}}{{\partial^2 F}\over {\partial x^2}};\quad V(\lambda ;x)= {{i}\over{2}}{{\partial F}\over {\partial x}}\eqno(5.35)$$
\noindent which have degrees $n$, $n+1$ and $n-1$ respectively, for typical $x$, and the coefficients are given by the $f_j$ and their partial derivatives with respect to $x$. Then we have a monic polynomial
$$p(\lambda )=U(\lambda ;x)W(\lambda ;x)+V(\lambda ;x)^2\eqno(5.36)$$
\noindent of degree $2n+1$, where the leading terms are 
$$U=\lambda^n +\lambda^{n-1}f_1+\dots ,\quad V=(i/2) f_1'\lambda^{n-1}+\dots , \quad W=\lambda^{n+1}+(f_1-u)\lambda^n+\dots .\eqno(5.37)$$
\noindent Then ${\cal E}$ gives a hyperelliptic curve over ${\bf C}$  of genus $g\leq n$. The Picard group ${\hbox{Pic}}({\cal E})$ is the set of divisors on ${\cal E}$, modulo linear equivalence as in [65]. Then the Jacobian $J$ is the subgroup of ${\hbox{Pic}}({\cal E})$ consisting of equivalence classes that have degree zero; this group $J$ gives an Abelian variety which may be determined algebraically by arguments presented in [57]. There is a Riemann theta function $\vartheta$ which is entire and quasi-periodic  on $J$, and has zero set  $\theta_0=\{ z\in J: \vartheta (z)=0 \}$. When we replace $u(x)$ by translation to $u(x+t)$, we make a corresponding automorphism $p(\lambda )\mapsto p(\lambda ;t)$, $U(\lambda ;x)\mapsto U(\lambda ;x,t)$, etc. One can then differentiate with respect to the parameter $t$, and thus introduce a vector field $D_\infty$ on $J$.\par
\indent One can also introduce a vector field $D_\theta$ on $J$ and the meromorphic function $\wp (z)=D^2_\infty \log \vartheta $ on $J\setminus\theta_0$. Comparing our calculations with those of Mumford, we obtain  $D_\infty =(i/2){{\partial}/{\partial x}}$  and $4\wp +d=u/2$ for some constant $d$. Thus we can regard $x$ as a coordinate for motion along the vector field determined by $D_\infty$ in the tangent space to $J$, and $u=8\wp +2d$ as an extension of $u$ to a function on $J\setminus \theta_0$. Likewise, the functions $f_j, g_k$ of (5.1) may be extended to complex functions on $J\setminus \theta_0$ give rise to an embedding $J\setminus\theta_0\rightarrow {\bf C}^{2g+1}$ by $z\mapsto (\wp, D_\infty\wp, \dots, D_\infty^{2n}\wp ).$ By analogy, in the context of Theorem 5.3 we can regard $(A, \partial A, \dots, \partial^{2n}A)\in {\cal A}_\Sigma^{2n+1}$ as coordinates for the Jacobian. See [49, article 9].\par
\rightline{$\square$}\par
\vskip.05in

\noindent {\bf Example 5.6.} Suppose that $A$ satisfies $\lfloor a_3A^3+a_1A\rfloor =0$, which is the first non trivial case of Proposition 5.4. Then $a_3f_2-a_1f_1=0$, so
$$a_3\bigl( -u''/16 +3u^2/16+c_1u/8+c_2/4-c_1/16\bigr) -(a_1/4)u-(a_1/4)c_1=0.\eqno(5.38)$$
\noindent Hence we can identify $V$ with a curve $X_1^2=p(X_0)$ where $p$ is a cubic; see [49].  In the context of Theorem 4.3, we have $u(x)\rightarrow 0$ as $x\rightarrow\infty$ along with $u', u''$, so the solution has the form
$$u(x)={{2a_1}\over{a_3}}{\hbox{sech}}^2\sqrt{{-a_1}\over{4a_3}} (x+\gamma ),\eqno(5.39)$$
\noindent with constant $\gamma\in {\bf C}$, as  is familiar from the theory of solitons. However, (5.37) also has solutions in terms of Weierstrass's elliptic function $\wp$, and in section 7 we construct linear systems with potential $\wp$.\par
\indent (ii) For $g>1$, does there exist a linear system such that the corresponding potential is the $\wp$ function for a hyperelliptic curve of genus $g$? We obtain a partial solution of this in section 7, by using the Schottky--Klein prime function.\par
\vskip.05in

\noindent {\bf Proposition 5.7.} {\sl Suppose that $g_0(x_j, \zeta_j)=0$ and ${{\partial g_0}\over{\partial \zeta }}(x_j, \zeta_j)\neq 0$  for some $x_j>x_0$ and nonzero}  $\zeta_j\in {\bf C}\setminus {\hbox{Spec}}(A^2)$. \par
\indent {\sl (i) Then there exist $\varepsilon >0$ and a differentiable family of solutions of (5.13)  which are\par
\noindent parametrized by an arc $\{ \mu_j(t) : x_j-\varepsilon <t<x_j+\varepsilon )\}$ passing through $\zeta_j$ such that $\psi (t,\mu_j(t))=0$;\par
\indent (ii)  Dubrovin's equation holds}  
$${{d\mu_j}\over{dx}}(t)={{\pm\sqrt{\zeta}}\over{ {{\partial g_0}\over{\partial \zeta }}(t, \mu_j(t))}}.\eqno(5.40)$$
\vskip.05in
\noindent {\bf Proof.}  (i) First we make an observation about the zeros of $g_0(x, \zeta )$.  The exponential matrix 
$$\exp \Bigl( t\left[\matrix{0&1&0\cr 0&0&1\cr 0&4\zeta &0\cr}\right]\Bigr)\eqno(5.41)$$
\noindent has entries that are entire functions of $\zeta$ of order $\rho$, where $\rho \leq 1/2$. We deduce that the general solution $g(x,\zeta )$ of (5.11) may be written as an entire function of $\zeta$ of order $\rho\leq 1/2$, so there exist functions $g_1(x)$ and $\mu_j(x)$ such that 
$$g(x, \zeta )=g_1(x)(\zeta -\mu_0(x))\prod_{j=1}^\infty\Bigl( 1-{{\zeta}\over{\mu_j(x)}}\Bigr).\eqno(5.42)$$ 
\noindent We can, of course, divide this by any function of $\zeta$ and still have a solution of the linear differential equation; in particular, we obtain $g_0(x, \zeta )$ in this way. \par
\indent Let $x_j>x_0$ and suppose that $g_0(x_j, \zeta_j)=0$ where $\zeta_j\in {\bf C}\setminus {\hbox{Spec}}(A^2)$ so that $\zeta \mapsto g_0(x_j, \zeta )$ is holomorphic. By hypothesis,   we have ${{\partial g_0}\over{\partial \zeta }}(x_j, \zeta_j)\neq 0$, which rules out the possibility of a multiple zero at $\zeta_j$ in the factorization (). Now we apply the implicit function theorem to the formula $g_0(t, \zeta )=0$, noting that $\zeta\mapsto g_0(x_j, \zeta )$ has a simple zero at $\zeta=\zeta_j$. By Rouch\'e's theorem and the calculus of residues, there exist $\varepsilon_1, \varepsilon_2>0$ such that the contour integral 
$$\mu_j(t)={{1}\over{2\pi i}}\int_{C(\zeta_j, \varepsilon_1)} {{z{{\partial g_0}\over{\partial \zeta }}(t,z) dz}\over{ g_0(t,z)}}\eqno(5.43)$$
\noindent determines a continuously differentiable  function $\mu_j: (x_j-\varepsilon_2, x_j+\varepsilon_2 )\rightarrow {\bf C}$ which satisfies $\mu_j(x_j)=\zeta_j$ and 
$${{\partial g_0}\over{\partial x }}(t, \mu_j(t))+{{\partial g_0}\over{\partial \zeta }}(t, \mu_j(t)) {{d\mu_j}\over {dx}}(t)=0.\eqno(5.44)$$   
\indent (ii) By (5.22), we have ${{\partial g_0}\over{\partial x }}(x_j, \zeta_j)^2=\zeta\neq 0$  hence Dubrovin's equation holds. Also $d\mu_j/dx\neq 0$, so $\mu_j$ determines a differentiable arc in ${\bf C}\setminus {\hbox{Spec}}(A^2)$. Then for $\zeta=\mu_j(t)$ on this arc,  $\psi_+(x, \mu_j(t))\psi_-(x,\mu_j(t))=g_0(x, \mu_j(t))$ vanishes at $x=t$, so one of the solutions $\psi_{\pm }(x;\mu_j(t) )$ satisfies $\psi_{\pm }(t;\mu_j(t) ))=0$ for all $t\in (x_j-\varepsilon_3,x_j+\varepsilon_3)$ for some $0<\varepsilon_3<\varepsilon_1$.\par
\indent If $\zeta_j\in {\bf R}$ and $g_0(x,\zeta )$ is real-valued for all $(x,\zeta )\in {\bf R}^2$ in a neighbourhood of $(x_j, \zeta_j)$, then we can choose the arc $\{ \mu_j (t):  t\in (x_j-\varepsilon_3,x_j+\varepsilon_3)\}$ to be an open  subinterval of  ${\bf R}$.\par

\rightline{$\square$}\par
\vskip.05in
\noindent {\bf Definition} {\sl (Characteristic function).} For $(-A,B,C)$ as in Theorem 2.2 let
$$\varphi (\lambda ;x)=\int_0^\lambda {{ig_0(x; -\zeta )}\over{\sqrt{-\zeta }}}\, d\zeta ,\eqno(5.45)$$
\noindent and define the characteristic function by
$$\Delta (\lambda ; x)=2\cos \varphi (\lambda ;x)\qquad (2\pi -2\theta <\arg \lambda <2\theta ).\eqno(5.46)$$
\vskip.05in
\indent  Gesztesy and Simon [30] have developed another approach to $G(x,x;\zeta )$ for self-adjoint Schr\"odinger operators which involves the xi function and Krein's spectral shift; see also [34]. 
We adopt their terminology from [30, (6.8)] in the following context. To understand the sign conventions, it helps to bear in mind that, for $L$ self-adjoint,  the diagonal Greens function $G(x,x;\zeta )$ is purely imaginary for $\zeta$ in the essential spectrum of $L$.  \par
\vskip.05in 
\indent Suppose that $(\lambda_j)$ is a non-zero complex sequence such that \par
\indent (i) $(\Re \lambda_j)_{j=0}^\infty$ is increasing;\par
\indent (ii) $\sum_{j=0}^\infty 1/(1+\vert \Re\lambda_j\vert^{\alpha})$ converges for all $\alpha>1/2$;\par
\indent (iii) $\sum_{n=1}^\infty \vert\lambda_{2n-1}-\lambda_{2n}\vert$ converges.\par
\noindent By (ii), we can introduce the functions 
$$F_e(\lambda )=\prod_{j=0}^\infty\Bigl( 1-{{\lambda}\over{\lambda_{2j}}}\Bigr),  \quad F_o(\lambda )=\prod_{j=1}^\infty\Bigl( 1-{{\lambda}\over{\lambda_{2j-1}}}\Bigr),\eqno(5.47)$$
\noindent which are entire and of order $\rho\leq 1/2$. Then $F_0(\lambda )/F_e(\lambda )$ is meromorphic with only simple zeros and poles, so 
$${\cal E}=\Bigl\{ (z,\lambda ):z^2={{F_o(\lambda )}\over{F_e(\lambda )}}\Bigr\}\eqno(5.48)$$
\noindent determines a hyperelliptic curve of genus $g\leq \infty$. There exists a homology basis of ${\cal E}$ which includes loops $\alpha_j$ around $[\lambda_{2j-1}, \lambda_{2j}]$. Let $\mu_j(x)$ lie on $\alpha_j$.\par
\vskip.05in
\noindent {\bf Definition} {\sl (Discretely dominated).}  Then $g_0$ is discretely dominated if there exist such data, with only the $\mu_j$ depending upon $x$ and 
$$g_0(x,-\zeta )={{1}\over{2\sqrt{\lambda_0-\zeta}}}\prod_{j=1}^g {{\mu_j(x)-\zeta }\over{\sqrt{(\lambda_{2j}-\zeta )(\lambda_{2j-1}-\zeta )}}}\qquad (\Im \zeta >0).\eqno(5.49)$$
\vskip.05in
\indent Gesztesy and Simon provide several examples of real potentials such that $g_0$ is discretely dominated; our constants are consistent with their  Example 3.2. The $\mu_j(x)$ are referred to as Dirichlet eigenvalues, or tied eigenvalues, while the differential equation (5.40) was derived by Dubrovin for real finite-gap periodic potentials. Given the general form (5.49), the further analysis reduces to the system of coupled differential equations for $d\mu_j(x)/dx.$  When $g_0$ is discretely dominated, its properties are most easily understood in terms of conformal mapping.\par
\indent In the following, we define the square root function by $s(\lambda )=\sqrt{\lambda} =\vert\lambda\vert^{1/2}\exp (i\arg (\lambda )/2)$, where $-\pi <\arg \lambda\leq\pi$ is the principal value of the argument, and recall that $h^*(\lambda ,x)=\overline{h(\bar\lambda ,\bar x)}.$  \par
\vskip.05in
\noindent {\bf Proposition 5.8.} {\sl (i) Suppose that $(\lfloor A^{2j-1}\rfloor )_{j=1}^\infty$ is a sequence of real functions of $x\in {\bf R}$. Then $\varphi^*  (\lambda ;x)=\varphi (\lambda ;x)$.\par  
(ii) Suppose that $g_0$ is discretely dominated with $\lambda_0=0$ and $\lambda_j$ real, and suppose that that $\Delta (\lambda ;x)$ is real for all $\lambda\in {\bf R}$. Then $\lambda \mapsto \varphi (\lambda ;x)$ gives a conformal map of the upper half plane onto a slit domain in the first quadrant with possible vertical slits at $n\pi$, for $n\in {\bf N}$, and $\lambda\mapsto \Delta (\lambda ;x)$ is entire.} 
\vskip.05in
\noindent {\bf Proof.}  (i) By the uniqueness of asymptotic expansions in Theorem 5.2(i), we have $g_0^*(x;\zeta )=g_0(x;\zeta )$, and since $s^*(\zeta )=s(\zeta )$, we deduce that $\varphi^*  (\lambda ;x)=\varphi (\lambda ;x)$ for $x\in {\bf R}$ and $\lambda \in {\bf C}\setminus {\bf R}$. At this stage, we do not claim that $\varphi$ is continuous across the real axis.\par
\indent (ii) Suppose that $\lambda_0<\lambda_1<\lambda_2<\dots <\lambda_{2g}$ and $\lambda_{2j-1}<\mu_j<\lambda_{2j}$ for $j=1, \dots, g$; then for $c_0\geq 0$, let
$$\varphi (\lambda)=\int_0^\lambda {{1}\over{2\sqrt{\zeta}}}\prod_{j=1}^g {{\mu_j-\zeta }\over{\sqrt{(\lambda_{2j}-\zeta )(\lambda_{2j-1}-\zeta )}}}d\zeta +ic_0.\eqno(5.50)$$
\noindent Then $\varphi $ is holomorphic for $\{\lambda\in {\bf C}:\Im \lambda >0\}$ with $\varphi (\lambda )/\sqrt\lambda \rightarrow 1$ as $\vert\lambda\vert\rightarrow\infty$; the image of $\lambda\in (-\infty ,0)$ is the subinterval $(ic_0, i\infty )$ of the imaginary axis; whereas  the image of $(0, \infty )$ consists of horizontal line segments running from left to right, interspersed by vertical lines running upwards then downwards, and the horizontal line segment together run towards $\infty$. For all sufficiently large $c_0>0$, the image of  $\{\lambda\in {\bf C}:\Im \lambda >0\}$  is a domain $\Omega$ in the first quadrant, with boundary consisting of horizontal and vertical line segments. This is a degenerate case of the Schwarz--Christoffel map, in which  the triple $\lambda_{2j-1}<\mu_j<\lambda_{2j}$ corresponds to a degenerate triangle on the edge of $\Omega$; two of the vertices of the degenerate triangle may coincide, giving a slit.\par
\indent Now suppose $c_0=0$, and observe that $\Delta (\lambda ;x)$ is real, if and only if either $\varphi (\lambda ;x)\in {\bf R}$ or there exists $n\in {\bf Z}$ such that $\varphi (\lambda ;x)-n\pi \in i{\bf R}$. We deduce that the only possible slits are at $\Re \varphi (\lambda ;\pi )=n\pi $ for some $n\in {\bf N}$, that all slits are vertical, and each starts and finishes at the same point on the real axis. \par   
\indent We apply the Schwarz reflection principle. By (i), $\Delta^* (\lambda ;x)=\Delta (\lambda ;x)$, so we need to check continuity across $\lambda\in {\bf R}$ Whereas $\sqrt\zeta $ is discontinuous across $(-\infty ,0)$, we only need deal with $\cos\sqrt\zeta $, which is continuous. More precisely, by (i), the function $\varphi (\lambda ;x)$ satisfies $\varphi (\lambda ;x)=\varphi^* (\lambda ;x)$, and $g(x;-\zeta )$ and $s(-\zeta )^2$ are continuous across $(-\infty ,0)$, hence $\varphi (\lambda ;x)^2$ is holomorphic across $(-\infty ,0)$, so $\Delta (\lambda ;x)$ is holomorphic across $(-\infty ;0)$; likewise,  $\Delta (\lambda ;x)$ is holomorphic across the spectral bands $[\lambda_{2j}, \lambda_{2j+1}]$. As $\lambda$ approaches a spectral gap $(\lambda_{2j-1}, \lambda_{2j})$ the image $\varphi (\lambda ;x)$ approaches a slit on the boundary of $\Omega$; now $\cos (k\pi \pm iy )=(-1)^k\cosh y$, so $\Delta (\lambda ;x)$ takes the same value, irrespective of which side $\lambda$ approaches from, hence $\Delta (\lambda ;x)$ is continuous across ${\bf R}$ and defines an entire function.\par    
 \rightline{$\square$}\par

\vskip.05in
\noindent {\bf 6. The differential ring of a periodic
linear system}\par

\vskip.05in
\indent In this section we obtain analogues of Theorem 4.3 for periodic
groups. First we formulate the notion of a periodic linear system, where we take section 4 as our guide. We show that the corresponding $\tau$ functions have properties analogous to those in Theorem 4.3. For periodic and meromorphic $u$, the differential equation
$-\psi''+u\psi =\lambda\psi$ is known as the complex Hill's equation. We show how periodic linear systems appear in the Floquet solutions, and obtain a counterpart of Proposition 5.5. Previous authors [22, 28 37, 73] explored the connection between Hill's equation and scattering solutions of Schr\"odinger's equation on the line. In the current paper, we show how Lyapunov's equation is the basis for some analogies. In section 7 we consider particular examples, which exhibit subtle effects.\par  
\par 
\indent For periodic linear
systems, the defining integral for $R_x$ in Proposition 2.1 does
not converge, and the contour integral for $R_0$ in Proposition 4.5 is
inapplicable; nevertheless, we can adapt a result of 
Bhatia, Dacis and McIntosh discussed in [8] and otherwise construct $R_x$
satisfying (1.3).\par
\vskip.05in
\noindent {\bf Lemma 6.1.} {\sl Let $B\in {\cal L}^1(H)$ and $C\in {\cal L}(H)$, and let $(e^{-tA})_{t\in {\bf R}}$ be a bounded and strongly continuous group of
operators on $H$.\par
\indent (i) The space ${\cal D}(A^\infty )$ is dense.\par
\indent (ii) Suppose that the spectrum of $A$ does not intersect the
spectrum of $-A$. Then there exists $E\in {\cal L}^1(H)$ such that $R_x=e^{xA}Ee^{xA}$ gives a solution to the Lyapunov equation 
$-{{d}\over{dx}}R_x=AR_x+R_xA$ such that $AR_0+R_0A=BC$ and $R_x$ is trace class for all
$x\in {\bf R}$. \par}
\indent {\sl (iii)  Suppose that the range of $E$ is contained in ${\cal D}(A^\infty )$ and $A^kE\in {\cal L}^1(H)$ for all $k$. Then
$$\tau_\zeta (x)=\det \bigl( I+(\zeta I+A)(\zeta I-A)^{-1} e^{xA}Ee^{xA}\bigr)\qquad (x\in {\bf R})\eqno(6.1)$$
\noindent has an asymptotic expansion in powers of $\zeta^{-j}$ as $\zeta \rightarrow\pm \infty$.}\par
\indent {\sl (iv) Suppose further that $(e^{xA})$ is periodic with period $2\pi$. Then the spectrum of $A$  is contained in $i{\bf Z}$ and the coefficients in the asymptotic expansion are periodic with period $\pi$.}\par  
\vskip.05in
\noindent {\bf Proof.} (i) By standard results [21], $A^2$ generates an analytic semigroup
$$e^{tA^2}={{1}\over{\sqrt{4\pi t}}}\int_{-\infty}^\infty e^{-s^2/4t} e^{sA}\, ds\qquad (t>0).\eqno(6.2)$$
\noindent The domains of the powers of $A$  satisfy ${\cal D}(A)\supseteq {\cal D}(A^2)$, and ${\cal D}(A^\infty )$ is dense.\par
\indent (ii) The main problem is to find $E$ such
that $EA+AE=BC$. By a theorem of Sz.-Nagy, the group $(e^{-tA})$ is similar to a group
of unitaries, so there exists an invertible operator $S$ and a unitary
group $(U_t)_{t\in {\bf R}}$ such that $e^{-tA}=SU_tS^{-1}$. Hence the
spectrum of $A$ lies on $i{\bf R}$ and is a closed subset. By
hypothesis, there exists $\delta >0$ such that the spectra of $A$ and
$-A$ are separated by $\delta$ and $\sigma (A)\cup \sigma (-A)$ does not
intersect $[-i\delta , i\delta ]$. By Plancherel's theorem, we can
construct $f\in L^1({\bf R}; {\bf C})$ such that $\hat f(\xi )=1/\xi$ for
all $\xi\in {\bf R}$ such that $\vert \xi\vert \geq \delta$. Then the
integral
$$E=\int_{-\infty}^\infty e^{-xA}BCe^{-xA} f(x)dx\eqno(6.3)$$
\noindent has a weakly continuous integrand in ${\cal L}^1(H)$, and is absolutely convergent with
$$\Vert E\Vert_{{\cal L}^1(H)}\leq \int_{-\infty}^\infty \Vert
B\Vert_{{\cal L}^1(H)}\Vert C\Vert_{{\cal L}(H)} M^2 \vert f(x)\vert\, dx\eqno(6.4)$$
\noindent hence $E$ is trace class. Using the spectral representation of $U_t$, one can show that
$AE+EA=BC$. Next we introduce $R_x=e^{-xA}Ee^{-xA}$ which gives a
one parameter family of trace class operators such that
$-{{dR_x}\over{dx}}=AR_x+R_xA.$ On verifies this identity on 
${\cal D}(A)$ and then observes that both sides are trace class.\par
\indent (iii) For all $\zeta\in {\bf R}$ we can invert $\zeta I-A$, and there is an asymptotic expansion 
$$ (\zeta I+A)(\zeta I-A)^{-1}Ee^{2xA}=\Bigl( I+{{2A}\over{\zeta}} +{{2A^2}\over{\zeta^2}}+\dots \Bigr)Ee^{2xA}\eqno(6.5)$$
\noindent valid as $\zeta\rightarrow\pm \infty$. The result follows.\par
\indent (iv) If $(e^{xA})$ is periodic with period $2\pi$, then $(e^{xA})_{x\in {\bf R}}$ is bounded, and the spectrum of $A$ is contained in $i{\bf Z}$.\par 
\rightline{$\square$}\par
\vskip.05in
\noindent {\bf Example 6.2.} Let $H=L^2({\bf R}/2\pi {\bf Z}; d\theta/(2\pi ))$ and $A:e^{in\theta}\mapsto  i(1+\vert n\vert )e^{in\theta}$, so $(e^{tA})$ is a $2\pi$ periodic strongly continuous unitary group. Also, $(e^{itA})_{t>0}$ gives a strongly continuous contraction semigroup $e^{itA}: \sum_{n\in {\bf Z}} a_n e^{in\theta}\mapsto e^{-t}\sum_{n\in {\bf Z}} a_n e^{-\vert n\vert t}e^{in\theta}$. For comparison, the Poisson semigroup is $P_r : \sum_{n\in {\bf Z}} a_n e^{in\theta}\mapsto \sum_{n\in {\bf Z}} a_nr^{\vert n\vert}e^{in\theta}$, so $e^{itA}=e^{-t}P_{e^{-t}}$ for $t>0$.\par
\vskip.05in

\noindent {\bf Definition} {\sl (Periodic linear system).} (i)  Let $(e^{-xA})_{x\in {\bf R}}$ be a 
strongly continuous group of operators on $H$ such that $e^{2\pi A}=I$ and $A$ is
invertible. 
Suppose further that $E$ is trace class
operators on $H$, and
that $B:H_0\rightarrow H$ and $C:H\rightarrow H_0$ are bounded linear operators, 
such that $AE+EA=BC$ and either $B$ is trace class, or $B$ and $C$ are Hilbert--Schmidt. Then $\Sigma_\infty =(-A,B,C;E)$ is a
periodic linear system with input and output space $H_0$ and state space $H$. (Unlike in Theorem 2.2, we generally take $H=H_0$ so the linear system has infinitely many inputs and outputs.)\par 
\indent (ii)  Moreover, if $(e^{-xA})_{x\in {\bf R}}$ is uniformly continuous, or equivalently $A\in {\cal L}(H)$, we say that $(-A,B,C;E)$ is a  uniform periodic linear system.\par
\indent (iii)  The $\tau$ function of $\Sigma_\infty$ is $\tau (x)=\det (I+e^{-xA}Ee^{-xA})$; then let 
$u(x)=-2{{d^2}\over{dx^2}}\log\tau_\infty (x)$ be the potential. \par
\indent (iv) Let $\Phi(x)=Ce^{-xA}B$ be the operator scattering function so that $\phi (x)={\hbox{trace}}\,\Phi
(x)$ is the (scalar) scattering function.\par
\indent (v) Let $R_x=e^{-xA}Ee^{-xA}$, then we introduce
$F_x=(I+e^{-xA}Ee^{-xA})^{-1}$.\par 
\indent (vi) Let
${\hbox{Spec}}(A)$ be the spectrum of $A$ as an operator in $H$, let ${\bf P}={\bf C}\cup\{ \infty \}$ be the Riemann sphere and
introduce
the periodic linear system 
$$\Sigma_\lambda=\bigl(-A,(\lambda I+A)(\lambda I-A)^{-1}B, C; 
(\lambda I+A)(\lambda I-A)^{-1}E\bigr)\qquad (\lambda \in {\bf
P}\setminus {\hbox{Spec }}(A))\eqno(6.6)$$ 
\noindent and its accompanying tau function $\tau_\lambda$.\par
\indent (vii)  We also introduce
the (non commutative) algebra ${\cal S}={\bf C}\{I, A, BC, F_x\}$, 
and then let ${\cal A}$ be the subring of ${\cal S}$ spanned by $A^{n_1}$
and by the ordered products $A^{n_1}FA^{n_2}\dots FA^{n_r}$ for $n_j\in {\bf N}$.
\vskip.05in

\noindent {\bf Definition} {\sl (Bracket and $*$ product).} (i) As in Lemma 4.1, we introduce on  ${\cal S}$ the product
$\ast$ and derivation $\partial$ by
$$P\ast Q=P(AF+FA-2FAF)Q, \qquad \partial
P=A(I-2F)P+{{dP}\over{dx}}+P(I-2F)A.\eqno(6.7)$$
 \indent (ii) We also introduce the bracket  
 $\lfloor\, .\, \rfloor:{\cal S}\rightarrow {\cal M} ({\cal L}^1(H)):$ 
$$\lfloor P\rfloor=Ce^{-xA}FPFe^{-xA}B.\eqno(6.8)$$ 
\indent (iii) Let $\lfloor {\cal A}\rfloor =\{ \lfloor P\rfloor : P\in {\cal A}\}$ and 
 $\Theta =\{ P\in {\cal A}: {\hbox{trace}}\lfloor P\rfloor=0\}$, which is a linear subspace, and not necessarily a ring. Then 
${\cal A}/\Theta =\{ {\hbox{trace}}\lfloor P\rfloor :P\in {\cal A}\}$, so that ${\cal A}/\Theta$ is analogous to the differential ring generated by the potential
$u$.\par
\vskip.05in
\noindent {\bf Theorem 6.3.} {\sl Let $(-A, B,C;E)$ be a uniform periodic linear system.\par
\indent (i) Then  $\tau_\lambda (x)$ is holomorphic except at finitely many singularities; so $x\mapsto \tau_\lambda (x)$ is entire, while
$\lambda\mapsto \tau_\lambda (x)$ is holomorphic on} 
${\bf P}\setminus{\hbox{Spec}}(A)$ ;\par
\indent {\sl (ii) there is a homomorphism of complex differential rings 
$({\cal A}, \ast, \partial )\rightarrow {\bf M}_{\bf C}({\cal L}(H_0))$ given by $X\mapsto \lfloor X\rfloor$;}\par
\indent {\sl (iii) the potential $u$ is meromorphic and $\pi$-periodic on  ${\bf C}$ and belongs to ${\cal A}/\Theta .$}\par 
\indent {\sl (iv) Also, let $T(x,y)=-Ce^{-xA}F_xe^{-yA}B$. Then
$${{\partial^2}\over{\partial x^2}}T(x,y)-
{{\partial^2}\over{\partial
y^2}}T(x,y)=-2\Bigl({{d}\over{dx}}T(x,x)\Bigr) T(x,y),\eqno(6.9)$$
\noindent and } $u(x)=-2{{d}\over{dx}}{\hbox{trace}}\, T(x,x)$. \par

\vskip.05in

\noindent {\bf Proof.} (i) First we show that $A$ is an algebraic operator. By periodicity, the group $(e^{-xA})_{x\in
{\bf R}}$ is bounded and hence by Sz.-Nagy's theorem, $e^{xA}$ is similar to a unitary group on 
$H$, so $A$ is similar to a skew symmetric operator. 
By uniform continuity, $A$ is bounded, and hence has spectrum contained
 in  $\{ -iN, \dots , iN\}$ for some integer $N$; see [21]. Consequently, there exists a monic polynomial $p$ such that $p(A)=0$. As in Proposition 4.5, $Ce^{-xA}B$ satisfies a linear differential equation with constant coefficients.\par
\indent Hence $A$ is an invertible algebraic operator, so 
as in (4.21),  $A^{-1}$ is a polynomial in $A$ and 
$(\lambda I+A)(\lambda I-A)^{-1}\in {\cal S}$
for all $\lambda$ in the resolvent set of $A$. Observe that 
$(\lambda I+A)(\lambda I-A)^{-1}$ is a
polynomial in $A$ with coefficients that are rational
functions of $\lambda$, and holomorphic except when $\lambda$ is in
the spectrum of $A$; in particular it is holomorphic on $\{ \lambda :
\vert \lambda\vert <1\} \cup\{ \lambda : \vert \lambda \vert >\Vert A\Vert\}$.\par
\indent  We also introduce 
polynomials $p_j$ 
for each point in the spectrum of $A$ such that $p_j(ik)=\delta_{jk}$ for $k=-N, \dots ,N$, 
and since $A$ is similar to a skew operator, we deduce that 
$$e^{-xA}=\sum_{j=-N; j\neq 0}^Np_j(A)e^{-ijx}. \eqno(6.10)$$

 Hence $\tau_\lambda$ is a holomorphic function of $\lambda$,
except at $\lambda\in {\hbox{Spec}}(A)$, which is
a finite set.\par
\indent (ii) First we check that $({\cal S}, \ast, \partial )$ is a complex differential ring for $(-A, B,C; E)$
and for $\Sigma_\lambda$. By (4.19), the operator $E$ belongs to ${\cal S}$ and hence by (6.10)
$R_x=e^{-xA}Ee^{-xA}$ also belongs to ${\cal S}$. Hence we have
$${{d}\over{dx}}R_x=-e^{-xA}AEe^{-xA}-e^{-xA}EAe^{-xA}=
-e^{-xA}BCe^{-xA}.\eqno(6.11)$$
By the Riesz
theory of compact operators, $F_x$ is a meromorphic operator-valued function of
$x$, and so $AF+FA-2FAF=Fe^{-xA}BCe^{-xA}F$, hence
$${{dF}\over{dx}}=AF+FA-2FAF;\eqno(6.12)$$
\noindent so that 
$$\partial (P\ast Q)=(\partial P)\ast Q+P\ast (\partial Q);\eqno(6.13)$$
\noindent thus $({\cal S}, \ast , \partial )$ is a differential ring.
Moreover $\lfloor\, .\,  \rfloor :({\cal S}, \ast , \partial )\rightarrow 
({\cal S}, \cdot , d/dx )$ is a homomorphism of differential rings, in the
sense that 
$${{d}\over{dx}}\lfloor P\rfloor =\lfloor \partial P\rfloor, \quad 
\lfloor P\rfloor \lfloor Q\rfloor=\lfloor P\ast Q\rfloor .\eqno(6.14)$$
\noindent Note that ${\cal A}$ is a subring of ${\cal S}$, and hence
$\lfloor {\cal A}\rfloor$ is also a differential ring.\par 
\indent  (iii) Since $e^{-xA}$ is an entire operator function, we deduce that $\tau_\infty$ is 
entire, and $\pi$ periodic 
since $\tau_\infty (x)=\det (I+e^{2xA}E)$ and $e^{2\pi A}=I$. When 
$\tau_\infty (x)\neq 0$, we have
$$\eqalignno{{{d}\over{dx}}\log\det (I+e^{-xA}Ee^{-xA})&
=-{\hbox{trace}}\bigl((I+e^{-xA}Ee^{-xA})^{-1}e^{-xA}(AE+EA)e^{-xA}
\bigr)\cr
&=-{\hbox{trace}}\bigl((I+e^{-xA}Ee^{-xA})^{-1}e^{-xA}BCe^{-xA}\bigr)\cr
&=-{\hbox{trace}}\bigl(Ce^{-xA}(I+e^{-xA}Ee^{-xA})^{-1}e^{-xA}B\bigr)\cr
&=-{\hbox{trace}}\bigl( Ce^{-xA}Fe^{-xA}B\bigr),&(6.15)\cr}$$
\noindent and hence
$$\eqalignno{u&=-2{{d^2}\over{dx^2}}\log\det
(I+e^{-xA}Ee^{-xA})\cr
&=-4{\hbox{trace}}\, Ce^{-xA}FAFe^{-xA}B\cr
&=-4{\hbox{trace}}\, 
\lfloor A\rfloor ;&(6.16)\cr}$$
\noindent so $u$ belongs to ${\cal A}_0=\{ {\hbox{trace}}\lfloor
P\rfloor :P\in {\cal A}\}$. Likewise, the derivatives $u^{(j)}$ belong 
to ${\bf A}_0$
since $\lfloor{\cal A}\rfloor$ is a differential
ring.\par
\indent (iv) By repeated differentiation, we have 
$${{\partial^2}\over{\partial x^2}}T(x,y)=-Ce^{-xA}\bigl( FA^2-4FAFA-4FA^2F+8FAFAF\bigr) e^{-yA}B,$$
$${{\partial^2}\over{\partial
y^2}}T(x,y)=-Ce^{-xA}\bigl( FA^2\bigr) e^{-yA}B\eqno(6.17)$$
\noindent while the Lyapunov equation gives
$$\eqalignno{-2\Bigl({{d}\over{dx}}T(x,x)\Bigr) T(x,y)&=4Ce^{-xA}FAFc^{-xA}BCe^{-xA}Fe^{-yA}B\cr
&= 4Ce^{-xA}FA(AF+FA-2FAF)e^{-yA}B, &(6.18)\cr}$$
\noindent hence the result. In [12, Lemma 7.2] we obtained a variant of this formula via an integral equation in the style of Gelfand--Levitan, under special commutativity conditions.  \par
\rightline{$\square$}

\noindent {\bf Definition} {\sl (Solutions of Hill's equation).}  (i) Let $u\in C^2({\bf R}; {\bf C})$ be $\pi$-periodic, and consider Hill's equation in the form
$${{d}\over{dx}}\left[\matrix{\psi\cr \psi'\cr}\right] =\left[\matrix{0&1\cr u(x)-\lambda&0\cr}\right] \left[\matrix{\psi\cr \psi'\cr}\right].\eqno(6.19)$$
\noindent Let $F_\lambda (x)$ be the fundamental solution matrix, and $\Delta (\lambda )={\hbox{trace}}(F_\lambda (\pi ))$ the discriminant.\par
\indent (ii) The multiplier curve is ${\cal E}=\{{\bf p}=
(\lambda , \rho ) : \rho^2-\Delta (\lambda )\rho +1=0\}$, as in [22]. \par 
\indent (iii) A Floquet solution consists of a 
nonzero function $f$ such that $-f''(x)+q(x)f(x)=\lambda f(x)$  and $f(x+\pi )=\rho f(x)$ for some $\lambda$. We call $\lambda$ the eigenvalue and $\rho$ the multiplier, and $(\lambda , \rho )$ lies on ${\cal E}$.\par
\indent (iv) In particular, when $\rho =1$, we say that $f$ is periodic, and when $\rho=-1$, we say that $f$ is anti periodic. Then $(\lambda , \pm 1)$ is a branch point on ${\cal E}$. \par
\indent (v) Suppose that $u$ is real-valued; then $L=-d^2/dx^2 +u$ is essentially self-adjoint. The Bloch spectrum of Schr\"odinger's operator consists of those $\lambda\in {\bf C}$ such that there exists a nontrivial and bounded solution 
$f$ of $-f''+uf=\lambda f$.  Real potentials are said
to belong to the same spectral equivalence class if their multiplier
curves are equal. See [47, 50]. \par
\indent (vi)  The Dirichlet eigenvalues $(\mu_j)_{j=1}^\infty $ are the $\mu$ such that
$$-y''(x;\mu )+u(x)y(x,\mu )=\mu y(x; \mu )$$
$$y(0;\mu )=0=y(\pi ;\mu )\eqno(6.20)$$
\noindent has a nontrivial solution. Replacing $u(x)$ by $u(x+t)$ we obtain the Dirichlet eigenvalues $(\mu_j (t))_{j=1}^\infty.$   
\vskip.05in
\indent The following result is related to results from Brett's thesis [13], and relates to the case in which the Gelfand--Levitan equation is  scalar-valued, as in Proposition 3.1 and [27]. Given a periodic linear system $(-A,B,C)$, we introduce
$$S_x=\int_0^x e^{-tA}BCe^{tA}\, dt, \qquad V_x=\int_0^x e^{tA}BCe^{tA}\, dt\eqno(6.21)$$
\noindent which satisfy the Lyapunov equation
$${{d}\over{dx}}\left[\matrix{V&0\cr 0&S\cr}\right] =\left[\matrix{A&0\cr 0&-A\cr}\right] \left[\matrix{V&0\cr 0&S\cr}\right]+\left[\matrix{V&0\cr 0&S\cr}\right] \left[\matrix{A&0\cr 0&A\cr}\right]+\left[\matrix{BC&0\cr 0&BC\cr}\right].\eqno(6.22)$$
\noindent Let  $W_x=V_x-S_x$  and
$$L(x,y)=-Ce^{xA}(I+W_x)^{-1}(e^{yA}-e^{-yA})B;\eqno(6.23)$$
\noindent then let $Z_x=V_x+S_x$ and
$$K(x,y) =-Ce^{xA}(I+Z_x)^{-1}(e^{yA}+e^{-yA})B.\eqno(6.24)$$
\vskip.05in
\noindent {\bf Proposition 6.4.} {\sl Let $(-A,B,C)$ be a periodic linear system with $H_0={\bf C}$ such that $\phi(x)=Ce^{xA}B$ is even and suppose there exist $E, E_-\in {\cal L}^1(H)$ such that $AE+EA=BC$ and $-AE_-+E_-A=BC$.\par
\indent (i)  Then the potential $w(x)=-2{{d^2}\over{dx^2}}\log \det (I+W_x)$ is periodic, and 
$$\varphi (x)={{\sin kx}\over{k}}+\int_0^x L(x,y){{\sin ky}\over{k}}dy\eqno(6.25)$$
\noindent satisfies
$$\eqalignno{ -\varphi''(x)+w(x)\varphi (x)&=k^2\varphi (x)\cr
\varphi (0)=0,\quad &\quad \varphi'(0)=1.&(6.26)}$$ 
\noindent  (ii) Also the potential  $u(x)=-2{{d^2}\over{dx^2}}\log \det (I+Z_x)$ is periodic and 
$$\psi (x)=\cos kx+\int_0^x K(x,y)\cos ky\, dy\eqno(6.27)$$
\noindent satisfies}
$$\eqalignno{ -\psi''(x)+u(x)\psi (x)&=k^2\psi (x)\cr
\psi (0)=1,\quad &\quad \psi'(0)=-2CB.&(6.28)}$$ 
\vskip.05in
\noindent {\bf Proof.} (i) We have
$$\det (I+W_x)=\det (I+e^{xA}Ee^{xA}-e^{-xA}E_-e^{xA}-E+E_-)\eqno(6.29)$$
\noindent which is periodic. Also,  $L(x,0)=0$ and 
$$L(x,x)=-{{d}\over{dx}}\log\det (I+W_x).\eqno(6.30)$$
\noindent One can verify that
$$\phi (x+y)-\phi (x-y)+L(x,y)+\int_0^x L(x,t)\bigl( \phi (t+y)-\phi (t-y)\bigr)\, dt=0.\eqno(6.31)$$
\noindent and with  $w(x)=2{{d}\over{dx}}L(x,x)$, we deduce that
$${{\partial^2L}\over{\partial x^2}}-{{\partial^2L}\over{\partial y^2}}=w(x)L(x,y)\eqno(6.32)$$
\noindent and by manipulating the Gelfand--Levitan equation, one deduces the differential equation for $\varphi$.\par
\noindent (iii) We have
$$\det (I+Z_x)=\det (I+e^{xA}Ee^{xA}+e^{-xA}E_-e^{xA}-E-E_-),\eqno(6.33)$$ 
\noindent which is periodic, with ${{\partial}\over{\partial y}}K(x,0)=0$ and 
$$K(x,x)=-{{d}\over{dx}}\log\det (I+Z_x).\eqno(6.34)$$
Then we verify that
$$\phi (x+y)+\phi (x-y)+K(x,y)+\int_0^x K(x,t)\bigl( \phi (t+y)+\phi (t-y)\bigr)\, dt=0.\eqno(6.35)$$
\noindent and  with $u(x)=2{{d}\over{dx}}K(x,x)$, we have
$${{\partial^2K}\over{\partial x^2}}-{{\partial^2K}\over{\partial y^2}}=u(x)K(x,y);\eqno(6.36)$$
\noindent then  by manipulating the Gelfand--Levitan equation, one deduces the differential equation for $\psi$.\par
\rightline{$\square$}\par
\vskip.05in
\indent Given Proposition 6.4, it is tempting to seek a version of Proposition 6.1 for periodic linear systems, and try to express the general solution of Hill's equation in terms of quotients of tau functions of periodic linear systems. The following result indicates the restriction that must be imposed upon $u$ for such a representation to be valid.\par 
\vskip.05in
\noindent {\bf Proposition 6.5.} {\sl Suppose that $\Sigma =(-A,B,C;E)$ is a uniformly periodic linear system with $H_0=H$ such that for all but finitely many $\lambda\in {\bf C}$, Hill's equation  has a pair of linearly 
independent  solutions of the form
$$\psi_\lambda (x)=e^{\nu x}\prod_{j=1}^n {{\tau_{\zeta_j}(x-a_j)}\over{\tau_{\eta_j}(x-b_j)}}\eqno(6.37)$$
\noindent for some $\nu, a_j,b_j,\zeta_j,\eta_k$ depending upon $\lambda$ .\par
\indent (i) Then $u$ is a Picard potential, in the sense that Hill's equation has a meromorphic general solution for all but finitely many $\lambda\in {\bf C}$ ;\par
\indent (ii)  $u$ is finite gap as a potential for Hill's equation, and there exists a differential operator $P_{2g+1}$ of order $2g+1$ such that 
$$P_{2g+1}^2=\prod_{j =0}^{2g} (L-\lambda_j)\eqno(6.38)$$
\noindent for some  $\lambda_j\in {\bf C}$ such that $\Delta^2(\lambda_j)=4$, and the diagonal Greens function $g_0$ is discretely dominated and satisfies (5.49).\par 
\indent (iii) Suppose further that $u$ is real-valued. Then the Bloch spectrum of $L$ is associated with a hyperelliptic algebraic complex curve}
$${\cal E}_0=\Bigl\{ (\mu  , \lambda )\in {\bf C}^2: \mu^2=\prod_{j=0}^{2g} (\lambda
-\lambda^o_j)\Bigr\} \cup \bigl\{ (\infty , \infty )\bigr\}\eqno(6.39)$$
\noindent {\sl where $\lambda_j^o$ a simple real zero of $\Delta (\lambda )^2-4=0$, and such that the $\lambda\in \sigma_B$  give real
points $(\mu,  \lambda  )$ on ${\cal E}_0$. }\par
\vskip.05in
\noindent {\bf Proof.} (i) In this case $x\mapsto \det (I+e^{-xA}Ee^{-xA})$ is entire and hence $\psi_{\lambda}$ is meromorphic; hence there exists a meromorphic fundamental system of solutions on ${\bf C}$. 
Since $\tau$ is periodic, we have $\phi_\lambda (x+\pi )=e^{\nu \pi}\phi_\lambda (x)$, so that $\psi_\lambda$ is a Floquet solution with Floquet exponent $\pi\nu$. By choosing a pair of solutions with Floquet exponents $\pm \pi\nu (\lambda )$, we obtain a fundamental solution matrix such that $\Delta (\lambda ) =2\cosh \pi\nu (\lambda )$.  \par
\indent (ii) This is a particular case of the Burchnall--Chaundry theorem, and the required $P_{2g+1}$ is given in the proof of Proposition 5.5. See [29, Theorem 4.1] for details.
Also, for periodic $u\in C^2({\bf R}, {\bf C})$, Gesztesy and Weikard [29] obtain the formula (4.31) where the $\lambda_j$ are the periodic eigenvalues which correspond to periodic solutions of (6.19). The $\mu_j(x)$ are referred to as Dirichlet eigenvalues, or tied eigenvalues.\par
\indent (iii)   The set of all $\lambda$ that give $\rho=\pm 1$ is called the periodic spectrum. Such a  $u$ has periodic spectrum $\lambda_0<\lambda_1\leq\lambda_2<\lambda_3\leq\dots ,$ where for $\lambda_{2j}$
there exists a periodic eigenfunction, while for $\lambda_{2j+1}$ there exists an anti periodic eigenfunction. We have $\sigma_B=\cup_{k=0}^\infty [\lambda_{2k}, \lambda_{2k+1}]$, where many of these intervals are abutting.
 For real $x$, we deduce that just as $\psi_\lambda (x)$ in (6.19) has parameters $\nu ,a_j,\eta_j, b_j$ and $\eta_j$ for $j=1, \dots ,n$ and gives a solution of Hill's equation with $\lambda$, the complex conjugate 
$\overline{\psi_\lambda (x)}$ has parameters $\bar \nu, \bar a_j, \bar\zeta_j, \bar b_j$ and 
$\bar \eta_j$ for $j=1, \dots ,n$ and hence gives a solution of Hill's equation with $\bar\lambda$. 
We select $\psi_\lambda$ with Floquet exponent $\pi\nu$, and $\phi_{\bar\lambda}$ with Floquet exponent 
$\pi\bar \nu$, and observe that $\overline {\psi_{\bar\lambda }}(x)=\psi_\lambda (x)$ as they both solve Hill's equation with $\lambda$ and have the same Floquet exponent $\pi\nu$. \par
\indent Note that $\tau_\infty  (x)$ is real and non zero for real $x$, since  $x\mapsto u_\infty (x):{\bf R}\rightarrow 
{\bf R}$ is twice continuously differentiable, 
and such that Hill's equation has two linearly independent Floquet solutions $\psi_{\lambda }(x)$ for all but finitely many
$\zeta$. Then Gesztesy and Weikard proved [31, Theorem 4.1] that $u_\infty (x)$ is
finite gap, and hence that the discriminant equation $\Delta (\lambda )^2-4=0$ has only finitely many zeros that are simple, namely $\lambda_0^o< \dots < \lambda^o_{2g}$. Hence the Bloch spectrum 
is $\cup_{j=0}^{g-1}[\lambda_{2j}^0, \lambda_{2j+1}^o]\cup [\lambda_{2g}^o, \infty )$ and each point in the interior of this union of real intervals corresponds to a pair of real points on the curve ${\cal E}$.\par
\rightline{$\square$}\par
\vskip.05in
\noindent {\bf Example 6.6.} Let $u\in C^2({\bf R}; {\bf R})$ be a finite gap Hill's potential.  Here the fundamental solution matrix is real for all $\lambda\in {\bf R}$, hence $\Delta (\lambda )\in {\bf R}$ for all $\lambda\in {\bf R}$. Hochstadt observed that the simple periodic spectrum $(\lambda_j^o)$ determines the double periodic spectrum and the nontrivial roots $\lambda_j'$ of $\Delta'(\lambda )=0,$ so we can write
$${{\Delta'(\lambda )}\over{\sqrt{4-\Delta (\lambda )^2}}}=c_3{{\prod_{j=1}^g (\lambda -\lambda_j')}\over{\sqrt{ \prod_{j=0}^{2g}(\lambda -\lambda_j^o)}}}.\eqno(6.40)$$
\noindent McKean and van Moerbecke [50] proved that the set of potentials with simple periodic spectrum $\lambda_0^0<\dots <\lambda^o_{2g}$ gives a $2^g$ sheeted cover of the cell $[\lambda^0_1, \lambda^o_2]\times \dots \times [\lambda^o_{2g-1}, \lambda^o_{2g}]$ in which the Dirichlet spectrum lives. We select one such potential $q$ such that 
$\mu_j=\lambda_j'$.  Then 
$$\varphi (\lambda )=\cos^{-1} {{\Delta (\lambda )}\over{2}}=-\int_0^\lambda {{\Delta'(\zeta )d\zeta}\over{\sqrt{4-\Delta (\zeta )^2}}}\eqno(6.41)$$
\noindent gives a conformal mapping from the upper half plane $\{ \lambda : \Im \lambda >0\}$ to a slit domain, as in Proposition 5.8. \par
\vskip.05in
\noindent {\bf 7. Linear systems on the complex torus, and the hyperelliptic prime function}\par
\vskip.05in
\noindent In this section, we show that there are significant examples of periodic linear systems which satisfy the hypotheses of the previous section. In order to give explicit formulas, we do not seek the greatest generality in the presentation. We start with genus one, for which our results are most complete, and then progress to hyperelliptic cases.  \par
\vskip.05in

\noindent {\bf Definition} {\sl (Elliptic functions).}  (i) For $\omega_1,\omega_2\in {\bf C}\setminus \{0\}$ with $\Im (\omega_2/\omega_1)>0$ uppose let  $\Lambda ={\bf Z}2\omega_1+{\bf Z}2w_2$ be a lattice, and let
${\cal T}={\bf C}/\Lambda$ be the corresponding torus. A
meromorphic function on ${\bf C}$ is elliptic (of the first kind) if it
is doubly periodic with respect to $\Lambda$.  A meromorphic function is elliptic of
the second kind if there exist multipliers $\rho_j\in {\bf C}$ such that
$f(z+2\omega_j )=\rho_jf(z)$; so that $f$ is quasi-periodic with
respect to $\Lambda$.  A meromorphic function is elliptic of the third kind if there exist $a_j, b_j\in {\bf C}$ for $j=1,2$ such that $f(z+2\omega_j)=e^{a_jz+b_j}f(z)$.  See [43, 49].\par
\indent (ii) Let $\omega\in {\bf C}$ have $\Im \omega >0$; then Jacobi's elliptic theta function is
$$\vartheta _1(x\mid\omega
)=i\sum_{n=-\infty}^\infty (-1)^n e^{(2n-1)\pi ix+(n+1/2)^2\pi
i\omega}\qquad (x\in {\bf C}),\eqno(7.1)$$ 
\noindent which is elliptic of the third kind on with respect to ${\bf Z}+\omega {\bf Z}$.\par
\vskip.05in
\noindent {\bf Lemma 7.1.} {\sl (i) Let $\delta$ be a positive divisor on  ${\bf C}/\Lambda$. Then there exists a uniformly periodic linear system with tau function $\tau$, where $\tau$ is elliptic of the third kind with zero divisor $\delta$.\par
\indent (ii) Let $\eta$ be a divisor on ${\bf C}/\Lambda$ of degree zero. Then there exists a pair of uniformly periodic linear systems with tau functions $\tau_0$ and $\tau_1$ such that $\tau_1/\tau_0$ is elliptic of the second kind with $\eta$ the divisor of its poles and zeros.\par
\indent (iii) Let $\gamma=\sum_{j=1}^n \delta_{a_j}-\sum_{j=1}^n\delta_{b_j}$ be a divisor of degree zero, where $\sum_{j=1}^n (a_j-b_j)\in \Lambda$. Then there exists a pair of uniformly periodic linear systems with tau functions $\tau_1$ and $\tau_0$ such that $\tau_1/\tau_0$ is elliptic of the first kind with $\gamma$ the divisor of the poles and zeros.}\par
\vskip.05in

\noindent {\bf Proof} (i) We introduce ${\hbox{q}}=e^{\pi i\omega}$ and introduce Jacobi's elliptic function of the third kind by the product
$$\vartheta_3(x)=\prod_{n=1}^\infty (1-{\hbox{q}}^{2n})\prod_{n=1}^\infty (1+2{\hbox{q}}^{2n-1}\cos 2\pi x +{\hbox{q}}^{4n-2}).\eqno(7.2)$$
\noindent Observing that 
$$\det\Biggl( \left[\matrix{1&0\cr 0&1\cr}\right]  +{\hbox{q}}^{2n-1}\left[\matrix{\cos 2\pi x&\sin 2\pi x\cr -\sin 2\pi x& \cos 2\pi x\cr}\right]\Biggr)=1+2{\hbox{q}}^{2n-1} \cos 2\pi x +{\hbox{q}}^{4n-2},\eqno(7.3)$$
\noindent one introduces
$$A_n=C_n= \left[\matrix{0&\pi \cr -\pi &0\cr}\right],\quad B_n=2E_n=\left[\matrix{{\hbox{q}}^{2n-1}&0\cr 0 &{\hbox{q}}^{2n-1}\cr}\right]\eqno(7.4)$$
\noindent so that $A_nE_n+E_nA_n=B_nC_n$. By forming the direct sum $\bigoplus_{n=1}^\infty (A_n, B_n, C_n; E_n)$  
can easily construct a linear system $(-A, B, C; E)$ with $H=H_0=\oplus_{n-1}^\infty {\bf C}^{2\times 1}$ that generates a uniform periodic semigroup and with tau function $\vartheta_3(x)/c$ for $c=\prod_{n=1}^\infty (1-{\hbox{q}}^{2n})$.\par
\indent We then introduce $\vartheta_1(x)= -ie^{\pi ix}{\hbox{q}}^{1/4}\vartheta_3 (x+\omega /2+1/2)$.  By simple manipulations of (7.22), one shows that $\vartheta_1$ is entire and elliptic of the third kind, and from the product formula it is evident that $\vartheta_1$ has a simple zero at $x=0$ and no others in the fundamental cell of ${\bf C}/\Lambda$.
In section 7 of [11] we likewise obtained a uniform periodic linear system with tau 
function $\vartheta_1(x)$.\par
\indent Given a positive divisor $\delta=\sum_{j=1}^n \delta_{a_j}$, we can form block diagonal sums, 
and obtain a periodic linear system with tau function
 $\tau_0(x)=\prod_{j=1}^n\vartheta_1(x-a_j)$, so $\tau_0$ is elliptic of the third kind with simple zeros at $a_1, \dots, a_j$ and the points congruent to these with respect to $\Lambda.$\par
\indent (ii) Given $\eta=\sum_j\delta_{a_j}-\sum_k\delta_{b_k}$, we introduce periodic linear systems with tau functions $\tau_1$ and $\tau_0$ as in (i) so that
$${{\tau_0(x)}\over{\tau_1(x)}}= \prod_{j=1}^n{{\vartheta_1(x-a_j)}\over{\vartheta_1(x-b_j)}}\eqno(7.5)$$
\noindent is elliptic of the second kind with  zeros $a_1, \dots, a_n$ and poles  $b_1, \dots, b_n$ listed according to multiplicity modulo $\Lambda$.\par
\indent (iii) In the particular  case of (ii) in which $\sum_{j=1}^n (a_j-b_j)\in \Lambda$, then $\tau_0/\tau_1$ is elliptic of the first kind, as in Abel's theorem.\par
 \rightline{$\square$}\par
\vskip.05in

\noindent {\bf Proposition 7.2.} {\sl Suppose that $u$ is elliptic of the first kind, and that Hill's
equation $-\psi''(x)+u(x)\psi (x)=\lambda \psi (x)$ has a meromorphic fundamental system of solutions for some $\lambda\in {\bf C}$.  Then Hill's equation has a nontrivial elliptic solution of the second kind,
which may be expressed as a quotient of tau functions that arise from uniformly periodic linear
systems and systems with finite-dimensional state space.}\par
\vskip.05in
\noindent {\bf Proof.} The first part is conventionally attributed to Picard. 
 Let $V_\lambda$ be the vector space of meromorphic solutions of Hill's equation, and suppose that
$V_\lambda$ has dimension two. Observe
that the monodromy operators $T_j:\psi (z)\mapsto
\psi (z+2\omega_j)$ are commuting
operators such that $T_j(V_\lambda )\subseteq V_\lambda $ for $j=1,2$  
since $u$ is
elliptic, so we can take $\Lambda$ to be the group generated by
$T_1$ and $T_2$. Then $T_1$ and $T_2$ have a common eigenvector, which gives an elliptic solution of 
the second kind.
 (Furthermore, if $T_1$ or $T_2$ has distinct
eigenvalues as an operator on $V_\lambda$, then there exists a fundamental system of
elliptic functions of the second kind; in particular, this happens for $T_1$ when
$\Delta (\lambda )^2-4\neq 0$.)\par
\indent Let $\psi$ be a solution that is elliptic of the second kind
 hence has the
form 
$$\psi (x)=e^{bx+c}{{\prod_{j=1}^n
\vartheta_1(x-a_j)}\over{\prod_{j=1}^n\vartheta_1(x-b_j)}}.\eqno(7.6)$$
\noindent The initial exponential factor $e^{bx+c}$ is the tau function of a  linear system with state space ${\bf C}$, while the quotient of $\vartheta_1$ functions can be realized from a pair of uniformly periodic 
linear systems by Lemma 7.1.
 \par
\rightline{$\square$}\par

\vskip.05in
\noindent {\bf Example 7.3.} {\sl (Lam\'e's equation).} (i) By [49, page 132] there exists a constant $e_2$ such that 
$$\wp (x)=e_2+\Bigl( {{\vartheta_1'(0)\vartheta_3(x)}\over {\vartheta_1(x)\vartheta_3(0)}}\Bigr)^2,\eqno(7.7)$$
\noindent where Weierstrass's function $\wp$ is meromorphic and doubly periodic with
respect to ${\bf Z}+\omega {\bf Z}$.\par  
\indent  The fundamental example of a finite gap elliptic differential
 equation is Lam\'e's equation. Let $(X,Z)=(\wp (x), \wp'(x))$ and $(Y,W)=(\wp (y), \wp'(y))$ be 
points on the elliptic curve ${\cal T}=\{ (X,Z): Z^2=4X^3-g_2X-g_3\}\cup \{(\infty
,\infty )\}$, where $g_2^3-27g_3^2\neq 0$, with Klein's
invariant $J=g_2^3/(g_2^3-27g_3^2)$, and ${\bf K}={\bf C}(X)[Z]$ is the elliptic function field. 
 Then Lam\'e's equation is
$$-{{d^2}\over{dx^2}}\psi (x)+\ell (\ell +1)\wp (x)\psi (x) =\lambda \psi (x)\eqno(7.8)$$ 
\noindent and we write $\Psi (X)=\psi (x)$ to convert between coordinates on the 
curve and the torus. See [44] for a detailed discussion of various forms of the solutions.\par
\indent  To solve the case $\ell =1$, we introduce 
$$\psi_2(x,\alpha )=-2q^{1/4}e^{(\zeta (\alpha )-2\alpha\eta_1/\pi) x}
{{\vartheta_1(x-\alpha )}\over{\vartheta_1(\alpha )\vartheta_1(x )}}\prod_{n=1}^\infty (1-q^{2n})^{3} ,\eqno(7.9)$$ 
\noindent  which satisfies Lam\'e's equation with $\lambda =-\wp (\alpha )$ and is such that $\alpha\mapsto \psi_2(x,\alpha )$
 is doubly periodic and meromorphic, and $x\mapsto \psi_2(x,\alpha )$ is elliptic of the
 second kind; moreover $\psi_2(x,\alpha )\psi_2(-x,\alpha )
=\wp (\alpha )-\wp (x)$. By Lemma 7.1, $\psi_2(x,\alpha )$ can be 
expressed as a quotient of tau functions from periodic linear systems.\par
\indent (ii) The Lam\'e example is fundamental, since several elementary examples can be derived from it. In each of the following, $\gamma$ and the potential $u$ are meromorphic functions on 
a Riemann surface ${\cal E}$ and $\psi$ satisfies 
 the addition rule
$$\psi (x+y)={{\psi'(x)\psi (y)-\psi (x)\psi'(y)}\over{\gamma (x)-\gamma (y)}}.
\eqno(7.10)$$
\noindent and each entry can be obtained from periodic linear systems by taking a limit
of the real or imaginary period to infinity.\par
\vskip.05in
$$\matrix{ {{\cal E}}& u(x)&\psi (x)&\gamma (x)&\tau (x)\cr
{\bf P}^1& g(g+1)/x^2&(g+1)/x&-(g+1)/x^2&x^{g(g+1)/2}\cr
{\bf C}/\pi {\bf Z}&g(g+1){\hbox{cosec}}^2x& (g+1)\cot x&-(g+1){\hbox{cosec}}^2x & 
(\sin x)^{g(g+1)/2}\cr
{\bf C}/\pi i{\bf Z}&g(g+1){\hbox{cosech}}^2x&(g+1)\coth x&-(g+1){\hbox{cosech}}^2x&
(\sinh x)^{g(g+1)/2}\cr
{\bf C}/\Lambda & 2\wp (x\mid\Lambda )&\psi_2(x,\alpha )& -\wp (x\mid\Lambda ) 
&\vartheta_1(x\mid \Lambda )\cr}\eqno(7.11)$$
\vskip.05in
\indent Lemma 7.1 extends to symmetric compact Riemann surfaces of genus $g\geq 2$ via the Schottky--Klein function $\varpi$, which depends upon the Schottky's model for ${\cal E}$. In Schottky's model, ${\cal E}$ arises as the quotient of ${\bf C}_\infty={\bf C}\cup\{ \infty \}$ under the action of a discrete subgroup of $PSL(2, {\bf C})$. 
 A Riemann surface ${\cal E}$ is said to be symmetric if there exists an anti-conformal involution $\varphi :{\cal E}\rightarrow {\cal E}$. Whereas the defining formula for $\varpi$ does not seem well adapted for computation,  Crowdy {\sl et al} [5,17, 18,70] have recently devised algorithms which facilitate computing $\varpi$ in geometrical contexts which arise in applied mathematics.\par
\vskip.05in

\noindent {\bf Example 7.4} {\sl (Schottky's model).} To each  $T\in PSL(2, {\bf C})$, we associate the M\"obius transformation  
$$T=\left[\matrix{a&b\cr c&d\cr}\right]: z\mapsto {{az+b}\over{cz+d}}.\eqno(7.12)$$
\noindent In particular, each $T\in PSL(2, {\bf R})$ gives a M\"obius transformation which preserves ${\bf R}\cup\{ \infty \}$, and is said to be Fuchsian. Let $\hat\gamma_0$ be the unit circle $C(0,1)$, and let $\hat\gamma_j$ for $j=1, \dots , g$ be Jordan curves inside $\hat\gamma_0$ which are mutually exterior to one another.
 Let $\varphi$ be the anti conformal involution of ${\bf C}_\infty$ consisting of reflection in $\hat\gamma_0$, as in $\varphi (z)=1/\bar z$ then we obtain $\hat\gamma_{-j}$ by reflecting $\hat\gamma_j$ in $\hat\gamma_0$ via $\varphi$ so  $\hat\gamma_{-j}$ are Jordan curves in the exterior of $\hat\gamma_0$; then $\hat\gamma_{\pm 1}, \dots , \hat\gamma_{\pm g}$ are mutually exterior. The fundamental region is the set ${\cal F}$ which is exterior to all the $2g$ Jordan curves $\hat\gamma_{\pm j}$ for $j=1, \dots , g$, so ${\cal F}$ has $2g$ holes, and comes equipped with the anti conformal involution $\varphi$. Let $\hat T_j\in PSL(2, {\bf C})$ map $\gamma_j$ to $\hat\gamma_{-j}$ and the exterior of $\hat\gamma_j$ with respect to ${\bf C}_\infty$ onto the interior of $\hat\gamma_{-j}$. Then we introduce the group $\Gamma$ which is generated by $\hat T_j$ for $j=1, \dots ,g$. Let $D_0=D(0,1)$ be the interior of $\hat\gamma_0$ and $D_j$ the interior of $\hat\gamma_j$ for $j-0, \dots, g;$ then let $D_\Gamma=D_0\setminus \cup_{j=1}^gD_j$, so that $D_\Gamma$ is a set with $g$ holes. Then we can form the Schottky double of $D_\Gamma$ by attaching another copy of $D_\Gamma$ to itself along the curves $\gamma_j$ for $j=0, \dots ,g$ to form a symmetric Riemann surface ${\cal E}$ with genus $g$. \par
\vskip.05in   
\indent Now we separate off the identity transformation, and pair up all other transformations with their inverses, so there is a partition
$\Gamma =\{ I\} \cup\bigcup_{k=1}^\infty \{ T_k, T_k^{-1}\}$, where  
$T_k(z)=(a_kz+b_k)/(c_kz+d_k)$.  By classical results summarized in [3], there exists a singular set $E_0$ such that $\Gamma$ acts properly discontinuously on ${\bf C}_\infty\setminus E_0$, the set $E_0$ is perfect and nowhere dense, and $E_0$ has Lebesgue area measure zero. Given any relatively compact domain $D'$ contained in ${\bf C}_\infty \setminus E_0$, there exist only finitely many points $-d_k/c_k\in D'$; hence there exists a nonempty domain $D''\subset D'$ such that $D''$ does not contain $-d_k/c_k$ or $\infty$. \par
\vskip.05in
\noindent {\bf Definition} {\sl (Prime function).}  Suppose momentarily that the following product converges for some $z, \zeta\in D''$; then the product defines the Schottky--Klein prime function of ${\cal E}$ by  
$$\varpi (z,\zeta )=(z-\zeta ) \prod_{k=1}^\infty {{z-T_k(\zeta )}\over{ z-T_k(z)}} {{\zeta-T_k(z)}\over{ \zeta -T_k(\zeta)}}.\eqno(7.13)$$
\vskip.05in
\noindent Each factor in this product is a cross-ratio, hence does not change when $T_k$ is replaced by $T_k^{-1}$.  We write $T_{-k}=T_k^{-1}$ for $k=1, 2, \dots $. See [7, Chapter 7] for other functional relations.\par
\vskip.05in
\noindent {\bf Proposition 7.5.} {\sl Suppose that $\Gamma$ is a non elementary Fuchsian group such that the limit set is a proper subset of ${\bf R}$, and let $\zeta_0\in {\cal F}$ have $\Im \zeta_0>0$. 
Then there exists a sequence of uniformly periodic linear systems $\Sigma_k$ with tau functions $\tau_k$, and a constant $c_1$  such that} 
$$ {{\varpi (e^{i\theta}, \zeta_0)}\over{\varpi (e^{i\phi}, \zeta_0)}}=c_1\prod_{k=-\infty }^\infty {{\tau_k(\theta )}\over{\tau_k(\phi )}}.\eqno(7.14)$$ 
\vskip.05in

\noindent {\bf Proof.} Burnside [15] showed that for $\Gamma$ such a  Fuchsian group,  the series $\sum_k\vert c_kz+d_k\vert^{-2}$ converges for all $\Im z>0$; he referred to such groups as groups of the first class, a terminology which is no longer current.  The  set of limit points of  $(T_k(\zeta_0))_{k=-\infty}^\infty$ is contained in the limit set of $\Gamma$, hence is a proper subset of ${\bf R}$. Akaza considered the property
$$ \sum_{k=1}^\infty {{1}\over{\vert c_k\vert^2}}<\infty\eqno(7.15)$$
\noindent and showed that this is equivalent to the absolute and uniform convergence of various Poincare series of dimension $(-2)$, and this property holds in the present context. Hence the series
$$\sum_k \Bigl( {{1}\over{T_k(s )-z }}-{{1}\over{T_k(s)-\zeta}}\Bigr) T_k'(s )\eqno(7.16)$$
\noindent converges for all $s\in D''$, and $\zeta, z\in \hat\gamma_0$, as in Theorem A of [3]. Upon integrating (7.16), and noting Appendix A of [5] we obtain
$$\prod_{T\in \Gamma} {{T(\zeta_0 )-z }\over{T(\zeta_0)-\zeta}} =C {{\varpi (z, \zeta_0 )}\over{\varpi (\zeta, \zeta_0 )}},\eqno(7.17)$$
\noindent for some constant $C$; we write this as a convergent product 
$$   C {{\varpi (z, \zeta_0 )}\over{\varpi (\zeta, \zeta_0 )}}=\prod_{k=1}^\infty {{ (a_k\zeta_0 +b_k)-(c_k\zeta_0+d_k)z}\over {
 (a_k\zeta_0 +b_k)-(c_k\zeta_0+d_k)\zeta }}.\eqno(7.18)$$
\indent Now we introduce a periodic linear system for each $k$, as specified by the matrices 
$$A_k=\left[\matrix {0&1\cr -1&0\cr}\right] , \qquad E_k =\left[\matrix {\alpha_k&\beta_k\cr \gamma_k&\delta_k\cr}\right],\eqno(7.19)$$
\noindent with coefficients to be determined. We have
$$\tau_k(\theta) =\det \bigl(I+e^{\theta A} E\bigr)=1+(\alpha_k\beta_k-\gamma_k\delta_k)+(\alpha_k+\delta_k)\cos\theta +(\gamma_k-\beta_k)\sin\theta,\eqno(7.20)$$
\noindent so we choose complex $\alpha_k$ and $\beta_k$ to solve the quadratic equation
$${{a_k\zeta_0+b_k}\over{c_k\zeta_0+d_k}}={{1-\alpha_k^2-\beta_k^2}\over{c_k\zeta_0+d_k}}-\alpha_k+i\beta_k,\eqno(7.21)$$
\noindent (for instance, one can choose $\beta_k=i\alpha_k)$, then let 
$$\delta_k=-\alpha_k-(c_k\zeta_0+d_k), \qquad \gamma_k=\beta_k-i(c_k\zeta_0+d_k).\eqno(7.22)$$
\noindent Let $z=e^{i\theta}$ and $\zeta =e^{i\phi}$ so $z$ and $\zeta$ lie on the circle $\hat\gamma_0$, which lies inside ${\cal F}$, and the product (7.18) converges. The identity (7.14) follows from our choice of $A_k$ and $E_k$.\par 
\rightline{$\square$}\par

\vskip.05in
\noindent {\bf Remarks 7.6} (i) Lemma 7.1 is a particular case of Proposition 7.5. When $g=1$ and  $\hat\gamma_1=C(0, {\hbox{q}} )$ for some $0<{\hbox{q}}<1$, the concentric circles $\hat\gamma_1$ and $\hat\gamma_{-1}$ bound an annulus, from which one can construct a model of the torus ${\bf C}/({\bf Z}2\omega_1+{\bf Z}2\omega_2)$ where $\omega_1=\log 1/{\hbox{q}}$ and $\omega_2=\pi i$, by identifying points on the inner and outer circles as in [48, p 48]. Here the Schottky--Klein function for the complex torus is
$$\varpi (\zeta ,z)=(\zeta -z)\prod_{n=1}^\infty {{(1-{\hbox{q}}^2(\zeta /z+z/\zeta )+{\hbox{q}}^{4n})}\over{(1-{\hbox{q}}^{2n})^2}},\eqno(7.23)$$
\noindent which strongly resembles the formula for $\vartheta_3$ as given in (7.2), and the analogy with the product for $\vartheta_1$ is closer still; see [48 p135], and [18]. By contrast, for $g\geq 2$, the universal cover of ${\cal E}$ is the hyperbolic upper half plane ${\bf C}_+$ rather than ${\bf C}$, and $\Gamma$ is a non abelian group.\par
\indent (ii)  The function $z\mapsto \varpi (z, \zeta )$ is holomorphic and has a first order zero at $z=\zeta$, and no other zeros in ${\cal F}$; hence $\varpi (z,\zeta )$ can be used to build meromorphic functions on ${\cal E}$ with zeros and poles, subject to Abel's theorem. Baker [7] describes the functional equation of $\varpi$, viewing $z\mapsto \varpi (z, \zeta )$  as an automorphic function with respect to the action of $\Gamma$.\par 
\indent (iii) By Koebe's retrosection theorem, any compact Riemann surface can be uniformized by a Schottky group, as in the preceding discussion, where $\Gamma$ is some discrete subgroup of $PSL(2, {\bf C})$.
 In Proposition 7.5, we have chosen $\Gamma$ to be a Fuchsian group so as to simplify the geometry and to identify a class of Schottky groups for which the product (7.13) converges; in this case, one can choose the $\hat\gamma_j$ to be circles with centres on ${\bf R}$.  The Riemann surface ${\cal E}$ is not to be confused with ${\bf C}_+/\Gamma$. There are alternative definitions of $\varpi$ which avoid the infinite product (7.13); see [5] for an existence theorem for $\varpi$ which is based on potential theory.\par
\indent (iv)  Fay and Mumford [23, 57] introduce $\varpi$ via Fay's prime form, and their theory leads naturally to theta functions on the Jacobi variety of ${\cal E}$.
 Let ${\cal E}$ be a symmetric compact Riemann surface of genus $g$. Then ${\cal E}$ has a prime form $E$, and we suppose that $\varphi :{\cal E}\rightarrow {\cal E}$ satisfies
$$E(\varphi (\zeta ), \varphi (z))=\overline {E(\zeta , z)}.\eqno(7.24)$$
\noindent Then the Schottky--Klein function $\varpi : {\cal E}\times {\cal E}\rightarrow {\bf C}$ satisfies 
$$E(\zeta ,z)={{\varpi (\zeta ,z)}\over{\sqrt{d\zeta dz}}},\eqno(7.25)$$
\noindent and the functional relations
$$\varpi (\zeta ,z)=-\varpi (z, \zeta ), \quad \overline {\varpi (1/\bar \zeta , 1/\bar z)}={{-\varpi (\zeta ,z)}\over{\zeta z}}.\eqno(7.26)$$
\vskip.05in
\noindent {\bf Example 7.7.} The following configuration was studied by Burnside [16]. Let $1<e_1<e_2<\dots <e_{2g-1}<e_{2g}$, and let $\hat\gamma_{-j}$ be the circle with centre $(e_{2j-1}+ e_{2j})/2$ and radius $(	e_{2j}-e_{2j-1})/2$; then we let $\hat\gamma_0$ be the circle with centre $0$ and radius $1$, and let $\hat\gamma_j$ be the inversion of $\hat \gamma_{-j}$ in $\hat\gamma_0$ for $j=1, \dots, g$. This gives us $2g+1$ circles with centres on ${\bf R}$ such that $\hat \gamma_{\pm j}$ for $j=1, \dots ,g$ are $2g$ circles with mutually disjoint interiors, and we can construct the Schottky surface ${\cal E}$ as above. Baker [7] shows that the functions
$$X (\zeta )=\Bigl( {{\varpi (\zeta ,-1)\varpi (1,\zeta_0)}\over{\varpi (\zeta, 1)\varpi (-1, \zeta_0)}}\Bigr)^2,$$
$$Y(\zeta )={{\varpi (\zeta ,-1)\varpi (1,\zeta_0)}\over{\varpi (\zeta, 1)\varpi (-1, \zeta_0)}}\prod_{j=1}^g {{\varpi (\zeta ,e_{2j-1})\varpi (e_{2j},\zeta_0)}\over{\varpi (\zeta, e_{2j})\varpi (e_{2j-1}, \zeta_0)}}\eqno(7.27)$$
\noindent are meromorphic and invariant under the action of $\Gamma$, so $Y(T(\zeta ))=Y(\zeta )$ for all $T\in\Gamma$.  Hence  $X$ defines a rational function of degree two on ${\cal E}\rightarrow {\bf C}_\infty$, with double pole at $\zeta=1$ and a double zero at $\zeta =-1$, hence ${\cal E}$ is hyperelliptic. To identify the corresponding curve,  let $\lambda_j=X(e_j)$, then $Z=X\prod_{j=1}^g  (X-\lambda_{2j-1}).$ Baker [7] shows that ${\cal E}$ is the closure of the complex curve
$$\Bigl\{ (X,Z)\in {\bf C}^2: Z^2=X\prod_{j=1}^{2g}(X-\lambda_j)\Bigr\}.\eqno(7.28)$$
\noindent Renumbering the $\lambda_j$ if needs be, we can consider $(\lambda_{2j-1}, \lambda_{2j})$ as consecutive gaps. By Abel's theorem, there exists a discrete subgroup $\Lambda$ of ${\bf C}^g$ of full rank such that the map
$$J:\delta_{\zeta}-\delta_{\zeta_0}\mapsto \Bigl( \int_{\zeta_0}^\zeta {{X^{j-1}dX}\over{Z}}\Bigr)_{j=1}^g\in {\bf C}^g\qquad ({\hbox{mod}}\quad \Lambda )\eqno(7.29)$$
\noindent extends via ${\bf Z}$-linear combinations and defines a surjective group homomorphism from the divisors on ${\cal E}$ of degree zero onto the quotient group ${\bf C}^g/\Lambda $. Thus one can identify the Picard variety with ${\bf C}^g/\Lambda$ via the Jacobian map $J$, which is an isomorphism of abelian groups.\par
\vskip.05in
\indent Marcenko showed that the set of finite gap potentials is norm dense in $L^2([0, \pi ]; {\bf R})$; see[25] for more on spectral gaps for $L^2$ potentials.\par
\vskip.05in
\noindent {\bf 8.  Kadomtsev--Petviashvili differential equations}\par
\vskip.05in
\noindent  For a meromorphic complex function $u(x,y,s)$, the $KP$ equation is 
$${{\partial}\over{\partial x}}\Bigl( {{\partial^3 u}\over{\partial
x^3}}-6u{{\partial u}\over{\partial x}}+4\lambda {{\partial u}\over{\partial x}}
+4\alpha {{\partial u}\over{\partial s}}\Bigr)
 +3\beta^2{{\partial^2u}\over{\partial y^2}}=0,\eqno(8.1)$$
\noindent where the $\alpha, \beta, \lambda\in {\bf C}$ are parameters.  
 In this section, we use $(2,2)$-admissible linear systems, one can 
produce solutions to the $KP$ equations via scattering functions and 
the Gelfand--Levitan equation. Zakharov and Shabat [72] considered the
associated scattering function $\phi$, which we take to be a meromorphic complex function
$\phi (x,y,z,t)$ that satisfies the linear $KP$
equations 
$$\alpha{{\partial \phi}\over{\partial t}} +{{\partial^3\phi}\over{\partial
x^3}}+{{\partial^3\phi}\over{\partial z^3}}+\lambda \Bigl( 
{{\partial \phi}\over{\partial x}}+
{{\partial \phi}\over{\partial z}}\Bigr)=0\eqno(8.2)$$
\noindent and
$$\beta {{\partial \phi}\over{\partial y}}+{{\partial^2\phi}\over{\partial
x^2}}-{{\partial^2\phi}\over{\partial z^2}}=0.\eqno(8.3)$$

\vskip.05in

\noindent {\bf Definition} {\sl ($GL$ equation for $KP$).} The appropriate version of the 
Gelfand--Levitan equation for the linear $KP$ equations is 
$$\phi (x,z;y,t)+K(x,z;y,t)+\int_x^\infty K(x,s;y,t)\phi (s,z;y,t)ds
=0\qquad (x<z).\eqno(8.4)$$
\noindent  For a solution $K(x,z;y,t)$, define the potential 
by 
$$u(x;y,t)=-2{{d}\over{dx}}K(x,x;y,t).\eqno(8.5)$$
\vskip.05in
\noindent  {\bf Remark.} Note that in comparison with [72], our
potential has an extra minus sign. We regard $(x,z)$ as the main variables, $(x,x)$ as the diagonal and $(y,t)$ as parameters which describe the deformation of solutions of the integral equation. We write ${{\partial}\over{\partial x}}$ to indicate differentiation with respect to the first variable, ${{\partial}\over{\partial z}}$ to indicate differentiation with respect to the second variable, and ${{d}\over{dx}}$ to indicate differentiation along the diagonal.\par
\vskip.05in
\indent We shall introduce a suitable family of admissible linear systems such that  $u$ arises as their tau function, and solve the Gelfand--Levitan equation in a similar way to section 2. \par
\vskip.05in 

\noindent {\bf Definition.} Given $(2,2)$ admissible linear systems $(-A_1, B_0, C_0)$ and $(-A_2, B_0, C_0)$ as in Theorem 2.2 with $A_1,A_2\in {\cal L}(H)$. 
We write $(-A_1, -A_2; B_0, C_0)$ for brevity, which is not to be confused with the notation $(A,B,C,D)$ which is used as shorthand  for the colligation matrix of a linear system.  Let
$$C(y;t)=C_0e^{t(A_1^3+\lambda A_1)/\alpha -yA_1^2/\beta }\eqno(8.6)$$
$$B(y;t) =e^{t(A_2^3+\lambda A_2)/\alpha +yA_2^2/\beta}B_0.\eqno(8.7)$$
Then let 
$$\phi (x,z;y;t)=C(y;t)e^{-xA_1}e^{-zA_2}B(y;t),\eqno(8.8)$$
$$R_x=R_x(y,t)=\int_x^\infty e^{-A_2s}B(y;t)C(y;t)e^{-A_1 s}ds,\eqno(8.9)$$
\noindent and 
$$K(x,z;y;t)=-C(y;t)e^{-xA_1}(I+R_x)^{-1}e^{-zA_2}B(y;t).\eqno(8.10)$$
\vskip.05in
\noindent The signs before $A_1^2$ and $A_2^2$ in (8.6) and (8.7) are purposefully different. We do not assume that $A_1$ and $A_2$ commute, so that $\phi$ really does depend upon $y$ in general. When $A_1=A_2$, the formula (8.9) reduces to our usual $R$ operator in the style (1.7), Proposition 8.1 reduces to Proposition 2.4, and $\psi$ reduces to a
Hankel type kernel $\psi (x,z;y,t)=\phi (x+z;t)$, independent of $y$.\par 
\indent  Zakharov and Shabat used this method for finite rank $A_1$ and $A_2$ to 
produce soliton solutions of $KP$. For solutions in the style of Proposition 4.5, see [35, Proposition 14.12]. In order to ensure that various products and brackets  are well defined, we have imposed the condition $A_1,A_2\in {\cal L}(H)$; some of the results hold under less stringent conditions.\par
\vskip.05in 

\noindent {\bf Proposition 8.1.} {\sl  (i) Then $\phi (x,z;y;t)$ satisfies the scattering equations (8.2) and (8.3) for the $KP$ equation.\par
\indent (ii)  $K(x,z;y;t)$ satisfies the integral equation (8.4) and}
$$K(x,x;y;t)={{d}\over{dx}}\log\det (I+R_x).\eqno(8.11)$$
\indent {\sl (iii) there exists $x_0$ such that $K(x,z;y;t)$ and the corresponding $u$ from (8.5) satisfy} 
$${{\partial^2K}\over{\partial x^2}}-{{\partial^2K}\over{\partial z^2}}
+\beta{{\partial K}\over{\partial y}}=u(x;y;t)K(x,z;y;t)\qquad 
(x_0<x<z).\eqno(8.12)$$
\vskip.05in
\noindent {\bf Proof.} (i) Since the operators are all bounded, the functions are differentiable and one can verify the differential equations, without assuming that $A_1$ and $A_2$ commute.\par
\indent (ii) The linear system
$$\hat \Sigma=\Bigl(\left[\matrix{-A_1&0\cr 0&-A_2\cr}\right],
\left[\matrix{B_0&0\cr 0&B(y;t)\cr}\right],
\left[\matrix{0&-C_0\cr C(y;t)&0\cr}\right]\Bigr)\eqno(8.13)$$
\noindent is $(2,2)$ admissible and by Theorem 2.2 the corresponding $\hat R_x$ operator
$$\hat R_x=\left[\matrix{ 0&-\int_x^\infty e^{-sA_1}B_0C_0e^{-sA_2}ds\cr 
\int_x^\infty
e^{-sA_2}B(y;t)C(y;t)e^{-sA_1}ds&0\cr}\right],\eqno(8.14)$$
\noindent is trace class as in Proposition 2.4. Hence 
$$\eqalignno{\phi (x,z;y,t)&+K(x,z;y,t)+\int_x^\infty K(x,s;y,t)\phi (s,z;y,t)ds\cr
&=C(y;t)e^{-xA_1}e^{-zA_2}B(y;t) -C(y;t)e^{-xA_1}(I+R_x)^{-1}e^{-zA_2}B(y;t)\cr
&\quad -C(y;t)e^{-xA_1}(I+R_x)^{-1}\int_x^\infty 
e^{-sA_2}B(y;t)C(y;t)e^{-sA_1}e^{-zA_2}B(y;t)ds\cr
&=C(y;t)e^{-xA_1}\Bigl(I-(I+R_x)^{-1}-(I+R_x)^{-1}R_x\Bigr)
e^{-zA_2}B(y;t)\cr
&=0,&(8.15)\cr}$$
\noindent as in the proof of Proposition 2.4. 
One then verifies the determinant identity (8.11), which  involves the bottom left entry of $\hat R_x$ satisfying the asymmetric Lyapunov equation 
$${{d}\over{dx}}R_x=-A_2R_x-R_xA_1=-e^{-xA_2}B(y;t)C(y;t)e^{-xA_1}
\qquad (x>0).\eqno(8.16)$$
\indent (iii) The solution of the integral equation is unique for large enough $x$ since 
$\Vert e^{-xA_1}\Vert \rightarrow 0$ and $\Vert e^{-xA_2}\Vert \rightarrow 0$ 
exponentially fast as $x\rightarrow\infty$; hence $\Psi (x;z;y;t)\rightarrow 0$ 
exponentially fast as $x\rightarrow\infty$. Using the scattering equation (8.2), 
one shows by differentiating (8.4) repeatedly that 
$${{\partial^2}\over{\partial x^2}}K(x,z;y;t)-{{\partial^2}\over{\partial
z^2}}K(x,z;y;t)+\beta 
{{\partial }\over{\partial y}}K(x,z;y;t)$$ 
\noindent and $u(x;y;t)K(x,z;y;t)$ both satisfy the equation which appears when (8.4) is multiplied by $u(x;y;t)$, and so by uniqueness are equal.\par 

\rightline{$\square$}\par

\vskip.05in
\noindent {\bf Theorem 8.2.} {\sl The potential $u$ associated with $\tau (x,y,t)$ from $(-A_1, -A_2; B(y;t),C(y;t))$  satisfies the $KP$ equation (8.1).}
\vskip.05in
\indent The proof involves a calculation which extends the results of sections 4 and 5, and we split this into two Lemmas.
\vskip.05in
\noindent {\bf Definition} {\sl (Product and bracket).} In the notation of Proposition 8.1, let $F_x=(I+R_x)^{-1}$.  Let ${\cal B}$ be any 
differential ring of
functions from  $(0,\infty)\rightarrow {\cal L}(H_0)$, then let 
$${\cal A}={\hbox{span}}_{\bf C}\bigl\{A_1^{n_1}A_2^{m_1}A_1^{p_1}, A_1^{n_1}A_2^{m_1} A_1^{p_1}F_xA_1^{n_2}\dots F_xA_1^{n_r}A_2^{m_r}A_1^{p_r}: n_j, m_j,p_j\in {\bf
Z}_+\bigr\}\eqno(8.17)$$
\noindent be the algebra generated by $I,A_1, A_2$ and $F$. On ${\cal A}$ we introduce the associative product $\ast$ by 
$$P\ast Q=P(A_1F+FA_2-F(A_1+A_2)F)Q,\eqno(8.18)$$
which is distributive over the standard addition, and the derivation ${\partial }: {\cal A}\rightarrow {\cal A}$ by
$$\partial P=(A_2-(A_1+A_2)F)P+{{dP}\over{dx}}+P(A_1-F(A_1+A_2)),\eqno(8.19)$$
\noindent Then let the bracket $\lfloor\,\cdot\,\rfloor :{\cal A}\rightarrow {\cal B}$ be 
the linear map
$$\lfloor Y\rfloor =Ce^{-xA_1}F_xYF_xe^{-xA_2}B \qquad (Y\in {\cal A}).\eqno(8.20)$$
\vskip.05in
\noindent {\bf Lemma  8.3.} {\sl Then $({\cal A}, \ast , \partial )$ is a
differential ring,  and the bracket gives a 
homomorphism of differential rings 
$\lfloor\,\cdot\,\rfloor: ({\cal A}, \ast , \partial )\rightarrow ({\cal B}, \cdot , d/dx)$.}\par
\vskip.05in

\noindent {\bf Proof.}  The basic observation is that
$dF/dx=A_1F+FA_2-F(A_1+A_2)F$, so one can check that
$${{d}\over{dx}}\bigl( A_1F+FA_2-F(A_1+A_2)F\bigr) =\bigl(A_1-F(A_1+A_2)\bigr)(\bigl( A_1F+FA_2-F(A_1+A_2)F\bigr)$$
$$+\bigl( A_1F+FA_2-F(A_1+A_2)F\bigr)\bigl( A_2-(A_1+A_2)F\bigr)\eqno(8.21)$$
\noindent so that
$$\eqalignno{ {{d}\over{dx}}P\bigl( A_1F+FA_2&-F(A_1+A_2)F\bigr)Q\cr
&= {{dP}\over{dx}}\bigl( A_1F+FA_2-F(A_1+A_2)F\bigr)Q\cr
&\quad+P\bigl( A_1-F(A_1+A_2)\bigr)\bigl( A_1F+FA_2-F(A_1+A_2)F\bigr)Q\cr
&\quad +P\bigl( A_1F+FA_2-F(A_1+A_2)F\bigr)\bigl( A_2-(A_1+A_2)F\bigr)Q\cr
&\quad+P\bigl( A_1F+FA_2-F(A_1+A_2)F\bigr){{dQ}\over{dx}},&(8.22)\cr}$$
\noindent so by adding the terms at either end of this expression, one shows that 
$$\partial (P\ast Q)=(\partial P)\ast Q+P\ast (\partial Q).\eqno(8.23)$$
\indent Now we consider the bracket, and find from Lyapunov's equation that
$$\eqalignno{\lfloor P\rfloor \lfloor Q\rfloor &=Ce^{-xA_1}FPFe^{-xA_2}BCe^{-xA_1}FQFe^{-xA_2}B \cr
&=Ce^{-xA_1}FPF(A_2S+SA_1)FQFe^{-xA_2}B\cr
&=Ce^{-xA_1}FP(A_1F+FA_2)-F(A_1+A_2)F) QFe^{-xA_2}B\cr
&=\lfloor P\ast Q\rfloor;&(8.24)\cr}$$
\noindent and 
$$\eqalignno{{{d}\over{dx}}\lfloor P\rfloor&= {{d}\over{dx}} Ce^{-xA_1}FPFe^{-xA_2}B\cr
&= Ce^{-xA_1}F\Bigl(\bigl(  A_2-(A_1+A_2)F\bigr)P+{{dP}\over{dx}}+P\bigl(A_1-F(A_1+A_2)\bigr)\Bigr)Fe^{-xA_2}B\cr
&=\lfloor \partial P\rfloor.&(8.25)\cr}$$ 
\rightline{$\square$}\par
\vskip.05in
\indent With the usual operator multiplication, let ${\cal A}_0$ be the subalgebra of ${\cal A}$ that is generated by $I,A_1$, and $A_2$. Let ${\bf J}$ be the ideal in ${\cal A}$ generated by $F$; then the powers ${\bf J}^n$ give a decreasing chain of ideals such that $\cap_{n=1}^\infty {\bf J}^n=\{0\}$; any product including $n$ factors of $F$ belongs to ${\bf J}^n$. Now ${\cal A}/{\bf J}$ is isomorphic as an algebra to ${\cal A}_0$, and $\partial ({\bf J}^n)\subseteq {\bf J}^n$, so there are induced maps 
$\tilde \partial :({\bf J}^n/ {\bf J}^{n+1})\rightarrow ({\bf J}^n/ {\bf J}^{n+1})$; in particular $\tilde \partial :{\cal A}/{\bf J}\rightarrow {\cal A}/{\bf J}$ may be identified with $P\mapsto A_2P+PA_1$ on  ${\cal A}_0$. Thus we may regard ${\cal A}$ as a graded algebra consisting of polynomials in $F$ with a (noncommutative) algebra ${\cal A}_0$ of coefficients. We regard  $\partial$ as the sum of the derivation $d/dx$, the multiplications $X\mapsto A_2X+XA_1$ which typically preserve the degree, and the multiplications $X\mapsto -(A_1+A_2)FX$ and $X\mapsto -XF(A_1+A_2)$ which can raise the degree by one at most.\par

\vskip.05in

\noindent {\bf Lemma 8.4.} {\sl The function 
$$w=4\alpha {{\partial u}\over{\partial t}} +4\lambda {{\partial u}\over{\partial x}}+ 
{{\partial^3 u}\over{\partial x^3}}-6u{{\partial u}\over{\partial x}}\eqno(8.26)$$
\noindent is the image under the bracket $\lfloor\cdot\rfloor$ of the sum of terms} 
$$  6A_1^4-12A_2^2+6A_2^4$$
$$ +12 (A_1+A_2)FA_2(A_1-A_2)^2-12 (A_1^2-A_2^2)A_1F(A_1+A_2)$$
$$+ 6(A_1^2+2A_2A_1)+A_2)^2F(A_1^2-A_2^2)-6(A_1^2-A_2^2)F(A_1^2+2A_2A_1+A_2^2)$$
$$+12(A_1^2-A_2^2)F(A_1+A_2)F(A_1+A_2)-12(A_1+A_2)F(A_1+A_2)F(A_1^2-A_2^2).\eqno(8.27)$$
\vskip.05in 
\noindent {\bf Proof.} We introduce the ordered products, in which powers of $A_2$ powers occur before powers of $A_1$, 
$$A^{(1)}=A_1+A_2,\quad A^{(2)}=A_1^2+2A_2A_1+A_2^2,\quad A^{(3)}=A_1^3+3A_2A_1^2+3A_2^2A_1+A_2^3, \dots \eqno(8.28)$$
\noindent and the coefficients are as in Pascal's triangle.  It suffices to compute
$$W=4\alpha \partial_tU+4\lambda \partial U +\partial^3 U-3U\ast \partial U-3\partial U\ast U,\eqno(8.29)$$
\noindent since $w=\lfloor W\rfloor$ by Lemma 8.3. The terms from (8.29) are given in the following multiplication table (8.30). Starting with $U=-2A^{(1)},$ we have derivatives
$$\eqalignno{\partial U&= -2 A^{(2)}+4A^{(1)}FA^{(1)},\cr
\partial^2 U&=-2 A^{(3)}+6A^{(2)}FA^{(1)}+6A^{(1)}FA^{(2)} -12 A^{(1)}FA^{(1)}FA^{(1)},\cr
\partial^3U&=-2A^{(4)}+8 A^{(3)}FA^{(1)}+8A^{(1)}FA^{(3)} +12A^{(2)}FA^{(2)}\cr
&\quad  -24 A^{(1)}FA^{(1)}FA^{(2)}-24 A^{(1)}FA^{(2)}FA^{(1)}-24 A^{(2)}FA^{(1)}FA^{(1)}\cr
&\quad +48A^{(1)}FA^{(1)}FA^{(1)}FA^{(1)};&(8.30)\cr}$$
\noindent these products exhibit a high degree of symmetry. The only proof known to the authors is applying $\partial$ repeatedly, then  patiently multiplying out and gathering the various products. To respect the symmetry of terms, we use
$$\eqalignno{U\ast \partial U+\partial U \ast U&= 4A^{(1)}FA_2A^{(2)}+ A^{(2)} A_1FA^{(1)}+A^{(2)}FA_2A^{(1)} +A^{(1)}A_1FA^{(2)}\cr
&\quad-8A^{(2)}FA^{(1)}FA^{(1)}-8A^{(1)}FA^{(2)}FA^{(1)}-8A^{(1)}FA^{(1)}FA^{(2)}\cr
&\quad-4(A_1^2-A_2^2)FA^{(1)}FA^{(1)} +4A^{(1)}FA^{(1)}F(A_1^2-A_2^2)\cr
&\quad+16^{(1)}FA^{(1)}FA^{(1)}FA^{(1)}.&(8.31)\cr}$$
We likewise introduce
$$\eqalignno{\alpha\partial_tU&=2(A_2^3+\lambda A_2)(A_1+A_2)+2(A_1+A_2)(A_1^3+\lambda A_1)\cr
& -2(A_1^3+A_2^3+\lambda A_1+\lambda A_2)F(A_1+A_2)\cr
&\quad -2(A_1+A_2)F(A_1^3+A_2^3+\lambda A_1+\lambda A_2)&(8.32)\cr}$$
\noindent and 
$$\eqalignno{\beta^2\partial_y^2 U&=-2(A_1^5+A_2A_1^4-2A_2^2A_1^3-2A_2^3A_2^2+A_2^4A_1+A_2^5) \cr
& \quad +2(A_1+A_2)F(A_1^4-2A_1^2A_2^2+A_2^4)+2(A_1^4-2A_1^2A_2^2+A_2^4)F(A_1+A_2)\cr
&\quad+4(A_1^2-A_2^2)F (A_1^3+A_2A_1^2-A_2^2A_1-A_2^3)\cr
&\quad+4(A_1^3+A_2A_1^2-A_2^2A_1-A_2^3)F(A_1^2-A_2^2) \cr
&\quad-4 (A_1^2-A_2^2)F(A_1+A_2)F(A_1^2-A_2^2)\cr
& \quad-4(A_1+A_2)F(A_1^2-A_2^2)F(A_1^2-A_2^2) \cr
&\quad-4 (A_1^2-A_2^2)F(A_1^2-A_2^2)F(A_1+A_2).&(8.33)\cr}$$ 
\noindent Then one checks that $W$ reduces to the combination (8.29).\par
\rightline{$\square$}\par
\vskip.05in
\noindent {\bf Proof of Theorem 8.2} Let $\Theta =\{ X\in {\cal A}: \lfloor X\rfloor =0\}$ and observe that $\Theta$ contains the commutator subspace spanned by $Q\ast P-P\ast Q$. 
The final two terms in $W$ are of degree two in $F$, which would give terms of degree three in  $\partial W$, which do not appear in the formula for $\partial_y^2U$. Hence we replace them by terms of degree one, before differentiating; or equivalently, we show that $4\beta^2\partial_y^2U+\partial W$ belongs to $\Theta$. We have
$$\eqalignno{0&=\lfloor A_1^2-A_2^2\rfloor \lfloor A_1+A_2\rfloor -\lfloor A_1+A_2\rfloor\lfloor A_1^2-A_2^2\rfloor\cr
&= \lfloor (A_1^2-A_2^2)(A_1F+FA_2-F(A_1+A_2)F) (A_1+A_2) \rfloor\cr
&\quad - \lfloor (A_1+A_2)(A_1F+FA_1-F(A_1+A_2)F(A_1^2-A_2^2)\rfloor\cr
&= \lfloor (A_1^2-A_2^2)(A_1F+FA_2)(A_1+A_2) -(A_1+A_2)(A_1F+FA_2)(A_1^2-A_2^2)\rfloor\cr
&\quad  -\lfloor (A_1^2-A_2^2)F(A_1+A_2)F(A_1+A_2) -(A_1+A_2)F(A_1+A_2)F(A_1^2-A_2^2)\rfloor.&(8.34)\cr}$$
\noindent So when we replace the final two terms of degree two by terms such as $(A_1^2-A_2^2)(A_1F+FA_2)(A_1+A_2)$, we obtain the following collection of terms of degree one 
$$12(A_1+A_2)FA_2(A_1^2-A_2^2)-12(A_1^2-A_2^2)A_1F(A_1+A_2) $$
$$+6(A_1^2+2A_2A_1+A_2^2)F(A_1^2-A_2^2) - 6(A_1^2-A_2^2)F(A_1^2+2A_2A_1+A_2^2)$$
$$+ 6(A_1^2+2A_2A_1)+A_2)^2F(A_1^2-A_2^2)-6(A_1^2-A_2^2)F(A_1^2+2A_2A_1+A_2^2)$$
$$+12 (A_1^2-A_2^2)A_1F(A_1+A_2)+12 (A_1^2-A_2^2)FA_2(A_1+A_2)$$
$$-12 (A_1+A_2)A_1F(A_1^2-A_2^2)-12(A_1+A_2)FA_2(A_1^2-A_2^2)$$
$$=-12 (A_1^2-A_2^2)F(A_1^2-A_2^2).\eqno(8.35)$$
Now we compute
$$\eqalignno{\partial \bigl( 6A_1^4-&12A_2^2+6A_2^4 -12 (A_1^2-A_2^2)F(A_1^2-A_2^2)\bigr)\cr
&=6A_2A_1^4-12 A_2^3A_1^2+6A_2^5+6A_1^5 -12 A_2^2A_1^3 +6A_2^4A_1\cr
&\quad -(A_1+A_2)F(6A_1^4-12A_2^2A_1^2+6A_2^2) -(6A_1^4-12 A_2^2A_1^2-A_2^4) F(A_1+A_2)\cr
& \quad -12 A_2(A_1^2-A_2^2)F(A_1^2-A_2^2) -12(A_1^2-A_2^2)F(A_1^2-A_2^2)A_1\cr
&\quad +12 (A_1+A_2)F(A_1^2-A_2^2)F(A_1^2-A_2^2)+12(A_1^2-A_2^2)F(A_1^2-A_2^2)F(A_1+A_2)\cr
&\quad -12 (A_1^2-A_2^2)(A_1F+FA_2-F(A_1+A_2)F)(A_1^2-A_2^2)\cr
&= 6A_1^5+6A_2A_1^4-12 A_2^2A_1^3-12 A_2^3A_1^3+6A_2A_1^4+6A_1^5\cr
&\quad -6(A_1+A_2)F(A_1^4-2A_2^2A_1^2+A_1^4) -6(A_1^4-2A_2^2A_1^2+A_1^4)F(A_1+A_2)\cr
&\quad+(-12(A_2A_1^2+12A_2^3-12A_1^3+12A_2^2A_1)F(A_1^2-A_2^2) \cr
&\quad+(A_1^2-A_2^2)F(-12 A_1^3+12A_2^2A_1-12 A_2A_1^2+12A_1^3)\cr
&\quad+12 (A_1+A_2)F(A_1^2-A_2^2)F(A_1^2-A_2^2) +12 (A_1^2-A_2^2) F(A_1^2-A_2^2)F(A_1+A_2)\cr
&\quad +12 (A_1^2-A_2^2) F(A_1+A_2)F(A_1^2-A_2^2).&(8.36)\cr}$$
\noindent By comparing this with (8.27) we obtain the result.\par
\rightline{$\square$}\par
\vskip.05in
\noindent {\bf 9. The Baker--Akhiezer function for $KP$}\par
\vskip.05in
\indent In this final section, we obtain solutions to the time dependent Schr\"odinger equation in the form of quotients of tau functions for a family of admissible linear systems. When $u$ satisfies $KP$, one can choose $w$ so that the operators 
$$\beta {{\partial }\over{\partial y}}+L=\beta {{\partial }\over{\partial y}}-{{\partial^2}\over{\partial x^2}}+u(x;y,t)\eqno(9.1)$$
\noindent and 
$${{\partial}\over{\partial t}}+ M={{\partial}\over{\partial t}}+ {{\partial^3}\over{\partial x^3}}-{{3}\over{2}}u (x;y,t){{\partial}\over{\partial x}}-{{3}\over{4}} {{\partial u}\over{\partial x}} -3\alpha w(x;y,t)\eqno(9.2)$$
\noindent commute. In the following result, we obtain an explicit form for a common eigenfunction for both these operators, so \par
$$\bigl(\beta {{\partial }\over{\partial y}}+L\bigr)\psi_\zeta =0=\bigl({{\partial}\over{\partial t}}+ M\bigr)\psi_\zeta.\eqno(9.3) $$
\noindent By analogy with (3.2), we call a particular family of solutions the Baker--Akhiezer function.\par
\vskip.05in
\noindent {\bf Definition} {\sl (Baker--Akhiezer function).} Consider the linear system $(-A_1, A_2; B(y,t), C(y,t))$ from Theorem 8.2, with spectral parameter $\zeta$, and $R_x=R_x(y,t)$ as in (8.9). Then the Baker--Akhiezer function is
$$\psi_\zeta (x; y,t)=e^{\zeta x-\zeta^2y/\beta -\zeta^3t}{{\det (I+R_x(\zeta I+A_1)(\zeta I-A_2)^{-1})}\over{\det (I+R_x)}},\eqno(9.4)$$
\noindent defined on $\zeta\in {\bf C}\setminus {\hbox{Spec}}(A_2)$.\par
\vskip.05in
\noindent {\bf Proposition 9.1.} {\sl Let $\alpha w(x;y,t)={{\partial K}\over{2\partial y}}(x,x;y,t)$.\par
\indent (i)  Then  $\beta {{\partial }\over{\partial y}}+L$ and ${{\partial}\over{\partial t}}+ M$ commute;\par
 \indent (ii) The Baker--Akhiezer function $\psi_\zeta$ is meromorphic in $(x,y,t)\in {\bf C}^3$;\par
\indent (iii) also $\psi_\zeta $ satisfies
$$-{{\partial^2\psi_\zeta (x;y,t)}\over{\partial x^2}}+u(x;y,t)\psi_\zeta (x;y,t)+\beta {{\partial\psi_\zeta (x;y,t)}\over{\partial y}}=0,\eqno(9.5)$$
\indent  (iv) and satisfies} 
$${{\partial\psi_\zeta}\over{\partial t}}+ {{\partial^3\psi_\zeta}\over{\partial x^3}}-{{3}\over{2}}u (x;y,t){{\partial\psi_\zeta}\over{\partial x}}-{{3}\over{4}} {{\partial u}\over{\partial x}}\psi_\zeta  (x;y,t) -3w(x;y,t)\psi_\zeta (x;y,t)=0.\eqno(9.6)$$
\vskip.05in
\noindent {\bf Proof}  (i) We have
$$\alpha {{\partial w}\over{\partial x}}(x;y,t)= {{1}\over{2}}{{\partial }\over{\partial y}}{{dK}\over{dx}} (x,x;y,t)={{1}\over{4}}{{\partial u}\over{\partial y}}(x;y,t),\eqno(9.7)$$
\noindent which is what one needs to make the operators commute.\par
 \indent (ii) The function $(x,y,t)\mapsto R_x$ is entire from ${\bf C}^3\rightarrow {\cal L}^1(H)$, hence 
$$\det (I+(\zeta I+A_1)R_x(\zeta I-A_2)^{-1})\eqno(9.8)$$
\noindent is entire and $\psi_\zeta$ is a quotient of entire functions, hence is meromorphic.\par
\indent (iii)  By some simple manipulations of the determinants, we have 
$$\eqalignno{\psi_\zeta (x; y,t)&=e^{\zeta x-\zeta^2y/\beta -\zeta^3t}{{\det (I+((\zeta I-A_2)R_x+(A_2R_x+R_xA_1) (\zeta I-A_2)^{-1})}\over{\det (I+R_x)}}\cr
&=e^{\zeta x-\zeta^2y/\beta -\zeta^3t}{{\det (I+R_x+(\zeta I-A_2)^{-1}(A_2R_x+R_xA_1) )}\over{\det (I+R_x)}}\cr
&=e^{\zeta x-\zeta^2y/\beta -\zeta^3t}\det (I+ (I+R_x)^{-1}(\zeta I-A_2)^{-1}(A_2R_x+R_xA_1) );&(9.9)\cr}$$
\noindent so by Lyapunov's equation, we have
$$\eqalignno{\psi_\zeta (x;y,t)&=e^{\zeta x-\zeta^2y/\beta -\zeta^3t}\det (I+ (I+R_x)^{-1}(\zeta I-A_2)^{-1}e^{-xA_2}B(y,t)C(y,t)e^{-xA_1} )\cr
&=e^{\zeta x-\zeta^2y/\beta -\zeta^3t}\Bigl( 1+{\hbox{trace}} (I+R_x)^{-1}(\zeta I-A_2)^{-1}e^{-xA_2}B(y,t)C(y,t)e^{-xA_1} \Bigr)\cr
&=e^{\zeta x-\zeta^2y/\beta -\zeta^3t}\Bigl( 1+C(y,t)e^{-xA_1}(I+R_x)^{-1}(\zeta I-A_2)^{-1}e^{-xA_2}B(y,t)\Bigr)&(9.10)\cr}$$
\noindent since $B(y,t)C(y,t)$ has rank one; then we write this as 
$$\eqalignno{\psi_\zeta (x;y,t)&=e^{\zeta x-\zeta^2y/\beta -\zeta^3t}\Bigl( 1-C(y,t)e^{-xA_1}(I+R_x)^{-1}\int_x^\infty e^{-zA_2} e^{\zeta (z-x)}\, dzB(y,t)\Bigr)\cr
&= e^{\zeta x-\zeta^2y/\beta -\zeta^3t}\Bigl( 1+\int_x^\infty K(x,z;y,t) e^{\zeta (z-x)}\, dz\Bigr).&(9.11)\cr}$$
\noindent Now we calculate
$${{\partial \psi_\zeta}\over{\partial y}}=-(\zeta^2/\beta )e^{\zeta x-\zeta^2y/\beta -\zeta^3t}+\int_x^\infty {{\partial }\over{\partial y}}K(x,z,y,t)e^{\zeta z-\zeta^2y/\beta -\zeta^3t}\, dz\eqno(9.12)$$
\noindent and 
$${{\partial \psi_\zeta}\over{\partial x}}=\zeta e^{\zeta x-\zeta^2y/\beta -\zeta^3t}-e^{\zeta x-\zeta^2y/\beta -\zeta^3t}K(x,x;y,t)+\int_x^\infty {{\partial }\over{\partial x}}K(x,z;y,t)e^{\zeta z-\zeta^2y/\beta -\zeta^3t}\, dz\eqno(9.13)$$
\noindent hence
$$\eqalignno{{{\partial^2 \psi_\zeta}\over{\partial x^2}}&=\zeta^2 e^{\zeta x-\zeta^2y/\beta -\zeta^3t}-\zeta e^{\zeta x-\zeta^2y/\beta -\zeta^3t} K(x,x;y,t)+e^{\zeta x-\zeta^2y/\beta -\zeta^3t}{{d}\over{dx}}K(x,x;y,t)\cr
&\quad -{{\partial }\over{\partial x}}K(x,x;y,t)e^{\zeta x-\zeta^2y/\beta -\zeta^3t}+  \int_x^\infty {{\partial^2 }\over{\partial x^2}}K(x,z;y,t)e^{\zeta x-\zeta^2y/\beta -\zeta^3t}\, dz.&(9.14)\cr}$$
\noindent Integrating by parts, we obtain
$$\eqalignno{\int_x^\infty&{{\partial^2}\over{\partial z^2}}K(x,z;y,t)e^{\zeta z-\zeta^2y/\beta -\zeta^3t}dz\cr
&=-{{\partial}\over{\partial z}}K(x,x;y,t)e^{\zeta x-\zeta^2y/\beta -\zeta^3t} -\zeta \int_x^\infty{{\partial }\over{\partial z}} K(x,z;y,t) e^{\zeta z-\zeta^2y/\beta -\zeta^3t}\, dz\cr
&=-{{\partial}\over{\partial z}}K(x,x;y,t)e^{\zeta z-\zeta^2y/\beta -\zeta^3t} +K(x,x;y,t)e^{\zeta x-\zeta^2y/\beta -\zeta^3t}\cr
&\qquad +\zeta^2 \int_x^\infty K(x,z;y,t) e^{\zeta x-\zeta^2y/\beta -\zeta^3t}\, dz;&(9.15)\cr}$$
\noindent recalling the  (8.10), and the definition of $u(x;y,t)$, we deduce the differential equation for $\psi_\zeta$.\par
\indent (iv) One starts with
$$\psi_\zeta (x,y,t)=e^{\zeta x-\zeta^2y/\beta -\zeta^3t}+\int_x^\infty e^{\zeta z-\zeta^2y/\beta -\zeta^3t}K(x,z;y,t)dz.\eqno(9.16)$$
\noindent Then by manipulating the Gelfand--Levitan equation, one deduces (9.6).\par

\rightline{$\square$}
\vskip.05in
\indent Let $(\zeta_j)$ be a sequence of distinct complex numbers and $(\psi_{\zeta_j})_{j=1}^\infty$ a corresponding sequence of distinct solutions of the pair of equations (9.4) and (9.5), where $u$ is as in (8.5) and fixed. Then, taking derivatives in the $x$-variable, one forms the Wronskian
$$\Delta_n={\hbox{Wr}}\bigl( \psi_{\zeta_1}, \dots , \psi_{\zeta_n}\bigr)\eqno(9.17)$$
\noindent and introduces a sequence of new potentials and new Baker--Akhiezer functions by 
$$u_n(x,y,t)=u(x,y,t)-2{{\partial^2}\over{\partial x^2}}\log   \Delta_n(x,y,t)\eqno(9.18)$$
$$\Psi_n(x,y,t)={{\Delta_{n+1}(x,y,t)}\over{\Delta_n(x,y,t)}}.\eqno(9.19)$$
\vskip.05in
\noindent {\bf Corollary 9.2.} (Matveev) {\sl  Then  $\Psi_n$ and $u_n$ satisfy (9.5), and the corresponding (9.6).}\par
\vskip.05in 
\noindent {\bf Proof.} Matveev [46] showed that this follows from Proposition 9.1 by direct calculation.\par
\rightline{$\square$}\par
\vskip.05in
\noindent {\bf Remarks 9.3} (i) Corollary 9.2 enables us to generate a sequence of solutions $(u_n)$ of the $KP$ equation. If all the  $(\psi_{\zeta_j})$  belong to a differential field $({\cal F}, {{\partial}\over{\partial x}}, {{\partial}\over{\partial y}}, {{\partial}\over{\partial t}})$, then $(u_n)$ and $(\Psi_n)_{n=1}^\infty$ also belong to ${\cal F}$. The case in which $u=0$ and $\psi_\zeta (x; y,t)=e^{\zeta x-\zeta^2y/\beta -\zeta^3t}$ gives soliton solutions to $KP$.\par
\indent (ii) If the potential $u$, which appears as a coefficient  of $L$ in Proposition 9.1 does not depend upon $y$, then we can reduce Proposition 9.1(i) to a Lax equation ${{\partial L}\over{\partial t}}=[L,M].$ The results of this section are applicable even when the determinant quotient (9.4) indeed depends upon $y$.\par   
\indent (iii) Krichever and Novikov [41, 42]  consider Baker--Akhiezer functions 
$\psi_\zeta (x,t)$ that are meromorphic with respect to the spectral parameter
 $\zeta$ and produce examples based upon quasi-periodic theta functions;  
in particular, the function  $\psi_2(x, \zeta )$ of section 7 is meromorphic; 
see [39, 32]. Proposition 9.1 does not assert that $\zeta\mapsto \tau_\zeta (x,t)$ is meromorphic 
on ${\bf C}$ and the case 
when the state space $H$ has infinite dimension is problematic. In the rest of this section we circumvent this problem by introducing infinitely many time variables, and acting on the linear systems with a type of infinite dimensional Lie group.
However, Segal and Wilson [64, Proposition 6.11]  
have identified $\tau_\zeta$ functions that are meromorphic. \par
\vskip.05in
\noindent {\bf Definition} (i) For fixed $A_1, A_2\in {\cal L}(H)$, let ${\bf \Sigma}_{A_1, A_2}$ be the set of $\Sigma= (-A_1, -A_2; B, C)$ that give a $(2,2)$ admissible linear system, where  $B:{\bf C}\rightarrow H$ and $C:H\rightarrow {\bf C}$ vary. Let ${\bf C}^\infty_\eta=\{ (a_j)_{j=0}^\infty\in {\bf C}^\infty : \lim\sup_{j\rightarrow\infty} \vert a_j\vert^{1/j}\leq\eta \}$ be the space of coefficients of complex power series with radius of convergence greater than or equal to $1/\eta$, which may be identified with the algebra of holomorphic functions $D(0, \eta )\rightarrow {\bf C}$.  Let 
$$V(t)=\exp\Bigl({\sum_{j=1}^\infty t_jA^j_2}\Bigr),\qquad W(t)= \exp\Bigl({-\sum_{j=1}^\infty t_j(-A_1)^j}\Bigr)\qquad (t=(t_j)\in {\bf C}^\infty_0).\eqno(9.20)$$
\noindent and extend to $t\in {\bf C}^\infty_\eta$ with $\eta>0$ when the series converge absolutely. There is an action $\rho$ of ${\bf C}^\infty_0$ on ${\bf \Sigma}_{A_1, A_2}$ which is given by 
$$ \rho (t):\bigl(-A_1, -A_2, B, C\bigr)\mapsto  \bigl(-A_1, -A_2, V(t) B, CW(t)\bigr)\quad (t\in {\bf C}^\infty_0).\eqno(9.21)$$
\indent (ii) {\sl (Tau function).} The tau function of the right-hand side of (9.21) is defined to be
$$\tau (x;t)=\det \Bigl(I+V(t)R_xW(t)\Bigr)\qquad (t=(t_j)_{j=1}^\infty )\eqno(9.22)$$
\noindent where $R_x=\int_x^\infty e^{-vA_2}BCe^{-vA_1}dv.$\par
\indent (iii)  {\sl (Spectral shift).} We introduce $[s]=(s^j/j)_{j=1}^\infty$ so that for sufficiently small $\vert s\vert$, $\rho$ extends to  
$$\rho ([s]): \bigl(-A_1, -A_2, B, C\bigr)\mapsto  \bigl(-A_1, -A_2, (I-sA_2)^{-1} B, C(I+sA_1)\bigr),\eqno(9.23)$$
\noindent then choose $\zeta =1/s$ with $\zeta\in {\bf C}\setminus {\hbox{Spec}}(A_2)$ so  that the spectral shift is 
$$ \bigl(-A_1, -A_2, B, C\bigr)\mapsto  \bigl(-A_1, -A_2, (I-A_2/\zeta )^{-1} B, C(I+A_1/\zeta )\bigr).\eqno(9.24)$$
\noindent Then we define the Baker--Akhiezer function by
$$\psi_\zeta (x;t)=\exp\Bigl({x\zeta +\sum_{j=1}^\infty \zeta^j t_j}\Bigr) {{\tau (x;t+[1/\zeta ])}\over{\tau (x;t)}}\qquad (t=(t_j)_{j=1}^\infty ).\eqno(9.25)$$
This matches with the definition used in Proposition 9.1 when we choose
$$(t_1, t_2, t_3, t_4,\dots )=(t\lambda /\alpha , y/\beta, t/\alpha, 0, \dots ).\eqno(9.26)$$
\indent (iv) {\sl (Sato's integral).} The Sato integral is
$$\int_{\vert \zeta \vert =r} \tau \bigl(x; t+[1/\zeta ]\bigr)\tau \bigl(x; t'-[1/\zeta ]\bigr)\exp\Bigl({\sum_{j=1}^\infty (t_j-t_j')\zeta^j}\Bigr)\, d\zeta\qquad (t=(t_j)_{j=1}^\infty ,t'=(t'_j)_{j=1}^\infty \in {\bf C}^\infty_0) .\eqno(9.27)$$
\vskip.05in

\noindent {\bf Theorem 9.5.} {\sl Suppose that $A_1, A_2\in {\cal L}(H)$ be as in Theorem 2.2. Then Sato's integral vanishes identically for all $r>\max\{ \Vert A_1\Vert , \Vert A_2\Vert\}.$} 
\vskip.05in
\noindent {\bf Proof} We consider the integral
$$S(t,y)=\int_{\vert \zeta \vert =r} \exp\Bigl({-2\sum_{j=1}^\infty \zeta^jy_j}\Bigr) {{\tau (x, t+y+[1/\zeta ])\tau (x,t-y-[1/\zeta ] )}\over{\tau (x,t+y)\tau (x,t-y)}}d\zeta\eqno(9.28)$$
\noindent which as in (9.8) we can write as
$$\eqalignno{{}\int_{\vert \zeta \vert =r}&\Bigl( 1+Ce^{-xA_2}W(t+y)(I+V(t+y)R_xW(t+y))^{-1}V(t+y)e^{-xA_2} (\zeta I-A_2)^{-1}B\Bigr)\cr
&\times  \Bigl(1-C(\zeta I+A_1)^{-1}e^{-xA_1}W(t-y)(I+V(t-y)R_xW(t-y))^{-1}V(t-y)e^{-xA_2}B\Bigr) \cr
&\times \exp\Bigl( -2\sum_{j=1}^\infty \zeta^jy_j\Bigr)\, d\zeta &(9.29)\cr}$$
\noindent which we split as a sum of four terms: first we have 
$$\int_{\vert \zeta \vert =r} \exp\Bigl( -2\sum_{j=1}^\infty \zeta^jy_j\Bigr)d\zeta =0,\eqno(9.30)$$
\noindent by Cauchy's theorem; the second is  
$$\eqalignno{\int_{\vert \zeta \vert =r}& \exp\Bigl(-2\sum_{j=1}^\infty \zeta^jy_j\Bigr)Ce^{-xA_2}W(t+y)(I+V(t+y)R_xW(t+y))^{-1}V(t+y)e^{-xA_2} (\zeta I-A_2)^{-1}Bd\zeta \cr
&=2\pi i Ce^{-xA_1}W(t+y))(I+V(t+y)R_xW(t+y))^{-1}V(t+y)e^{-xA_2}\exp\Bigl(-2\sum_{j=1}^\infty y_jA_2^j\Bigr)B\cr
&=2\pi i Ce^{-xA_1}W(t+y))(I+V(t+y)R_xW(t+y))^{-1}V(t-y)e^{-xA_2}B;&(9.31)\cr}$$
\noindent by the residue theorem;  the third is
$$\eqalignno{-\int_{\vert \zeta \vert =r}& \exp\Bigl(-2\sum_{j=1}^\infty \zeta^jy_j\Bigr)C(\zeta I+A_1)^{-1}e^{-xA_1}W(t-y)\cr
&\times (I+V(t-y)R_xW(t-y))^{-1}V(t-y)e^{-xA_2}B\, d\zeta \cr
 &=-Ce^{-xA_1}\exp\Bigl(-2\sum_{j=1}^\infty (-A_1)^jy_j\Bigr)W(t-y)(I+V(t-y)R_xW(t-y))^{-1}V(t-y)e^{-xA_2}B\cr
&=-Ce^{-xA_1}W(t+y)(I+V(t-y)R_xW(t-y))^{-1}V(t-y)e^{-xA_2}B&(9.32)\cr}$$
\noindent likewise; and finally
$$\eqalignno{-&\int_{\vert \zeta \vert =r} \exp\Bigl(-2\sum_{j=1}^\infty \zeta^jy_j\Bigr)Ce^{-xA_2}W(t+y)(I+V(t+y)R_xU(t+y))^{-1}V(t+y)e^{-xA_2} (\zeta I-A_2)^{-1}\cr &\times BC(\zeta I+A_1)^{-1}e^{-xA_1}W(t-y)(I+V(t-y)R_xW(t-y))^{-1}V(t-y)e^{-xA_2}B\, d\zeta&(9.33)\cr}$$
\noindent which involves
$$J(y)=\int_{\vert \zeta \vert =r}\exp\Bigl(-2\sum_{j=1}^\infty \zeta^jy_j\Bigr)(\zeta I-A_2)^{-1}BC(\zeta I+A_1)^{-1}d\zeta.\eqno(9.34)$$
\noindent This integral resembles (4.40), and likewise gives a solution to a type of Lyapunov equation. Now
$$\eqalignno{ -{{1}\over{2}}{{\partial J}\over{\partial y_1}}+JA_1&= \int_{\vert \zeta \vert =r}\exp\Bigl(-2\sum_{j=1}^\infty \zeta^jy_j\Bigr) (\zeta I-A_2)^{-1}BCd\zeta\cr
&=2\pi i V(-2y)BC&(9.35)\cr}$$
\noindent and 
$$\eqalignno{ -{{1}\over{2}}{{\partial J}\over{\partial y_1}}-A_2J&= \int_{\vert \zeta \vert =r}\exp\Bigl(-2\sum_{j=1}^\infty \zeta^jy_j\Bigr) (\zeta I-A_2)^{-1}BCd\zeta\cr
&=2\pi i BCW(2y)&(9.36)\cr}$$
\noindent so by subtracting and applying the residue theorem, we have
$$A_2J+JA_1=2\pi i (V(-2y)BC-BCW(2y)).\eqno(9.37)$$
\noindent Then we introduce $J_0=2\pi i ( V(-2y)R_0-R_0W(2y)),$ which satisfies
$$\eqalignno{A_2J_0+J_0A_1&=2\pi i V(-2y)(A_2R_0+R_0A_1)-2\pi i(A_2R_0+R_0A_1)W(2y)\cr
& =2\pi i (V(-2y)BC-BCW(2y))
=A_2J+JA_1,&(9.38)\cr}$$
\noindent by Lyapunov's equation. By the uniqueness of solution of this equation, we deduce that
$$J(y)=J_0= 2\pi i \bigl( V(-2y)R_0-R_0W(2y)\bigr).\eqno(9.39)$$
\noindent Then, combining the terms (9.30), (9.31), (9.32) and (9.33) via (9.39), we have 
$$\eqalignno{S(y,t)&=2\pi iCW(t+y)e^{-xA_1}\Bigl(  (I+V(t+y)R_xW(t+y))^{-1}  - (I+V(t-y)R_xW(t-y))^{-1}\cr
&\qquad- (I+V(t+y)R_xW(t+y))^{-1}V(t+y)e^{-xA_1}Je^{-xA_2}W(t-y)\cr
&\qquad \times (I+V(t-y)R_xW(t-y))^{-1}\Bigr) V(t-y)e^{-xA_2}B\cr
&= 2\pi iCW(t+y)e^{-xA_1}(I+V(t+y)R_xW(t+y))^{-1}\Bigl( (I+V(t-y)R_xW(t-y))\cr
&\qquad -(I+V(t-y)R_xW(t-y))+( V(t+y)R_xW(t+y)-V(t-y)R_xW(t-y))
\Bigr)\cr
&\qquad \times (I+V(t-y)R_xW(t-y))^{-1} V(t-y)e^{-xA_2}B\cr
&=0,&(9.40)}$$
\noindent as required.\par
\rightline{$\square$}\par
\vskip.05in

\noindent {\bf Corollary 9.6.}{\sl  (i) For $s_0, s_1, s_2, s_3\in {\bf C},$ let $\sigma_{jk}=(s_j-s_k)\tau (x;t+[s_j]+[s_k])$. Then Fay's identity holds
$$\sigma_{0,1}\sigma_{2,3}-\sigma_{0,2}\sigma_{1,3}+\sigma_{0,3}\sigma_{1,2}=0;\eqno(9.41)$$
\indent (ii) for the Wronskian with derivatives in the $x$-variable,  the differential form of Fay's identity holds}
$${\hbox{Wr}}\bigl(\tau \bigl(x;t+[s_1]\bigr), \tau \bigl(x;t+[s_2]\bigr)\bigr)={{s_1-s_2}\over{s_1s_2}}\det\left[\matrix{\tau (x;t)&\tau (x;t+[s_2])\cr \tau (x;t+[s_1])& \tau (x;t+[s_1]+[s_2])\cr}\right];
\eqno(9.42)$$
\indent {\sl (iii) the second-order differential identity holds} 
$${{\partial^2 }\over{\partial \zeta \partial x}}\log \tau \bigl(x;t+[1/\zeta ]\bigr)=1-{{\tau (x;t)\tau (x;t+2[1/\zeta ])}\over{\tau (x; t+[1/\zeta ])^2}}.\eqno(9.43)$$  
\noindent 
\noindent {\bf Proof.} (i) See [66] and [2]. We have written the result in the style of a Pl\"ucker relation.  \par
\indent (ii)  See [66] and [2]. This spectral addition rule has a similar style to the Toda equation (3.18).\par
\indent (iii) We divide  (ii) by $\tau (x;t+[s_1])\tau (x;t+[s_2])$  so as to obtain ${{\partial}\over{\partial x}}\log \tau (x;t+[s_2])/\tau (x;t+[s_1])$ on the left-hand side; then we differentiate 
 with respect to $s_2$, thus obtaining 
$$\Bigl({{\partial^2 }\over{\partial s_2 \partial x}}\Bigr)_{s_2=s_1}\log \tau (x;t+[s_2])=-{{1}\over{s_1^2}}\Bigl( 1-{{\tau (x;t)\tau (x;t+2[s_1])}\over{\tau (x; t+[s_1])^2}}\Bigr);\eqno(9.44)$$
\noindent then we change variables to $s_1=1/\zeta. $ This  resembles the proof in [57, II  3.124].\par

\rightline{$\square$}\par

\vskip.05in
\noindent {\bf Acknowledgements}\par
\noindent The work of the second-named author was supported by an EPSRC Research studentship. The authors are grateful to Prof H.P. McKean for helpful comments which motivated this work.\par

\vskip.05in

\noindent {\bf References}\par
\vskip.05in
\noindent [1] M. Adler, Integrable systems, random matrices and random processes, pp. 131-225, J. Harned edr. {\sl Random Matrices,
Random Processes and Integrable Systems} CRM (Springer, 2011).\par
\noindent [2] M. Adler and P. van Moerbeke, Birkhoff strata, B\"acklund
transformations, and regularization of isospectral operators, {\sl
Adv. Math.} {\bf 108} (1994), 140--204.\par
\noindent [3] T. Akaza, Poincar\'e theta series and singular sets of Schottky groups, {\sl Nagoya Math. J.} {\bf 24} (1964), 43--65.\par


\noindent [4] T. Aktosun, F. Demontis and C. van der Mee, Exact solutions
to the focusing nonlinear Schr\"odinger equation, {\sl Inverse
Problems} {\bf 23} (2007), 2171--2195.\par
\noindent [5] Y.A. Antipov and D.G. Crowdy, Riemann--Hilbert problem for automorphic functions and the Schottky--Klein prime function, {\sl Complex Anal. Oper. Theory } {\bf 1}  (2007), 317-334.\par

\noindent [6] M.F. Atiyah and I.G. Macdonald, {\sl 
Introduction to commutative
algebra}, (Addison Wesley, Reading, Massachusetts, 1969).\par   
\noindent [7] H.F. Baker {\sl Abelian functions: Abel's theorem and the
allied theory of theta functions}, Cambridge University Press,
edition of 1995.\par
\noindent [8] R. Bhatia and P. Rosenthal, How and why to solve the
operator equation $AX-XB=Y$, {\sl Bull. London Math. Soc.} {\bf 29}
(1997), 1--21.\par
\noindent [9] G. Blower, Linear systems and determinantal random point
fields, {\sl J. Math. Anal. Appl.} {\bf 335} (2009), 311--334.\par
\noindent [10] G. Blower, On linear systems and $\tau$ functions
associated with Lam\'e's equation and Painlev\'e's equation VI, {\sl J.
Math. Anal. Appl.} {\bf 376} (2011), 294--316.\par
\noindent [11] G. Blower, On tau functions for orthogonal polynomials
and matrix models, {\sl J. Phys. A} {\bf 44} (2011),
285202.\par
\noindent [12] A. Borodin and P. Deift, Fredholm determinants,
Jimbo--Miwa--Ueno $\tau$-functions and representation theory, {\sl
Comm. Pure Appl. Math.} {\bf 55} (2002), 1160--2030.\par
\noindent  [13] C. Brett, {\sl Sampling from the spectrum of Hill's equation} (Ph.D. Thesis, Lancaster, 2015).\par
\noindent [14] Y.V. Brezhnev, What does integrability of finite-gap or
soliton potentials mean? {\sl Philos. Trans. R. Soc. Lond. Ser. A
Math. Phys. Eng. Sci.} {\bf 366} (2008), 923--945.\par
\noindent [15] W. Burnside, On a class of automorphic functions,  {\sl Proc. London Math. Soc.} {\bf 23} (1891) 49-88.\par
\noindent [16] W. Burnside, Further note on automorphic functions,  {\sl Proc. London Math. Soc.} {\bf 23} (1891) 281-295.\par


\noindent [17] D.G. Crowdy, C.C. Green, E.H. Kropf and M.M.S. Nasser, the Schottky--Klein prime function: a theoretical and computational tool for applications, {\sl IMA J. Appl. Math.} {\bf 81} (2016), 589--628.\par 
\noindent [18] D. Crowdy and J. Marshall, Analytical formulae for the Kirchhoff--Routh path function in multiply connected domains, {\sl Proc. R. Soc. A} {\bf 461} (2005), 2477--2501.\par 
\noindent [19] E.B. Davies, {\sl Linear Operators and their Spectra}, (Cambridge University Press, 2007).\par
\noindent [20] F.J. Dyson, Fredholm determinants and inverse scattering problems, {\sl Commun. Math. Phys.} {\bf 47} (1976), 171--183.\par 
\noindent [21] K.-J. Engel and R. Nagel, {\sl One parameter semigroups for
linear evolution equations}, Springer--Verlag, New York, 2000.\par
\noindent [22] N. Ercolani and H.P. McKean, Geometry of KdV. IV: Abel sums,
Jacobi variety and theta function in the scattering case, {\sl Invent.
Math.} {\bf 99} (1990), 483--544.\par
\noindent [23] J.D. Fay, Theta functions on Riemann surfaces, Lecture
Notes in Mathematics 352, Springer, Berlin 1973.\par
\noindent [24] N.C. Freeman  and J.J.C. Nimmo, Soliton solutions of the Korteweg-- de Vries and Kadomtsev -- Petviashvili equations: the Wronskian 
technique, {\sl Proc. Roy. Soc. London Ser. A} {\bf 389} (1983), 319--329.\par 
\noindent [25] J. Garnett and E. Trubowitz, Gaps and bands of one dimensional periodic Schr\"odinger operators II, {\sl Comment. Math. Helv.} {\bf 62} (1987), 18-37.\par

\noindent [26] I.M. Gelfand and L.A. Dikii, Integrable nonlinear
equations and the Liouville theorem, {\sl Funct. Anal. Appl.} {\bf 13}
(1979), 6-15.\par 
\noindent [27] I.M. Gelfand and B.M. Levitan, On the determination of a differential equation from its spectral function, {\sl Izvestiya Akad. Nauk SSSR Ser. Mat.} {\bf 15} (1951), 309--360).\par
\noindent [28] F. Gesztesy, Integrable systems in the infinite genus limit, {\sl Differential and Integral Equations} {\bf 14} (2001), 671--700.\par
\noindent [29] F. Gesztesy and H. Holden, {\sl Soliton equations and
their algebro-geometric solutions Volume I: $(1+1)$-dimensional
continuous models}, (Cambridge University Press,
2003).\par
\noindent [30] F. Gesztesy and B. Simon, The xi function, {\sl Acta Math.} 176 (1996), 49-71.\par
\noindent [31] F. Gesztesy and R. Weikard, Picard potentials and Hill's
equation on a torus, {\sl Acta Math.} {\bf 176} (1996), 73--107.\par
\noindent [31'] W. Hayman, {\sl Meromorphic functions}, (Oxford University Press, 1964).\par 
\noindent [32] E. Hille, {\sl Lectures on Ordinary Differential Equations}, (Addison-Wesley, 1968).\par
\noindent [33] M. Jimbo and T. Miwa, Monodromy preserving
deformation of linear ordinary differential equations with rational
coefficients II {\sl Phys. D.} {\bf 2} (1981), 407--448.\par
\noindent [34] R. Johnson and J. Moser, The rotation number for almost periodic potentials, {\sl Comm. Math. Phys.} {\bf 84}, (1982), 403-438.\par
\noindent [35] V.G. Kac, {\sl Infinite dimensional Lie algebras}, (Cambridge University
Press, 1985).\par
\noindent [36] K. Kajiwara, T. Masuda, M. Noumi, Y. Ohta, and Y. Yamada, Determinant formulas for the Toda and discrete Toda equations,
{\sl Funkcial. Ekvac.} {\bf 44} (2001), 291--307.\par
\noindent [37] S. Kamvissis, Inverse scattering as an infinite
period limit, {\sl C. R. Acad. Sci. Paris S\'er I} {\bf 325} (1997)
969--974.\par  
\noindent [38] H. Kato and S. Saito, Generalization of the Boson--Fermion equivalence and Fay's addition theorem, {\sl Lett. Math. Phys.} {\bf 18} (1989), 177--183.\par
\noindent [39] I.M. Krichever, The integration of nonlinear equations
by the methods of algebraic geometry, {\sl Funct. Anal.
Appl.} {\bf 11} (1978), 12--26.\par
\noindent [40] I.M. Krichever, Elliptic solutions of the
Kadomcev--Petviashvili equations, and integrable systems of particles, {\sl Funct Anal. Appl.} 
{\bf 14} (1980), 282--290.\par
\noindent [41] I.M. Krichever, Vector bundles and Lax equations on algebraic curves, {\sl Comm. Math. Phys.} {\bf 229} (2002), 229--269.\par
\noindent [42] I.M. Krichever and S.P. Novikov, Holomorphic bundles
and nonlinear equations, {\sl Physica D} {\bf 3} (1981), 267--293.\par
\noindent [43] S. Lang, {\sl Introduction to algebraic and abelian
functions}, Second edition, Springer, New York, 1982.\par 
\noindent [44] R. S. Maier, Lam\'e polynomials, hyperelliptic reductions
and Lam\'e band structure, {\sl Philos. Trans. R. Soc. Lond. Ser. A
Math. Phys. Eng. Sci.} {\bf 366} (2008), 1115--1153.\par    
\noindent [45] S.V. Manakov, The inverse scattering transform for the time dependent Schr\"odinger equation and the Kadomtsev--Petviashvili equation, {\sl Physica} {\bf 3} D (1981), 420--427.\par
\noindent [46] V.B. Matveev, Darboux transformation and explicit solutions of the Kadomtcev-\par
\noindent -Petviaschvily equation, depending upon functional parameters, {\sl Lett. Math. Phys} {\bf 3} (1979), 213--216.\par

\noindent [47] H.P. McKean, Geometry of KDV. III: determinants and unimodular isospectral flows, 
{\sl Commun. Pure Appl. Math.} {\bf 45} (1992),  389-415.\par 
\noindent [48] H.P. McKean, Fredholm determinants, {\sl Cent. Eur. J. Math.} {\bf
9} (2011), 205--243.\par

\noindent [49] H. McKean and V. Moll, {\sl Elliptic curves: function
theory, geometry, arithmetic,} (Cambridge University Press, 1999).\par
\noindent [50] H.P. McKean and P. van Moerbeke, {The spectrum of Hill's
equation}, {\sl Invent. Math.} {\bf 30} (1975), 217--274.\par
\noindent [51] H.P. McKean and E. Trubowitz, Hill's operator and hyperelliptic function theory in the presence of
infinitely many branch points, {\sl Comm. Pure Appl. Nath.} {\bf 29} (1976), 143--226.\par
\noindent [52] H.P. McKean and E. Trubowitz, Hill's surfaces and their theta functions, {\sl Bull. Amer. Math. Soc.} {\bf 84} (1978), 1042--1085.\par 
\noindent [53] A.V. Megretskii, V.V. Peller and S.R. Treil, The inverse
spectral problem for self-adjoint Hankel operators, {\sl Acta Math.}
{\bf 174} (1995), 241--309.\par  
\noindent [54] M. L. Mehta, {\sl Random Matrices}, Second Edition, (Academic press, 1991).\par
\noindent [55] T. Miwa, M. Jimbo, and E. Date, {\sl Solitons:
Differential equations, symmetries, and infinite dimensional algebras,}
Cambridge University Press, 2000.\par
\noindent [56] M. Mulase, Cohomological structure in soliton
equations and Jacobian varieties, {\sl J. Differential Geom.}
{\bf 19} (1984), 403--430.\par
\noindent [57] D. Mumford, {\sl Tata Lectures on Theta II}, Birkh\"auser, Boston, 1984.\par
\noindent [58] N.K. Nikolski, {\sl Operators, Functions and Systems: An Easy Reading.
Volume 1}, (American Mathematical Society, 2002).\par

\noindent [59] S. Novikov, S.V. Manakov, L.P. Pitaevskii and V.E. Zakharov,
 {\sl Theory of Solitons, the Inverse
Scattering Method}, Consultants Bureau, New York and London, 1984).\par

\noindent [60] V.V. Peller, {\sl Hankel operators and their
applications}, (Springer--Verlag, New York, 2003).\par
\noindent [61] C. P\"oppe, The Fredholm determinant method for the KdV
equations, {\sl Physica} {\bf 13} D (1984), 137--160.\par
\noindent [62] C. P\"oppe and D.H. Sattinger, Fredholm determinants and
the $\tau$ function for the Kadomtsev --Petviashvili hierarchy, {\sl
Publ. Re. Inst. Math. Sci.} {\bf 24} (1988), 505--538.\par  
\noindent [63] M. van der Put and M.F. Singer, {\sl Galois theory of
linear differential equations}, (Springer, Berlin, 2003).\par
\noindent [64] G. Segal and G. Wilson,  Loop groups and equations of
KdV type, {\sl Inst. Hautes \'Etudes Sci. Publ. Math.} {\bf 61}
(1985), 5--65.\par  
\noindent [65] I.R. Shafarevich, {\sl Basic algebraic geometry},
(Springer,  Berlin, 1977).\par
\noindent [66] T. Shiota, Characterization of Jacobian varieties in
terms of soliton equations, {\sl Invent. Math.} {\bf 83} (1986),
333--382.\par
\noindent [67] B. Simon, {\sl Trace ideals and their applications}, second edition, (American Mathematical
Society, 2005).\par 
\noindent [68] C.A. Tracy and H. Widom, Fredholm determinants, differential equations and matrix models, {\sl Comm. Math. Phys.} 
{\bf 163} (1994), 33--72.\par
\noindent [69] M. Trubowitz, The inverse problem for periodic potentials, {\sl Commun. Pure Appl. Math.} {\bf 30} (1977), 321--337.\par

\noindent [70] G.L. Vasconcelos, J.S. Marshall and D.G. Crowdy, Secondary Schottky--Klein prime functions associated with multiply connected planar domains, {\sl Proc. Roy. Soc. London } A {\bf 471}, 2014, \par 
\noindent [71] V.E. Zakharov, `The inverse scattering method', pp 243--285 in {\sl Solitons}, R.K. Bullough and P.J. Caudrey, (Springer-Verlag, 1980).\par
\noindent [72] V.E. Zakharov and A. Shabat, A plan for integrating the
nonlinear equations of mathematical physics by the method of the
inverse scattering problem. I., {\sl Funct. Anal. Appl.} {\bf 8} (1975),
226--235.\par 
\noindent [73] T. Zhang and S. Venakides, Periodic limit of inverse
scattering, {\sl Comm. Pure Appl. Math.} {\bf 46} (1993), 819--865.\par
\vfill
\eject
\end

\indent As in de Branges's theory of entire functions, we write $f^*(z)=\overline {f(\bar z)}$ and consider the space of entire functions $E$ such that $E$ has no zeros in the closed  upper half plane, $\vert E(z)\vert \geq \vert E(\bar z)\vert$ for $\Im z\geq 0$; and introduce the Hilbert space  of entire functions ${\cal H}(E)$ consisting of entire functions $f$ such that $f/E$ and $f^*/E$ belong to $H^2$, with the inner product
$$\langle f,g\rangle_{{\cal H}(E)}=\int_{-\infty}^\infty {{f(x)\overline{g(x)}}\over{\vert E(x)\vert^2 }}dx.\eqno(10.2)$$
\noindent Then ${\cal H}(E)$ has the fundamental properties:\par
\indent (H1) $f(z)\mapsto f(z)(z-w)/(z-\bar w)$ is an isometry on ${\cal H}(E)$ for all $w\in {\bf C}$;\par
\indent (H2) for all $w\in {\bf C}\setminus {\bf R}$, the map $f(z)\mapsto f(w)$ is a continuous linear 
functional on ${\cal H}(E)$;\par 
\indent (H3) $f\mapsto f^*$ defines a real linear isometry on ${\cal H}(E)$.\par
\vskip.05in
\indent We also introduce the entire functions $F(z)=(E(z)+E^*(z))/2$ and $G(z)=i(E(z)-E^*(z))/2$, and have the reproducing kernel
$$h_w(z)={{G(z)\bar F(w)-F(z)\bar G(w)}\over{\pi (z-\bar w)}},\eqno(10.3)$$
\noindent so that $f(w)=\langle f, h_w\rangle_{{\cal H}(E)}.$\par
\indent In particular, we can take $E(z)=e^{-i\gamma z}$ for $\gamma >0$ and obtain the Paley--Wiener space
 $PW_\gamma $ which consists of entire functions $g$ of exponential type such that
 $\int_{-\infty}^\infty\vert g(x)\vert^2\, dx$ converges and 
$$\lim\sup_{y\rightarrow\pm\infty}{{\vert g(iy)\vert}\over{\vert y\vert}}\leq \gamma ;\eqno(10.4)$$
\noindent we write $PW$ for $PW_\pi$.\par
\vskip.05in
\noindent {\bf Definition} A phase function is a continuous and increasing real function that takes $0$ as a value.
 We say that phase functions $\nu_1$ and $\nu_2$ are equivalent if $\nu_1(x)=\nu_2(x)$ whenever $\nu_1(x)=k\pi$ or 
$\nu_2(x)=k\pi$ for some $k\in {\bf Z}$. For a Hilbert space of entire functions, the phase function is 
$E(x)=\vert E(x)\vert e^{-i\nu (x)/2}$.\par
\vskip.05in  
\indent We introduce a phase function for Hill's equation as follows. Let 
$$V(\lambda )=\alpha +\int_{\lambda_0}^\lambda {{\varepsilon (s)\Delta'(s) ds}\over{\sqrt{4-\Delta (s)^2}}}\eqno(10.5)$$
\noindent where $\varepsilon (s)=\pm 1$ is chosen to be constant on all open intervals of stability and instability such 
that $\nu (\lambda )=\Re V(\lambda )$ is a phase function. Then $V$ itself is known as a spike function since the graph 
$\{ V(\lambda ): \lambda >0\}$ consists of abutting line segments $[\beta +j\pi , \beta +(j+1)\pi ]$ 
parametrized by the intervals of stability $[\lambda_{2j}, \lambda_{2j+1}]$, and vertical segments above the 
$\beta +j\pi$ which are parametrized by the intervals of instability $[\lambda_{2j+1}, \lambda_{2j+2}].$  \par
\indent We show that the discriminant of Hill's equation gives an equivalence class of phase functions which 
determines a Hilbert space of entire functions essentially uniquely.\par
\vskip.05in

\noindent {\bf Lemma 10.1} {\sl
Suppose that $\sup\{ \vert \sqrt{\lambda_{2j}}-\pi j\vert \quad :\quad j=0, 1, \dots \} <\pi /4$. Then\par
\indent (i)  there exists a Hilbert space ${\cal H}={\cal H}(E)$ of entire functions with norm 
$$\Vert f\Vert_{\cal H}=\bigl(\sum_{j=0}^\infty \vert
f(\sqrt{\lambda_{2j}})\vert^2\bigr)^{1/2};\eqno(10.6)$$
\indent (ii) ${\cal H}$ has a phase function equivalent to $\nu$, and when $E$ has no real zeros, 
the phase function $\nu$ determines ${\cal H}$ up to a canonical isometric isomorphism.}\par
\vskip.05in

\noindent {\bf Proof.} (i) Let $t_n=\sqrt{\lambda_{2n-1}}$ for $n=1, 2, \dots $; $t_0=0$ and $t_{-n}=-t_n$ for $n=1, 2, \dots .$  By Kadec's $1/4$ theorem, the exponentials $(e^{it_nx})_{j=0}^\infty$ form a Riesz basis for $L^2(-1,1)$, so there exist constants $\kappa_3, \kappa_4>0$ such that 
$$\kappa_3\sum_{j=-\infty }^\infty \vert a_j\vert^2\leq \int_{-1}^1 \Bigl\vert\sum_{j=-\infty }^\infty a_je^{it_nx
}\Bigr\vert^2 dx\leq \kappa_4\sum_{j=-\infty}^\infty \vert a_j\vert^2\qquad ((a_j)_{j=-\infty }^\infty\in \ell^2).\eqno(10.7)$$
\indent Let ${\hbox{sinc}}(t)=(\sin t)/t$ be the unnormalized ${\hbox{sinc}}$ function and 
consider 
$$g(t)=\sum_{n=-\infty}^\infty
a_n {\hbox{sinc}}(\pi t-t_n);$$
\noindent then $g\in PW$. Moreover, for all $f\in PW$, we have
$$\int_{-\infty}^\infty f(t)\bar g(t)\, dt=\sum_{n=-\infty}^\infty a_nf(t_n)$$
by Plancherel's formula, so $(t_n)_{n=-\infty}^\infty$ is a sampling sequence for $PW$. Furthermore, there exist constants
$\kappa_5, \kappa_6>0$ such that 
$$\kappa_5\sum_{n=-\infty}^\infty \vert f(t_n)\vert^2\leq 
\int_{-\infty}^\infty \vert f(t)\vert^2\, dt\leq \kappa_6\sum_{n=-\infty}^\infty \vert f(t_n)\vert^2 ;$$
\noindent Hence we can introduce an equivalent Hilbert space norm on $PW$ by 
$\Vert f\Vert^2=\sum_{n=-\infty}^\infty \vert f(t_n)\vert^2$. Indeed, by de Branges's fundamental existence theorem, there exists an entire function $E$ as above such that 
${\cal H}(E)=PW$ as sets, and 
$$\sum_{n=-\infty}^\infty\vert f(t_n)\vert^2=\int_{-\infty}^\infty {{\vert f(t)\vert^2}\over{\vert E(t)\vert^2}}dt.\eqno(10.8)$$
\indent (ii) Furthermore, the phase function for $E$ satisfies $\nu (t_k)=\pi k$ for all 
$k\in {\bf Z}$, and $E$ has no real zeros. Evidently we can 
choose the constants $\alpha$ and $\beta$ so that $\nu$ is equivalent to the phase function 
$$\alpha+ {\hbox{sgn}}(t)\int_{{\lambda_0}}^{t^2} {{\varepsilon (s)\Delta'(s)ds}\over{\sqrt{4-\Delta (s)^2}}}.\eqno(10.9)$$
\noindent Clearly the points $(t_k)_{k=-\infty}^\infty$ determine the phase function up to equivalence. Suppose that 
${\cal H}(E_0)$ is another Hilbert space of entire functions that has phase function $\nu_0$ equivalent to $\nu$. Then by Theorem 24 of [dB], there exists an entire function $S(z)$ such that all zeros of $S$ are real and $S(x)$ is real for all $x\in {\bf R}$, and $f(z)\mapsto S(z)f(z)$ gives an isometric isomorphism of ${\cal H}(E_0)\rightarrow {\cal H}(E)$.
 Thus the  phase function determines ${\cal H}(E)$ up to a specific linear isometric equivalence. Indeed, when $E$ has no real zeros and $\vert E(x-iy)\vert <\vert E(x+iy)\vert$ for all $y>0$, there exists $p>0$ such that  phase function may be chosen to satisfy [dB , p63]
$${{\partial}\over{\partial y}}\log \vert E(x+iy)\vert =py+\int_{-\infty}^\infty {{d\varphi (t)}\over{(t-x)^2+y^2}}\qquad (y>0, x\in {\bf R}).$$
\rightline{$\square$}

\indent The elements of $PW$ satisfy$\int_{-\infty}^\infty \vert f(t)\vert^2dt<\infty$ and $\vert f(t+is)\vert \leq Me^{\pi\vert s\vert}$ for some $M$ and all $s,t\in {\bf R}$.\par
\indent in the context of Lemma 10.1, we can suppose that the phase function $\nu$ for $E(z)$ and the phase
function $\tilde \nu$ for $E^*(-z)$ are equivalent, and that $E(z)$ has no zeros on the real axis. Then
there exists an entire function $S(z)$ which is real on the real axis, with all of its zeros on the real
axis, such that 
$$\int_{-\infty}^\infty {{\vert S(t)f(t)\vert^2}\over{\vert E(-t)\vert^2}}dt=
\int_{-\infty}^\infty {{\vert f(t)\vert^2}\over{\vert E(t)\vert^2}}dt.$$
Let $f$ be an even function in ${\cal H}(E)$. Then there exists an entire function $g$ such that
$f(z)=g(z^2)$, so $\vert g(z)\vert \leq Me^{\kappa\sqrt\vert g\vert}$ and there exist universal
constants of eqivalence such that
$$\sum_{n=0}^\infty \vert g(\lambda_n )\vert^2=\int_{-\infty}^\infty {{\vert g(t^2)\vert^2}\over{\vert
E(t)\vert^2}}dt\sim \int_0^\infty {{\vert g(u)\vert^2}\over{\sqrt{u}}}du,$$
and \par

\vskip.05in
\noindent {\bf Definition} Let $\Sigma =(-A,B,C;E)$ be a periodic linear system such that
$u\in C^2({\bf R}; {\bf C})$.\par
\indent (i) Suppose that $u$ is real valued. Then we say that $\Sigma$ has the phase function $\nu$, and that $\Sigma$ has 
Hilbert space of entire functions ${\cal H}$, which is associated with the equivalence class of $\nu$ as in Lemma 10.1.\par
\indent (ii) (Picard system) Say that a periodic linear system 
$\Sigma$ is a Picard system if $-\psi''(x)+u(x)\psi (x)=\lambda\psi (x)$
has a meromorphic general solution $\psi $ for all but finitely many
$\lambda\in {\bf C}$. See [28].\par
\indent (iii) When $u$ is meromorphic, then we say that $\Sigma$ has
divisor $D$, the divisor of $u$. In particular, this occurs in case (ii).\par

\vskip.05in

\vskip.05in
\noindent {\bf Definition} (Divisors)
 Given a complex analytic
manifold, a Weyl divisor is a finite sum with integral coefficients
of irreducible subvarieties of codimension one. For a (possibly non compact) Riemann surface 
${\cal E}$ with meromorphic functions ${\bf M}({\cal E})$, the divisors are those functions
 $\nu:{\cal E}\rightarrow {\bf Z}$ such that the restriction $\nu\vert_K$ has finite support for all 
compact subsets $K$ of ${\cal E}$ and the divisors that have finite support on ${\cal E}$ are said to be finite. 
The degree of a finite divisor $\nu$ is $\sum_{x\in {\cal E}}\nu (x)$.  In particular, for each non-zero 
$g\in {\bf M}({\cal E})$, let $\nu (g,x)$ be the order of the zero $x$ of $g$, and let 
$(g)=\sum_{x\in {\cal E}} \nu (g,x)\delta_x-\sum_{x\in {\cal E}}\nu (1/g, x)\delta_x$ be the 
principal divisor associated with $g$. Let ${\hbox{Div}}_0({\cal E})$ be the space of finitely supported divisors 
that have degree zero. Each $D\in {\hbox{Div}}({\cal E})$ may be identified with a Radon
measure by associating $\delta_x$ with the point mass at $x$; thus ${\hbox{Div}}({\cal E})_0$
may be embedded into the complex vector space of bounded Radon measures on ${\cal E}$. When
$D(t)\in {\hbox{Div}}({\cal E})_0$ depends upon $t\in {\bf R}$, we can consider evolution
equations involving $D(t)$ as a family of measures. \par
\indent When ${\cal E}$ is compact, all  the principal divisors belong to ${\hbox{Div}}({\cal E})_0$ by Abel's theorem. 
One can identify the Picard group  ${\hbox{Pic}}({\cal E})$ with ${\hbox{Div}}({\cal E})_0$ modulo the principal divisors.\par
\indent Let $u$ be as in Proposition 10.2. Then there are only finitely many simple periodic eigenvalues, and we can replace 
the infinite series in () by finite sums, and the ${\hbox{Pic}}({\cal E})$ can be embedded into a complex Abelian manifold 
of dimension $g$.   For any polynomial function $R$ of degree $2g+2$,  surface the surface ${\cal E}=\{ (X,Y): Y^2=R(X)\}$ 
is hyperelliptic, and the Abel map ${\hbox{Pic}}({\cal E})\rightarrow {\bf C}^g$  is  
$$\sum_{j=1}^g(\delta_{p_j}-\delta_\infty )\mapsto \Bigl(\sum_{j=1}^g\int_\infty^{p_j}{{X^{k-1}dX}\over{\sqrt{R(X)}}}
\Bigr)_{k=1}^g.$$   
Abel calculated the canonical integral of the third kind $\Pi_{\alpha, \beta}^{x,z}$, as in Baker p. 195, then one can
 express any Abelian integral in terms of sums of $\Pi_{\alpha_j, \beta_j}^{x_j,z_j}$ and their derivatives and Abelian
 integrals of the first kind. Hence one can invert the Abel map without introducing the 
Riemann theta function, as in Baker p244.\par 
\indent Mumford shows that one can endow the Jacobi variety with the structure of an algebraic group, without using the 
transcendental analytic functions.
The Mumford shows [48, 3.155] that there exists a hyperelliptic function $\wp$ and a derivation $D_{\bf X}$ on the rational 
functions on ${\bf X}$ such that any rational function $u$ on ${\bf X}$ has the form 
$u=q(\wp ,D_{\bf X}\wp, \dots ,D^g_{\bf X}\wp )$ for some complex polynomial $q$.\par 
\indent In the remainder of this section, we follow through some of these ideas, without 
assuming in general that $u$ is finite gap. We introduce the transcendental hyperelliptic curve
$${\cal E}=\{ (\lambda , \mu ):\mu^2=4-\Delta (\lambda )^2\}$$
\noindent which has infinite genus.\par 
\vskip.05in

\noindent {\bf Proposition 10.3} {\sl Let $u$ be a real, $C^2$ and $\pi$-periodic potential. 
Let ${\cal E}$ be the multiplier curve for the corresponding Hill's equation, and let $\phi_{\bf p}$ be a Floquet solution 
that corresponds to ${\bf p}\in {\cal E}$. Then there exists a family of linear systems ${\Sigma}_{\bf p}$ with scattering
 functions $\phi_{\bf p}$ for ${\bf p}\in {\cal E}$. The linear systems are admissible and depend continuously on ${\bf p}$
 on the complement of the Bloch spectrum, and periodic on the periodic spectrum.}\par
\vskip.05in

\noindent {\bf Proof.} For $\vert \mu\vert \leq 1$, let $W_\mu =\{ f\in L^2_{loc}(0, \infty ): f(x+\pi )=\mu f(x); x>0\}$; 
then for $\vert \mu\vert <1,$ this $W_\mu$ is an eigenspace of $S:L^2(0, \infty )\rightarrow L^2(0, \infty )$, 
$S(f(x)=f(x+\pi )$ while for $\vert \mu\vert =1$, this $W_\mu$ gives approximate eigenvectors for $S$. We introduce 
$$\eqalignno{L_{\bf p}&=\{ f: -f''+uf=\lambda f; f\in W_\mu \}\qquad \vert \mu\vert \leq 1;\cr
&=\{ f: -f''+uf=\lambda f; f(-x)\in W_{1/\mu }\} \qquad \vert\mu\vert >1.&(10.18)\cr}$$
\noindent Loosely speaking, $L_{\bf p}$ gives a line bundle over ${\cal E}$, but we need to take care with 
the Bloch spectrum and the periodic eigenvalues in particular.\par

\indent The gaps in the Bloch spectrum are given by 
$\{ \lambda \in {\bf R}: \Delta (\lambda )^2>4\}$, which are known as intervals of instability since all nontrivial solutions of Hill's equation are unbounded 
on ${\bf R}$.  
\vskip.05in
\noindent {\bf Lemma 10.4} {\sl (i)  For all $\lambda\in {\bf R}$ such that $\Delta (\lambda )^2>4$,  there exists $\mu\in {\bf
C}$ such that $\vert \mu\vert <1$ and Hill's equation $-\psi''(x)+u(x)\psi (x)=\lambda\psi (x)$
\noindent has a nontrivial solution such that $\psi (x+\pi )=\mu \psi (x)$.\par
\indent (ii) There exists a linear system $(-A,B,C)_{\bf p}$ with bounded $A$, state space $W_\mu$ and input and output 
spaces ${\bf C}$, such that 
$\psi$ is the scattering function of $(-A,B,C)_{\bf p}$, and $\psi_{(\pi )}$ is the scattering function of
 $(-A,\mu B, \mu C)_{\bf p}$;\par
\indent  (iii) $\Gamma_\psi$ satisfies Corollary 5.6. }\par

\vskip.05in
\noindent {\bf Proof.} (i) This is immediate from Floquet's theorem.\par
\indent (ii) Since $\vert \mu\vert <1$, $x\vert\psi (x)\vert^2$ is integrable over $(0, \infty)$, so  the Hankel
 operator $\Gamma_\psi$ is  Hilbert--Schmidt; hence $\psi$ can be realised as the scattering function of such  a linear 
system by the Theorem 3 of [46]. Also, $\psi_{(\pi )}(x)=\mu^2\psi (x)$ so is realised as $\psi_{(\pi )}(x)=\mu^2Ce^{-xA}B$.\par
\indent (iii) This follows from Corollary 5.6.\par
\rightline{$\square$}\par

\noindent {\bf Proof of Proposition 10.3} We let $\psi_{\lambda }$ be a family of solutions to
 $-\psi_\lambda''(x)+u(x)\psi_\lambda (x) =\lambda\psi_\lambda (x)$ which depends analytically upon $\lambda$ and satisfies $\psi_\lambda (x+\pi )=\mu\psi_\lambda (x)$ where $\mu^2+\Delta (\lambda )\mu +1=0$.  In particular, for $\lambda\in (\lambda_{2j+1}, \lambda_{2j+2})$, we can suppose that either $\vert \mu\vert <1$, in which case we tke $\phi_\lambda (x)=\psi_\lambda (x)$; or $\vert\mu\vert >1$, in which case we take $\phi_\lambda (x)=\psi_\lambda (-x)$. \par
\indent Let $\mu =e^{-\eta \pi -i\beta\pi}$. By Corollary 5.6, the exponentials $(e^{-\eta x-i\beta x-2ijx})_{j=0}^\infty$ give a Riesz basis for a closed linear subspace $W_{e^{-\eta\pi}}$  of $L^2(0, \pi)$ and there is an admissible  linear system $(-A,B,C)$ with input and output space ${\bf C}$ such that $\phi_\lambda (x)=Ce^{-xA}B$, With $\phi_\lambda (x)=\sum_{j=-\infty}^\infty a_j(\lambda) e^{-\eta x-i\beta -2ijx}$, we let
$$\eqalignno{B:{\bf C}\rightarrow W_{e^{-\eta \pi}}:&\quad b\mapsto b\sum_{j=-\infty}^\infty \vert a_j(\lambda )\vert^{1/2} e^{-\eta x-i\beta x- 2ijx};\cr
A:{\cal D}(A)\rightarrow W_{e^{-\eta\pi}}:&\quad f(x)\mapsto -f'(x);\cr
C:W_{e^{-\eta \pi}}\rightarrow {\bf C}: &\quad \sum_{j=-\infty}^\infty b_je^{-\eta x-i\beta x-2ijx}\mapsto \sum_{j=-\infty}^\infty \vert a_j(\lambda )\vert^{-1/2}a_j(\lambda )b_j.&(10.19)\cr}$$   
\indent Then
$$R_x=\Bigl[ {{a_je^{-2\eta x-2i\beta x -2ij x-2ikx}}\over {2\eta +2i\beta +ij+ik}}\Bigr]_{j,k=1}^\infty.\eqno(10.20)$$
\noindent where $a_j=\int_0^\pi\phi_\lambda (x)e^{i\beta x+i\eta x}e^{2ijx}dx/\pi$  and $R_x$ depends continuously on $\lambda$ in the trace class norm. Indeed, the matrix
$$\Bigl[{{1}\over{j^2(i\eta +j+k)}}\Bigr]$$
\noindent is trace class, and $\lambda \mapsto (j^2a_j(\lambda ))_{j=-\infty}^\infty$ is continuous ${\bf C}\setminus \sigma_B\rightarrow \ell^\infty$. We deduce that 
$\tau (x;\lambda )$ is uniformly continuous on compact sets of $(0, \infty ) \times ({\bf C}\setminus \sigma_B).$\par 
\indent For $\lambda$ in an interval of stability, $-2<\Delta (\lambda )<2$, and ${\bf p}_{\pm}=(-\lambda ,\mu_{\pm})$ where $\mu_{\pm}= 2^{-1}(-\Delta (\lambda )\pm i\sqrt{4-\Delta (\lambda )^2}))=e^{\mp i\pi \beta}$ is associated with a bounded solution of Hill's equation. The space $L^2(0, \pi )$ has a natural grading $L^2=W_\mu \oplus W_\mu^\perp$, where we choose $W_\mu$ to be the closed linear span of $( e^{-i\beta x+i2jx})_{j=1}^\infty$; note that the classical Hardy space on the circle $H^2$ is the closed linear span of $e^{ijx}$ for $j=0, 1, \dots $ in $L^2(0, 2\pi )$.\par 

\indent Suppose that $\lambda_j$ is a periodic eigenvalue with corresponding periodic or anti periodic eigenfunction $\psi_{\lambda_j}$; after translating and scaling $\psi_{\lambda_j}$ as needs be we can suppose that  $\psi_{\lambda_j}(0)=1$.
\indent Letting $\lambda\rightarrow \lambda_j$, we obtain a periodic linear system with scattering function $\phi_{\lambda_j}$. The matrix
$$ R_x=\Bigl[ {{a_je^{-2ij x-2ikx}}\over {ij+ik}}\Bigr]_{j,k=1}^\infty\eqno(10.21)$$
\noindent arises from the block of indices $j,k\geq 1$ as $\eta\rightarrow 0+,$ and defines a trace class operator. \par

\indent There are various tau functions associated with Hill's equation. Let $S=\cup_{j=0}^N [\lambda_{2j}, \lambda_{2j+1}]$ be a finite union of intervals of stability, and introduce 
$$\tau_n(\lambda_0, \lambda_1, \dots , \lambda_{2N+1})=\det_{j,k=0, \dots , n-1}\Bigl[ \int_S {{x^{j+k}dx}\over{\sqrt{4-\Delta (x)^2}}}\Bigr]\eqno(10.23)$$
\noindent be the sequence of Hankel determinants which are regarded as functions of the endpoints of $S$. Then $\tau_n$ satisfy systems of differential equations which arise in the theory of orthogonal polynomials.\par

\indent We now show how to obtain the curve ${\cal E}$ via the differential rings of sections 3 and 5.\par
\vskip.05in
\indent  Suppose that $\Delta (\lambda )^2-4$ has only simple zeros, and consider the transcendental curve ${\cal E}=\{ (\mu , \lambda ): \mu^2=\Delta (\lambda )^2-4\}$ which we regard as hyperelliptic and of infinite genus. We extend the classical language of algebraic curves to deal with this context.  For a meromorphic function $f$, let $n(r,f)$ be the number of poles, counted according to multiplicity in $\{ z\in {\bf C}: \vert z\vert \leq r\}$, and introduce the Nevanlinna characteristic by
$$T(r,f)=\int_0^{2\pi} \log_+ \vert f(re^{i\theta})\vert\, {{d\theta}\over{2\pi}} +\int_0^r (n(t,f)-n(0,f)){{dt}\over{t}}+n(0,f)\log r.\eqno(10.33)$$
The set ${\bf K}$ of $f$ such that $T(r,f)\leq {\gamma \sqrt{r}}+M$ for some $, \gamma ,M>0$ and all $r>0$ forms a field under the standard pointwise multiplication. Also, ${\bf K}$ 
has a subring ${\bf A}$ consisting of those entire functions $f$ such that $\vert f(z)\vert\leq Me^{\gamma \sqrt{\vert z\vert }}$ for some $\gamma ,M>0$ and all $z\in {\bf C}$.  By Hadamard's factorization theorem,  ${\bf A}$ is a unique factorization domain with ${\bf K}$ as its field of fractions. Let ${\bf A}^{\gamma}$ be the linear subspace of ${\bf A}$ consisting  of entire functions such that 
$\vert f(z)\vert \leq Me^{\gamma \vert z\vert^{1/2}}$ for some $M\geq 0$ and all $z\in {\bf C}$. 
\vskip.05in
\noindent {\bf Lemma 10.5} {\sl (i) Suppose that $(-A,B,C;E)$ is a
uniformly periodic linear system such that $\tau (x)$ is real and
nonzero on ${\bf R}$.\par
\indent (i) Then  for all $\kappa >0$, there exists $M_\kappa $ such that
$\lambda_{2n+2}-\lambda_{2n+1}<M_\kappa n^{-\kappa}$ for all $n$.\par 
\indent (ii) Under the hypotheses of Lemma 10.1, $\Delta (\lambda )$ belongs to ${\bf A}^{(1)}$ and ${\cal H}\subset {\bf A}^{(1)}$. }\par
\vskip.05in
\noindent {\bf Proof.} (i) There exists $\varepsilon >0$ such that
$\tau (z)$ is  entire, $\pi$-periodic and non-zero on $\{ z: -\pi \leq
\Re z\leq \pi ; -\varepsilon \leq \Im z\leq \varepsilon \}$. Hence
$u(x)$ is periodic and  real analytic on ${\bf R}$, so we can invoke
Hochstadt's estimates Theorem 2.13 of [MW].\par  
\indent The eigenvalues satisfy $\lambda_{2j}, \lambda_{2j-1}=j^2\pi^2 +O(1)$ 
as $j\rightarrow\infty$, so $\sqrt{\lambda_{2j+2}}-\sqrt{\lambda_{2j}}=\pi +O(1/j)$ as $j\rightarrow\infty$. 
By comparing $\Delta (\lambda )$ with $(\sin{\sqrt{\lambda}})/\sqrt{\lambda}$ we see that $\vert\Delta (\lambda )\vert \leq Me^{\sqrt\vert\lambda\vert}$ for some $M>0$ and all $\lambda$, and this cannot be improved.\par
\rightline{$\square$}\par
\indent We  
let ${\bf K}_{\cal E}={\bf K}[\sqrt{\Delta (\lambda )^2-4}]$ be a field of meromorphic functions on ${\cal E}$ and 
${\bf A}_{\cal E}={\bf A}[\mu ]/(\mu^2+4-\Delta (\lambda )^2)$, which is a ring of coordinate functions on ${\cal E}$. By an
Abelian integral of the first kind, we mean
$$\int_{x_0}^z {{f(\lambda )d\lambda}\over{\sqrt{4-\Delta (\lambda )^2}}}$$
\noindent for some $f\in {\bf A}$. As in Lemma 10.1, suppose that ${\cal H}$ is a reproducing kernel Hilbert space of entire functions with reproducing kernels $h_\lambda$ so $f(\lambda )=\langle f, h_\lambda \rangle$, and suppose that for all $\mu_j\in [\lambda_{2j+1}, \lambda_{2j+2}]$, the sequence $(h_{\mu_j})_{j=0}^\infty$ gives a Riesz basis for ${\cal H}$. Then for $\nu_j\in [\lambda_{2j+1}, \lambda_{2j+2}]$ such that $\sum_{j=0}^\infty \lambda_{2j+1}\lambda_{2j+2}(\nu_j-\lambda_{2j+1})^2/(\lambda_{2j+2}-\lambda_{2j+1})^2$ converges,  the coefficients that are determined by 
$$x_j=\sum_{k=0}^\infty \int_{\lambda_{2k+1}}^{\nu_k} {{h_{\lambda_{2j}}(x)dx}\over{\sqrt{\Delta (x)^2-4}}}.\eqno(10.25)$$
form $(x_j)_{j=0}^\infty \in \ell^2$, and there is a bounded linear functional
$$\sum_{j=0}^\infty \int_{\lambda_{2j+1}}^{\nu_j} {{f(x)dx}\over{\sqrt{4-\Delta (x)^2}}}=\sum_{j=0}^\infty x_jf(\lambda_{2j})
\qquad (f\in  {\cal H}).\eqno(10.26)$$

\noindent Let $B_k$ be a contour in ${\bf C}\setminus\cup_{k=0}^\infty [\lambda_{2k}, \lambda_{2k+1}]$ which goes once round $[\lambda_{2\ell}, \lambda_{2\ell+1}]$ for $\ell =0, \dots , k-1$ once in the positive sense. Let $A_j$ be a contour that goes once round $[\lambda_{2j-1}, \lambda_{2j}]$ once in the negative sense.\par
\indent  Let $X=[x_{jk}]_{j,k=0}^\infty$ be the matrix determined by the coefficients
$$x_{jk}=\int_{\lambda_{2j+1}}^{\lambda_{2j+2}}{{h_{\lambda_{2k}}(x)dx}\over{\sqrt{4-\Delta (x)^2}}},\eqno(10.27)$$
\noindent and suppose that $X$ is invertible with inverse $Y=[y_{jk}]_{j,k=0}^\infty$ as operators on $\ell^2$. Then the elements of ${\cal H}$ given by 
$$f_k(x)=\sum_{j=0}^\infty y_{jk}h_{\lambda_{2j}}, \quad g_j=\sum_{k=0}^\infty \bar x_{jk}h_{\lambda_{2k}}\eqno(10.28)$$
\noindent are biorthogonal in the sense that $\langle f_k, g_\ell\rangle_{\cal H}=\delta_{k,\ell}$ and  hence
$$\int_{\lambda_{2j+1}}^{\lambda_{2j+2}}{{f_k(x)dx}\over{\sqrt{\Delta (x)^2-4}}}=\delta_{jk}\qquad (j,k=0,1, 2, \dots ).\eqno(10.29)$$
\vskip.05in
\noindent {\bf Definition}  The $(f_k)_{k=0}^\infty$ give the canonical basis for ${\cal H}(E)$ corresponding to the homology basis. For some choice of square roots, the corresponding period matrix is $[I; \Omega ]$, where the entries of $\Omega$ are
$$\Omega_{nk}=\int_{B_n}{{f_k(x)dx}\over{\sqrt{\Delta (x)^2-4}}}.\eqno(10.30)$$
\vskip.05in
\indent  Then for any base point $x_0\in {\cal E}$, any element of the fundamental group $\pi_1({\cal E}, x_0)$ is homotopic to some product of powers of elements of the fundamental homology basis given by $B_k$ and loops $A_k$ around the intervals of instability. Let ${\bf L}$ be the additive subgroup of ${\bf C}^\infty$ that is generated by 
$$\Bigl( \int_{A_k}{{f_j(x)dx}\over{\sqrt{4-\Delta (x)^2}}}\Bigr)_{j=0}^\infty,\quad 
\Bigl( \int_{B_k}{{f_j(x)dx}\over{\sqrt{4-\Delta (x)^2}}}\Bigr)_{j=0}^\infty, \qquad (k=1,2, \dots ).\eqno(10.31)$$
Hence there is a natural group homomorphism $\pi_1({\cal E}, x_0)\rightarrow {\bf L}$, and define 
$H^1({\cal E}, {\bf Z})$ to be the range of this map.\par

\indent The notion of a principal divisor on a noncompact Riemann surface is
 open to interpretation. A $1$-chain $\gamma=\oplus_{j=1}^n \gamma_j$ is a direct
sum of continuous $\gamma_j:[0,1]\rightarrow{\cal E}$ such that $\gamma_j\vert
(0,1)$ is continuously differentiable; the boundary of $\gamma$ is
$\partial\gamma=\sum_{j=1}^n (\delta_{\gamma_j(1)}-\delta_{\gamma_j(0)}).$ Let $D\in {\hbox{Div}}({\cal E})_0$ and suppose that
there exists a $1$-chain $\alpha$ on ${\cal E}$ with boundary $\partial \alpha
=D$ such that $\int_\alpha {{f_j(\lambda )d\lambda }\over{\sqrt{\Delta (\lambda )^2-4}}}=0$ for all
$j$. We say that $D\in {\hbox{Div}}({\cal E})_P$.\par

\vskip.05in
\noindent {\bf Lemma 10.5} {\sl There is an additive group homomorphism}
$J_0:{\hbox{Div}}({\cal E})_0\rightarrow {\bf C}^\infty
/{\bf L}$ {\sl given by 
$$J_0:D\mapsto \Bigl( \int_{\gamma}{{f_j(\lambda )d\lambda }\over{\sqrt{\Delta (\lambda )^2-4}}}\Bigr)_{j=0}^\infty$$
\noindent where $\gamma$ is any $1$-chain in ${\cal E}$ such that 
$\partial \gamma =D$. The nullspace of $J_0$ is}
${\hbox{Div}}({\cal E})_P$.\par
\vskip.05in 
\noindent {\bf Proof.} First we check that $J$ is well-defined, by considering any $1$-chains $\gamma_1$ and
$\gamma_2$ on ${\cal E}$ such that $\partial \gamma_1 =D=\partial\gamma_2$. Then we introduce $\sigma =\gamma_1
-\gamma_1$ which is a $1$-cycle since $\partial \sigma
=D-D=0$, so there exist finitely many non zero integers $n_k$ and $m_k$ so that $\sigma
=\sum_k(n_kA_k+m_kB_k)$ in $H^1({\cal E}, {\bf Z})$. 
\noindent  Then 
$$\Bigl(\int_\sigma {{f_j(\lambda )d\lambda }\over{\sqrt{\Delta (\lambda
)^2-4}}}\Bigr)_{j=0}^\infty =\sum_kn_k\Bigl(\int_{A_k} 
{{f_j(\lambda )d\lambda }\over{\sqrt{\Delta (\lambda
)^2-4}}}\Bigr)_{j=0}^\infty+m_k\Bigl(\int_{B_k}{{f_j(\lambda )d\lambda }
\over{\sqrt{\Delta (\lambda )^2-4}}}\Bigr)_{j=0}^\infty\in {\bf L}.$$
\noindent This proves that $J(D)$ is independent of the 
particular choice of chain $\gamma$. Clearly $J$ is additive.\par
\indent Now let $D\in {\hbox{Div}}({\cal E})_P$. Then we can choose a $1$-chain $\alpha$ on ${\cal E}$ such that 
$\partial \alpha=D$ and $\int_\alpha {{f_j(\lambda )d\lambda }\over{\sqrt{\Delta (\lambda )^2-4}}}=0$ for all $j$, so 
$J_0(D)=0.$ Conversely, if $J_0(D)=0$ for all $D\in
{\hbox{Div}}({\cal E})_0$, then $D\in {\hbox{Div}}({\cal E})_p$.\par
\rightline{$\square$}\par
\noindent {\bf Definition} We define the Picard group to be 
$${\hbox{Pic}}({\cal E})_0={\hbox{Div}}({\cal E})_0/{\hbox{Div}}({\cal E})_P,$$
\noindent and define the Jacobian to be the injective group homomorphism
 $J:{\hbox{Pic}}({\cal E})_0\rightarrow {\bf C}^\infty/{\bf L}$ given by Lemma 10.5. We say that elements $D_1$ and $D_2$ of ${\hbox{Div}}({\cal E})_0$  are congruent and write $D_1\sim D_2$ if
$D_1-D_2\in {\hbox{Div}}({\cal E})_P.$ 
\vskip.05in
\indent When we translate the linear system $(-A, B,C; E)$ to
$(-A, e^{-sA}B, Ce^{-sA};e^{-sA}Ee^{-sA})$, we thereby translate $\tau
(x)$ to $\tau (x+s)$, so that $u(x)\mapsto u(x+s)$. We introduce
the auxiliary spectrum $(\mu_j(s))_{j=0}^\infty$, namely the 
spectrum of $-f''(x)+u(x+s)f(x)=\lambda f(x)$ with boundary conditions $f(0)=f(\pi )=0$ for 
$f\in C^2({\bf R};{\bf C})$, with $\lambda_{2j+1}\leq \mu_j(s)\leq \lambda_{2j+2}$, so $\mu_j(s)$ typically lies 
in an interval of instability.\par
\vskip.05in 
\noindent {\bf Theorem 10.6} (Baker's formula) {\sl Let $(a_j)_{j=0}^\infty$ be the Taylor coefficients of 
$\Delta (\lambda )^2-4$, and $(d_j)_{j=0}^\infty $ be the Taylor coefficients about $\lambda =0$ of $\Delta'(\lambda )$. 
Then  all the $a_j$ and $d_k$ are invariant under all isospectral flows. 
The auxiliary spectrum determines the potential by}
$$u(s)-\int_0^\pi u(t){{dt}\over{\pi}}=-2{{d}\over{ds}}\sum_{j=0}^\infty \int_{\lambda_{2j+1}}^{\mu_j(s)}
{{\Delta'(\lambda )d\lambda}\over{\sqrt{\Delta (\lambda )^2-4}}}.$$
\vskip.05in
\noindent {\bf Proof.} The discriminant function for the auxiliary spectrum is the entire function of order $1/2$ given by the convergent product
$$\bigtriangledown_s(\lambda )=\prod_{j=0}^\infty {{\mu_j (s)-\lambda}\over{\pi^2j^2}}\eqno(10.56)$$
\noindent which has derivative satisfying 
$$\bigtriangledown_s'(\mu_j)={{-1}\over {\pi^2j^2}}\prod_{k=1; k\neq j}^\infty{{\mu_k-\mu_j}\over{\pi^2k^2}}.\eqno(10.57)$$ 
\noindent Trubowitz has shown that the operation of translation on the potential leads to the following coupled system of ordinary differential equations for the auxiliary spectrum
$${{d\mu_j(s)}\over{ds}}=-{{\sqrt{\Delta (\mu_j(s))^2-4}}\over{\bigtriangledown_s'
(\mu_j(s))}}\qquad (j=0, 1,2, \dots )\eqno(10.58)$$
\indent  and the initial condition $\mu_k(0)$ is taken as the $k^{th}$ tied
eigenvalue. For $u\in C^\infty$, the Fourier coefficients satisfy 
$j^k\hat u(j)\rightarrow 0$ as $j\rightarrow\pm\infty$ for all $k\in {\bf N}$, and hence
$j^k(\lambda_{2j+2}-\lambda_{2j+1})\rightarrow 0$ as $j\rightarrow\infty$ for all $k\in {\bf N}$. We let $(\mu_j)_{j=0}^\infty \in\times_{j=0}^\infty [\lambda_{2j+1}, \lambda_{2j+2}]$ and introduce 
the principal function
$$S((\mu_j)_{j=0}^\infty ; (a_j)_{j=0}^\infty) =\sum_{j=0}^\infty \int_{\lambda_{2j+1}}^{\mu_j(s)}\sqrt{\Delta (t_j)^2-4}\, 
dt_j\eqno(10.59)$$
\noindent where the series is dominated by $\sum_{j=0}^\infty
(\lambda_{2j+2}-\lambda_{2j+1})^2/\sqrt{\lambda_{2j+2}\lambda_{2j+1}},$ hence converges.\par 
\indent We introduce the space of real divisors 
$${\hbox{Div}}({\cal E})_R=\Bigl\{ \sum_{j=0}^\infty (\delta_{\mu_j}-
\delta_{\nu_j}): \nu_j, \mu_j\in [\lambda_{2j+1},
\lambda_{2j+2}]; \lim\sup_{j\rightarrow\infty}e^{\kappa j}\Bigl\vert
{{\mu_j-\nu_j}\over{\lambda_{2j+2}-\lambda_{2j+1}}}\Bigr\vert=0;\forall \kappa >0\Bigr\}.$$

\indent Then for $\sum_{j=0}^\infty (\delta_{\mu_j}-\delta_{\lambda_{2j+1}})\in
{\hbox{Div}}({\cal E})_R$, we define the 
angle variables $(\varphi_j)$ that are conjugate to $(a_j)$ by
$$\varphi_k={{\partial S}\over{\partial a_k}}=
{{1}\over{2}}\sum_{j=0}^\infty 
 \int_{\lambda_{2j+1}}^{\mu_j(s)}{{t_j^kdt_j}\over{\sqrt{\Delta (t_j)^2-4}}};\eqno(10.60)$$
\noindent where the series is dominated by 
$$\sum_{j=0}^\infty
\lambda_{2j+2}^k\sqrt{{{\lambda_{2j+2}\lambda_{2j+1}(\mu_j-\lambda_{2j+1})}
\over{(\lambda_{2j+2}-\lambda_{
2j+1})}}},$$
\noindent hence converges. We introduce the map $J_{R}:{\hbox{Div}}({\cal E})_R\rightarrow {\bf
C}^\infty$ by $\sum_{j=0}^\infty (\delta_{\mu_j}-\delta_{\lambda_{2j+1}})\mapsto
(2\varphi_k)_{k=0}^\infty$. For comparison, the classical Abel map 
$J:{\hbox{Div}}({\cal E})_0\rightarrow {\bf C}^\infty $ is given by
$J:\sum_{j=1}^n(\delta_{\mu_j}-\delta_{\lambda_{2j+1}})\mapsto (2\varphi_k)_{k=0}^\infty$.\par
\indent In the same spirit, we introduce the phase map $\Phi_R :{\hbox{Div}}({\cal E})_R\rightarrow {\bf C}$ by 
$$\Phi_R
:\sum_{j=0}^\infty (\delta_{\alpha_j}-\delta_{\beta_j})\mapsto \sum_{j=0}^\infty \int_{\beta_j}^{\alpha_j}{{\Delta'(\lambda
)d\lambda }\over{\sqrt{\Delta (\lambda )^2-4}}}.$$  
\indent To compare these maps, we note that $\Delta'\in {\bf A}$  and hence the Taylor coefficients satisfy
$\vert d_k\vert\leq Me^{2k-k\log (2k/\gamma )}$ for some $M, \gamma>0$ and all $k$.
So we can multiply the identity () by $d_k$ and sum over $k$ to obtain a convergent expression
$$\sum_{k=0}^\infty d_k\varphi_k ={{-1}\over{2}}\sum_{j=0}^\infty  \int_{\lambda_{2j+1}}^{\mu_j(s)}{{\Delta'(t_j) dt_j}\over{\sqrt{\Delta (t_j)^2-4}}};
\eqno(10.62)$$
\noindent so that the right-hand side is the sum of the heights of the spikes above the
 partial intervals of instability $\cup_{k=0}^\infty [\lambda_{2k+1},\mu_k]$. The principal function and the phase function are thus related by
$$-\sum_{k=0}^\infty d_k{{\partial}\over{\partial a_k}}S= {{1}\over{2}}\sum_{j=0}^\infty 
 \int_{\lambda_{2j+1}}^{\mu_j(s)}{{\Delta'(t_j) dt_j}\over{\sqrt{\Delta (t_j)^2-4}}}.\eqno(10.64)$$
\indent The phases have have derivatives 
$$\eqalignno{{{d\varphi_k}\over{ds}}&={{1}\over{2}}\sum_{j=0}^\infty {{\mu_j(s)^k}\over{\sqrt{\Delta (\mu_j(s))^2-4}}}{{d\mu_j'(s)}\over{ds}}\cr
&=-{{1}\over{2}}\sum_{j=0}^\infty {{\mu_j(s)^k}\over{\bigtriangledown_s'(\mu_j(s))}}.&(10.61)\cr}$$
In turn,  the derivative of (10.62) is
$$\eqalignno{\sum_{k=0}^\infty d_k{{d\varphi_k}\over{ds}}&=-{{1}\over{2}}\sum_{j=0}^\infty {{\Delta'(\mu_j(s))}\over{\bigtriangledown'_s(\mu_j(s))}},\cr
&=-{{1}\over{2}}\sum_{j=1}^\infty (\lambda'_j-\mu_j)&(10.63)\cr}$$
\noindent where we have used an argument of McKean and Trubowitz [Bull AMS ]at the last step. A standard trace formula gives
$$\eqalignno{u(s)&=\lambda_0+\sum_{j=1}^\infty (\lambda_{2j-1}+\lambda_{2j}-2\mu_j(s))\cr
&=\lambda_0+\sum_{j=1}^\infty (\lambda_{2j-1}+\lambda_{2j}-2\lambda'_j)+2\sum_{j=1}^\infty (\lambda'_j-\mu_j(s))\cr
&=\int_0^{\pi} u(t){{dt}\over{\pi}}-2{{d}\over{ds}}\sum_{k=0}^\infty \int_{\lambda_{2k-1}}^{\mu_k(s)}
 {{\Delta'(\lambda )d\lambda }\over{\sqrt{\Delta (\lambda )^2-4}}}.\cr}$$
\noindent This is a version of Baker's formula for the hyperelliptic transcendental curve ${\cal E}$.\par
\rightline{$\square$}\par
\indent By integrating, we obtain the inversion formula
$$\tau (x)=\tau (0)\exp\Bigl( x(\log \tau
)'(0)+\int_0^x\sum_{k=0}^\infty
\int_{\mu_k(0)}^{\mu_k(t)} {{\Delta'(\lambda )d\lambda}\over{\sqrt{\Delta (\lambda
)^2-4}}}dt\Bigr)$$
\noindent which expresses $\tau$ in terms of the phase function on the real divisor.

\indent The previous proof introduces a version of the Abel map $J_R$ and an
inversion formula,
without specifying the range of $J_R$. In the following result, we
construct an infinite torus which is a geometrical model for the Picard
variety.\par

\vskip.05in
\noindent {\bf Proposition 10.7} {\sl Suppose that $\Delta (\lambda )^2-4$ has only 
simple zeros and let ${\cal E}=\{ (\mu , \lambda ): \mu^2=\Delta (\lambda )^2-4\}$.
 For $n=1, 2, \dots $ there exist a quartic $\alpha_n$ in finitely many variables such that the system of equations 
$$\alpha_n (x)=0\qquad (n=1, 2, \dots )$$
\noindent in infinitely many variables specifies a Jacobian ${\bf X}$ for the transcendental curve ${\cal E}$. The operation of translating
the potential is given by the canonical equations of motion for a Hamiltonian system with 
canonical coordinates that determine 
${\bf X}$.}\par
\vskip.05in

 \noindent {\bf Proof.} We can now go from the periodic linear system to the isospectral torus of Hill's equation. We introduce $\ell^2 =\ell^2({\bf Z}_+; {\bf C})$ with the complex bilinear form
 $\langle .\, ,\, \rangle :\ell^2\times \ell^2\rightarrow {\bf C}:$ 
$\langle x,y\rangle=\sum_{j=0}^\infty x_jy_j$. We also write $x\cdot \nabla $ for the formal differential operator $\sum_{j=0}^\infty x_j (\partial /\partial y_j)$. Then the tangent bundle to the unit sphere in the real Hilbert space $\ell^2({\bf R})$ 
has complexification
$$TS^\infty =\{ (x,y)\in \ell^2\times \ell^2 : \langle x,x\rangle =1, \langle x,y\rangle =0\}.\eqno(10.43)$$
\noindent Let $H^1=\{ x\in \ell^2 : \sum_{j=0}^\infty (1+j^2)\vert x_j\vert^2<\infty \}$, and given a strictly  increasing  sequence $(\lambda_{2j})$  and $M_1, M_2>0$ such that $\vert \lambda_{2j}- M_1j^2\vert \leq M_2$ for all $j$, we let 
$\Lambda :H^1\rightarrow \ell^2$ be the linear operator $\Lambda: (x_j)\mapsto (\lambda_{2j}x_j)$. On $H^1$ there is defined the Hamiltonian
$$H=2^{-1}(\Lambda x,x)+2^{-1}\bigl( \langle x,x\rangle \langle y,y\rangle -\langle x,y\rangle^2\bigr) -2^{-1}\langle \Lambda x,x\rangle \bigl( \langle x,x\rangle -1\bigr)\eqno(10.44)$$
\noindent which is associated with the derivation
$$D_{\bf X}=y\cdot \nabla_x-(\Lambda x)\cdot \nabla_y+\bigl( \langle \Lambda x,x\rangle -\langle y,y\rangle \bigr) x\cdot \nabla_y\eqno(10.45)$$
\noindent so that the canonical equations of motion are $\dot {x}=D_{\bf X}x$ and $\dot{y}=D_{\bf X}y$. Neumann considered $H$ as a model of particles that are subject to a quadratic potential $(\Lambda x,x)$ and 
are constrained to lie on a sphere given by $\langle x,x\rangle =1$  so the Hamiltonian incorporates a Lagrange multiplier to deal with the constraint. One can check that 
$$G_j(x,y)=x_j^2+\sum_{k:k\neq j}{{(x_jy_k-x_ky_j)^2}\over{\lambda_{2j}-\lambda_{2k}}}\qquad (j=0,1,\dots )\eqno(10.46)$$
\noindent are also invariant with respect to $D_{\bf X}$, in the sense that $D_{\bf X}G_j=0$, and likewise $D_{\bf X}\langle x,x\rangle =0$ and $D_{\bf X}\langle x,y\rangle =0$ on $TS^\infty$. \par
\indent Recall that $u$ has periodic spectrum $\lambda_0<\lambda_1\leq\lambda_2<\lambda_3\leq\dots ,$ so we need to introduce the odd indexed periodic eigenvalues.  Suppose that there exist $M>\gamma_j>0$ such that the function
$$F(\lambda )=\prod_{j=0}^\infty \Bigl( 1-{{\lambda}\over{\lambda_{2j}}}\Bigr)^2\sum_{j=0}^\infty {{\gamma_j}\over{\lambda -\lambda_{2j}}}\eqno(10.47)$$
\noindent has simple zeros $\lambda_j$ for $j=0, 1, \dots $. Then $F(\lambda )=c(4-\Delta (\lambda )^2)$ for some nonzero $c$. Then we define
$${\bf X}=\bigl\{ (x,y)\in \ell^2\times \ell^2: \langle x,x\rangle =1,
 \langle x,y\rangle =0, G_j(x,y)=\gamma_j, j=0,1,2, \dots \bigr\}\eqno(10.48)$$
\noindent to be the complex abelian manifold associated with $u$, and we observe that $D_{\bf X}$ is a tangent vector to ${\bf X}$.\par

\indent There is an injective map from the space of $\pi$-periodic potentials $u\in C^\infty ({\bf R}; {\bf R})$ such that Hill's equation has periodic spectrum $\lambda_0<\lambda_1<\dots $ to the points on 
the torus ${\bf X}$, where ${\bf X}$ is specified by $(\lambda_{2j})_{j=0}^\infty$ and $(\gamma_j)_{j=0}^\infty$. 
. When $\lambda_{2n}$ is a simple periodic eigenvalue of Hill's equation, 
there exists a periodic solution $f_n$ such that [Tru, p329], $-f_n''+uf_n=\lambda_{2n}f_n$,  
$\int_0^\pi f_n(s)^2ds/\pi =1$, which is given in terms of the discriminant by  
$$f_n(s)^2=-{{\bigtriangledown_s (\lambda_{2n})}\over{\Delta'(\lambda_{2n})}}\eqno(10.49)$$ 
We note in passing the difference between this and (10.61), then recall that McKean and Trubowitz [45'] have shown that there exist $\varepsilon_n>0$ such that 
$x_n(s)=\varepsilon_nf_n(s)$ and $y_n(s)=\varepsilon_n f_n'(s)$ such that $x(s)=(x_j(s))_{j=0}^\infty$ and $ y(s)=(y_j(s))_{j=0}^\infty$ determines a 
point $(x(s),y(s))\in TS^\infty$ for all times $s$. Then we can recover the potential from $(x(s),y(s))$ by taking the bilinear product of $-x''+u x=\Lambda x$ with  $x$, and noting that
 $\langle x,x''\rangle +\langle y,y\rangle =0$; hence 
$$u(s)=\langle \Lambda x(s), x(s)\rangle -\langle y(s), y(s)\rangle .\eqno(10.50)$$
\indent Note that $t_1\mapsto (t_1, t_2, \dots )$ gives a line in ${\bf C}^\infty$ and that for any function $G:{\bf X}\rightarrow {\bf C}$ that depends upon only finitely many coordinates, we have ${{\partial }\over{\partial t_1}}G=D_{\bf X}G$. The 
function $t\mapsto u(t_1, t_2, \dots )$ is thus the restriction to a line of a function of 
rational character on ${\bf X}$. The divisor $D(t_2, t_3, \dots )$ of $t_1\mapsto 
u(t_1, t_2, \dots )$  may be treated as a function of the parameters $(t_2, t_3, \dots )$.\par
\indent Observe that $({\bf T}^\infty )_0$ acts on ${\bf X}$ by $(x(s), y(s))\mapsto (x(s+t),y(s+t))$ for $(x,y)\in {\bf X}$ and $s,t\in ({\bf T}^\infty)_0$.\par
\indent The definition of ${\bf X}$ involves infinite series involving the $(x_j,y_j)$. 
Next we show that ${\bf X}$ is also determined by an infinite sequence of quartic equations, each
of which involves only finitely many variables. We wish to replace $Z^2=4-\Delta^2$ by $Z^2=UW+V^2$, and to do so, we introduce the entire functions 
$$F_1(\lambda )=\prod_{j=0}^\infty \Bigl( 1-{{\lambda}\over{\lambda_{2j}}}\Bigr),\quad 
U(\lambda) =F_1(\lambda) \sum_{k=0}^\infty {{x_k^2}\over{\lambda -\lambda_{2k}}},
\eqno(10.51)$$
$$V(\lambda )=iF_1(\lambda )\sum_{k=0}^\infty {{x_ky_k}\over{ \lambda -\lambda_{2k}}},
\quad W(\lambda )=F_1(\lambda )\Bigl( 1+\sum_{k=0}^\infty {{y_k^2}\over{\lambda -\lambda_{2k}}}\Bigr).\eqno(10.52)$$
\noindent For $\gamma >0$, let ${\bf A}^{(\gamma )}$ be the space  of entire functions $f$ 
such that $\vert f(z)\vert\leq Me^{\gamma \vert z\vert^{1/2}}$ for all $z\in {\bf C}$ and 
some $ M\geq 0$. Then 
$U,V,W\in {\bf A}^{(1)}$, and by using the special choice of $(x,y)\in {\bf X}$ and considering
the cancellation of apparent poles at $\lambda =\lambda_{2j}$, one sees that the functions satisfy
$$U(\lambda )W(\lambda )+V(\lambda )^2=F(\lambda );$$
\noindent the details are provided by Mumford in 3.157-8 of [Mum]. Reversing this argument, one can recover the defining
equations for ${\bf X}$ from the identity (). Thus the coefficients in
the power series expressions for $U,V,W$ give a system of coordinates 
for ${\bf X}$. We introduce $U(\lambda )=\sum_{j=0}^\infty u_j\lambda^j$, $V(\lambda )=\sum_{j=0}^\infty v_j\lambda^j$ and $W(\lambda )=\sum_{j=0}^\infty w_j\lambda^j$ and consider 
$$\Delta (\lambda )^2-4-U(\lambda )W(\lambda )-V(\lambda )^2=\sum_{j=0}^\infty \alpha_j (U,V,W)\lambda^j.\eqno(10.53)$$
\noindent where
$$\alpha_0=u_0v_0+v_0^2+4-a_0\eqno(10.54)$$
$$\alpha_\ell =\sum_{j=0}^\ell (u_jw_{\ell -j}+v_jv_{\ell -j})-a_\ell\qquad (\ell=1, 2, \dots ).\eqno(10.55)$$ 
\noindent Then $\alpha_\ell =0$ for $\ell =0, 1, \dots $ is a system of quartic equations 
each of which involves only finitely many variables that characterizes 
$\Delta (\lambda )^2-4=U(\lambda )W(\lambda )+V(\lambda )^2$. \par
\vskip.05in

\noindent {\bf 8 Theta functions on Abelian manifolds}\par
\vskip.05in
\indent The notion of a tau function was introduced to generalize the notion of a theta 
function of an algebraic curve, so we now clarify how these are related.\par
\vskip.05in
\noindent {\bf Definition} (Riemann's theta function) Suppose that $\Omega_0$ and $\Omega_1$ are real symmetric
$g\times g$ matrices with $\Omega_1$ positive definite, and let $\Omega
=\Omega_0+i\Omega_1$; then let $\Lambda={\bf Z}^g+\Omega {\bf Z}^g$ be a lattice in ${\bf
C}^g$. Then
$$\theta (x\mid \Omega )=\sum_{m\in {\bf Z}^g} e^{2\pi im^Tx+\pi
im^T\Omega m}\eqno(8.1)$$
\noindent is Riemann's theta function for the Abelian manifold 
${\bf X}={\bf C}^g/\Lambda$ and $m^T$ denotes the row with integral entries.\par
\noindent Then $\theta (x\mid \Omega )$ is entire and quasi periodic with respect to $\Lambda$  in the sense that
$$\theta (x+m\mid \Omega )=\theta (x\mid \Omega ), \quad \theta (x+\Omega m\mid \Omega )=\theta (x\mid \Omega )e^{-2\pi im^Tx-\pi im^T\Omega m}\qquad (x\in {\bf C}^g, m\in {\bf Z}^g).$$
\vskip.05in
\noindent {\bf Proposition 8.1} {\sl Given a Riemann theta function $\theta (x\mid \Omega )$, there
exists a linear system such that}
$$\theta (x\mid \Omega )=\det (I+U(x)BC).$$
\vskip.05in
\noindent {\bf Proof.} Let $H=\ell^2({\bf Z}^g)$ and introduce the linear system $(0,B,C)$ with periodic deformation group $(U(x))_{x\in {\bf T}^g}$ by 
$$\eqalignno{ B:{\bf C}\rightarrow H:&\quad s\mapsto (e^{\pi im^T\Omega m/2}s)_{m\in {\bf Z}^g};\cr
U(x): H\rightarrow H: &\quad (b_m)_{m\in {\bf Z}^g} \mapsto (e^{2\pi im^Tx}b_m)_{m\in {\bf Z}^g};\cr
C:H\rightarrow {\bf C}:&\quad (b_m)_{m\in {\bf Z}^g}\mapsto \sum_{m\in {\bf Z}^g\setminus \{0\}}e^{\pi im^T\Omega m/2}b_m.\cr}$$
\noindent Then $B$ and $C$ have rank one, so we can realise Riemann's theta function as 
$$\theta (x\mid \Omega )=\det (I+U(x)BC).$$
\rightline{$\square$}
\indent For any compact Riemann surface
${\cal E}$ of genus $g$, one can define a homology basis and a
$g$-dimensional space of Abelian differentials of the first kind. Then
one defines a corresponding lattice $\Lambda$ of periods and a 
Jacobi variety ${\bf J}={\bf C}^g/\Lambda$ with a period matrix
$\Omega$, and hence  
Riemann's theta function $\theta (x\mid \Omega )$ by (3.20). Schottky [45, 54] asked how one can 
characterize the $\theta$ functions that arise from compact Riemann 
surfaces amongst all the possible functions $\theta (x\mid \Omega )$ on
Abelian varieties as in (5.1). \par
\indent We choose $\kappa ,\gamma, \delta ,\zeta\in {\bf R}^\infty$ and introduce the affine
subspace $\{t=\kappa x+\gamma y+\delta s=\zeta :x,y,z\in {\bf R}\}$ and thereby reduce to
functions of three variables $(x,y,s)$.\par
\vskip.05in

\noindent {\bf Proposition 8.2} (Shiota) {\sl  Suppose that $\theta (t\mid\Omega )$ is 
Riemann's theta function for some Abelian manifold ${\bf X}={\bf C}^g/\Lambda$ of dimension $g$; let 
$Q(x,y,s)$ be a quadratic form, let $\kappa, \gamma,\delta ,
\zeta\in {\bf C}^g$ with $\kappa\neq 0$, and for
$$\sigma (x,y,s;\zeta )=e^{Q(x,y,s)}\theta (\kappa x+\gamma y+\delta
s+\zeta \mid \Omega ),\eqno(10.1)$$
\noindent let $w(x,y,s; \zeta )=-2{{\partial^2}\over{\partial
x^2}}\log\sigma (x,y,s;\zeta ).$ Then the following two
conditions are equivalent:\par
\indent (i) ${\bf X}$ is isomorphic to the Jacobian variety of a
complete algebraic curve;\par
\indent (ii) the $\theta$ divisor is irreducible, and $w$ satisfies the $KP$
 equation  
$${{\partial}\over{\partial x}}\Bigl( {{\partial^3 w}\over{\partial
x^3}}-6w{{\partial w}\over{\partial x}}-4{{\partial w}\over{\partial s}}\Bigr)
 +3{{\partial^2w}\over{\partial y^2}}=0,\eqno(6.48)$$
\noindent for some $\kappa ,\gamma, \delta$ and all 
$\zeta\in {\bf C}^g$; moreover, $w$ is the restriction of a periodic
function to an affine space.}\par

\vskip.05in
\noindent {\bf Proof.} This is contained in the papers of Shiota [54] and Mulase [45]. Our statement
slightly different from these papers, so we note the following. 
Given a smooth algebraic curve of genus $g$, the Jacobian is a complex abelian group
with universal cover ${\bf C}^g$; the theta function is entire on ${\bf
C}^g$. There exists a $g\times \infty$ complex matrix $\bar A$ of rank $g$ such
that $\tau (t) =e^{Q(t)}\theta (\bar A t+\zeta )$ gives a solution to
(3.29) when $t=(x,y,s, 0, \dots )$. There exists a
$g$ dimensional vector subspace $V$ of ${\bf C}^\infty$ such that $\{ t_0\in V:
\bar At_0\in {\bf Z}^g+\Omega {\bf Z}^g\}$ is a lattice, and for such $t_0$  
there exists $\rho_{t_0,
\zeta}\neq 0$ such that $\tau (t+t_0)=\rho_{t_0, \zeta }\tau (t)$ for all $t$. 
Hence $w$ restricted to $t\in V+\zeta $ is periodic.\par 

\rightline{$\square$}\par

\vskip.05in
\indent We proceed to introduce the notion of a theta function on an
infinite-dimensional Abelian manifold, taking (6.3) as our . Let $\ell^2({\bf C})=\{ x=(x_\nu
)_{\nu\in {\bf Z}^\infty}: x_\nu \in {\bf C}: \sum_{\nu\in {\bf
Z}^\infty}\vert x_\nu\vert^2<\infty\}$ be the Hilbert space with 
inner product $\langle x,y\rangle =\sum_{\nu\in {\bf Z}^\infty}x_\nu \bar
y_\nu$, and let $\ell^2$ the the real part of $\ell^2({\bf C})$.  
Let $({\bf Z}^\infty
)_0$ be the space of finitely supported integral sequences, and note that $({\bf Z}^\infty )_0
=\ell^2\cap {\bf Z}^\infty$.
Let $Q$ be a self-adjoint  and positive operator on
$\ell^2({\bf C})$ that is trace class, injective and maps $\ell^2$ into
$\ell^2$. As in the spectral theorem, let $(e_\nu
)_{\nu\in {\bf Z}^\infty}$ be the standard complete orthonormal basis of
$\ell^2({\bf C})$ so that $Qe_j=q_je_j$, and map $({\bf Z}^\infty)_0\rightarrow
\ell^2$ by $(\nu_j)\mapsto \sum_{j=1}^\infty \nu_je_j$.\par
\indent Let ${\bf T}={\bf R}/{\bf Z}$ be the real torus, which we identify with the circle $\{ e^{2\pi it }:
0\leq t <1\}$ by the map $t\mapsto e^{2\pi it}$.  
A lattice is a countable additive subgroup of a Hilbert space; in particular, $\Lambda_0=({\bf Z}^\infty)_0
+iQ({\bf Z}^\infty)_0$ and $\Lambda_1=({\bf Z}^\infty)_0
+iQ^{-1}({\bf Z}^\infty)_0$ are lattices in $\ell^2({\bf C})$. \par
\vskip.05in
\noindent {\bf Definition} We define ${\bf X}_0=\ell^2({\bf C})/\Lambda_0$
and ${\bf X}_1=\ell^2({\bf C})/\Lambda_1$ to be the dual pair of Abelian manifolds associated with $Q$ for
the real bilinear pairing $(\xi ,\eta )\mapsto \Re \langle \xi ,\eta \rangle$ on $\ell^2({\bf
C})\times \ell^2({\bf C})$. The real parts are $\Re {\bf
X}_0=\Re {\bf
X}_1=\ell^2/({\bf Z}^\infty )_0$, and the imaginary parts are $\Im {\bf
X}_0=\ell^2/Q({\bf Z}^\infty )_0$ and  $\Im {\bf
X}_1=\ell^2/Q^{-1}({\bf Z}^\infty )_0$.\par
\par 
\vskip.05in
\noindent {\bf Lemma 8.2.} {\sl (i) $\Re {\bf X}_0$ is subgroup of 
the compact, metrizable and abelian group ${\bf T}^\infty$.\par 
\indent (ii) Let $\Lambda_0$ have the discrete topology. Then 
 ${\hbox{Hom}}(\Lambda_0; {\bf T})$ is compact and metrizable abelian group which
contains ${\bf X}_1$ as a subgroup.\par
\indent (iii) Let $\Lambda_1$ have the discrete topology. Then 
 ${\hbox{Hom}}(\Lambda_1; {\bf T})$ is compact and metrizable abelian group which
contains ${\bf X}_0$ as a subgroup.}\par      
\vskip.05in 
\noindent {\bf Proof.} (i) The lattice $({\bf Z}^\infty)_0$ is
countable, so ${\bf T}^{({\bf Z}^\infty)_0}$ is metrizable, and compact for the
Tychonov topology.\par
(ii) We note that $\Lambda_1=\{ \eta\in \ell^2 ({\bf C}):
\Re\langle \eta , \xi\rangle \in {\bf Z}, \forall \xi\in \Lambda_0\}$ and likewise 
$\Lambda_0=\{ \eta\in \ell^2 ({\bf C}):
\Re\langle \eta , \xi\rangle \in {\bf Z}, \forall x\in \Lambda_1\}$.\par
Let ${\hbox{Hom}}(\Lambda_0; {\bf T})$ be the space of
characters on $\Lambda_0$. We associate $\eta\in \ell^2({\bf C})$ with the
character $\varphi_\eta :\xi\mapsto e^{2\pi i\Re\langle \xi , \eta\rangle}$, and deduce that 
$\ell^2({\bf C})/\Lambda_1$ is a subgroup of ${\hbox{Hom}}(\Lambda_0; {\bf T})$ which 
is compact 
for the Tychonov product topology and metrizable, as in (i).\par   
\indent (iii) As in (ii).\par 
\rightline{$\square$}\par
\vskip.05in

\indent We introduce the Hilbert space $H_0=\{ t\in \ell^2: Q^{-1/2}t\in \ell^2\}$ with the inner
product $\langle t,s\rangle_{H_0}=\langle Q^{-1}t,s\rangle_{\ell^2}$;
then let $H_1$ be the completion of $\ell^2$ for the norm associated
with the pre inner product $\langle t, s\rangle_{H_1}=\langle
Q^{1/2}t,Q^{1/2}s\rangle_{\ell^2}$ for $s,t\in \ell^2$. This gives us Hilbert spaces
with dense linear inclusions $H_0\subset \ell^2\subseteq H_1$.
\vskip.05in

\noindent {\bf Definition} (Theta functions) A type is a pair 
$(L, J)$ such that 
$L:\Lambda_1\rightarrow \ell^2({\bf C})$ and $J:\Lambda_1\rightarrow
{\bf C}$. A function $\theta :H_0+iH_0\rightarrow {\bf C}$ is entire if
$z\mapsto \theta (\xi +z\eta )$ is entire for all $\xi, \eta \in H_0+iH_0.$ A
quotient $\theta$ of nonzero entire functions on $H_0+iH_0$ is said to
be a theta function if there exists a type $(L, J)$ such that 
$\theta (z+\gamma )=e^{2\pi i(\langle z, L(\gamma )\rangle +J(\gamma
))}\theta (z)$ for all $\gamma\in \Lambda_1$ and $z\in H_0$.\par
\vskip.05in
\indent  Given $\psi\in H_1$, the function 
$f(\xi )=e^{2\pi i\langle Q\xi \, \bar xi \rangle -\langle \xi
\psi\rangle -c}$ give a Gaussian theta function, which is regarded as trivial in the present
discussion.\par
\vskip.05in
\noindent {\bf Proposition 8.3} {\sl Let $Q$ be a positive and trace
class operator on $\ell^2$ such that $Q\leq I$.\par
\indent (i)  Then $Q$ is the correlation function of a
Gaussian measure on $\ell^2$ that induces a function $\theta_0$ on
${\bf T}^\infty$.\par
\indent (ii) There exists an entire theta function
$\theta_1$ 
with respect
to $\Lambda_1 =({\bf Z}^\infty)_0 +iQ^{-1}({\bf Z}^\infty)_0/2\pi$ 
with type $(L, J)$ where
$$L(\alpha +i\beta ) =-4\pi^2iQ\beta ,\qquad  J(e_j+iQ^{-1}e_k/2\pi )
=2^{-1}\langle Q^{-1}e_k, e_k\rangle,\eqno(9.66)$$
\noindent so that $\theta_1(x)/\theta_1(x-i\beta )$ is a meromorphic
theta function for ${\bf X}_1$ with type $(0, L(i\beta ) )$, where
$\beta\mapsto L(i\beta )/i$  gives a natural
isomorphism}
$$\Im {\bf X}_1\sim {\hbox{Hom}}(\Im \Lambda_1 , {\bf T}).\eqno(9.67)$$    
\indent {\sl (iii) There exists a periodic 
linear
system $(-A, B,C;E)$ that has corresponding tau function $\theta_1 (x)=\tau (x)$.}\par
\indent {\sl (iv) Suppose that $\Omega $ from (9.23) satisfies $\Omega =iQ^{-1}$. Then ${\bf L}=\Lambda_1$.}\par
\vskip.05in
\noindent {\bf Proof.} (i) There
exists a Gaussian measure $P_0$ on this triple $H_0\subset \ell^2\subset H_1$ such that the
characteristic function satisfies
$$\int_{H_1}e^{i\langle t,\xi\rangle_{\ell^2}} P_0(dt)=\exp \bigl(-
\langle Q\xi ,\xi\rangle_{\ell^2}/2\bigr)\qquad (\xi\in \ell^2).\eqno(9.68)$$
\noindent By a standard result on Gaussian measures, there exists $\varepsilon >0$ such that $\int
e^{2\pi^2\varepsilon \langle Q x,x\rangle } P_0(dx)$ is finite. \par
\indent Next we regard $\ell^2({\bf R})$ as a subset of ${\bf
R}^\infty$, and consider the integer lattice ${\bf Z}^\infty$ in 
${\bf R}^\infty$. The group ${\bf T}^\infty ={\bf R}^\infty /{\bf
Z}^\infty$ is compact for the Tychonov product topology and hence has a unique Haar probability measure $\Theta$.
Let $\varphi :{\bf R}^\infty\rightarrow {\bf T}^\infty$
$:(t_j)_{j=1}^\infty\mapsto (\{ t_j\})_{j=1}^\infty$ be
the quotient map with respect to ${\bf Z}^\infty$, and $\{ x\}$ is the fractional part of a
real $x$. Note also that ${\bf Z}^\infty
\cap \ell^\infty$ is a subgroup of $\ell^\infty$ and that ${\bf T}^\infty
=\ell^\infty /(\ell^\infty \cap {\bf Z}^\infty )$. Then $\varphi$ induces
a probability measure $P_0'$ on ${\bf T}^\infty$, so that $\int_{{\bf
T}^\infty} f(t)P_0'(dt)=\int_{{\bf R}^\infty}f(\varphi (y))P_0(dy)$ for all
bounded and continuous function $f$ that depend on only finitely many
coordinates. We define $\theta_0 (t)$ to be the Radon--Nikodym derivative
$dP_0'/d\Theta$.\par 
\indent  Let
$\nu=(\nu_j)_{j=1}^\infty$ be a finitely supported sequence of
integers and $t=(t_j)_{j=1}^\infty \in {\bf R}^\infty$. Then $\langle
t, \nu\rangle =\sum_{j=1}^\infty t_j\nu_j$, and $e^{2\pi i\langle t,\nu
\rangle}$ depends only $\varphi (t)$, hence $t\mapsto e^{2\pi i\langle t,\nu \rangle}$
 defines a character
of ${\bf T}^\infty$. Then
$$\int_{{\bf T}^\infty} e^{-2\pi  i\langle t, \nu \rangle} \theta (t)\Theta
(dt)=\int_{{\bf R}^\infty} e^{-2\pi i\langle y, \nu \rangle}P(dy)
=e^{-2\pi^2\langle Q\nu , \nu\rangle}\eqno(9.69)$$
\noindent so that there is a Fourier series expansion 
$$\theta_0 (t)\sim \sum_{\nu\in ({\bf Z}^\infty )_0} e^{2\pi i\langle t,
\nu\rangle -2\pi^2\langle Q\nu , \nu
\rangle}\eqno(9.37)$$
\noindent where the sum is over all finitely supported $\nu$ in the
lattice. For $\xi\in H_0$, the translation $t\mapsto t+\xi$ induces a
Gaussian probability measure $P_\xi$ with mean $\xi$ from $P_0$.
Moreover, $P_\xi$ and $P_0$ are equivalent in the sense that they have
the same collection of Borel null sets. Thus one can derive $\theta
(t+\xi ).$ Certainly $Qe_j\in H_0$, so by manipulating the Fourier series, one can verify that the
standard orthonormal basis $(e_j)$ of $\ell^2$ satisfies
$$\theta_0 (x-2\pi i Qe_j)=\exp (2\pi i \langle e_j, x\rangle
+2\pi^2\langle Qe_j, e_j\rangle)\theta_0 (x).\eqno(9.70)$$
\indent (ii) Now $f\mapsto \int f(Q^{-1}x)P(dx)$ is a positive linear functional on
$C_b(H_1; {\bf R})$ and in particular, we have
$$\int e^{i\langle \xi , Q^{-1}x\rangle }P_0(dx)=e^{-\langle Q^{-1}\xi ,
\xi\rangle}\qquad (\xi \in H_1).\eqno(9.71)$$

\noindent There exists a finitely additive set function $\hat P_1$ on $H_1$ such that 
$\int e^{i\langle \xi , y\rangle} \hat P_1(dy)=e^{-\langle Q^{-1}\xi , \xi \rangle}$ for
all $\xi\in H_1$; we do not claim that $\hat P_1$ is a Radon (inner regular) probability measure. By analogy with the foregoing, we define
$$\theta_1(t) =\sum_{\nu\in ({\bf Z}^\infty)_0} e^{2\pi i\langle \nu ,t\rangle
-\langle Q^{-1}\nu , \nu \rangle/2}\eqno(9.72)$$
\noindent which we interpret as the density of the finitely additive set function on
${\bf T}^\infty$ that
is induced by $\varphi$ from $\hat P_1$.\par 
\indent Let $Q$ have eigenvalues $q_j$, listed according to multiplicity, and
observe that $\langle Q^{-1}\nu ,\nu \rangle \geq \sum_{j=1}^\infty
q_j^{-1}\nu_j^2$, so 
$$\eqalignno{\sum_{\nu \in ({\bf Z}^\infty )_0} \exp \bigl(-\alpha \langle Q^{-1}\nu ,
\nu\rangle\bigr)&\leq  \sum_{\nu \in ({\bf Z}^\infty )_0} \exp
\bigl(-\alpha\sum_{j=1}^\infty q_j^{-1}\nu_j^2\bigr)\cr
&\leq \prod_{j=1}^\infty \sum_{\nu_j=-\infty}^\infty \exp \bigl( -\alpha
\nu_j^2q_j^{-1}\bigr)\cr
&\leq \prod_{j=1}^\infty \Bigl(
{{1+e^{-\alpha/ q_j}}\over{1-e^{-\alpha /q_j}}}\Bigr),&(9.73)\cr}$$
\noindent for some $C(\alpha )>0$, where the product converges since $Q$ is
trace class. Furthermore,  
$$\vert \exp\bigl( {2\pi i\langle t,\nu \rangle }\bigr\vert \leq \exp\bigl(\pi^2\Vert
Q^{1/2}t\Vert^2/\varepsilon\bigr)\exp\bigl(\varepsilon\Vert Q^{-1/2}\nu\Vert^2\bigr),\eqno(9.74)$$
 and $\sum_{\nu\in ({\bf Z}^\infty )_0}e^{2\pi^2(\varepsilon -1)\langle
Q^{-1}\nu ,\nu \rangle}$ converges,  so we have convergence for
all $t\in H_1$ and $\theta_1 (t)$ is entire.

Clearly
$$\theta_1 (t+e_j)=\theta_1(e_j);$$
$$\theta_1\Bigl( t+{{iQ^{-1}e_k}\over{2\pi}}\Bigr)=e^{\langle
Q^{-1}e_k,e_k\rangle/2+2\pi i\langle e_k, t\rangle}\theta_1(t),\eqno(9.75)$$
\noindent so $(L_1, J_1)$ is a well defined type for the lattice 
$\Lambda_1=({\bf Z}^\infty )_0+iQ^{-1}({\bf Z}^\infty )_0$ in $\ell^2+iH_1$,  
, and $\theta_1 $ has type $(L, J)$. \par
\indent (iii) Let $Z$ be the set of finitely supported sequences of integers
$\nu =(\nu_j)_{j=0}^\infty$ that are not identically zero, so $Z$ may
be identified with the nonzero integral polynomials. Then let $H$ be
the Hilbert space $\{ (a_\nu)_{\nu\in Z}:a_\nu\in {\bf C}; \sum_{\nu\in
Z}\vert a_\nu\vert^2<\infty\}$ with the inner product $\langle
 (a_\nu )_{\nu\in Z},(b_\nu )_{\nu\in Z}\rangle_H=\sum_{\nu\in
Z}a_\nu \bar b_\nu $. Next we introduce the linear system
$$\eqalignno{B:{\bf C}\rightarrow H:&\quad  \beta\mapsto (  
e^{-\langle Q^{-1}\nu , \nu
\rangle/4})_{\nu\in Z}\beta ;\cr
U(t):H\rightarrow H:&\quad (a_\nu)_{\nu \in Z}\mapsto (e^{2\pi i\langle t,
\nu\rangle }a_\nu )_{\nu\in Z};\cr
C:H\rightarrow {\bf C} :&\quad (a_\nu)_{\nu\in Z}\mapsto 
\sum_{\nu\in Z}e^{-\langle Q^{-1}\nu , \nu
\rangle/4}a_\nu;&(9.76)}$$
\noindent then $BC$ has rank one, so from the Fourier expansion, we have
$$\theta_1 (t)=\det (I+U(t)BC).\eqno(9.77)$$
\indent We fix $t_2, t_3,\dots$ and consider the multiplier $A:(a_\nu )_\nu \mapsto (-2\pi
i\nu_1a_\nu )_\nu$ which generates a strongly continuous and unitary group of operators on $H$. 
Then $R_{t_1}=U(t/2)BCU(t/2)$ satisfies ${{\partial}\over{\partial t_1}}R_{t_1}
=-AR_{t_1}-R_{t_1}A$ and hence we have a solution of Lyapunov's 
equation, and can introduce the uniform periodic linear system 
$$(-A,  U(0, t_2/2,
t_3/2, \dots )B, CU(0, t_2/2, t_3/2, \dots ); U(0, t_2/2, t_3/2, \dots )
BCU(0, t_2/2, t_3/2, \dots )),\eqno(9.78)$$ 
so that $\theta (x,t_2, t_3, \dots )=\det (I+e^{-xA}U(0, t_2,
t_3, \dots )BC)$ is the corresponding tau function.\par  

\rightline{$\square$}\par

\indent The addition rule makes (9.21) an obvious analogue of the Tracy--Widom integrable operator in these special cases. 
Suppose that $q$ is an elliptic function of the first kind with fundamental periods $2\omega_1$
and $2\omega_2$, where $\Re \omega_2/\omega_1\geq 0$; usually we take $\pi =2\omega_1$, and assume
that $q$ has no poles on the real line. We consider the differential equation $-\psi''(x)+q(x)\psi(x)=\lambda \psi(x)$ and proceed to introduce an appropriate class of integrable
operators for the complex torus.\par

\vskip.05in

\noindent {\bf Theorem  8.4} {\sl Suppose that $q$ is elliptic, that $q$ is real-valued with no 
poles on ${\bf R}$, and that $q(x)$ is a
finite gap potential for Hill's equation; suppose that $\Delta (\lambda )^2>4$. 
for some $\lambda \in {\bf R}$.\par
\indent (i) Then there exists a nontrivial solution
$\psi$ that is elliptic of the second kind and 
$\psi 
(x+\pi)=\mu\psi (x)$ for some $\vert \mu\vert <1$, so for some $c\in {\bf
C}$ the kernel
$$K(x,y) ={{\theta_1(x+y)e^{c(x+y)}}\over{\theta_1(x)^2\theta_1(y)^2}}{{\psi(x)\psi'(y)-\psi'(x)\psi(y)}\over{\wp
(x)-\wp (y)}},\eqno(8.13)$$
\noindent is bounded and Hilbert--Schmidt. For $0<\varepsilon <1$, there exist a
Hilbert space $H$ and functions $F_\varepsilon\in L^2((0, \infty ); H)$ 
and $G_\varepsilon\in L^2((0,
\infty ); H)$ such that the operator $K_{\varepsilon}$ with kernel
$K(x+i\varepsilon , y+i\varepsilon )$ satisfies $K_\varepsilon
=\Gamma_{F_\varepsilon}\Gamma_{G_\varepsilon}$.\par
\indent (ii) Suppose further that $q$ is even. Then there exists a
Hilbert space $H$, and $F\in L^2((0, \infty ); H)$ and $G\in L^2((0,
\infty );H')$ such that $K=\Gamma_F\Gamma_G$.}\par
\indent {\sl (iii) In the high density limit, as $i\omega_2\rightarrow -\infty$, the kernel reduces to
a trigonometric kernel} 
$$K(x,y)\rightarrow {{\psi(x)\psi'(y)-\psi'(x)\psi(y)}\over{16\pi^2\sin \pi (y-x)}}.
\eqno(8.14)$$
\vskip.05in
\noindent {\bf Proof.} (i) Gesztesy and Weikard [28] considered 
$-\psi''+q\psi =\lambda\psi$ for $q\in {\bf K}_{\cal
T}^1$, and
showed that $q$ is finite gap if and only if $q$ is a Picard potential. 

When $\lambda$ lies in a real interval of instability, there 
exists a Floquet multiplier $\mu_1$ for $T_1$ such
that $\vert \mu_1 \vert <1$; otherwise, all of the solutions to Hill's
equation would be bounded, contrary to the choice of $\lambda$. Combining these results, we obtain a
non zero meromorphic function $\psi$ that satisfies Hill's equation and 
$$\psi (z+2)=\mu_1\psi (z); \quad \psi (z+\omega_2) =\mu_2\psi (z); \quad \vert
\mu_1\vert <1 \qquad (z\in {\bf C}).\eqno(8.15)$$
Observe that
$$K(z,w) ={{\theta_1(z+w)e^{c(z+w)}}\over{\theta_1(z)^2\theta_1(w)^2}}
\Bigl({{\psi(z)\psi'(w)-\psi'(z)\psi(w)}\over{\psi (z+w)(\wp
(z)-\wp (w))}}\Bigr) \psi (z+w)\eqno(8.16)$$
\noindent where the first factor is elliptic of the third kind in $z$ and $w$ , and 
independent of $q$; the factor in parentheses is elliptic of the first kind in $z$ and $w$; the
final factor is elliptic of the second kind and gives a kernel in Hankel
 form. Also, we have $K(z+2,w)=\mu_1K(z,w)=K(z,w+2)$. \par
\indent  Next we check that $K(x,y)$ is bounded for $x,y\in [0,2]$,
 despite the first factor in the denominator which is problematic since $\theta_1(0)=0$. Let 
$$\Delta =16(e_1-e_2)^2(e_2-e_3)^2(e_1-e_3)^2\eqno(8.17)$$
\noindent be the
discriminant, and then recall that $\Delta=16\pi^4\theta_1'(0)^8$, so
$$K(z,w)\rightarrow {{2\pi
(\psi(0)\psi'(w)-\psi'(0)\psi (w))}\over{\Delta^{1/4}\theta_1(w)}} \qquad
(z\rightarrow 0);\eqno(8.18)$$
\noindent we define this right-hand side to be $K(0,w)$, and note further
that
$$K(0,w)\rightarrow {{(2\pi)^{3/2}(\psi
(0)\psi''(0)-\psi'(0)^2)}\over{\Delta^{3/8}}}\qquad (w\rightarrow 0).\eqno(8.19)$$
\noindent Hence $K(x,y)$ is bounded and continuous on $[0,2]\times
[0,2]$; so $K(x,y)$
defines a Hilbert--Schmidt operator on $L^2(0, \infty )$, and
$K(x,y)\rightarrow 0$ as $x\rightarrow\infty$ or $y\rightarrow\infty$.\par 
\indent Jacobi's theta function and $\wp$ function are related by 
$$\wp (x)=-\bigl( \log\theta_1(x)\bigr)''+e_1+\bigl(
\log\theta_1\bigr)''(1/2),\eqno(8.20)$$
\noindent while the addition rule for the cubic leads to 
$${{\wp'(x)-\wp'(y)}\over{\wp (x)-\wp (y)}}+2\int_{1/2}^{x+y}\wp (w)dw  
-2\int_{1/2}^{x}\wp (w)dw- 2\int_{1/2}^{y}\wp (w)dw-\int_{i\omega/2}^{i\omega/2+1}\wp
 (w)dw=0;\eqno(8.21)$$
\noindent see []. The function $\zeta (z)=-\int_{1/2}^z \wp (x)\, dx$ is additively
quasi-periodic with periods $1$ and $\omega_2$, so there exists $\alpha$ such
that $\zeta (x+i\varepsilon
)-\alpha x$ is infinitely differentiable and $1$-periodic on ${\bf R}$.
Hence there exists a linear system $(-A_\varepsilon , B_\varepsilon , C_\varepsilon
)$ such that $\zeta (x+i\varepsilon )=C_\varepsilon e^{-xA_\varepsilon}B_\varepsilon
$, and $\Vert e^{-xA_\varepsilon}\Vert \leq M(1+\vert x\vert )$ for some $M$ and all $x\in
{\bf R}$.\par
\indent There exist complex rational functions $u_0$ and $u_1$ such that
$$q(x)=u_0(\wp (x))+u_1(\wp (x))\wp'(x),\eqno(8.22)$$
\noindent so we can write
$$\eqalignno{{{q(x)-q(y)}\over{\wp (x)-\wp (y)}}&={{u_0(\wp (x))-u_0(\wp (y))}\over{\wp
(x)-\wp (y)}}+u_1(\wp (x)){{\wp'(x)-\wp' (y)}\over{\wp (x)-\wp (y)}}
+{{u_1(\wp
(x))-u_1(\wp (y))}\over{\wp (x)-\wp (y)}}\wp' (y)\cr
&=\sum_{j=1}^N v_j(x)w_j(y)+2u_1(\wp (x))\Bigl(-\zeta (x+y)+\zeta (x)+\zeta (y)\Bigr) 
&(8.23)\cr} $$
\noindent for some elliptic functions $v_j(x)$ and $w_j(y)$. 
Consider the derivatives
$$\eqalignno{ \Bigl( {{\partial}\over{\partial x}}+
&{{\partial}\over{\partial y}}\Bigr)K(x,y)&(8.24)\cr
&=\Bigl(
{{2\theta_1'(x+y)}\over{\theta_1(x+y)}}-{{2\theta_1'(x)}\over{\theta_1(x)}}-{{2\theta_1'(y)}
\over{\theta_1(y)}}+2c\Bigr) K(x,y)-\Bigl({{\wp'(x)-\wp'(y)}\over{\wp (x)-\wp (y)}}\Bigr) K(x,y)\cr
&\quad +{{\theta_1(x+y)e^{c(x+y)}}\over{\theta_1(x)^2\theta_1(y)^2}}
\Bigl({{q(y))-q(x))}\over{\wp (x)-\wp (y)}}\Bigr) \psi(x)\psi(y)\cr
&= {{\theta_1(x+y)e^{c(x+y)}}\over{\theta_1(x)^2\theta_1(y)^2}}
\Bigl( \sum_{j=1}^N v_j(x)w_j(y)+2u_1(\wp (x))\bigl(-\zeta (x+y)+\zeta (x)+\zeta (y)\bigr)\Bigr)\psi(x)\psi(y),\cr}$$
\noindent for some elliptic functions $v_j$ and $w_j$ for $j=1, \dots ,
N$. 
\indent By Proposition 7.2, there exists a periodic linear system $(-A_1,
B_1,C_1;E_1)$ such 
that $\theta_1(x)=C_1e^{-xA_1}B_1$.
 Then we can introduce $F_\varepsilon\in
L^2((0, \infty ); H)$ and $G_\varepsilon\in L^2((0, \infty ); H')$ such that 
$$\Bigl( {{\partial}\over{\partial x}}+
{{\partial}\over{\partial y}}\Bigr)K(x+i\varepsilon ,y+i\varepsilon )
=-\langle F_\varepsilon (x),G_\varepsilon (y)\rangle_H;\eqno(8.25)$$
\noindent hence
$$K(x+i\varepsilon ,y+i\varepsilon )=\int_0^\infty \langle F_\varepsilon (x+s),G_\varepsilon (s+y)\rangle_H \, 
ds.$$  
\indent (ii) In this case we use the double pole of $\wp (x)$ to cancel the
double zero of $\theta_1 (x)$ at $x=0$, and we do not have the
additional complication of the $\zeta$ function. We can write
$${{\theta_1(x+y)}\over{\theta_1(x)^2\theta_1(y)^2}}{{u_0(\wp
(x))-u_0(\wp (y))}\over{\wp (x)-\wp
(y)}}={{-\theta_1(x+y)}\over{\theta_1(x)^2\wp (x)\theta_1(y)^2\wp
(y)}}{{u_0(\wp (x))-u_0(\wp (y))}\over{1/\wp (x)-1/\wp (y)}},\eqno(8.26)$$
\noindent where $u_0(\wp (x))$ is bounded on ${\bf R}$, so there exist
$a_j\neq 0$ and integers $n_j,N\geq 0$ such that 
$$u_0(\wp (x))=\sum_j {{c_j}\over{(1/\wp
(x)-a_j)^{n_j}}}+\sum_{k=0}^N{{d_k}\over{\wp (x)^k}}\eqno(8.27)$$
\noindent for some constants $c_j$ and $d_k$. Hence   
$${{u_0(\wp (x))-u_0(\wp (y))}\over{1/\wp (x)-1/\wp (y)}}=\sum_j
v_j(\wp (x))w_j(\wp (y))\eqno(8.28)$$
\noindent where $v_j$ and $w_j$ are rational functions such that
$v_j(\wp (x))$ and $w_j(\wp (y))$ are bounded on $[0,2]$. We can
conclude the proof as in (i).\par 
\indent Let ${\cal D}(A_0)=\{ f\in L^2((0,\infty ); H): f'\in
L^2((0,\infty );H)\}$ and then introduce
$$\eqalignno{ C_0:{\cal D}(A_0)\rightarrow H: &\quad U(x)\mapsto U(0);\cr
A_0:{\cal D}(A_0)\rightarrow L^2((0, \infty );H):&\quad f(x)\mapsto
-f'(x);\cr
B_1:{\bf C}\rightarrow {\cal D}(A_0):&\quad  b\mapsto F(x)b;\cr
B_2:{\bf C}\rightarrow {\cal D}(A_0):&\quad b\mapsto G(x)b;&(8.29)\cr}$$
\noindent and then let
$$(-A, B,C)=\Bigl( -\left[\matrix{A_0&0\cr 0&A_0\cr}\right]
,\left[\matrix{B_1&0\cr 0&B_2\cr}\right],\left[\matrix{0&C_0\cr
C_0&0\cr}\right]\Bigr).\eqno(8.30)$$
\noindent Then we have
$$\left[\matrix{0&\Gamma_F\cr \Gamma_G&0\cr}\right]=\int_0^\infty
e^{-tA}BCe^{-tA}\, dt=R_0,\eqno(8.31)$$     
\noindent so $\det (I-K)=\det (I-R_0)$.\par
\indent (iii) We introduce the elliptic nome ${\hbox{q}}=e^{i\omega_2}$, and observe that
${\hbox{q}}\rightarrow 0$. Then from the produce formula for $\theta_1$, we have $\theta_1(x\mid
{\hbox{q}})\rightarrow 2\sin \pi x$ and $\wp (x\mid {\hbox{q}})\rightarrow
\pi^2{\hbox{cosec}}^2x-2\pi^2/3$, so
after some trigonometric reduction we find 
$${{\theta_1 (x+y\mid {\hbox{q}})}\over{ \theta_1 (x\mid {\hbox{q}})^2\theta_1 (y\mid
{\hbox{q}})^2(\wp (x\mid {\hbox{q}})-\wp (y\mid {\hbox{q}}))}}\rightarrow {{1}\over{16\pi^2 \sin \pi
(y-x)}},\eqno(8.32)$$
\noindent which implies the stated result.\par

\rightline{$\square$}\par
\vskip.05in
\noindent {\bf Corollary 8.5} {\sl Let} $\psi (x)={\hbox{sech}}^2x$
{\sl and } $\phi (x)={\hbox{sech}}^3x$.
{\sl Then the operator with kernel }
$${{\psi' (x)\psi (y)-\psi (x)\psi'(y)}\over{\sinh
(x-y)}}=\int_0^\infty \phi (x+t)\phi (t+y)\, dt\eqno(8.33)$$
\noindent {\sl defines a trace class integrable operator on $L^2(0,
\infty )$, which is positive.}\par
\vskip.05in
\noindent {\bf Proof.} One can differentiate both sides with respect to $\partial/\partial x+\partial /\partial y$, and then proceed as in [60]. Alternatively, one can deduce the result from a degenerate case of Theorem 7.4. Lam\'e's equation with $\ell=2 $ and
$\lambda^2=3g_2$ has solution $\wp (x\mid g_2, g_3)$. By taking
$i\omega \rightarrow -\infty$, and replacing $x$ by $i(x+\pi /2)$ 
we obtain $\psi (x)$.\par
\rightline{$\square$}\par

\noindent {\bf Remarks 8.7.}  (i)  Whereas Proposition 7.2  applies
to Lam\'e's equation, Theorem 7.4 does not since $\wp$ has a pole
at $x=0$.\par
\
\vskip.1in

\noindent {\bf 9 Tracy--Widom kernels and semi-classical weights}\par 

\vskip.05in   
\indent In [19], Deift, Its and Zhou observe that the correlation functions 
for many exactly solvable models admit determinant representations involving
 integrable operators, although only a few of them are of convolution operators. 
In random matrix theory, one often encounters integrable kernels
that are the products of Hankel integral operators on $L^2(0,
\infty )$; see [58, 59] and section 6 below for examples. 
The most significant examples in random matrix theory are associated 
with linear systems of differential equations with rational matrix 
coefficients, as in Proposition 6.4 and Remark 6.5. In
this section we realise such examples from a more general context.
 we introduce $(2,2)$ admissible linear systems with vectorial input and output spaces,
 so that we can
realise integrable operators as products of Hankel operators with matrix symbols. Then
 we consider the properties of the tau functions that arise from these 
Hankel operators, and derive Schlesinger's system of differential equations. \par
\indent Let ${\cal R}$ be the ring of bounded linear operators on
$L^2(a,b)$ that have continuous kernel, and let ${\cal F}$ the
ideal in ${\cal R}$ of finite rank operators; let $\partial_M$ be 
the derivation $\partial_M (K)=KM-MK$, where $Mh(x)=xh(x)$ for $h\in L^2(a, b)$.
Now ${\cal R}_{\cal F}$ is a ring which contains
 a subring of operators known as the integrable operators.
\vskip.05in
\noindent {\bf Definition} (Integrable operators) [19]  
An integrable kernel has
the form 
$$K(x,y)={{\sum_{j=1}^n f_j(x)g_j(y)}\over{x-y}},\eqno(9.1)$$ 
where $f_j,g_j$ are continuous and bounded functions on $(0, \infty )$,
and we suppose
further that $\sum_{j=1}^n f_j(x)g_j(x)=0,$ so $K$ is nonsingular on 
$x=y$.\par
\vskip.05in
\noindent {\bf Definition}(Basic equation) Let the basic equation be the differential equations
$$J{{d}\over{dx}}\left[\matrix{f\cr g\cr}\right]=\Omega (x)\left[\matrix{f\cr
g\cr}\right],\qquad\Omega (x)= \left[\matrix{\gamma &\alpha \cr \alpha
&\beta \cr}\right],\qquad J=\left[\matrix{0&-1\cr
1&0\cr}\right],\eqno(9.2)$$
\noindent with $\alpha, \beta$ and $\gamma$ rational functions. Let
${\bf C}(x)$ be the differential field of rational functions,
let $W_F$ be the fundamental solution matrix of (9.2) on some neighbourhood
of $x_0$, and let ${\bf K}_{PV}$ be Picard--Vessiot field, namely the field of fractions generated
by ${\bf C}(x)$ and the entries of $W_F$. In the following calculations, a
crucial feature of the differential equation is that the entries of $(\Omega
(x)-\Omega (y))/(x-y)$ all belong to ${\bf C}(x)\otimes {\bf C}(y).$\par

\indent As in Tracy and Widom's theory of matrix models, we introduce the kernel 
$$K_{(z)}(x,y)={{f(x+2z)g(y+2z)-f(y+2z)g(x+2z)}\over{x-y}},\eqno(9.3)$$
\noindent and $\tau (x;\lambda )=\det (I+\lambda K_{(x)})$.
Further, suppose that $L_{(z)}$ by $(I-L_{(z)})(I+K_{(z)})=I$; see [58, 59]. \par
\vskip.05in
\noindent {\bf Theorem 8.1} {\sl Suppose that $\alpha, \beta$ and $\gamma$
are proper rational functions with $n$ poles of order less than or equal
to $p$, for some $p\in {\bf N}$, and all poles are in ${\bf
C}\setminus [0, \infty )$; suppose that
$f,g\in L^2(0, \infty )$ are solutions of (6.2) and that $f(x),g(x)\rightarrow 0$ as
$x\rightarrow\infty$.\par
\indent (i) Then there exist Hilbert--Schmidt Hankel operators $\Gamma_\Phi$
and $\Gamma_\Psi$ with $2np^2\times 2np^2$ matrix symbols $\Phi$ and $\Psi$ such that} 
$$\det (I+\lambda K_{(z)})=\det (I+\lambda
\Gamma_{\Phi_{(z)}}\Gamma_{\Psi_{(z)}}
).\eqno(9.4)$$
\noindent {\sl The scattering functions are given by
finite matrices with entries in the differential field ${\bf
C}(x,f,g)$.}\par
\indent {\sl (ii) There exists $x_0$ such that $L_{(z)}$ is a bounded
integrable operator for all $z\geq x_0$.}\par  
\indent {\sl (iii) Suppose further that $e^{2\varepsilon x}
f(x)\rightarrow 0$ and $e^{2\varepsilon x}g(x)\rightarrow 0$ as
$x\rightarrow \infty$ for some $\varepsilon >0$. Then $\Phi$ and $\Psi$ are realised by $(2,2)$
admissible linear systems.}\par
\vskip.05in
\noindent {\bf Proof.} (i) We can write
$$\Omega
(x)=E_0+\sum_{k=1}^n\sum_{\ell
=1}^{p_k}{{E_{k,\ell}}\over{(x-a_k)^{\ell}}},\eqno(9.5)$$
\noindent where the $E_0$ and $E_{k, \ell}$ for $\ell =1, \dots ,p_k$ with $p_k\leq p$ and $k=1, \dots , n$ are symmetric 
$2\times 2$ matrices and the poles $a_j$ lie in ${\bf C}\setminus [0, \infty )$. From the differential equation, we have 
$$\eqalignno{ \Bigl({{\partial}\over{\partial
x}}+{{\partial}\over{\partial y}}\Bigr){{f(x)g(y)-f(y)g(x)}\over{x-y}}&=
\Bigl\langle {{\Omega (x)-\Omega (y)}\over{x-y}}\left[\matrix{f(x)\cr
g(x)\cr}\right], \left[\matrix{f(y)\cr g(y)\cr}\right]\Bigr\rangle&(9.6)\cr
&=-\sum_{k=1}^n \sum_{\ell =1}^{p_k}\sum_{\nu=0}^{\ell}\Bigl\langle
{{E_{k,\ell}}\over{(x-a_k)^{\ell-\nu}}}\left[\matrix{f(x)\cr
g(x)\cr}\right], {{1}\over{(y-a_k)^{\nu +1}}}\left[\matrix{f(y)\cr g(y)\cr}\right]
\Bigr\rangle ,\cr}$$
\noindent where we have used the real inner product. Noting that
$E_{k, \ell}$ has rank less than or equal to two, let $N=2np^2$ and introduce scalar-valued functions
$\phi_j(x)$ and $\psi_j(y)$ such
that the previous sum  equals $-\sum_{j=1}^{N} \phi_j(x)\psi_j (y)$, and since the poles are off
$(0, \infty )$, we can ensure that $\int_0^\infty x(\vert\phi_j
(x)\vert^2+\vert\psi_j(x)\vert^2)dx$ is finite, so $\phi_j$ and
$\psi_j$ give the symbols of Hilbert--Schmidt Hankel operators on $L^2(0, \infty )$. Then one verifies the
identity
$${{f(x)g(y)-f(y)g(x)}\over{x-y}}
=\int_0^\infty \sum_{j=1}^{N}\phi_j(x+s)\psi_j(s+y)\, ds;\eqno(9.7)$$
\noindent indeed by the preceding calculation, the difference between the two sides of (6.7) is a 
function of $x+y$, which goes to zero as $x\rightarrow \infty $ or $y\rightarrow \infty$ in any way. and hence vanishes identically. 
Finally, we build the
$N\times N$ matrices
$$\Phi (x)=\left[\matrix{\phi_1(x)&\phi_2(x)&\dots &\phi_{N}(x)\cr 0&0&\dots  &0\cr
\vdots&\vdots &\ddots &\vdots\cr 0&0&\dots &0\cr}\right], 
\Psi (y)=\left[\matrix{\psi_1(y)&0&\dots &0\cr
                          \psi_2(y)&0&\dots  &0\cr
                          \vdots&\vdots&\ddots &\vdots\cr
                            \psi_{N}(y)&0&\dots
&0\cr}\right]\eqno(9.8)$$
\noindent so that $\Gamma_\Phi$ and $\Gamma_\Psi$ are Hilbert--Schmidt
matrix operators, 
and with $\phi_{j,(z)}(x)=\phi_j(x+2z)$ {\it etc} we have
$$\det \bigl(I+\lambda K_{(z)}\bigr)=\det \Bigl(I+\lambda \sum_{j=1}^{N} 
\Gamma_{\phi_{j,(z)}}\Gamma_{\psi_{j,(z)}}\Bigr)=\det
\bigl(I+\lambda \Gamma_{\Phi_{(z)}}\Gamma_{\Psi_{(z)}}\bigr).\eqno(9.9)$$
Hence
$$\Xi (x)=\left[\matrix{0&\lambda \Phi (x)\cr -\Psi (x)&0\cr}\right]$$
\noindent has $\det (I+\lambda K_{(z)})=\det (I+\Gamma_{\Xi_{(z)}}).$\par
\indent We can assume that the solutions $f$ and $g$ are meromorphic on
some nonempty open subset of ${\bf C}$, so from the linear base
equation we deduce that ${\bf C}(x,f,g)$ is a differential field. The
definition of $\phi_j$ and $\psi_j$ shows that these functions also
belong to ${\bf C}(x,f,g)$.\par  
\indent (ii) We can define $L_{(z)}=K_{(z)}(I+K_{(z)})^{-1}$ for
all $z$ such that $\Vert K_{(z)}\Vert <1$.  Deift et al [] have shown that, when $I+K$
is an invertible integrable operator, $L=(I+K)^{-1}-I$ is also integrable, and 
for any derivation $\partial :{\cal R}\rightarrow {\cal R}$ 
$$\partial L=(I+K)^{-1}(\partial K)(I+K)^{-1}.\eqno(9.10)$$
\noindent In short, we obtain
$L$ from the kernel
$${{F(x)G(y)-F(y)G(x)}\over{x-y}},\qquad {\hbox{where}}\qquad
\left[\matrix{F\cr G\cr}\right]=\left[\matrix{(I+K)^{-1}f\cr 
(I+K)^{-1}g\cr}\right].\eqno(9.11)$$ 
\noindent Moreover, $\partial K=[d/dx, K]$ is the finite rank
integral operator that is represented by the kernel (8.6), so $\partial
L$ is also finite rank.\par
\indent (iii) Given that $f$ and $g$ are of exponential decay, the
integral $\int_0^\infty xe^{2\varepsilon x}\vert\phi_j (x)\vert^2dx$
converges, and hence the Hankel operator $\Gamma_j$ with symbol $e^{\varepsilon
x}\phi_j(x)$ is bounded. We decompose $\phi_j=\Re\phi_j+i\Im \phi_j$
so that we can work with the self-adjoint Hankel operators 
$\Gamma_{\Re \phi_j}$ and $\Gamma_{\Im \phi_j}$; so by theorem 2.1 of 
[44, p.257], there exist linear systems $(-A'_j, B'_j, C'_j)$ and $(-A_j'',
B_j'', C_j'')$ with input and output
spaces ${\bf C}$, and state space $H$, and all operators bounded, such that $e^{\varepsilon x}\Re \phi_j(x)=C'_je^{-xA'_j}B'_j$
and $e^{\varepsilon x}\Im \phi_j(x)=C_j''e^{-xA_j''}B_j''$; then we introduce the linear system
$$\bigl(-A_j,B_j, C_j\bigr)=\Bigl(-\left[\matrix{ A_j'&0\cr
0&A_j''\cr}\right], \left[\matrix{ B_j'\cr
B_j''\cr}\right], \left[\matrix{
C_j'&iC_j''\cr}\right]\Bigr),\eqno(9.12)$$
\noindent which has scattering function $e^{\varepsilon x}\phi_j(x)= C_je^{-xA_j}B_j$. Hence
we can introduce the linear system
$$\bigl(-A,B,C\bigr)=\Bigl(-\left[\matrix{\varepsilon I+A_1&\dots &0\cr 0&\ddots&\vdots \cr
0&\dots & \varepsilon I+A_N\cr}\right],  
\left[\matrix{B_1&\dots &0\cr 
0&\ddots&\vdots \cr
0&\dots & B_N\cr}\right], \left[\matrix{C_1&\dots &C_N\cr 0&\ddots&\vdots \cr
0&\dots &0\cr}\right]\Bigr)\eqno(9.13)$$
\noindent where $A:H^{2N}\rightarrow  H^{2N}$, $B:{\bf C}^N\rightarrow 
H^{2N}$ and $C:H^{2N}\rightarrow
{\bf C}^N$ are bounded linear operators. Since $\Re \langle A\xi, \xi \rangle_{H^N} \geq \varepsilon
\langle \xi, \xi\rangle_{H^N}$ for all $\xi\in H^N$, Lemma 2.1 and Proposition 2.2 show that $(-A,B,C)$ is a
$(2,2)$ admissible linear system. Evidently $(-A,B,C)$ realises 
$\Phi$, and we can likewise
realise $\Psi$ by a $(2,2)$ admissible linear system.\par 
\rightline{$\square$} 
 
\vskip.05in
\indent By taking $\alpha =0$, $\gamma $ to be a negative proper rational
function and $1/\beta $ to be a positive polynomial on $(0, \infty )$, one
can produce solutions of (6.2) that satisfy the hypotheses of Theorem
6.1(ii).
\par

\noindent {\bf Remarks} (i) Whereas Theorem 6.1 does not give an
explicit form for the admissible linear system $(-A,B,C)$, we can
produce one explicitly in several important cases; 
see  [9, 10, 11]. \par
\noindent (ii) An important special case is the system
$$\left[\matrix{0&-1\cr 1&0\cr}\right]{{d}\over{dx}}\left[\matrix{f\cr g\cr}\right]+{{1}\over{2}}
\left[\matrix{p+q&i(p-q)\cr i(p-q)&-(p+q)\cr}\right]
\left[\matrix{f\cr g\cr}\right]=E\left[\matrix{f\cr g\cr}\right]\eqno(9.14)$$
\noindent for $p,q\in {\bf C}(x)$. When $p=1$, this reduces to Schr\"odinger's
equation $-g''+qg=E^2g$.\par
\vskip.05in
\noindent {\bf Theorem 9.2.} {\sl  For $W_j\in M_2({\bf C})$, let  $W_\infty=-\sum_{j=0}^{N}W_j$ and suppose that $W_\infty$ has eigenvalues
 $\pm\theta_\infty/2$, where $\theta_\infty$ is not a positive integer. 
Then the differential equation 
$${{d}\over {dx}}\Phi =\sum_{j=0}^N {{W_j(t)}\over{x-t_j}}\Phi\eqno(9.15)$$
\noindent has general  solution $\Phi (x)$ such that $\Phi (x)x^{W_\infty}$ is holomorphic on a neighbourhood $\{ x\in {\bf C}: \vert x\vert >\kappa_0\}$; 
the Laurent coefficients $C_j(t)$ and $W_\infty$ are determined by a recurrence relation and determine a linear system $(-A,B,C(t))$ with state space 
$L^2((0, \infty ); \ell^2)$ and scattering function $\Phi (x).$ The operator $R_x$ is given by an infinite matrix kernel}
$$R_x\leftrightarrow \Bigl[ \kappa_0^{j-k}\Gamma 
(jI+W_\infty )^{-1}{{e^{-x(u+v)}v^{j-1}v^{W_\infty}C_k}\over{u+v}}\Bigr]_{j,k=0}^\infty
\qquad (u,v>0; x>\kappa_0).\eqno(9.16)$$  
\vskip.05in
\noindent {\bf Proof.} Then for any $\Phi_0\in M_2({\bf C})$ there exists a solution of the form 
$$\Phi (x)=\Bigl( I+\sum_{\ell =1}^\infty {{C_\ell (t)}\over{x^\ell}}\Bigr)x^{-W_\infty}
\Phi_0\eqno(9.17)$$
\noindent where $C_k(t)$ are $2\times 2$ matrices for $k=1,2, \dots$ and $C_0=I_2$.
 Indeed, the matrices $W_\infty $ and $W_\infty+kI$ have no eigenvalues in common,
 and hence by a theorem of Sylvester,
the recurrence relation
$$W_\infty C_k(t)-C_k(t)(W_\infty +kI)=\sum_{j=0}^N
 \sum_{\nu=1}^kt_j^\nu W_jC_{k-\nu }(t)\eqno(9.18)$$
\noindent has a unique solution for $C_k (t)$ given $C_0, C_1(t), \dots ,C_{k-1}(t)$; 
more explicitly, we can write
$D_{k-1}=\sum_{j=0}^N \sum_{\nu=1}^k t_j^kW_jC_{k-\nu}$ and then choose
$$C_k=-\int_0^\infty e^{sW_\infty}D_{k-1}e^{-sW_\infty} e^{-sk}ds,\eqno(9.19)$$
\noindent which evidently satisfies the relation. We can select positive constants such that
 $T\geq \max\{ 1, \vert t_j\vert; j=0, \dots ,N\}$, $M_0\geq \Vert W_j(t)\Vert$ ,
 $\Vert e^{-sW_\infty}\Vert \leq M_\infty e^{-s\theta_\infty}$ and 
$M_2>\max\{1,  k/\vert k-\theta_\infty\vert ;k=1,2,\dots\}.$ Then from the integral we have
$$\Vert C_k\Vert\leq {{M_\infty^2\Vert D_{k-1}\Vert}\over{\vert k-\theta_\infty \vert}},
\eqno(8.20)$$
\noindent where $\Vert D_k\Vert \leq \sum_{j=0}^N \sum_{\nu=1}^kT^\nu M_0\Vert C_{k-\nu}
\Vert$, hence
$$\Vert C_k\Vert\leq {{M_\infty^2M_1(N+1)}\over{\vert k-\theta_\infty\vert}}
\sum_{\mu=0}^{k-1}T^{k-\mu}\Vert C_\mu\Vert,\eqno(9.21)$$
\noindent and we can deduce by induction that 
$$\Vert C_k\Vert\leq (M_0M_2M_\infty^2(N+1)T)^k\Vert C_0\Vert .\eqno(9.22)$$
 Hence  with $\kappa_0=M_0M_2M_\infty^2(N+1)T,$  the series () converges for all $x\in {\bf C}$
 such that $\vert x\vert >\kappa_0.$ 

Now let $H=L^2((0,\infty );\ell^2)$, and ${\cal D}(A)=\{ f\in H: sf(s)\in H\}$; then we introduce linear operators by 
$$b_j(s)=\Gamma (jI+W_\infty )^{-1}s^{(j-1)I+W_\infty}\qquad (j=0,1,\dots ),\eqno(9.23)$$
\noindent then introduce the linear maps
$$\eqalignno{ A:{\cal D}(A)\rightarrow H:& f(s)\mapsto sf(s);\cr
e^{-xA}B: M_2({\bf C})\rightarrow H:& \beta \mapsto 
\bigl(\kappa_0^jb_j(s)e^{-xs}\beta\bigr)_{j=0}^\infty;\cr
C(t):H\rightarrow M_2({\bf C}):& (f_j)_{j=0}^\infty \mapsto
 \sum_{j=0}^\infty \int_0^\infty \kappa_0^{-j} C_j f_j(s)ds.&(9.24)\cr}$$
\noindent Then $(e^{-sA})_{s>0}$ is a strongly continuous and bounded semigroup on $H$, and $e^{-xA}B:M_2({\bf C})\rightarrow H$ is bounded for $x>\kappa_0$ and $\Phi (x)=C(t)e^{-xA}B$ for $x>\kappa_0$.\par
\rightline{$\square$}\par


\indent In Theorem 9.1, we showed that $\tau (z,t)$ is given by the Fredholm
determinant of a product of Hankel operators, and in Proposition 2.4, we expressed ${{\partial}\over{\partial z}}\log\det
(I+\Gamma_{\phi_{1, (z)}}\Gamma_{\phi_{2, (z)}})$ in terms of the
solution of a Gelfand--Levitan equation; thus 
${{\partial }\over{\partial z}}\log \tau (z,\lambda )$ is given
in terms of the solution of a Gelfand--Levitan equation.\par

\indent The tau function of $\Sigma (t)=(-A, B, C(t))$ may be expressed in terms of $\tau (x; t_0,
\dots t_N)$, where the $t_j$ are treated as parameters. Let ${\bf K}_0$ be the space of meromorphic complex function of $t=(t_0, \dots , t_N)$ for $t\in \Omega$, where $\Omega$ is a domain in ${\bf C}^{N+1}$. We now compare this 
with the tau functions introduced by Jimbo, Miwa and Ueno to describe the isomonodromy 
of rational differential equations. In the particular case of Theorem 8.2, their tau function is
locally meromorphic in the $t_j$ and satisfies
$$\log \tau (t_0, \dots ,t_N)={{1}\over{2}}\sum_{j,k:j\neq k} {\hbox{trace}}
(W_jW_k)d\log (t_j-t_k).\eqno(9.25)$$
We can express this more explicitly via determinants in a particular case.\par

Let $\tau_0=1$ and let $(\tau_n)_{n=1}^\infty$ be a sequence of
positive terms; let $\tau_{-n}=\tau_n$, and then let
$a_n=\tau_{n+1}\tau_{n-1}/\tau_n^2$ for $n\in {\bf Z}$. For any
complex sequence $(b_n)$, we introduce 
$$V_n(x)=\left[\matrix{0&1\cr -a_n&x+b_n\cr}\right],\qquad
Y_{-1}=\left[\matrix{0\cr 1\cr}\right],\qquad   
Y_{0}=\left[\matrix{1\cr x+b_1\cr}\right],\eqno(9.26)$$
\noindent so that the recurrence relation 
$$Y_{n+1}(x)=V_n(x)Y_n(x)\qquad (n=0, 1, \dots )\eqno(9.27)$$
\noindent gives $Y_n(x)=\left[\matrix{p_n(x)\cr p_{n+1}(x)\cr}\right]$, where
$p_n(x)$ is a monic polynomial of degree $n$. By a Theorem of Favard,
given $(\tau_n)_{n=1}^\infty$, there exists a real sequence $(b_n)$
and a bounded and positive measure $\omega$ on $(-\infty ,\infty )$
such that $(p_n(x))_{n=1}^\infty$ is orthogonal with respect to the
usual inner product on $L^2(\omega )$, and 
$$\tau_n=\det \Bigl[ \int_{-\infty}^\infty x^{j+k}\omega
(dx)\Bigr]_{j,k=0}^n\qquad (n=1, 2, \dots ).\eqno(9.28)$$
\vskip.05in 
\noindent {\bf Definition} (Semi-classical weights)  Let $S$ be a finite union of real intervals
and $w$ a real functions such that\par
\indent  (i) $w$ is positive and Lebesgue
integrable on $S$; \par
\indent (ii) $w$ is differentiable on $S$, and $w(x)\rightarrow 0$
$x\in S$ tends to an endpoint of $S$, including infinity in the case of
$S$ unbounded;\par
\indent (iii) $w'(x)/w(x)$ is rational.\par
\noindent Then $w$ is said to be a semi classical weight.\par

\indent Let $N=2g+1$ and for $t=(t_0, \dots ,t_{N})$ where 
$t_0<t_1<\dots <t_N$, let $S= \cup(t_{2j},
t_{2j+1})$, where the union is taken over some subset of $\{
0,\dots ,g\}$. Then we introduce the nonnegative weight 
$$w_t(x)=(t_N-x)^{1/2}\prod_{j=0}^{g}\bigl((t_{2j+1}-x)(x-t_{2j})\bigr)^{1/2}{\bf I}_{S}
(x).\eqno(9.29)$$
\noindent We introduce the orthogonal projection $K(n,t):L^2(S,
w_t)\rightarrow {\hbox{span}}\{ 1, x, \dots , x^{n-1}\}$, which by the
Christoffel--Darboux formula may be expressed as the integral operator
on $L^2(S, w_t)$ with kernel
$$K(y,z;n,t)=-{{\gamma_{n-1}(t)^2}}{{\langle JY_{n-1}(y,t),
Y_{n-1}(z,t)\rangle}\over{y-z}}\eqno(9.30)$$
\noindent where $\gamma_{n-1}(t)^{-2}=\int p_{n-1}(x,t)^2w_t(x)dx$. For
semi classical weights as in (), the $Y_n(x,t)$ satisfy a linear
differential equation with rational coefficients, by results of Magnus
or Tracy and Widom. The crucial point is that
$${{w'(x)/w(x)-w'(y)/w(y)}\over{x-y}}\eqno(9.31)$$
\noindent is an integrable operator of finite rank. 
There exist $2\times 2$ complex matrices $W_j(n,t)$
for $j=0, \dots ,N,\infty$ such that
$${{d}\over{dx}}Y_n(x,t)=\Bigl(W_\infty (n,t)+\sum_{j=0}^N
{{W_j(n,t)}\over{x-t_j}}\Bigr) Y_n(x,t).\eqno(9.32)$$
Next we let $K_x(n,t)$ be the compression of $K(n,t)$ to $L^2(S\cap (-\infty ,x];w_t)$ so that $K_x(n,t)$ has kernel 
$$K_x(y,z;n,t)={\bf I}_{(-\infty ,x]}(y) K(y,z;n,t){\bf I}_{(-\infty ,x]}(z),\eqno(9.33)$$
\noindent and the associated determinant
$$\tau_n(x,t)=\det \Bigl[ \int_{-\infty}^x t^{j+k}\omega_t
(dt)\Bigr]_{j,k=0}^n\qquad (n=1, 2, \dots ).\eqno(9.34)$$
\noindent We also introduce
$$L_j={{\partial}\over{\partial t_j}}+{{W_j(n,t)}\over{x-s_j}}\qquad (s_j=t_j).\eqno(9.35)$$
Schlesinger's equations are
$$\eqalignno{{{\partial}\over{\partial t_k}}W_j(n,t)&={{W_j(n,t)W_k(n,t)-W_kW_j(n,t)}\over{t_j-t_k}}\qquad (j\neq k; j,k=0, 
\dots ,N);&(8.36)\cr {{\partial}\over{\partial t_j}}W_j(n,t)&=-\sum_{k=0; k\neq j}^{N}{{W_j(n,t)W_k(n,t)-W_k(n,t)W_j(n,t)}\over{t_j-t_k}}\cr
&\quad -(W_j(n,t)W_\infty (n,t)-W_\infty (n,t)W_j(n,t))
\qquad (j=0, \dots ,N).&(9.37)\cr}$$
\noindent Let
$$A(\mu, x,t)=\mu I_2-\sum_{j=0}^N{{W_j(n,t)}\over{x-t_j}}-W_\infty
(n;t)\eqno(9.39)$$

\vskip.05in 

\noindent {\bf Proposition 9.3} {\sl For a semi classical weight $\omega_t$
with zeros $t=(t_0, \dots , t_N)$, the sequence $(\tau_n(\infty ,
t))_{n=0}^\infty$. is given by the Hankel determinants of the moments
of $\omega_t$ on $S$. The monic orthogonal polynomials $(p_n,
p_{n+1})$ satisfy a
rational differential equation with singular points at the $t_j$ and
corresponding $2\times 2$ residue matrices $W_j(n,t)$. Suppose that Schlesinger's equations hold.\par 
\indent (i) Then for all $x\in S$, the determinant $\det (I+\lambda
K_x(n,t))$ is a polynomial in $\lambda$ of degree $n$, with leading
coefficient $\tau_n(x,t)/\tau_n(\infty ,t)$.\par
\indent (ii) Then the $L_j$ commute and there is an exact differential form}
$$d\log \tau_n(\infty ,t)=\sum_{0\leq j\neq k\leq
N}{{{\hbox{trace}}(W_j(n,t)W_k(n,t))}\over{t_j-t_k}}dt_j\eqno(9.38)$$
\noindent {\sl so that $\tau_n(\infty ,t)$ is the Jimbo--Miwa--Ueno $\tau $
function for the algebraic family ().}\par
\indent {\sl (iii) The operators $(\partial/\partial t_j)K(x,y;n,t)$ are of finite rank, and the Laurent coefficients $C_k(t)$ of Theorem 8.2 satisfy a system of linear partial differential equations that is determined recursively by Schlesinger's equations.}\par
\indent {\sl (iv) The curves 
$${\cal E}_n=\bigl\{ (x,\mu ):\det A(\mu ,t,x)=0\bigr\}$$
\noindent are invariant under the flows generated by the $L_j$ for all $n$. Each simple zero of $\mu\mapsto \det A(\mu ,x,t)$ has an eigenvector $\Psi (x,t)$ such that 
the entries of $\Psi $ are meromorphic on ${\cal E}_n$ and there exists $\xi_j\in {\bf K}_0[\mu ]$ such that $L_j\Psi (x,t)=\xi_j(x,t)\Psi (x,t)$ for all $j=1, \dots ,n$.}

\vskip.05in
\noindent {\bf Proof.} (i) We have  Fredholm expansion
$$\det (I+\lambda K_x(n,t))=\sum_{j=0}^n
{{\lambda^j}\over{j!}}\int_{(-\infty ,x)^j}\det [K(y_k, y_\ell
;n,t)]_{j,\ell =1}^j w_t(y_1)\dots w_t(y_j)dy_1\dots dy_j\eqno(9.40)$$
\noindent in which the series terminates since $K_x$ has rank $n$.
Indeed, $(\gamma_k(t)p_k(y,t))_{k=0}^\infty$ is an orthonormal basis
for $L^2(w_t)$, so
$K(y,z;n,t)=\sum_{j=0}^{n-1}\gamma_k(t)^2p_k(y,t)p_k(z,t)$, hence
$$\det [K(y_k, y_\ell ;n,t)]_{k,\ell =1}^n =\det
[\gamma_{m-1}(t)p_{m-1}(y_k,t)]_{k,m=1}^n  
   \det [\gamma_{m-1}(t)p_{m-1}(y_\ell ,t)]_{\ell ,m=1}^n,\eqno(9.41)$$
\noindent so we can express the leading terms as 
$$\eqalignno{
{{1}\over{n!}}\prod_{m=0}^{n-1}\gamma_m(t)^2\int_{(-\infty ,x)^n}&\det
[y_k^{\ell-1}]_{k, \ell =1, \dots ,n}w_t(y_1)\dots w_t(y_n)dy_1\dots
dy_n\cr
&=\prod_{m=0}^{n-1}\gamma_m(t)^2\int_{(-\infty ,x)}\det [y^{k+\ell
-2}]_{k, \ell =1,\dots ,n}w_t(y)dy\cr
&={{\tau_n(x,t)}\over{\tau_n(\infty ,t)}},&(9.42)\cr}$$
\noindent where the final step follows from a standard formula.\par
\indent (ii) Schlesinger [51] observed that the system (8.26) is
consistent if and only if the family of solutions 
satisfies an isomonodromy
condition with respect to infinitesimal deformation, or
equivalently that the family of differential operators
commutes. When forming products $L_jL_k$, one evaluates at $t_j=s_j$ after differentiation, and
by hypothesis $L_jL_k-L_kL_j=0$. 
The variational equation for the pole $t_j$ is
$${{\partial}\over{\partial
t_j}}Y_n(x;t)={{-W_j(n;t)}\over{x-t_j}}Y_n(x;t).\eqno(9.43)$$

Then the poles $t_j$ are independent of $n$, and the residue
matrices satisfy the recurrence relation
$$ W_j(t; n+1)=V_n(a_j;t)W_j(t;n)V_n(a_j;t)^{-1}\qquad (n=1, 2, 
\dots ; t\in {\cal T})\eqno(9.44)$$
\noindent The recurrence relation () and the variational equation lead to the
consistency condition
$${{\partial}\over{\partial
t_j}}V_n(x;t)-{{V_n(x;t)W_j(n;t)}\over{x-t_j}}=-{{W_j(n+1;t)V_n(x;t)}
\over{x-t_j}}.\eqno(9.45)$$
Given (i), and these equations, we can deduce the conclusion from Chen and Its
paper [].\par
\indent (iii) From the recurrence relation, we see that $w_j={\hbox{trace}}W_j(t;n)/2$ does not depend upon $n$, while Schlesinger's equations imply that $w_j$ does not depend upon $t$. Note that $(W_j(n;t)-w_jI)^TJ=-J(W_j(n;t)-w_jI)$; then from the variational equation (9.33) and (9.28) we obtain
$${{\partial }\over{\partial t_j}}\Bigl({{(x-t_j)^{w_j}(y-t_j)^{w_j}}\over{\gamma_{n-1}(t)^{2}}}K(x,y;n,t)\Bigr)=-\Bigl\langle
{{J(W_j(n,t)-w_jI)Y_n(x;t)}\over{x-t_j}},
{{Y_n(y;t)}\over{y-t_j}}\Bigr\rangle,\eqno(9.46)$$
\noindent which gives an operator on $L^2(0, \infty )$ of finite rank. Since $K(x,y;n,t)$ has rank $n$, we deduce that $(\partial/\partial t_j)K(x,y;n,t)$ gives a finite rank operator.\par

\indent By partially differentiating (9.17) and (9.18) with respect to $t_\ell$, we obtain expressions for 
$$\eqalignno{W_\infty {{\partial C_k(t)}\over{\partial t_\ell}}-{{\partial C_k(t)}\over{\partial t_\ell}}(W_\infty+kI)&=-{{\partial W_\infty }\over{\partial t_\ell}}C_k(t)
+C_k(t){{\partial W_\infty}\over{\partial t_\ell}}+\sum_{\nu =1}^k \nu t_\ell^{\nu -1}W_\ell C_{j-\nu} (t)\cr
&\quad +\sum_{j=0}^N \sum_{\nu =1}^k t_j^\nu {{\partial W_j}\over{\partial t_\ell}} 
C_{k-\nu }(t) +\sum_{j=0}^N \sum_{\nu =1}^k t_j^\nu W_j
{{\partial C_{k-\nu}(t)}\over{\partial t_\ell}}&(9.47)\cr}$$ 
\noindent where the right-hand side $D_k(t)$ is a linear combination of $C_k(t)$, the $C_\nu$ and ${{\partial}\over{\partial t_\ell}}C_\nu (t)$ for $\nu =0, \dots ,k-1$. As in the proof of Theorem 8.2, $D_k(t)$  determines ${{\partial}\over{\partial t_\ell}}C_k(t)$ as 
$${{\partial C_k(t)}\over{\partial t_\ell}}=-\int_0^\infty e^{sW_\infty}
D_k(t)e^{-(kI+W_\infty )s}\, ds.\eqno(9.48)$$
\indent Hence we obtain the coefficients that arise in the formula
$${{\partial}\over {\partial t_\ell}}\log \tau_n (x;t)=
{\hbox{trace}}\int_x^\infty {{\partial C(t)}\over{\partial t_\ell}}
e^{-sA}(I+R_x)^{-1}e^{-sA}B\, ds.\eqno(9.49)$$ 
\indent (iv) 
Moreover, the $L_j$ form an algebraic family in the sense of
Krichever and Novikov, since 
$$\Bigl[ L_\ell, W_\infty (n,t)+\sum_{j=0}^N
{{W_j(n,t)}\over{x-s_j}}\Bigr]=0, \qquad (s_j=t_j; j,\ell =0, \dots , N).$$
The coefficients of the polynomial  $\mu\mapsto \det A(\mu, x,t)$ 
\noindent are invariant under the flows associated with the $L_j$. Hence the curve ${\cal E}_n$  is invariant under the flows generated by the $L_j$. \par
\indent Let ${\bf K}_0$ be the field of meromorphic functions in $t=(t_1, \dots ,t_n)$. Moreover, each 
$(x, \mu )$ on ${\cal E}_n$ is associated with a maximal ideal 
${\bf p}$ of ${\bf K}_0[\mu ]$; furthermore, ${\bf K}_0[\kappa ]/{\bf p}$ is 
a finite algebraic extension of ${\bf K}_0$ by the weak Nullstellensatz [5].\par 
\indent For each simple zero of $\mu \mapsto \det A(x,\mu ,t)$  We choose $(x, \mu )$ on ${\cal E}_n$ and then by Cramer's rule select 
$\Psi$ such that
$A(x,\kappa ,t)\Psi=0$, where the entries of  $\Psi$ are rational expressions in $\mu$ and the entries of $W_j(t)/(x-t_j)$ for $j=1, \dots ,n$ and $W_\infty (n,t)$. In 
particular, $t\mapsto \Psi_t (t, x )$ is meromorphic and the equation
$$\Bigl( W_\infty (n,t)+\sum_{k=0}^N {{W_k(t)}\over{x-t_k}}\Bigr)\Psi =\mu \Psi\eqno(9.50)$$
\noindent is also satisfied by $L_j\Psi$ for all $j$; hence 
by uniqueness of eigenvectors corresponding to simple eigenvalues, for each $(x,t)$ we obtain $\xi_{j}(x,t)\in {\bf C}$ such that $L_j\Psi
=\xi_{j}(x,t)\Psi$, and by considering a non-zero element of $\Psi$, we deduce that $\xi\in {\bf K}_0[\mu ]$.\par

\rightline{$\square$}\par
\vskip.05in

\vskip.05in
\vskip.05in


\noindent {\bf 10. The Kadomtsev--Petviashvili equations}\par

\indent In this section we describe the theta functions which may be
expressed in terms of elliptic linear systems as in section 6. Quasi-periodic theta functions  satisfy additional identities which reflect
their geometrical origins

\vskip.05in

One can thus view the tau function of Theorem 3.1 as a generalized
theta function, in the spirit of [31]. In the matrix model for two dimensional quantum gravity, 
a $\tau$ function of the KP hierarchy appears as the square
 root of a partition function.

 The $KP$ differential equations reduce to KdV equations in specific cases, and the KdV
hierarchy is specifically associated with hyperelliptic curves. 
A solution to $KP$ that is given by a rational curve with only 
ordinary double points is known as a soliton solution; see [45].\par
\indent Dynamical systems associated with the $KP$ differential equation
can be described in terms of commutative algebras. Mulase showed that any
finite dimensional orbit of a KP dynamical system is an open subset of the
Picard variety of a complete algebraic curve, and the orbit determines a
linear flow on this variety. Shabat observed that the scattering theory
for the KP equation, one can work with matrix coefficients, thus extending
the theory in part to noncommutative rings.\par

\indent (iii) Suppose that $ABC=BCA$, and let $E=2^{-1}A^{-1}BC$, and $D=d/dx$. Then the ring of formal power series ${\bf B}={\bf C}[[x,E,A]]$ is commutative and Noetherian, and a differential ring which is
 closed under indefinite integration and forming exponentials $f\mapsto e^f$. Under these hypotheses, Mulase considers the formal expression $U(t)=\exp (\sum_{j=1}^\infty t_jD^j)U(0)$ where $t=(t_j)_{j=1}^\infty$ is a complex sequence and $U(0)$ is some initial value in ${\bf B}$.\par

\vskip.1in

\indent We now consider elliptic solutions of the KP equations.  Let $\Lambda_0$
 be a lattice in ${\bf C}$ and $\Lambda_1$ a lattice
in ${\bf C}^g$. Let $\theta$ be Riemann's theta function for the Abelian manifold 
${\bf X}_1={\bf C}^g/\Lambda_1$, and suppose that $u,v,w\in {\bf C}^g$ have the property
that $\nu\Lambda_0\subset \Lambda_1$. Then $x\mapsto \theta (ux+vt+w)$ is a theta
function for ${\bf X}_0={\bf C}/\Lambda_0$. We suppose further that 
$\theta$
degenerates into a product
$$\theta (\nu x+vt+w)=\prod_{j=1}^n \sigma (x-x_j(t)),\eqno(10.23)$$
\noindent where $\sigma$ is a theta function for ${\bf X}_0$. The corresponding
potential
$$u(x,t)=-2{{\partial^2}\over{\partial x^2}}\log\theta (\nu
x+vt+w).\eqno(10.24)$$
\noindent Whereas the motion in ${\bf C}^g$ is rectilinear, we seek to express the
solutions in terms of elliptic functions on ${\bf X}_0$. The
idea is that we take $x$ to be the primary variable of an elliptic
function, with poles $x_j(t)$ that move subject to a classical
dynamical system. Thus we reduce the partial differential equation to a
system of coupled ordinary differential equations.\par
\vskip.05in
\noindent {\bf Lemma 10.4} {\sl Let $(x_j, p_j)$ for $j=1, \dots , n$ be the canonical
 coordinates
for the classical Hamiltonian dynamical system with Hamiltonian
$$H={{1}\over{2}}\sum_{j=1}^n p_j^2-4\sum_{j,k:j<k}\wp (x_j-x_k).\eqno(10.25)$$
There exists a periodic linear system that has tau function $\theta (\nu
x+vt+w)$ so that $x\mapsto u(x,t)$ has poles at the $x_j(t)$.}\par

\vskip.05in
\noindent {\bf Proof.} By Lemma 6.1, exists a linear system $(-A,B,C;E)$ that has $\sigma$ as it $\tau$
function. Hence the linear system
$$\bigoplus_{j=1}^n \Bigl( -A, e^{x_j(t)A/2}B, Ce^{x_j(t)A/2};
e^{x_j(t)A/2}e^{x_j(t)A/2}\Bigr)\eqno(10.26)$$
\noindent has $\tau$ function $\prod_{j=1}^n \sigma (x-x_j(t)).$\par
\vskip.05in

\indent Krichever [36] shows that the system  (10.26) is integrable in the sense that 
there exists a compact Riemann surface ${\cal Y}_N$ which covers the elliptic curve $N$-fold,
 and the solution of the Hamiltonian dynamical system can be expressed in action-angle variables 
with the angles in the Jacobi variety of ${\cal Y}_N$.\par 
\indent There exists a complex differential field $({\bf K}_0, d/dx)$ that is
generated by $\sigma$ and its derivatives. Indeed, we take the elliptic function
field ${\bf C}(\wp )[\wp ]$ and adjoin $\zeta$ such that $\zeta'=-\wp$ and then
adjoin $\sigma$ such that $\sigma'/\sigma =\zeta$. By repeating this construction in
the variables, we can introduce a differential field ${\bf K}_n$ for the derivatives
${{\partial}\over{\partial x_1}}, \dots , {{\partial}\over{\partial x_n}}$ that is
generated by $\sigma (x_1), \dots , \sigma (x_n), \sigma (x_1+\dots +x_n)$ and their
derivatives.\par
\indent The potential of this Hamiltonian belongs to ${\bf K}_n$. Indeed, 
$$2\sum_{j,k:j<k}\wp (x_j-x_k)=\sum_{\ell =1}^n {{\partial^2}\over{\partial
x_\ell^2}}\log \prod_{j,k:j<k}\sigma (x_j-x_k),\eqno(10.27)$$
\noindent and the right-hand side belongs to ${\bf K}_n$ on account of the classical
formula
$$\eqalignno{\det\left[\matrix{1&\wp (x_1)&\dots &\wp^{(n-2)}(x_1)\cr
\vdots &{}&\ddots &\vdots\cr
1&\wp (x_n)&\dots &\wp^{(n-2)}(x_n)\cr}\right] &\cr
=(-1)^{(n-1)(n-2)/2}\prod_{k=1}^n
k! &{{\sigma (x_1+\dots +x_n)}\over{\prod_{j=1}^n \sigma
(x_j)^n}}\prod_{j,k:j<k}\sigma (x_j-x_k).&(10.28)\cr}$$
\indent Next we prove that $u$ satisfies $KP$. We introduce the elliptic function of the second kind $\psi_2(x,\alpha )$ as in section 6.
Then we introduce
$$\Psi (x,y,t)=\sum_{j=1}^n a_j(y,t)\psi_2 (x-x_j(y,t), \alpha )e^{kx+k^2y+k^3t}\eqno(10.29)$$
\noindent such that
$${{\partial}\over{\partial y}}\Psi =\Bigl( {{\partial^2}\over{\partial
x^2}}-2\sum_{j=1}^n\wp (x-x_j(y,t))\Bigr) \Psi\eqno(10.30)$$
\noindent and 
$${{\partial}\over{\partial t}}\Psi =
 \Bigl( {{\partial^3}\over{\partial
x^3}}-3\sum_{j=1}^n\wp (x-x_j(y,t)){{\partial}\over{\partial x}}
+{{3}\over{2}}\sum_j \wp (x-x_j){{\partial x_j}\over{\partial y}}
-{{3}\over{2}}\sum_j\wp'(x-x_j)\Bigr) \Psi\eqno(10.31)$$
\noindent Note that $x\mapsto \Psi (x,y,t)$ is elliptic of the second
kind.\par
\vskip.05in
\noindent {\bf Proposition 10.4} (Krichever) {\sl The most general elliptic solution of
$KP$ is given by 
$$u(x,y,t)=C+2\sum_{j=1}^n \wp (x-x_j(y,t))\eqno(10.32)$$ 
\noindent where the canonical equations of motion are used to determine the dynamics 
with respect to time $t$ of the poles
of the functions $x_1(t,y), \dots , x_n(t,y)$ with respect to $y$.}\par
\vskip.05in

\vskip.05in
\noindent {\bf Remark.} We leave it as an open problem to characterize all finite gap cases of Hill's equation in terms of periodic linear systems.\par
\vskip.05in

\indent The $KP$ equation (10.4) is the first of a hierarchy of 
differential equations.\par
\vskip.1in

\vskip.1in

\noindent {\bf 11. The differential ring associated with the
Painlev\'e II equation}\par
\vskip.05in

\noindent In this section we consider a linear system which is
important in random matrix theory. Whereas the state ring ${\bf S}$
is finitely generated, the linear system is not integrable in the sense
that $\tau$ does not emerge from ${\bf C}(x)$ by successive Abel integrations.
Let $H(w,v;x)$ be a Hamiltonian which is rational in
the canonical variables
$(w,v)$ and a meromorphic function of time $x$, and let $(w(s), v(s))$ be 
solutions of the
canonical equations of motion, and suppose momentarily that these are
meromorphic functions of $s$. Then the corresponding tau function is
$$\tau (x)=\exp \int_0^x H(w(s), v(s);s)\, ds, \eqno(11.1)$$
\noindent where the integral is taken along an orbit in phase space; so
 the value of $\tau$ is
locally independent of the path of integration, provided the path
avoids poles.\par
\indent The Hamiltonians which arise on random matrix theory 
have additional properties which are described in the following result,
which is a variant of  Theorem 1 in Okamoto's paper [46].\par 
\vskip.05in
\noindent {\bf Proposition 11.1} {\sl Suppose that the Hamiltonian $H(w,v;x)$ is rational in
$x$, a polynomial in $v$, and a quadratic polynomial in $w$, let $u$
be the potential that corresponds to $\tau$; let ${\bf A}_0={\bf C}(x)[w,v]$. Then
${\bf A}_0$ is a differential ring containing $u$ such that for all maximal ideals
$J_0$ of ${\bf A}_0$, there exists a finite algebraic cover of the Riemann
sphere.}\par
\vskip.05in
\noindent {\bf Proof.} We write $H=A(v,x)w^2+B(v,x)w+C(v,x)$. Then the canonical
equations are ${{dv}\over{dx}}={{\partial H}\over{\partial w}}$ and
 ${{dw}\over{dx}}=-{{\partial H}\over{\partial v}}$.
 Hence ${\bf A}_0={\bf C}(x)[w,v]$ is a
commutative differential ring for the derivative
${{d}\over{dt}}={{\partial}\over{\partial x}}+
{{\partial H}\over{\partial w}}{{\partial}\over{\partial v}}-
{{\partial H}\over{\partial v}}{{\partial}\over{\partial w}}$.
Evidently ${\bf A}_0$ is finitely generated over the field ${\bf
C}(x)$, hence is Noetherian.\par
\indent Note also that 
$$u(x)=-2{{\partial H}\over{\partial x}}=-2{{\partial A}\over{\partial
x}}w^2-2{{\partial B}\over{\partial x}}w-2{{\partial C}\over{\partial x}}\eqno(11.2)$$
\noindent belongs to ${\bf A}_0$.\par
\indent For each maximal ideal $J_0$ of ${\bf A}_0$, we have by the weak 
Nullstellensatz, a finite algebraic extension ${\bf A}_0/J_0$ of ${\bf C}(x)$,
which may be realised as an algebraic function. Hence there exists a Riemann
surface ${\bf X}$ and a holomorphic and surjective map $\pi :{\bf X}\rightarrow {\bf
P}^1$ which gives a cover of the Riemann sphere with finitely many sheets.\par
\rightline{$\square$}\par
\vskip.05in  
Suppose furthermore that $K$ defines a trace class operator on 
$L^2(a,b )$ for some real interval $(a,b)$; hence the Fredholm determinant $\det (I-K)$ is defined. In particular, one can consider $K:L^2(0, \infty )\rightarrow L^2(0, \infty )$ that is trace class and
such that $0\leq K\leq I$, so there exists  a determinantal  random point field on
$(0, \infty )$, and $\det (I-K{\bf I}_{(s, \infty )})$
is the probability that all random points are in $(0, s)$. In applications
to random matrix theory, the random points are eigenvalues of Hermitian matrices with random entries. For significant cases of (1.2),
such as the Airy kernel or Bessel kernel [59, 60],
there exists a Hankel integral operator $\Gamma_\phi$ such that
$\Gamma_\phi^2 =K$; hence one can describe $\det (I-K)$ in terms of
$\tau (x ,\mu )$. In [10] we showed how one can realise $\Gamma_\phi$ by
means of linear systems. In the present paper, we take linear systems
as the starting point and show how general properties of the linear
system are reflected in the $\tau$ functions and systems of
differential equations so produced.
 By making appropriate substitutions, one can reduce each of the main models in [59] to a $2\times 2$
rational differential equation and a family of subintervals of the real line. 

\indent Let ${\bf K}_0={\bf C}(x)$, and suppose that $H(\lambda , \mu ;x)$ is a rational function 
which we regard as a Hamiltonian with canonical equations
$$\lambda'=-{{\partial H}\over{\partial \mu}}, \qquad \mu'={{\partial H}\over{\partial
\lambda}},\qquad {}'={{d}\over{dx}}.$$
\noindent Then ${\bf C}(\lambda , \mu ,x)$ is a differential field. Okamoto showed that each of
the Painleve transcendental differential equations $P_{I}, \dots , P_{VI}$ could be expressed as
canonical equations of motion for suitable Hamiltonians $H_1, \dots ,H_{VI}$, and that the corresponding tau functions
are
$$\tau_J (x)=\exp \Bigl( \int_0^x H_J(\lambda (t), \mu (t); t)\, dt\Bigr).\eqno(11.3)$$
\noindent We can therefore introduce differential fields ${\bf K}_2(P_J)={\bf C}(\lambda, \mu ,x,
\tau_J)$ for $J=I, \dots ,VI$.  
$$\matrix{ {\hbox{field}}& {\hbox{soft edge}}& {\hbox{hard edge}}& {\hbox{Hermite}}&{\hbox{Laguerre}}& {\hbox{Lam\'e}}\cr
            {\bf K}_0& {\bf C}(x)& {\bf C}(x)&{\bf C}(x)& {\bf C}(x)&{\bf C}({\cal T})\cr
             {\bf K}_1& {\bf K}_0({\hbox{Ai}}(x), {\hbox{Ai}}'(x))& {\bf K}_0({\hbox{J}}_\alpha (\sqrt{x}),
\sqrt{x}{\hbox{J}}_\alpha'(\sqrt{x}))& {\bf K}_0(e^{-x^2/2})& {\bf K}_0(x^{\alpha /2}e^{-x/2})& {\bf C}(\Gamma_\ell )\cr
{\bf K}_2& P_{II}& P_{III}& P_{IV}& P_V&{}\cr}$$


\indent Okamoto [46] has shown that
each of the Painlev\'e transcendental differential equations $P_{I}, \dots ,
P_{VI}$ arises from a Hamiltonian as in Proposition 10.1, and $\tau$ is meromorphic on a suitable covering surface. In particular,  for $x\in {\bf C}$  and a complex constant $\alpha$, let 
$$H_{II}(w,v;x)={{1}\over{2}}\Bigl(w-{{x}\over{2}}\Bigr)^2+\Bigl( 
v^2+{{x}\over{2}}\Bigr)\Bigl(
w-{{x}\over{2}}\Bigr) -\alpha v +{{x^2}\over{8}}.\eqno(10.4)$$ 
\vskip.05in
\noindent {\bf Proposition 10.2} {\sl Under the canonical equations of
motion with Hamiltonian $H_{II}$, \par
\indent (i) $v$ satisfies 
$P_{II}
:v''=xv+2v^3+\alpha $ and the corresponding $\tau$ function is
$$\tau (x)=\exp\Bigl( -{{1}\over{2}}\int_x^\infty (s-x)v(s)^2\,
ds\Bigr);\eqno(10.5)$$
\indent (ii) 
$w$ satisfies $w'''+6ww'-(2w+xw')=0$ and 
$U(x,t)=(3t)^{-2/3}w(3^{-1}t^{-1/3}x)$ satisfies}
$${{\partial^3U}\over {\partial
x^3}}+{{2U}\over{3^{1/3}}}{{\partial U}\over{\partial
x}}-{{1}\over{9}}{{\partial U}\over{\partial t}}=0.\eqno(10.6)$$ 

\vskip.05in
\noindent {\bf Proof.} (i) Proposition 10.1 applies due to the special form of the Hamiltonian. 
The equation (10.2) reduced to $u=v^2$, and hence we obtain the
expression for $\tau$.\par
\indent (ii) Now $w$ satisfies    
$$K_2:\qquad w''+2w^2-xw+{{\alpha (\alpha
+1)+w'-(w')^2}\over{2w-x}}=0.\eqno(10.7)$$
\indent One can then verify that $U$ satisfies KdV; see [2] for
further discussion. \par
\rightline{$\square$}\par
\vskip.05in
\indent Proposition 11.2 indicates that we have potentials linked by the Miura transform, as in (2.37), and the next step is to introduce suitable linear systems and operators to solve the differential equation $P_{II}$. We introduce
Airy's function ${\hbox{Ai}}(x)=\int_{-\infty}^\infty e^{i\xi
x+i\xi^3/3}d\xi /(2\pi )$, which satisfies ${\hbox{Ai}}''(x)=
x{\hbox{Ai}}(x)$. Let 
$\phi (x)={\hbox{Ai}}(x)$ and let $\zeta =\phi'/\phi$; then
${\bf S}={\bf C}[x,\phi (x) ,\zeta (x)]$ is a differential ring with
respect to $d/dx$. In the context of Theorem 10.3(iii) below the
integrable kernel
$$R_0^2(x,y)={{{\hbox{Ai}}(x){\hbox{Ai}}'(y)-{\hbox{Ai}}'(x){\hbox{Ai}}
(y) }\over{x-y}}\eqno(11.9)$$
\noindent is known as Airy's kernel, which is associated with soft
edges of eigenvalue distributions [59]. The Fredholm determinants of $R_0^2$ 
 lead to a solution of the Painlev\'e II nonlinear differential 

\vskip.05in
\indent The results of this section extend to differential equations which arise in unitary matrix models. For each $\ell =1, 2,\dots,$ there exists a real solution ${\hbox{A}}_\ell  (x)$ to the generalized Airy equation 
$${\hbox{A}}_\ell^{(2\ell )}(x)+(-1)^\ell x{\hbox{A}}_\ell (x)=0
\eqno(11.45)$$
\noindent such that ${\hbox{A}}_\ell (x)\rightarrow 0$ exponentially as $x\rightarrow 0$ and ${\hbox{A}}_\ell (x)$ is algebraically damped and oscillating as $x\rightarrow -\infty$.
We introduce a formal group of operators $U(t)$ 
$$U(t)f(x)={{1}\over{2\pi}}\int_{-\infty}^\infty \exp\Bigl({\sum_{j=0}^\infty (-1)^jit_{2j+1}\xi^{2j+1} +i\xi x}\Bigr) \hat f(\xi )\, d\xi\qquad  (t=(t_{2j+1})_{j=0}^\infty )$$
\noindent for real sequences $t$ that have only finitely many nonzero terms.

In [16], Crukovic, Douglas and Moore select a suitable curve $C$ in the integral
$${\hbox{A}}_\ell\ (x)
={{1}\over{\pi}}\int_C \exp i\Bigl( wx+{{w^{2\ell +1}}\over{2\ell +1}}
\Bigr) dx\eqno(11.46)$$
\noindent and use the stationary phase method to prove the existence and uniqueness of $\phi$. \par
\vskip.05in
\noindent {\bf Definition} (Airy equation) (i) Let $\phi (x)={\hbox{A}}_\ell (x/2)$, and let the generalized Airy kernel be 
$$K_{(x)}(t,z)=2^{2\ell +1}
\sum_{k=1}^\ell (-1)^{k-1-\ell} {{\phi^{(2\ell -k)}_{(x)}
(t)\phi^{(k-1)}_{(x)}(z)
-\phi^{(k-1)}_{(x)}(t)\phi^{(2\ell
-k)}_{(x)}(z)}\over{t-z}}\eqno(11.47)$$
\noindent which reduces in the case $\ell =1$ and $x=0$ to the 
standard Airy kernel as in [59].\par
\indent (ii)  Let $J$ be a finite union of open sub intervals of $(0,\infty )$ and let $P_J$ be the orthogonal projection on $L^2(0, \infty )$ given by multiplying functions by the indicator function of $J$. A point configuration on $J$ is a function $\nu :J\rightarrow \{ 0, 1, 2, \dots \}$ such that $\{ x\in S: \nu (x)>0\}$ is finite for all compact subsets $S$ of $J$; write $\nu (S)=\sum_{x\in S}\nu (x)$. Let ${\hbox{Conf}}(J)$ be the set of all point configurations on $J$ with the $\sigma$-algebra generated by the functions $\nu$, so $\nu\mapsto \nu (S)$ is a random variable on ${\hbox{Conf}}(J)$. We regard $\nu (S)$ as the number of random points in $S$. \par
\indent (iii) For the remainder of this section, we temporarily 
change notation, so that we conform with the standard references. 
Gelfand and Dikii [26] introduced a sequence of polynomials $R_j[u]\in {\bf R}[u,u',u'', \dots ]$ by the recurrence relation
$$R'_{j+1}=4^{-1}R_j'''-uR_j'-2^{-1} u'R_j\eqno(11.48)$$
\noindent and the sequence begins with
$$R_0=2^{-1},\quad R_1=-4^{-1}u,\quad
R_2=16^{-1}(3u^3-u'').\eqno(11.49)$$

\vskip.05in 
\noindent {\bf Proposition 11.5} {\sl (i) Then $K_{(x)}=\Gamma_{\phi_{(x)}}^2$, so $K_{(x)}$ defines a trace class operator on $L^2(0,\infty )$ such that $0\leq K_{(x)}$ and $K_{(x)}\rightarrow 0$ as $x\rightarrow\infty $.}\par
\indent {\sl (ii) Let $v$ be the differentiated quotient of determinants
$$v(x)={{d}\over{dx}}{{1}\over{2}}\log 
{{\det (I+\Gamma_{\phi_{(x)}})}\over{\det (I-\Gamma_{\phi_{(x)}})}};
\eqno(11.50)$$
\noindent then let $u(x)=v'(x)+v(x)^2$, so that $u$ is the Miura transform of $v$. Then $v$ satisfies the string equation
$$P_{II}^{(\ell )}:\qquad {{d}\over{dx}}
R_\ell [u]-vR_\ell [u]+(-1)^{\ell +1}2^{-2\ell }xv(x)=0,\eqno(11.51)$$
\noindent which reduces in the case $\ell=1$ to $P_{II}$.}\par
\indent {\sl (iii) Then there exists a probability measure on} ${\hbox{Conf}}((0, \infty ))$ {\sl such that 
$$E(n;J)={{(-1)^n}\over{n!}}{{d^n}\over{d\lambda^n}}
\Bigl\vert_{\lambda =1}\det (I-\lambda K_{(x_0)}P_J)\eqno(11.52)$$
\noindent equals the probability of $\{ \nu :\nu (J)=n\}$ for all $J$ and all $n=0, 1, \dots.$ }\par

\vskip.05in
\noindent {\bf Proof.} (i) Since $\phi$ is real and rapidly decreasing, the Hankel operator $\Gamma_{\phi_{(x)}}$ is self-adjoint and Hilbert--Schmidt, hence $\Gamma_{\phi_{(x)}}^2\geq 0$ is trace class. To prove that $K_{(x)}=\Gamma_{\phi_{(x)}}^2$, we multiply the right-hand side of the claimed identity by $(t-z)$ and use the differential equation to obtain
$${{(-1)^{\ell +1}(t-z)}\over{2^{2\ell +1}}}\int_{2x}^\infty \phi 
(t+s)\phi (s+z)ds=\int_{2x}^\infty \bigl( \phi^{(2\ell )}(t+s)\phi 
(s+z)-\phi (t+s)\phi^{(2\ell )}(s+z)\bigr)ds.\eqno(11.53)$$
\noindent After integrating by parts $\ell$ times, we obtain the formula
$$K_{(x)}(t,z)=\int_{2x}^\infty \phi (t+s)\phi (s+z)ds.\eqno(11.54)$$
\indent (ii) This follows from Theorem 2.6. Zakharov and Shabat
[62] show that the Gelfand--Levitan equation leads directly to a system of differential equations, which in this case gives the generalized Painlev\'e equations
From the recurrence relation, we obtain the sequence of differential equations that starts with
$$\eqalignno{P_{II}^{(1)}:&\qquad v''+2v^3-xv=0,\cr
P_{II}^{(2)}:&\qquad v^{(4)}-10v''v^2-10(v')^2v+6v^5+xv=0.&(11.55)\cr}$$
\indent (iii)
For all $x_0>0$ sufficiently large, the kernel satisfies 
$0\leq K_{(x)}\leq I$ for all $x\geq x_0$ and by Mercer's theorem 
$P_JKP_J$ is a trace class kernel for any compact subset $J$ of 
$(0, \infty )$ as in [55]. Hence $K_{(x)}$ defines a determinantal random point 
field on $(x_0, \infty )$ by Soshnikov's theorem [Theorem 3, 53].
We refer the reader to this paper for more precise definitions of the construction of the 
measure such that $E(n,J)$ is the probability of $\{ \nu \in {\hbox{Conf}}((x_0,\infty )): \nu (J)=n\}$.\par
\rightline{$\square$}\par
\vskip.05in 
\indent 
The string equation (11.51) arises from a double scaling of a 
probabilistic model of unitary random matrices, so one would expect
 $K_{(x)}$ to have a probabilistic interpretation in the style of 
Tracy and Widom. The system of differential equations (11.55) is associated 
with scaling solutions of the Zakharov--Shabat hierarchy, and
 is sometimes known as the Painlev\'e II hierarchy; although there are competing notions of this concept. One can also introduce a hierarchy of differential equations from the system from the dispersive water wave equations, the first terms of which reduces to the system of linear equations which that Jimbo and Miwa [30] use to linearize $P_{II}$; but at the time of writing it is not clear as to whether the various approaches lead to equivalent hierarchies. 
\par

 \vskip.05in

\vskip.05in

\vskip.05in
\noindent {\bf Acknowledgements} The EPSRC supported SLN's research via a research studentship. GB carried out part of this work
during visits to the Isaac Newton Institute in Cambridge and the University of Glasgow. The authors are grateful to Prof H.P.  McKean for his
comments on [41, 43], to Dr M. Mazzocco and Prof S.Power for helpful remarks.\par 
\vskip.05in

\end
The curve ${\cal E}$ and the differential equation (7.2)
 depend upon
the branch points $a_1, \dots ,a_{g+1},$ $b_1, \dots , b_g$. To investigate
this further, we introduce $2\times 2$ matrices $W_k(t)$ 
which are meromporphic functions of $t=(t_1, \dots , t_n)$ and consider the
differential equation

$${{d}\over {dx}}\Phi =\sum_{j=0}^N {{W_j(t)}\over{x-t_j}}\Phi\eqno(7.9)$$
\noindent 
Krichever and 
Novikov introduced the
notion of spectral data for a family of commuting differential
 operators [38]. We give the simplified version which was previously considered by
Garnier.\par
\vskip.05in 
\noindent {\bf Definition} (Algebraic family) Let ${\bf K}_0$ be the 
field of meromorphic functions in $t=(t_0, \dots , t_N)$, let ${\bf K}_1={\bf K}_0(x)$ so that $({\bf K}_1, d/dx)$ is a differential field with constants ${\bf K}_0$. Let $W_j(t)$ be $2\times 2$ matrices with entries in ${\bf K}_0$, and such that ${\hbox{trace}}
W_j(t)=0$ for all $t$
and $j=0, \dots , N$; then let
 $$L_j={{\partial 
}\over{\partial t_j}}+
{{W_j(t)}\over{x -t_j}}\qquad (j=0, \dots, N),\eqno(7.10)$$
Suppose that [25] $[L_j, L_k]=0$ for all $j,k$ and that 
$$\Bigl[ L_j,\sum_{k=0}^N {{W_k}\over {x -t_k}}\Bigr]=0
\qquad (j=0, \dots, N).\eqno(7.11)$$
\noindent Let ${\bf K}_1$ be the space of meromorphic functions in $x$ and
$t$. Let $W_k(t)$ be meromorphic functions of $t=(t_0, \dots ,t_N)$ and
suppose that $\det W_k(t)\neq 0$ whereas ${\hbox{trace}}W_k(t)=0$ for all
$k$.  We introduce the $2\times 2$ polynomial matrix function
$$A(x, \kappa ,t)=I_2\prod_{j=0}^N (x-t_j) -\kappa \sum_{j=0}^N
W_j(t)\prod_{k=0; k\neq j}^N(x-t_k),\eqno(7.12)$$
\noindent which is monic in the sense that the leading term is
$x^{N+1}I_2$. and introduce the curve
$${\cal F}=\Bigl\{ (\kappa , x ): \det A(x,\kappa ,t)=0\Bigr\}.\eqno(7.13)$$
\vskip.05in
Let ${\bf K}$ be a field and $A(z)=z^{N+1}I_2+\dots +A_0$ a
monic $2\times 2$ matrix polynomial of degree $N+1$ with coefficients
$A_j$ that have entries in ${\bf K}$. Say that $zI-B$ is a right
divisor if there exists a $2\times 2$ matrix polynomial $C(z)$ such that
$A(z)=C(z)(zI-B)$. The matrix polynomial $A(z)$ is said to be generic
if there exist $2N+2$ linearly independent eigenvectors $\Psi_0,
\dots , \Psi_{2N+1}$.
\vskip.05in
\noindent {\bf Proposition 7.3} {\sl (i) The curve ${\cal F}$ is independent of $t$ and 
gives a hyperelliptic curve of genus $g$ where $g\leq N$.\par
\indent (ii) Suppose that $x\mapsto \det A(x,\kappa ,t)$ has simple zeros
for some $\kappa$. Then for each zero, there exists a corresponding eigenvector
$\Psi (t, x )$ such that 
$$L_\ell \Psi (t,
x )=\xi_\ell (x ,t)\Psi (x , t)\eqno(7.14)$$
\noindent  for some $\xi_\ell (x,\kappa ,t)\in {\bf C}$. These
eigenvectors determine the right divisors of $A(x,\kappa ,t)$.}\par 
\vskip.05in
\noindent {\bf Proof.} (i) Suppose that $\det A(x,\kappa ,t)=0$, and note
that $x\neq t_k$ for all $k$, so we deduce that $\det (I-\kappa \sum_{k=0}^N
W_k(t)/(x-t_k))=0$. We note that ${{\partial }/{\partial t_j}}
\det (I-\kappa \sum_{k=0}^N W_j(t)/(x-t_k))=0$, so 
$\{ x:\det (A(x,\kappa ,t))=0\}$ is
independent of $t$. Also $\det A(x,\kappa ,t)$
is quadratic in $\kappa$ and a polynomial in $x$ of degree
less than or equal to $2N+2$. Hence ${\cal F}$ gives a algebraic curve 
of genus $g\leq N$. 

\indent We introduce
$2\times 2$ matrices $B_\ell =B_\ell (\mu ,t )$ by 
$$B_\ell\Psi_{2\ell}=x_{2\ell}\Psi_{2\ell},
B_\ell\Psi_{2\ell+1}=x_{2\ell +1}\Psi_{2\ell+1}\qquad (\ell =0, \dots ,
N).\eqno(7.16)$$
\noindent Then for each $k=0, \dots, N$, there exists a $2\times 2$
matrix $C_k(x,\kappa ,t)$ such that $x\mapsto C_k(x,\kappa ,t)$ is a
polynomial of degree $N$, and 
$$A(x,\kappa ,t)=C_k(x, \kappa ,t)(xI_2-B_k(\mu ,t))\eqno(7.17)$$
\noindent as a polynomial identity; thus $xI_2-B_k(\mu ,t)$ is a right
divisor of $A(x, \kappa ,t)$. \par
\indent To conclude the proof, we recall the following result from
[4'].There is a birational correspondence between
generic monic matrix polynomials and their right divisors 
$$A(z)\leftrightarrow (zI_1-B_0, \dots , zI_2-B_N)\eqno(7.18)$$
 \noindent such that the entries of the coefficient matrices $A_j$ a
rational functions of the $B_k$ and the eigenvalues of the $B_k$.\par
\rightline{$\square$}\par
\vskip.05in
\vskip.05in

\noindent {\bf 10. Discrete integrable operators}\par

\indent In the literature of group representation theory, $\tau$ functions of this form 
appear implicitly. Let $\lambda$ be a partition of $N$, and recall that the conjugate 
$\lambda'$ is the partition such that the Ferrers diagrams corresponding to $\lambda$ and $\lambda'$ are mirror images upon reflecting about the main diagonal. 
There exists an irreducible representation of the symmetric group on $N$ letters on a Specht module $H_\lambda$; the dimension of $H_\lambda$ appears in the following $\tau$ function. We introduce
$$\eqalignno{ B:\ell^2\rightarrow \ell^2:& \qquad \beta\mapsto (2i\beta /j!)_{j=0}^\infty ;
\cr
   A:{\cal D}(A)\rightarrow \ell^2:& \qquad (\xi_j)_{j=0}^\infty \mapsto
 (i(2j+1)\xi_j)_{j=0}^\infty ;\cr
C:\ell^2\rightarrow {\bf C}:& \qquad (\xi_j)_{j=0}^\infty \mapsto \sum_{j=0}^\infty \xi_j/j^!;\cr
E:\ell^2\rightarrow\ell^2:& \qquad E=\Bigl[{{1}\over{(j+k+1)\, j!k!}}\Bigr]_{j,k=0}^\infty.
&(10.3)\cr}$$
\vskip.05in
\noindent {\bf Proposition  10.1} {\sl Then $(-A,B,C;E)$  gives a periodic linear system, with corresponding $\tau$ function}
$$\tau (x)=\sum_\lambda e^{-2ix\vert \lambda\vert}{{{\hbox{dim}} 
(H_\lambda)}\over{\vert\lambda \vert!}}\eqno(10.4)$$
\noindent {\sl where the sum is over all partitions $\lambda$ such that $\lambda'=\lambda$.}\par
\vskip.05in
\noindent {\bf Proof.}  The hypotheses of Lemma 8.1(ii) are satisfied, since $(e^{-xA})_{x\in {\bf R}}$ is unitary and periodic as a group on $\ell^2$, and it is straightforward to check that $AE+EA=BC$. Let $\lambda$ be any partition such that $\lambda=\lambda'$ , and let $p_j$ be the number of boxes in the $j^{th}$ row of the Ferrers diagram to the right of the main diagonal. Then Frobenius's coordinates satisfy  $\lambda =(p_1, \dots , p_d; p_1, \dots ,p_d)$ where $p_1>p_2>\dots >p_d$ are positive integers such that $2\sum_{j=1}^d p_j=N-d$. We introduce
$$L_\lambda =\Bigl[ {{1}\over{(p_j+p_k+1)\, p_j!p_k!}}\Bigr]_{j,k=1, \dots ,d}\eqno(10.5)$$
\noindent which is a diagonal submatrix of $E$. We can express the $\tau$ function as a series of determinants of diagonal blocks, and write 
$$\det (I+e^{-xA}Ee^{-xA})=\sum_\lambda e^{-2ix\vert\lambda\vert}\det L_\lambda ,\eqno(10.6)$$ 
\noindent where the value of $\det L_\lambda$ is given by a version of the hook length formula.\par
\rightline{$\square$}\par

\noindent {\bf Definition} (Discrete integrable operators ) (i) Let $H_0$, $H_1$ and $H=\sum_{j=-N}^NH_1$ be Hilbert spaces, and let 
let $B_j:H_0\rightarrow H_1$ and $C_j:H_1\rightarrow H_0$ be Hilbert--Schmidt operators such that $B_0=0=C_0$ and $B_jC_{-j}=0$ for all $j$. Then
we say that $E:H\rightarrow H$ 
$$E=\Bigl[ {{iB_jC_k}\over{j+k}}\Bigr]_{j,k=-N, \dots ,-1,1, \dots , N}\eqno(10.7)$$
\noindent is a discrete integrable operator.\par
\indent (ii) For  $H_0={\bf C}$ and  $H_1={\bf C}^{2\times 1}$, let  $Y_0\in H_0$, and $V_n\in M_2({\bf C})$ such that $\det V_n=1$  and suppose that $n\mapsto V_n$ is a rational function;  then  define $Y_n\in H_1$ by the recurrence relations
$$Y_{n+1}=V_nY_n,\qquad  Y_{-n-1}=V_n^{-1}Y_{-n}\qquad (n=0, 1, 2, \dots);\eqno(10.8)$$
\noindent then let $E$ be the linear operator with matrix
$$[E_{n,k}]=\Bigl[ {{\langle JY_n, Y_{-k}\rangle}\over {n+k}}
 \Bigr]_{n,k=-N, \dots ,-1,1, \dots , N},\quad J=\left[\matrix{0&-1\cr
1&0\cr}\right],\eqno(10.9)$$
\noindent where $\langle,\rangle: H_1^2\rightarrow H_0$ is the usual bilinear map. The recurrence relations and the matrix $E$ define a discrete Tracy--Widom system.  \par
\vskip.05in
\indent These concepts have been used implicitly in several papers, and the analogy with rational differential equations used to advantage in [Borodin, etc]. An important example in random matrix theory is known as the discrete Bessel kernel, in which the Bessel function is
$$J_n(z)={{1}\over{2\pi}}\int_{-\pi}^\pi e^{-in s+iz\sin s}ds\eqno(10.10)$$
\noindent and $Y_n(\theta )={\hbox{column}}[\sqrt{\theta}J_{n}(2\sqrt{\theta}), J_{n+1}(2\sqrt{\theta})]$. 
More generally, we have the following result.\par 
\vskip.05in
\noindent {\bf Proposition 10.2} {\sl (i) Given a discrete integrable operator, there exists a uniform periodic linear system $(-A,B,C;E)$ with corresponding $E$.}\par 
\indent {\sl (ii) For any discrete Tracy--Widom system, the discrete integrable operator $E$ satisfies}
$$E_{m+1, -k-1}=E_{m,-k}+\Bigl\langle J{{V_m-V_k}\over{m-k}}Y_m,
Y_{k+1}\Bigr\rangle.\eqno(10.11)$$

\vskip.05in
\noindent {\bf Proof.} (i) We write $H$ as a column of copies of $H_1$, and introduce 
$$\eqalignno{A:H\rightarrow H:&\qquad {\hbox{column}}[\xi_j ]\mapsto {\hbox{column}}[-ij\xi_j]\qquad (\xi_j\in H_1);\cr
B:H_0\rightarrow H:   &\qquad \xi\mapsto {\hbox{column}} [B_j\xi]\qquad (\xi\in H_0);\cr
C:H\rightarrow H_0:&\qquad  {\hbox{column}}[\xi_j]\mapsto \sum_j C_j\xi_j\qquad (\xi_j\in H_1)
&(10.12)\cr}$$
\noindent gives a periodic linear system $(-A, B, C;E)$ where $e^{-xA}$ is a uniformly continuous periodic group and $AE+EA=BC.$
The operator $A$ generates rotations of the circle, and hence gives a unitary operator on $L^2{(\bf T})$.\par
\indent (ii) Note that $\langle JY_n,Y_n\rangle =0$, so we have a discrete integrable operator. From the recurrence relation, we have
$$\langle JY_{m+1}, Y_{n+1}\rangle =\langle V_n^TJV_nY_m,Y_n\rangle +\langle J(V_m-V_n)Y_m,V_nY_n\rangle ,$$
\noindent and the identity follows by direct calculation, since $V_n^TJV_n=J.$ \par
\rightline{$\square$}\par

\vskip.1in

\indent In the following result, we consider $\phi$ that is a bounded
and holomorphic function on the upper half plane, and satisfying some
additional decay conditions. The corresponding $\tau$ function from () is
shown to be a scattering function of a linear system on a large
Hilbert space. This is consistent with the definition of the tau
function on a Fock space given in [Miwa].\par 
\vskip.05in
\noindent {\bf Theorem 4.5} {\sl 
 Let $\phi$ be a scattering function that arises as $\phi
(x)=(2\pi)^{-1}\int_{0}^\infty e^{i\sigma x} f(\sigma )d\sigma$, where
$f(\sigma )$ and $f(\sigma ) /\sigma\in L^1(0, \infty )$, and suppose
that $\vert \phi (x)\vert \leq M/(1+x)^2$ for all $x>0$ and some $M$. 
Then there exists a linear system $(-A,B,C)$ on $H$, and $x_0\geq 0$, such that}
$$\det (I+\Gamma_{\phi_{(x)}})=Ce^{-xA}B\qquad (x>x_0)\eqno(4.16)$$
\noindent {\sl and there exists a solution to the Lyapunov equation ().}
\vskip.05in
\noindent {\bf Proof.} By considering each summand in the determinant
expansion, we deduce that 
$$\bigl\vert \det [\phi_{(x)} (x_j+x_k)]_{j,k=1}^m\bigr\vert\leq
n!M^n\prod_{j=1}^n{{1}\over{(1+2x+x_j)^2}},\eqno(4.17)$$
\noindent and hence 
$${{1}\over{n!}}
\int_{(0,\infty )^n} \bigl\vert 
\det [\phi_{(x)} (x_j+x_k)]_{j,k=1}^m\bigr\vert dx_1\dots dx_n\leq
{{M^n}\over{(1+2x)^n}}.\eqno(4.18)$$
\noindent This estimate ensures convergence of the Fredholm determinant
expansion, as follows. Denoting the exterior power by a wedge, we have
$$\det (I+\Gamma_{\phi_{(x)}})=1+\sum_{n=1}^\infty
{{1}\over{n!}}{\hbox{trace}}(\Gamma_{\phi_{(x)}})^{ \wedge n},\eqno(4.19)$$
\noindent where the trace in the $n^{th}$ summand is 
$$\eqalignno{\int_{(0, \infty )^n}& \det \bigl[ \phi_{(x)}
(x_j+x_k)\bigr]_{j,k=1}^ndx_1\dots dx_n\cr
&={{1}\over{2\pi}}\int_{(0, \infty )^n}\int_{(0, \infty )^n}e^{2ix\sum_{j=1}^n \sigma_j} f(\sigma_1)\dots
f(\sigma_n)\det [e^{i\sigma_j(x_j+x_k)}]_{j,k=1}^n d\sigma_1\dots d\sigma_n dx_1\dots dx_n\cr
&={{1}\over{2\pi}}\int_{(0, \infty )^n}\int_{(0, \infty )^n}e^{2ix\sum_{j=1}^n \sigma_j} f(\sigma_1)\dots
f(\sigma_n)\det\Bigl[{{i}\over{\sigma_j+\sigma_k}}\Bigr]_{j,k=1}^nd\sigma_1\dots d\sigma_n\cr 
&={{i^n}\over{(4\pi
)^n}}\int_{(0,\infty )^n}e^{2ix\sum_{j=1}^n\sigma_j}\prod_{j=1}^n 
{{f(\sigma_j)}\over{\sigma_j}}\prod_{j<k}
\Bigl({{\sigma_j-\sigma_k}\over{\sigma_j+\sigma_k}}\Bigr)^2d\sigma_1
\dots d\sigma_n,&(4.20)\cr}$$ 
\noindent where the final step follows from the Cauchy--Binet formula. We introduce a linear system $(-A,B,C)$ by
$$\eqalignno{B_n:{\bf C}\rightarrow L^2((0, \infty )^n):& \beta\mapsto
{{i^n}\over{(2\pi )^n}}\prod_{j=1}^n \Bigl\vert {{f(\sigma_j)}\over{2\sigma_j}}\Bigr\vert^{1/2}\beta;\cr
e^{-xA_n}:L^2((0, \infty )^n)\rightarrow L^2((0, \infty )^n):& g(\sigma
)\mapsto e^{2ix\sum_{j=1}^n \sigma_j} g(\sigma );\cr
C_n:L^2((0, \infty )^n)\rightarrow {\bf C}:& g\mapsto \int_{(0, \infty )^n}
g(\sigma )\prod_{j=1}^n {{f(\sigma_j)}\over{\vert
2\sigma_jf(\sigma_j)\vert^{1/2}}}\prod_{j<k}\Bigl(
{{\sigma_j-\sigma_k}\over{\sigma_j+\sigma_k}}\Bigr)^2 d\sigma_1\dots
d\sigma_n&(4.21)\cr}$$
\noindent which satisfies
$$C_ne^{-xA_n}B_n={\hbox{trace}}(\Gamma_{\phi_{(x)}})^{ \wedge n}.\eqno(4.21)$$
\noindent Also let, $A_0=0$, $B_0=1=C_0$. On the Hilbert space $H={\bf C}\oplus\bigoplus_{n=1}^\infty L^2((0,
\infty )^n)$, the linear system 
$$(-A,B,C)=\Bigl( -\bigoplus_{n=0}^\infty A_n,  
\bigoplus_{n=0}^\infty {{B_n}\over{\sqrt{ n!}}},  
\bigoplus_{n=0}^\infty {{C_n}\over{\sqrt{n!}}}\Bigr),\eqno(4.22)$$
\noindent where $A$ is a sum of diagonal blocks, $B$ is a column and $C$ is
a row, has the required properties. Let $\nu_j=\sum_{\ell
=1}^j\sigma_\ell$. The Lyapunov equation has solution
$$R_x=\Bigl[{{e^{2ix\nu_j}B_jC_ke^{2ix\nu_k}}\over{2i\sqrt{j!k!}(\nu_j+\nu_k)}}\Bigr]_{j,k=0}^\infty .$$
\indent Observe that the formula () extends to define a bounded linear
operator $e^{-(x+iy)A_n}$ for $x\in {\bf R}$ and $y>0$, so we can
interpret $Ce^{-xA}B$ as the boundary values of a bounded and
holomorphic function on the upper half plane.\par
\rightline{$\square$}\par
\vskip.05in

\noindent {\bf Remarks 4.6} (i) Ercolani and McKean interpret
the $\sigma_j$ as random points on the line, and deduce an analogy
between the determinant () and classical theta functions on ${\bf C}^n$
with respect to the lattice ${\bf Z}^n+i{\bf Z}^n$. The probabilistic
interpretation can be expressed in a form which is familiar from
electrostatics. Suppose that $f(\sigma )=e^{i\theta (\sigma )} \vert f(\sigma
)\vert$ is the polar decomposition of $f$, and let $w(\sigma
)=\log \sigma -\log\vert
f(\sigma )\vert.$  As in the method of images, let a unit of positive charge be
distributed on $(0, \infty )$ according to the probability measure
$\mu$, while a unit of negative charge is distributed on $(-\infty
,0)$ symmetrically with respect to reflection in the origin. The term
$\log \vert s+t\vert$ arises from the mutual attraction of charge of
unlike polarity on opposite sides of the origin; whereas 
$-\log \vert s-t\vert$ arises from the repulsion of charges of like
polarity on the same side of the origin. Also, $w$ represents an
electric field. Then 
$$E_n (\mu )=n\int_0^\infty w(t)\mu (dt)+n^2\int\int_{s\neq
t}\log\Bigl\vert {{s+t}\over{s-t}}\Bigr\vert \mu (ds )\mu (dt)\eqno(4.23)$$
\noindent represents the total energy associated with $n\mu$ on $(0, \infty
)$. 

\indent In particular, 
let $\mu^\sigma={{1}\over{n}}\sum_{j=1}^n \delta_{\sigma_j}$ be the
empirical distribution of the unordered points $\sigma
=(\sigma_j)_{j=1}^n$ on $[0,\infty
)$. Then there exists $Z_n>0$ such that 
$dP=Z_n^{-1}\exp (-E_n(\mu^\sigma ))d\sigma_1\dots
d\sigma_n$ gives a probability measure on $(0, \infty )^n$, which
gives the distribution of $\sigma$. Then $\sigma\mapsto\mu^\sigma$
defines a random measure on $(0, \infty )$, which has characteristic 
functional $\Xi (g)=\int \exp (i\int
g(t)\mu^\sigma (dt)) dP$ for any continuous and bounded real function
$g$ on $(0, \infty )$. With $g(t)=nxt+n\theta (t)$, the integrand of () becomes
$$\Xi (g)=Z_n^{-1}\int_{(0,\infty )^n}\exp \Bigl( ix\sum_{j=1}^n
\sigma_j+i\sum_{j=1}^n\theta (\sigma_j)\Bigr) \prod_{j=1}^n {{\vert
f(\sigma_j)\vert}\over{\sigma_j}} \prod_{j<k}\Bigl(
{{\sigma_j-\sigma_k}\over{\sigma_j+\sigma_k}}\Bigr)^2 d\sigma_1\dots
d\sigma_n.\eqno(4.24)$$
\indent (ii) The preceding results of this section apply in
particular when $A$ is a finite matrix such that all the
eigenvalues have $\Re \lambda_j>0$.\par

\noindent In the statement of the following result, we use De Branges's terminology
from [16'], in which (H1),(H2) and (H3) characterize Hilbert spaces of entire
functions.\par  
\vskip.05in
\noindent {\bf Proposition 4.7} {\sl There exists a Hilbert space ${\cal H}(F^*)$ of entire functions such
that 
$$\pi^{-1}K(z,w)={{f(z)g(\bar w)-f(\bar w)g(z)}\over{\pi (z-\bar
w)}}\eqno(4.27)$$
\noindent is the reproducing kernel;\par
\indent (H1) there is an isometry from $\{ k\in {\cal H}(F^*): k(w)=0\}$ to
${\cal H}(F^*): k(z)\mapsto k(z)(z-\bar w)/(z-w)$ for all $w\in {\bf C}\setminus
{\bf R};$\par
\indent (H2) for all $w\in {\bf C}\setminus {\bf R}$, there is a bounded linear
functional $k\mapsto k(w)$; \par
\indent (H3) the map $k\mapsto k^*$ is an isometry on ${\cal
H}(F^*)$.}\par
\vskip.05in

\noindent {\bf Proof.} Let
$F^*(z)=\overline{F(\bar z)}$, so that $f(z)=(F(z)+F^*(z))/2$ and
$g(z)=(F(z)-F^*(z))/(2i),$ and hence the integrable kernel of () is given
by 
$${{f(z)g(\bar w)-f(\bar w)g(z)}\over{\pi (z-\bar
w)}}={{F^*(z)F(\bar w)-F(z)F^*(\bar w)}\over{2\pi i(z-\bar w)}}.\eqno(4.28)$$
\noindent
The weighted Hardy space $H^2(F^*)=\{ k=F^*h: h\in H^2\}$ is defined to be the range of multiplication by $F^*$ on
$H^2$, and the inner product is $\langle k_1,k_2\rangle =\int_{-\infty}^\infty
k_1(x)\overline{k_2(x)} dx/\vert F^*(x)\vert^2 .$
Then $\theta (s)=F(s)/F^*(s)$ is a meromorphic function with zeros $(\sigma_n)$ and poles $(\bar \sigma_n)$, 
and $\theta$ is an inner function; so that, $\vert \theta (z)\vert <1$ for all $\Im z>0$ and 
 $\vert \theta (x)\vert
=1$ for all almost $x\in {\bf R}$ with respect to Lebesgue measure.

 The multiplication $k\mapsto \theta k$ is an isometry on $H^2(F^*)$ and we let ${\cal H}(F^*)$
be the orthogonal complement of the range $\theta H^2(F^*)$ in $H^2(F^*).$ Then ${\cal H}(F^*)$
is a reproducing kernel Hilbert space with reproducing kernel $\pi^{-1}K(z, \zeta )$,
so $k(\zeta )=\pi^{-1}\langle k(z), K(z, \zeta )\rangle$; and one can extend the formula to
define an entire function of $z$.\par
\rightline{$\square$}

\indent Let $a(x)$ and $b(x)$ be meromorphic functions on ${\bf C}$. 
In [] we considered the functional equation
$$\alpha (x+y)={{a'(x)-a'(y)}\over{a(x)-a(y)}},\qquad \beta (x+y)={{b(x)-b(y)}\over{a(x)-a(y)}},$$
\noindent which implies $a'''(x)=Aa'(x)$ for some constant $A$, hence $a$ is a quadratic polynomial, trigonometric, or hyperbolic. The differential equation 
$$a(x){{d}\over{dx}}\left[\matrix{f(x)\cr g(x)\cr}\right] =\left[\matrix{0&1\cr \lambda a(x)-a(x)b(x)&0\cr}\right] \left[\matrix{f(x)\cr g(x)\cr}\right]$$
\noindent is of Tracy--Widom type, for a suitable change of variables, 
and is equivalent to $Lf=0$ where 
$Lf=-(d/dx)a(x)(d/dx)f(x)+b(x)f(x)+\lambda f(x)$. 
There is an associated differential equation 
$\phi''(x)+\alpha (x)\phi' (x)-\beta (x)\phi (x)=0$ such that the 
Hankel operator $\Gamma_\phi$ commutes with $L$. In each of these cases, 
() reduces to the 
hypergeometric equation or the confluent hypergeometric equation.
 Whereas this analysis introduces Hankel operators that arise in
 random matrix theory, it does not deal with the elliptic case, 
and involves rather contrived operators in the trigonometric case.
 In the present paper, we address this issue by introducing periodic
 linear systems in section 8, including elliptic systems in section 11.\par

\indent Suppose that $a:(0, \infty )\rightarrow {\bf R}$ is strictly increasing and twice continuously differentiable, and $b:(0, \infty )\rightarrow {\bf R}$ is continuous, and that there exist real functions $\alpha$ and $\beta$ such that 
$$\alpha (x+y)={{a'(x)-a'(y)}\over{a(x)-a(y)}}, \qquad \beta (x+y)={{b(x)-b(y)}\over{a(x)-a(y)}}.\eqno(5.12)$$
\noindent In [], we investigated this condition, and considered the important examples $a(x)=e^{-\kappa x}$ and $a(x)=a_2x^2+a_1x+a_0$, which arise in random matrix theory. Now consider the system
$$a(x){{d}\over{dx}}\left[\matrix {f(x)\cr g(x)\cr}\right]=\left[\matrix{0&1\cr a(x)b(x)-\lambda a(x)& 0\cr}\right]\left[\matrix{f(x)\cr g(x)\cr}\right],\eqno(5.13)$$
\noindent which is equivalent to the second order equation $Lf=\lambda f$, where $Lf(x)=-{{d}\over{dx}}a(x){{df}\over{dx}}+b(x)f(x).$ We introduce the kernel
$$K(x,y)=a'((x+y)/2){{f(x)g(y)-f(y)g(x)}\over{a(x)-a(y)}}\qquad (x,y>0).\eqno(5.14)$$
\vskip.05in
\noindent {\bf Theorem 5.6} {\sl Suppose that\par
\indent (i)  $a$ and $\beta$ are linear combinations of $x^je^{-\kappa x}$ for $j=0, 1, 2, \dots $
 and $c\in {\bf R}$;\par 
\indent (ii) $f$ and $g$ satisfy (), and\par
\indent (iii) there exists $\varepsilon >0$ such that 
$$e^{\varepsilon x}a'(x/2)f(x),\quad  e^{\varepsilon x}\beta (x)a'(x/2)f(x),\quad {\hbox{ and}} \quad
e^{\varepsilon x}a'(x/2)g(x)/a(x)$$
\noindent are all in $L^2(0, \infty)$ and $K(x,y)\rightarrow 0$ as $x,y\rightarrow\infty$.\par
\noindent Then there exist a finite-dimensional Hilbert space $H$ and $F\in L^2((0,\infty ); H)$ and $G\in L^2((0,\infty ); H')$ such that} $$K=\Gamma_F\Gamma_G.\eqno(5.15)$$
\vskip.05in
\noindent {\bf Proof.} We observe that from the differential equation
$$\eqalignno{\Bigl( {{\partial}\over{\partial x}}+{{\partial}\over{\partial y}}\Bigr) K(x,y)&={{a''((x+y)/2)}\over{a'((x+y)/2)}}K(x,y)-{{a'(x)-a'(y)}\over{a(x)-a(y)}}K(x,y)\cr
&+a'((x+y)/2){{b(y)-b(x)}\over{a(x)-a(y)}}f(x)f(y)-a'((x+y)/2){{g(x)g(y)}\over{a(x)a(y)}}\cr
&=-a'((x+y)/2)\beta (x+y)f(x)f(y) -a'((x+y)/2){{g(x)g(y)}\over{a(x)a(y)}},&(5.16)\cr}$$
\noindent where we have used the special properties of () at the last step, especially $\alpha (2x)=a''(x)/a'(x)$.\par 
\indent Suppose that $\phi :{\bf R}\rightarrow {\bf R}$ is a linear combination of $x^je^{-\kappa x}$ for $j=0, 1, \dots$ and $\kappa \in {\bf R}$. Then there exists a linear system $(-A_0, B_0, C_0)$ with finite-dimensional state space such that $\phi (x)=C_0e^{-A_0x}B_0$; indeed, $\phi (x)$ satisfies some nontrivial linear differential equation $\sum_{j=0}^N c_j \phi^{(j)}(x)=0$ with $c_j\in {\bf R}$, so one can find $(-A_0, B_0, C_0)$ by considering the Jordan canonical form of those matrices that satisfy $\sum_{j=0}^N c_j(-A_0)^j=0$.  Observe that the eigenvalues of $A$ can be chosen from the values of $\kappa$. In particular, we introduce $(-A,B,C)$ with state space $H_0$ and $(-A',B',C')$ with state space $H$ such that $\beta (x)=Ce^{-xA}B$ and $a'(x/2)=C'e^{-xA'}B'$, where the eigenvalues of $A$ and $A'$ are chosen from the $\kappa$ where $x^je^{-\kappa x}$ is a  summand on $\beta (x)$ or $a'(x/2)$; hence we have $\Vert e^{-xA}\Vert =O(e^{\varepsilon x/4}\vert \beta (x)\vert )$ as $x\rightarrow\infty$, while $\Vert e^{-A'x}\Vert =O(e^{\varepsilon x/4}a'(x/2))$ as $x\rightarrow\infty$ . We then introduce the Hilbert space $H=H_1\oplus (H_0\otimes H_1)$ and
 $F\in L^2((0,\infty ); H)$ and $G\in L^2((0,\infty ); H')$ by
$$F(x)=\left[\matrix{ {{g(x)}\over{a(x)}}C'e^{-xA'} & f(x)Ce^{-xA}\otimes C'e^{-xA'}\cr}\right]\eqno(5.17)$$
\noindent and  
$$ G(x)=\left[\matrix{ e^{-xA'}B'{{g(x)}\over{a(x)}}\cr e^{-xA}B\otimes e^{-xA'}B'f(x)\cr}\right],\eqno(5.18)$$
\noindent so that 
$$\eqalignno{F(x)G(y)&=C'e^{-(x+y)A'}B'{{g(x)g(y)}\over{a(x)a(y)}}+Ce^{-(x+y)A}BC'e^{-(x+y)A'}B'f(x)f(y)\cr
&=a'((x+y)/2){{g(x)g(y)}\over{a(x)a(y)}}+a'((x+y)/2)\beta (x+y)f(x)f(y).&(5.19)\cr}$$
We deduce that 
$$a'((x+y)/2){{f(x)g(y)-f(y)g(x)}\over{a(x)-a(y)}}=\int_0^\infty F(x+s)G(s+y)ds;\eqno(5.20)$$
\noindent indeed, both sides of this claimed equation converge to zero as $x\rightarrow\infty$ or $y\rightarrow\infty$, and have equal values when we apply ${{\partial}\over{\partial x}}+{{\partial }\over{\partial y}}.$
\par
\rightline{$\square$}\par
\vskip.05in
\indent In the following section, we obtain a more 
refined version of Theorem 5.6 for rational 
differential equations. Whereas the conditions (5.14) have 
trigonometric solutions, the statement of Theorem 5.6 needs modifying before one can 
adapt if to periodic systems of functions, as we see in sections 8 and 9. In section 12,
we obtain a version of Theorem 5.6 for elliptic functions.\par

\noindent {\bf 9. Hill's determinant}\par
\vskip.05in

\indent In his famous paper on periodic differential equations, Hill introduced various determinants $D, \nabla$ and $\square$, and here we reconcile $\nabla$ with the $\tau$ of a specially chosen periodic system. Let $H=L^2({\bf T}; {\bf C})$ have orthonormal basis $(e^{ijx})_{j=-\infty}^\infty$, and the standard convolution operation $\ast$, and let ${\cal D}(A) =\{ f\in H: f'\in H\}$ be the usual H\"older space, then let $H_0=L^2({\bf T}; {\bf C}^{3\times 1})$. Suppose that $q\in H$ is real-valued and $\pi$ periodic with $q'\in H$, and has Fourier coefficients $(q_j)_{j=-\infty}^\infty$. Hill's operator is $-d^2/dx^2+q(x)$ on $\{f\in H; {\bf C}):f''\in H\}$  and the spectrum is the periodic spectrum of Hill's equation $-f''+qf=\lambda f$. The periodic spectrum consists of the principal series for which the eigenfunctions are $\pi$-periodic, and the complementary series for which the eigenfunctions anti periodic in the sense that $f(x+\pi )=-f(x)$. \par
\indent To express this in terms of linear systems,  we introduce $b_{\pm}\in H$ by
$$b_+(x)=\sum_{j=-\infty}^\infty {{je^{ijx}}\over{4j^2-1}},\qquad b_-(x)=\sum_{j=-\infty}^\infty {{e^{ijx}}\over{4j^2-1}};\eqno(8.25)$$
\noindent then  we introduce a linear system $(-A,B_{\rho ,\lambda},C_t;E_{\rho,\lambda ,t})$ which satisfies some of the properties of a periodic linear system, where the operators are
$$\eqalignno{ A:{\cal D}(A)\rightarrow H: & f(x)\mapsto -f'(x);\cr
B_{\rho, \lambda }: H_0\rightarrow H:& \left[\matrix{f\cr g\cr h\cr}\right]\mapsto 2(4\rho^2+1-4\lambda )b_+\ast f+ 4\rho b_-\ast f+8b_+\ast g-4b_-\ast h, \qquad (\lambda \in {\bf C})\cr
C_t:H\rightarrow H_0:&f(x)\mapsto \left[\matrix{if(x)\cr iq(x+t)f(x)\cr iq'(x+t)f(x)\cr}\right], \qquad (t\in {\bf R}).&(8.26)\cr}$$
\noindent Let $E_{\rho, \lambda ,t }:H\rightarrow H$ be the operator with matrix
$$[E_{\rho, \lambda ,t}(j,k)]_{j,k=-\infty}^\infty =\Bigl[ {{4q_{j-k}e^{i(j-k)t}+(8j\rho+4\rho^2+1-4\lambda )\delta_{jk}}\over{4j^2-1}}\Bigr]_{j,k=-\infty}^\infty .\eqno(8.27)$$
\noindent{\bf Definition} Suppose that $E$ is a Hilbert--Schmidt operators with eigenvalues $\mu_j$, listed according to multiplicity; then Carleman's determinant is $\det_2(I+E)=\prod_{j=1}^\infty (1+\mu_j)e^{-\mu_j}$.  Refining Hill's definition [], we introduce 
$\nabla_2(\rho ,\lambda ,t)=\det_2(I+E_{\rho, \lambda ,t}).$\par
\vskip.05in
\noindent {\bf Proposition 9.1} {\sl Consider Hill's equation with potential $q(x+t)$.\par
 \indent (i) Then $e^{2\pi i\rho}$ is a Floquet multiplier, if and only if $\nabla_2(\rho ,\lambda ,t)=0.$}\par 
\indent {\sl (ii) The operators $C_t$ are uniformly bounded for $t\in {\bf R}$ , $B_{\rho,\lambda}$ and $E_{\rho,\lambda ,t}$ are Hilbert--Schmidt and 
$$AE_{\rho,\lambda ,t} +E_{\rho, \lambda ,t}A=B_{\rho, \lambda} C_t+4i\rho I.\eqno(8.28)$$
\indent (iii) Suppose $\rho=0$. Then $E_{0,\lambda ,t}$ is trace class,  $\tau (x,t,\lambda )=\det (I+e^{-xA}E_{0,\lambda ,t}e^{-xA})$ is $\pi$-periodic in $x$ and $t$, and defines an entire function of $\lambda$; the set of zeros of $\tau (0,t,\lambda )$ coincides with the  periodic spectrum.}\par
\vskip.05in
\noindent {\bf Proof.}  (i) Consider $g(x)=e^{i\rho x}f(x)$, where $f(x)=\sum_{j=\infty}^\infty f_je^{ijx}$ is a $2\pi$ -periodic function. Then Hill's equation for $g$ reduces to the system of linear equations for the Fourier coefficients of $f$, namely 
$$f_j+{{8j\rho+4\rho^2+1-4\lambda }\over{4j^2-1}}f_j+4\sum_{k=-\infty}^\infty{{q_{j-k}e^{i(j-k)t}}\over{4j^2-1}}f_k=0;\qquad (j\in {\bf Z});\eqno(8.29)$$
\noindent introducing $F=(f_k)_{k=\infty}^\infty$, we write this as $(I+E_{\rho, \lambda ,\rho ,t})F=0$. We observe  $\sum_{j=-\infty}^\infty \vert q_j\vert $ converges, so that $E_{\rho ,\lambda ,t}$ is a Hilbert--Schmidt operator, hence $I+E_{\rho, \lambda ,0}$ is invertible, unless $\nabla_2 (\lambda ,\rho )=0$.  In which case, there exists $F_\lambda \in H$ non zero and such that $(I+E_{\lambda ,0} )F_\lambda =0$, hence the system of equations ()  has a solution with $G\in \ell^2$, so $g(x)$ gives a solution of Hill's equation such that $g(x+2\pi )=e^{2\pi i\rho}g(x).$\par
 \indent  (ii) The functions $q$ and $q'$ are bounded, so $C_t$ is bounded. Next one computes the Fourier components of $AE_{0,\lambda ,t}+E_{0,\lambda ,t}A$; a key step in the computation is the formula
$$\eqalignno{(j+k)E_{\rho,\lambda ,t}(j,k)&=4\rho \delta_{jk}+ {{4\rho\delta_{jk}}\over{ 4j^2-1}}+ {{8jq_{j-k}e^{i(j-k)t}}\over{4j^2-1}}\cr
&\quad -{{4(j-k)q_{j-k}e^{i(j-k)t}}\over{4j^2-1}}+{{(8\rho^2+2-8\lambda )j\delta_{jk}}\over{4j^2-1}}&(8.30)\cr}$$
\noindent which gives the decomposition corresponding with $q_t(x)=q(x+t)$ to
$$B_{\rho,\lambda}C_tf=4i\rho b_-\ast f+ 8ib_+\ast (q_tf)- 4ib_-\ast (q'_tf)+ 2i(4\rho^2+1-4\lambda) b_+\ast f.\eqno(8.31)$$ 
\indent (iii) Suppose $\rho=0$. Then the diagonal entries of $E_{0,\lambda ,t}$ are now summable so \par
\noindent $\sum_{j,k=-\infty}^\infty\vert E_{0,\lambda ,t}(j,k)\vert $ converges, and hence  $E_{0,\lambda ,t}$ is trace-class. Consequently, \par
\noindent $\det (I+e^{-xA}E_{0,\lambda ,t}e^{-xA})$ is well defined, and defines an entire 
 function of $\lambda$. Suppose that $\det (I+E_{0,\lambda ,t})=0$ for 
some $\lambda\in {\bf C}$. Then the system of equations () has a 
nonzero solution given by  $F_\lambda =(f_j)_{j=-\infty}^\infty
\in\ell^2({\bf Z}).$  We deduce that $f_\lambda (x )=\sum_{j=-\infty}^\infty f_je^{ijx}$ satisfies $-f_\lambda''(x)+q(x+t)f_\lambda (x)=\lambda f_\lambda (x)$, so $f_\lambda$ is a periodic eigenfunction, so $\lambda$ is in the periodic spectrum of Hill's equation.\par
\rightline{$\square$}

\indent Let $U(x,t,\lambda )$ be the fundamental solution of 
$${{d}\over{dx}}U=\left[\matrix{ 0&1\cr q(x+t)-\lambda & 0\cr}\right]
U; U(0,t,\lambda )=I_2.$$
\noindent Let $y_2(x,t,\lambda )$ be the top right entry of
$U(x,t,\lambda 0$, which satisfies 
$$y_2(x,t,\lambda )={{\sin \sqrt{\lambda}
x}\over{\sqrt{\lambda}}}+\int_0^x {{1}\over{\sqrt{\lambda}}}\sin
\sqrt{\lambda}(x-s)q(s+t) y_2(s,t,\lambda )\, ds.$$
We introduce the linear system with state space $H$ and ${\cal D}(A)=\{ f\in L^2[0, \pi
]: f'\in L^2[0, \pi ]\}$ given by 
$$\eqalignno{ A:{\cal D}(A)\rightarrow H:& f(x)\mapsto -f'(x);\cr
              B:{\bf C}\rightarrow {\cal D}(A): & b\mapsto
y_2(x,t,\lambda )b;\cr
C:{\cal D}(A)\rightarrow {\bf C}:& f\mapsto f(0).\cr}$$
\noindent Then $\phi (x)=Ce^{-xA}B$ reduces to $\phi
(x)=y_2(x,t,\lambda )$ and $R=\int_0^\pi e^{-sA}BCe^{-sA}\, ds$ reduces
to
$$Rg(x)=\int_0^\pi y_2(x+s,t,\lambda )g(s)\, ds\qquad (g\in L^2[0, \pi
]).$$
\noindent We also introduce ${\cal D}(Q)=\{f\in L^2[0, \pi ]: f',
f''\in L^2; f(0)=f(\pi )=0\}$ and consider $Q^{(t)} f(x)=-f''(x)+q(x+t
)f(x)$. Note that the periodic eigenvalues of Hill's equation are
fixed irrespective of $t$, whereas the spectrum of $Q^{(t)}$
depends upon $t$. Indeed, $Q^{(t)}$ is self-adjoint in ${\cal D}(Q)$ and the
eigenvalues $\mu^{(t)}_1<\mu^{(t)}_1<\dots$ give the tied spectrum
namely the double periodic eigenvalues $\lambda_{2j-1}=\lambda_{2j}$
together with the auxiliary spectrum with terms lying between simple
periodic eigenvalues $\lambda_{2j-1}\leq\mu^{(t)}_j\leq\lambda_{2j}$.
As $t$ ranges over $[0, \pi ]$, the eigenvalue $\mu_j^{(t)}$ describes
$[\lambda_{2j-1}, \lambda_{2j}]$ $j$ times. 
Observe that $y_2(\pi,t,\lambda )=0$ if and only if $\lambda$ belongs
to the tied spectrum, so we can introduce
$$m(t, \lambda )={{y_2(\pi ,t, \lambda )}\over{y_2(\pi ,0,\lambda
)}}=\prod_{j=1}^\infty {{(1-\lambda/\mu_j^{(t)})}\over{(1-\lambda
/\mu_j^{(0)})}}.$$
\vskip.05in
\noindent {\bf Proposition}{\sl (i) Let $\lambda =\mu^{(t)}$ belong to
the tied spectrum of Hill's equation. Then $R$ is a periodic Hankel
operator.\par
\indent (ii) Let $\lambda$ be a simple periodic eigenvalue of Hill's
equation with eigenfunction $f_\lambda (x)$. Then $m(x,\lambda )$ is
proportional to $f_\lambda (x)^2$.}\par
    
{\bf Proof.} (i) In this case, $y_2(x,t, \mu^{(t)})$ is periodic with
period $\pi$ or $2\pi$.\par
\indent (ii) This follows from Lemma  of [McKvM].\par

\indent Let $(\lambda_n)_{n=0}^\infty$ be the periodic spectrum of Hill's equation, with no repeats, and for $\eta >0$ consider 
$$\sigma_n=\cases{ i\eta +\sqrt{\lambda_{n-1}}, & for $n=1,2, \dots$\cr
                               i\eta, & for $n=0$;\cr
                             i\eta-\sqrt{\lambda_{-n-1}},& for $n=-1,-2,\dots $\cr}$$

Then $(e^{i\sigma_{2n}x})$ and $(e^{i\sigma_{2n+1}x})$ gives Riesz sequences in $L^2(0,
2\pi )$. Whereas $(\sigma_n)$ does not satisfy condition $(C)$, there
is a one sided bound
$$\bigl\Vert\sum_{j=-\infty}^\infty a_ne^{i\sigma_nx}\bigr\Vert_{L^2(),
2\pi )}\leq C\Vert (a_n)\Vert_{\ell^2}\qquad ((a_n)\in \ell^2).$$

\indent For a periodic potential, we consider the Tracy--Widom kernel of exponential type [p46,of CMP163] , and obtain a result analogous to Theorem 6.1.\par
\vskip.05in
\noindent {\bf Theorem 9.2} {\sl Let the potential of Hill's equation be a rational function $u$ of $e^{2ix}$, which is real and with no poles for $x\in {\bf R}$; suppose that the discriminant satisfies $\Delta (\lambda )^2>4$ for some $\lambda\in {\bf R}$. Then there exists $-1<\mu <1$ and a non zero solution $g(x)=e^{i\rho x}f(x)$ of Hill's equation, where $e^{2i\pi \rho}=\mu$ and $f$ is $2\pi$-periodic. The kernel
$$K_\lambda (x,y)=e^{i\rho (x+y)} {{f(x)f'(y)-f'(x)f(y)}\over{\sin (x-y)}}$$
\noindent defines a trace class operator on $L^2(0,\infty )$, and there exist $N\times N$ matrices $\Phi$ and $\Psi$ such that:\par
\indent (i) the entries of $\Phi$ and $\Psi$ belong to ${\bf C}(\sin x,\cos x, g(x),g'(x))$  and satisfy $\Phi (x+2\pi )=\mu \Phi (x)$, $\Psi (x+2\pi )=\mu\Psi (x)$;\par
\indent (ii) the Laplace transform of $g$ is meromorphic with only possible poles at $s=2ij+\pi^{-1}\log\mu$ for $j\in {\bf Z}$;\par
\indent (ii) $\Gamma_{\Phi}$ and $\Gamma_{\Psi}$ are Hilbert--Schmidt operators on $L^2((0,\infty ); {\bf C}^N)$ and such that}
$$\det (I+K_\lambda )=\det (I+\Gamma_{\Phi}\Gamma_{\Psi }).$$
\vskip.05in
\noindent {\bf Proof.} (i)  By Floquet's theorem, there exists a solution $g$ of this form for each $\lambda$ in an interval of instability. Observe that $g(x+2\pi )=\mu g(x)$, so that $g\in L^2(0,\infty )$ and $K_\lambda (x+2\pi ,y)=K_\lambda (x,y+2\pi )=\mu K_\lambda (x,y)$, so $K_\lambda (x,y)\rightarrow 0$ as $x\rightarrow\infty$ or $y\rightarrow\infty$. \par
\indent (ii) Let $s_j=\pi ij+\pi^{-1}\log \mu$, so that $e^{\pi s_j}=\mu$ for all $j\in {\bf Z}$. Then using the fact that $g$ is a Floquet solution, we compute the Laplace transform
$$\eqalignno{{\cal L}(g; s)&=\int_0^\infty e^{-sx}g(x)\, dx\cr
&=\sum_{k=0}^\infty\int_{k\pi}^{(k+1)\pi}e^{-sx}g(x)\,dx\cr 
&={{\int_0^{\pi}e^{-s(x-\pi /2)}g(x)dx}\over { 2\sqrt{\mu}\sinh ((s-s_0)\pi/2)}}\cr}$$
\noindent where the numerator is an entire function of exponential type and growth rate $\pi/2$, and the denominator is entire and vanishes only at the points $s_j$, all of which lie on a bilateral arithmetic progression in the open left half plane.\par 
\indent (iii) By a simple calculation, we have  
$$\Bigl({{\partial }\over{\partial x}}+{{\partial }\over{\partial y}}\Bigr)K_\lambda (x,y)={{u(e^{2ix})-u(e^{2 iy})}\over{\sin (x-y)}}e^{i\rho x}f(x)e^{i\rho y}f(y),$$
\noindent where $u(e^{2ix})=q(x)/p(x)$ where $p(x)$ and $q(x)$ are trigonometric polynomials, so 
$${{u(e^{2ix})-u(e^{2iy})}\over{\sin (x-y)}}={{q(x)-q(y)}\over{p(x)\sin (x-y)}}-{{(p(x)-p(y))q(y)}\over{p(x)p(y)\sin (x-y)}}.$$
\noindent Suppose momentarily that $p(x)=1$, so that only the first difference quotient appears. Then writing $q(x)=\sum_{j=-N}^Na_je^{2ijx}$ for some $a_j\in {\bf C}$, we can introduce $\phi_j(x)=e^{i\rho x+ijx}f(x)$  and write this as 
$$\eqalignno{{{q(x)-q(y))}\over{\sin (x-y)}}e^{i\rho x}f(x)e^{i\rho y}f(y)&=\sum_{j=1}^N\sum_{k=0}^{j-1}2ia_j\phi_{2k+1}(x)\phi_{2(j-k)-1}(y)\cr
&\quad +\sum_{j=-N}^{-1}\sum_{k=0}^{-j-1}2ia_j\phi_{-2k+1}(x)\phi_{2j+2k+1}(y)\cr}$$
\noindent We note that this equals
$$\eqalignno{-\Bigl( {{\partial }\over{\partial x}}+{{\partial }\over{\partial y}}\Bigr)\int_0^\infty \Bigl(& \sum_{j=1}^N\sum_{k=0}^{j-1}2ia_{j}\phi_{2k+1}(x+s)\phi_{2j-2k-1}(y+s)\cr
&+\sum_{j=-N}^{-1}2ia_j\phi_{-2k+1}(x+s)\phi_{2j+2k+1}(y+s)\Bigr)ds\cr}$$
\noindent where $\phi_j(x)=e^{i\rho x+ijx}f(x)$ has $\int_0^\infty x\vert \phi_j(x)\vert^2dx$ convergent. We deduce that 
$$\eqalignno{K_\lambda (x,y)&=-\sum_{j=1}^N\sum_{k=0}^{j-1}\int_0^\infty 2ia_j\phi_{2k+1}(x+s)\phi_{2j-2k-1}(y+s)ds\cr
&+\sum_{j=-N}^{-1}\sum_{k=0}^{-j-1}\int_0^\infty 2ia_j\phi_{-2k-1}(x+s)\phi_{2j+2k+1}(s+y)ds.\cr}$$
\noindent since the left and side and the right-hand side differ by a function $h(x-y)$, and both the  left-hand side and right-hand side converge to zero as $x\rightarrow\infty$ or $y\rightarrow\infty$.
Finally, we build the
$N\times N$ matrices
$$\eqalignno{\Phi (x)&=\left[\matrix{\phi_{-2N+1}(x)&\phi_{-2N+3}(x)&\dots &\phi_{2N-1}(x)\cr 0&0&\dots  &0\cr
\vdots&\vdots &\ddots &\vdots\cr 0&0&\dots &0\cr}\right],\cr 
\Psi (y)&=\left[\matrix{2ia_{-2N+1}\phi_{-2N+1}(y)&0&\dots &0\cr
                          2ia_{-2N+3}\phi_{-2N+3}(y)&0&\dots  &0\cr
                          \vdots&\vdots&\ddots &\vdots\cr
                            2ia_{2N-1}\phi_{2N-1}(y)&0&\dots
&0\cr}\right]&(8.8)\cr}$$
\noindent so that $\Gamma_\Phi$ and $\Gamma_\Psi$ are Hilbert--Schmidt
matrix operators, so that $K_\lambda (x,y)$ is the top left entry of $\Gamma_\Phi\Gamma_\Psi$, and all other entries in the matrix kernel are zero. \par
\indent The general case of $u=q/p$ is similar, with a modified definition of $\Phi$ and $\Psi$.\par
\rightline{$\square$}\par 

\vskip.05in
   
\vskip.1in

\noindent In this section we consider a differential ring ${\bf S}$ associated with a linear system $(-A,B,C)$, and obtain solutions of the Korteweg-- de Vries equation and its higher order counterparts. First we introduce an infinite collection of time variables $t_j$ and consider solutions which are parametrized by the infinite torus.\par
\indent We introduce the infinite real torus
$${\bf T}^\infty =
\bigl\{ (t_{2j-1})_{j=1}^\infty\in {\bf R}^\infty : 
\lim\sup\vert t_{2j-1}\vert^{1/j}\rightarrow 0,\quad j\rightarrow\infty 
\bigr\}\eqno(9.1)$$
\noindent which forms a vector space under coordinate wise addition and multiplication by scalars. We regard the space variable $x$ as a time parameter $t_{-1}$.
 By means of  the map $(t_{2j-1})_{j=1}^\infty \mapsto \sum_{j=1}^\infty t_{2j-1}z^{j-1}$, we embed ${\bf T}^\infty$ as a linear subspace of the space ${\bf H}_{\bf C}$ of odd entire functions that are real on the real axis, where ${\bf H}_{\bf C}$  has the topology of uniform convergence on compact subsets of ${\bf C}$. 
The linear functionals $\ell_{2j-1}(f)=f^{(2j-1)}(0)/(2j-1)!$ are 
continuous, and any linear map 
$\alpha :{\bf C }^g\rightarrow {\hbox{span}}\{ \ell_j: j=1, 3, \dots \}$
 is continuous and has a continuous transpose 
$\alpha^t:{\bf H}_{\bf C}\rightarrow {\bf C}^g$. The matrix that 
represents $\alpha^t$, with respect to $(z^{2j-1})_{j=1}^\infty$ and 
the standard basis of ${\bf C}^g$, has finitely many non-zero entries.\par 
We regard $t_{2j-1}$ as deformation parameters, and ${\bf T}^\infty$ as a subgroup of the formal
 Lie group $\{ \exp (\sum_{j=1}^\infty t_{2j-1}\partial_{2j-1}): t_{2j-1}\in {\bf R}\}$  with 
$\partial_{2j-1}=\partial/\partial t_{2j-1}$. \par
As in [31], for $z\in {\bf C}$ with $\vert z\vert >1$ we introduce $\lfloor z\rfloor =(2z^{-2j+1}/(2j-1))_{2j-1}^\infty$, so $t\mapsto t+\lfloor z\rfloor $ is the Sato's shift on odd coordinates.\par 
For $v,t_0\in {\bf T}^\infty$, the equation $t=vs+t_0$ gives rectilinear motion.\par 
\indent Now we construct representations of ${\bf T}^\infty$. For a linear system 
$\Sigma =(-A,B,C)$, with $A$ as in the following result, we write 
$(t\cdot \tilde A)=\sum_{j=1}^\infty t_{2j-1}A^{2j-1}$ for 
$t=(t_{2j-1})_{j=1}^\infty \in {\bf T}^\infty$.\par  
\vskip.05in 
\noindent {\bf Lemma 9.1 } {\sl Suppose that either\par
\indent  (i) $A\in {\cal L}(H)$, or\par
\indent  (ii) $(e^{-xA})_{x\in {\bf R}}$ is a bounded and strongly continuous group of operators on $H$.\par
\noindent Then $U(t)=\exp (t\cdot \tilde A)$ for  $t\in {\bf T}^\infty$ defines a family of bounded linear operators on $H$ such that $U(t+s)=U(t)U(s)$ and 
$$U(\lfloor \zeta \rfloor) =(\zeta I+A)(\zeta I-A)^{-1}\qquad 
(\Re \zeta >x_0).\eqno(9.2)$$ 
\noindent As $t$ undergoes a rectilinear motion with velocity $v$, the resolvent operator satisfies}
$${{\partial R}\over{\partial v}}=(v\cdot \tilde A)R+R(v\cdot \tilde A).\eqno(9.3)$$
\vskip.05in  
\noindent {\bf Proof.}  (i) When $A\in {\cal L}(H)$, the series $\sum_{j=1}^\infty t_{2j-1}A^{2j-1}$ converges, and the addition rule is clear.
The elementary identity,
$$ (\zeta I+A)(\zeta I-A)^{-1}=\exp \Bigl(\sum_{j=1}^\infty
{{2A^{2j+1}}\over{ (2j+1)\zeta^{2j+1}}}\Bigr)=\exp\Bigl( 2\lfloor 1/\zeta \rfloor \cdot \tilde A\Bigr) ,\eqno(9.4)$$
\noindent where the series
$\sum_{j=1}^\infty A^{2j+1}/(2j+1)\zeta^{2j+1}$ converges 
for $\vert
\zeta\vert >\Vert A\Vert$, so we can use the left-hand side as a definition of the
right-hand side for all $\zeta$ outside the spectrum of $A$.  
\par
\indent (ii) By a theorem of Sz.-Nagy,  $(e^{-xA})_{x\in {\bf R}}$ is similar to a unitary group on $H$, so $A$ is similar to a skew symmetric operator $\hat A$. Let $\sigma$ be the spectrum of $i\hat A$ and $\rho$ the scalar spectral 
measure of $\hat A$, and $H_\xi$ a measurable family of Hilbert spaces so that $H=\int_\sigma H_\xi\rho (d\xi )$ is the direct integral decomposition given by the spectral theorem; and for $t=(t_{2j-1})_{j=1}^\infty \in {\bf T}$, let $U(t)=\exp({\sum_{j=1}^\infty t_{2j-1}A^{2j-1}})$  be the bounded linear operator that operates on $H_\xi$ as multiplication by $\exp({\sum_{j=1}^\infty t_{2j-1}(-i\xi )^{2j-1}})$; then $U(t)$ is similar to a unitary and $U(t)U(s)=U(t+s)$ for all $s,t\in {\bf T}.$ One can obtain the identity (9.1) from the spectral theorem.\par
\indent In either case, we note that $U(t)R_xU(t)$ gives the resolvent operator for $\Sigma (t)$, where $R_x$ is the resolvent operator for $\Sigma$. The differential equation (9.2) follows directly.\par
\rightline{$\square$}\par
\indent Given a representation $t\mapsto U(t)$ of ${\bf T}^\infty$ as bounded linear operators on $H$, and a 
uniformly periodic linear system $(-A, B,C; E)$, we introduce 
$$\Sigma (t)=\bigl( -U(t)AU(-t), U(t)B, CU(t); U(t)EU(t)\bigr)$$
\noindent which is also a periodic linear system with $\tau$ function $\tau (x,t)$.
 Then $s:\Sigma (t)\mapsto \Sigma (s+t)$ is a group action of ${\bf T}^\infty$ on the space of linear systems, and $\tau (x,t)\mapsto\tau (x, s+t)$ the flow on the corresponding $\tau$ functions.

\noindent {\bf Definition} Let $U(t)$ be as in Lemma 11.1 Then $U(t)$ is the odd deformation group of 
the linear system $\Sigma (0)=(-A, B_0, C_0)$ which is associated with the odd deformation parameters $t_{2j-1}$. We introduce the family of linear systems
$$\Sigma_\zeta (t)=(-A, U(\lfloor \zeta \rfloor )U(t)B, CU(t))
\qquad (t\in {\bf T}^\infty).\eqno(9.5)$$
Then the corresponding $\tau$ function is given explicitly by
$$\tau_\zeta  (x,t)=\det \Bigl(I+(\zeta I+A)(\zeta I-A)^{-1}\exp ({2xA+\sum_{j=1}^\infty
2t_{2j+1}A^{2j+1}})E\Bigr)\qquad (t\in {\bf T}^\infty )$$
is the tau function that arises as in Theorem 8.2. In contrast to Proposition 2.3(iii), $\tau_\zeta  (x,t)$ does depend upon $t$.
Given a linear map $\rho :{\bf C}^g\rightarrow {\bf
C}^\infty$ of rank $g$ such that $\rho (e_j)\in {\bf Z}^\infty$ has
only finitely many non-zero entries with respect to the standard bases,
then $\rho^t:{\bf C}^\infty \rightarrow {\bf C}^g$ satisfies
$\rho^t({\bf Z}^\infty )\subseteq {\bf Z}^g$. Then 
$\tau \circ \rho :{\bf C}^g\rightarrow {\bf C}$ is entire and periodic with
respect to ${\bf Z}^g$. \par
\vskip.05in 
\indent  We allow $C:H\rightarrow {\bf C}$ and 
$B:{\bf C}\rightarrow H$ to
evolve through time so that $C=C_0U(t)$ and $B=U(t)B_0$ for some initial
$C_0:H\rightarrow {\bf C}$ and $B_0:{\bf C}\rightarrow H$, and correspondingly
$R(x,t)=U(t)R_xU(t)$. The formulas involving 
$C,B$ and $R$ are symmetrical
with respect to time evolution, since $B$ and $C$ both evolve under 
the same group. In contrast to Corollary 9.3, we do not assume that $A$
commutes with $BC$; that $BC$ here will typically have rank one, whereas $A$
will have infinite rank.\par

\vskip.05in

\noindent {\bf Proposition 9.2} {\sl Suppose that $A$ is bounded and let 
$F_x=(I+R)^{-1}$. Then}
$$\eqalignno{{\cal A}&={\hbox{span}}_{\bf C}\Bigl\{A^{n_1}, 
A^{n_1}F_xA^{n_2}\dots F_xA^{n_r}:n_1, n_2, 
\dots , n_r\in {\bf
N}\Bigr\}&(9.7)\cr}$$
\noindent {\sl is a differential subring of} $C^\infty ((0, \infty )^2;
{\cal L}(H))$, {\sl and the map $\lfloor \, .\,\rfloor :{\bf
A}\rightarrow C^\infty ((0, \infty )^2; {\bf C})$ 
$$\bigl\lfloor P\bigr\rfloor=Ce^{-xA}F_xPF_xe^{-xA}B\eqno(9.8)$$
\noindent has range $\lfloor{\cal A}\rfloor$, where $\lfloor{\cal A}\rfloor$ is a differential ring with pointwise multiplication and
derivatives $\partial /\partial x$ and $\partial /\partial t_1$.}\par
\vskip.05in
\noindent {\bf Proof.} As in Lemma 3.2, the basic relations are
$$\eqalignno{{{\partial}\over{\partial
x}}\bigl\lfloor P\bigr\rfloor&=\Bigl\lfloor A(I-2F_x)P+{{\partial }\over{\partial
x}}P+P(I-2F_x)A\Bigr\rfloor ,&(9.9)\cr
{{\partial}\over{\partial t_3}}\bigl\lfloor P\bigr\rfloor&=\Bigl\lfloor A^3(I-2F_x)P+
{{\partial }\over{\partial
t_3}}P+P(I-2F_x)A^3\Bigr\rfloor ,\cr
\bigl\lfloor P\bigr\rfloor\bigl\lfloor Q\bigr\rfloor&=\bigl\lfloor
P(AF_x+F_xA-2F_xAF_x)Q\bigr\rfloor.&(9.10)\cr}$$
\indent Indeed it follows from the Lyapunov equation (1.10) that $F_x$ 
satisfies the differential equations
$${{\partial F_x}\over{\partial x}}=AF_x+F_xA-2F_xAF_x,\eqno(9.11)$$
$${{\partial F_x}\over{\partial
t_3}}=A^3F_x+F_xA^3-2F_xA^3F_x\eqno(9.12)$$
\noindent and hence the derivatives from the first and last
factors in (9.8) satisfy
$${{\partial }\over{\partial x}}Ce^{-xA}F_x=Ce^{-xA}F_xA(I-2F_x), 
\quad {{\partial }\over{\partial
x}}F_xe^{-xA}B=(I-2F_x)AF_xe^{-xA}B;\eqno(9.13)$$
$${{\partial }\over{\partial
t_3}}Ce^{-xA}F_x=Ce^{-xA}F_xA^3(I-2F_x);\quad {{\partial }\over{\partial
t_3}}F_xe^{-xA}B=(I-2F_x)A^3F_xe^{-xA}B.\eqno(9.14)$$
\noindent By applying Leibniz's rule, we deduce that 
$\lfloor{\cal A}\rfloor$ is closed under ${\partial }/{\partial x}$
and ${\partial }/{\partial t_3}.$ Furthermore
$$F_xe^{-xA}BCe^{-xA}F_x=AF_x+F_xA-2F_xAF_x,\eqno(9.15)$$
\noindent so $\lfloor{\cal A}\rfloor$ is closed under
multiplication, and the product rule (9.10) holds. 
\rightline{$\square$}\par
\vskip.05in

\vskip.05in
\indent We now point out some
particular solutions which are realised via Lemma 9.3, some of
which were also noted by P\"oppe [49]. Let $\lambda_j$ be distinct complex numbers 
for $j=1, \dots ,m$, such that $\Re\lambda_j >0$, and let 
$H={\hbox{span}}\{ x^{j}e^{-\lambda_\ell x}: j=0, \dots , 
n_\ell -1; \ell =1, \dots ,m\}$, which we view as a subspace of
 $L^2(0, \infty )$, and let $A=-{{d}\over{dx}}$ on $H$.\par 
\vskip.05in
\noindent {\bf Corollary 9.4} (Solitons) {\sl (i) Then $(e^{-sA})_{s\in {\bf R}}$ 
defines a $C_0$ group of operators on $H$ such that 
$\Vert e^{-sA}\Vert <1$ for $s>0$, and $\phi (x;t_3)$
\noindent satisfies} ${{\partial \phi}\over{\partial
t_3}}=2{{\partial^3\phi}\over{\partial
x^3}},$ {\sl and $u(x;t_3)\in {\bf C}(x,t_3,
e^{-\lambda_jx}, e^{-2\lambda_j^3t_3})$ satisfies $KdV$}.\par
{\sl (ii) In particular, suppose that $A$ has distinct and simple
eigenvalues, and that} $B_0=(b_j)\in {\bf C}^{n\times 1}$ {\sl and}
$C_0=(c_j)\in {\bf C}^{1\times n}$. {\sl Then} 
$$\det (I+\mu R_x)=\sum_{m=0}^N\mu^m\sum_{\sigma \subseteq \{ 1, 
\dots, N\}, \sharp
\sigma =m}\prod_{j\in \sigma } b_jc_je^{-2\lambda_j^3t_3-2\lambda_jx}\,
\prod_{j,k\in \sigma :j\neq k}
{{\lambda_j-\lambda_k}\over{\lambda_j+\lambda_k}}\eqno(9.24)$$
\vskip.05in
\noindent {\bf Proof.} (i) The group $e^{-sA}$ operates as translations
$e^{-sA}f(x)=f(x+s)$, and hence $e^{-sA}$ is a strict contraction
on the finite dimensional space $H={\bf C}^n$ for $s>0$. In effect, we have returned to
the setting of Proposition 5.4. The generator is $-A=d/dx$, and can introduce 
$A^3=-d^3/dx^3$ and the group $e^{-t_3A^3}$ which is associated 
with the linearized Korteweg--de Vries equation. By Lemma 9.3, $u$
satisfies the KdV equation (9.23), and by Theorem 5.2, $u$ is rational in
the basic variables.\par 
\indent (ii) Apply Proposition 5.4 and Lemma 9.3.\par
\rightline{$\square$}\par
\vskip.05in
\noindent Let $H=L^2(-\infty,
\infty )$ and as in section 6, we can take $Af(x)=-f'(x)$  and we note that
$e^{-tA^3}$ is the Airy group\par
$$e^{-tA^3}f(y)={{1}\over{2\pi}}\int_{-\infty}^\infty \hat f(\xi
)e^{-it\xi^3+iy\xi}\, d\xi .\eqno(9.25)$$
\noindent Then with $g\in {\cal D}(A^4)$ we choose 
$B_0:\alpha\mapsto g(y)\alpha$ and
$C_0:f\mapsto f(0)$, and let
$$\gamma_n=(-1)^n\int_{-\infty}^\infty \hat g(\xi ){{(i\xi
+1)^n}\over{(-i\xi +1)^{n+1}}}{{d\xi}\over{\pi}}.\eqno(9.26)$$ 
\vskip.05in
\noindent {\bf Corollary 9.5} (Non solitons) {\sl Let $g\in {\cal
D}(A^4)$ have $\sum_{n=0}^\infty (1+n)\vert \gamma_n\vert<1$. Then $\phi 
(x;t)$ satisfies (9.16) and 
$u(x;t)$ satisfies
(9.17).}\par
\vskip.05in
\noindent {\bf Proof.} By
Plancherel's formula we identify 
${\cal D}(A^4)=\{ g\in L^2: \int_{-\infty}^\infty
(1+\xi^8)\vert\hat g(\xi )\vert^2d\xi <\infty\}$, so the maps are well
defined. By a simple calculation, have 
$$R_xf(y)=\int_x^\infty g(y+s)f(s)\, ds,\qquad (f\in L^2(0, \infty
))\eqno(9.27)$$
\noindent so in particular $R_0$ is the Hankel integral operator with
kernel $g(y+s)$. Hence $R_0$ is unitarily equivalent to
$[\gamma_{j+k}]_{j,k=0}^\infty$ on $\ell^2$, which by the hypotheses is a
trace class operator; likewise, $R_x$ is trace class. Furthermore, 
$I+R_x$ is invertible, and the inverse $F$ is given by a Neumann series.
Given these facts, we can apply Lemma 9.3.\par
\rightline{$\square$}\par
\indent For $m\geq 4$, We can choose 
$g(x)={\bf I}_{(0, \infty )}(x)x^me^{-x}$ in
Corollary 9.4. Whereas the choice of $g(y)=\delta_0$ is technically inadmissible,
the resulting expression
$\phi (x;t)=t^{-1/3}{\hbox{Ai}}(-x/(6t)^{1/3})$ does give a
solution of (9.16).\par
\vskip.05in
\noindent {\bf Proposition 9.5} {\sl Suppose that 
$C_0A^6:H\rightarrow {\bf
C}$ and $A^5B_0:{\bf C}\rightarrow H$ are bounded.\par
\indent (i) Then the scattering function $\phi
(x;t_2)=C_0e^{-2t_2A^5-xA}B_0$ satisfies 
$${{\partial \phi}\over{\partial t_2}}=2{{\partial^5 \phi}\over{\partial
x^5}}.\eqno(9.28)$$
\indent (ii) Let $v(x)=T(x,x)$, so that
$$v(x,t)=-C_0e^{-xA-t_5A^5}(I+R)^{-1}e^{-xA-t_5A^5}B_0.\eqno(9.29)$$
\noindent Then $u(x,t_5)={{\partial v}\over{\partial x}}$ 
satisfies the KdV(5) equation}
$$16{{\partial u}\over{\partial t_5}}={{\partial^5u}\over{\partial
x^5}}+10u{{\partial^3u}\over{\partial x^3}}+20 {{\partial
u}\over{\partial x}}{{\partial^2
u}\over{\partial x^2}}+30u^2{{\partial
u}\over{\partial x}}.\eqno(9.30)$$
\vskip.05in
\noindent {\bf Proof.} We shall prove that
$$16{{\partial
v}\over{\partial t_2}}={{\partial^5
v}\over{\partial x^5}}+10{{\partial^3
v}\over{\partial x^3}}{{\partial
v}\over{\partial x}}+5\Bigl({{\partial^2
v}\over{\partial x^2}}\Bigr)^2+20\Bigl( {{\partial
v}\over{\partial x}}\Bigr)^6.\eqno(9.31)$$
\indent The basic identities required follow from (9.20), namely
$$\eqalignno{{{\partial^4 v}\over{\partial
x^4}}&=16\lfloor A^4\rfloor-64\lfloor A^3F_xA+AF_xA^3\rfloor-
96\lfloor A^2F_xA^2\rfloor\cr
&\quad +112\lfloor A^2F_xAF_xA+AF_xA^2F_xA+AF_xAF_xA^2\rfloor-
384\lfloor AF_xAF_xAF_xA\rfloor;&(9.32)\cr
{{\partial^5 v}\over{\partial
x^5}}&=32\lfloor A^5\rfloor-160\lfloor A^4F_xA+AF_xA^4\rfloor-320\lfloor A^3F_xA^2+A^2F_xA^3\rfloor\cr
&\quad +640
\lfloor A^3F_xAF_xA+AF_xA^3F_xA+AF_xAF_xA^3\rfloor\cr
&\quad +960\lfloor A^2F_xA^2F_xA+A^2F_xAF_xA^2+AF_xA^2F_xA^2\rfloor\cr
&\quad -1920\lfloor A^2F_xAF_xAF_xA+AF_xA^2F_xAF_xA+AF_xAF_xA^2F_xA+AF_xAF_xAF_xA^2\rfloor\cr
&\quad +3840
\lfloor AF_xAF_xAF_xAF_xA\rfloor.&(9.33)\cr}$$ 
\noindent Using these, one checks that (9.30) holds.\par
\rightline{$\square$}\par

\vskip.1in

\vskip.05in 
\noindent {\bf 11. The  Baker--Akhiezer function}
\vskip.05in
\indent Krichever [35] defines the Baker--Akhiezer function $\psi (x; {\bf p} )$ for spectral data $({\cal E}, {\bf p}_0, D,x_0)$ as follows. Given a compact Riemann surface ${\cal E}$ of genus $g$ with distinguished point ${\bf p}_0=\infty$,  let the local parameter at ${\bf p}_0$ be $z=\lambda^{-1}$ and suppose that ${\bf p}\mapsto \psi (x;{\bf p})$ is meromorphic on ${\cal E}\setminus \{{\bf p}_0\}$ with poles $\gamma_1, \dots , \gamma_g$ which are independent of $x=(x_1, \dots ,x_N)$; suppose also that there exist $\xi_j(x)$ for $j=1,2, \dots $, which are holomorphic for $x$ in a neighbourhood of $(x_{1,0}, \dots ,x_{N,0})$, such that there is a convergent series expansion
$$\psi (x,{\bf p})=\Bigl(1+\sum_{j=1}^\infty {{\xi_j(x)}\over{\lambda^j}}\Bigr)
\exp\Bigl( \sum_{k=1}^N(x_k-x_{k,0})\lambda^k\Bigr)
\qquad ({\bf p}\rightarrow {\bf p}_0)\eqno(11.1)$$ 
\noindent see [35, 54]. Given $x, x_0$ and the positive divisor $D=\sum_{j=1}^g k_j\delta_{{\gamma}_j}$, there exists a unique $\psi (x, {\bf p})$ on ${\cal E}\setminus \{ {\bf p}_0\}$ with these properties, up to a scalar multiple, and one can construct such a
function from quotients of Riemann's theta functions on the Jacobi variety for ${\cal E}$. 
Krichever showed that potentials which are finite gap or algebro-geometric on a complete algebraic curve give rise to an algebraic family in the sense of Lemma 5.3, and he constructed the meromorphic matrix function $W$ from the Baker--Akhiezer function and its derivatives. 
To deal with
commuting families of differential operators of rank greater than
one, he introduces matricial $\psi (x, \lambda )$ and matricial divisors in [36].\par
\indent In contrast, we define our $\psi_{BA}$ for linear systems, irrespective of whether there exists a suitable ${\cal E}$, and then seek to recover the properties of Krichever's definition. We start by looking at admissible linear systems, and then consider periodic linear systems. Then we introduce more parameters into the discussion, and consider $\psi_{BA}(x,t,\lambda )$ with $t$ in an infinite dimensional torus. Finally, we recover the case in which the linear system is associated with a compact algebraic spectral curve. \par

\indent In this section we consider scattering functions and the multiplication rule for tau functions  which is
analogous to the addition rule for positive divisors divisors on an 
algebraic curve. The multiplication $B\mapsto (\lambda
I-A)(\lambda I+A)^{-1}B$ is associated with adding the
divisor associated with a pole on the spectral curve. There is a
consequent formula for addition of divisors, which Ercolani and McKean  [21] 
credit to Darboux, as in Proposition 5.2. To subtract divisors, one considers the quotient of tau functions, as in the definition of the Baker--Akhiezer function.\par

\vskip.05in

\noindent {\bf Definition} (Baker--Akhiezer function) Given an admissible linear system
 $\Sigma_\infty =(-A,B,C)$
with tau function $\tau_\infty (x)=\det (I+\Gamma_{\phi_{(x)}})$ 
as in Proposition 2.2, we introduce 
$$\Sigma_\lambda =\Bigl(-A,(\lambda I+A 
)(\lambda I-A )^{-1}B, C\Bigr)\qquad (\Re \lambda >0)\eqno(11.2)$$
\noindent with tau function $\tau_\lambda (x)$, and the Baker--Akhiezer
function 
$$\psi_{BA}(x;\lambda )=\exp \bigl({\lambda x}\bigr){{\tau_\lambda (x)}\over 
{\tau_\infty (x)}}.\eqno(11.3)$$
\vskip.1in
\vskip.05in
\indent Let $C_0^\infty
({\bf R};{\bf R})$ denote the space of infinitely differentiable
functions $f:{\bf R}\rightarrow {\bf R}$ such that $\vert
x\vert^jf^{(k)}(x)\rightarrow 0$ as $x\rightarrow\pm\infty$, and
suppose that $u\in C_0^\infty ({\bf R}; {\bf R})$. Then with $\lambda
=k^2$, let $s(k)$ be the scattering matrix, which depends analytically
upon $k$, and let $s_{21}(k)$ be the bottom left entry, which satisfies
$s_{21}\in C_0^\infty ({\bf R}; {\bf R})$ and
$\overline{s_{21}(k)}=s_{21}(-k)$. There exists a solution
$-f''_-(x;k)+u(x)f_-(x;k)=k^2f_-(x;k)$ such that 
$$f_-(x;k)\asymp e^{-ikx}+s_{21}(k)e^{ikx}\qquad (x\rightarrow\infty )$$
\noindent so it is natural to consider the Fourier integral
$\int_{-\infty}^\infty (f_-(x;k)-e^{-ikx})e^{iky}dk$, and its asymptotic
approximation
$$\phi (x+y)={{1}\over{2\pi}}\int_{-\infty}^\infty e^{ik(x+y)}s_{21}(k)\,
dk\eqno(11.4)$$
\noindent gives a real function. Dyson inverted the scattering map 
$q\mapsto s_{21}$ by the formula (1.12).\par

\vskip.05in
\noindent {\bf 11.1 Admissible linear systems in scattering} As in [9, Theorem 4.2] and
[21, p. 486] we can introduce a linear system and Hankel determinant to realise scattering
functions. The following formulas are similar, but slightly
different from those in [21]. Let $H=L^2({\bf R}; {\bf C})$ and let $b_1,b_2:{\bf
R}\rightarrow {\bf C}$ be smooth
functions of compact support such that $b_1(-k)=\overline{b_1(k)}$,
$b_2(-k)=\overline{b_2(k)}$ and $\vert b_1(k)\vert=\vert b_2(k)\vert$
for all $k\in {\bf R}$, and let
$$\eqalignno{B:{\bf C}\rightarrow H:&\,\alpha\mapsto b_1(k)\alpha ;\cr
             e^{-xA}:H\rightarrow H:&\,f(k)\mapsto e^{ixk}f(k);\cr    
               C:H\rightarrow {\bf C}: &\,f(k)\mapsto
{{1}\over{2\pi}}\int_{-\infty}^\infty f(k)b_2(k)\, dk.&(11.5)\cr}$$
\noindent The potential $u$ is in $C^\infty_0({\bf R}; {\bf R})$, and
we assume that there are no bound states, so we are in the scattering
case of Schr\"odinger's equation. Then $(-A,B,C)$ has scattering 
function $\phi
(x)=\int_{-\infty}^\infty e^{ixk}b(k)dk/2\pi$, while
$\Sigma_{i\kappa}=(-A,(i\kappa I -A)
(i\kappa I+A)^{-1}B,C )$ has scattering function $\phi_{i\kappa} (x)=
\int_{-\infty}^\infty e^{ixk}b(k)(\kappa +k)(\kappa -k)^{-1}dk/2\pi$,
which is unambiguously defined for real $\kappa$ since the
Hilbert transform is bounded on $H$;
the corresponding potential is
$u_{i\kappa}(x)=-2{{d^2}\over{dx^2}}\log \tau_{i\kappa}(x)$.\par 

\indent The pole divisor ${\bf z}_\lambda (t)$ is determined by $\{
z_n(t): \psi_{BA}(z_n(t), t; \lambda )=0\}$ and is associated with the
potential $u_\lambda (x;t)=-2{{d^2}\over{dx^2}}\log\tau_\lambda (x,t)$. In
this section, we introduce dynamical systems on ${\cal T}_{\bf
R}^\infty$ such that $u_\lambda (x,t)$ undergoes the nonlinear evolution
associated with the KdV hierarchy. To obtain $KdV(2n+1)$, we vary
$t_{2n+1}$ while fixing $t_{2j+1}$ for $j\neq n$.\par
\vskip.05in

\noindent {\bf Definition} (Torus) Let $(\tau_\infty)=\{ p_n: n=1, 2,\dots \}$ be the zeros of $\tau_\infty $ and $O_n$ be the real oval in
$\cup_{\lambda\in (-\infty , 
\infty]}(\tau_\lambda )$  that is based upon $p_n$. Then let
${\cal T}^\infty_{\bf R}=\prod_{n=1}^\infty O_n$ and consider 
${\bf z}_\lambda=\{ z_n:n=1,2, \dots
\}=(\tau_\lambda )$ with $z_n\in O_n$. Then ${\bf z}_\lambda\in {\cal
T}^\infty_{\bf R}$ is the
pole divisor of $\psi_{BA}(x,\lambda )$ in the infinite real torus 
${\cal T}_{\bf R}^\infty$.\par   
\vskip.05in
\noindent {\bf Proposition 11.1} {\sl (i) The Baker--Akhiezer function $\lambda \mapsto \psi_{BA}(x,\lambda )$ is holomorphic on} ${\bf C}\setminus {\hbox{Spec}}(A).$ 
{\sl If $E$ has finite rank, then
$\lambda\mapsto \psi_{BA}(x,\lambda )$ is meromorphic on
${\bf C}$ with the only possible poles being
 on the spectrum of $A$.}\par
\indent {\sl (ii) $\psi_{BA}(x;\lambda )$ belongs
to a Liouvillian extension of the field of fractions of ${\bf
A}_0$ and satisfies, in the notation of Proposition 2.2,} 
$$\psi_{BA}(x,\lambda )=e^{\lambda x}\det \Bigl(I-\int_x^\infty
T(x,y)e^{\lambda (y-x)}dy\Bigr)\qquad (\Re \lambda <0).\eqno(11.20)$$
\indent {\sl (iii) Suppose that $\Sigma$ is a block diagonal direct sum
$\oplus_{j=1}^\infty\Sigma_j$, where $\Sigma_j$ is a periodic linear
system with $T_j$ as in Theorem 8.1. Then}
$$\psi_{BA}(x,\lambda )=e^{\lambda x}\prod_{j=1}^\infty\det \Bigl(I-\int_x^\infty
T_j(x,y)e^{\lambda (y-x)}dy\Bigr)\qquad (\Re \lambda <0).\eqno(11.21)$$ 
\indent {\sl (iv) Suppose that $B$ and $C$ have rank one, and $A$ is bounded. Then the Baker--Akhiezer function satisfies $\psi_{BA}(x,ik)-e^{ikx}\rightarrow
0$ as $x\rightarrow\infty$ and}
$$-\psi_{BA}''(x,ik)+u(x)\psi_{BA}(x,ik)=k^2\psi_{BA}(x,ik)\qquad (x\in
{\bf R}).\eqno(11.9)$$
\noindent {\sl There exists $r_0>0$ such that} 
$$\int_{\vert \lambda\vert =r}\psi_{BA}(x, \lambda )\psi_{BA}(y, -\lambda )d\lambda =0\qquad (x,y\in {\bf R}, r>r_0).$$
\indent {\sl (v) If $\tau_\lambda$ has only simple 
zeros, then each zero of $\psi_{BA}(z,\lambda)$ in
$(\tau_\lambda )$ processes in a real oval based at a pole of $\psi_{BA}(z,
\lambda )$ in $(\tau_\infty )$ as $\lambda$ describes $(-\infty
,\infty )$. The pole divisor defines a map $\Sigma_\lambda\mapsto
{\bf z}_\lambda$ from the periodic linear system to the real torus ${\cal
T}^\infty_{\bf R}$.}\par
\vskip.05in
\noindent {\bf Proof.}  (i) The first part follows from Theorem 8.2. Suppose then that $E$ has finite rank, and note that 
$(\lambda I+A)(\lambda I-A)^{-1}E$ is a rational function with values in
the space of operators on a finite-dimensional Hilbert space.
Hence the determinant $\tau_\lambda$ is meromorphic as a function of
$\lambda$ on ${\bf P}^1$, and the only possible poles are at the $\lambda\in {\hbox{Spec}}(A)$. Hence $\lambda \mapsto \psi_{BA}(x,\lambda )$ is meromorphic on ${\bf C}$ and has poles independent of $x$, as in Krichever's definition of the Baker--Akhiezer function.\par

 \indent (ii) We have $u(x,\lambda )=-2(\log
\tau_\lambda )''$ in ${\cal A}_0$ by Theorem 8.2, hence $\psi''_{BA}(x,
\lambda )$ belongs to ${\cal A}_0$; we integrate this to obtain
$\psi_{BA}$ in some Liouville extension. By some simple manipulations, we have
$$\det \Bigl( I+R_x(\lambda I+A)(\lambda I-A)^{-1}\Bigr)=\det
(I+R_x)\det \Bigl( I+(\lambda
I-A)^{-1}(AR_x+R_xA)(I+R_x)^{-1}\Bigr)\eqno(11.23)$$
\noindent where $AR_x+R_xA=e^{-xA}BCe^{-xA}$, and hence

$${{\det \Bigl( I+R_x(\lambda I+A)(\lambda I-A)^{-1}\Bigr)}\over{\det
(I+R_x)}}=\det \Bigl( I+Ce^{-xA}(I+R_x)^{-1}(\lambda
I-A)^{-1}e^{-xA}B\Bigr),\eqno(11.24)$$
\noindent and $\int_x^\infty e^{\lambda (y-x)}e^{-yA}=-(\lambda
I-A)^{-1}e^{-xA}$, which leads to the stated identity. 
Moreover,
the right-hand side is analytic in $\lambda$  when $\vert
\lambda\vert >\Vert A\Vert$, and $\psi_{BA}(x, \lambda )=e^{\lambda
x}(1+O(\lambda^{-1}))$ as $\vert\lambda\vert\rightarrow\infty$.\par
By the proof of Theorem
8.2, ${{d^2}\over{dx^2}}\log \psi_{BA} (x; \lambda )$ belongs to
${\cal A}_0$.\par 
\indent (iii) 
With $T_{i\kappa}$ and
the corresponding potential 
$u_{i\kappa}(x)=-2{{d^2}\over{dx^2}}\log \tau_{i\kappa}(x)$ defined for the linear system
$\Sigma_{i\kappa}$, we introduce 
$$f_{i\kappa}(x,k)=e^{ikx}+\int_x^\infty T_{i\kappa } (x,y)e^{iky}\,
dy.\eqno(11.11)$$
\noindent By repeated integration by parts, one verifies  
$-f_{i\kappa}''(x, \pm
k)+u_{i\kappa }(x)f_{i\kappa}(x,\pm k)=k^2f_{i\kappa}(x,\pm k)$ and
$f_{i\kappa}(x,\pm k)-e^{\pm ikx}\rightarrow 0$ as $x\rightarrow\pm
\infty $. In particular, with $i\kappa =\infty$ we can express
$$\eqalignno{f_\infty (x,k)&=e^{ikx}-Ce^{-xA}(I+R_x)^{-1}\int_x^\infty
e^{-yA}Be^{iky}dy\cr
&=e^{ikx}\bigl((1+Ce^{-xA}(I+R_x)^{-1}(ikI-A)^{-1}e^{-xA}B)\cr
&=e^{ikx}\det
\bigl(I+(ikI-A)^{-1}e^{-xA}BCe^{-xA}(I+R_x)^{-1}\bigr)&(11.12)\cr}$$
\noindent where we have used a simple identity for rank-one
operators, hence
$$f_\infty (x,k)=e^{ikx}{{\det
(I+R_x+(ikI-A)^{-1}e^{-xA}BCe^{-xA})}\over{\det (I+R_x)}},\eqno(11.13)$$
\noindent and we can finish by using Lyapunov's equation
$$f_\infty (x,k)=e^{ikx}{{\det (I+R_x-(ikI-A)^{-1}R_x')}\over{\det
(I+R_x)}},\eqno(11.14)$$
\noindent where the determinant on the numerator is
$$\det \bigl(I+R_x+(ikI-A)^{- 1}(AR_x+R_xA)\bigr)=
\det\bigl(I+R_x(ikI+A)(ikI-A)^{-1}\bigr).\eqno(11.15)$$ 

This follows from (i), since
$$\psi_{BA}(x, \lambda )=e^{-\lambda x}-\int_x^\infty T(x,y)e^{\lambda
y}\, dy$$
\noindent and one can consider the differential equation for $T(x,y)$ to
verify that $\psi_{BA}(x., \lambda )$ satisfies ().\par
\indent (iv) We reduce to the case of the admissible linear
system $(-A-\varepsilon I, B,C)$, which has input and output space
${\bf C}$, as in Proposition 2.5. For $\varepsilon >0$, let
$R^{(\varepsilon )}_x =e^{-2\varepsilon x}e^{-xA}Ee^{-xA}$, so that
$R^{(\varepsilon )}_x\rightarrow 0$ exponentially fast as
$x\rightarrow\infty$, and $R^{(\varepsilon )}_x$ satisfies the Lyapunov
equations
$$-{{d}\over{dx}}R^{(\varepsilon )}_x=(A+\varepsilon I)
R^{(\varepsilon )}_x+R^{(\varepsilon )}_x(A+\varepsilon
I),\eqno(11.25)$$
\noindent with 
$$-{{d}\over{dx}}R^{(\varepsilon )}_x\vert_{x=0}=BC+2\varepsilon
E.\eqno(11.26)$$
\noindent Since $BC$ and $E$ are trace class, we can introduce
$\tau_\infty^{(\varepsilon )}(x)=\det (I+R^{(\varepsilon )}_x)$ and 
$$\tau_\lambda^{(\varepsilon )}(x)=\det (I+R^{(\varepsilon )}_x(\lambda
I+\varepsilon I+A)(\lambda I-\varepsilon I-A)^{-1}),\eqno(11.27)$$
\noindent whenever $\lambda -\varepsilon$ is in the resolvent set of
$A$; likewise we can introduce $u^{(\varepsilon
)}(x)=-2{{d^2}\over{dx^2}}\log \tau_\infty^{(\varepsilon )}(x)$. Now
the Baker--Akhiezer function 
$$f^{(\varepsilon )}(x, k)=e^{ikx}{{\tau_{ik}^\varepsilon
(x)}\over{\tau_\infty^{(\varepsilon )}(x)}}\eqno(11.28)$$
\noindent satisfies
$$-{{d^2}\over{dx^2}}f^{(\varepsilon )}(x)+u^{(\varepsilon
)}(x)f^{(\varepsilon )}(x,k)=k^2f^{(\varepsilon )}(x,k);\eqno(11.29)$$
\noindent letting $\varepsilon\rightarrow 0$, we obtain
$$-{{d^2}\over{dx^2}}f(x)+u(x)f(x,k)=k^2f(x,k);\eqno(11.30)$$  
\noindent as required.\par
\indent We consider the integral
$$K(x,y)=\int_{\vert\lambda\vert = r}e^{\lambda (x-y)} {{\det (I+(\lambda +A)(\lambda
-A)^{-1}R_x) 
\det (I-(\lambda -A)(\lambda +A)^{-1}R_y)}\over {\det (I+R_x)\det (I+R_y)}}d\lambda,$$
\noindent and prove that
$$\Bigl(-{{\partial^2}\over{\partial x^2}}+u(x)\Bigr) K(x,y)
=\Bigl(-{{\partial^2}\over{\partial y^2}}+u(y)\Bigr) K(x,y)$$
\noindent with the initial conditions
$$\eqalignno{ K(x,x)&=0&(11.65)\cr
{{\partial }\over{\partial x}}K(x,x)&=0.&(11.67)\cr}$$
\noindent The partial differential equation follows immediately from
$-\psi_{BA}(x, \lambda )+u(x)\psi_{BA}(x, \lambda )=\lambda^2\psi_{BA}(x, \lambda )$ and
its companion identity with $-\lambda$ instead of $\lambda$. The initial condition () follows
from Cauchy's Theorem, as applied to the integrand in a neighbourhood of $\lambda=\infty$. The
initial condition on the derivative involves the integral
$$\eqalignno{ {{\partial }\over{\partial x}}&K(x,x)\cr
&=\int_{\vert\lambda \vert =r}\Bigl( \lambda +{\hbox{trace}}\bigl(
(I+(\lambda I+A)(\lambda I-A)^{-1}R_x)^{-1}(\lambda I+A)(\lambda I-A)^{-1}R_x'
-(I+R_x)^{-1}R_x'\bigr)\Bigr)\cr
&\qquad\det \bigl( I+2A(\lambda I -A)^{-1}R_x(I+R_x)^{-1}\bigr) \det\bigl( I-2A(\lambda I
+A)^{-1}R_x(I+R_x)^{-1}\Bigr)d\lambda . \cr}$$
\noindent We expand the various terms as series in $\lambda^{-1}$ and evaluate the integral by the calculus of residues. After applying Lyapunov's equation to
simplify the terms in $R_x'$, we are left with a multiple of the term
$${\hbox{trace}}\bigl( R_xF_xA R_xF_xA 
\bigr) +2{\hbox{trace}}\Bigl(\bigl(AR_xF_x\bigr)^{\wedge 2}\Bigr)
-\Bigl({\hbox{trace}}\bigl(AR_xF_x\bigr)\Bigr)^2,$$
\noindent where $F_x=(I+R_x)^{-1}$; by considering the eigenvalues of $AR_x(I+R_x)^{-1}$, we see that this vanishes. Hence
$K(x,y)$ satisfies a hyperbolic differential equation with null initial conditions, and hence
vanishes.\par 
    

\indent (v) Clearly the poles of $\psi_{BA}(z,\lambda )$
occur at the zeros of $\tau_\infty (z)$, and hence form the set
$(\tau_\infty )$, for all $\lambda$. The zeros of
$\psi_{BA}(z, \lambda )$ form the set $(\tau_\lambda )$, which does vary
with $\lambda$. The subset $\{ 
(\lambda I-A)(\lambda I+A)^{-1}: \lambda \in {\bf R}\}$ of
${\cal L}(H)$ is
compact in the norm topology since the spectrum of $A$ is separated
from ${\bf R}$; hence $\tau_\lambda $ gives a compact family of
holomorphic functions for the topology of uniform convergence on
compact sets, with $\tau_{-\infty }(z)=\tau_\infty (z)$. For each
bounded open subset $\Omega$ of ${\bf C}$, the set $\{ z\in \Omega :
\tau_\lambda (z)=0\}$ has a uniformly bounded number of terms for $-\infty \leq \lambda\leq
\infty$, by Jensen's formula and the Proposition 8.4. Each zero 
depends continuously upon $\lambda$ by the inverse function theorem, and 
describes an oval for 
$-\infty \leq \lambda\leq\infty$.\par    

\rightline{$\square$}\par

\noindent Given a tau function from a periodic linear system
$(-A,B,C;E)$  where the
input and output space are both ${\bf C}$, we
consider the conditions under which $\tau$ arises from the theta
functions on a compact algebraic curve. First we consider families of linear systems as in
Theorem 11.2, with common $A$, which are parametrized by $\lambda\in
{\bf P}^1\setminus {\hbox{Spec}}(A)$ and time parameters $(t_3, t_5,
\dots )$, giving tau functions $\tau_\lambda (x,t)$. Initially $x$ and
$t_{2j+1}$ are real, and $\tau_\lambda (x,t)$ is $\pi$ periodic in each
variable, hence $\tau_\lambda (x,t)$ gives a periodic function on the infinite real torus
${\bf R}^\infty/\pi {\bf Z}^\infty$; then we extend to complex $x$ and
$t_{2j+1}$, so that $\tau_\lambda (x,t)$ is entire. By forming quotients of
such functions, we aim to realise typical tau functions.\par
\indent To introduce the required linear systems, we let
$${\bf T}=\{(x,t_3, t_5, \dots )\in {\bf R}^\infty : 
\lim\sup_{j\rightarrow\infty}\vert
t_{2j+1}\vert^{1/j}=0\}\eqno(11.28)$$ 
\noindent which gives an
Abelian group under addition, and for
$(x,t) \in {\bf T}$, let\par
\noindent $U(t)=\exp (-\sum_{j=1}^\infty t_{2j+1}A^{2j+1})$, which gives a
multi parameter group of operators such that $U(s+t)=U(s)U(t)$.
 Then we replace
$\Sigma_\infty (0)=(-A,B,C;E)$ of Theorem 11.2 by 
$$\Sigma_\lambda (t)=\Bigl(-A,(\lambda I+A 
)(\lambda I-A )^{-1}U(t)B, CU(t), (\lambda I+A
)(\lambda I-A)^{-1}U(t)EU(t)\Bigr)\eqno(11.29)$$
\noindent for $\lambda\in {\bf P}^1\setminus {\hbox{Spec}}(A)$. Each $\Sigma_\lambda (t)$ gives a space ${\cal A}_0(t, \lambda
)$ of potentials as in Theorem 11.2(iii), while $\lambda$ is a spectral
parameter as in Proposition 11.5. Let $({\bf A_0}, d/dx)$ be the differential
ring generated by $\Sigma$ as in Theorem 11.3(ii), and let $({\bf
A}_\infty , \partial /\partial x, \partial /\partial t_{2j+1} )$ be the
differential ring generated by all the $\Sigma_\lambda (t)$; then ${\bf
A}_0\subseteq {\bf A}_\infty$, and the inclusion splits by mapping
$t_{2j+1}\mapsto 0$ for all $j=1, 2, \dots .$ \par
\vskip.05in
\noindent {\bf Definition} (Baker--Akhiezer function) Let 
$\xi (x,t,\lambda )= {{x}{\lambda}}+\sum_{j=1}^\infty
{{t_{2j+1}}{\lambda^{2j+1}}}$ and define the
quotient 
$$\psi_{BA}(x,t;\lambda )=e^{\xi (x,t,\lambda )}{{\tau_\infty \Bigl(x-{{1}\over{\lambda}} , t_3-{{1}\over{3\lambda^3}},
t_5-{{1}\over{5\lambda^5}}, \dots \Bigr)}\over 
{\tau_\infty (x,t_3, t_5, \dots )}}\eqno(11.30)$$
\noindent to be the Baker--Akhiezer function of the periodic linear
system $(-A,B,C;E)$ under $U(t)$. Let $[\alpha ]=(\alpha, \alpha^3/3, \alpha^5/5, \dots ).$\par
\indent (ii) Say that the family of linear systems $\Sigma_\zeta (t)$ satisfies $KdV$ if for all $(x,t), (y,t)\in {\bf T}^\infty$ there exists $r_0>0$ such that
$$\int_{\vert \zeta\vert =r}\psi_{BA}(x,t; \zeta )\psi_{BA}(y,s;-\zeta )\, d\zeta =0\qquad (r>r_0).$$ 
\vskip.05in
\indent This definition is consistent with earlier in this section, 
and with Shiota's definition [54]. However, we
cannot expect a precise analogue of Proposition 11.2, which
expresses eigenfunctions in terms of $\psi_{BA}$.\par
\vskip.05in
\noindent {\bf Lemma 11.3} {\sl (i) The scattering function
$\Phi_\lambda (x,y)=CU(t)e^{-xA}(\lambda I+A)(\lambda I-A)^{-1}U(t)B$ for
$\Sigma_\lambda (t)$ satisfies
$${{\partial^{2j+1}}\over{\partial
x^{2j+1}}}\Phi_\lambda (x,t)+{{\partial}\over{\partial t_{2j+1}}}\Phi_\lambda 
(x,t)=0.\eqno(11.31)$$
\indent (ii) $\tau_\lambda (x,t)$ is holomorphic for} 
$(x,t, \lambda )\in {\cal C}\times {\cal C}^\infty \times ({\bf P}^1\setminus
{\hbox{Spec}}(A)),$ {\sl where  ${\cal C}={\bf C}/\pi{\bf Z}$ is the 
complex cylinder.\par
\indent (iii) $\lambda\mapsto \psi_{BA}(x,t, \lambda )$ is
meromorphic on} ${\bf C}\setminus {\hbox{Spec}}(A),$
{\sl while
 $(x,t)\mapsto\psi_{BA}(x,t, \lambda )$ is meromorphic and quasiperiodic with respect
to the lattice $\pi{\bf Z}^\infty$ in ${\bf C}^\infty$.}\par 
\vskip.05in
\noindent {\bf Proof.} (i) Since $U(t)$ is actually
analytic in each $t_j$ this is a straightforward computation.\par
\indent (ii) First we observe that
$$\tau_\lambda (x,t)=\det 
\Bigl(I+(\lambda I+A)(\lambda I-A)^{-1}U(2t)e^{-2xA}E\Bigr).
\eqno(11.32)$$ 
\noindent as in Lemma 8.1. Hence $\lambda
\mapsto \tau_\lambda (x,t)$ is holomorphic on ${\bf P}^1\setminus
{\hbox{Spec}}(A)$, and $(x,t)\mapsto \tau_\lambda (x,t)$ is entire in
each variable since $A$ is bounded. The spectrum of $A^{2j+1}$ is contained in\par
\noindent $\{ -iN^{2j+1}, -i
(N-1)^{2j+1}, \dots , iN^{2j+1}\}$, so $e^{2\pi A^{2j+1}}=I$, and by
Theorem 11.2 $\tau_\lambda (x+\pi ,t)=\tau_\lambda (x,t)$; likewise
$\tau_\lambda (x,t)$ is unchanged by adding $\pi$ to $t_j$; so
$\tau_\lambda (x,t)$ is periodic with respect to 
$\pi {\bf Z}^\infty$
in ${\bf C}^\infty$.\par
\indent (iii) Observe that by the choice of $(x,t)$, the function
$\lambda\mapsto \xi (x,t, \lambda )$ is entire. Then by calculations as 
in Proposition 11.1, we have the expression 
$$\eqalignno{\psi_{BA}(x,t, &\lambda )\cr
&=e^{\xi (x,t,\lambda )}\Bigl( 1+{\hbox{trace}}\, CU(t)e^{-xA}\bigl(
I+e^{-xA}U(t)EU(t)e^{-xA}\bigr)^{-1}(\lambda
I-A)^{-1}e^{-xA}U(t)B\Bigr)\cr}$$ 
\noindent which is evidently holomorphic on ${\bf C}\setminus
{\hbox{Spec}}(A)$. If $(\lambda I-A)^{-1}$ has a pole at $\lambda_j\in 
{\hbox{Spec}}(A)$ with residue $A_j$, then $\psi_{BA}(x,t, \lambda )$
possibly has a possible pole at $\lambda_j$ with residue
$$ e^{\xi (x,t,\lambda_j )}{\hbox{trace}}\, CU(t)e^{-xA}\bigl(
I+e^{-xA}U(t)EU(t)e^{-xA}\bigr)^{-1}A_je^{-xA}U(t)B.$$
 With $(e_j)_{j=0}^\infty$ the standard unit vector basis in
${\bf T}^\infty,$ we deduce from (ii) that $\psi_{BA}(x, t+\pi e_j,
\lambda )=e^{2\pi \lambda^{2j+1}}\psi_{BA}(x,t, \lambda )$, and 
$(x,t)\mapsto \psi_{BA}(x,t,\lambda )$ is meromorphic.\par 
\rightline{$\square$}\par
\indent For any uniform periodic linear system as in Theorem 8.2, we let 
$$\Phi (x,t;\lambda ;x_0)=\Bigl[{{\partial^{j-1}\psi_{k}(x,t;\lambda 
;x_0)}\over{\partial x^{j-1}}}\Bigr]_{j,k=1}^N;\qquad \Phi
(x_0,t;\lambda ;x_0)=I_N\eqno(11.33)$$ 
\noindent be the Wronskian matrix associated with some collection of functions, and suppose that $\Phi$ is invertible; then with $U_j(x,t;\lambda;x_0 )={{\partial \Phi}\over{\partial t_{2j-1}}}\Phi^{-1}$, the operators $L_j={{\partial}\over {\partial t_{2j-1}}}-U_j$  give a commuting family of matrix differential operators. \par
\vskip.05in
\noindent {\bf Proposition 11.4} {\sl Let $\Sigma_\zeta (t)$ be a uniformly
periodic linear system with tau function $\tau_\infty (x,t)$, and for $\alpha, \beta, \gamma\in {\bf C}
\setminus \{0\}$ such that $\alpha^{-1}+\beta^{-1}+\gamma^{-1}=0$, let 
$(y,s)-(x,t)=2[\alpha ]+2[\beta ]+2[\gamma ]$. Then for all $r>\max\{ 1/\vert\alpha\vert,
1/\vert\beta\vert, 1/\vert\gamma\vert, \Vert A\Vert\}$, }
$${{1}\over{2\pi i}}\int_{\vert \zeta\vert =r}\psi (x,t,\zeta )\psi
(y,s,-\zeta )\, d\zeta =2{{\partial}\over{\partial x}}\log{{\tau_\infty
(x,t)}\over{\tau_\infty ((x,t)+2[\alpha ]+2[\beta] +2[\gamma ])}}.$$
\vskip.05in
\noindent {\bf Proof.} The right-hand side has the same form as the
expression in Corollary 5.3(ii). The left-hand side is 
$${{1}\over{2\pi i}}\int_{\vert \zeta\vert =r}{{(1+\alpha
\zeta )(1+\beta \zeta )(1+\gamma \zeta )}\over{(1-\alpha\zeta )(1-\beta \zeta )(1-\gamma
\zeta )}}{{\tau_\zeta
(x,t)\tau_{-\zeta}(y,s)}\over{\tau_\infty (x,t)\tau_\infty
(y,s)}}d\zeta,$$
\noindent which by previous calculations gives
$${{1}\over{2\pi i}}\int_{\vert \zeta\vert
=r}{{(1+\alpha
\zeta )(1+\beta \zeta )(1+\gamma \zeta )}\over{(1-\alpha\zeta )(1-\beta \zeta )(1-\gamma
\zeta )}}\det \bigl( I+2(\zeta
I-A)^{-1}R'_xF_x\bigr) \det \bigl( I-2(\zeta
I+A)^{-1}R_y'F_y\bigr) d\zeta,$$
\noindent where $R_x'$ and $F_x$ are evaluated at $(x,t)$, while $R_y'$ and $F_y$ at $(y,s)$
for $\Sigma_\zeta (t)$. The values of the integral is given  by the calculus of residues to be
$$-{{2}\over{\alpha}}-{{2}\over{\beta}}-{{2}\over{\gamma }}-2{\hbox{trace}}\,\bigl( R_x'(I+R_x)^{-1}\bigr)
+2{\hbox{trace}}\,\bigl( R'_y (I+R_y)^{-1}\bigr),$$
\noindent which implies the result.\par

\noindent {\bf Corollary 11.4} {\sl Suppose further that there exists a compact Riemann surface ${\bf P}$ with marked points $\{ {\bf p}_1,\dots , {\bf p}_m\}$ and field ${\bf K}$ of meromorphic functions on ${\bf P}\setminus \{ {\bf p}_1, \dots ,{\bf p}_m\}$ and that $W(t,{\bf p})$ is a $M_N$-valued function on ${\bf T}^\infty \times {\bf P}$ such that the entries of ${\bf p}\mapsto W(t,{\bf p})$ belong to ${\bf K}$, and that  $[L_j, W]=0$. Then  the family $L_j$ is algebraic and has a spectral curve ${\cal E}$, with marked points $\{ {\bf q}_1,\dots ,{\bf q}_{nm}\}$,  that gives a finite cover of ${\bf P}$.} \par
\vskip.05in
\noindent {\bf Proof.} This is immediate from Lemma 5.3.\par
\rightline{$\square$}\par

\vskip.05in

\noindent We now consider how the results of section 7 relate to the notions of
Liouville integrability and finite gap integration. The results of
this section are essentially corollaries of some subtle results
proved elsewhere, and the most interesting relate to elliptic
potentials.\par
\vskip.05in 

\vskip.05in

\indent The solutions of (11.1) turn out to be complicated polynomials in
 $u$ and its derivatives, as one can prove by induction. Nevertheless,
we can express a solution $g_m$ simply in terms of
$\lfloor{\bf A}\rfloor$. The following proposition is a compilation of known
results, and included for completeness.\par
\vskip.05in

\noindent {\bf Proposition 11.5} {\sl Let $(-A,B,C;E)$ be a periodic linear system as in Theorem 8.2.}\par
\indent {\sl (i) The complex vector space spanned by the $\lfloor A^{2m-1}\rfloor$ for $m=1, 2, \dots $ is finite dimensional.}\par 
\indent {\sl (ii) Suppose further that ${\bf S}$ is commutative, as in Corollary 8.4, and that the operators in $\lfloor {\bf A}\rfloor$ all have a common eigenvector. Then there is a homomorphism of differential rings $\rho:\lfloor {\bf A}\rfloor\rightarrow C(\Omega ; {\cal L}(H))$ such that $g_m=\rho(\lfloor A^{2m+1}\rfloor)$ gives solutions to the the stationary KdV equations.}\par
\indent {\sl (iii) Suppose that $u$ has finite gap, so that $g_m=0$ for some $m$. 
Then there exists a hyperelliptic spectral curve ${\cal E}$, 
with a branch point $p_0$, a meromorphic function $\lambda$ on ${\cal E}$, and a
pair of distinct points $p_j, q_j\in {\cal E}$ for each point} $ij\in
{\hbox{Spec}}(A)$, {\sl all independent of $x$, such that
 $\lambda\mapsto\psi_{BA}(x, \lambda )$ is holomorphic on}
 ${\cal E}\setminus \{ p_j, q_j: j=0; ij\in {\hbox{Spec}}(A)\}$. \par

\vskip.05in
\vskip.05in
\noindent {\bf Proof.}  (i) Let $m$ be the minimal polynomial of degree $N$ for the algebraic operator $A$. Then for
each entire function $f$, either $f(A)=0$ or there exists a polynomial $r$ of degree less than or equal to $N$
such that $f(A)=r(A).$ Hence 
the span of the $A^{2m-1}$ for $m=1,2, \dots $ is finite-dimensional, and hence its image under $\lfloor \, .\, \rfloor$ is also finite-dimensional.\par
\indent (ii) By Theorem 8.2(ii) and Corollary 8.4, $\lfloor {\bf A}\rfloor$ is a commutative differential ring of operators on $H$ over $\Omega$. Let $\psi$ be a common eigenvector
and let $\rho (\lfloor P\rfloor )=\lambda (P)$ where $\lfloor P\rfloor \psi=\lambda (P)\psi$, so that $\rho$ defines homomorphism of commutative rings. By Proposition 3.4. the images $g_m=\rho (\lfloor A^{2m+1}\rfloor)$ satisfies the recurrence relation for the KdV hierarchy.\par
\indent (iii) From the 
recurrence relation (3.15), we deduce that $g_n=0$ for all $n\geq m$, so $u$ is
 finite gap and ${\bf C}[\lambda , u,u',
u'', \dots ]={\bf C}[\lambda , u,u',
u'', \dots , u^{(m+1)}]$ is a differential ring. Any solution of the stationary KdV 
equations is meromorphic on 
${\bf C}$ [52, 6.10]. 
Let $\lambda_0<\lambda_1<\dots <\lambda_{2g}$ be the simple
zeros of $4-\Delta (\lambda )^2=0$, and introduce the spectral curve
$${\cal E}=\Bigl\{ (z,w):w^2=\prod_{j=0}^{2g}(z-\lambda_{j})\Bigr\}\cup\{
(\infty ,\infty )\},\eqno(11.35)$$
\noindent Now there exists a solution $\rho (x,\lambda )$ to Drach's
equation (8.2) 
$$\mu^2=-{{1}\over{2}}\rho (x, \lambda )\rho''(x, \lambda
)+{{1}\over{4}}\rho'(x, \lambda )^2 +(u(x)+\lambda )\rho (x, \lambda
)^2\eqno(11.36)$$
\noindent such that $\mu (\lambda )$ is independent of $x$ and
$\lambda\mapsto \rho (x, \lambda )$ is a polynomial, which we factor as
$\rho (x, \lambda )=\prod_{j=1}^g (\lambda -\gamma_j(x))$. Brezhnev [13] gives the
solution 
$$\psi_{\pm} (x)=\exp\Bigl( \sum_{j=1}^g \int^{\gamma_j(x)}
{{(w\pm\mu )dz}\over{(z-\lambda )w}}\Bigr),\eqno(11.37)$$
\noindent where the integral is taken along ${\cal E}$. Here $u$ and its derivatives are rational functions on
${\cal E}$; see [34, 50]. For such a potential $u$, the functions
$\psi_{\pm}$ of (11.37) give locally meromorphic solutions to 
Schr\"odinger's equation.\par
\indent Suppose that
${\cal E}$ has genus $g\geq 2$, and  choose $p_0$ to be one of
the $2g+2$ branch points of the holomorphic two sheeted cover ${\cal
E}\rightarrow {\bf P}$, and then observe that there exists a
meromorphic function $\lambda$ on ${\cal E}$ such has precisely one
pole, namely a double pole at $p_0$, and hence has degree two 
(When $g=1$, we can use $\lambda (p)=\wp (p-p_0)$).\par
\indent The exponential $e^{x\lambda}$ gives an
essential singularity in the variable $\lambda$ for ${\bf p}$ close to
${\bf p}_0$. Also, $\lambda\mapsto (\lambda I+A)(\lambda I-A)^{-1}E$ is a rational 
function, with trace class values, and the only possible poles are on
the spectrum of $A$; hence ${\bf p}\mapsto \tau_\lambda$ gives a holomorphic function,
except at finitely many points of ${\cal E}$, which we list as
${\bf p}_j,{\bf q}_j$ for $ij$ in the spectrum of $A$. (These can give essential singularities when $E$ has infinite rank.)\par
\rightline{$\square$}\par

bounded interval $[a,b]$, we let  
$\omega (f,g)=\int_a^b f(x)g'(x)dx.$ In particular, if $g(x)$ is nowhere zero on 
$[a,b]$, then $\omega (g^{-1}, g)$ gives the variation of $\log g$ over $[a,b]$.\par

\indent We can view $\Omega^1C^\infty$ as a space of integrable
operators, via the following construction. Let $M:L^2[a,b]\rightarrow L^2[a,b]:$ $f(x)\mapsto xf(x)$ be the operator
of multiplication by $x$, and let ${\cal F}=C^\infty\otimes C^\infty$ be the 
space of finite-rank integral operators on $L^2[a,b]$ that have smooth kernels. Then 
${\cal R}=\{ K\in B(L^2[a,b]):KM-MK\in {\cal F}\}$ is the space of integral operators which have 
kernels of the form\indent Let ${\cal R}$ be a unital algebra over ${\bf C}$ with ideal ${\cal
F}$, and suppose that $\delta :{\cal R}\rightarrow {\cal F}$ is a
derivation.  The standard multiplication ${\cal R}\times
{\cal R}\rightarrow {\cal R}:(P,Q)\mapsto PQ$ gives a short exact sequence 
$0\rightarrow \Omega^1{\cal R}\rightarrow {\cal R}\otimes {\cal R}\rightarrow {\cal
R}\rightarrow 0$ of linear maps. We write $PdQ=P\otimes Q-PQ\otimes 1$, which
clearly belongs to $\Omega^1{\cal R}$, and introduce the right module multiplication 
$(PdQ)R=Pd(QR)-PQdR$, so that $\Omega^1{\cal R}$ is the linear span of 
$PdQ$ for
$P,Q\in {\cal R}$ is an ${\cal R}$ bimodule.\par 
\vskip.05in
\noindent {\bf Lemma} {\sl Suppose further that 
$\tau_0:{\cal F}\rightarrow {\bf C}$ is a trace on the bimodule ${\cal F}$ with respect
to the algebra ${\cal R}$. Then the bilinear form 
$\omega :{\cal R}\times
{\cal R}\rightarrow {\bf C}$ defined by $\omega (P,Q)=\tau_0(P\delta Q),$ gives a trace} 
$$\Omega^1{\cal R}\rightarrow {\bf C}:PdQ\mapsto \tau_0(P\delta Q).$$ 
\vskip.05in
\indent The proof is immediate from the universal property of the derivation $d$, as 
in as Lemma 2.6.2 of [Loday].
There is an ${\cal R}$ bimodule map
$\Omega^1{\cal R}\rightarrow {\cal F}: PdQ\mapsto P\delta Q,$ and we compose this with $\tau_0$.

\vskip.05in 
\noindent {\bf Definition} For a fixed $\delta F\in {\cal F}$ introduce a linear functional on ${\cal
R}$ by 
$\lfloor P\rfloor =\tau_0(P\delta F)$.\par
\vskip.05in
\indent We introduce the main examples, beginning with linear systems.\par
 (i) Suppose that $(-A,B,C)$ is a linear system on state space $H$ such that $(e^{-tA})$
defines a holomorphic semigroup. In section 2, we formulate conditions under which $R_z$ is a
holomorphic solution of Lyapunov's equation with values in the trace class operators, and 
introduce a domain $\Omega$ on which $\tau (z)=\det (I+R_z)$ is holomorphic and non-zero. Then
$F_z=(I+R_z)^{-1}$ is a bounded linear operator, with derivative in the trace class.
In section 4, we introduce a differential ring 
${\bf S}$ of holomorphic functions $\Omega\rightarrow {\cal L}(H)$, which contains $A, BC, R_z$ and
$F_z$, so that we can solve (1.10) and (1.11) inside ${\bf S}$. Then we define $\omega (P,Q)={\hbox{trace}}(PQ')$ for $P,Q\in
{\bf S}$, and introduce  $\lfloor P\rfloor={\hbox{trace}}(PF')$. In particular, we have 
$$\lfloor F^{-1}\rfloor =-{{d}\over{dx}}\log\det F=\tau'(z)/\tau (z),$$
\noindent while
$$u(z)=-4\lfloor A\rfloor.$$ 
\noindent If we 
can choose ${\bf S}$ to be a right Noetherian ring, then
we say that $(-A,B,C)$ is finitely generated. 
Thus $\lfloor
\, .\,\rfloor $ linearizes the determinant formula in (1.12).\par
\indent  In section 5 we show that if $A$ is a finite
matrix with eigenvalues $\lambda_j$ such that $\Re\lambda_j>0$, then
 $(-A,B,C)$ is finitely generated. We also recover some
determinant formulas from the theory of solitons.\par 
\indent P\"oppe [41, 49, 50] proved some remarkable product formulas for products and traces of Hankel integral operators,
and applied them to the scattering theory of the Schr\"odinger equation with localized smooth potential. Using the $\lfloor\,
\rfloor$ functional on ${\bf S}$, we obtain a more general version of his calculus, which extends to encompass 
periodic linear systems.\par    
\indent (ii) First we introduce $C^\infty =C^\infty ({\bf R}; {\bf C})$ with the
derivation $\delta :C^\infty\rightarrow C^\infty : f\mapsto f'$; then for any
$$K(x,y)=\sum_{j=1}^n {{f_j(x)g_j(y)}\over{x-y}}$$
\noindent where $f_j, g_j\in C^\infty$. In particular, the subspace 
${\cal R}_0=\{ K\in {\cal R}: KM-MK\in \Omega^1C^\infty\}$ consists of those $K$ such that 
$\sum_{j=1}^n f_j(x)g_j(x)=0$ for all $x\in [a,b]$. Then ${\cal R}_0$ forms an algebra, and 
the elements of ${\cal R}_0$ are said to be integrable operators. Then 
$$\sum_{j=1}^n \omega (f_j,g_j)=\sum_{j=1}^n \int_a^b f_j(x)g'_j(x)\,
dx=-\int_a^b K(x,x)\, dx$$
\noindent is minus the trace of $K$. \par
\noindent {\bf Corollary 7.3} {\sl Let $q$ be as in Proposition 7.2.
Then there exists a hyperelliptic $\wp$ function and a derivation $D_{\bf X}$
on the rational functions on ${\bf X}$ such that $q=Q(\wp, D_{\bf X}\wp , \dots
, D_{\bf X}^{g}\wp )$ for some complex polynomial $Q$.}\par
\vskip.05in

\noindent {\bf Proof.} Let $q$ be a finite gap Hill's potential. Then $q$, when regarded as a
periodic function on the real torus ${\bf T}^g$, can be extended to a periodic
function on some Abelian manifold ${\bf C}^g/\Lambda .$  
\indent Consider positive real points $a_1<b_1<a_2<\dots <b_g<a_{g+1}$ and the
Hamiltonian
$$H={{1}\over{2}}\sum_{j=1}^{g+1} a_jx_j^2+{{1}\over{2}}\Bigl( 
\sum_{j=1}^{g+1}x_j^2\sum_{j=1}^{g+1} y_j^2-\Bigl( \sum_{j=1}^{g+1}
x_jy_j\bigr)^2\Bigr) -\mu\Bigl( \sum_{j=1}^{g+1} x_j^2-1\Bigr)\eqno(7.5)$$
\noindent in the canonical variables $(x_j,y_j)$ 
where $\mu ={{1}\over{2}}\sum_{j=1}^{g+1} a_jx_j^2$. This system was
introduced by Neumann, and represents the motion of a particle that is
constrained to lie on the unit sphere in ${\bf R}^{g+1}$. On the
tangent bundle $\{ (x,y)\in S^{g}\times {\bf R}^{g+1}:\Vert x\Vert =1,
\langle x,y\rangle =0\}$, we introduce the vector field
$$D=\sum_{j=1}^{g+1} y_k{{\partial}\over{\partial x_k}}-
\sum_{j=1}^{g+1}a_kx_k{{\partial }\over{\partial y_k}}+\Bigl( 
\sum_{j=1}^{g+1}a_jx_j^2-\sum_{j=1}^{g+1} y_j^2\Bigr) \sum_{j=1}^{g+1}
x_j{{\partial}\over{\partial y_j}}$$
\noindent so that the canonical equations of motion are 
$${{dx_j}\over{dt}}=Dx_j, \quad  {{dy_j}\over{dt}}=Dy_j \qquad
(j=1, \dots ,g+1).\eqno(7.6)$$
\noindent From the canonical equations, one easily can check that 
$$F_j(x,y)=x_j^2+\sum_{k:k\neq j}{{(y_jx_k-x_jy_k)^2}\over{a_j-a_k}}
\qquad (j=1, \dots , g+1)\eqno(7.7)$$
\noindent give a family of integrals of the motion, and hence that
$F_j(x,y)=c_j$ for some constants $c_j$. Now suppose that the real
polynomial
$$\prod_{j=1}^{g+1}
(s-a_j)\sum_{k=1}^{g+1}{{c_j}\over{s-a_j}}\eqno(7.8)$$
\noindent has real zeros $b_1, \dots , b_g$, and consider the curve
$${\cal E}=\bigl\{ (z,w)\in {\bf C}^2: w^2=\prod_{j=1}^{g+1}
(z-a_j)\prod_{j=1}^g(z-b_j)\bigr\}\cup \{ (\infty , \infty )\}.$$
\noindent The curve ${\cal E}$ is typically hyperelliptic and hence has
a Jacobi variety ${\bf X}$ which is a complex torus of dimension $g$.
The solution curves to the canonical equations are associated with
straight lines in ${\bf C}^g$ of the form $\{ z_0+te: t\in {\bf C}\}$.
Now let $\theta$ be Riemann's theta function for ${\bf X}$, and let
$\Theta =\{ z\in {\bf X}: \theta (z)=0\}$. The vector field $D$, when
complexified, lifts to a derivation $D_{\bf X}$ on the meromorphic
functions on ${\bf X}$; in particular, ${\wp}(z)=D^2_{\bf X}\log \theta
(z)$ is rational on ${\bf X}$ with no poles on ${\bf X}\setminus
\Theta$. Mumford shows that 
$$\wp (z)={{1}\over{8}}\bigl(
\sum_{j=1}^{g+1}a_jx_j^2-y_j^2\Bigr)+\kappa_2$$
for some constant $\kappa_2$. With ${\wp}$ we form a system of coordinates 
for ${\bf X}$; the
map 
$$\phi_n: {\bf X}\setminus \Theta \rightarrow {\bf C}^{2g}:z\mapsto ({\wp}(z),
 D_{\bf X}{\wp}(z) , \dots ,D_{\bf
X}^{2g-1}{\wp}(z))$$
\noindent is an
embedding and the derivatives $D_{\bf X}^j{\wp}(z)$ for $j=0, \dots , 2g-1$
generate the affine ring of ${\bf X}\setminus \Theta$.\par
\rightline{$\square$}\par

\noindent We begin with simple existence result, showing how linear systems in
continuous time give rise to Hankel matrices. Subsequent results will introduce
stronger hypotheses to ensure the existence of Fredholm determinants, and hence $\tau$ functions.\par
\noindent {\bf Proposition 2.1} {\sl Suppose that $H$ is a separable Hilbert
space, and that\par
\indent (i) $C:H\rightarrow {\bf C}$ and $B:{\bf C}\rightarrow H$ are
bounded linear operators;\par
\indent (ii) $A$ is a densely defined linear operator in $H$;\par
\indent (iii) $A$ is accretive, so $\Re \langle Af,f\rangle \geq 0$ for
all $f\in {\cal D}(A)$;\par
\indent (iv) $\lambda I+A$ is invertible for some $\lambda >0$.\par
\noindent Then $(e^{-tA})_{t>0}$ is a $C_0$ contraction semigroup on
$H$, so
$\phi_{(x)} (s)=Ce^{-(2x+s)A}B$ is bounded and continuous on $(0, \infty )$; the
cogenerator $V=(A-I)(A+I)^{-1}$ satisfies $\Vert V\Vert\leq 1$ as an
operator on $H$, and
there is a unitary equivalence between $\Gamma_{\phi_{(x)}}$ on 
$L^2(0, \infty )$ and the Hankel matrix on
$\ell^2({\bf N}\cup \{0\})$ that is given by} 
$$\Gamma_{\phi_{(x)}} \leftrightarrow
\Bigl[\sqrt{2}Ce^{-2xA}V^{n+m}(I+A)^{-1}B\Bigr]_{n, m=0}^\infty .\eqno(2.1)$$ 
\vskip.05in
\noindent {\bf Proof.} (i) By the Lumer--Phillips theorem [21], $-A$ generates a
$C_0$ contraction semigroup. Directly from the definition (iii) of an
accretive operator and hypothesis (iv), one proves that $\Vert V\Vert \leq 1.$ \par
\indent We introduce the Laguerre polynomials of
order zero $L_n^{(0)}(s)=(n!)^{-1}e^{s}(d/ds)^n s^ne^{-s}$ and then the
functions $h_n(s)=\sqrt{2}e^{-s}L_n^{(0)}(2s)$, so that
$(h_n)_{n=0}^\infty$ gives a complete orthonormal basis of $L^2(0,
\infty )$. By integrating by parts, one can verify that
$$\eqalignno{\int_0^\infty \phi
(s)h_n(s)\,ds&={{1}\over{\sqrt{2}n!}}\int_0^\infty
Ce^{-2xA}e^{-(A-I)s/2}B{{d^n}\over{ds^n}}\Bigl( s^ne^{-s}\Bigr)\, ds\cr
&= \sqrt{2} Ce^{-2xA}(A-I)^n(A+I)^{-n-1}B.&(2.2)\cr}$$ 
\noindent Peller [50, p.233] shows that $\Gamma_\phi$ is unitarily equivalent to
the Hankel matrix under the unitary correspondence
$(h_n)_{n=0}^\infty\leftrightarrow (e_j)_{j=0}^\infty$, where $(e_j)$
is the standard orthonormal basis of $\ell^2$.\par
\rightline{$\square$}\par

\indent Now we write
 $$v_{\zeta }(x,t)={{d}\over{dx}}\log\Bigl( {{e^{x\zeta} \tau_{\zeta}(x,t)}\over{\tau_{\infty  }(x,t)}}\Bigr),\quad \psi_\zeta (x,t)={{e^{\zeta x}\tau_\zeta (x,t)}\over{\tau_\infty (x,t)}} .\eqno(5.21)$$ 
\noindent  In the following result we describe the spectral shift on linear systems, and derive an expression for $v_{1/\alpha }$ in terms of Sato's integral.\par
\vskip.05in
\noindent {\bf Proposition 5.2} {\sl (i) Then $\Sigma_\zeta (t)$ satisfies the
hypotheses of Lemma 3.1, and gives rise to a differential ring 
${\bf A}_t$ for the derivatives ${{\partial}\over{\partial x}}$ and
${{\partial }\over{\partial t_j}}$ .}\par
\indent {\sl (ii) The function $\zeta\mapsto v_\zeta (x,t)$ is holomorphic 
except for possible singular points 
for $\zeta$ in the spectrum of $A$ as an operator on} ${\cal L}(H)$. 
{\sl Miura's transformation rule $v_\zeta'+v_\zeta^2=u_\infty+\zeta^2$ is equivalent to $-\psi_\zeta''+u_\infty \psi_\zeta =-\zeta^2 \psi_\zeta$.}\par
\indent {\sl (iii) If $A$ is bounded, then $\Sigma_{1/\alpha}(t)$ with $\vert \alpha\vert <\Vert A\Vert^{-1}$ gives a deformation family such that $\Sigma_{1/\alpha }(t)=\Sigma_\infty (t+\{\alpha  \})$, and Sato's integral satisfies the addition rule}
$$S((x,t)+\{\alpha \}; (x,t)-\{\alpha \})=-{{2}\over{\alpha }}-{{\tau'_{1/\alpha } (x,t)}\over{\tau_{1/\alpha} (x,t)}}+{{\tau'_{-1/\alpha } (x,t)}\over{\tau_{-1/\alpha} (x,t)}}.\eqno(5.22)$$

\vskip.05in
\noindent {\bf Proof.} (i) See Theorem 3.3.\par
\indent (ii) We observe that $\zeta \mapsto \det (I+(\zeta I+A)(\zeta I-A)^{-1}e^{-xA}U(t)EU(t)e^{-xA})$  is differentiable, provided that $(\zeta I-A$ is invertible and for typical $(x,t)$, the derivative will be nonzero. For such $\zeta $ and $(x,t)$, the function $\zeta\mapsto v_\zeta (x,t)$ is differentiable; we do not claim that it is meromorphic when $E$ has infinite rank.  The identity $-\psi_\zeta''+u_\infty \psi_\zeta=-\zeta^2\psi_\zeta$ is equivalent to 
$${{\tau_\zeta''}\over{\tau_\zeta}}+{{\tau_\infty''}\over{\tau_\infty}}+2\zeta \Bigl(
{{\tau'_\zeta}\over{\tau_\zeta}}-{{\tau_\infty'}\over{\tau_\infty}}\Bigr)=0,\eqno(5.23)$$
\noindent which implies
$$\eqalignno{ v_\zeta'+v_\zeta^2&={{\tau_\zeta''}\over{\tau_\zeta}}-{{\tau_\infty''}\over{\tau_\infty }}
+2\Bigl({{\tau_\infty'}\over{\tau_\infty}}\Bigr)^2+\zeta^2+2\zeta\Bigl( {{\tau_\zeta'}\over{\tau_\zeta}}-{{\tau_\infty'}\over{\tau_\infty}}\Bigr)\cr
&=\zeta^2-2{{\tau_\infty''}\over{\tau_\infty}}+2\Bigl({{\tau_\infty'}\over{\tau_\infty}}\Bigr)^2\cr
&=u_\infty +\zeta^2.&(5.24)\cr}$$
\noindent Also note that 
$$v_\zeta (x,t)=\zeta +2{\hbox{trace}}\Bigl( Ce^{-xA}\bigl(I+(\zeta I+A)(\zeta I-A)^{-1}R_x\bigr)^{-1}A(\zeta I-A)^{-1}(I+R_x)^{-1}e^{-xA}B\Bigr).\eqno(5.25)$$
\indent (iv) Suppose that $\tau_\infty (x,t)\neq 0$, and observe that $\zeta\mapsto \psi_\zeta (x,t)$ is holomorphic, except where $\zeta I-A$ is not invertible. For $\vert\zeta \vert >\Vert A\Vert$, we observe that $\zeta I-A$ is invertible and we can choose
$\alpha =1/\zeta$ with $\vert \alpha \vert
<\Vert A\Vert^{-1}$ as the deformation parameter. The Laurent series for
$\log (\zeta I+A)(\zeta I-A)^{-1}$ involves only odd powers $1/(2j+1)\zeta^{2j+1}$. This tau is associated with the Hirota transformation, and also Darboux's
spectral addition rule considered in [22, p. 484].\par
\indent By choosing $(x,t)-(y,s)=2\{ \alpha \}$, we obtain a meromorphic function $\zeta\mapsto \psi_\zeta ((x,t)+\{\alpha \})\psi_{-\zeta}((x,t)-\{ \alpha\})$ on $\Vert A\Vert <\vert\zeta \vert< \vert\alpha\vert^{-1}$, which we continue analytically to $\Vert A\Vert<\vert\zeta\vert$ and have a holomorphic function for $1/\vert \alpha \vert <\vert\zeta\vert$. Hence for $1/\vert\alpha\vert<r$ we obtain the integral
$$S((x,t)+\{\alpha \};(x,t)-\{\alpha\})=\int_{\vert \zeta \vert =r} {{(1+\alpha \zeta )}\over{(1-\alpha \zeta )}}{{\tau_\zeta  ((x,t)+\{\alpha \})\tau_{-\zeta }((x,t)-\{\alpha \})}\over 
{\tau_\infty ((x,t)+\{\alpha\})\tau_{\infty }((x,t)-\{\alpha \})}}{{d\zeta}\over{2\pi i}}\eqno(5.26)$$
\noindent and for $\Sigma_\infty (t)$ we have
$$\eqalignno{ {{\tau_\zeta (x,t)}\over{\tau_\infty (x,t)}}&=\det (I+(\zeta I-A)^{-1}(AR+RA)(I+R)^{-1})\cr
&= 1+\zeta^{-1}{\hbox{trace}}\bigl( (AR+RA)(I+R)^{-1}\bigr) +O(1/\zeta^2 )\cr
&=1-{{1}\over{\zeta}}{{\tau_\infty'(x,t)}\over{\tau_\infty (x,t)}}+O(1/\zeta^2)&(5.27)\cr}$$
\noindent as $\vert \zeta\vert\rightarrow\infty$, where the final step follows from Liapunov's equation. Hence we can compute the residue of the integral at $\zeta=\infty$, and deduce that
$$S((x,t)+\{\alpha\}; (x,t)-\{ \alpha\}) ={{-2}\over{\alpha }}-{{\tau'_\infty ((x,t)}+\{\alpha \})\over{\tau_\infty( (x,t)}+\{\alpha\})}+{{\tau_\infty'((x,t)-\{ \alpha\})}\over{\tau_\infty ((x,t)-\{ \alpha \})}},\eqno(5.28)$$
\noindent and we conclude the proof by using the identity of the linear
 systems $\Sigma_{1/\alpha }(t)=\Sigma_\infty (t+\{ \alpha \})$. We recognise 
the right-hand side as $v_{-1/\alpha }(x,t)-v_{1/\alpha }(x,t)$.\par

\indent Zakharov and Shabat used this method for finite rank $A_1$ and $A_2$ to produce soliton solutions of $KP$.\par 
simplify the expression in ()  by taking 

In Proposition 8.4(iv) we use a variant on their proof to produce meromorphic eigenfunctions arising in the case when $H$ is finite-dimensional, so we can write 
$$\mu I+(\zeta I+A)(\zeta I-A)^{-1}E=\mu I+ E_0+\sum_{k=1}^n \sum_{\ell =1}^{p_k}{{E_{k,\ell}}\over{(\zeta -a_k)^\ell}}\eqno(7.12))$$
\noindent where $a_k$ are the eigenvalues of $A$ and $E_{j,k}$ are matrices on $H$. \par
\indent In the section 8, we consider some linear systems associated with the curve ${\cal E}$ of
(6.2), especially when we deform the $a_j$.\par
\indent {\sl (iii) Suppose that there exist constants $a_j$ such that $a_{2g}\neq
0$ and
$$-24A(A-2FA)^2+8A^3-\sum_{j=0}^{2g}a_jA(AFA+FA^2-2FAFA)^{j}\in
\Theta\eqno(9.32)$$
\noindent Then there
exist $\beta, \nu_j\in {\bf C}$ such that $\sum_{j=0}^{2g+2}1/\nu_j=0$ and} 
$$(u')^2=\beta \prod_{j=1}^{2g+2}(u-\nu_j).\eqno(9.33)$$
\indent (iii) From Lemma 3.2(i), we have the multiplication rules
$$\eqalignno{u(x)^j&=(-4)^{j}\bigl\lfloor
A\bigl((AF+FA-2FAF)A\bigr)^{j-1}\bigr\rfloor,\cr
u'(x)&=-8\bigl\lfloor A(I-2F)A\bigr\rfloor ,\cr
u''(x)&=-24\bigl\lfloor
A\bigl((I-2F)A\bigr)^2\bigr\rfloor+8\bigl\lfloor A^3\bigr\rfloor
,&(9.36)}$$ 
\noindent which in this case give rise to the identity
$$u''+\sum_{j=0}^{2g}{{a_ju^{j+1}}\over{(-4)^{j+1}}}=0.\eqno(9.37)$$
\noindent The polynomial in identity (9.28) has no
constant term, so when we multiply by $u'$ and integrate the resulting polynomial has zero as the 
coefficient of $u$, hence the special form for the roots. \par
\noindent {\bf Definition} (Toda's equation) The sequence $(1, \sigma_1, \sigma_2,
\dots )$ in ${\bf A}$ satisfies Toda's equation if
$${{d^2}\over{d x^2}}\log\sigma_n(x)={{\sigma_{n+1}(x)\sigma_{n-1}(x)}
\over{\sigma_n(x)^2}}\qquad (n=1,2, \dots ).\eqno(4.23)$$ 
\noindent The left-hand side is a multiple of a potential associated with a tau function, whereas the logarithm of right-hand side is the second order difference sequence of $(\log \sigma_n(x))$.\par 

\vskip.05in
\noindent {\bf Proposition 4.5} {\sl Let $\sigma_0=1$ and $P\in
{\bf A}$. Then the sequence given by}
$$\sigma_n(x)={\hbox{trace}}_H \bigl( I_H\ast
\det^*[\partial^{j+k-2}P]_{j,k=1}^n\bigr)\eqno(4.24)$$
\noindent {\sl gives a solution of Toda's equation.}\par  
 \vskip.05in
\noindent {\bf Proof.} Darboux observed that a solution is given by the determinants
$$\sigma_n(x)
=\det\Bigl[ {{d^{j+k-2} \psi (x)}\over{dx^{j+k-2}}}\Bigr]_{j,k=1}^n.
\eqno(4.25)$$
\noindent The entries of the determinant belong to ${\bf S}$, and hence the expansion $\sigma_n$ also belongs
to ${\bf S}$. Letting $\psi =\lfloor P\rfloor$, we use Lemma 3.2, to
find the entries of the determinant and then express them as
$$\det [\psi^{(j+k-2)}]_{j,k=1}^n=\det \bigl[\, \lfloor \partial^{j+k-2}P
\rfloor\bigr]_{j,k=1}^n.\eqno(4.26)$$ 
\par
\rightline {$\square$}\par

\noindent {\bf Proposition 5.4} {\sl (i) Let $(-A_1,B,C)$ and $(-A_2,B,C)$ be $(2,2)$ admissible linear systems as in Lemma 2.2 with input and output spaces ${\bf C}$, where $A_1$ and $A_2$ are bounded on a common state space $H$. Let $C(y;t)=Ce^{t(A_1^3+\lambda A_1)/\alpha -yA_1^2/\beta}$ and $B(y;t)=e^{t(A_2^3+\lambda A_2)/\alpha +yA_2^2/\beta}B$.\par
\indent (i) Then $\psi (x,z;y;t)=C(y;t)e^{-xA_1}e^{-zA_2}B(y;t)$ satisfies the scattering equations (5.30) and (5.31) for the KP equation.\par
\indent (ii) With $S_x=\int_x^\infty e^{-A_2s}B(y;t)C(y;t)e^{-A_1 s}ds,$ let
$$K(x,z;y;t)=-C(y;t)e^{-xA_1}(I+S_x)^{-1}e^{-zA_2}B(y;t).\eqno(5.36)$$
\noindent Then $K(x,z;y;t)$ satisfies the integral equation (5.35) and }
$$K(x,x;y;t)={{d}\over{dx}}\log\det (I+S_x).\eqno(5.37)$$
\noindent {\sl Moreover, there exists $x_0$ such that $K(x,z;y;t)$ satisfies} 
$${{\partial^2K}\over{\partial x^2}}-{{\partial^2K}\over{\partial z^2}}
+\beta{{\partial K}\over{\partial y}}=u(x;y;t)K(x,z;y;t)\qquad 
(x_0<x<z)\eqno(5.38)$$
\indent {\sl (iii) Also, $u(x,y,t)$ satisfies $KP$ when $\beta =i$.}\par
\indent {\sl (iii) Suppose furthermore that $A_1=A_2=A$. Then the 
scattering function satisfies $\psi (x,z;y; t)=\phi (x+z;t)$ with 
$\phi (x;t)=Ce^{2t (A^3+\lambda A)/\alpha -xA}B$, and the
corresponding potential is 
given in terms of a Hankel determinant by}
$$u(x;t)=-2{{d^2}\over{dx^2}}\log\det
(I+\Gamma_{\phi_{(x)}}).\eqno(10.12)$$
\vskip.05in
\noindent {\bf Proof.} (i) Since the operators are all bounded, the functions are differentiable and one can verify the differential equations, without assuming that $A_1$ and $A_2$ commute.\par
\indent (ii) The linear system
$$\Bigl(\left[\matrix{-A_1&0\cr 0&-A_2\cr}\right],
\left[\matrix{B&0\cr 0&B(y;t)\cr}\right],
\left[\matrix{0&-C\cr C(y;t)&0\cr}\right]\Bigr)\eqno(5.39)$$
\noindent is $(2,2)$ admissible and by Lemma 2.2 the resolvent operator
$$\left[\matrix{ 0&-\int_x^\infty e^{-sA_1}BCe^{-sA_2}ds\cr 
\int_x^\infty
e^{-sA_2}B(y;t)C(y;t)e^{-sA_1}ds&0\cr}\right],\eqno(5.40)$$
\noindent is trace class as in Proposition 2.5. Hence one can
verify the integral equation as in the proof of Proposition 2.5. 
One then verifies the determinant identity (5.19), which  involves the bottom left entry $S_x$ satisfying  the asymmetric Lyapunov equation 
$${{d}\over{dx}}S_x=-A_2S_x-S_xA_1=-e^{-xA_2}B(y;t)C(y;t)e^{-xA_1}
\qquad (x>0).\eqno(5.41)$$
\indent The solution of the integral equation is unique for large enough $x$ since $\Vert e^{-xA_1}\Vert \rightarrow 0$ and $\Vert e^{-xA_2}\Vert \rightarrow 0$ exponentially fast as $x\rightarrow\infty$; hence $\Psi (x;z;y;t)\rightarrow 0$  exponentially fast as $x\rightarrow\infty$. Using the scattering equation (5.30), one shows by differentiating (4.35) repeatedly that ${{\partial^2}\over{\partial x^2}}K(x,z;y;t)-{{\partial^2}\over{\partial z^2}}K(x,z;y;t)-\beta 
{{\partial }\over{\partial y}}K(x,z;y;t)$  and $u(x;y;t)K(x,z;y;t)$ both satisfy (3.5) multiplied by $u(x;y;t)$, and so by uniqueness are equal.\par 
\noindent (iii) We introduce the function
$$\phi (x,y,t)=e^{ikt-ik^2y}+\int_x^\infty K(x,z,t;t)e^{ikz-ik^2y}\, dz\eqno(4.42)$$
\noindent and show by integrating by parts that (4.38) leads to
$$i{{\partial \phi}\over{\partial y}}+{{\partial^2\phi}\over{\partial x^2}}=u\phi.\eqno(4.43)$$
\noindent Thus we identify $K$ with the triangular factor of the scattering operator, as considered by Manakov. Hence the diagonal $K(x,x,y,t)$ has derivative which gives $u(x,y,t)$, namely a solution of $KP$.\par
\indent (iii) The determinant $\det (I+S_x)$ of (ii) is not generally a tau function in the strict sense of Theorem 1.1, but in the following case where $A_1=A_2$  we do obtain a genuine tau function from a linear system with input and output space 
$$R_x=\int_x^\infty  e^{t(A^3+\lambda A)/\alpha}e^{-sA}BCe^{-sA}
e^{t(A^3+\lambda A)/\alpha}ds,\eqno(10.16)$$
\noindent hence $\det (I+S_x)=\det (I+\Gamma_{\phi_{(x)}})$ by 
Proposition 2.4.\par
\noindent {\bf Proposition 5.3} {\sl Let 
$\Sigma_j(t)=(-A_j, U_j(t/2)B_j, C_jU_j(t/2))$ be a $(2,2)$ admissible linear 
system as in Proposition 3.6 with scattering function
 $\psi_j(x,t)$, and let $\tau_n(x,t)
=Wr(\psi_1,\psi_2, \dots ,\psi_n$; let 
$u(x,t)=-2{{\partial^2}\over{\partial x^2}}\log \tau_n(x,t)$. 
Then $u(x,t)$ satisfies the $KP$ equation.}\par 
\vskip.05in
\noindent {\bf Proof.} The scattering functions $\psi_k(x;t)$ satisfy
$${{\partial \psi_k}\over {\partial t_3}}=
{{\partial^3 \psi_k}\over {\partial t_1}};\qquad
{{\partial \psi_k}\over {\partial t_2}}=
-{{\partial^2 \psi_k}\over {\partial x^2}}\eqno(5.33)$$
so by a result of Freeman and Nimmo [24], $\tau_n$ satisfies Hirota's 
bilinear form of the KP equation
$$\eqalignno{-3\Bigl({{\partial \tau_n}\over{\partial y}}\Bigr)^2
+3\Bigl({{\partial^2 \tau_n}\over{\partial t_1^2}}\Bigr)^2&
+3\tau_n\Bigl({{\partial^2 \tau_n}\over{\partial t_3^2}}\Bigr)
+4\Bigl({{\partial \tau_n}\over{\partial t_2}}\Bigr)\Bigl({{\partial 
\tau_n}\over{\partial x}}\Bigr)\cr
-4\tau_n\Bigl({{\partial^2 \tau_n}\over{\partial t_1\partial t_2}}\Bigr)^2
&-4\Bigl({{\partial \tau_n}\over{\partial t_1}}\Bigr)\Bigl(
 {{\partial^3\tau_n}\over{\partial t_1^3}}\Bigr)+
\tau_n\Bigl({{\partial^4\tau_n}\over{\partial
t_1^4}}\Bigr)=0&(5.34)\cr}$$  
\noindent and hence $u$ satisfies $KP$.\par

\indent For the linear system $\Sigma_1=(-A,B, Ce^{t(A^3+\lambda
A)/\alpha +yA^2\beta})$ we introduce the scatttering function
$\phi_1 (x;y,t)$ and for $\Sigma_1=(-A,B, Ce^{t(A^3+\lambda
A)/\alpha -yA^2\beta})$ we introduce $\phi_2(x;y,t);$ then we introduce
$$\psi (x,z;y,t)=\int_0^\infty \phi_2(x+s; y,t)\phi_1(s+z;y,t)\, ds;$$
\noindent note that $\psi $ is not a function of $x+z$. Next let
$$(-\hat A, \hat B, \hat C)=\Bigl( \left[\matrix{-A& 0\cr
0&-A\cr}\right],\left[\matrix{B& 0\cr
0&B\cr}\right], \left[\matrix{0& Ce^{t(A^3+\lambda A)/\alpha
-yA^2/\beta }\cr
-Ce^{t(A^3+\lambda A)/\alpha +yA^2/\beta}&0\cr}\right]\Bigr);$$
\noindent then let $\Phi (x)=\hat C e^{-\hat Ax}\hat B$ and $\hat
R=\int_x^\infty e^{-s\hat A}\hat B\hat C e^{-s\hat A}\, ds,$ and 
$$\hat T (x,z)=-\hat C E^{-x\hat A}(I+\hat R_x)^{-1} e^{-\hat A z} \hat
C,$$
\noindent which we write as
$$\hat T (x,z)=\left[\matrix{K(x,z;y,t)& L(x,z;y,t)\cr M(x,z;y,t)&N(x,z;y,t)\cr}\right].$$
\noindent {\bf Proposition} {\sl (i) Then $\psi(x,z;y,t)$ satisfies
the linear differential equations () and ();\par
\indent (ii) $\hat T(x,z;y,t)$ satisfies the Gelfand-Levitan equation ();\par
\indent (iii) the integral operator $\Psi_{(s)}$ with kernel
$\psi (x+s,z+s;y,t)$ on $L^2(0, \infty )$ satisfies }
$${{d}\over{ds}} \log\det (I+\Psi_{(s)})={\hbox{trace}}\, T(s,s;y,t);$$
\indent (iv) {\sl $u(x,y,t)=-2{{\partial}\over{\partial x}}K(x,x;y,t)$ satisfies (KP).}\par  
\vskip.05in
\noindent {\bf Proof.} (i) This is a simple calculation.\par
\indent (ii) This is similar to Theorem 3.1.\par
\indent Note that 
$$ \Phi (x)=\left[\matrix{0& \phi_2(x)\cr
-\phi_1 (x)&0\cr}\right]$$
\noindent so that $\Psi_{(s)} =\Gamma_{\phi_{2, (s)}}\Gamma_{\phi_{1,
(s)}}$ and 
$$\eqalignno{\det (I+\hat R_s)&=\det
(I+\Gamma_{\Phi_{(s)}})\cr
&=\det\Bigl(\left[\matrix{I& \Gamma_{\phi_{2,
(s)}}\cr
-\Gamma_{\phi_{1, (s)}}&I\cr}\right]\Bigr)\cr 
&=\det (I+\Gamma_{\phi_{2, (s)}}\Gamma_{\phi_{1, (s)}})\cr
&=\det (I+\Psi_{(s)}).\cr}$$
\indent (iv) Let $K(x,z;y,t)$ and $L(x,z;y,t)$ be the top row of $\hat T(x,z)$, which satisfy the
coupled integral equations
$$K(x,z;y,t)-\int_x^\infty L(x,w)\phi_1(w+z;y,t)\, dw=0$$
$$\phi_2(x+z;y,t)+L(x,z;y,t)+\int_x^\infty K(x,w;y,t)\phi_2(w+z;y,t)\, dw=0$$
\noindent which combine to give
$$\eqalignno{K(x,z;y,t)+&\int_x^\infty \phi_2 (x+w;y,t)\phi_1(w+z;y,t)\, dw\cr
&+\int_x^\infty K(x,w;y,t)\int_x^\infty
\phi_2(w+s;y,t)\phi_1(s+z;y,t) \, dsdw=0,\cr}$$
\noindent or
$$0=K(x,z;y,t)+\psi (x,z;y,t)+\int_x^\infty K(x,w;y,t)\psi (w,z; y,t)\, dw=0.$$
\noindent Then $u(x,y,t)=-2{{\partial }\over{\partial x}} K(x,x;y,t)$ satisfies KP by the 
general case (37) of
Zahkarov and Shabat's paper.\par

 \noindent {\bf Definition} (Deformation of differential equations) We use $d$ to denote the exterior complex  derivative on differential forms on $({\bf T}^\infty)_0$ that depend on only finitely many coordinates. Let $B$ be an $n\times n$ matrix function with entries that are meromorphic functions of $(x,t)$; let $W_j$ be $n\times n$ matrices with entries that are meromorphic functions of $(x,t)$ for $j=1, \dots ,N$ and such that the differential form $\Omega =\sum_{j=1}^N W_jdt_j$ satisfies $d\Omega =\Omega\wedge \Omega.$ Then the differential equation ${{\partial }\over{\partial x}}Y=BY$ has a system of deformations $dY=\Omega Y$. In this case, the pair of differential equations are said to be consistent. \par
\vskip.05in

\noindent {\bf Proposition 4.7} {\sl Let $\tau (z;t)$ be the tau function of $\Sigma (t)$.\par
\indent (i)  Then the simple zeros of $z\mapsto \tau (z;t)$ depend
holomorphically upon each $t_j$.\par 
\indent (ii) Suppose that $x_0$ is a regular point of
$u(x,0)$. Then there exists a matrix  $W_j$ such that the pair of ordinary differential equations
$${{d}\over{dx}}Y=\left[\matrix{0&1\cr u(x,t)-\lambda &0\cr}\right]Y
,\quad {{d}\over{dt_j}}Y=W_j Y\eqno(4.30)$$
\noindent is consistent and have holomorphic solutions on $\{ (x,t_j): \vert x-x_0\vert , \vert
t_j\vert <\delta_1\}$ for some $\delta_1>0$.}\par
\indent {\sl (iii) There exists a nonconstant meromorphic function $w$ such that the system (3.22) with $\Omega =W_1dt_1$ is consistent, where}
$$W_1=\left[\matrix{ {{-1}\over{2}}
{{\partial w}\over{\partial x}}& w\cr 
{{-1}\over{2}}{{\partial^2w}\over{\partial x^2}} 
+w(u-\lambda )& {{1}\over{2}}
{{\partial w}\over{\partial x}}\cr}\right].\eqno(4.31)$$ 
\vskip.05in

\noindent{\bf Proof.} (ii) By Proposition
2.4(iii), $\tau (x,t)$ depends uniformly continuously upon $t_j$ on compact subsets of $x\in {\bf C}$. So by Hurwitz's theorem, given an isolated zero
$z_1$ of $\tau (z;0)$ and suitably small $\varepsilon >0$, there exists
$\delta_0>0$ such that $\tau (z;t_j)$ and $\tau (z;0)$ have the same
number of zeros in $\{ z: \vert z-z_1\vert <\varepsilon\}$ for all 
$\vert t_j\vert <\delta_0$. For a simple zero $z_1(t_j)$, we can write
$$z_1(t_j)={{1}\over{2\pi i}}
\int_{\vert z-z_1(0)\vert=\varepsilon /2} {{z\tau'(z;t)dz}\over{\tau
(z;t)}}.\eqno(4.32)$$
\indent (ii) Thus the positions of the poles of $u(x,t)$
depend continuously on $t_j$. Hence we can consider $\delta_1 >0$ such that $x\mapsto u(x,t)$ has no poles on 
 $\{ (x,t): \vert x-x_0\vert , \vert
t\vert <\delta_1\}$, and solve the first differential equation to obtain an invertible matrix function
$Y(x,t)$ which is holomorphic in $(x,t_j)$. Then we define
 $W (x,t)=({{\partial Y}/{\partial t_j}})Y^{-1}$
on  $\{ (x,t): \vert x-x_0\vert , \vert
t_j\vert <\delta_1\}$.\par 
\indent (iii) The consistency condition reduces to the linear partial differential equation
$${{\partial u}\over{\partial t_1}}=w{{\partial u}\over{\partial x}}+2(u-\lambda ){{\partial w}\over{\partial x}}-{{1}\over{2}}{{\partial^3w}\over{\partial x^3}},$$ 
\noindent which was previously derived by L. and R. Fuchs. For each fixed $t$, this is a linear ordinary differential equation for $w$, and lifts to the differential equation
$$-8A(I-2F)A=-(8A+2\lambda I)\ast (\partial W)-8A(I-2F)A\ast W-2^{-1}\partial^3 W$$
\noindent in $({\bf A}, \ast, \partial )$. On a neighbourhood on which $F$ is holomorphic, one can integrate this differential equation and obtain a holomorphic solution by Picard iteration.\par

\par
\rightline{${\square}$}\par

\indent From $\phi (t)=Ce^{-tA}B$ in $L^2(0, \infty )$, we introduce
$$s_{21}(k)=\int_0^\infty e^{ikt}\phi (t)\, dt$$
\noindent and normalize it so that $\vert s_{21}(k)\vert\leq 1.$ Then $s_{12}$
determines a function in $H^2$ of the upper half plane, by the Paley--Wiener
theorem. We introduce an outer function $s_{11}$ be defining $\vert
s_{11}(k)\vert=(1-\vert s_{21}(k)\vert^2)^{1/2}$, and then letting
$$s_{11}(k)=\exp\Bigl(\log \vert s_{11}(k)+{{i}\over{\pi}}{PV}
\int_{-\infty}^\infty {{\log \vert s_{11}(t)\vert \, dt}\over{k-t}}\Bigr)$$
\noindent which externds to define a function in $1+H^2 $. \par
\indent Next we introduce $L>0$ and build the periodic function
$q(x)=\sum_{j=-\infty}^\infty u(x+jL)$. Then the system of differential equations
$${{d}\over{dx}}\left[\matrix{f\cr g\cr}\right]=\left[\matrix{0&1\cr
q(x)-\lambda &1\cr}\right] \left[\matrix{f\cr g\cr}\right]$$
\noindent has a fundamental solution matrix $F_\lambda $ such that
$F_\lambda(0)=I$. The
discriminant is $\Delta_L (\lambda )={\hbox{trace}}\, F_\lambda (L)$. 
$$\Delta_L (\lambda )={{e^{ikL}}\over{\bar
s_{11}(k)}}+{{e^{-ikL}}\over{s_{11}(k)}}.$$

\noindent {\bf Proposition 4.4} {\sl (i) Let ${\bf A}_u$ be the complex subalgebra of $\lfloor {\bf
A}\rfloor$ generated by $u$ and its
derivatives, and suppose that $u^{(k)}=P_k(u,u', \dots , u^{(k-1)})$ for some complex
polynomial $P_k$. Then there exists an irreducible affine variety $V$ in ${\bf
C}^{k+1}$ such that ${\bf A}_u$ is naturally isomorphic to the coordinate ring ${\bf
C}[V]$.}\par
\indent {\sl (ii) Suppose that  $A$ is an algebraic operator. Then $\phi$ satisfies
a linear differential equation with constant coefficients and belongs
to some differential field ${\bf K}_1$ that satisfies an addition rule.}\par 
\vskip.05in
\noindent {\bf proof.} (i) Since $u$ is meromorphic on a connected set $\Omega$, the algebra ${\bf
A}_u$ is an integral domain. By mapping $X_j\mapsto u^{(j-i)}$ for $j=1, \dots , k+1$,
we obtain a short exact sequence of algebra homomorphisms $0\rightarrow J\rightarrow
{\bf C}[X_1, \dots , X_{k+1}]\rightarrow {\bf A}_u\rightarrow 0$, where $J$ is a
nontrivial prime ideal, and hence associated with the irreducible 
affine variety $V=\{
z\in {\bf C}^{k+1}: f(z)=0, \forall f\in J\}$. \par
\indent (ii) Let $a_j$ be complex numbers such that
$A^n+a_{n-1}A^{n-1}+\dots +a_0I=0.$ Then by a simple calculation, we
have
$$\phi^{(n)}(x)-a_{n-1}\phi^{(n-1)}(x)+\dots +(-1)^{n-1}a_0\phi
(x)=0.\eqno(4.22)$$
\noindent Note that $R_z=e^{-zA}R_0e^{-zA}$, where $z^{-zA}$ reduces to a polynomial expression in $A$ with coefficients that are entire functions of $z$. We choose a largest $N$ such that $\phi, \phi', \dots
,\phi^{(N-1)}$ are linearly independent, and build the matrix $\Phi
(x)=[\phi^{(j+k-2)}(x)]_{j,k=1,\dots , N}$ which satisfies $\Phi'
(x)=A_0\Phi (x)$ for some constant $N\times N$ matrix $A_0$, and has $\Phi
(0)$ invertible. Then we write $\Phi (x)=e^{xA_0}\Phi (0)$, so 
$\Phi (x+y)=e^{xA_0}e^{yA_0}\Phi (0)=\Phi (x)\Phi (0)^{-1}\Phi (y)$. Thus
the 
differential field ${\bf K}_1={\bf C}(\phi, \phi', \dots ,\phi^{(N-1)})$
satisfies an addition rule.\par 
\rightline{$\square$}\par
\noindent {\bf Remark} (i) The conclusions of Proposition 4.4 are satisfied in the context of certain rational differential
 equations, according to a theorem of Halphen; see [Ince 15.5]. Suppose that 
$${{d^n}\over{dx^n}}y+a_{n-1}(x){{d^{n-1}}\over{dx^{n-1}}}y+\dots +a_1(x){{d}\over{dx}}y+a_0(x)y=0,\eqno(3.17)$$
\noindent where (a) the $a_j$ are proper complex rational functions;\par
\indent (b) $a_j(x)$ has poles of order less than or equal to $n-j$ at all $x\in {\bf C}$, so finite singular points of the 
differential equation are regular;\par
\indent (c) the general solution $y(x;c)$ is meromorphic for all constants $c\in {\bf C}^n$, so (3.17) is a Picard system.\par
\noindent Then there exists a nonzero polynomial $p$ such that $\phi (x;c)=p(x)y(x;c)$ is entire for all $c$ and satisfies a 
nontrivial linear differential equation with constant coefficients, and ${\bf K}_1$ is a subfield of 
${\bf C}(x, e^{\lambda_1x}, \dots ,e^{\lambda_nx})$ for some $\lambda_j\in {\bf C}$. See [Singer, p53]. 
We consider this situation in detail in section 4.\par 

\indent  (ii) Suppose that ${\bf A}_0={\bf A}/\Theta$ is a commutative algebra, as in Lemma 3.2(iii). Let ${P}$ be the set of prime ideals of ${\bf A}_0$,   For $E\subseteq {P}$, let $V(E)=\{ {\bf p}\in {P}:E\subset  {\bf p}\}$, and let $\{ {P}\setminus V(E): E\subseteq {P}\}$ be the collection of open sets 
associated with the Zariski topology of ${\bf A}_0$ on ${P}$.  In particular, suppose that ${\bf A}_0$ is an algebra of bounded holomorphic functions on some domain $\Omega$ as in Theorem 3.3(iii),
and for each $\omega \in \Omega$  let ${\bf p}(\omega )$ be the prime ideal given by the null space of the homomorphism  $g\mapsto g(\omega ):$ ${\bf A}_{0}\rightarrow {\bf C}$; thus ${\Omega}$ may be regarded as a subset of ${P}$. For any 
net $(\omega_\alpha)$ in $\Omega$;  then ${\bf p}(f)=\lim_\alpha f(\omega _\alpha )$ determines a homomorphism ${\bf A}_0\rightarrow {\bf C}$ such that the null space is either a maximal ideal or all of ${\bf A}_0$. \par 
\indent (iii) In an important special case, the maximal ideal space of ${\bf A}_0$ contains the open right half plane.
 Let $A=A^\dagger\geq 0$ and $B=C^\dagger $ give a linear system with state space $H$ and input space ${\bf C}$, and 
suppose that $AR_0+R_0A=BC$ has a solution $R_0\in {\cal L}(H)$. Then by Theorem 3.2 of [PMT], 
 the integral $\int_0^\infty e^{-tA}BCe^{-tA}\, dt$ is weakly convergent, and equals $R_0$; furthermore, 
 we obtain a solution to Lyapunov's equation from $R_z=e^{zA}R_0e^{-zA}$ such that $F_z=(I+R_z)^{-1}$ is bounded and 
holomorphic on $\{ z\in {\bf C}: \Re z> 0\}$.\par

\vskip.1in
\indent (vi) First we introduce an operator-valued  counterpart of the Weierstrass
zeta function by considering the integral equation
$$Z_x=-\int_0^x F_t^{-1} dt-\int_0^x (A(I-2F_t)Z_t+Z_t(I-2F_t)A)\,
dt,\eqno(9.22)$$
\noindent and observing that by Gronwall's inequality, any continuous
solution $Z:[0, \pi]\rightarrow {\cal L}(H)$ satisfies the inequality
$$\Vert Z_x\Vert_{B(H)}\leq \Bigl( \int_0^x \Vert
F_t^{-1}\Vert_{B(H)}\, dt\Bigl) \exp \Bigl( 
2\Vert A\Vert_{B(H)}\int_0^x \Vert I-2F_t\Vert_{B(H)}\, dt\Bigr).\eqno(9.23)$$
\noindent Using Banach's fixed point theorem, one can obtain a solution
of this integral equation in $C([0, 2\pi ]; {\cal L}(H))$. Noting that
$\partial$, $\ast$ and $\lfloor\, ,\lfloor$ extend to 
$C([0, \pi ]; {\cal L}(H))$, one checks that $\partial Z=-F^{-1},$
so $(d/dx)\lfloor Z\rfloor =\tau_\infty'/\tau_\infty$. Now we introduce the
modified exponential series $\Upsilon =Z+Z\ast
Z/2!+Z\ast Z\ast Z/3!+\dots, $ which converges in $C([0,
\pi ]; {\cal L}(H))$, and satisfies $\Upsilon_0=Z_0=0$ and 
$\lfloor \Upsilon\rfloor =\exp \lfloor Z\rfloor -1=(\tau_\infty (x)-\tau_\infty (0)/\tau_\infty
(0).$ We do not assert that $Z$ or $\Upsilon$ is periodic;
however, $\log \tau_\infty (x)$ is $\pi$-periodic, so $\lfloor Z_{x+\pi
}-Z_x\rfloor =0$.\par

indent {\sl (v) Suppose that $F_x$ is continuous on $[0, \pi ]$.
Then there exists a bounded and continuous function
$\Upsilon:[0, \pi ]\rightarrow {\cal L}(H)$ such that} 
$$\lfloor \Upsilon\rfloor ={{\tau_\infty (x)-\tau_\infty (0)}\over{\tau_\infty (0)}}.\eqno(9.7)$$

\indent {\sl (ii) the potential of $\Sigma$ is given by a second order difference}
$$-2{{d^2}\over{dx^2}}\log \tau_\infty (x;t)=-2\lim_{\zeta\rightarrow\infty}\zeta^2\Bigl( 
{{\tau_\infty (x;t)\tau_{\zeta ; \zeta}(x;t)-\tau_\zeta (x,t)^2}\over{\tau_\infty (x;t)^2}}\Bigr).\eqno(9.6)$$
\indent (ii)  To deal with $\zeta\rightarrow \infty$, we introduce
$$X_\infty =\bigl(AR+RA\bigr)\bigl(I+R)^{-1},\eqno(9.9)$$
$$X_\zeta =\bigl( A(\zeta I+A)(\zeta I-A)^{-1}R+(\zeta I+A)(\zeta I-A)^{-1}RA\bigr)
\bigl( I+(\zeta I+A)(\zeta I-A)^{-1}R\bigr)^{-1}\eqno(9.10)$$
\noindent such that $\zeta (X_\zeta -X_\infty )$ is a bounded family of trace class operators as $\zeta\rightarrow\infty$. Now    
from (4.27) and obtain
$$\eqalignno{ {{\tau_\zeta (x,t)}\over{\tau_\infty (x,t)}}&=\det (I+(\zeta I-A)^{-1}(AR+RA)(I+R)^{-1})\cr
&=\det (I+\zeta^{-1}X_\infty )\det (I +\zeta^{-2} AX_\infty )+O(1/\zeta^3)\cr
&=1+\zeta^{-1}{\hbox{trace}} X_\infty+\zeta^{-2}{\hbox{trace}}(X_\infty^{\wedge 2})+\zeta^{-2}{\hbox{trace}}(AX_\infty )+O(\zeta^{-3})\cr
&=1-{{1}\over{\zeta}}{{\tau_\infty'(x,t)}\over{\tau_\infty (x,t)}}+\zeta^{-2}{\hbox{trace}}(X_\infty^{\wedge 2})+
\zeta^{-2}{\hbox{trace}}(AX_\infty )+O(\zeta^{-3}).&(9.11)\cr}$$
\par
\noindent and likewise
$$ {{\tau_{\zeta,\zeta } (x,t)}\over{\tau_\zeta (x,t)}}=1-{{1}\over{\zeta}}
{{\tau_\infty'(x,t)}\over{\tau_\infty (x,t)}}+\zeta^{-2}{\hbox{trace}}
(X_\zeta^{\wedge 2})+\zeta^{-2}{\hbox{trace}}(AX_\zeta )+O(\zeta^{-3}),\eqno(9.12)$$
\noindent in which we can replace the terms in $X_\zeta$ by $X_\infty$  at the cost of adjusting the term in $O(1/\zeta^2)$. Hence we obtain 
$$\eqalignno{-{{1}\over{\zeta }}{{d}\over{dx}}{{\tau'_\infty}\over{\tau_\infty}}&={{\tau_\zeta}\over{\tau_\infty}}\Bigl( {{\tau'_\zeta}\over{\tau_\zeta}}-{{\tau'_\infty}\over{\tau_\infty}}\Bigr)+O(1/\zeta^2)\cr
&=\zeta {{\tau_\zeta}\over{\tau_\infty}}\Bigl( {{\tau_\zeta}\over{\tau_\infty}}-
{{\tau_{\zeta ,\zeta}}\over{\tau_\zeta}}+O(1/\zeta^3)\Bigr) +O(1/\zeta^2 ).&(9.13)\cr}$$

\noindent {\bf Example} Suppose that $FA=AF$ and $A^2=-I$. It follows that $e^{-\pi A/2}=A$, so we
can apply Theorem 9.3 (v). Then
$v=\lfloor A\rfloor$ satisfies
$$v^2 =-2\lfloor AF(I-F)\rfloor,\quad v^3 =4\lfloor
AF^2(I-F)^2\rfloor,$$
$$v'=-2\lfloor I-2F\rfloor, (v')^2 =8\lfloor AF(I-F)(I-2F)^2\rfloor,$$
\noindent so 
$$(v')^2=-4v^2-8v^3;$$
\noindent let $u=-2v$, so $(u')^2=u^3-u^2$. The scalar solution of this
differential equation is found to be 
$$u(x)={\hbox{cosec}}^2\, (x+c).$$
\noindent Now 
$$\sum_{n=-\infty}^\infty \pi {\hbox{cosec}}^2\pi
(z-in)=\pi^2/3+\sum_{n=\infty ; n\neq 0}^\infty \pi {\hbox{cosec}}^2\,\pi
in+\wp (z)$$,
\noindent for the lattice $\pi{\bf Z}+\pi \pi {\bf Z}$, where 
$$(\wp' (z))^2=(\wp (z) -e_3)(\wp (z)-e_2)(\wp (z)-e_1).$$
\vskip.05in


\indent We introduce the linear systems
$$\Sigma (t) =(-A, U(t)B,CU(t);U(t)EU(t) ),\quad \Sigma^\dagger (t) =(-A^\dagger, U(\bar t)^\dagger C^\dagger, B^\dagger U(\bar t)^\dagger; U(\bar t)^\dagger E^\dagger U(\bar t)^\dagger ).\eqno(7.11)$$ 
\vskip.05in
\indent (ii) Suppose further that the tau functions of $\Sigma (t)$ and 
$\Sigma^\dagger (t)$ are equal and  $\tau_\infty (x,t)$ has no zeros on ${\bf R}$ for real $t$.
 Then $\tau_{\bar\zeta}(\bar x, \bar t)=\overline{\tau_\zeta (x,t)}$ and $u_\infty (x,t)$ 
$$\overline{\tau_{\bar\zeta}(\bar x,\bar t)}
=\det (I+(\zeta I +A^\dagger)(\zeta I-A^\dagger )^{-1}U(\bar t)^\dagger e^{-xA^\dagger}
E^\dagger U(\bar t)^\dagger e^{-xA^\dagger} )=\tau_\zeta (x,t)\eqno(7.15)$$ 
\indent Let $t=(t_1,t_2, \dots )\in ({\bf T}^\infty)_0$ and let $(-A,B,C;E)$ be a periodic linear system such that\par
\noindent $e^{-\pi A/2}Ee^{-\pi A/2}=-E$; let $U(t) =\exp (-\sum_{j=1}^\infty t_jA^{2j-1})$ and let 
$$\Sigma (t)=(-A, U(t)B, CU(t); U(t)EU(t))\eqno(7.34)$$ 
have tau function  
 $$\tau_0 (t)=\det (I+U(t)EU(t))\eqno(7.35)$$ 
and $\tau_1(t_1, t_2, \dots )=\tau_0(t_1+\pi /2, t_2, \dots );$ then let 
$$v(t)={{\partial}\over{\partial t_1}}\log {{\tau_1(t)}\over{\tau_0(t)}};\quad u(t)=-2{{\partial^2}\over{\partial t_1^2}}\log\tau_0(t);\eqno(7.36)$$
\noindent and consider Hill's equation
$$-{{\partial^2}\over{\partial t_1^2}}f(t;\zeta )+u(t)f(t;\zeta )=\zeta^2f(t;\zeta ).\eqno(7.37)$$
\vskip.05in
\noindent {\bf Definition} If  () has characteristic function $\Delta (\lambda )$ which is independent of $t$, then say that the family $\Sigma (t)$ of linear systems is isospectral.\par
\vskip.05in
Let $G_1=A$ and suppose that the $G_m$ satisfy the recurrence relation 
$$\partial G_{j+1}= 8(G_1\ast \partial G_m)+8\partial (G_1\ast G_m)+\partial^3G_m.\eqno(7.38)$$
\vskip.05in
\noindent {\bf Proposition 7.7} {\sl There exists a sequence of universal noncommutative polynomials $(Q_j)_{j=0}^\infty $ such that $Q_j= 8\partial_jA+\partial G_{j+1}$ and 
$\Sigma (t)$ is isospectral if and only if $Q_j(A,F)\in \Theta $ for all $j$.}\par
\vskip.05in

\noindent {\bf Proof.} We introduce the differential operators  $L=-{{\partial^2}\over{\partial x^2}}+u(t)$ and 
$$P_{2n+1}=\sum_{j=0}^n \Bigl( g_{n-j}{{\partial}\over{\partial x}}-2^{-1}{{\partial g_{n-j}}\over{\partial x}}\Bigr) (-L)^j.\eqno(7.39)$$
\noindent Then the isospectral condition for deformation with respect to parameter $t_j$ is
$$[{{\partial}\over{\partial t_j}}+P_{2n+1}, L]=0,\eqno(7.40)$$ 
\noindent which reduces to 
$${{\partial u}\over{\partial t_j}}=2{{\partial g_{j+1}}\over{\partial x}}.\eqno(7.41)$$
In terms of the algebra ${\bf A}$, we have $u=-4{\hbox{trace}}\lfloor A\rfloor $ and $ {{\partial g_{j+1}}\over{\partial x}}={\hbox{trace}}\lfloor \partial G_{j+1}\rfloor .$ Now
$\partial_j A=-A^j(I-2F)A-A(I-2F)A^j$, so there exists a universal noncommutative polynomial $Q_j$ in two variables such that $Q_j(A,F)=8\partial_jA+G_{j+1}$, and () holds if and only if 
${\hbox{trace}}\lfloor Q_j(A,F)\rfloor =0$, that is $Q_j(A,F)\in \Theta.$ 
\par
\rightline{$\square$}\par
\noindent See [Gesztesy] for discussion.\par

\vskip.05in

\noindent {\bf Proposition 9.7} {\sl Let $(-A,B,C;E)$ be a 
uniform periodic linear system as in Theorem 4.2, and let $\tau_\lambda$ be the
tau function of 
$\Sigma_\lambda$.}\par
\indent {\sl (i) 
$\tau_\infty\in {\bf H}_{\cal C}$ satisfies $\log_+\log_+
\vert\tau_\infty (z)\vert\leq 2N\vert z\vert +c_1$ for some $
c_1>0$
and all $z$, where $N$ is the spectral radius of $A$.}\par
\indent {\sl (ii)  Let $(\tau_\lambda) 
=\{ z\in {\bf C}: \tau_\lambda (z)=0\}$ for all $\lambda\in (-\infty ,\infty
)\cup\{\pm\infty\}$, which is either empty or countably infinite.
Every zero of $\tau_\lambda$ gives rise to a double
pole of $u_\lambda =-2(\log \tau_\lambda )''.$}\par
\indent {\sl (iii) If $E$ has finite rank, then $\tau_\infty$ is of
exponential type and in ${\bf C}_{\cal C}$. Conversely, if 
$\tau_\infty$ is of exponential type, then there 
exist $\alpha_j\in {\cal C}$, $\alpha\in {\bf Z}$ and $\beta\in {\bf C}$
such that
$$\tau_\infty (z)=e^{2i\alpha z+\beta}\prod_{j=1}^m \sin 2(z-\alpha_j)
\eqno(9.39)$$
\noindent and} 
$$u(z)=\sum_{j=1}^m {{8}\over{\sin^22(z-\alpha_j)}}.\eqno(9.40)$$

\vskip.05in
\noindent {\bf Proof.} (i) Let $a_j(V)$ be the approximation numbers of a compact operator $V$. Then we have
$a_n(e^{-zA}Ee^{-zA})\leq \Vert e^{-zA}\Vert^2a_n(E)$ and hence by a standard
bound on the determinant 
$$\log \vert\det (I+e^{-zA}Ee^{-zA})\vert \leq c_0e^{2N\vert
z\vert}\sum_{j=1}^\infty a_j(E).\eqno(9.41)$$
\indent (ii) If $\tau_\lambda (z)=0$, then $\tau_\lambda (z+k\pi )=0$
for all $k\in {\bf Z}$. \par
\indent By computing $u=-2(\log \tau_\infty )''$, we obtain a potential as
in (9.10), which is a rational
function of $e^{ix}$ and $e^{-ix}$. 
\vskip.05in
\noindent We now introduce the analogue for periodic linear systems 
of the Darboux transform. Let $(-A,B,C;E)$ be a periodic linear system as in
Theorem 5
.2, and let
$$T(x,y)=\left[\matrix { W(x,y)&V(x,y)\cr
V(x,y)&W(x,y)\cr}\right],\quad  \Phi (x)=\left[\matrix{0&\phi (x)\cr \phi (x)&0\cr}\right]\eqno(9.43)$$
\noindent  where 
$$V(x,y)=-Ce^{-xA}(I-(e^{-xA}Ee^{-xA})^2)^{-1}e^{-yA}B,\eqno(9.44)$$
$$W(x,y)=Ce^{-xA}(I-(e^{-xA}Ee^{-xA})^2)^{-1}e^{-xA}Ee^{-xA}e^{-yA}B\eqno(9.26)$$
\noindent and $\phi (x)=Ce^{-xA}B$. Then we let $\tau_0(x)=\det (I+e^{-xA}Ee^{-xA})$ and $\tau_1(x)=\tau_0(x+\pi /2)$.
\vskip.05in
\noindent {\bf Corollary 9.9} {\sl Suppose that $e^{-\pi A/2}Ee^{-\pi A/2}=-E$. Then\par
\indent (i) $T$ satisfies the Gelfand--Levitan equation
$$T(x,y)+\Phi (x+y)+{{1}\over{2}}\int_{x}^{x+\pi /2} T(x,z)\Phi (z+y)\, dz=0.\eqno(9.45)$$
\indent (ii) The entries of $T(x,x)$ satisfy} $v(x)={\hbox{trace}}V(x,v)$ {\sl and} $w(x)={\hbox{trace}}\, W(x,x)$ {\sl where}
$$v(x)={{1}\over{2}}{{d}\over{dx}}\log {{\tau_1(x)}\over{\tau_0(x)}}, \qquad w(x)={{1}\over{2}}{{d}\over{dx}}\log\tau_0(x)\tau_1(x).\eqno(9.27)$$
\indent {\sl (iii) If $\tau_0$ has all its zeros simple and $\tau_0(z_0)=0$ implies $\tau_0(z_0+\pi /2)\neq 0$, then $v(x)^2+w'(x)$ has poles of order less than or equal to one.}\par
\vskip.05in
\noindent {\bf Proof.} (i) This follows from Theorem 4.5(iv). We introduce the periodic linear system
$$\Bigl( \left[\matrix{-A&0\cr 0&-A\cr}\right], \left[\matrix{B&0\cr 0&B\cr}\right], 
\left[\matrix{0&C\cr C&0\cr}\right]; \left[\matrix{0&E\cr E&-0\cr}\right]\Bigr)\eqno(9.46)$$
\noindent and follow the proof of Proposition 2.5.\par
\indent (ii) This follows from a similar calculation to Theorem 2.9.\par
\indent (iii) This is a purely local calculation. The poles of $v$ and $w$ are all simple and occur at the zeros of $\tau_0(z)$ and the zeros of $\tau_0(z+\pi /2)$, with no cancellation between these; however, the second order poles of $v^2$ cancel all the second order poles of $w'$.\par
\rightline{$\square$}\par
\vskip.05in

\noindent {\bf Remarks 9.10} (i) Airault, McKean and Moser [3] consider 
the cases of Theorem 3.3(iv) given by $u'''=12uu'$ for $u$
rational, trigonometric and elliptic. These cases can be identified
directly in terms of the state ring via Proposition 9.2(iii), without recourse to the general
theory of finite gap potentials.\par
\vskip.1in

\indent (ii) Let $\psi$ be an elliptic solution of the second kind such that 
$\psi (x+\pi )=\mu \psi (x)$ where $\vert \mu\vert <1$ and $\psi$ has no real poles. 
Then the Hankel operator $\Gamma_\psi$ satisfies Proposition 4.9.

\indent (ii) This follows from Corollary 5.8.\par

\noindent {\bf Proposition 9.2} {\sl Let $\phi (x)=\sum_{n=-\infty}^\infty a_n e^{inx}$, where
$\sum_{n=-\infty }^\infty \vert
a_n\vert <\infty$ and $a_0=0$. then there exists a periodic linear system $(-A,B,C;E)$ such that} 
$\phi (x)={\hbox{trace}}\,
Ce^{-xA}B$.\par
\vskip.05in
\noindent {\bf Proof.} We consider $b_n=\vert a_n\vert^{1/2}$ and $c_n=a_n\vert a_n\vert^{-1/2}$ with $c_n=0$ when
$a_n=0$. We form the $\ell^2$ sequences $b_+={\hbox{column}}(b_n)_{n=1}^\infty$ 
 and $b_-={\hbox{column}}(b_n)_{n=-\infty}^{-1}$; likewise 
we introduce $c_+={\hbox{column}}(c_n)_{n=1}^\infty$ 
 and $c_-={\hbox{column}}(c_n)_{n=-\infty}^{-1}$; then we build infinite square matrices 
$$B_{\pm }=\left[\matrix{ b_\pm &0\cr}\right] ,\qquad 
 C_{\pm }=\left[\matrix{ c_\pm\cr  0\cr}\right] \eqno(9.3)$$
\noindent and introduce the
infinite diagonal matrices $A_+={\hbox{diagonal}}(-in)_{n=1}^\infty$ and 
$A_-={\hbox{diagonal}}(in)_{n=1}^\infty$. By Lemma 4.1(ii), there exist matrices $E_{\pm}$ such that $(-A_{\pm}, 
B_{\pm},
C_{\pm}; E_{\pm})$ are periodic linear systems. Then we introduce  
$$\Bigl( -\left[\matrix {A_+&0\cr 0&A_+\cr}\right] , 
\left[\matrix {B_+&0\cr 0&B_-\cr}\right], \left[\matrix {C_+&0\cr 0&C_-\cr}\right]; \left[\matrix {E_+&0\cr
0&E_-\cr}\right]\Bigr)\eqno(9.4)$$
\noindent which gives a periodic linear system such that 
$\phi (x) =c_+e^{-xA_+}b_++c_-e^{-xA_-}b_-$, which
reduces to ${\hbox{trace}}\,
Ce^{-xA}B$.\par
\rightline{$\square$}\par

 phase function, as above,
$$\nu (\lambda )=\pm \int_0^\lambda {{\Delta'(s)}\over{\sqrt{\Delta (s)^2-4}}}{{ds}\over{\pi}}+\nu_0.\eqno(10.13)$$
\noindent {\bf Proposition 8.2} {\sl For each $a\in {\bf C}\setminus
{\Lambda}$ there exists linear systems that have tau functions $\tau_0(x)$ and
$\tau_1(x)=\tau_0(x-a)$ such that
$$v(x)={{d}\over{dx}}\log {{\tau_0(x-a)}\over{\tau_0(x)}}\eqno(8.4)$$
\noindent is an elliptic function with Miura transform $u(x)=-2\wp (x)+c$ for
some constant $c$. Moreover, $v''=Q(v,v')$ for some complex rational function $Q$.}\par
\vskip.05in
\noindent {\bf Proof.} There exists a periodic linear system with $\tau$
function $\theta_1(x)$, and so for any $a\in {\bf C}$ we can introduce
linear systems with $\tau$ functions 
$\tau_0(x)=e^{\alpha x^2+\beta x}\theta_1(x)$ and $\tau_1(x)=\tau_0(x-a)$ such
that 
$$v(x)=\zeta (x)+\zeta (a-x)-\zeta (a)+{{1}\over{2}}\int_{i\omega /2}^{i\omega
+1} \wp (z)dz.\eqno(8.5)$$
By the addition rule for elliptic functions, we can write
$$2v(x)={{\wp'(x)-\wp'(a-x)}\over{\wp (x)-\wp (a-x)}}\eqno(8.6)$$
\noindent so $v''+2vv'=-2\wp'(x)$. By integrating this equation, we obtain the
stated result.\par
\indent By a classical result of function theory, there exists a nonzero complex polynomial $P$ such
that $P(v,v')=0$, hence by differentiation we obtain $v''=Q(v,v')$.\par 
\rightline{$\square$}\par
\indent Hochstadt observed that these simple zeros of $\Delta (s)^2-4$ are interlaced with zeros $\lambda'_j$ of $\Delta'(s)$, such that
$${{\Delta'(s)}\over {\sqrt{\Delta (s)^2-4}}}=c{{\prod_{j=1}^{2g}(s-\lambda_j')}\over{\sqrt{\prod_{j=0}^{2g}(s -\lambda_j^o)}}},\eqno(10.16)$$
\noindent and hence we have an Abelian integral on ${\cal E}$
$$\eqalignno{\nu (\lambda)&={{c}\over{\pi}} \int_0^\lambda {{\prod_{j=1}^{2g}(s-\lambda_j')}\over{\sqrt{\prod_{j=0}^{2g}(s -\lambda_j^o)}}}ds+\nu_0\cr
&={{c}\over{\pi}}\int^\lambda_{\lambda_0^o} \prod_{j=1}^{2g}(s-\lambda_j') {{ds}\over{\mu (s)}}+\nu_0&(10.17)\cr}$$

 We obtain a characterization of the elliptic potentials that are finite gap in terms 
of Hill's equation. All elliptic potentials can be realised as quotients of $\tau$ 
functions from periodic linear systems.\par
\indent Krichever and Novikov [39]
formulated the notion of an algebro-geometric system in terms of a family of matrix differential operators $L_j={{\partial }\over{\partial t_j}}+W_j(t)/(x-s_j)$ which commute as $s_j=t_j$.  In particular, this applies to
finite gap Schr\"odinger equations, where the spectral parameter
may be chosen to be a meromorphic function on a hyperelliptic
Riemann surface. \par

\indent In section 8 we show how to 
realise kernels of the form (1.15) from linear systems by means of 
products of Hankel operators with matricial symbols. The system of differential equations (1.14) depends upon the poles $t=(t_0, \dots ,t_N)$
of $\alpha, \beta$ and $\gamma$, and assuming that the $L_j$ satisfy equations (ii) associated with the classical concept of isomonodromy, and we prove the following.\par
\vskip.05in
\noindent {\bf Theorem 1.5} {\sl Let $W_j(t)$ and $W_\infty (t)$ be $2\times 2$ complex analytic matrices such that:\par
\indent (i) $W_0+\dots +W_N=-W_\infty$;\par
\indent (ii) the $W_j$ satisfy Schlesinger's equations as functions of $t=(t_0, \dots ,t_N)$;\par
\indent (iii) the eigenvalues of $W_\infty$ are $\pm\theta_\infty/2$, where $\theta_\infty$ is not an integer.\par
\noindent Let $\Phi (x,t)$ be a solution of 
$${{d}\over{dx}}\Phi (x,t)=\Bigl( \sum_{j=0}^N {{W_j(t)}\over{x-t_j}}+W_\infty \Bigr) \Phi (x,t)\eqno(1.23)$$
\noindent of the form $\Phi (x,y)=(I+\sum_{\ell=1}^\infty x^{-\ell }C_\ell (t))x^{-W_\infty }\Phi_0$. \par
\indent Then there exists a family of linear systems $(-A, B, C(t))$ with $H_0={\bf C}^2$ and state space $H$ such that $\Phi (x,t)=C(t)e^{-xA}B$, 
where $C(t)$ is determined by the Laurent coefficients $C_\ell (t)$, which satisfy a system of linear partial differential equations.}\par
\vskip.05in
\indent  For rational matrix linear differential equations that correspond to semiclassical weights on the real line, we use Schlesinger's equations in the style of Chen and 
Its [15] to show that $\tau (x, t_0, \dots ,t_N)$ from $(-A, B, C(t))$ coincides with the Jimbo--Ueno--Miwa tau functions [53, 24, 32].\par 

\indent (iii) Conversely, if $S(x,t;y,t)=0$, then} 
$$T(x,w)\partial_2 T(w,w)=\partial_2T(x,w)T(w,w).\eqno(7.4)$$
\indent (iii) Let 
$$f(x)=e^{\zeta x} +\int_x^\infty e^{\zeta y}T(y,x)\, dx +\int_x^\infty e^{\zeta y}\int_x^y
T(x,w)T(y,w)\, dw dy,\eqno(7.19)$$
\noindent and $u(x)=-2{{d}\over{dx}}T(x,x).$
\noindent In the case (ii) , one can show by integration by parts that $f$ satisfies Schr\"odinger's
equation. We now drop the commutativity assumption of (ii) and investigate the converse. By direct differentiation, we have
$$\eqalignno{f''(x)&=\zeta^2 e^{\zeta x} -\zeta e^{\zeta x}T(x,x)-e^{\zeta x}
{{d}\over{dx}}T(x,x)-e^{\zeta x} \partial_2T(x,x)\cr
& +\int_x^\infty e^{\zeta y}\partial_2^2
T(y,x)\, dy-\Bigl({{d}\over{dx}}T(x,x)\Bigr)
\int_x^\infty e^{\zeta y}T(y,x)\, dy
+T(x,x)^2e^{\zeta x}\cr
& -T(x,x) \int_x^\infty e^{\zeta y}\partial_2 T(y,x)\, dy -\partial_1T(x,x)\int_x^\infty e^{\zeta y} T(y,x)\, dy\cr
&+\int_x^\infty \partial_1^2 T(x,w)\int_w^\infty e^{\zeta y} T(y,w) dy
dw.\cr}$$
\noindent Now $(\partial_1^2-\partial_2^2)T(x,y)=u(x)T(x,y)$, so
$$\eqalignno{\int_x^\infty e^{\zeta y}\partial_2^2 T(y,x)\,
dy&=-\int_x^\infty
u(y)e^{\zeta y}T(y,x)\, dy -e^{\zeta x}\partial_1T(x,x)\cr
& +\zeta e^{\zeta x}T(x,x)+\zeta^2\int_x^\infty e^{\zeta y} T(y,x)\, dy&(7.20)\cr}$$    
We likewise deal with the final integral in () by repeated integration
by parts, obtaining 
$$\eqalignno{\int_x^\infty \partial_1^2 T(x,w)\int_w^\infty e^{\zeta y} T(y,w) dy
dw&=u(x)\int_x^\infty T(x,w)\int_w^\infty e^{\zeta
y}T(y,w)dydw\cr
&-\partial_2 T(x,x)\int_x^\infty e^{\zeta y}
T(y,x)dy+T(x,x)\int_x^\infty e^{\zeta y}\partial_2T(y,x) \, dy\cr
&+\int_x^\infty T(x,w)\int_w^\infty e^{\zeta y}
\partial_2^2 T(y,w)dydw;&(7.21)\cr}$$
\noindent where 
$$\eqalignno{\int_x^\infty T(x,w)\int_w^\infty e^{\zeta y}\partial_2^2
T(y,w)dydw&=-\int_x^\infty
T(x,w)\int_w^\infty u(y)e^{\zeta y}T(y,w)dydw\cr
&-\int_x^\infty T(x,w)e^{\zeta w}\partial _1T(w,w)dw+\zeta\int_x^\infty
T(x,w)e^{\zeta}T(w,w)dw\cr
&+\zeta^2\int_x^\infty T(x,w)\int_w^\infty e^{\zeta y}T(y,w)dydw.&(7.22)\cr}$$

Under the general hypotheses above, one can thus show by repeated integration
by parts that 
$$\eqalignno{f''(x)-u(x)f(x)-\zeta^2 f(x)&=
\int_x^\infty e^{\zeta y} u(y)\Bigl( T(x,y)-T(y,x)-\int_x^y
T(x,w)T(y,w)\, dw\Bigr)dy\cr
&+\int_x^\infty e^{\zeta y} \bigl( T(x,y)\partial_2T(y,y)
-\partial_2 T(x,y) T(y,y)\bigr)\,
dy.&(7.23)\cr}$$ 
\noindent If $S(x,t;y,t)=0$, then $f$ satisfies Schr\"odinger's equation and 
$$\int_x^\infty e^{\zeta w}
\bigl(T(x,w)\partial_2T(w,w)-\partial_2T(x,w)T(w,w)\bigr)dw=0;\eqno(7.24)$$
\noindent hence by Lerch's theorem on uniqueness of Laplace transforms, the integrand vanishes.\par

\indent Generally, we can assume that the
elements of ${\bf K}_2$ are meromorphic on a region $\Omega_1$. In particular, we can select an entry
$\phi\in {\bf K}_2$ called the scattering function from the fundamental solution matrix $Y$, and  we consider a differential field ${\bf
K}_2$ that contains suitable quotients of $\tau$ functions associated with the scattering function 
$\phi$. Whereas meromorphic solutions are essential in the current paper, Jimbo and Miwa [32] develop the theory of tau functions from formal power series solutions of differential equations, and  differential 
Galois theory also uses formal power series [61].\par
\indent The differential field ${\bf K}_2={\bf K}_0(u,u', u'', \dots )$ may be regarded as a counterpart of the Picard--Vessiot field.\par
\vskip.05in

\indent (1) We recall an example from [10, Example 5.5] which
relates to the logistic distribution. For ${\bf K}_0={\bf C}$ and $\phi (x)=2e^{-x}$, we can take  
$\tau (x) =1+e^{-2x}$ and 
$u(x)=-2{\hbox{sech}}^2x$ and $v(x)=-{\hbox{cosech}}\, 2x$ so that $u$
and $v$ are related by the Miura
transformation $v'-2v^2=-u/2,$ and satisfy the differential equations with polynomial right-hand
sides
$u''=3u^2+4u$ and $v''=4v+8v^3$. In this case ${\bf K}_1={\bf K}_2={\bf C}(e^{-x})$. In section 4 we consider in more detail the function fields generated by exponential functions and the associated 
Hankel operators.\par
\indent (ii) In section 6 we show that $\phi (x)=\sin x$ is associated with Weierstrass's elliptic function $u(x)=\wp (x)$ via a periodic linear system, so have ${\bf K}_0={\bf C}$, ${\bf K}_1={\bf C}(\sin x, \cos x)$ and ${\bf K}_2={\bf C}(\wp )[\wp']$. One can say that a multiple of the sine function is the scattering function associated with $2\wp$ as the potential in Lam\'e's equation, although this suppresses most of the information about the linear systems.\par
\indent (iii) For ${\bf K}_0={\bf C}(x)$ we can consider the Airy equation $\phi''(x)=x\phi (x)$, which gives a simple example of a system that is not Picard. Nevertheless, there does exist an entire solution $\phi$, so ${\bf K}_1={\bf C}(x, \phi ,\phi')$, and $u$  and satisfies Painleve equations $P_{II}$ $u''=2u^2+xu+\alpha$.\par

\indent In particular, this applies to the kernels that are used in
Tracy and Widom's theory of matrix models, as in [59]. 
One starts with a system of
differential equations 
$$J{{d}\over{dx}}\left[\matrix{f\cr g\cr}\right]=
\left[\matrix{\gamma &\alpha \cr \alpha
&\beta \cr}\right]\left[\matrix{f\cr
g\cr}\right],\qquad J=\left[\matrix{0&-1\cr 1&0\cr}\right],\eqno(1.21)$$
\noindent with $\alpha, \beta$ and $\gamma$ rational functions, then one
introduces a kernel 
$$K(x,y)={{f(x)g(y)-f(y)g(x)}\over{x-y}},\eqno(1.22)$$
\noindent which is known as an integrable operator. 
Tracy and Widom deal with the cases relating to classical differential equations over the rational or exponential field. In section 6, we extend their results to deal with differential equations over the elliptic function field. We show 
how several examples of kernels arising from integrable systems with trigonometric or rational base field ${\bf K}_0$ can be realised as degenerate cases of elliptic systems.\par 
\par  
\indent An advantage of this sufficient condition for finite gap potentials is that there is an explicit formula for the solution, whereas Drach's criterion [14] involves an unknown function. While the The proof uses results of Gesztesy and Weikard [29], so extending their terminology, we identify a particular class of linear linear systems.\par
\vskip.05in 
\noindent {\bf Definition} Let $(-A,B,C)$ be a linear system with meromorphic potential $u$ on $\Omega$, and let $S$ be a domain in ${\bf C}$. Say that a linear system $(-A_2,B_2,C_2)$ with tau function $\tau^{(2)}$ on $\Omega$ covers $(-A,B,C)$ finitely over $S$  if the general solution of $-\psi''+u\psi =\lambda \psi$ for all but finitely many $\lambda\in  S$ is given by a quotient of translated tau functions $\tau_\zeta^{(2)}(x-a)$ .\par

\noindent We prove the following.\par
\indent (i) Let $(-A,B,C)$ be a $(2,2)$-admissible linear system with bounded $A$, and input and output space ${\bf C}$. Then $(-A,B,C)$ covers itself finitely over ${\bf C}\setminus D(0, \Vert A\Vert )$.\par
\indent (ii) Let $(-A,B,C;E)$ be a periodic linear system with real $C^2$ potential $u$, and suppose that that a periodic linear system $(-A_2, B_2,C_2;E_2)$ covers $(-A,B,C)$ finitely over ${\bf C}$. Then $u$ is finite gap, and has a compact algebraic spectral curve ${\cal E}$ which has a degree two covering map of ${\bf P}^1$.\par

The term finite cover is suggested by the following example.  Lam\'e's equation 
$-y''+\ell (\ell +1)\wp y=\lambda y$ has coefficients which are meromorphic on the complex torus, while the general solution is meromorphic on the Jacobian of a hyperelliptic curve of genus $\ell =1,2,\dots $.\par

\indent  The following classical results of function theory suggest how one can classify some linear systems. Let $\phi $ be a meromorphic function that satisfies the addition rule
$$\phi (x+y)=Q(\phi (x), \phi'(x), \phi (y), \phi'(y))\qquad (x,y\in {\bf C})\eqno(1.20)$$
\noindent for some rational function $Q$. Then ${\bf K}_0={\bf C}(\phi (x), \phi' (x))$ 
is a differential field, and 
by classical function theory ${\bf K}_0$ is either the field ${\bf C}(x)$ of rational functions, 
the field ${\bf C}(e^{\lambda x})$ for some $\lambda$ with elements which are periodic, or the elliptic
function field ${\bf C}(\wp )[\wp']$, with elements that are doubly periodic.\par
\indent Suppose that $\Sigma$ has an entire $\tau$ function, so that $u$ is meromorphic on 
${\bf C}$.
 Let ${\bf K}_0$ be a subfield of the field of meromorphic functions on ${\bf C}$ such that 
${\bf A}_n={\bf K}_0[u,u', \dots ]$ is finitely generated 
over ${\bf K}_0$; let ${\bf M}$ be a maximal ideal of ${\bf A}_u$. The following gives a basic description of the Noetherian rings over meromorphic function fields that satisfy an addition rule. \par
\vskip.05in
\noindent {\bf Proposition 1.4} {\sl  (1) Suppose that ${\bf K}_0={\bf C}$. 
Then there exists an irreducible algebraic variety $V$ in ${\bf C}^n$ for some $n$ such that 
${\bf A}_u\sim {\bf C}[V].$\par
\indent (2) Suppose that ${\bf K}_0={\bf C}(t)$. Then ${\bf A}_u/{\bf M}$ 
is a finite algebraic extension of ${\bf C}(t)$, namely a compact Riemann surface.\par
\indent (3) Suppose that ${\bf K}_0={\bf C}(\wp )[\wp']$, the elliptic function field. 
Then ${\bf A}_u/{\bf M}$ is a finite algebraic extension of the elliptic function field, 
namely a Lam\'e  cover.}\par
 
\vskip.05in
\noindent There results follow from the Nullstellensatz for Noetherian rings, 
as in [Atiyah].\par 
In subsequent sections, we formulate conditions on $\Sigma$ which ensure that ${\bf A}_u$ is indeed finitely generated. \par

\noindent {\bf Definition} (Picard Systems) Let ${\bf K}_1$ be a differential field of meromorphic functions containing the field of constants ${\bf C}$ on a domain $\Omega$ in ${\bf C}$, and given a $n\times n$ matrix $W(x)$ 
with entries in ${\bf K}_0$ let the base equation be the differential equation $(d/dx)Y=W(x)Y(x)$.\par
\indent (i)  Then there exists a differential field ${\bf K}_2$ that
contains the entries of the fundamental solution matrix $Y$, called the Picard--Vessiot field.\par 
\indent (ii) If the entries of $Y$ are meromorphic on ${\bf C}$, then the system is called a Picard system [29].\par
\vskip.05in

\indent The following theorem summarizes the main results of sections 3 and 4, justifies the definitions and terminology. \par
In section 3 we conclude the proof of Theorem 1.1 by introducing the  algebra $({\bf A}, \ast
, \partial )$ and a suitable operation $\lfloor\, .\, \rfloor$.\par 
\indent A trivial solution of (1.10) is given by $F=I$, and this coincides with the case in which the
multiplication $\ast$ on ${\bf A}$ is trivial; indeed, we show that $u$ reduces to a multiple of the 
scattering  function in this case. 
Then we consider a scattering
function $\phi$ that arises from an admissible linear system, and show how the quotient and product of tau functions,
$${{\det (I+2^{-1}\Gamma_{\phi_{(x)}})}\over {\det (I-2^{-1}\Gamma_{\phi_{(x)}})}} \quad {\hbox{and}}\quad 
\det (I-4^{-1}\Gamma_{\phi_{(x)}}^2)\eqno(1.15)$$
\noindent are related by a common Gelfand--Levitan equation with matrix entries. 
This gives rise to an ordinary differential equation,
which is essentially the Darboux transformation. In section 2 we produce a Darboux transformation for
admissible linear systems, and in section 5 we do likewise for periodic linear systems.\par

\indent In [12] we introduced examples of periodic linear systems as in (ii), and here 
develop a systematic theory which shares some common elements of scattering theory from case (i). Traditionally, scattering theory has been considered in the context of Schr\"odinger's
equation $-f''+uf=\lambda f$  when $u(x)$ decays rapidly as $x\rightarrow\pm \infty$ and $\phi$ describes the asymptotic behaviour of solutions as $x\rightarrow\pm\infty$. Nevertheless, Ercolani and McKean
developed a scattering theory for the case of periodic $u$. A significant application of $R_x$ and ${\bf S}$ is to spectral theory of periodic linear systems, which seem
to lie outside the scope of [43, 51].  The following result gives a criterion for a tau function from a periodic linear system to arise from the theta function of a Jacobi variety.\par
\vskip.05in

\vskip.05in
\indent Then $\psi:V\mapsto
{\hbox{trace}}(VR)$ is a linear functional on ${\bf S}$ and 
$$ \lfloor V\rfloor= {\hbox{trace}}( (FVFA+AFVF)R),\eqno(4.15)$$
so the range of $\psi$ contains  ${\bf S}/\Theta_0$. Conversely, for all 
$V\in {\bf S}$, there exists $Y\in {\bf S}$ such that ${\hbox{trace}}(V_xR_x)=\lfloor Y\rfloor.$ To prove this, we let
$$W_x=\int_0^\infty  e^{-sA} V_xe^{-sA}\, ds,\eqno(4.16)$$
so that $AW_x+W_xA=V_x$, and then we let $Y_x=(I+R_x)W_x(I+R_x)$ so that
$$\eqalignno{{\hbox{trace}}(V_xR_x)&={\hbox{trace}}\bigl(W_x(AR_x+R_xA)\bigr)\cr
&=Ce^{-xA}W_xe^{-xA}B\cr
&= Ce^{-xA}(I+R_x)^{-1}Y_x(I+R_x)^{-1}e^{-xA}B\cr
&=\lfloor Y\rfloor.&(4.17)\cr}$$ 
\noindent The ring $({\bf A}_\Sigma, \ast, \partial )$ has an exterior product 
$$\Psi_1\wedge \Psi_2\wedge \dots \wedge \Psi_n=\sum_{\sigma\in S_n}{\hbox{sign}}(\sigma )\partial^{\sigma (1)-1}\Psi_1\ast \partial^{\sigma (2)-1}\Psi_2\ast \dots\ast \partial^{\sigma (n)-1}\Psi_n.$$
\noindent where the sum is taken over the group of permutations of $\{ 1, \dots , n\}$.\par
\indent {\sl (iii) Let 
 $\psi_j=\lfloor \Psi_j\rfloor$; then their usual Wronskian is given by}
$${\hbox{Wr}}(\psi_1, \psi_2, \dots , \psi_n) =\lfloor \Psi_1\wedge \Psi_2\wedge \dots \wedge
\Psi_n\rfloor\eqno(4.14)$$
\indent (iii) By applying the homomorphism $\lfloor \, .\, \rfloor$ to the ordered product,
 we obtain the Wronskian.\par
\indent Suppose that ${\bf L}={\bf C}(\phi_1, \dots , \phi_N)$, where each
$\phi_j$ is meromorphic on $\Omega$ and arises through some Abel
extension; suppose further that there exists $F_j\in {\bf C}(x_1, \dots ,
x_{2N})$ such that 
$$\phi_j(z+w)=F_j(\phi_1(z), \dots , \phi_N(z), \phi_1(w), \dots ,
\phi_N(w)).\eqno(4.18)$$
\noindent Then we say that ${\bf L}$ satisfies an addition rule.\par
;\par 
\indent (v) the translated potential $u_\infty (x-\alpha )$ arises from the limit
of $(M^{n/\alpha })^n$ acting on $u_\infty (x)$ as $n\rightarrow\infty.$
\indent (v) By repeated multiplication, we have
$$\bigl( M^{n/\alpha }\bigr)^n:(-A, B, C)\mapsto \Bigl( -A, (1+\alpha
A/n)^n(1-\alpha A/n)^{-1}B, C\Bigr),$$
\noindent which converges to $(-A, e^{2\alpha A}B, C)$ as $n\rightarrow\infty$.
Hence the effect is to replace $\tau_\infty (x)$ by $\tau_\infty (x-\alpha )$.\par
\rightline{$\square$}\par
\vskip.05in
\noindent By a  Pl\"ucker relation on the $2\times 2$ minors of the $2\times 4$ matrix 
$$\left[\matrix {\psi_{\zeta_1 }&\psi_{\zeta_2}&\psi_{\zeta_3}&\psi_{\zeta_4}\cr
\psi'_{\zeta_1}&\psi'_{\zeta_2 }&\psi'_{\zeta_3}&\psi'_{\zeta_4}\cr}\right]\eqno(6.38)$$
\noindent we have an identity for Wronskians 
$$W(\zeta_1,\zeta_2)W(\zeta_3,\zeta_4)=\det\left[ \matrix{ W(\zeta_1,\zeta_3)&W(\zeta_1,\zeta_4)\cr W(\zeta_2,\zeta_3
)&W(\zeta_2,\zeta_4)\cr}\right].\eqno(6.39)$$ 
\indent The infinitesimal translation operation is
$${{M^{\zeta+\varepsilon}M^{-\zeta} -I}\over{\varepsilon}}
\bigl(-A,B,C\bigr)\rightarrow 
\bigl( -A, -2A (\zeta^2I-A^2)^{-1}B, C\bigr)\qquad
(\varepsilon\rightarrow 0).\eqno(6.40)$$
\indent The functions $\psi_{\zeta}$ and $\psi_{-\zeta}$ both satisfy $\psi''=(u_\infty
+\zeta^2)\psi$, so $Wr(\psi_\zeta , \psi_{-\zeta})$ is constant, and by scaling if
necessary we can suppose that it equals one. Then we write
$\psi_{\zeta+\varepsilon}=\psi_\zeta +\varepsilon\varphi_\zeta +O(\varepsilon^2)$,
where $\varphi_\zeta$ satisfies the variational equation
$\varphi_\zeta''=2\zeta \psi_\zeta +(u_\infty +\zeta^2)\varphi_\zeta;$ hence 
$$Wr( \psi_{\zeta +\varepsilon }, \psi_{-\zeta})=1+\varepsilon 
Wr( \varphi_{\zeta}, \psi_{-\zeta})+O(\varepsilon^2)\eqno(6.41)$$
hence 
$$\eqalignno{{{d}\over{dx}}\log Wr( \psi_{\zeta +\varepsilon }, \psi_{-\zeta})
&=\varepsilon \left|\matrix{
\varphi_\zeta& \psi_{-\zeta}\cr \varphi_\zeta''&
\psi_{-\zeta}''\cr}\right|+O(\varepsilon^2)\cr
&=-2\zeta\varepsilon \psi_\zeta \psi_{-\zeta}+O(\varepsilon^2);&(6.42)\cr}$$
so 
$$-2{{d^2}\over{dx^2}}\log Wr( \psi_{\zeta +\varepsilon },
\psi_{-\zeta})=4\zeta\varepsilon {{d}\over{dx}}\psi_\zeta
\psi_{-\zeta}+O(\varepsilon^2).\eqno(6.43)$$

\vskip.05in
\noindent {\bf Corollary 6.4} {\sl Let ${\bf K}_1$ be the
complex differential
field generated by $u_\infty$, ${\bf K}_2$ be the differential field
generated by the solutions of $-\psi_\zeta''+u_\infty \psi_\zeta
=\zeta^2\psi_\zeta$; let ${\bf K}_3$ be the differential field generated
by the $\tau_\zeta$, and the $e^{\zeta x}$. Then ${\bf K}_1\subseteq 
{\bf K}_2\subseteq {\bf K}_3$.}\par 
\vskip.05in

\indent The following result is a companion to Theorems 3.1 and Theorem
6.1, except that we allow the input and output space to be $H$ instead
of ${\bf C}$, and we impose the hypothesis that $A$ commutes with
$R_x$. The result gives a generalized Darboux addition rule for
infinite-dimensional input spaces.\par
\vskip.05in   
\noindent {\bf Proposition 6.5} {\sl Consider the linear systems 
$$\Sigma_\zeta =(-A, (\zeta I+A)(\zeta I-A)^{-1}B,C)\eqno(6.44)$$
\noindent where\par
\indent (i) the operators have
$A, B,C\in {\cal L}(H)$ and $BC\in {\cal
L}^1(H)$;\par
\indent (ii) $A$ commutes with $BC$.\par
\noindent Then $\psi_\zeta =e^{\zeta
x}\tau_\zeta/\tau_\infty (x)$ satisfies}
$$\psi_\zeta ''(x)=\bigl(\zeta^2+u_\infty+q_\zeta (x)\bigr)\psi_\zeta (x).$$
\noindent {\sl where}
$$u_\infty (x)= -2\Bigl(
{{\tau_\infty '(x)}\over{\tau_\infty (x)}}\Bigr)',\eqno(6.45)$$
$$q_\zeta (x)=\bigl( {\hbox{trace}}
(T_\infty (x,x)-T_\zeta (x,x))\bigr)^2-{\hbox{trace}}((T_\infty
(x,x)-T_\zeta (x,x))^2\bigr)\eqno(6.46)$$
\noindent {\sl satisfies} $q_\zeta (x)=O(\zeta^{-2})$ {\sl as} $\zeta\rightarrow\infty$.
\vskip.05in
\noindent {\bf Proof.} By simple calculations, we find
$$\eqalignno{{{\psi_\zeta ''(x)}\over{\psi_\zeta (x)}}&=
\zeta^2-2\Bigl(
{{\tau_\infty '(x)}\over{\tau_\infty (x)}}\Bigr)' +
{{\tau_\zeta''(x)}\over{\tau_\zeta (x)}}+
{{\tau_\infty''(x)}\over{\tau_\infty (x)}}-
\Bigl( {{\tau_\zeta'(x)}\over{\tau_\zeta (x)}}\Bigr)^2- 
\Bigl( {{\tau_\zeta'(x)}\over{\tau_\zeta (x)}}\Bigr)^2\cr
&+2\zeta\Bigl( {{\tau_\zeta'(x)}\over{\tau_\zeta (x)}}-
{{\tau_\infty'(x)}\over{\tau_\infty (x)}}\Bigr)+
\Bigl( {{\tau_\zeta'(x)}\over{\tau_\zeta (x)}}-
{{\tau_\infty'(x)}\over{\tau_\infty (x)}}\Bigr)^2&(6.47)\cr}$$ 

\noindent Now we have
$$T_\infty (x,x)=-Ce^{-xA}\bigl(
I+R_x\bigr)^{-1}e^{-xA}B\eqno(6.48)$$
\noindent and
$$T_\zeta (x,x)=-Ce^{-xA}\bigl(
I+(\zeta I-A)(\zeta I-A)^{-1}R_x\bigr)^{-1}(\zeta I+A)(\zeta
I-A)^{-1}e^{-xA}B\eqno(6.49)$$
\noindent which implies by commutativity that
$$T_\infty (x,x)-T_\zeta (x,x)=2Ce^{-xA}(I+R_x)^{-1}\bigl(
I+(\zeta I-A)(\zeta I-A)^{-1}R_x\bigr)^{-1}A(\zeta
I-A)^{-1}e^{-xA}B,\eqno(6.50)$$
\noindent which is $O(\zeta^{-1})$ as $\zeta\rightarrow\infty$. Also

$$\eqalignno{{{d}\over{dx}}& T_\infty (x,x)+ {{d}\over{dx}} T_\zeta (x,x)
+(T_\infty (x,x)-T_\zeta (x,x))^2 -2\zeta (T_\infty (x,x)-T_\zeta
(x,x))\cr
&=Ce^{-xA}\Bigl\{ (I+R_x)^{-2}(I+(\zeta I+A)
(\zeta I-A)^{-1}R_x)^{-2}\cr
&\quad\Bigl( 2A\bigl(I+(\zeta I+A)(\zeta I-A)^{-1}R_x\bigr)^2 +2A 
(\zeta I+A)(\zeta
I-A)^{-1}(I+R)^2\cr
&\qquad +8A^3R_x(\zeta I-A)^{-2} -4\zeta A(\zeta I-A)^{-1}
(I+R_x)\bigl(I+(\zeta I+A)(\zeta I-A)^{-1}R\bigr)\Bigr)\Bigr\} e^{-xA}B\cr 
&=0.&(6.51)\cr}$$
We deduce that 
$$\eqalignno{{{\tau_\zeta''(x)}\over{\tau_\zeta (x)}}&+
{{\tau_\infty''(x)}\over{\tau_\infty (x)}}-
\Bigl( {{\tau_\zeta'(x)}\over{\tau_\zeta (x)}}\Bigr)^2- 
\Bigl( {{\tau_\infty'(x)}\over{\tau_\infty (x)}}\Bigr)^2+2\zeta\Bigl( {{\tau_\zeta'(x)}\over{\tau_\zeta (x)}}-
{{\tau_\infty'(x)}\over{\tau_\infty (x)}}\Bigr)\cr
&={\hbox{trace}}\Bigl( {{d}\over{dx}} T_\infty (x,x)+ {{d}\over{dx}} T_\zeta (x,x)
-2\zeta (T_\infty (x,x)-T_\zeta
(x,x)\Bigr)\cr
&=-{\hbox{trace}}\bigl( (T_\infty (x,x)-T_\zeta (x,x))^2\bigr).&(6.52)\cr}$$

\noindent This gives the stated result. In the scattering case, we have $q_\zeta (x)=0$, 
and () reduces to Schr\"odinger's equation; but we cannot generally assume this here.\par

\noindent {\bf Corollary 7.2} {\sl Suppose that \par
\indent (i) $A$ is bounded on $H$ with $\Vert e^{-t_0A}\Vert <1$ for some
$t_0>0$;\par
\indent (ii) $B: {\bf C}\rightarrow H$ and $C: H\rightarrow 
{\bf C}$ are bounded with $\Vert B\Vert \Vert C\Vert <1$;\par
\indent (iii) $A$ commutes with $BC$.\par
\indent Then $\Sigma(t)= (-A, U(t)B, CU(t))$ gives a $KdV$ system.}\par
\vskip.05in
\noindent {\bf Proof.} Evidently $A$ and $U(t)R_x$ commute under multiplication, and for some $x_0$ 
we can form $F_x=(I+U(t)R_xU(t))^{-1}$  
which also commutes with $A$ for all $x>x_0$. Then
${\cal A}={\hbox{span}}(A^jF_x^k;x>x_0, j,k=0, 1, 2, \dots \}$ forms a commutative algebra
for the usual operator multiplication and for ${}\ast{}$. So by
Proposition 7.1(ii) $S(x,t;y,t)=0$ for all $x,y>0$ which implies (7.25) by (7.2).\par
\rightline{$\square$}\par
\vskip.05in

\noindent {\bf Definition} Let ${\bf \Sigma}_{2,2}$ be the space of $(2,2)$ admissible linear systems that have input and output space $H_0$ and state space $H$. Then ${\bf \Sigma}_{2,2}$ has a natural topology determined by the following collection of basic open sets. For $j=1, \dots ,N$ with any $N\in {\bf N}$,  let $(-A_j, B_j,C_j)\in {\bf \Sigma}_{2,2}$ have tau functions $\tau_j$, let $\varepsilon_j>0$ and let $K_j$ be compact subsets of $(0, \infty )$; then the corresponding basic open set is 
$$\{ (-A,B,C)\in {\bf \Sigma}_{2,2}: \vert \tau (x)-\tau_j(x)\vert<\varepsilon,\quad  \forall x\in K_j, \quad\forall j;\quad  \tau \quad {\hbox{tau function of}}\quad (-A,B,C)\}.$$
\noindent We can extend this in the obvious way when the tau functions are holomorphic on some domain $\Omega$ of ${\bf C}$. This gives the topology of uniform convergence of tau functions.\par
\indent There is a stricter topology given by the open sets
$$\{ (-A,B,C)\in {\bf \Sigma}_{2,2}: 
\Vert R_x-R_x^{(j)}\Vert_{{\cal L}^1}<\varepsilon_j\quad  \forall x\in (0, \infty );\quad  
R^{(j)}_x \quad {\hbox{from}}\quad (-A_j,B_j,C_j); \quad\forall j\}.$$
\noindent We refer to this as the norm topology on 
${\bf \Sigma}_{2,2}$ as it coincides with the topology on the $R_x$ given by 
the supremum norm on $C_b((0, \infty ); {\cal L}^1(H))$, the space of continuous and
bounded functions $(0, \infty )\rightarrow {\cal L}^1(H)$.\par
\vskip.05in
\noindent {\bf Corollary 2.3} {\sl Suppose that $A$ has
bounded imaginary powers so $(A^{it})_{t\in {\bf R}}$ forms a strongly
continuous group of bounded linear operators on $H$ such that $\Vert
A^{it}\Vert_{{\cal L}(H)}\leq \kappa e^{\mu \vert t\vert}$ for all $t\in
{\bf R}$ and some $\mu <\pi /2$ $\kappa>0$. Then the map
$$H^\infty \times {\cal L}^2\times {\cal L}^2\rightarrow {\bf
\Sigma}_{2,2}: (r, B,C)\rightarrow (-A, r(A)B, C)$$
\noindent is continuous for the norm topology.}\par
\vskip.05in
\noindent {\bf Proof.}
Then by
McIntosh's theorem, there is a bounded linear map
$\Phi_A:H^\infty\rightarrow
{\cal L}(H)$ $f\mapsto \Phi_A(f)$ such that
$\Phi_A(fg)=\Phi_A(f)\Phi_A(g)$ and $\Phi_A(r)=r(A)$ for all stable
rational functions. So we can introduce a $(2,2)$ admissible
linear system $(-A, \Phi_A(r)B,C)$ for all $r\in H^\infty$. Then
$$(r,B,C)\mapsto r(A)\int_x^\infty e^{-tA}BCe^{-tA}\, dt$$
\noindent is continuous from $H^\infty \times {\cal L}^2
\times {\cal L}^2$ with the product norm topology to ${\cal L}^1.$\par
\par
equation. Ablowitz and Segur solved $P_{II}$ by a slightly different 
method, Borodin and Deift [12] obtained a solution by considering a
matrix Riemann--Hilbert problem involving (11.9) and we
include the proof of Theorem 11.4(iii) to illustrate the general theory of
linear systems.\par 

\vskip.1in
\noindent {\bf 10. Solutions to Painleve II and XXXIV}\par
\vskip.05in
\indent In previous sections we started from an admissible linear system
and produced a Hankel integral operator $\Gamma_\phi$. In this section
 we begin with a technical result
which realises a typical Hilbert--Schmidt Hankel operator 
$\Gamma_\phi$ from an explicit linear system $(-A,B,C)$ chosen for 
$\phi$. Here $A$ is defined on ${\cal D}(A)=
\{ f\in L^2(0, \infty ); f'\in
L^2(0, \infty )\}$ and $C$ is bounded on ${\cal D}(A)$. Suppose that $\phi$ and $\psi$ are continuous 
functions on ${\bf R}$ such that
$\int_0^\infty (1+t)(\vert \phi (t)\vert^2 +\vert\psi (t)\vert^2)\, dt<\infty$. Then we let
$H=L^2(0, \infty )$ and introduce the operators
$$\eqalignno{ A: f(x)&\mapsto -f'(x)\qquad f\in {\cal D}(A);\cr
              B:\beta&\mapsto \phi (x)\beta ;\cr
              E:\beta&\mapsto \psi (x)\beta ;\cr 
              C: g(x)&\mapsto g(0)\qquad (g\in {\cal D}(A))
,&(10.1)\cr}$$
\noindent so that $\phi (x)=Ce^{-xA}B$ and $\psi (x)=Ce^{-xA}E$. We introduce the operators on $H$ given by 
$R_x=\int_x^\infty e^{-tA}BCe^{-tA}\,dt$ and 
$S_x=\int_x^\infty e^{-tA}ECe^{-tA}\, dt$. We also introduce the
observability Gramian $Q_x=\int_x^\infty e^{-tA^\dagger}C^\dagger
Ce^{-tA}\, dt$ and we observe that $Q_x$ is 
the orthogonal projection
$Q_x:L^2(0, \infty )\rightarrow L^2(x,\infty )$. We consider the 
Gelfand--Levitan integral equation (2.8)
where $T(x,y)$ and $\Phi (x+y)$ are $2\times 2$
matrices, and 
$$\Phi (x)=\left[\matrix{0&\psi (x)\cr \phi (x)&0\cr}\right].
\eqno(10.2)$$
\vskip.05in
\noindent {\bf Lemma 10.1} {\sl (i) There exists $\delta >0$ such
that for all $\mu\in {\bf C}$ such that $\vert \mu\vert <\delta$, the operator $I-\mu^2R_xS_x$ has 
inverse} $G_x\in {\cal L}(H)$ {\sl and the matrix function
$$\hat T(x,y)=\left[\matrix {\mu Ce^{-xA}G_xS_xe^{-yA}B 
&-Ce^{-xA}G_xe^{-yA}E\cr
-Ce^{-xA-yA}B-\mu^2 Ce^{-xA}R_xG_xS_xe^{-yA}B& \mu Ce^{-xA}
R_xG_xe^{-yA}E\cr}\right]\eqno(10.3)$$
\noindent satisfies the Gelfand--Levitan equation (2.8).}\par\indent {\sl (ii) The determinants satisfy}
$$\det (I-\mu^2 R_xS_x)=\det (I-\mu^2
\Gamma_{\psi_{(x)}}\Gamma_{\phi_{(x)}}).\eqno(10.4)$$
\noindent {\sl and}
$${\hbox{trace}}\, \hat T(x,x)={{1}\over{\mu}}{{d}\over{dx}}\log\det
(I-\mu^2\Gamma_{\phi_{(x)}}\Gamma_{\psi_{(x)}}).\eqno(10.5)$$
\vskip.05in
\noindent {\bf Proof.} (i) We introduce
$$\hat A=\left[\matrix{A&0\cr 0&A\cr}\right],\quad  
\hat B=\left[\matrix{B&0\cr 0&E\cr}\right], \quad 
\hat C=\left[\matrix{0&C\cr C&0\cr}\right]\eqno(10.6)$$
\noindent and follow the computations of Proposition 2.3(i) to find
$T$.\par
\indent (ii) We observe that $R_xf(z)=\int_x^\infty \phi
(z+u) f(u)\, du$, so $R_x$ is a Hilbert--Schmidt operator, and $R_0$ is the Hankel
operator $\Gamma_\phi$; likewise $S_x$ is Hilbert--Schmidt; 
hence $R_xS_x$ is trace class. The identity (10.4) follows
from Proposition 2.3.\par
\indent Whereas $R_xS_x$ is not a differentiable function of
$x$, and we cannot adopt the direct approach of Proposition 2.3(iii), 
we can differentiate the Hankel product 
$${{d}\over{dx}}\int_0^\infty \phi (2x+s+u)\psi (2x+t+u)\, du=-2\phi
(2x+s)\psi (2x+t),\eqno(10.7)$$
\noindent so
$${{d}\over{dx}}\Gamma_{\phi_{(x)}}\Gamma_{\psi_{(x)}}=-2e^{-2xA}B
Ce^{-2xA}S_0,\eqno(10.8)$$
\noindent where the right-hand side is a rank one and bounded linear
operator. We recall from [48] the following identities regarding the shift and
Hankel operators 
$$e^{-xA^\dagger }e^{-xA}=Q_x,\qquad e^{-xA}e^{-xA^\dagger}=I,\qquad  
R_0e^{-xA^\dagger}=e^{-xA}R_0\eqno(10.9)$$
\noindent and the following special identities which may be checked by
looking at the kernels 
$$R_x=R_0Q_x,\quad  e^{-xA}R_0=R_xe^{-xA^\dagger},\quad 
\Gamma_{\phi_{(x)}}=e^{-xA}R_0e^{-xA^\dagger}.\eqno(10.10)$$
\indent Hence we can differentiate using (10.15), obtaining
$$\eqalignno{ {{d}\over{dx}}&\log\det
\bigl(I-\mu^2\Gamma_{\phi_{(x)}}\Gamma_{\psi_{(x)}}\bigr)\cr
&=2\mu^2{\hbox{trace}}\Bigl(
\bigl(I-\mu^2\Gamma_{\phi_{(x)}}\Gamma_{\psi_{(x)}}\bigr)^{-1}
e^{-2xA}BCe^{-2xA}S_0\Bigr)\cr
&=2\mu^2Ce^{-2xA}S_0\bigl(I-\mu^2e^{-xA}R_0e^{-xA^\dagger}
e^{-xA}S_0e^{-xA^\dagger}\bigr)^{-1}e^{-2xA}B\cr
&=2\mu^2Ce^{-xA}S_xe^{-xA^\dagger}\bigl(I-\mu^2e^{-xA}
R_0e^{-xA^\dagger}e^{-xA}S_0e^{-xA^\dagger}\bigr)^{-1}e^{-2xA}B.&(10.11)\cr}$$
\noindent We now use the identity $K(I+LK)^{-1}=(I+KL)^{-1}K$ to
shuffle terms around, and obtain 
$$\eqalignno{{}&=2\mu^2Ce^{-xA}S_x(I-\mu^2e^{-xA^\dagger}
e^{-xA}R_0e^{-xA^\dagger}e^{-xA}S_0e
^{-xA^\dagger})^{-1}e^{-xA^\dagger}e^{-xA}e^{-xA}B\cr
&=2\mu^2Ce^{-xA}S_x(I-\mu^2Q_xR_0Q_xS_0e^{-xA^\dagger})^{-1}Q_xe^{-xA}B\cr
&=2\mu^2 Ce^{-xA}S_x(I-\mu^2R_xS_x)^{-1}e^{-xA}B\cr
&=2\mu^2Ce^{-xA}(I-\mu^2S_xR_x)^{-1}S_xe^{-xA}B;&(10.12)\cr}$$ 
\noindent which is a multiple of the top left entry of $T(x,x),$ and likewise 
$$2\mu^2 Ce^{-xA}R_x(I-\mu^2S_xR_x)^{-1}e^{-xA}E=2\mu^2Ce^{-xA}
R_xG_xe^{-yA}E;\eqno(10.13)$$ 
as in the bottom left entry of $T(x,x)$ so we obtain the expected result
$${{d}\over{dx}}\log\det (I-\mu^2 R_x^2)=\mu
{\hbox{trace}}\, T(x,x).\eqno(10.14)$$ 
\rightline{$\square$}\par

\indent Let $(-A,B,C)$ be as in Lemma 10.3, where $T(x,y)$ and $\Phi (x+y)$ are $2\times 2$
matrices, and 
$$T(x,y)=\left[\matrix{W(x,y)&V(x,y)\cr V(x,y)&W(x,y)\cr}\right]
,\qquad \Phi (x)=\left[\matrix{0&\phi (x)\cr \phi (x)&0\cr}\right].
\eqno(10.15)$$
\noindent where the entries of $T$ are given by
$$W(x,y)=i\mu Ce^{-xA}(I+\mu^2R_x^2)^{-1}R_xe^{-yA}B\eqno(10.16)$$
\noindent and
$$V(x,y)=-Ce^{-xA}(I+\mu^2R_x^2)^{-1}e^{-yA}B.\eqno(10.17)$$
\vskip.05in
\noindent {\bf Theorem 10.4} {\sl (i) There exists $\delta >0$ such that 
for all $\mu\in {\bf C}$ such that $\vert\mu\vert <\delta$,
the operator $I+\mu^2R_x^2$ is invertible, and $T$ satisfies
$$T(x,y)+\Phi ( x+y)+i\mu\int_x^\infty T(x,z)\Phi
(z+y)dz=0.\eqno(10.18)$$
\indent {\sl (ii) Suppose further that $(-A,B,C)$ is $(2,2)$ admissible. Then the determinant satisfies}
$$V(x,x)={{1}\over{2i\mu}}{{d}\over{dx}}\log{{\det
(I+i\mu\Gamma_{\phi_{(x)}})}\over{\det I-i\mu\Gamma_{\phi_{(x)}})}},
\eqno(10.19)$$
$$W(x,x)={{1}\over{2i\mu}}{{d}\over{dx}}\log\det 
(I+\mu^2\Gamma_{\phi_{(x)}}^2),\eqno(10.20)$$
\noindent hence}
$${{1}\over{2i\mu}}{{d}\over{dx}}W(x,x)=-V(x,x)^2.\eqno(10.21)$$
\indent {\sl (iii) In particular, let} 
$\phi (x)={\hbox{Ai}}(x/2)$. {\sl Then $v(x)=V(x,x)$
satisfies Painlev\'e's equation 
$$P_{II}\qquad 
{{d^2}\over{dx^2}}v(x)=xv(x)-8\mu^2v(x)^3\eqno(10.22)$$
\noindent and} $v(x)\asymp -{\hbox{Ai}}(x)$ {\sl as} 
$x\rightarrow\infty $.\par
\indent {\sl (iv) The entries of $T(x,x)$ all lie in the differential ring
${\bf C}[x,v,v']$, and the potential is} 
$$u(x)=-8\mu^2v(x)^2.\eqno(10.23)$$
\indent {\sl (v) The cumulative distribution function of the Tracy--Widom
distribution [56] satisfies} 
$$F_2(x)=\det ((I-\Gamma_{\phi_{(x)}}^2/4).\eqno(10.24)$$ 
\vskip.05in 
\noindent {\bf Proof.}(i) This follows from Proposition 2.3.\par  
\indent (ii) The first part follows from Proposition 10.2(ii), while the
 second part is as in Theorem 2.6.\par
\indent (iii) First, note that $V(x,x)=-Ce^{-xA}(I+\mu^2R_x^2)^{-1}e^{-xA}B$ where
${\hbox{Ai}}(x/2)=Ce^{-xA}B$, so $V(x,x)$ is asymptotic to
$-{\hbox{Ai}}(x)$ as $x\rightarrow\infty$.\par
It follows from the Gelfand--Levitan equation that 
$$V(x,y)+\phi (x+y)+\mu^2\int_x^\infty\int_x^\infty V(x,z)\phi (z+s)\phi
(s+y)\, dzds=0.\eqno(10.25)$$
\noindent Let $L=({{\partial}\over{\partial x}}+
{{\partial}\over{\partial y}})^2-{{x+y}\over{2}}$ and $\phi
(x)={\hbox{Ai}}(x/2)$, so that $L\phi (x+y)=0$. Also from Airy's equation,
we obtain
$${{y-z}\over{2}}\int_x^\infty \phi (z+t)\phi (t+y)\, dt=4\bigl(
\phi'(z+x)\phi (x+y)-\phi (z+x)\phi'(x+y)\bigr),\eqno(10.26)$$
\noindent and by repeatedly integrating by parts, we can reduce (10.27) to the
expression
$$\eqalignno{LV(x,y)-&4\mu^2 \Bigl({{d}\over{dx}}\int_x^\infty V(x,z)\phi (z,x)\,
dz\Bigr)\phi (x+y)\cr
&+\mu^2\int_x^\infty\int_x^\infty  LV(x,z)\phi(z+s)\phi (s+y)\,
dzds=0.&(10.27)\cr}$$
\noindent Now we have 
$$W(x,y)+i\mu\int_x^\infty V(x,z)\phi (z+y)dz=0,\eqno(10.28)$$
\noindent so 
$$\eqalignno{{{d}\over{dx}}\int_x^\infty V(x,z)\phi (z+x)dz&={{-1}\over{i\mu}}{{d}\over{dx}}W(x,x)\cr
&=2V(x,x)^2.&(10.29)\cr}$$
On multiplying (10.38) by
$-8\mu^2V(x,x)^2$ and using uniqueness, we deduce that 
$$LV(x,y)=-8\mu^2 V(x,x)^2 V(x,y),\eqno(10.30)$$
\noindent and on the diagonal we have
$$P_{II}\qquad {{d^2}\over{dx^2}}V(x,x)-xV(x,x)=-8\mu^2 V(x,x)^3.
\eqno(10.31)$$   
\indent (iv) From (10.32), we see that ${\bf C}[x,v,v']$ is a differential ring. Note that $w'(x)=-2i\mu v(x)^2$ so the
differential equation gives
$$w(x)=2i\mu (-xv(x)^2+v'(x)^2+4\mu^2 v(x)^4),\eqno(10.33)$$
\noindent which are all elements of ${\bf C}[x,v,v']$, hence the entries
of $T(x,x)$ are all in ${\bf C}[x,v,v']$.\par 
\indent (v) With $\mu =i/2$, the potential gives rise to the Tracy Widom
distribution function
$$F_2(x)=\exp\Bigl(-2^{-1}\int_x^\infty (s-x)u(s)\,
ds\Bigr)\eqno(10.34)$$
\noindent that is associated with the soft spectral edge of the Gaussian
unitary ensemble; see [55, 56 (1.17)].\par
\rightline{$\square$}\par

\noindent {\bf Remark} The main application of the Miura transformation is in the context of the (modified) 
Korteweg--de Vries equation. The mKdV equation is associated with a pair of 
$\tau$ functions, namely $\tau_{0}$ and $\tau_{\infty}$ which satisfy (2.29) and give 
rise to a solution $v$. Then $u$ satisfies KdV, and is 
associated with $\tau_{0}$ only, so the relationship between KdV and mKdV is not 
symmetrical. When describing matrix models for quantum field theory [30], 
one says that the partition function for KdV is the square of one 
$\tau$ functions; whereas the partition function for mKdV is the product of two $\tau$ functions
 which are linked by (2.37) as in [].\par
\indent Given a solution of 
$$P_{II}:\qquad {{d^2u}\over{dz^2}}=2u^3 -4zu +4(\alpha +1/4),\eqno(10.35)$$
\noindent we introduce the new dependent variable
$$w(z)= -{{du}\over{dz}} -u^2(z) +2z\eqno(10.37)$$
\noindent which satisfies the Painlev\'e equation in Gambier's listing
$$P_{XXXIV}:\qquad 2w{{d^2w}\over{dz^2}}-\Bigl({{dw}\over{dx}}\Bigr)^2 +4w^3-8z w^2+16\alpha^2=0,\eqno(10.38)$$
\noindent as one can check by direct computation. Then one can introduce variables $(x,t)$ and scale the variable to $z=x/(3t)^{1/3}$ and write
$$W(x,t)={{1}\over{ (3t)^{2/3}}}\bigl( w(z)-2z\bigr)\eqno(10.39)$$
\noindent which satisfies the $KdV$ equation
$$4{{\partial W}\over{\partial t}}={{\partial^3W}\over{\partial x^3}} +6W{{\partial W}\over{\partial x}}.\eqno(9.40)$$

\indent The scattering functions have a simple determinantal structure. As in Corollary 2.3, let $(-A_j, B_j,C_j)$ be linear systems with scattering functions $\phi_j(x)$ for $j=1, \dots ,n)$. Then there exists a linear system with scattering function $\psi (x; s_1, \dots , s_n)=\phi_1(x+s_1)\phi_2(x+s_2)\dots \phi_n(x+s_n)$. Now we from the Wronskian of these functions 
$$\eqalignno{{\hbox{Wr}}(\phi_1(x+s_1), \dots , \phi_n(x+s_n)&=\det[ \phi_j^{(k-1)}(x+s_j)]_{j,k=1}^n \cr
&=\prod_{1\leq j<k\leq n}\Bigl( {{\partial }\over{\partial s_j}}-{{\partial }\over{\partial s_k}}\Bigr)\psi (x; s_1, \dots , s_n).\cr}$$ 
 Let $\phi (x)=Ce^{-xA}B$ and form the sequence $\sigma_0=1,$ 
$$\sigma_n(x)=\det[ \phi^{(j+k-2)}(x)]_{j.k=1, \dots ,n},$$
\noindent which is the determinant of a finite Hankel matrix, sometimes called a double Wronskian. Then by a result of Darboux [33], 
$${{d^2}\over{dx^2}} \log \sigma_n(x)= {{\sigma_{n+1}(x)\sigma_{n-1}(x)}\over{\sigma_n(x)^2}}.$$
\noindent The right-hand side is the ratio of terms in the sequence $(\sigma_n/\sigma_{n-1})$, and the formula suggests a linkage between derivatives and ratios of tau functions might hold. 
\noindent {\bf 8. The right shift semigroup, and exponential bases}\par
\vskip.05in
\indent In this section, we are concerned with the right shift semigroup on
$L^2(0, \infty )$  given by $e^{-tA}f(x)=f(x-t){\bf I}_{(0, \infty )}(x-t)$ for
$f\in L^2(0, \infty )$, as discussed by Partington and Weiss [PW]. Here
$(e^{-tA})_{t\geq 0}$ is isometric, completely non unitary and 
weakly asymptotically null, so does not satisfy Theorem 2.1(iii). In addition to its 
widespead applications, this semigroups is
fundamental to the Sz.-Nagy--Foias model theory of semigroups. We begin with a
simple sufficient condition for $W_o$ to be trace class, which applies to the standard Lyapunov
equation and to the right-shift semigroup with suitable $C$, possibly vector-valued.\par
\par

\vskip.05in
\noindent {\bf Proposition 8.1} {\sl Let $(e^{-tA})_{t\geq 0}$ be the right-shift
semigroup on $H=L^2(0, \infty )$, and suppose that $X\in {\cal L}^1(H)$, so that
$X=\sum_{j=1}^\infty B_jC_j$ where $B_j, C_j\in H$ satisfy $\sum_{j=1}^\infty 
\Vert B_j\Vert_H\Vert C_j\Vert_H<\infty$. Suppose further that\par
\indent (i) $B={\hbox{row}}(B_j)_{j=1}^\infty$ belongs to $L^2((0, \infty );
\ell^2; xdx )$;\par
\indent (ii) 
$C={\hbox{column}}(B_j)_{j=1}^\infty$ belongs to $L^2((0, \infty );
\ell^2; xdx )$;\par \par
\noindent Then $W=\int_0^\infty e^{-tA^\dagger}Xe^{-tA}\, dt$ is a trace class
linear operator on $H$, such that $W=\Gamma_B\Gamma_C$ and} 
$$A^\dagger W+WA=X.$$
\vskip.05in
\noindent {\bf Proof.} For $f,g\in H$, we have
$$C_je^{-tA}f=\int_0^\infty C_j(x)f(x-t)\, dx=\Gamma_{C_j}f(t)\eqno(8.1)$$
\noindent where $\Gamma_{C_j}$ is the Hankel operator on $H$ with scattering
function $C_j$ and
Then we build the vector Hankel operator $\Gamma_C :L^2(0, \infty
)\rightarrow L^2((0, \infty ), \ell^2)$ by 
$f\mapsto (\Gamma_{C_j}f)_{j=1}^\infty$. Now $\int_0^\infty \Vert
C(x)\Vert^2_{\ell^2}xdx$ converges  by hypothesis, so the vector Hankel operator 
$\Gamma_C: L^2(0,\infty )\rightarrow 
L^2((0,\infty ), \ell^2)$ is Hilbert--Schmidt. \par
\indent Likewise we have 
$$B_je^{-tA}g=\int_0^\infty B_j(x)g(x-t)\, dx=\Gamma_{B_j}g(t)\eqno(8.2)$$
\noindent where $\Gamma_{B_j}$ is the Hankel operator on $H$ with scattering
function $B_j$ and
Then we build the vector Hankel operator $\Gamma_B :L^2(0, \infty
)\rightarrow L^2((0, \infty ), \ell^2)$ by 
$f\mapsto (\Gamma_{B_j}f)_{j=1}^\infty$, except that we write $B$ as a row rather
than a column. Now $\int_0^\infty \Vert
B(x)\Vert^2_{\ell^2}xdx$ converges  by hypothesis, so the vector Hankel operator 
$\Gamma_B: L^2(0,\infty )\rightarrow 
L^2((0,\infty ), \ell^2)$ is Hilbert--Schmidt. Hence

$$\eqalignno{ \bigl\langle g, Wf\bigr\rangle &=\bigl\langle g, \int_0^\infty
e^{-tA^\dagger }Xe^{-tA} f\,dt\bigr\rangle_H\cr
&=\sum_{j=1}^\infty \int_0^\infty \Gamma_{B_j}g(t) \Gamma_{C_j}f(t)\, dt\cr
&\leq \Bigl(\sum_{j=1}^\infty \Vert \Gamma _{B_j}g\Vert_H^2\Bigr)^{1/2} 
\Bigl(\sum_{j=1}^\infty \Vert \Gamma_{C_j}f\Vert_H^2\Bigr)^{1/2}.&(8.3)\cr}$$
\noindent By the preceding estimates, we can use
the pairing of row and column versions of $\ell^2$
to write $W=\Gamma_B\Gamma_C$, so that  defines a trace class operator on
$H$.\par
\indent We check that $e^{-tA^\dagger} \Gamma_B\Gamma_C e^{-tA}=\Gamma_BP_{(t,
\infty )}\Gamma_C,$ and taking the derivative at $t=0+$, we obtain 
$A^\dagger W+WA=X.$
 
\rightline{$\square$}\par

\indent Let $(-A,B,C)$ be a linear system with input and output
space $H_0={\bf C}$, and state space a separable Hilbert space $H$. 
Suppose that $(e^{-tA})_{t\geq 0}$ is a strongly continuous
semigroup of contractions on $H$ so that $\Vert e^{-tA}\Vert_{{\cal
L}(H)}\leq 1$, and that there exists a constant $M$ such that 
$$\eqalignno{(\Re s)\Vert C(sI+A)^{-1}\Vert_{H}^2&\leq M \cr
(\Re s)\Vert B^\dagger (sI+A^\dagger )^{-1}\Vert^2_{H}&\leq M,\qquad (\Re s>0).&(8.4)\cr}$$
Then Jacob and Partington have shown that $W_c$ and $W_o$ are bounded
linear operators on $H$. The notion of a $(2,2)$ admissible linear
system requires that $W_c$ and $W_o$ are trace class, so we need more
stringent hypotheses on $(-A,B,C)$.\par
 The conditions (8.4) 
are related to the notion of a Carelson measure. \par   
\vskip.05in
\noindent {\bf Definition} Let $\mu$ be a positive Radon measure on ${\bf C}_+=\{ z\in {\bf C}:\Im z>0\}$. We say
that $\mu$ is a Carleson measure if there exists $\kappa>0$ such that for all bounded intervals
$I=(a,b)$ and corresponding squares $S_I=\{ x+iy: x\in (a,b), y\in (0, b-a)\}$ the measure
satisfies $\mu (S_I)\leq \kappa (b-a).$ 
\vskip.05in
\noindent {\bf Example 8.2} We recall the example of a linear systems for which the generator $(-A)$
acts as a diagonal operator on some orthonormal basis, an in [JacobPart]. Let $H=\ell^2$, and let $(\sigma_j)_{j=1}^\infty$ be a squence in ${\bf C}_+$. Then ${\cal
D}(A)=\{ (a_j)\in\ell^2: \sum_{j=1}^\infty (1+\vert \sigma_j\vert^2)\vert a_j\vert^2<\infty\}$
is a dense linear subspace of $\ell^2$, and $(e^{-tA})_{t\geq 0}$ defines a strongly continuous
semigroup on $H$, where $e^{-tA}: (a_j)\mapsto (e^{-i\bar\sigma_j t} a_j)$ with generator
$(-A): {\cal D}(A)\rightarrow H$. Let $(c_j)$ be a complex sequence such that $\sum_{j=1}^\infty
\vert c_j\vert^2/(1+\vert\sigma_j\vert^2)<\infty$; then $C:{\cal D}(A)\rightarrow {\bf C}$
defines a bounded  linear functional, where $C:(a_j)\mapsto \sum_{j=1}^\infty a_jc_j$, and
${\cal D}(A)$ is taken to have the squared graph norm $\Vert a\Vert^2_{{\cal D}(A)}=\Vert
a\Vert_H^2+\Vert Aa\Vert^2_{H}.$ We also introduce the positive measure $\mu$ on ${\bf C}_+$
 by
$\mu=\sum_{j=1}^\infty \vert c_j\vert^2\delta_{\sigma_j}.$ \par
\indent One can check that 
$$\Vert C(sI+A)^{-1}\Vert^2_H=\sum_{j=1}^\infty {{\vert c_j\vert^2}\over{\vert
s+i\bar\sigma_j\vert^2}}=\int\!\!\int_{{\bf C}_+} {{\mu (dz)}\over{\vert i\bar s+z\vert^2}};
\eqno(8.6)$$
and there exists $\kappa>0$ such that 
$$\int\!\!\int_{{\bf C}_+} {{\mu (dz)}\over{\vert i\bar s+z\vert^2}}\leq {{\kappa}\over {\Re
s}}\qquad (\Re s>0)\eqno(8.7)$$
if and only if $\mu$ is a Carleson measure. 
Partington and Jacob have shown that the controllability operator $W_o=\int_0^\infty
e^{-tA^\dagger }C^\dagger Ce^{-tA}\, dt$ is a bounded linear operator on $H$ if and only if
$\mu$ is a Carleson measure.\par
\vskip.05in
 Let $H^2$ be the usual Hardy space of holomorphic functions $h$ on $\{ z\in {\bf C}:
\Im z>0\}$ such that $\sup_{y>0} \int_{-\infty}^\infty \vert h(x+iy)\vert^2\, dx $ is 
finite. Let $\Gamma_\phi$ be a bounded Hankel integral operator on $L^2(0, \infty )$, 
with nullspace $N=\{ f\in L^2(0, \infty ): \Gamma_\phi f=0\}$. We recall how $L^2\ominus
N$ contains a natural sequence of exponential functions. Let
$T_s:L^2(0, \infty )\rightarrow L^2(0, \infty )$ be $f(x)\mapsto f(x-s)$ for $x,s\geq 0$ and 
$M_s:L^2(-\infty , \infty )\rightarrow L^2(-\infty , \infty )$ $f(x)\mapsto e^{isx} f(x)$
for $x\in (-\infty ,\infty ), s\geq 0$. Then $\Gamma_\phi T_s=T_s^\dagger \Gamma_\phi$, so
$N$ is a closed linear subspace of
$L^2(0, \infty )$ such that $T_s(N)\subseteq N$. By the Paley--Wiener theorem, the inverse
Fourier transform ${\cal F}^\dagger :L^2(0, \infty )\rightarrow H^2$ is a unitary
isomorphism, such that ${\cal F}^\dagger T_s=M_s{\cal F}^\dagger.$ It follows that ${\cal
F}^\dagger N$ is a closed linear subspace of $H^2$ such that $M_s{\cal F}^\dagger
(N)\subseteq {\cal F}^\dagger (N)$ for all $s\geq 0$, so ${\cal F}^\dagger N$ is invariant
under the semigroup of shift operators on $H^2$. By Beurling's theorem,  
${\cal F}^\dagger N$ is either zero, or infinite-dimensional; in the latter case, 
there exists an inner function $\theta$ such that 
$\theta H^2={\cal F}^\dagger (N).$ \par
\indent Now space $K_\theta =H^2\ominus \theta H^2$ has $L^2\ominus N= {\cal F}K_\theta$, 
and one can show that $K_\theta$ is a reproducing kernel Hilbert space  on $\{
\zeta \in {\bf C}: \Im \zeta >0\}$ with reproducing  kernels 
$$h^\theta_\zeta (z)={{1-\bar\theta (\zeta)\theta (z)}\over{(2\pi i)( z-\bar \zeta
)}}\qquad (\Im \zeta >0)$$
\noindent so $\langle f, h^\theta_\zeta\rangle =f(\zeta )$ for all $\zeta\in {\bf C}$ such
that $\Im \zeta >0$. The inner function $\theta$ may be written as $\theta
=u\theta_s\beta$, where $u$ is a unimodular constant, $\theta_s$ is a singular
inner function and $\beta$ is a Blaschke product with zeros $\sigma_n$. Since $\theta
(\sigma_n)=0$, one can easily show that ${\cal F}h^\theta_{\sigma_n}
=e^{-i\bar\sigma_n x}/\sqrt{2\pi}$, so $e^{-i\bar\sigma_n x}/\sqrt{2\pi}\in L^2\ominus N$. It is of interest to determine when
$(e^{-i\bar\sigma_n x})$ gives a Riesz basis. \par
\indent 
 We now extend Proposition 5.4 to the context of infinite exponential bases.
 In [], we considered exponentials $(e^{-x\lambda_j})_{j=1}^\infty$ in $L^2(0, \infty )$, 
and formulated conditions on $\phi (x)=\sum_{j=1}^\infty \xi_j e^{-\lambda_j x}$ 
to give a trace class Hankel operator. By the Paley--Wiener theorem,
 $L^2(0, \infty )$ corresponds to $H^2$  under the Fourier transform. 
In the current paper, we consider exponential bases of $L^2(0, 2\pi )$, 
or equivalently, shift invariant subspaces $\theta H^2$.\par 
\vskip.05in

\indent 
 We consider stronger conditions, which a view to having trace
class operators.\par

\vskip.05in
\noindent {\bf Definition} $(C)$ We say that a sequence $(\sigma_j)_{j=-\infty}^\infty$ in $\{z: \Im z>0\}$ 
satisfies Carleson's interpolation condition if there exists $\delta >0$ such that
$$\prod_{j=-\infty; j\neq k}^\infty \Bigl\vert{{\sigma_j-\sigma_k}\over{\sigma_j-\bar \sigma_k}}\Bigr\vert\geq 
\delta \qquad (k\in {\bf Z}).\eqno(8.8)$$

\indent In the next result we establish the properties of the following linear system, where 
$H$ is the closed linear span in $L^2(0, \infty )$ of 
$(e^{-i\bar\sigma_jx})_{j=-\infty}^\infty$ and  
$${\cal D}(A)=\Bigl\{ \sum_{j=-\infty}^\infty g_je^{-i\bar\sigma_jx}: 
\sum_j (1+\vert \sigma_j\vert^2)\vert g_j\vert^2/\Im \sigma_j<\infty \Bigr\}:\eqno(8.9)$$
$$\eqalignno{ A: {\cal D}(A)\subset H\rightarrow H: &\quad \sum_{j}g_je^{-i\bar\sigma_j x}
\mapsto \sum_j i\bar\sigma_j g_je^{-i\bar\sigma_jx} ;\cr
B: {\bf C}\rightarrow H:&\quad  b\mapsto b\sum_{j}a_je^{-i\bar\sigma_jx}\qquad (b\in {\bf C});\cr
C:{\cal D}(A)\rightarrow {\bf C}:& \quad \sum_{j}g_je^{-i\bar\sigma_jx}\mapsto \sum_{j}g_j.&(8.10)\cr}$$
\vskip.05in

\noindent {\bf Theorem 8.3} {\sl Suppose that $(C)$ holds, and the following series converge: $$\sum_j\vert a_j\vert^2/
(\Im \sigma_j),\leqno(1)$$
$$\sum_j\vert a_j\vert^2/
(\Im \sigma_j)^2.\leqno(2)$$
\indent (i) Then $(-A,B,C)$ is an linear system with scattering function
  $\phi (x)=\sum_j a_j e^{-i\bar\sigma_j x}$, where $\Gamma_\phi\in {\cal L}(L^2(0,
\infty ))$;\par
\indent (ii) $\Gamma_{\phi}$ is
a Hilbert--Schmidt operator on $L^2(0, \infty )$  with null space} 
$$N=\{ f\in L^2(0, \infty ): \hat f(\bar\sigma_j)=0; j\in {\bf Z}\}.\eqno(8.11)$$ 
\indent {\sl (iii) $R_x$ defines a Hilbert--Schmidt
operator on $H$;}\par

\indent {\sl (iv) Suppose further that $\sum_{j=-\infty}^\infty \vert a_j\vert/\Im\sigma_j$ converges. 
Then $\Gamma_{\phi_{(x)}}$ and $R_x$ are
trace class and the corresponding tau function is}
$$\det (I+\lambda \Gamma_{\phi_{(x)}})=\sum_{m=0}^\infty \lambda^m\sum_{S:S\subset {\bf Z}; \sharp S=m} \Delta_S,$$ 
\noindent {\sl where for each finite subset $S\subset {\bf Z}$ the
corresponding minor is} 
$$\Delta_S=\prod_{k\in S}{{a_ke^{-2i\bar\sigma_kx}}\over{2i\bar\sigma_k}}\prod_{j,k\in
S:j\neq k}{{\bar \sigma_j-\bar\sigma_k}\over{\bar \sigma_j+\bar\sigma_k}}.\eqno(8.12)$$
\vskip.05in
\noindent {\bf Proof.} (i) We repeat some arguments from [NikII], as need to extend them later. By 
[Garnett, p314], a complex sequence $(\sigma_j)\subset {\bf C}_+$ is
a finite union of interpolating sequences $(C)$, if and only if $\mu =\sum_j \Im \sigma_j
\delta_{\sigma_j}$ is a Carleson measure. Let $g(x)=\sum \sqrt{\Im \sigma_j }e^{-i\bar \sigma_jx} g_j$
 and let $f\in L^2(0, \infty )$. Then by the Paley--Wiener theorem, $\check f$ is a typical $H^2$ function
 and 
$$\eqalignno{\Bigl\vert\int_0^\infty \bar g(x) f(x)\, dx\Bigr\vert &=\Bigl\vert \sum_{j=-\infty}^\infty 
\bar g_j\check  f(\sigma_j ) \sqrt{\Im \sigma_j}\Bigr\vert\cr
&\leq \Bigl( \sum_{j=-\infty}^\infty \vert g_j\vert^2\Bigr)^{1/2}\Bigl(\sum_{j=-\infty}^\infty 
\Im \sigma_j \vert \check
f(\sigma_j)\vert^2\Bigr)^{1/2}\cr
&\leq \kappa \Bigl( \sum_{j=-\infty}^\infty \vert g_j\vert^2\Bigr)^{1/2}\Vert f\Vert_{L^2(0, \infty
)}&(8.13)\cr}$$
\noindent where $\kappa$ is the Carleson measure constant of $\mu$. Conversely, given $(g_j)\in \ell^2$, 
we can by the
condition $(C)$ and the Paley--Wiener Theorem, choose $f\in L^2(0, \infty )$ such that 
$\check f(\sigma_j)\Im \sigma_j
=g_j $ and $\Vert f\Vert_{L^2(0, \infty )}\leq \kappa_2 \Vert (g_j)\Vert_{\ell^2}$, where $\kappa_2$ 
depends only opon $\delta$. We deduce that there exist constants $c_1, c_2>0$ such that  
$$c_2\Bigl(\sum_{j=-\infty}^\infty \vert g_j\vert^2\Bigr)^{1/2}\leq \Bigl\Vert \sum_{j=-\infty}^\infty g_j 
e^{-i\bar\sigma_jx}\sqrt {\Im \sigma_j}\Bigr\Vert_{L^2}\leq
c_2\Bigl(\sum_{j=-\infty}^\infty \vert g_j\vert^2\Bigr)^{1/2}\eqno(8.14)$$
\noindent for all $(g_j)\in \ell^2$; hence $( e^{-i\bar\sigma_jx}\sqrt {\Im \sigma_j})$ gives a Riesz basis 
for its closed linear span in $L^2(0,
\infty )$, and 

$$H=\{ \sum_{j=-\infty}^\infty g_j e^{-i\bar\sigma_jx}\sqrt {\Im \sigma_j}: (g_j)\in
\ell^2\}.\eqno(8.15)$$ 

Hence $\sum_j a_je^{-i\bar\sigma_j x}$ belongs to $H$, so $B$ is well defined. (Note that when
$(\Im\sigma_j)_{j\in {\bf Z}}$ is bounded above, convergence of series (2) implies convergence of (1).) The semigroup $(e^{-tA})$ 
is similar to the
diagonal semigroup from the previous example, hence is similar to a strongly 
continuous semigroup of contractions on Hilbert
space. 
Note also that $e^{-tA}f(x)=f(x+t)$ for $f\in H$, where the right shift semigroup is strongly continuous 
on $L^2(0, \infty )$
by standard arguments. Hence the generator $-A$ has domain 
${\cal D}(A)$ in $H$ as stated, and ${\hbox{Spec}}(A)$ 
is the closure of $\{ i\bar\sigma_j: j\in {\bf Z}\}$. Also, the right shift semigroup is (strongly)
asymptotically stable in the sense that $\Vert e^{-tA}f\Vert_{L^2}\rightarrow 0$ as $t\rightarrow\infty$ 
for all $f\in
H$.\par
\indent To check that $C:{\cal D}(A)\rightarrow {\bf C}$ is bounded, we use the Cauchy--Schwarz inequality 
to show that  
$$ \Bigl(\sum_{j\in {\bf Z}}\vert g_j\vert\Bigr)^2
\leq\sum_ {j\in {\bf Z}}
{{\vert g_j\vert^2\vert \sigma_j+i\vert^2}\over{\Im \sigma_j }}\sum_{j\in {\bf Z}}
{{\Im \sigma_j}\over{\vert \sigma_j+i\vert^2}};\eqno(8.16)$$
\noindent and this sum is finite by the Carleson's condition $(C)$. Hence $\sum_jg_j$ converges 
absolutely for any $\sum_j g_je^{-i\bar\sigma_j x}\in {\cal D}(A).$\par
\indent By the assumptions on $(a_j)$, we have $M$ such a that $\vert a_j\vert\leq M\Im
\sigma_j$, and $\mu=\sum_j\Im \sigma_j \delta_{\sigma_j}$ is a Carleson measure by
$(C)$, so $\sum_j \vert a_j\vert \delta_{\sigma_j}$ is also Carleson measure. Then
$\phi (x)=\sum_{j} a_je^{-i\bar\sigma_jx}$ is the scattering function of a bounded
Hankel operator on $L^2(0, \infty )$. Indeed, for $f,g\in L^2(0, \infty )$, we have
$$\eqalignno{ \Bigl\vert\int_0^\infty \!\!\!\int_0^\infty \bar\phi (x+y)f(x)g(y)\,
dxdy\Bigr\vert &=\Bigl\vert\sum_j
a_j\int_0^\infty e^{i\sigma_jx}f(x)\, dx\int_0^\infty e^{i\sigma_jy} g(y)\,
dy\Bigr\vert\cr
&=\Bigl\vert\sum_j a_j \check f(\sigma_j)\check g(\sigma_j)\Bigr\vert\cr
&\leq \Bigl(\sum_j \vert a_j\vert \vert \check f(\sigma_j)\vert^2\Bigr)^{1/2}  
\Bigl(\sum_j \vert a_j\vert \vert \check g(\sigma_j)\vert^2\Bigr)^{1/2}.&(8.17)\cr}$$
which is bounded by $M\kappa \Vert f\Vert_{L^2(0, \infty )}\Vert g\Vert_{L^2(0, \infty
)}.$\par
\indent Now suppose that $\sum_j\vert a_j\vert /\Im \sigma_j \leq M_1$. Let $\Gamma_{(j)}$ be the Hankel 
operator on $L^2(0, \infty )$ with scattering function $e^{-i\bar\sigma_j x}$, so that $\Gamma_{(j)}$
 has rank one and operator norm $\Vert \Gamma_{(j)}\Vert_{{\cal L}(L^2)}\leq \kappa \vert a_j\vert/\Im\sigma_j$; 
so $\Gamma_\phi=\sum_j \Gamma_{(j)}$ is trace class with norm 
$$\Vert \Gamma_\phi\Vert_{{\cal L}^1}\leq\sum_j \Vert \Gamma_{(j)}\Vert_{{\cal L}^1}\leq \kappa\sum_j 
\vert a_j\vert/\Im\sigma_j.\eqno(8.18)$$
We deduce that the map $(a_j/\Im \sigma_j)\mapsto \phi (x)=\sum_j a_j e^{-i\bar\sigma_j x}$ gives rise to 
bounded linear operator $((a_j/\Im\sigma_j)\mapsto \Gamma_\phi$ which is bounded $\ell^\infty\rightarrow
 {\cal L}(L^2)$ and $\ell^1\rightarrow {\cal L}^1(L^2);$ 
we deduce by interpolation Theorme 2.10 of [Simon] that for $1\leq p\leq \infty$ there exists 
$\kappa_p$ such that 
$$\Vert \Gamma_{\phi }\Vert_{{\cal L}^p(L^2)}\leq \kappa_p 
\Bigl( \sum_{j} {{\vert a_j\vert^p}\over{(\Im\sigma_j)^p}}\Bigr)^{1/p}.$$ 
In particular, for $p=2$ we see that $\Gamma_\phi$ is a Hilbert--Schmidt
operator. 
For $f\in L^2(0, \infty )$, we introduce $\check f(\xi )=\int_0^\infty e^{i\xi y} f(y)dy/2\pi$, so that $\check f\in H^2$ for the upper half plane. With 
$\phi (x)=\sum_{j\in {\bf Z}}a_j e^{-i\bar\sigma_j x}$, we obtain
$$\eqalignno{\Gamma_\phi f(x)&=\int_0^\infty \phi (x+y)f(y)dy\cr
&=\sum_{j\in {\bf Z}}a_j e^{-i\bar\sigma_j x} 2\pi \check f(-\bar \sigma_j)&(8.19)\cr}$$
\noindent so that $\Gamma_\phi f=0$ if and only if $\check f(-\bar \sigma_j)=0$ 
for all $j\in {\bf Z}$; that is $\hat f(\bar\sigma_j)=0$ for all $j\in {\bf Z}$. The closed linear span of $(e^{-is\bar\sigma_j})_{j=-\infty}^\infty$ in $L^2(0, \infty )$
is invariant under the operation of $\Gamma_{\phi_{(x)}}$ for all $x>0$. Hence 
$\Gamma_{\phi_{(x)}}$ is represented by the matrix
$$\Bigl[
{{-ia_je^{-2i(\bar\sigma_j +\bar\sigma_k) x}\sqrt{\Im \sigma_k}}\over{(\bar\sigma_j+\bar\sigma_k)\sqrt{\Im\sigma_j}}}
\Bigr]_{j,k=-\infty}^\infty\eqno(8.20)$$
\noindent with respect to the Riesz basis. This matrix determines a 
Hilbert--Schmidt operator.\par
\indent (iii) The operator $R_x=\int_x^\infty e^{-tA}BCe^{-tA}\, dt$ acts on the terms in the 
Riesz  basis by 
$$R_x: \sqrt{\Im \sigma_j} e^{-i\bar\sigma_js}\mapsto 
\sum_{k=-\infty}^\infty  {{a_ke^{-i(\bar\sigma_k+\bar\sigma_j
)x-i\bar\sigma_js}\sqrt{\Im \sigma_k}\sqrt{\Im \sigma_j}}
\over{i\sqrt{\Im \sigma_k}(\bar\sigma_j+\bar\sigma_k)}}\qquad (j\in {\bf Z}),\eqno(8.21)$$ 
\noindent so $R_x$ is represented by the infinite matrix 
$$R_x=\Bigl[ {{a_ke^{-i(\bar\sigma_k+\bar\sigma_j
)x}\sqrt{\Im
\sigma_j}}\over{i\sqrt{\Im \sigma_k}(\bar\sigma_j+\bar\sigma_k)}}\Bigr]_{j,k=-\infty}^\infty 
.\eqno(8.22)$$ 
\noindent with respect to the Riesz basis of $H$. Let $U={\hbox{diagonal}}(e^{i\bar\sigma_k x})$,
and note that $U^{-1}R_x U$ has the same spectrum as $R_x$. The matrix of  $U^{-1}R_x U$ is equal to the matrix that represents
$\Gamma_{\phi_{(x)}}$, so both $R_x$ and  $\Gamma_{\phi_{(x)}}$ are Hilbert--Schmidt.
By Lemma 2.1 the integral
$\int_x^\infty e^{-tA}BCe^{-tA}dt$ is weakly convergent on $H$, so $R_x$ unambiguously defines
a Hilbert--Schmidt operator.\par 
\par 

\indent (iv) Under the stronger hypothesis on $(a_j)$, we have seen that $\Gamma_{\phi_{(x)}}$ 
is trace class. Since $\Gamma_{\phi_{(x)}}$ and $R_x$ are both represented
by the same matrix with respect to some Riesz basis, $R_x$ is also trace class, and their Fredholm
determinants are equal, so
$$\tau (x;\lambda )=\det (I+\lambda\Gamma_{\phi_{(x)}})=\det (I+\lambda R_x).\eqno(7.28)$$ 
\noindent To compute this determinant, we use the matrix ()and the Cauchy--Binet formula.\par
\rightline{$\square$}\par
\vskip.05in

corresponds to $H^2\ominus uH^2$. We can interpret $uH^2$ as the set of
$f\in H^2$ such that $f(i\bar\lambda_j)=0$ for all $j\in {\bf Z}$.\par
\rightline{$\square$}

\noindent {\bf Remark} (i) Let $\phi (x)={\hbox{sech}}^kx$ for $k=1, 2, \dots,$ as in the theory of solitons. Then $\Gamma_\phi$ is represented with respect to $(e^{-(k+2\ell )x})_{\ell =0}^\infty$ by the matrix
$$\Bigl[ {{2^{k-1}(-1)^j(k+j-1)!}\over{j!(k-1)! (j+k+\ell )}}
\Bigr]_{j,\ell =0}^\infty .\eqno(8.37)$$
\indent (ii) In Corollary 5.6, the semigroup $e^{-tA}$ is of exponential decay as $t\rightarrow\infty$. the next section, we consider the case in which $e^{-tA}$ is a periodic group. In section 7, we obtain examples of functions such as $\psi$ from Hill's equation via periodic groups.\par 

\vskip.05in

\indent If we aim to let $n\rightarrow\infty$ in Example 7.5, then we soon encounter technical issue.
For instance, for typical real sequences $(\lambda_j)_{j=1}^\infty$ the algebra 
${\bf C}(e^{-t\lambda_1}, e^{-t\lambda_2}, \dots )$ has infinite transcendence degree over ${\bf C}$
and is not Noetheian. Using potentail theory, one can make 
some staements about the tau function from Proposition 7.4.\par

\noindent {\bf Corollary 7.6} {\sl 
Let $\lambda_j>0$, and suppose that $x/n\rightarrow \kappa >0$ as $n\rightarrow\infty$ for some
$\kappa>0$. Then}
$$\lim\sup_{n\rightarrow\infty}{{1}\over{n^2}}\log\Bigl(\prod_{1\leq
j<k\leq n}{{(\lambda_j-\lambda_k)^2}\over{ 
(\lambda_j+\lambda_k)^2}}{e^{-2\sum_{j=1}^n\lambda_jx}}\Bigr)\leq {{-\pi^2}\over{4}},\eqno(7.13)$$
\noindent {\sl where the constant is best possible.}\par
\vskip.05in
\noindent {\bf Proof.} Let $V$ be a positive and continuous real function such that
$V(\lambda)/\lambda^\delta \geq c_0$ for all $\lambda\geq c_1$ for some $\delta, c_0, c_1>0$. The Greens function of
$\{ z: \Re z>0\}$ is 
$$G(z,w)=\log\Bigl|{{z+\bar w}\over{z-w}}\Bigr|\qquad (z\neq w, -\bar w).\eqno(7.14)$$
For compactly supported probability measures on $(0, \infty )$, we introduce the
energy of $\rho$ associated with the external field $V$ by  
$$E_V(\rho )=\int_0^\infty V(\lambda) \rho (d\lambda)+\int\!\!\!\int_{[\lambda\neq
y]} G(\lambda,\mu)
\rho (d\lambda)\rho (d\mu);\eqno(7.15)$$  
\noindent Then we let the minimum energy of the field with external potential $V$ be
$$E_V=\inf_\rho\{  E_V(\rho )\}.\eqno(7.16)$$
\noindent There exists a compactly supported $\rho_0$ that attains the minimum, 
and satisfies 
$$V(\lambda)=-2\int_S G(\lambda,\mu)\rho_0 (d\mu)+c\eqno(7.17)$$ 
\noindent for all $\lambda$ is the support $S$ of $\rho_0$, and some $c$. 
 For $0<\lambda_1\leq \dots\leq \lambda_n$, we introduce the empirical
distribution $\rho_n(d\lambda )=n^{-1}\sum_{j=1}^n\delta_{\lambda_j}(d\lambda ).$ 
Then we observe that 
$$\eqalignno{\prod_{j=1}^n e^{-nV(\lambda_j)}
\prod_{1\leq j<k\leq n}{{(\lambda_j-\lambda_k)^2}\over{ 
(\lambda_j+\lambda_k)^2}}&=\exp\Bigl[-n^2\Bigl(\int V(\lambda )\rho_n(d\lambda )+\int\int
G(\lambda,\mu) \rho_n(d\lambda)\rho_n(d\mu)\Bigr)\Bigr]\cr
&=\exp\Bigl(  -n^2E_V(\rho_n)\Bigr).&(7.18)\cr}$$  
In particular, let
$V(\lambda )=2x\lambda /n$ and $b=\pi n/x$. Then 
$$\rho (\lambda )={{\lambda}\over{b\sqrt{b^2-\lambda^2}}}\qquad (0<\lambda <b)\eqno(7.19)$$
\noindent satisfies 
$$\int_0^{b} \rho (\lambda )d\lambda =1,$$
\noindent and
$$\eqalignno{ {{d}\over{d\lambda }}\int_0^{b} 2\log \Bigl|{{\lambda -\mu}\over{\lambda
+\mu}}\Bigr|\rho (\mu )d\mu&={\hbox{PV}}\,4\int_0^{b}{{\mu \rho (\mu
)d\mu}\over{\lambda^2-\mu^2}}\cr
&={\hbox{PV}}\,{{2}\over{b}}\int_0^{b^2} \sqrt{ {{\xi}\over{b^2-\xi}}}
{{d\xi}\over{\lambda^2-\xi}}\cr
&={{2\pi}\over{b}}={{2x}\over{n}}\cr
&=V'(\lambda ).\cr}$$ 
\noindent Hence $\rho_0(d\lambda )=\rho(\lambda )d\lambda$ gives the minimizer. By considering the identity
$$V(\lambda )=c-2\int_0^b G(\lambda , \mu )\rho (\mu )d\mu \qquad
(0<\lambda <b)$$
\noindent as $\lambda\rightarrow 0+$, we deduce that $c=0$. Hence
$$\eqalignno{E_V(\rho )&
= {{1}\over{2}}\int_0^b V(\lambda )\rho (\lambda )d\lambda\cr
&={{1}\over{2}}
\int_0^b {{2x\lambda^2d\lambda}\over{nb\sqrt{b^2-\lambda^2}}}\cr 
&={{\pi^2}\over{4}}.&(7.20)\cr}$$

\indent  Kronecker's theorem asserts that a bounded Hankel integral
operator has finite rank if and only if the transfer function $\hat
\phi$ is a rational function with all its poles in $\{ z\in {\bf C}:
\Re z<0\}$. Such rational functions are known as stable. In [23], the authors consider factorization of the transfer
function in $M_{n\times n}({\bf C}(\lambda ))$ and the subring of
stable matrix rational functions. Their results describe the
properties of $\hat {\bf S}$ rather than ${\bf S}$ itself.\par 
\vskip.05in

\indent (ii) Suppose that 
$e^{-tA}$ commutes with $R_y$; then }
$$S(x,t;y,t)=0.\eqno(7.3)$$
\indent (ii) Under the commutativity assumption, the previous
integral from (7.16) equals
$$\eqalignno{Ce^{-xA} (I+R_x)^{-1}& \int_y^x
e^{-zA}BCe^{-zA}(I+R_y)^{-1}e^{-yA}B \,dz\cr
&=
Ce^{-xA} (I+R_x)^{-1} \int_y^x R_z'\, 
dz (I+R_y)^{-1}e^{-yA}B\cr
&= Ce^{-xA} (I+R_x)^{-1} (R_x-R_y)(I+R_y)^{-1}e^{-yA}B\cr
&=-Ce^{-yA}(I+R_y)^{-1}e^{-xA}B+
Ce^{-xA}(I+R_x)^{-1}e^{-yA}B\cr
&=T(z,x)-T(x,y);&(7.17)\cr}$$
\noindent so 
$$S(x,t;y,t)={{1}\over{2\pi i}}\int_{\vert \zeta\vert =r}e^{\zeta
(x-y)} (1+V_\zeta (x))(1+V_{-\zeta }(y))\, d\zeta =0.\eqno(7.18)$$
\noindent We can express this condition in terms of derivatives by writing
$$\eqalignno{T(&x,w)\partial_2 T(w,w)-\partial_2T(x,w)T(w,w)\cr
&= Ce^{-xA}(I+R_x)^{-1}\bigl( 
e^{-wA}BCe^{-wA}(I+R_w)^{-1} (-A) +Ae^{-wA}BCe^{-wA}(I+R_w)^{-1}\bigr) e^{-wA}B\cr
&= -Ce^{-xA}(I+R_x)^{-1} \bigl( R_w'(I+R_w)^{-1}A-AR_w'(I+R_w)^{-1}\bigr) e^{-wA}B\cr
&=0.\cr}$$

\rightline{$\square$}\par  

\noindent {\bf Example 4.6} The following example appears implicitly in the theory of solitons. Let $[a_{jk}]_{j,k=1}^n$ satisfy $a_{jk}=a_{kj}$
with $a_{jj}>0$ and $a_{jk}\leq 0$ for all $j\neq k$, and $\sum_{k=1}^n
a_{jk}=0$ for all $j$. Then let $A=[a_{jk}]+\varepsilon I_n$ for
$\varepsilon >0$, $B={\hbox{column}}\,[1]$, and $C={\hbox{row}}\,[1]$. We obtain a
semigroup $(e^{-tA})$ of strict contractions on ${\bf C}^n$ such that
$e^{-tA}BC=e^{-t\varepsilon }BC$ and $BCe^{-tA}=BCe^{-t\varepsilon}$.\par 
\indent Let $\lambda_j$ $(j=1, \dots, n)$ be the eigenvalues of $A$,
let ${\bf K}_0={\bf C}(e^{-\lambda_1 x}, \dots , e^{-\lambda_nt})$, 
and introduce an orthogonal matrix $S$ such that $D=S^{-1}AS$ is
diagonal. Then we work with $(-D, S^{-1}B, CS)$ instead of $(-A,B,C)$,
where both linear systems have the same tau function. In particular, 
let $R_0=\int_0^\infty e^{-sD}S^{-1}BCSe^{-sD}$, so that $R_0, D$ and
$S^{-1}BCS$ all commute. Hence ${\cal A}={\bf K}_0[I, D, R_0, S^{-1}BCS]$ is a
commutative algebra which contains $e^{-xD}$ and hence
$R_x=e^{-xD}R_0e^{-xD}$. By Hilbert's basis theorem, ${\cal A}$
 is a Noetherian ring. Let ${\cal M}$ be a maximal ideal in ${\cal
A}$. Then the field ${\cal A}/{\cal M}$ is a finite algebraic
extension of ${\bf K}_0$ by weak Nullstellensatz. For instance, ${\bf
K}_0[I,D]$ is such a field.\par
\vskip.05in

We can write
$$\eqalignno{\varphi_k(x)&=\cos kx\Bigl( 1-8C\cosh (xA)((I+S_x)^{-1}(k^2I+A^2)^{-1}\sinh (xA)B\Bigr)\cr
&\quad -8k\sin kx\Bigl( C\cosh (xA)(I+S_x)^{-1}(k^2I+A^2)^{-1}\cosh (xA)B\Bigr),&(9.40)\cr}$$
\noindent where the first factor involves 

 $$\eqalignno{1-8{\hbox{trace}}\bigl(C\cosh (xA)&((I+S_x)^{-1}(k^2I+A^2)^{-1}A\sinh (xA)B\bigr)\cr
&= 1-8{\hbox{trace}}\bigl(k^2I+A^2)^{-1}A\sinh (xA)BC\cosh (xA)((I+S_x)^{-1}\bigr)\cr
&=\det \Bigl( I-(k^2I+A^2)^{-1}A(V_xA+AS_x)(I+S_x)^{-1}\Bigr)\cr
&={{\det\bigl( I+S_x-(k^2I+A^2)^{-1}(A^2S_x+AV_xA)\bigr)}\over{\det (I+S_x)}}\cr
&={{\det\bigl( I+(k^2I+A^2)^{-1}(k^2S_x-AV_xA)\bigr)}\over{\det (I+S_x)}}\cr
&={{\det\bigl( k^2I+A^2+ k^2S_x-AV_xA\bigr)}\over{\det (k^2I+A^2)\det (I+S_x)}},&(9.41)\cr}$$
\noindent where the final line is a notationally suggestive rewriting of the preceding line.\par
\vskip.05in
\indent For more general real bounded potentials,  $u\in C_b^2({\bf R}; {\bf R})$, we can seek to produce scattering functions via the matrix spectral measure. For $\lambda\in {\bf C}$ with $\Im \lambda >0$, we solve the system 
$${{d}\over{dx}}\left[\matrix{ h_+&h_-\cr h_+'&h_-'\cr}\right] = \left[\matrix{ 0&1\cr u+\lambda &0\cr}\right] \left[\matrix{ h_+&h_-\cr h_+'&h_-'\cr}\right]\eqno(2.35)$$
\noindent so that $h_+\in L^2(0, \infty )$, $h_-\in L^2(-\infty ,0)$ and $h_+h_-'-h_+'h_-=1$; the existence follows from a theorem of Weyl; see Theorem 10.1.4 of [27].\par
\vskip0.05in
\noindent {\bf Definition} (i) The Greens function 
$$G(x,y;\lambda )=\cases {h_+(x)h_-(y) & for $y<x$;\cr 
h_+(y)h_-(x)& for $y>x$;}\eqno(2.36)$$
\noindent gives an integral operator on $L^2({\bf R}; {\bf C})$ such that $g(x)=\int G(x,y;\lambda )f(y)\, dy$ satisfies $\lambda g-g''+ug=f$.  See [27, section 10.7].\par
\indent (ii) The spectral matrix is
$$M(\lambda ; x_0)=  \left[\matrix{ 2h_-h_+&h_+'h_-+h_+h_-'\cr h_+h_-'+h_+h_-'&2h_+'h_-'\cr}\right]\qquad (\Im \lambda >0),\eqno(2.37)$$
\noindent with all entries evaluated at $x_0\in {\bf R}$. Note that the top left entry determines $G(x_0, x_0; \lambda )$. \par
\indent (iii) The function $\lambda\mapsto -M(\lambda ;x_0)$ is of positive type on the upper half plane $\{ \lambda: \Im \lambda >0\}$ , so $(-M(\lambda ; x_0)+M(\lambda ; x_0)^\dagger )/(2i)$ is positive definite for $\Im \lambda >0$, and $M(\lambda ;x_0)$ extends to a holomorphic matrix function on ${\bf C}\setminus (-\infty ,x_1]$ for some $x_\in {\bf R}$; also, there exists a $2\times 2$ positive  matrix-valued measure $\omega (ds;x_0)$ on 
$[-x_1, \infty )$, called the spectral weight, such that 
$$M(\lambda; x_0)=\Re M(-i; x_0)+\int_{-x_1}^\infty \Bigl({{1}\over{s+\lambda}}- {{s}\over{1+s^2}}\Bigr)\omega (ds; x_0)\qquad (\Im\lambda >0).\eqno(2.38)$$
\vskip.05in
\indent Given a positive matrix measure $\omega$, we can introduce $M$ via (2.38),  and  write
$$  M(\lambda ; x_0)=\left[\matrix{\gamma (\lambda )& \alpha (\lambda )\cr \alpha^*(\lambda )&\beta (\lambda )\cr}\right]\eqno(2.39)$$
\noindent where $\alpha (\lambda )^*=\overline{\alpha (\bar\lambda )}$; then we consider the differential equation
 
$${{d}\over{dx}}\left[\matrix{ f\cr g\cr}\right]=\left[\matrix{\alpha (x)&\beta (x)\cr -\gamma (x)&-\alpha (x)\cr}\right] \left[\matrix{ f\cr g\cr}\right].\eqno(2.40)$$
\noindent This has the style of the fundamental differential equation from [62] except that Tracy and Widom restrict to rational functions as coefficients. The basic kernel is 
$$K(x,y)={{f(x)g(y)-f(y)g(x)}\over{x-y}},\eqno(2.41)$$
\noindent  and it is of interest to compute $\det (I-K)$. In [9, Theorem 1.1] we showed that if $f,g\in L^2((0, \infty ; {\bf R})$ are continuous and $f(x), g(x)\rightarrow 0$ as $x\rightarrow\infty$, then 
 there exists a linear system with input and output space $H_0$ and scattering function $\phi \in L^2((0, \infty ); H_0)$ such that $K=\Gamma_\phi^\dagger\Gamma_\phi. $ In the case in which $\phi \in L^2((0, \infty ); {\bf R})$, we can apply Theorem 2.1. The paper [10] analyses cases of (2.38) which produce determinantal random point fields. In the next section, we show how to compute $\det (I-\Gamma_\phi^2);$  in Proposition 5.5, we obtain a formula for the moments of the top left entry of $\omega (ds; x_0)$.\par

indent {\sl (iv) For any stable rational function $r(z)$, the linear system $(-A, r(A)B,
C)$ is also $(2,2)$ admissible, and $r(A)R_x$ satisfies Lyapunov's equation.}\par
\indent (iv) For any $X\in {\cal L}(H)$ such that $XA=AX$, the system $(-A,XB,C)$ is
$(2,2)$ admissible and it is clear from commutativity that $XR_x$ gives a solution of
Lyapunov's equation with initial condition $-XBC$.  Note that for 
$\zeta\in {\bf C}$ with $\Re \zeta >0$, the operator
$(\zeta I+A)^{-1}$ is bounded; likewise, for any stable rational function $r(z)$
the partial fraction decomposition shows that $r(A)\in {\cal L}(H)$ by addition and
multipliction. Further, we note that when $\phi$ is the scattering
function for $(-A,B,C)$, the function 
$$\int_t^\infty e^{-\zeta u}\phi (u)\, du$$
\noindent is the scattering function for $(-A, (\zeta I+A)^{-1}B, C)$.\par
We introduce
$$V(x)=\int_x^\infty T(x,y)e^{ik(y-x)}\, dy,\eqno(6.14)$$
\noindent so that
$$V'(x)=-T(x,x)+\int_x^\infty {{\partial T}\over{\partial
x}}(x,y)e^{ik(y-x)}\, dy -ik\int_x^\infty T(x,y)e^{ik(y-x)}\, dy\eqno(6.15)$$
\noindent hence
$$\eqalignno{V''(x)&=-{{d}\over{dx}}T(x,x)
-{{\partial T}\over{\partial x}}(x,x)+
\int_x^\infty {{\partial^2T}\over{\partial x^2}}(x,y)e^{ik(y-x)}\, dy\cr
&\quad -2ik\int_x^\infty {{\partial}\over{\partial x}}T(x,y)e^{ik(y-x)}\, dy+ikT(x,x)-k^2\int_x^\infty
T(x,y)e^{-ik(y-x)}\, dy;&(6.16)\cr}$$
\noindent here by the Gelfand--Levitan equation
$$\eqalignno{\int_x^\infty {{\partial^2T}\over{\partial x^2}}
(x,y)e^{ik(y-x)}\,
dy&={{-\partial T(x,x)}\over{\partial y}}+ikT(x,x)-k^2\int_x^\infty
 T(x,y)e^{ik(y-x)}\, dy\cr
&-2{{d}\over{dx}}T(x,x)\int_x^\infty T(x,y)e^{ik(y-x)}\, dy&(6.17)\cr}$$
\noindent hence 
$$V''(x)=-2{{d}\over{dx}}T(x,x)\bigl( I+V(x)\bigr)-2ikV'(x).\eqno(6.18)$$
\noindent We deduce that
$$\eqalignno{\psi_{ik}(x)&=e^{ikx}{{\det
(I+{{ikI+A}\over{ikI-A}}R_x)}\over{\det (I+R_x)}}\cr
&=e^{ikx}\det (I+V(x))&(6.19)\cr}$$
\noindent satisfies
$${{\psi'_{ik}(x)}\over{\psi_{ik}(x)}}=ik+{\hbox{trace}}\Bigl(
(I+V(x))^{-1}V'(x)\Bigr),\eqno(6.20)$$
\noindent so 
$$\eqalignno{{{\psi''_{ik}(x)}\over{\psi_{ik}(x)}}-
\Bigl({{\psi'_{ik}(x)}\over{\psi_{ik}(x)}}\Bigr)^2
&={\hbox{trace}}\Bigl(
(I+V(x))^{-1}V''(x)\Bigr)\cr
&\qquad -{\hbox{trace}}\Bigl(
(I+V(x))^{-1}V'(x)(I+V(x))^{-1}V'(x)\Bigr)\cr
&=-2{\hbox{trace}}{{d}\over{dx}}T(x,x)-2ik{\hbox{trace}}\Bigl(
(I+V(x))^{-1}V'(x)\Bigr)\cr
&\qquad -{\hbox{trace}}\Bigl(
(I+V(x))^{-1}V'(x)(I+V(x))^{-1}V'(x)\Bigr)&(6.21)\cr}$$
so 
$$\eqalignno{ {{\psi''_{ik}(x)}\over{\psi_{ik}(x)}}&
=-k^2-2{\hbox{trace}}{{d}\over{dx}}T(x,x)+\Bigl({\hbox{trace}}\Bigl(
(I+V(x))^{-1}V'(x)\Bigr)\Bigr)^2\cr
&-{\hbox{trace}}\Bigl(
(I+V(x))^{-1}V'(x)(I+V(x))^{-1}V'(x)\Bigr).&(6.22)\cr}$$
\noindent When $v'(x)=0$, the final terms cancel, so $\psi_{ik}(x)$
gives a solution of Schr\"odinger's equation.\par

\noindent {\bf Definition} (Sato's bilinear integral) Suppose that there exists $r_0$ such that the function  $\zeta\mapsto \psi_\zeta (x,t)\psi_{-\zeta }(y,s)$ is holomorphic on $\{ \zeta\in {\bf C}: 
\vert\zeta\vert>r_0\}$, and define Sato's bilinear integral by
$$S(x,t;y,s)={{1}\over{2\pi i}}\int_{\vert \zeta\vert =r}
\psi_\zeta (x,t)\psi_{-\zeta }(y,s)\, d\zeta \qquad (r>r_0)\eqno(7.11)$$
\noindent  If $S(x,t;y,t)=0$ for all $(x,t), (y,t)\in ({\bf T}^\infty )_0$, then 
we say that the family of tau functions $\tau_\zeta (x,t)$ satisfies
the KdV hierarchy.\par
\vskip.05in
\noindent {\bf Proposition 7.2} {\sl (i) Then $S(x,t;y;t)$ satisfies the hyperbolic differential
equation $$(\partial_x^2-\partial_y^2)S(x,t;y,t)
=(u_\infty (x)-u_\infty (y))S(x,t;y,t)\eqno(7.12)$$
\indent (ii)  and} 
$$S(x,t;y,t)=T(x,y)-T(y,x)-\int_x^y T(x,z)T(y,z)\, dz\qquad
(x,y>0; t\in {\bf T}_0^\infty).\eqno(7.13)$$

\vskip.05in
\noindent {\bf Proof.} (i) By Schr\"odinger's equation for $(-A, U(t)B, C)$, we have
$$\Bigl(-{{\partial^2}\over{\partial x^2}}+u_\infty (x;t)+{{\partial^2}\over{\partial
y^2}}-u_\infty (y;t)\Bigr) S(x,y;t,t)=0.\eqno(7.14)$$
\indent (ii) Now we consider
$$S(x,t;x,t)={{1}\over{2\pi i}}\int_{\vert \zeta\vert =r} {{\tau_\zeta
(x;t)\tau_{-\zeta}(x;t)}\over {\tau_\infty (x;t)^2}}d\zeta .\eqno(7.15)$$
\noindent With $C_j=CU(t)e^{-xA}(I+R_x)^{-1}A^je^{-xA}U(t)B$, we have
$$\eqalignno{ {{\tau_\zeta (x;t)}\over{\tau_\infty
(x;t)}}&=\det \bigl(I+\sum_{j=0}^\infty \zeta^{-j-1} C_j\bigr)&(7.16)\cr
&=1+\zeta^{-1}{\hbox{trace}}(C_0)+\zeta^{-2}\bigl(2^{-1}
{\hbox{trace}}(C_0)^2-2^{-1}{\hbox{trace}}(C_0^2)+{\hbox{trace}}(C_1)
\bigr)+O(\zeta^{-3}).\cr}$$
Hence
$$\eqalignno{S(x,t;x,t)&={{1}\over{2\pi i}}\int_{\vert \zeta\vert
=r}\bigl(1+\zeta^{-1}C_0+\zeta^{-2}C_1+O(\zeta^{-3})\bigr)
\bigl(1-\zeta^{-1}C_0+\zeta^{-2}C_1+O(\zeta^{-3})\bigr)d\zeta\cr
&=0&(7.17)\cr}$$
\noindent by Cauchy's integral formula, since the residue at infinity
 vanishes. Now we write
$$S(x,t;y,t)={{1}\over{2\pi i}}\int_{\vert \zeta\vert =r}e^{\zeta
(x-y)} (1+V_\zeta (x))(1+V_{-\zeta }(y))\, d\zeta\eqno(7.18)$$
\noindent By functional calculus, for $r>\Vert A\Vert$, we have
$$\eqalignno{{{1}\over{2\pi i}}\int_{\vert \zeta\vert =r} e^{\zeta
(x-y)}V_\zeta (x)\, d\zeta &= \int_{\vert \zeta\vert =r} e^{\zeta (x-y)}Ce^{-xA}(I+R_x)^{-1}(\zeta
I-A)^{-1} e^{-xA}B{{d\zeta}\over{2\pi i}}\cr
&=Ce^{-xA}(I+R_x)^{-1}e^{-yA}B\cr
&=-T(x,y);&(7.19)\cr}$$
\noindent and likewise
$$\eqalignno{{{1}\over{2\pi i}}\int_{\vert \zeta\vert =r} e^{\zeta
(x-y)}V_{-\zeta} (y)\, d\zeta &=
\int_{\vert \zeta\vert =r} e^{\zeta (x-y)}Ce^{-yA}(I+R_y)^{-1}(\zeta
I+A)^{-1} e^{-yA}B{{d\zeta}\over{2\pi i}}\cr
&=
Ce^{-yA}(I+R_y)^{-1}e^{-xA}B\cr
&=-T(y,x).&(7.20)\cr}$$
The other contribution to $S(x,t,y,t)$ is given by 
$$-\int_{\vert \zeta\vert =r} e^{\zeta (x-y)} V_\zeta (x)V_{-\zeta }(y)
{{d\zeta}\over{2\pi i}}\eqno(7.21)$$
$$=Ce^{-xA} (I+R_x)^{-1} \int_{\vert \zeta\vert =r}
e^{\zeta (x-y)}e^{-xA}(\zeta
I-A)^{-1}BCe^{-yA}(I+R_y)^{-1}(\zeta I+A)^{-1} e^{-yA}B
{{d\zeta}\over{2\pi i}}.\eqno(7.22)$$
\noindent We focus attention on the first two occurrences of $\zeta$,
and consider the family of integrals 
$$\eqalignno{{{d}\over{dz}}\int_{\vert \zeta\vert =r}&
e^{\zeta (z-y)}e^{-zA}(\zeta
I-A)^{-1}BCe^{-yA}(I+R_y)^{-1}(\zeta I+A)^{-1} e^{-yA}B
{{d\zeta}\over{2\pi i}}\cr
&=\int_{\vert \zeta\vert =r}
e^{\zeta (z-y)}e^{-zA}BCe^{-yA}(I+R_y)^{-1}(\zeta I+A)^{-1} e^{-yA}B
{{d\zeta}\over{2\pi i}}\cr
&=e^{-zA}BCe^{-yA}(I+R_y)^{-1}e^{-zA+yA}&(7.23)\cr}$$
\noindent where by the calculus of residues
$$\int_{\vert \zeta\vert =r}e^{-yA}(\zeta
I-A)^{-1}BCe^{-yA}(I+R_y)^{-1}(\zeta I+A)^{-1} e^{-yA}B
{{d\zeta}\over{2\pi i}}=0.\eqno(7.24)$$
\noindent So when we multiply (7.14) by $Ca^{-xA}(I+R_x)^{-1}$ and integrate with respect to $z$ from $z=y$ to $z=x$, we
find that
$$\eqalignno{Ce^{-xA} (I+R_x)^{-1}
 \int_{\vert \zeta\vert =r}&
e^{\zeta (x-y)}e^{-xA}(\zeta
I-A)^{-1}BCe^{-yA}(I+R_y)^{-1}(\zeta I+A)^{-1} e^{-yA}B
{{d\zeta}\over{2\pi i}}\cr
&=Ce^{-xA} (I+R_x)^{-1} \int_y^x
e^{-zA}BCe^{-yA}(I+R_y)^{-1}e^{-zA}B \,dz\cr
&=\int_y^x T(x,z)T(y,z)\, dz,&(7.25)\cr}$$
\noindent as required to prove (7.12).\par
\rightline{$\square$}\par

\vskip.05in
\noindent {\bf Definition} ($KdV$ system) A $KdV$ system consists of a family of $(2,2)$
admissible linear systems $\Sigma(t)= (-A,U(t)B, CU(t))$ with $A$ bounded, input and output space ${\bf
C}$ for all
$t\in {\bf T}^\infty_0$ where $U(t)=\exp (-\sum_{j=1}^\infty t_{j-1}A^{2j-1})$ such
that the corresponding $T$ function for $\Sigma (t)$ satisfies
$$T(x,y;t)=T(y,x;t)+\int_x^y T(x,w;t)T(y,w;t)\, dw\qquad
(x,y>0).\eqno(7.26)$$
\vskip.05in

\vskip.05in
\noindent {\bf Proposition 7.3} {\sl Suppose that $\Sigma (t)$ is a $KdV$ system. Then
$\tau (x;t)=\det (I+U(t)R_xU(t))$ satisfies the partial differential equation}
$${{\partial^3 u}\over{\partial x^3}}-4{{\partial u}\over{\partial s}}-6u{{\partial
u}\over{\partial x}}=0.\eqno(7.27)$$

\vskip.05in
\noindent {\bf Proof.} By Proposition 7.2, for a $KdV$ system the bilinear integral $S(x,t;y;t)$ vanishes. By Hirota's residue identity,
this is equivalent to the Kadomtsev--Petviashvili system of partial differential equations, starting with 
$$\Bigl( {{\partial^4}\over{\partial t_1^4}} +3{{\partial^2}\over{\partial t_2^2}}-4{{\partial^2}\over{\partial t_1\partial
t_3}}\Bigr)\log \tau (x;t)+6\Bigl( 
 {{\partial^2}\over{\partial t_1^2}}\log\tau (x;t)\Bigr)^2=0.\eqno(7.28)$$
\noindent see Appendix A to [] concerning the Hirota symbol residue identity. The proof involves computing the residue of the Sato integral at
$\zeta =\infty$ as a function of $t$. By letting $t_1=x,$ fixing $t_2=0$ and letting $t_3=s$, 
we can reduce to
$\tau_\zeta (x,s)$ for the system 
$(-A, (\zeta I+ A)(\zeta I-A)^{-1}e^{-s A^3},
e^{-sA^3}C$. Then $u_\infty (x,s)=-2{{\partial^2}\over{\partial x^2}}\log \tau_\infty
(x,s)$ evolves according to () and $\psi_\zeta (x)=e^{\zeta x}\tau_\zeta (x,s)/\tau_\infty (x,s)$ gives a
family of solutions of $\psi''_\zeta (x,s)=(\zeta^2+u_\infty (x,s))\psi_\zeta
(x,s)$.\par \par
\rightline{$\square$}\par
\vskip.05in
\indent We will return to Sato's integral in Theorem 9.3, and give conditions under which (7.26) holds.

\noindent In this section we show that $\tau$ functions for admissible linear systems have some of the main properties of Riemann theta functions on algebraic curves. The fundamental concept is that of a wave function which satisfies a Schr\"odinger equation, and such that the wave function and the potential are subject to a 
deformation by an abelian group of linear transformations. The main theorem of this section proves this, and also contains a version of Fay's identities for finite sums of pairwise products of shifted tau functions. The deformations are isomonodromic in the sense that they preserve local holomorphy, whereas they can change the spectrum and hence  are generally not isospectral. 

\noindent {\bf Definition} (Deformation) Let $({\bf T}^\infty)_0=\{ t=(t_j)_{j=1}^\infty \in {\bf
C}^\infty : \lim\sup_{j\rightarrow\infty} \vert t_j\vert^{1/j}=0\}$ be
the group of rapidly convergent complex sequences with the pointwise 
addition operation. Alternatively, one can view $({\bf T}^\infty )_0$ as the coefficients of entire functions $\sum_{j=0}^\infty t_jz^j$. Let ${\bf U}=\{ U(t): t\in ({\bf T}^\infty)_0\}$ be a group of operators on $H$, let 
$${\bf \Sigma}=\{ \Sigma (t)=(-A, U(t)B, CU(t)): t\in({\bf T}^\infty)_0\}\eqno(6.20)$$
\noindent be a set of admissible 
linear systems with corresponding algebras $\{ ({\bf A}_t, \ast, \partial ):t\in ({\bf T}^\infty)_0 \}$. A deformation of 
${\bf \Sigma}$ consists of the group action of ${\bf U}$ on ${\bf \Sigma}$ and the corresponding 
homomorphisms of algebras.\par
As in section 3, let
 $({\bf T}^\infty)_0=\{ t=(t_j)_{j=1}^\infty \in {\bf
C}^\infty : \lim\sup_{j\rightarrow\infty} \vert t_j\vert^{1/j}=0\}$ be
the group of rapidly convergent complex sequences with the pointwise 
addition operation, and for $(-A,B,C)$ as in Lemma 3.1, let
$U(t)=\exp ({\sum_{j=1}^\infty t_jA^{j}})$ and 
$$\Sigma_\zeta (t) =\bigl(-A,(\zeta I+A)(\zeta I-A)^{-1}U(t)B,
CU(t))\eqno(6.25)$$ 
\noindent with tau function $\tau_\zeta (x;t)$. For $\alpha\in {\bf C}$ such that $\vert \alpha\vert<1$, let $\{ \alpha \}=(\alpha , 0, \alpha^3/3, 0, \alpha^5/5, \dots )$. We write $\tau'_\zeta (x,t)=(\partial /\partial x) \tau_\zeta (x,t)$ for short.\par  
\vskip.05in
\noindent {\bf Definition} (Baker--Akhiezer function) The 
Baker--Akhiezer (or wave) function for the family
of linear systems $\Sigma_\zeta (t)$ is
$$\psi_\zeta (x;t)={{\tau_\zeta (x;t)}\over{\tau_\infty
(x;t)}}\exp\bigl({\zeta x+\sum_{j=1}^\infty
t_j\zeta^{j}}\bigr).\eqno(6.26)$$ 
\indent This definition is consistent with that of Ercolani and
McKean [19, p487], who use a Hankel operator in place of $R_x$. In Proposition 6.3, we give a criterion for $\psi_\zeta (x,t)$ to be a Baker--Akhiezer function in Krichever's sense. First we deal with the admissible case, in which $u(x)\rightarrow0$ as $x\rightarrow\infty$.\par
\indent Let $\partial =({{\partial}\over{\partial t_1}},{{\partial}\over{3\partial t_3}},{{\partial}\over{5\partial t_5}},\dots )$ and let $\eta (t, \zeta, )=\sum_{j=0}^\infty t_{2j+1}\zeta^{2j+1}$. According to Date, the wave function of the BKP 
hierarchy is 
$${{e^{\eta (t,\zeta )}e^{-2\eta (\partial ,\zeta^{-1})}\tau (t)}\over{\tau (t)}}={{\tau (x, t_1-{{2}\over{\zeta}} , t_3-{{2}\over{3\zeta^3}}, t_5-{{2}\over{5\zeta^5}}, \dots )}\over{\tau (t_1,t_3, \dots )}}\exp \bigl(\sum_{j=0}^\infty t_{2j+1}\zeta^{2j+1} \bigr),\eqno(6.27)$$
\noindent which is a variant on the above wave function.\par
\noindent {\bf Remark} Consider the differential equation $-f''+uf=\lambda f$ where one of the following situations occurs:\par
\indent (i) $u\in C^2({\bf R}; {\bf R})$ is periodic with period $\pi$ and the characteristic equation $\Delta (\lambda )=4$ has only finitely many simple solutions [45];\par
\indent (ii) $u$ is elliptic and the Bloch spectrum has only finitely many bands [26];\par
\indent (iii) $u\in C^\infty ({\bf R}; {\bf C})$ and the $n^{th}$ order stationary KdV equation is satisfied [25]. \par
\noindent Then there exists a hyperelliptic algebraic curve with Jacobi variety $V$ with Riemann theta function $\vartheta$ such that Schr\"odinger's equation has general solution given by an exponential factor times the quotient of $\vartheta (x-a_j)$. Consequently, the Baker--Akhiezer function can be expressed in terms of $\vartheta$. \par 
\vskip.05in

\indent Any transformation of the form 
$$ (-A,B,C)\mapsto (-A, (\beta I\pm A)(\alpha I-A)^{-1}B, C)\eqno(6.16)$$
\noindent is a composition of finitely many transformations of the form
$$ M^\alpha :(-A,B,C)\mapsto (-A, (\alpha I+A)(\alpha I-A)^{-1}B, C),\eqno(6.17)$$
$$S^\gamma :(-A,B,C)\mapsto (-A-\gamma I, B, C).\eqno(6.18)$$
\noindent The effect of $S^\gamma$ is to take $R_x\mapsto e^{-2\gamma x}R_x$.
We investigate $M^\alpha$ further.
\vskip.05in 
\noindent {\bf Definition} (Quasi periods)  Let $f:{\bf C}^\infty_0\rightarrow {\bf C}$ be a non-constant meromorphic function and let 
$$\Lambda=\{t_0\in {\bf C}^\infty_0: f(t+t_0)=f(t)e^{2\pi i(L(t,t_0)+J(t_0))}; \forall t\in {\bf C}^{\infty}_0\}\eqno(9.52)$$
\noindent where $t\mapsto L(t, t_0)$ is some linear function and $J(t_0)\in {\bf C}$ for all $t_)\in {\bf C}^\infty_0$. Then $\Lambda$ is a subgroup of ${\bf C}^\infty_0$ under addition, and we say that $t_0$ is a quasiperiod and that $f$ is quasiperiodic with respect to $\Lambda$. 
In particular the set $\Lambda_0=\{ t_0\in {\bf C}^\infty_0: f(t+t_0)=f(t);  \forall t\in {\bf C}^{\infty}_0\}$ is a subgroup of $\Lambda$, and we say that $t_0$ is a period and $f$ is periodic with respect to $\Lambda_0$. We identify ${\bf C}^g$ with $\{ (t_1, \dots, t_{g}, 0, \dots )\in {\bf C}^\infty_0\}$ and consider quasiperiods for functions $f_g: {\bf C}^g\rightarrow {\bf C}$ with the obvious changes. (The properties of $L$ and $J$ are described in Chapter VI of [39].) \par
\vskip.05in
\indent Let $V$ be a complex vector space of dimension $g$, choose a basis of $V$ as a real vector space, and let $\Lambda$ be the group generated by this basis of $V$. Then $V/\Lambda$ is a complex torus. A polarization is a postive definite bilinear form $H:V\times V\rightarrow {\bf C}$ with imaginary part $E$ such that $E$ restricted to $\Lambda\times \Lambda$ takes values in ${\bf Z}$. We extend this concept to $g=\infty$ and establish that theta functions exist for infinite dimensional Abelian varieties. Let $H_0=H_1=\ell^2({\bf N}; {\bf R})$ and let $V=H_0+iH_1=\ell^2({\bf C})$. We introduce the discrete group $\Lambda_0=H_0\cap {\bf Z}^\infty$ and observe that $H_0/\Lambda_0$ is a subgroup of the infinite real torus ${\bf T}^\infty$. Now ${\bf T}^\infty$ is a compact metrizable space, and by the Stone--Weierstrass theorem the algebra $A_0={\hbox{span}}\{ e^{2\pi \langle x, k\rangle}; k\in \Lambda_0\}$ is dense subpace of $C({\bf T}^\infty; {\bf C})$ for the topology associated with the supremum norm. Let $\pi_0: H_0\rightarrow  H_0/\Lambda_0$ be the quotient map. Now let $\Omega_1:{\cal D}(\Omega_1)\subset H_1\rightarrow H_1$ be a linear operator such that $\Lambda_0\subset {\cal D}(\Omega_1)$ and $\Omega_1$ is invertible with $\Omega_1^{-1}$ a positive and self-adjoint operator. Then let $\Omega_0:{\cal D}(\Omega_1)\rightarrow H_0$ be a linear operator which is relatively bounded with respect to $\Omega_1$ and suppose further that $\Omega_0=\sqrt{\Omega_1}Z\sqrt{\Omega_1}$ where $Z\in {\cal L}^1(H_0)$. Then we introduce $\Omega=\Omega_0+i\Omega_1:{\cal D}(\Omega_1)\rightarrow V$,  then let $\Lambda =\Lambda_0+\Omega\Lambda_0$.\par
We now intepret 
$$\vartheta (z)=\sum_{n\in \Lambda_0}e^{2\pi i\langle z,n\rangle +\pi i\langle \Omega n,n\rangle}\eqno(9.53)$$
\noindent and establish its basic properties.\par

\vskip.05in
\noindent {\bf Lemma 9.5} {\sl (i)  Then   $L(z,w)=-\langle z,\Omega_1^{-1}\Im w\rangle$ is a bounded real bilinear form $L: V\times V\rightarrow {\bf C}$ which is complex linear in the first argument and satisfies $L(z, m+\Omega n)=-\langle z,n\rangle$ for all $z\in V$ and $n\in \Lambda_0$. \par 
\indent (ii)  $E(z,w)=L(z,w)-L(w,z)$  is a bounded and real-bilinear form  $E:V\times V \rightarrow {\bf R}$ which is skew symmetric, and $E$ restricted to $\Lambda\times \Lambda$ takes values in ${\bf Z}$;\par
\indent (iii)   $H(z,w)=E(iz,w)+E(w,iz)$  is a bounded and real-bilinear form $H:V\times V\rightarrow {\bf C}$ which is Hermitian and non-negative, and has imaginary part $E$;\par
\indent (iv) $\vartheta$ is quasiperiodic with respect to $\Lambda$ and satisfies
 $$\vartheta (z+k+\Omega m)=e^{2\pi iL(z, k+\Omega m)-\pi i\langle \Omega m,m\rangle}\vartheta (z);$$
\indent (v) For each $z\in H_0$, there exists a bounded measure $\varpi_x$ on ${\bf T}^\infty$ that has Fourier expansion}
$$\sum_{k\in \Lambda_0} e^{2\pi i\langle x,k\rangle -\pi \langle \Omega_1^{-1}(k-z), k-z\rangle}\qquad (x\in {\bf T}^\infty ).\eqno(9.54)$$
\vskip.05in
\noindent {\bf Proof} (ii) We have $E(x+iy,\xi +i\eta )=-\langle \Omega_1^{-1}x,\eta \rangle +\langle \Omega_1^{-1}\xi ,y\rangle$ for all $x,y,\eta,\xi\in H_0$. 
 It follows from the symmetry of $\Omega_0$ and $\Omega_1$ that $E(m_1+\Omega n_1, m_2+\Omega n_2)=-\langle m_1,n_2\rangle +\langle m_2, n_1\rangle$ for all $n_1,n_2,m_1,m_2\in \Lambda_0$.\par
\indent (iii) With the usual inner product on $V$, we can write $H(z,w)=\langle z, \Omega_1^{-1}w\rangle.$\par
\indent (iv)  By standard results, including Sazanov's theorem, there exists a Radon probability measure $\mu_0$ on $H_0$ such that the Fourier transform satisfies 
$$e^{-\pi\langle \Omega_1^{-1}w,w\rangle}=\int_{H_0} e^{-2\pi i\langle w, \zeta\rangle} \mu_0(d\zeta ).$$
\noindent Hence there exists a bounded Radon measure $$d\mu =\det (I-iZ)^{-1/2}e^{i\langle \Omega_0\xi,\xi\rangle}d\mu_0\eqno(9.55)$$
\noindent on $H_0$; likewise, we can consider $e^{2\pi i\langle z,\xi\rangle} d\mu$ for all $z\in V$.  By standard properties of Gaussian measures, there exists $\varepsilon >0$ such that $\int_{H_0}e^{\varepsilon \Vert \zeta\Vert} \mu_0(d\zeta )$ is finite, so the series converges absolutely, and we can 
rearrange terms in the sum, and verify the stated result.\par
\indent  (v)  The reader is invited to think of $\mu$ as $\sqrt{\det (\Omega_1)}e^{-\pi\langle \Omega \xi, \xi\rangle} d\xi$ as motivation for the following computation.
 The quotient map $\pi_0: H_0\rightarrow {\bf T}^\infty$ induces a probability measure $\nu_0=\pi_0\sharp\mu_0$ on ${\bf T}^\infty$ such that 
$$\int_{{\bf T}^\infty} f(x)\nu_0(dx)=\int_{H_0} f(\pi_0(\zeta ))\mu_0(d\zeta )\eqno(9.56)$$
\noindent for all $f\in {C({\bf T}^\infty }; {\bf C});$ likewise, $\pi_0$ induces a bounded measure $\omega_z$ on ${\bf T}^\infty$ from $e^{2\pi i\langle z, \xi\rangle}d\mu_0$. In particular, we can choose $f(x)=e^{2\pi i\langle x,k\rangle}$ for each $k\in \Lambda_0$  to obtain $\hat \nu (k)=e^{-\pi\langle \Omega_1^{-1}(k-z),k-z\rangle}$. \par

 \vskip.05in
\noindent {\bf Proposition 9.6} {\sl Suppose that $A_1, A_2\in {\cal L}(H)$ are invertible algebraic operators, so that $A_1$ and $A_2$ both satisfy nontrivial polynomial equations.\par
\indent (i)  Then there exists a finite-dimensional Lie group ${\bf C}^g$ with an action $\tilde\rho$ on ${\bf {\Sigma}}_{A_1,A_2}$ such that the orbits of $\tilde\rho$ coincide with the orbits of $\rho$, and $\tau$ is periodic with respect to a nontrivial discrete subgroup of  ${\bf C}^g$.}\par
\indent {\sl (ii) $t\mapsto \tau (t)$ is entire and $t\mapsto u(t)$ is meromorphic, and every quasiperiod of $\tau$ is a period of $u$.}\par
\indent {\sl (iii) $\zeta\mapsto \psi_\zeta (x,t)$ is holomorphic except at finitely many points in ${\bf C}$.}\par
\vskip.05in

\noindent {\bf Proof.} (i) Under these hypotheses, $\{ p(X)\in {\bf C}[X]: p(-A_1)=p(A_2)=0\}$ is a nontrivial ideal in ${\bf C}[X]$; for instance, it contains the product of the minimal polynomials of $-A_1$ and $A_2$; hence this ideal is generated by 
 a unique monic polynomial of degree $g$.   Since $A_1$ and $A_2$ are invertible, we can assume that $p(x)=x^g+a_{g,g-1}x^{g-1}a_{g,g-2}x^{g-2}+\dots+ a_{g,0}$ where $a_{g,0}\neq 0$. Then by induction 
one shows that there exists $M>0$ such that for all $k\geq 1$ and $j=0, \dots , g-1$  there exist $a_{k,j}\in {\bf C}$ with $\vert a_{k,j}\vert \leq M^k$ and $A_2^k=\sum_{j=0}^{g-1} a_{k,j} A_2^j$ and $(-A_1)^k=\sum_{j=0}^{g-1}a_{k,j}(-A_1)^j$. Let $\alpha_j=(a_{k,j})_{k=1}^\infty$ and write $\langle \sigma, \alpha_j\rangle =\sum_{k=1}^\infty s_ka_{k,j}$; then 
$$\sum_{k=1}^\infty s_kA_2^k=\sum_{j=0}^{g-1} \langle \sigma, \alpha_j\rangle A_2^j, \quad -\sum_{k=1}^\infty s_k(-A_1)^k=-\sum_{j=0}^{g-1} \langle \sigma, \alpha_j\rangle (-A_1)^j\qquad (t\in {\bf C}^\infty _0),\eqno(9.57)$$ 
\noindent and $\sum_{j=0}^{g-1}\langle \sigma ,\alpha_j\rangle A_2^j=0=\sum_{j=0}^{g-1}\langle \sigma ,\alpha_j\rangle (-A_1)^j$  implies $\langle \sigma, \alpha_j\rangle =0$ for all $j=0, \dots, g-1$. Hence the map $Q: {\bf C}^\infty_0\rightarrow {\bf C}^g$ $\sigma\mapsto (\langle \sigma, \alpha_j\rangle )_{j=0}^{g-1}$ is surjective. We can introduce $\sigma_\ell\in {\bf C}^\infty_0$ such that $\langle \sigma_\ell, \alpha_j\rangle =\delta_{\ell, j}$ for $j, \ell =0, \dots, g-1$; then let $J:{\bf C}^g\rightarrow {\bf C}^\infty_0$ $(c_\ell )_{\ell =0}^{g-1}\mapsto \sum_{\ell =0}^{g-1} c_\ell \sigma_\ell$ so $QJ=I$. Then $\tilde \rho =\rho \circ J$ gives an action of ${\bf C}^g$ on ${\bf {\Sigma}}_{A_1,A_2}$ . Finally, observe that 
$$V(J (c_\ell )_{\ell =0}^{g-1})=\exp \Bigl({\sum_{\ell =0}^{g-1} c_\ell A_2^\ell}\Bigr), \quad W(J (c_\ell )_{\ell =0}^{g-1})=\exp\Bigl({-\sum_{\ell =0}^{g-1} c_\ell (-A_1)^\ell}\Bigr),\eqno(9.58)$$
\noindent so $\rho (\sigma )=\tilde \rho (Q\sigma )$ for all $\sigma\in {\bf C}^\infty _0$.\par 
\indent We identify ${\bf Z}^g$ with the diagonal $g\times g$ matrices with entries in ${\bf Z}$ and consider 
$$\Lambda_0 =\bigl\{ (p_1, \dots , p_g)\in {\bf C}^\infty_0: \sum_{j=1}^g p_jA_2^j, \sum_{j=1}^g p_j(-A_1)^j\in 2\pi i {\bf Z}^g\bigr\}\eqno(9.59)$$,
and observe that $\Lambda_0$ is an additive subgroup of ${\bf C}^g$. Also, $$-2\pi i (A_2^g+a_{g,g-1}A_2^{g-1}+\dots +a_{g,1}A_2)=2\pi i I,\eqno(9.60)$$
\noindent so $\Lambda_0 $ contains non-zero elements, and by the minimality of $p$, it is clear that $\Lambda_0$ is a 
discrete subset of ${\bf C}^g$. Identifying ${\bf C}^g$ with the subspace of ${\bf C}^\infty_0$ formed from the first $g$ coordinates, we can write $V(p)=W(p)=I$ for all $p\in \Lambda_0$, hence $\tau (x, t+p)=\tau (x,t)$, and $\tau$ is periodic with respect to $\Lambda_0$.\par
\indent (ii) Here $t\mapsto V(t)S_0W(t)$ is an entire function with values in ${\cal L}^1(H)$, so $\tau (t)= \det (I+V(t)S_0W(t))$ is entire on ${\bf C}^n$. Hence $u(x,t)=-2{{\partial^2}\over{\partial x^2}}\tau (x;t)$ is a quotient of entire functions, hence meromorphic. By taking the logarithmic derivative of  (9.32),  with $f=\tau$, we see that $u$ is periodic.\par 
\indent (iii) The function $\psi_\zeta$ is holomorphic on ${\bf C}$ except at the points in ${\hbox{Spec}}(A_2)$, which is a finite set.\par 
\vskip.05in

\indent A compact and connected Riemann surface is equivalent to a complete algebraic curve  ${\cal C}$ over ${\bf C}$. Suppose that  ${\cal C}$ cannot be parametrized by rational functions 
on ${\bf C}$, hence has genus $g\geq 1$. Then by Riemann's theory, one can choose a basis $\{ \omega_1, \dots, \omega_g\}$ of the space of holomorphic differential one-forms on ${\cal C}$ and the homology basis $\{ A_1, \dots ,A_g; B_1, \dots ,B_g\}$ of 
$H_1({\cal C}; {\bf Z})$ consisting of one cycles on ${\cal C}$ so that 
$$\int_{A_j}\omega_k=\delta_{jk};\eqno(9.61)$$
$$\Bigl[\int_{B_k}\omega_j\Bigr]=\Omega\eqno(9.62)$$
\noindent where $\Omega$ is symmetric with $\Im \Omega$ positive definite. Then the ${\bf Z}$-combinations of the columns of the $g\times (2g)$ matrix $[I\quad \Omega ]$ generate a lattice 
$\Lambda ={\bf Z}^g+\Omega {\bf Z}^g$ in $V={\bf C}^g$; see Corollary 2.3 of Mumford [51;  I]. Then $(V, \Lambda )$ is an Abelian variety known as the Jacobian of ${\cal C}$. One can construct suitable entire functions such as 
$$\vartheta (z)=\sum_{n\in {\bf Z}^g} e^{2\pi i\langle z, n\rangle +\pi i\langle \Omega n,n\rangle}\eqno(9.63)$$
\noindent known as theta functions, so that $\vartheta$ is quasiperiodic on $\Lambda$ with $H=(\Im \Omega )^{-1}$, and purely periodic for the subgroup ${\bf Z}^g$ of $\Lambda$,  as in [Mumford theta I]. \par 
\vskip.05in

\vskip.05in
\noindent {\bf Corollary 9.7} {\sl Let $(V,\Lambda )$ be a finite-dimensional abelian variety with theta function $\vartheta$, let ${\bf C}^N$ be the subspace of ${\bf C}^\infty_0$ generated by the first $N$ standard basis vectors and suppose that there exists a surjective linear map $\phi :{\bf C}^N\rightarrow V$ such that 
$$\tau (0;t)=\vartheta (\phi (t))\qquad (t\in {\bf C}^N).\eqno(9.64)$$
\noindent Then there exists a compact and connected Riemann surface such that $(V, \Lambda )$ is the Jacobian of ${\cal C}$.}\par
\vskip.05in
\noindent {\bf Proof.} By Theorem 9.3 and [1], $\tau$ satisfies the $KP$ hierarchy of differential equations.
The solutions to $KP$ apparently depend upon infinitely many variables
$t=(t_1, \dots )$. Shiota
[59] introduced the notion of a finite-dimensional solution $L$ to the $KP$
equations, and introduced a complex abelian Lie group ${\cal T}_L$ which effectively
meromorphically parametrizes $L$. If ${\cal T}_L$ is compact, then there exists a
smooth, complete, reduced and irreducible curve ${\cal C}$ over ${\bf C}$ such that
${\cal T}_L$ is isomorphic to the Jacobian of ${\cal C}$.  We only need finitely many differential equations to characterize the Jacobians of compact algebraic curves, hence can work on some subspace ${\bf C}^N$ of ${\bf C}^\infty_0$. Shiota and Mulase have shown that $\vartheta \circ\phi$ arises as the theta function on a  Jacobian. \par 
\rightline{$\square$}\par

\indent Suppose that $A_1=A_2=A$ is algebraic, and let $U(t)=\exp (\sum_{j=1}^\infty t_{2j-1}A^{2j-1})$ be the KdV group acting on a $(2,2)$ admissible linear system $(-A,B,C)$. Then there exists $m$ such that ${\bf A}={\bf C}[u, \partial u/\partial x, \dots, \partial^mu/\partial x^m]$ is a differential ring for $\partial /\partial x$ and the usual multiplication. Then $\tau (x,t)$ is an entire function such that 
$t\mapsto \tau (x,t)$ has quasi-periods in ${\bf C}^g$, and $t\mapsto u(x,t)$ is meromorphic and has a nontrivial period $p_0$. Suppose further that $(\pi, 0, \dots ,0)$ is a period for $u$, and consider Hill's equation
  $-\psi''(x)+u(x;t)\psi (z)=\lambda\psi (x)$. Let  $\Delta (\lambda ;t)$ be the discriminant of Hill's equation, and ${\cal E}= \{ (X,Z): Z^2=4-\Delta (X)^2\}$ the corresponding spectral curve. Then the Bloch spectrum is invariant under the KdV flow associated with $t\mapsto u(x,t)$. If ${\cal E}$ is finite gap, then  $\psi_\zeta (x;t)=e^{\zeta x}\tau (x;t+[1/\zeta ])/\tau (x; t)$ can be expressed as a quotient of Riemann theta functions for a hyperelliptic algebraic curve. \par

\noindent {\bf Corollary 10.3.} {\sl (v) Suppose that $e^{-A\pi /2}Ee^{-A\pi/2}=-E$. Then}
$$\Phi (x+y)+T(x,y)-{{1}\over{2}}\int_x^{x+\pi /2} T(x,z)\Phi (z+y)\,
dz=0.\eqno(10.15)$$ 
\indent {\sl (ii) If $E$ has finite rank, then $\tau_\infty$ is of
exponential type and in ${\bf C}_{\cal C}$. Conversely, if 
$\tau_\infty$ is of exponential type, then there 
exist $\alpha_j\in {\cal C}$, $\alpha\in {\bf Z}$ and $\beta\in {\bf C}$, not necessarily distinct, 
such that
$$\tau_\infty (z)=e^{2i\alpha z+\beta}\prod_{j=1}^m \sin 2(z-\alpha_j)
\eqno(10.16)$$
\noindent and} 
$$u(z)=\sum_{j=1}^m {{8}\over{\sin^22(z-\alpha_j)}}.\eqno(10.17)$$
\vskip.05in
\noindent {\bf Proof.} (i) One can verify this by direct computation, and the crucial
 identity is
$$\int_x^{x+\pi /2} e^{-zA}BCe^{-zA}\, dz=\bigl[
-e^{-zA}Ee^{-zA}\bigr]_{x}^{x+\pi /2} =2e^{-xA}Ee^{-xA}.\eqno(10.)$$
\indent (ii) There exists a projection $P$ of finite 
rank $\rho$ such that $PEP=E$ and hence $\tau_\infty (z)=\det (I+PEPe^{-2zA}P)$,
where $Pe^{-2zA}P$ is a finite matrix with entries that are in ${\bf
C}_{\cal C}$; in particular, the entries are functions of exponential growth. Hence from the expansion of this
determinant, we deduce that there exist $c_1, c_2>0$ such that $\vert
\tau_\infty (z)\vert \leq c_1e^{2\rho N\vert z\vert +c_2}$ for all $z$.\par

\indent Suppose conversely that $\tau$ is of exponential type, and recall a standard argument of function theory. By
Jensen's formula, the number of zeros of $\tau_\infty$ inside a circle of
radius $r$ grows like $c_3r+c_4$ for some $c_3, c_4>0$. Likewise, the number of zeros inside $\{ z\in {\bf C}: -r<\Re z<r; -r<\Im z<\Im r\}$ is less than $2c_3+c_4$ and since
$\tau_\infty$ is also $\pi$-periodic, we deduce that there exists $m<\infty$ 
such that the only zeros of $\tau_\infty$ in $\{ z: -\pi /2 <\Re
z\leq \pi/2\}$ are  $\alpha_1, \dots , \alpha_m$; hence $\tau_\infty (x)/\prod_{j=1}^m \sin 2(z-\alpha )$ is $\pi$-periodic, entire and without any zeros. Hence
there there exists an entire function $g$ such that 
$$\tau_\infty (z) =e^{g(z)}\prod_{j=1}^m \sin 2(z-\alpha_j),\eqno(9.26)$$
\noindent where  $g(z+\pi )-g(z)=2\pi i\ell$ for some $\ell\in {\bf Z}$.
Since $\vert \sin (x+iy)\vert\rightarrow \infty$ as
$y\rightarrow\infty$, we deduce that $\vert g(z)\vert\leq c_5\vert z\vert
+c_6$ for some $c_5,c_6>0$, and we finally obtain $g(z)=2i\alpha z+\beta$
where $\alpha\in {\bf Z}$.\par

\rightline{$\square$}\par
\noindent {\bf Remark 9.4}  The potential (5.19) can be interpreted in 
terms of a simple model in
electrodynamics, considered by Sutherland [61]. Consider $m$ fixed unit charges placed at points $e^{i\alpha_j}$ on a 
circular ring, and a further unit charge which has variable position $e^{ix}$ on the ring.
 Then the electrostatic energy of the moving charge is $u$. 
In particular, when $m=1$ we have 
$u(z)=8/\sin^22(z-\alpha_1)$, so we can
rescale this to the familiar case of $C{\hbox{sech}}^2z$ for some $C$.\par
In section 10, we show how this can otherwise be realised as a limiting case of periodic 
linear systems with elliptic potentials.\par
\vskip.05in
\noindent If we assume more commutativity, the proofs simplify and the
results become stronger.\par
\vskip.05in

\noindent {\bf Proposition 9.5} {\sl Let $(-A,B,C;E)$ be a uniformly periodic system such that $AE=EA$, where} ${\hbox{Spec}}(E)=\{ \varepsilon_j:j\in {\bf N}\}$. \par
\indent  (i) Then there exist $n(j,m)\in {\bf N}\cup \{ 0\}$ such that 
$$\tau (x)=\prod_{m=-N}^N \prod_{j=1}^\infty (1+\varepsilon_j e^{2imx})^{n(j,m)}\eqno(9.27)$$
\noindent \noindent {\sl with uniform convergence on compact subsets of ${\bf C}$.}\par
\indent {\sl (ii) $(-A,B,C)$ is finitely generated,
 since the algebra ${\bf S}$ is commutative and Noetherian, and is 
complex state ring for $(-A,B,C)$ on ${\cal C}'={\bf C}/2\pi {\bf Z}$.}\par
\indent {\sl (iii)  Let ${\bf K}={\bf C}(e^{ix})$ and let $\varepsilon$ be a nonzero element of} ${\hbox{Spec}}(E)$. {\sl Then there exists} $\alpha\in {\hbox{Spec}}(A)$ {\sl and a nonzero $\xi\in H$ such that $A\xi=\alpha \xi$ and 
$F\xi =\nu\xi$ where $\nu =(1+e^{-2\alpha x}\varepsilon)^{-1}\in {\bf K}$.  Then
${\bf p}= (\alpha , \nu )$ determines a maximal ideal in ${\bf K}[A,F]$ and determines a finite algebraic cover of the
Riemann sphere.}\par

\vskip.05in
\noindent {\bf Proof.} (i) By Sz.-Nagy's theorem, $A$ is similar to a skew-symmetric operator with ${\hbox{Spec}}(A)\subseteq\{ -N, \dots ,N\}$; hence there exist 
$P_m\in {\cal L}(H)$ such that $P_m^2=P_m$, $AP_m=imP_m$, $I=\sum_{m=-N}^N P_m$ and $P_mP_r=0$ for all $m\neq r$. Furthermore, $EP_m=P_mE$ for all $m$, and each $EP_m$ has a 
Riesz decomposition, where each non-zero point in the spectrum is an eigenvalue with ${\hbox{Spec}}(E)=\cup_{m=-N}^N {\hbox{Spec}}(P_mE)$. 
Hence there is a determinant given by a product
$$\eqalignno{\tau (x)&=\det (I+e^{-xA}Ee^{-xA})\cr
&=\prod_{m=-N}^N \det (I+e^{-2im x}P_mE)\cr}$$  
\noindent and each factor is a product over the eigenvalues of $P_mE$.\par

\indent (ii) Here  $e^{-xA}$ is a polynomial 
in $A$, $e^{ix}$ and $e^{-ix}$,
hence $e^{-xA}$ and likewise $R_x$ belong to ${\bf K}_{{\cal C}'} [I,E,A]$.
Observe that the set $S=\{ (I+e^{-xA}Ee^{-xA})^n: n=0, 1,
\dots \}$ is multiplicatively closed and does not contain $0$ since $
I+e^{-xA}Ee^{-xA}$ is invertible in the Calkin algebra of
${\cal L}(H)$ modulo the compact operators on $H$. Hence we can
identify ${\bf S}$ with the ring of fractions of ${\bf K}_{{\cal
C}'}[A, BC]$ modulo $S$. There is a natural surjective ring homomorphism 
${\bf K}_{{\cal C}'} [X_1, X_2,X_3]\rightarrow {\bf S}$ given 
by $X_1\mapsto A$, $X_2\mapsto BC$, $X_3\rightarrow F_x$, so 
by Hilbert's basis theorem, ${\bf S}$ is Noetherian as a
commutative ring.\par 
\indent (iii) For each
non zero $\varepsilon$ in the spectrum of the compact operator $E$, the operator $A$ has a finite dimensional invariant subspace $\{\xi \in H: E\xi =\varepsilon \xi\}$ and hence 
there exists a common eigenvector $\xi\in H$ such that $A\xi=\alpha\xi$ and
$E\xi=\varepsilon \xi$, so $F\xi =(1+e^{-2\alpha x}\varepsilon
)^{-1}\xi$. Now observe that  ${\bf K}[A,F]$ is a finitely generated commutative algebra over the field ${\bf K}$, and that $Q(A,F)\xi =Q(\alpha,\nu )\xi$, so 
${\bf M}=\{ Q(A,F): Q(\alpha ,\nu )=0; Q\in {\bf K}[X,Y]\}$ is a maximal ideal. By 
the weak Nullstellensatz, ${\bf K}[A,F]/{\bf M}$ is a finite algebraic extension of ${\bf K}$. \par
\indent Observe that ${\bf K}$ may be transformed to ${\bf C}(t)$ by the rational substitution $t=\tan (x/2)$, so the finite algebraic extension ${\bf K}[A,F]/{\bf M}$ of ${\bf C}(t)$,
which may be realised as an algebraic function. Thus there exists a Riemann
surface ${\bf X}$ and a holomorphic and surjective map $\pi :{\bf X}\rightarrow {\bf
P}^1$ which gives a cover of the Riemann sphere with finitely many 
sheets.\par
\rightline{$\square$}\par
\vskip.05in
  \indent The following was considered by Brett in her PhD thesis. Let $H=L^2({\bf R}/\pi {\bf Z}; d\theta/\pi )$ be the space of $\pi$-periodic $L^2$ functions on the circle, and 
$${\cal D}(A)=H^1=\{ f\sim\sum_{n=-\infty}^\infty b_ne^{2in\theta }:b_n\in {\bf C};\sum_{n=-\infty}^\infty (1+n^2)\vert b_n\vert^2<\infty\};$$
\noindent  let $0<\varepsilon <1/2$ be rational. We consider the linear system 
$$\eqalignno{ e^{-tA}:f(x)\mapsto e^{-i\varepsilon t}f(x-t)&\qquad (f\in H);\cr
                       B: b\mapsto \phi (x)b/2&\qquad b\in {\bf C};\cr
                       C: f(x)\mapsto f(0)&\qquad  f\in {\cal D}(A),&(9.28)\cr}$$
\noindent where $\phi\in {\cal D}(A)$ is even and $\pi$-periodic. We regard $\phi$ as a scattering function, and in section 10 justify this terminology. Then we introduce the operators
$$S_x=\int_{-x}^x (e^{tA}+e^{-tA})BC(e^{tA}+e^{-tA})\, dt\in {\cal L}^2(H)\eqno(9.29)$$
$$V_x=\int_{-x}^x (e^{tA}-e^{-tA})BC(e^{tA}-e^{-tA})\, dt\in {\cal L}^2(H);\eqno(9.30)$$
\noindent which may be expressed as integral operators
$$S_xf(z)={{1}\over{2}}\int_{-x}^x \bigl(e^{i\varepsilon t} \phi (z+t)+e^{-i\varepsilon t}\phi (z-t)\bigr)\bigl(e^{i\varepsilon t} f(t)+e^{-i\varepsilon t}f(-t)\bigr)\, dt\qquad (f\in H),\eqno(9.31)$$
$$V_xf(z)={{1}\over{2}}\int_{-x}^x \bigl(e^{i\varepsilon t} \phi (z+t)-e^{-i\varepsilon t}\phi (z-t)\bigr)\bigl(e^{i\varepsilon t} f(t)-e^{-i\varepsilon t}f(-t)\bigr)\, dt\qquad (f\in H).\eqno(9.32)$$
\noindent Then we let
$$T(x,y)=-C(e^{xA}+e^{-xA})(I+S_x)^{-1}(e^{yA}+e^{-yA})B,\eqno(9.33)$$
$$M(x,y)=-C(e^{xA}-e^{-xA})(I+V_x)^{-1}(e^{yA}-e^{-yA})B;$$
\noindent so
$$T(x,x)=-{{d}\over{dx}}\log \det (I+S_x)$$
$$M(x,x)=-{{d}\over{dx}}\log\det (I+V_x).$$
\indent For $0<\varepsilon <1/2$ rational,   let $E$ be the operator on $H$ which has matrix
$$E=\Bigl[{{\hat \phi (j)}\over{4i(\varepsilon +j+k)}}\Bigr]_{j,k\in {\bf Z}}\eqno(9.35)$$
\noindent with respect to the orthonormal basis $(e^{2ij\theta })_{j\in {\bf Z}})$. \par
\vskip.05in
\noindent {\bf Lemma 9.6} {\sl (i) The group  $(e^{xA})_{x\in {\bf R}}$ is strongly continuous, isometric and periodic.\par
\indent (ii) For $0<\varepsilon <1/2$ rational,  $E$ is trace class, $AE+EA=BC, $  and the function $\tau (x)=\det (I+e^{xA}Ee^{xA})$ is continuous and periodic on ${\bf R}$. \par
 \indent (ii) Let $\varepsilon =0$.  Then $T$ satisfies $T(0, y)=-2\phi (y)$ and the integral equation}
$$\phi (x+y)+\phi (x-y)+T(x,y)+\int_{-x}^x T(x,z)\bigl( \phi (z+y)+\phi (z-y)\bigr)\, dz=0.\eqno(9.36)$$
\vskip.05in
\noindent {\bf Proof.}  (i) Strong continuity follows from a standard result on rotation of the torus, as in [18].  One can check that $e^{\pi A}=e^{\pi i\varepsilon }I$, so $e^{\pi NA}=I$ for some $N\in {\bf N}$.\par
\indent (ii) When $\phi\in {\cal D}(A)$, the series $\sum_{j\in {\bf Z}}\vert \hat \phi (j)\vert$ converges. Also $\sum_{k\in {\bf Z}}1/(\varepsilon +j+k)^2$ is finite and independent of $j$. 
Hence we can express $E$ as a convex combination of  bounded linear operators on $\ell^2$, all of rank one; hence $E$ is trace class. It is clear that $(AE+EA)$ and $BC$ map $e^{ij\theta}$ to the same element.\par
\indent (iii) We have $\cos \varepsilon y\phi (y)=C(e^{yA}+e^{-yA})B$, so  $\phi (y)= C(e^{yA}+e^{-yA})B$ when $\varepsilon =0$, and $T(0,y)=-2C(e^{yA}+e^{-yA})B$ since $S_0=0$. The integral equation involves a restriction of the kernel, and here $-x<z<x$.\par
\rightline{$\square$}\par
\indent We now derive some consequences when $\varepsilon =0$. By repeatedly differentiating and integrating by parts, we deduce that  
$${{\partial^2T}\over{\partial x^2}}-{{\partial^2T}\over{\partial y^2}}=u(x)T(x,y)\eqno(9.37)$$
\noindent where $u(x)=2{{d}\over{dx}}(T(x,x)+T(x,-x))$, and $u(x)=-2{{d^2}\over{dx^2}}\log\det (I+S_x)$.\par
\indent Then we introduce $h$ such that $T(x,x)+T(x,-x)=h+(1/2)\int_0^xu(t)\, dt$.
Then one introduces 
$$\varphi_\zeta (x)=\cos \zeta x +\int_{-x}^x T(x,y) e^{i\zeta y}\, dy\eqno(9.38)$$
\noindent which satisfies the initial conditions $\varphi_\zeta (0)=1$ and  $\varphi_\zeta '(0)=h$, and the differential equation
$$-{{d^2\varphi_k}\over{dx^2}}+u(x)\varphi_k(x) =k^2\varphi_k(x).\eqno(9.39)$$

\indent The operators are linked by a Lyapunov equation
$${{d}\over{dx}}\left[\matrix{W_x&Z_x\cr Z_x&W_x}\right] = \left[\matrix{0&A\cr A&0}\right]  \left[\matrix{W_x&Z_x\cr Z_x&W_x}\right]+  \left[\matrix{W_x&Z_x\cr Z_x&W_x}\right] \left[\matrix{A&0\cr 0&A}\right]+\left[\matrix{0&2BC\cr 2BC&0}\right] .$$
For $\varepsilon>0$, the Sylvester equation
$$\left[\matrix{\varepsilon I&A\cr A&\varepsilon I\cr}\right] \left[\matrix{ E&F\cr F&E\cr}\right] 
+ \left[\matrix{ E&F\cr F&E\cr}\right]\left[\matrix{\varepsilon I&A\cr A&\varepsilon I\cr}\right] = \left[\matrix{ BC&0\cr 0&BC\cr}\right]$$
\noindent has a solution since
$$\exp\Bigl( -t\left[\matrix{\varepsilon I&A\cr A&\varepsilon I\cr}\right]\Bigr) =e^{-t\varepsilon} \left[\matrix{\cosh (tA)&\sinh (tA)\cr \sinh (tA)&\cosh (tA)\cr}\right]$$
\noindent is an admissible semigroup. To solve this in the particular case of interest, let $F=[F_{mn}]$ and $E=[E_{mn}]$ be the matrices with respect to the orthonormal basis $(e^{2in\theta})_{n=-\infty}^\infty$ of $L^2([0, \pi ]; d\theta/\pi )$ so that 
$$\eqalignno{ 2\varepsilon E_{mn}+(2im+2in)F_{mn}& = \hat\phi (m)\cr
2\varepsilon F_{mn}+(2im+2in)E_{mn}&=0\cr}$$ 
\noindent with solution
$$E=\Bigl[{ {\varepsilon \hat \phi (m)}\over{2\varepsilon^2+2(m+n)^2}}\Bigr], \quad F=\Bigl[ {{-i(m+n)\hat\phi (m)}\over{ 2\varepsilon^2+2(m+n)^2}}\Bigr].$$
\noindent The operators $E$ and $F$ are trace class for all $\phi\in {\cal D}(A)$, and 
$$\int_0^x e^{tA}BCe^{tA}\, dt=e^{xA}Fe^{xA}-F+2\varepsilon \int_0^x e^{tA}Ee^{tA}\, dt.$$
\vskip.05in
\noindent {\bf Proposition} {\sl Let $h=2CB$. Then the eigenvalue problem
$$\eqalignno{-\psi''+w\psi&\lambda \psi\cr
h\psi (0)+\psi'(0)&=0\cr
h\psi (\pi )+\psi'(\pi )&=0\cr}$$
\noindent has a sequence of eigenvalues $\lambda=\nu_j^2$ such that\par
\indent (i) $(e^{\pm i\nu_jx})_{j=0}^\infty$ is a Riesz basis for $L^2[-\pi ,\pi ]$;\par
\indent (ii) $\psi (x; \nu_j)$ is an eigenfunction with $\psi (0; \nu_j)=1$;\par
\indent (iii) there exist $(a_j)\in \ell^2$ such that $\phi (x)=\sum_{j=0}^\infty a_j\cos \nu_j x$ and }
$$K(0, y)=\sum_{j=0}^\infty -2a_j \cos \nu_j y, \qquad K(\pi , y)=\sum_{j=0}^\infty -2a_j\psi (\pi , \nu_j) \cos \nu_j y.$$ 
\vskip.05in
\noindent {\bf Proof}. The following observation is based upon Gelfand and Levitan's paper. As in (i), we can write $\phi (x)=\sum_{j=1}^\infty a_n\cos\sqrt{\nu_n}x$. Then $T(0,y)=-2\phi (y)$ and 
$$T(\pi ,y)=-2C(I+S_\pi )^{-1}(e^{yA}+e^{-yA})B=-2\sum_n b_n\cos\sqrt{\nu_n}y,\eqno(10.10)$$
\noindent where from the Gelfand--Levitan equation
$$\eqalignno{ 0&=K(\pi , y)+2\sum_n a_n\cos\sqrt{\nu_n}\pi\cos\sqrt{\nu_n}y +2\int_{-\pi}^\pi K(\pi ,z)\sum_n a_n\cos\sqrt{\nu_n}z\cos \sqrt{\nu_n}y\, dz\cr
&=K(\pi, y)+2\sum_n a_n\psi_{\nu_n}(\pi)\cos\sqrt{\nu_n}y.(10.11)\cr}$$
\rightline{$\square$}\par
\vskip.05in
\noindent {\bf Corollary} {\sl Suppose further that $w$ is $\pi$-periodic. Then $\psi (x; \nu_k)$ is a Floquet solution with multiplier}
$$\rho (\nu_k^2)={{\Delta (\nu_k^2)\pm \sqrt{\Delta (\nu_k^2)^2-4}}\over{2}}.$$ 
\vskip.05in
\noindent {\bf Proof.} Since $w$ is periodic, $\psi (x+\pi ; \nu_k^2)$ satisfies the same boundary value problem, and hence is a multiple of  $\psi (x; \nu_k^2)$. Then ${\hbox{col}}[\psi (0; \nu_k^2), \psi' (0; \nu_k^2)$ is an eigenvector of the monodromy matrix, and the eigenvalues are given by the above formula.\par
\rightline{$\square$}\par
\vskip.05in

\noindent {\bf 10.  Hill's equation on the torus}\par
\vskip.05in
\indent In this section we periodic linear systems in more detail and introduce the phase
function which has its origins in Baker's formula for algebraic curves [6, p 232].

It is convenient
to express the results in terms of the spectral curve ${\cal E}$ of Hill's differential equation. 
Ercolani and
McKean [19] have also considered phases and theta functions for Hill's curve, from a 
different starting
point. In our treatment, we take the $\tau$ function for the linear system as the starting
point, and derive the other functions from there.\par

\indent Consider Hill's equation in the form
$$L\Psi (x,\lambda )=0 , \qquad L={{d}\over{dx}}-
\left[\matrix{0&1\cr u-\lambda&0\cr}\right] \qquad 
(-\infty <x<\infty )\eqno(10.1)$$
\noindent where $u$ is continuous, real-valued and $\pi$-periodic
on ${\bf R}$, and let $U_\lambda (x)$ be the  fundamental solution matrix so that 
$LU_\lambda (x)=0$ and $U_\lambda (0)=I_2$. 
The monodromy operator $S:\Psi (x)\mapsto \Psi (x+\pi )$, which commutes with $L$. Now $H_\lambda =\{ \Psi : L\Psi =0\}$ is a two-dimensional complex vector space, hence the restriction of the monodromy operator to $H_\lambda$ may be represented by the monodromy matrix $U_\lambda (\pi )$, which has determinant one and the discriminant is $\Delta (\lambda )={\hbox{trace}}\, U_\lambda (\pi )$. \par

Let ${\cal F}\psi (z)=\int_{-\infty}^\infty e^{-izx}\psi (x)dx/(2\pi )$ be the complex Fourier transform; note that $\overline{{\cal F}\psi (\bar z)}={\cal F}\psi (z)$ if and only if $\overline{\psi (-x)}=\psi (x)$. Let $PW_\pi$ be the Paley--Wiener space of entire functions $f$ of exponential type such that 
$$\lim\sup_{y\rightarrow\pm\infty} {{\log \vert f(iy)\vert}\over{\vert y\vert}}\leq \pi,\eqno(10.2)$$
\noindent  such that $\Vert f\Vert_{L^2}=( \int_{-\infty}^\infty \vert f(x)\vert^2\, dx)^{1/2}$ is finite; the real linear subspace of $f\in PW_\pi$ such that $\overline {f(\bar z)}=f(z)$ is denoted $RPW_\pi$ . Then ${\cal F}$ restricts to a unitary map $L^2([-\pi ,\pi ]; dx; {\bf C})\rightarrow PW_{\pi}$. We introduce a real Hilbert space
$$V=\{ f\in RPW_\pi: z^2 f(z)\in RPW_\pi \};\qquad  \Vert f\Vert_V^2=\int_{-\infty}^\infty f(x)^2 dx+\int_{-\infty}^\infty x^4 f(x)^2 dx.\eqno(10.3)$$ 
The space ${\bf X}$ is known as the isospectral torus, as every  potential $u$ with prescribed (anti)-preriodic spectrum gives rise to a point in ${\bf X}$. One can in particular associate each $u$ with the $(\mu_j(y))\in {\bf X}$ where $(\mu_j(y))$ is the spectrum associated to the Dirichlet boundary value problem on $[0, \pi ]$ with potential $u(x+y)$. Then $V$ convenient space of functions on $V$ and the real Jacobian of ${\cal E}$ is to be interpreted as a space of bounded linear functionals on $V$. In Proposition 10.1, 
we show how all scattering functions $\phi$ give bounded linear functionals on $V$ and thus produce funcions of the spectral data; this justifies the terminology scattering function. Then we briefly show how these coincide with the functionals which arise from divisors and Abelian integrals on ${\cal E}$, as in [46]. \par 
\vskip.05in

\noindent {\bf Proposition 10.1} {\sl let $\Omega ={\bf C}\setminus \{ \lambda\in {\bf C}: \Delta (\lambda )\in [-2,2]\}$.\par
\indent (i) $\phi (x)=\sum_{k=0}^\infty a_k\cos\sqrt{\xi_k} x$ determines a bounded linear functional on $PW_\pi$ via sampling;
$$f\mapsto {{a_0}\over{2}}+{{1}\over{2}}\sum_{n=1}^\infty a_n \bigl(f(\xi_n)+f(-\xi_n)\bigr).$$
\indent (ii) Then there exists a holomorphic function $\Omega \rightarrow L^2((0, \infty ); {\bf C})$ such that $\psi_\lambda $ is a Floquet solution of Hill's equation with $\psi_\lambda (0)=1$ and multiplier $\rho (\lambda );$\par
\indent (iii) $\psi_\lambda'(0)=-m_+(\lambda )$, where $m_+:\Omega\rightarrow {\bf C}$ is a holomorphic function such that $m_+:{\bf C}_+\rightarrow {\bf C}_+$.\par
\indent (iv) Let $\lambda\in \cup_{j=1}^\infty (\lambda_{2j-1}, \lambda_{2j})$. Then for all $k$ there exist a Floquet solution with eigenvalue $\nu_k\in  [\lambda_{2k-1}, \lambda_{2k}]$ and multiplier $\rho (\nu_k)$ such that $m_+(\lambda )\psi (0)+\psi '(0)=0;$\par
\indent (v) There exists $(a_j)\in \ell^2$ such that}
 $$T(0,y)=\sum_{j=1}^\infty -2a_j\cos\sqrt{\nu_j} y, \qquad 
T(\pi ,y)=\sum_{j=1}^\infty -2\rho (\nu_j)a_j\cos\sqrt{\nu_j} y.\eqno(10.4)$$ 
\vskip.05in
\noindent {\bf Proof.} (i)  By Kadec's $1/4$ theorem, each real sequence $(\xi_n)_{n=0}^\infty $ with  $(\xi^2_n)_{n=0}^\infty \in {\bf X}$ gives a Riesz basis $(e^{\pm i\xi_nx})_{n=0}^\infty$ for $L^2[-\pi, \pi ]$. Suppose  that 
$\phi \in L^2[-\pi ,\pi ]$ is absolutely continuous with $\phi'\in L^2[\pi ,\pi ]$ and $\phi$ is even with $\phi (x)=\phi (-x)$; then $\phi (x)=\sum_{j=0}^\infty a_j\cos\xi_j x$ where $\sum_{j=0}^\infty (1+j^2)\vert a\vert^2<\infty$. Also
$${\cal F}\phi (z)={{1}\over{2}}a_0{\hbox{sinc}}(z)+ {{1}\over{2}}\sum_{n=1}^\infty a_n\bigl( {\hbox{sinc}}(z-\xi_j)+ {\hbox{sinc}}(z+\xi_j)\bigr),\eqno(10.5)$$
\noindent so $\phi$ determines a bounded linear functional on $PW_\pi$ via
$$f\mapsto \int_{-\infty}^\infty f(x)\overline{{\cal F}\phi}(x) dx={{1}\over{2}}a_0f(0)+{{1}\over{2}}\sum_{n=1}^\infty a_n\bigl( f(\xi_n )+f(-\xi_n )\bigr).\eqno(10.6)$$
\indent We show that every scattering function $\phi\in H^1$ determines a bounded linear functional on $PW_\pi$.\par
\indent (ii) By combining Floquet theory with Weyl's results, one shows that there is a holomorphic function ${\cal E}\rightarrow L^2(0, \infty )$ such that $\psi_\lambda $ gives a solution of Hill's equation $-\psi_\lambda''(x)+u(x)\psi_\lambda (x)=\lambda \psi_\lambda (x)$ with $\psi_\lambda (0)=1.$\par
\indent Suppose that $\Delta (\lambda )\in {\bf C}\setminus [-2,2]$; then there exist roots $\rho_{\pm}$ of $\rho^2-\Delta (\lambda )\rho +1=0$ such that $\vert \rho_+\vert <1<\vert \rho_-\vert$. Let 
$$\left[\matrix{1\cr -m_+(\lambda )\cr}\right], \qquad\left[\matrix{1\cr- m_-(\lambda )\cr}\right]\eqno(10.7)$$
\noindent  be an eigenvectors of $U_\lambda (\pi )$ corresponding to $\rho_+$ and  $\rho_-$ respectively. Then 
$$\Psi_+(x;\lambda )=U_\lambda (x) \left[\matrix{1\cr -m_+(\lambda )\cr}\right], \Psi_-(x; \lambda )= U_\lambda (x)\left[\matrix{1\cr -m_+(\lambda )\cr}\right]\eqno(10.8)$$ 
\noindent give solutions of Hill's equation such that $\Psi_+(\pi +x);\lambda )=\rho_+\Psi_+(x;\lambda )$ and $\Psi_-(\pi +x;\lambda )=\rho_-\Psi_-(x;\lambda )$; hence $\Psi_+\in L^2((0, \infty ); {\bf C}^{2\times 1})$; whereas 
  $\Psi_-\in L^2((-\infty ,0); {\bf C}^{2\times 1})$. \par
\indent (iii) We can choose $\lambda \mapsto m_+(\lambda )$ to be holomorphic from the upper half plane  $\{ \lambda\in {\bf C}: \Im \lambda >0\}$ to itself,  so that $\lambda\mapsto  \Psi_+(x);\lambda )$ is holomorphic from the upper half plane into
 $L^2((0, \infty ); {\bf C}^{2\times 1}).$ Since $u$ is real, we can take $\Delta (\bar\lambda )=\overline{\Delta (\lambda )}$, then  $m_+(\bar\lambda )=\overline{m_+(\lambda )}$, $\rho_+(\bar\lambda )=\overline{\rho_+(\lambda )}$  and $\psi_{\bar\lambda }(x)=\overline{\psi_{\lambda}(x)}$. \par 
\indent (iv) By standard results, $(a,b)\subset {\bf R}$ does not intersect the spectrum if and only if $\Im m_+(\lambda )\rightarrow 0$  as $\lambda \rightarrow c$ for all $c\in (a,b)$ and $m_+ :(a,b)\rightarrow {\bf R}$ defines a continuous function. Then we can continue $m_+(\lambda )$ through the spectral gaps to an analytic function by $m_+(\bar\lambda )=\overline{m_+(\lambda )}$. Suppose that $\lambda \in (\lambda_{2k-1}, \lambda_{2k})$, so $\lambda$ belongs to the $k^{th}$ spectral gap. Then there exists a Floquet multiplier $\rho_+\in (-1,1)$ and $U_\lambda (\pi )$ has 
$\rho_+(\lambda )$ as an eigenvalue. Hence the regular Sturm--Liouville problem
$$\eqalignno{ -\psi''+u\psi &=\nu\psi ,\cr
m_+(\lambda )\psi (0)+\psi'(0)&=0,\cr
m_+(\lambda )\psi (\pi )+\psi'(\pi)&=0,&(10.9)\cr}$$
\noindent has eigenvalues $\nu_0<\nu_1<\dots $.  The eigenfunctions $\psi$ are Floquet solutions since $\psi (x+\pi)$ is a real solution that satisfies the same boundary conditions, hence $\psi (x+\pi )=\rho (\nu_j)\psi (x)$ where $\rho(\nu_j )$ is real, so $\nu_j\cup_{\ell =1}^\infty [\lambda_{2\ell -1}, \lambda_{2\ell}]$, and the sequence of monodromy matrices $(U_{\nu_j}(\pi))_{j=0}^\infty $ have a common eigenvector $\left[\matrix{1\cr -m_+(\lambda )\cr}\right].$ The sequences $(\lambda_{2j})$ and $(\nu_j )$ grow at the rame rate, by Borg's estimates,
 so we can choose $\nu_j\in [\lambda_{2j -1}, \lambda_{2j}]$. Thus the boundary value problem () determines a point $(\nu_j)\in {\bf X}$. \par

\rightline{$\square$}\par
\indent Given ${\bf X}$, let ${\cal M}$ be the space of $\pi$ periodic $q\in C^2({\bf R}; {\bf R})$ that have (anti)-periodic spectrum $\{ \lambda_j ;j=0, 1, \dots \}$ and hence share a common discriminant $\Delta$. We regard ${\cal M}$ as a real torus based upon ${\bf X}$, or equivalently as the product ${\bf X}\times \{ (\varepsilon_j)_{j=1}^\infty ;\varepsilon_j =\pm 1\}$ where the signs $\varepsilon_j$ refer to the sign of $\sqrt{\Delta (\mu_j)^2-4}$, and may be chosen independently. The variation of $\lambda_{2j}$ with respect to changes in $q$ was computed by McKean and van Moebecke [45, p. 232] to be 
${{\partial \lambda_{2j}}\over{\partial u}}=\varphi (x;\lambda_{2j})^2$, where $\varphi (x;\lambda_{2j})$ is the periodic eigenfunction for eigenvalue $\lambda_{2j}$, normalized to that $\int_0^\pi  \varphi (x;\lambda_{2j})^2dx=1$. From the product formula
$$4-\Delta (\lambda )^2=c(\lambda -\lambda_0)\prod_{j=1}^\infty \Bigl( 1-{{\lambda}\over {\lambda_{2j-1}}}\Bigr)\Bigl( 1-{{\lambda}\over{\lambda_{2j}}}\Bigr),\eqno(10.17)$$
\noindent we deduce that 
$${{\partial \Delta}\over{\partial q}}(\lambda_{2j}) =-\Delta' (\lambda_{2j}){{\partial \lambda_{2j}}\over{\partial q}}=-\Delta'(\lambda_{2j})\varphi (x;\lambda_{2j})^2.$$
\indent let $f(x;\lambda )$ be the solution of 
$$\eqalignno{-f''(x; \lambda )+u(x)f(x;\lambda )&=\lambda f(x;\lambda )\cr
f(0; \lambda )=1,&\quad f'(0; \lambda )=0.\cr}$$
 \noindent Then the residual spectrum is the sequence $(\nu_j)_{j=0}^\infty$ such that $f'(\pi ;\nu_j)=0;$ we have $\nu_0\leq \lambda_0$ and $\lambda_{2j-1}\leq \nu_j\leq \lambda_{2j}$ for $j=1, 2, \dots $. By manipulating the differential equation, one proves that 
$$\int_0^\pi f(x; \nu_j)^2dx=-f(\pi ; \nu_j){{\partial^2f}\over{\partial x\partial\lambda }}(\pi, \nu_j).$$
\noindent Borg showed that the periodic spectrum, tothether with the residual spectrum and the norming constants together determine $q$. We can write
$$q=-2\sum_{j=1}^\infty {{\lambda_{2j}-\lambda_{0}}\over{ {{\partial^2f}\over{\partial x\partial\lambda }}(\pi, \nu_j)}} {{\partial\Delta}\over{\partial q}}(\lambda_{2j}) +\lambda_0-\sum_{j=1}^\infty (\lambda_{2j}-\lambda_{2j-1}).$$

differentiating along ${{\partial}\over{\partial x_j}}$ corresponds to the vector field
$$X_ju=-\Delta'(\lambda_{2j} ){{\partial }\over{\partial x }} \varphi (x;\lambda_{2j})^2.\eqno(10.18)$$
\noindent Note also that $\varphi (x;\lambda_{2j})^2$ satisfies Drach's equation () .\par

 Given $u\in {\cal M}$,  there is a preferred element $q_0$ such that the tied spectrum is $(\lambda_{2j-1})_{j=1}^\infty$ and the signs are the same as for $u$.\par

\indent We now show how Proposition 10.1 matches with the theory developed in [46] and [8]. We now consider differentials of the first kind, which are of the form 
$$\omega ={{h (\lambda )d\lambda}\over{\sqrt{\Delta (\lambda )^2-4}}}\eqno(10.12)$$
\noindent where $h$ is entire and of order $1/2$. In particular, we can take  $h=\Delta'$ for the special differential $\omega_\infty =\Delta' (\lambda )d\lambda /  \sqrt{\Delta (\lambda )^2-4},$ which satisfies
$$\exp \Bigl( \int_{\lambda_{2j-1}}^\lambda \omega_\infty \Bigr) ={{\rho (\lambda )}\over{\rho (\lambda_{2j-1})}},\eqno(10.13)$$
\noindent where $\rho (\lambda_{2j})=\rho (\lambda_{2j-1})=\pm 1.$  For every $(p_j)\in {\bf X}$, there is a divisor $\delta_{\bf p}=\sum_j \pm (\delta_{p_j}-\delta_{\lambda_{2j-1}})$ of degree zero, and hence a linear functional
$$\delta_{\bf p}:h \mapsto 2\sum_j\int_{\lambda_{2j-1}}^{p_j} {{h (x)dx}\over{\sqrt{\Delta (x)^2-4}}}=\sum_j x_j h (\lambda_{2j}).\eqno(10.14)$$
\noindent The real periods are defined to be the linear functionals 
$$\alpha_j:h\mapsto 2\int_{\lambda_{2j-1}}^{\lambda_{2j}} {{h (x)dx}\over{\sqrt{\Delta (x)^2-4}}},\eqno(10.15)$$
\noindent and the need for initial factor of $2$ arises when we think of a contour encircling the interval $[\lambda_{2j-1}, \lambda_{2j}]$ and the square root taking opposite signs on the top and bottom of  $[\lambda_{2j-1}, \lambda_{2j}]$ .  
\vskip.05in
\noindent {\bf Corollary 10.2} {\sl (i)  The functionals given by $\delta_{\bf p}$ and  $2\psi (x)=2\sum_{j=1}^\infty x_j\cos \sqrt{\lambda_{2j}}x$ coincide on all 
 $h(z)=f(\sqrt{z})+f(-\sqrt{z})$ where $f\in V$.\par
\indent (ii) The coefficient $x_j$ depends upon the variation of $\lambda_{2j}$ with respect to $u$.}\par
\vskip.05in
\noindent {\bf Proof.} (i) Such an $h$ is entire and of order $1/2$ with $\overline{h(\bar z)}=h(z)$ and $\int_0^\infty x^{3/2}h(x)^2dx<\infty$. We have
$$2\sum_{j=1}^\infty \int_{\lambda_{2j-1}}^{p_j} {{(f (\sqrt{x})+f(-\sqrt{x}))dx}\over{\sqrt{\Delta (x)^2-4}}}=\sum_{j=1}^\infty x_j(f(\sqrt{\lambda_{2j}})+f(-\sqrt{\lambda_{2j}})),\eqno(10.16)$$
which is the functional that arises from $2\psi (x)=2\sum_j x_j\cos \sqrt{\lambda_{2j}}x$ as in Propostion 10.1 with $(\nu_j)=(\lambda_{2j})$.\par
\indent (ii) To find the dependence of $x_j$ upon $u$, we recall that 
\rightline{$\square$}\par
\indent Finally, we realize ${\bf X}$ as the real part of the Jacobian of the curve ${\cal E}$. One can introduce a basis of differentials $\omega_k$ such that $\alpha_j(\omega_k)=\delta_{j,k}$ for $j,k=1, \dots;$ then there is a map 
$$(\omega_1, \omega_2, \dots )\mapsto (\delta_{\bf p}(\omega_1), \delta_{\bf p}(\omega_2), \dots )\in {\bf R}^\infty/{\bf Z}^\infty\eqno(10.19)$$
\noindent so we can regard the divisor $\delta_{\bf p}$ as a point in a real vector space modulo the lattice generated by $(\alpha_j(\omega_1, \alpha_j (\omega_2), \dots )$.\par
\vskip.05in

\vskip.05in

\indent Originally, $\tau$ functions were introduced in [29] to
generalize the classical notion of a theta function, as we now
recall. In this section we show that uniformly periodic linear systems
generalize the notion of a theta function on the complex torus
of genus one, essentially the elliptic curve. \par
\vskip.05in

\indent Let $(-A,B,C)$ be as in Theorem 5.3, let $U:{\bf C}^{N-1}\rightarrow {\cal L}(H)$ be holomorphic, such that
$U(t)A=AU(t)$. Let 
$\Sigma (t)=(-A, U(t)B, CU(t); U(t)EU(t))$ for $t\in {\bf C}^{N-1}$ have corresponding tau functions 
$\tau (x,t)$.\par
\vskip.05in
   \noindent {\bf  11 The sinh-Gordon equation}\par
\vskip.05in
\indent We introduce
$$U=\left[\matrix{i\lambda &iq\cr ir&-i\lambda\cr}\right],\quad   
V={{1}\over{i\lambda}}\left[\matrix{ \cosh u&-i\sinh u\cr -i\sinh u&-\cosh
u\cr}\right],\eqno(11.1)$$
\noindent and consider the consistency condition for the pair of differential
equations
$$\eqalignno{ {{d}\over{dx}}\Psi &=U\Psi \cr
{{d}\over{dt}}\Psi &=V\Psi .&(11.2)\cr}$$
\noindent We require
$${{dU}\over{dt}}-{{dV}\over{dx}}+[U,V]=0,\eqno(11.3)$$  
\noindent which with $r=-q=-u_x/2,$ reduces to $u_{xt}=2\sinh u.$\par
\indent We now show how to generate a family of solutions of this PDE. Let
$(-A,B,C_0)$ be the linear system introduced in the proof of Theorem 3.4. Now consider the family of linear systems 
$$(-A, B, C(\mu ,t)), \quad C(\mu ,s)=\mu C_0e^{-sA^{-1}},
\eqno(11.4)$$
\noindent which have scattering function $\phi_{(x)}(t;\mu ,s)=C(\mu ,s )e^{-(t+2x)A}B,$
which satisfies
$${{\partial^2}\over{\partial s\partial x}} \phi_{(x)}(t;\mu ,s)
=2\phi_{(x)}(t;\mu ,s).\eqno(11.5)$$
We suppress $t$ in the notation, and regard $(x,s)$ as the primary variables; 
then we introduce $\tau_1(x,s;\mu )=\det (I+R_x),$ and  $\tau_{0}(x,s;\mu )=\det (I-R_x)$,
where as usual 
$$R_x=\int_x^\infty e^{-tA}BC(\mu ,s)e^{-tA}\, dt.\eqno(11.6)$$
\vskip.05in
\noindent {\bf Proposition 11.4} (sinh-Gordon) {\sl (i) The tau functions satisfy
$${{\partial^2}\over{\partial s\partial x}}\log \tau_1={{1}\over{2}}\Bigl(
{{\tau_1^2}\over{\tau_{0}^2}}-1\Bigr);\eqno(11.7)$$
\indent (ii) and 
$$u=4\tanh^{-1}\Bigl({{\tau_1-\tau_{0}}\over{\tau_1+\tau_{0}}}\Bigr)\eqno(11.8)$$
satisfies}
$$ {{\partial^2 u}\over{\partial s\partial x}}=2\sinh u.\eqno(11.9)$$
\vskip.05in
\noindent {\bf Proof} (i) Let 
$$z= {{\partial^2}\over{\partial s\partial x}}\log\tau_1,
\qquad y={{\tau_1^2}\over{\tau_{0}^2}}, \eqno(11.10)$$
which we regard as functions of $\mu$. We have $z(0)=0$ and $y(0)=1$, and we aim to prove
that $z=(1/2)(y-1)$; to accomplish this, we prove that 
$$\mu {{dz}\over{d\mu }}= w(z+(1/2)),\eqno(11.11)$$
where $w=\mu {{d}\over{d\mu }}\log y$. Let $F=(I+R_x)^{-1}$ and introduce
$$U={{1}\over{4}}\Bigl( (2F-I)^{-1}A^{-1}+A^{-1}(2F-I)^{-1}\Bigr);\eqno(11.12)$$
then one checks that 
$$\lfloor U\rfloor ={\hbox{trace}}\bigl( (I-R^2)^{-1}R\bigr).
\eqno(11.13)$$
Also, $\mu {{d}\over{d\mu}}R =R,$ so one can readily check that
$${{\mu}\over{4}} {{d}\over{d\mu}}\log y={{\mu}\over{2}}{{d}\over{d\mu}}\log {{\det (I+R)}\over{\det
(I-R)}}={\hbox{trace}}\bigl( (I-R^2)^{-1}R\bigr).\eqno(11.14)$$
Also $z=-\lfloor A^{-1}\rfloor$ satisfies
$$\mu {{d}\over{d\mu}} z=\lfloor A^{-1}-FA^{-1}-A^{-1}F\rfloor;
\eqno(11.15)$$
\noindent so the claimed differential equation is
$$\lfloor A^{-1}-FA^{-1}-A^{-1}F\rfloor=-4\lfloor U\rfloor\bigl( \lfloor A^{-1}\rfloor
-1/2\bigr)\eqno(11.16)$$
or more symmetrically
$$\lfloor A^{-1}-FA^{-1}-A^{-1}F\rfloor=-2\bigl( \lfloor U\rfloor\lfloor A^{-1}\rfloor
+\lfloor A^{-1}\rfloor\lfloor U\rfloor-\lfloor U\rfloor\bigr).
\eqno(11.17$$
\noindent One proves this by repeatedly using the trace identity.\par
\indent (ii) Now 
$$\eqalignno{ {{\partial^2 u}\over{4\partial s\partial x}}&={{\partial^2}\over{2\partial s\partial
x}}\bigl( \log \tau_1-\log\tau_{0}\bigr)\cr
&={{1}\over{4}}\Bigl(
{{\tau_1^2}\over{\tau_{0}^2}}-{{\tau_{0}^2}\over{\tau_{1}^2}}
\Bigr)&(11.18)\cr}$$
\noindent while
$$\eqalignno{{{1}\over{2}}\sinh u&={{2\tanh (u/4) (1+\tanh^2(u/4))}\over{ (1-\tanh^2(u/4))^2}}\cr
&={{1}\over{4}}\Bigl(
{{\tau_1^2}\over{\tau_{0}^2}}-{{\tau_{0}^2}\over{\tau_{1}^2}}
\Bigr);&(11.19)\cr}$$    
\noindent hence the sinh-Gordon equation is satisfied.\par 
\vskip.05in

\indent (ii) In [40], Maier uses Lam\'e's equation to produce explicit covering maps 
${\cal Y}\rightarrow {\cal T}$ of the elliptic curve by hyperelliptic curves of arbitrary genus. For generic values of $J$ and $\ell=1,2, \dots $, 
there exists a hyperelliptic curve ${\cal Y}_\ell$ of genus $\ell$ and a 
holomorphic covering
map ${\cal Y}_\ell\rightarrow {\cal T}$ of degree $\ell (\ell +1)/2$; 
in special cases, one can reduce hyperelliptic integrals to elliptic integrals.\par
 \indent Having constructed the potential $\wp$ from a periodic 
linear system, we can produce a family of Hankel kernels and potentials from
 standard limiting arguments which are associated with exactly solvable problems 
in quantum mechanics. Consider an interacting system of $N$ identical particles at
 positions $x_j$ on the real line which interact only pairwise, and where the 
strength of the mutual separation.

. We consider
$$\tau (x)=e^{ax^2+bx+c}{{\prod_{j=1}^n
\vartheta_1(x-a_j)}\over{\prod_{j=1}^m\vartheta_1(x-b_j)}},\eqno(10.23)$$
\noindent which is meromorphic with divisor 
$\delta_{\tau)}=\sum_j\delta_{a_j}-\sum_k\delta_{b_k}$ on some cell of the quotient space ${\bf
C}/\Lambda$ so ${\hbox{deg}}(\tau )=n-m$. If ${\hbox{deg}}(\tau )=0$, $\sum_{j=1}^n(a_j-b_j)\in \Lambda$
and $a=b=0$, then by Abel's theorem $\tau$ is elliptic of the first kind. In all cases, $\tau $ is elliptic of the third kind, and $u(x)=-2(d^2/dx^2)\log\tau$ is elliptic of the first kind with possible poles at the $a_j$, $b_j$  and congruent points with respect to the lattice $\Lambda$. \par
Let $\tau_1(x)=e^{bx+c}\prod_{j=1}^n\vartheta_1(x-a_j)$ and 
$\tau_0(x)=\prod_{j=1}^n \vartheta_1(x-b_j)$, and  further suppose that the 
$a_j$ and $b_j$ all give distinct points in ${\bf C}/\Lambda $; 
next let $v(x)=(d/dx)\log (\tau_1(x)/\tau_0(x))$ and $w(x)=(d/dx)\log (\tau_1(x)\tau_0(x)).$ \par
\vskip.05in
 Any linear system of rational matrix ordinary differential equations 
gives rise to an integrable operator $K$ as in Tracy and Widom's theory of matrix models. Under general conditions on existence of solutions, it is shown that there exist Hankel operators $\Gamma_\Phi$ and $\Gamma_\Psi$  with matrix symbols such that $\det (I+\mu K)=\det (I+\mu \Gamma_\Phi\Gamma_\Psi )$. The paper derives differential equations 
for $\tau$  in terms of the singular points of the differential equation. Any semi classical weight gives rise to a system of ordinary differential equations of Fuchsian type with coefficients that are $2\times 2$ rational matrices. The paper produces a linear system such that the scattering function satisfies this differential equation, and gives the associated tau function. This paper also introduces an admissible linear system with tau function which gives a solution of Painlev\'e's equation $P_{II}$.\par
\indent (ii)  The periodic linear system $\Sigma$ has a tau function $\tau$ and a periodic potential $u$. Hence $\Sigma$ is associated with Hill's discriminant $\Delta (\lambda )$ and a spectral curve ${\cal E}$, which is typically a transcendental hyperelliptic curve of infinite genus. The Jacobi variety ${\bf X}$ of ${\cal E}$ is then an infinite dimensional complex torus. Various meromorphic functions on ${\cal E}$ and ${\bf X}$ can be realised from operator rings on the state space of $\Sigma$. 
Periodic linear systems give periodic potentials as in Hill's equation $-\psi''+u\psi =\lambda\psi$. If $u$ is real $C^2$ and periodic, and Hill's equation has independent Floquet solutions for all but finitely many $\lambda$, then $u$ is finite gap and $\tau$ is the restriction of a theta function to a straight line in the Jacobian of a hyperelliptic curve. The paper deals with the case of elliptic  $u$, so $u$ is expressed as a quotient of tau functions from periodic linear systems. 
The paper produces integrable operators that arise from 
linear $2\times 2$ matrix differential equations with
 elliptic coefficients, extending examples given by Tracy and Widom. 
The paper exhibits the most general algebraic tau functions that can be produced from elliptic linear systems. \par

\noindent In our previous paper [12], we showed that $\vartheta_1(x)$ arises as the tau function of a periodic linear system. In the current article, we develop the analogy between tau functions and the theta functions which are associated with 
Jacobians of algebraic curves.\par

\vskip.05in
\noindent {\bf Theorem 1.3} {\sl Suppose that $\Sigma =(-A,B,C)$ is a linear system with
input and output space $H$, and $(e^{-x A})$ is a
uniformly continuous and $\pi$-periodic group on $H$. Suppose that
there exists $E\in {\cal L}^1$ such that
$AE+EA=BC$, and let ${\bf B}$ be the space of $\pi$-periodic and meromorphic functions with values in
${\cal L}^1$.\par
\indent (i)  Then there
exists a solution to (1.10) and (1.11) such that $\tau (x)=\det (I+R_x)$ is
entire and $\pi$-periodic, and $u$ is meromorphic and $\pi$-periodic.\par
\indent (ii) There exists an algebra $({\bf A}_\Sigma, \ast , \partial )$ containing
 $F=(I+R)^{-1}$ and a ring homomorphism 
$\lfloor \, .\, \rfloor :{\bf A}_\Sigma\rightarrow {\bf B}$ such that} 
$u(x)=-4{\hbox{trace}}\lfloor A\rfloor$.\par
\indent {\sl (iii) Let $\tau_\zeta (x)=\det (I+(\zeta I+A)(\zeta I-A)^{-1}R_x)$. Suppose that 
$\psi_{\pm \zeta} (x)=e^{\pm \zeta x}\tau_{\pm\zeta} (x)/\tau (x)$ give a pair of linearly independent solutions of Hill's equation
$$-{{d^2}\over{dx^2}}\psi_\zeta (x)+u(x)\psi_\zeta (x)=-\zeta^2\psi_\zeta (x)\eqno(1.19)$$
\noindent such that $\zeta\mapsto \psi_\zeta (x)$ and $x\mapsto \psi_\zeta (x)$ are meromorphic for all but finitely many $\zeta , x\in {\bf C}$. Then 
$u$ is a finite gap potential for Hill's equation, and the multiplier curve is hyperelliptic.}\par   
\vskip.05in
Krichever has determined the most general elliptic
solution of $KP$, and described this in terms of the classical Hamiltonian system
$$H={{1}\over{2}}\sum_{j=1}^n p_j^2 -2\sum_{j,k:j<k}\wp (x_j-x_k).\eqno(1.22)$$
\noindent We realise this general elliptic solution in terms of elliptic
linear systems.\par

\indent $(\tau Q)$ The general solution of Schr\"odinger's equation is given by a quotient of translations of $\tau$ functions for all but finitely many $\lambda\in {\bf C}$.\par
\indent $(Picard)$ The general solution to Schr\"odinger's equation is meromorphic on ${\bf C}$ for all but finitely many $\lambda\in {\bf C}$. \par
\indent $(FG)$ The Bloch spectrum of Schr\"odinger's equation has only finitely many gaps. \par
\indent $(KP)$ $\tau$ satisfies the Kadomtsev--Petviashvilli hierarchy of differential equations.\par
\noindent For $u$ periodic, real and $C^2$, $(\tau Q)\Rightarrow (Picard) \Rightarrow (FG)\Rightarrow (KP)$. \par
\noindent For $u$ elliptic, $(\tau Q)\Leftrightarrow (Picard) \Leftrightarrow (FG)\Rightarrow (KP)$, and all the elliptic solutions of $KP$ are determined for elliptic linear systems.\par
A finite gap potential is associated with a hyperelliptic spectral curve
${\cal E}=\Bigl\{ (\mu  , \lambda )\in {\bf C}^2: \mu^2=\prod_{j=0}^{2g} (\lambda
-\lambda^o_j)\Bigr\} \cup \bigl\{ (\infty , \infty )\bigr\}\eqno(10.21)$$
\noindent {\sl with $\lambda_j^o$ real and such that the $\lambda\in \sigma_B$  give real
points $(\mu  , \lambda )$ on ${\cal E}$. }\par

\indent The notion of a tau function of a linear system generalizes the classical concept of a theta function for an abelian variety.
 Let $\Lambda$ be a discrete and countable subgroup of ${\bf C}^g$ such that 
${\bf C}^g/\Lambda $ is compact, and hence gives a complex torus. Let
 $L: {\bf C}^g\times \Lambda \rightarrow {\bf C}$ be a function such that $x\mapsto
L(x,\gamma )$ is linear for all $\gamma\in \Lambda$ and $J:\Lambda \rightarrow
{\bf C}$ be any function; then $(L, J)$ is a type.
Let $\theta$ be a quotient of entire functions on ${\bf C}^g$ 
such that $\theta (z+\gamma )=e^{2\pi i(L(z,\gamma )+J(\gamma ))}\theta (z)$ for all $z\in {\bf C}^g$; then $\theta$ is a theta function of 
type $(L, J)$. An abelian
function is a meromorphic and periodic function on ${\bf C}^g$ with respect to $\Lambda$, or equivalently a meromorphic theta function of type $(0,0)$.
 A fundamental result is that any 
abelian function can be expressed as the quotient of entire theta functions of the same 
type.\par 

\vskip.05in  
\noindent {\bf Definition} (Complex differential rings and state rings) Let 
$\Omega$ be a domain in ${\bf C}$ such that $z+w\in \Omega$ whenever $z,w\in\Omega$. Let $X$ be a complex 
Banach algebra  and let 
${\cal M}_\Omega (X)$ be the meromorphic functions from $\Omega$ to $X$. 
Suppose that ${\bf S}$ as above is also a differential subring of ${\cal M}_\Omega (X)$, for
the standard complex derivative $d/dx$. Then we say
that ${\bf S}$ is a complex state ring for $(-A,B,C)$ on
$\Omega$. \par
\vskip.05in

\noindent {\bf Definition} (i) (Scattering field) Let ${\bf K}_0$ be a complex differential field on
$\Omega$, and let ${\bf K}_1$ be the smallest complex differential
field containing ${\bf K}_0$ that also contains $\phi (z)=Ce^{-zA}B$.\par
\indent (ii) Let ${\bf K}_2$ be the smallest complex differential field that contains ${\bf K}_0$ and the
potential $u$.\par
\indent (ii) (Abelian extension) Let ${\bf F}_j$ be a field of complex functions with
differential $\partial$, and adjoin a complex function $h$ to 
${\bf F}_j$ to form ${\bf F}_{j+1}$ where
either:\par
\indent (A1) $h=\int g$ for some $g\in {\bf F}_j$, so that $\partial
h=g$;\par
\indent (A2) $h=\exp \int g$ for some $g\in {\bf F}_j$;\par
\indent (A3) $h$ is algebraic over ${\bf F}_j$; \par
\indent (A4) or adjoin functions $h_j$ for $j=0, \dots ,g$ by substituting in Abelian
functions. Let $f_1, \dots , f_g$ be transcendental over ${\bf C}$
and $f_0$ algebraic over ${\bf C}(f_1, \dots , f_g)$ so that ${\bf
C}(f_0, \dots ,f_g)$ is an Abelian function field on some torus ${\bf
C}^g/\Lambda$. Then for some $u\in {\bf F}_j$, let $h_j(z)=f_j(u(z))$ and
form ${\bf F}_{j+1}={\bf F}_j(h_0, \dots ,h_g)$.\par
\noindent Then ${\bf F}_{j+1}$ is an Abelian extension of ${\bf F}_j$ as in [65,66].
More generally, a field ${\bf L}$ is a Abelian extension of ${\bf
K}$ if there exist differential fields ${\bf F}_j$ such that ${\bf
K}={\bf F}_0\subset {\bf F}_1\subset \dots \subset {\bf F}_n={\bf L}$,
and each ${\bf F}_j$ arises from ${\bf F}_{j-1}$ by applying (A1), (A2),
(A3) or (A4).\par

\vskip.05in

\noindent {\bf Corollary 4.4} {\sl Let $(-A,B,C)$ be a $(2,2)$ admissible linear system
as in Theorem 2.2, and suppose furthermore that $A$ is bounded and
$H_0={\bf C}$.\par
\indent (i) Then $(-A,B,C)$ has a complex state ring ${\bf S}$
on ${\bf C}$ on which
$R_z$ is unique, and $R_{z+t}\in {\bf S}$ for all $t\in {\bf C}$.\par
\indent (ii) The scattering function is 
$\phi (2x)=\bigl\lfloor F_x^{-2}\bigr\rfloor$ and the potential
$u(x)=-4\bigl\lfloor A\bigr\rfloor$ has Miura transform $v=2\lfloor
(2F-1)^{-1}\rfloor$.\par
\indent (iii) There exists an open sector $S$ containing $(x_0, \infty )$ such that the range ${\bf A}_u=\lfloor{\bf A}_\Sigma\rfloor $ is a differential rings, and every elements $w$ of
 ${\bf A}_u$ is a bounded and holomorphic functions on  $\Omega$ such that $w(x)\rightarrow 0$ as $x\rightarrow\infty$ .\par
\indent (iv) The field of
fractions ${\bf K}$ of ${\bf A}_u$ is a differential
field, and $\tau (x)=1/\det F_x$ is entire and belongs to a Abel extension ${\bf
L}$ of ${\bf K}$.}\par

\vskip.05in
\noindent {\bf Proof.} (i) Mainly this follows from Lemma
2.1 and Proposition 2.3. Observe that $x\mapsto R_x$ is entire, and there exists $x_0$ such that $\Vert
R_x\Vert_{{\cal L}^1}<1$ for all $x\in {\bf C}$ such that $\Re x>x_0$, so we can easily carry out calculations on
$\{ x: \Re x>x_0\}$ and then extend identities by analytic continuation. Every point in the spectrum of $R_x$ other than $\lambda =0$ gives a pole of $(\lambda I-R_x)^{-1}$, so we can define
$$F_x={{1}\over{2\pi i}}\int_{\Gamma} (\lambda I-R_x)^{-1}{{d\lambda}\over{1+\lambda }}\eqno(4.19)$$
\noindent for a suitable contour $\Gamma$ that surrounds the spectrum of $R_x$. 
 By Riesz's theory of compact operators, the
$F_x=(I+R_x)^{-1}$ defines a meromorphic operator valued function on
${\bf C}$.  Hence we can select ${\bf S}$ to be the subring of 
meromorphic functions from $\Omega$
to ${\cal L}(H)$ generated by $I, A, BC, 
R_x, e^{-xA}$ and $F_x$. On $\{ x: R_x+R_x^\dagger >-2I\}$, the
function $F_x$ is holomorphic and satisfies $F_x'=FA+AF-2FAF$. By Lemma 2.2, there exist $x_0>0$ and $\varepsilon, \delta >0$ such that $\Vert R_z\Vert <1-\delta$ on 
the sector $S=\{ z\in {\bf C}: -\varepsilon <\arg (z-x_0)<\varepsilon\}$, so $F$ is bounded and holomorphic on $S$. Hence all the elements of $\lfloor {\bf A}_\Sigma\rfloor$ are bounded and holomorphic on $\Omega$.\par 
\indent (ii)  Evidently 
$$\bigl\lfloor F^{-2}\bigr\rfloor =Ce^{-2xA}B=\phi (2x),\eqno(4.20)$$
while we can write (1.12) as $u(x)=-2(d/dx)^2\log\det (I+R_x)$ where  $(d/dx)\log\det (I+R_x)=\bigl\lfloor
F^{-1}\bigr\rfloor$ and differentiate using (2.2). Also, we have
$$2\partial (2F-I)^{-1}=-4A -4(2F-1)^{-1}\ast (2F-I)^{-1},\eqno(4.21)$$
\noindent hence $v= 2\lfloor (2F-I)^{-1}\rfloor $ satisfies $v'+v^2=u$.\par
\indent (iii) By Lemma 3.2(iii), the image of $\lfloor\, .\,\rfloor$ is a
commutative differential ring. One can introduce $x_0\geq 0$ such that 
$\Vert R_{x_0}\Vert <1$, then use the sector $\Omega=\{ x=x_0+z \in {\bf C}:
 \Re z\langle Af,f\rangle \geq 0,  \forall f\in {\cal D}(A)\}$; observe that $z+t\in \Omega$ for all $z\in\Omega$ and $t>0$. Then that $F_x$ is
 holomorphic on $\Omega$ and $F_x\rightarrow I$ as $x\rightarrow\infty $ along $(x_0, \infty )$  in $\Omega$, and $\Vert e^{-xA}\Vert_{{\cal L}(H)}\rightarrow 0$; hence
$w(x)\rightarrow 0$ as $x\rightarrow\infty$ for all $w\in {\bf A}_u$.\par 
\indent (iv) Note also that  ${\bf A}_u$ is a ring of holomorphic functions on a connected set $\Omega$ and  hence ${\bf C}1+{\bf A}_u$ is an integral domain with a field of
fractions ${\bf K}$, which are meromorphic functions on $\Omega$. We have $2(d/dx)^2\log \det F_x=u(x)\in {\bf K}$, 
so we can recover $\det F_x$ by integration and exponential
integration. By (2.3) and Morera's theorem, $R_x$ is an entire
${\cal L}^1$-valued function, hence $\det (I+R_x)$ is entire.\par 

\vskip.05in
We write $f'=df/dx$., so for all $\lambda\in {\bf C}$ such that $\Im \lambda>0$, there exist nontrivial solutions 
$\psi_{\pm}(x; \lambda )$ to $-\psi_{\pm}''(x; \lambda )+u(x)\psi_{\pm} (x;\lambda )=\zeta \psi_{\pm} (x;\lambda )$ such that 
$\psi_+(x; \lambda )\in L^2(0, \infty )$ and  $\psi_-(x; \lambda )\in L^2(-\infty ,0)$. Then we let 
$$m_+(x; \lambda )={{\psi'_+(x; \lambda )}\over{\psi_+(x; \lambda)}}, \quad m_-(x; \lambda )={{\psi'_-(x; \lambda )}\over{\psi_-(x; \lambda )}}.\eqno(5.24)$$
\noinden This $m_+(x;\lambda )$ is the Weyl $m$-function, although the definition is usually written differently; our sign convention is such that $\lambda\mapsto m_+(x; \lambda )$,  $\lambda\mapsto -m_-(x; \lambda )$ and hence $\lambda\mapsto g(x;-\lambda )$ take the upper half plane to itself. \par
\indent These definitions give at once the formulas
$$m_+(x; \lambda )-m_-(x; \lambda )={{-1}\over {g(x; -\lambda )}},\eqno(5.26)$$
$$m_+(x; \lambda  )+m_-(x; \lambda )={{g'(x; -\lambda)}\over {g(x; -\lambda )}};\eqno(5.27)$$
\noindent hence $g$ and $g'$ together determine $m_+$ and $m_-$, and conversely. Also, we have the Riccati equations
$$m_+'={{\psi''_+}\over{\psi_+}}-\Bigl(  {{\psi'_+}\over{\psi_+}}\Bigr)^2=u-\lambda -m_+^2,\eqno(5.28)$$
\noindent and likewise for $m_-$. The Riccati equation reduces to Drach's equation (5.3) with $-\lambda$ in place of $\zeta$. One can check that $-i\sqrt{\lambda}g(x; -\lambda )\rightarrow 1/2$ as $\lambda\rightarrow\infty$.\par
\indent (ii) For $x_0\in {\bf R}$, let $L_{x_0}$ be $-{{d^2}\over{dx^2}}+u(x)$ with the Dirichlet boundary condition $f(x_0)=0.$. Then $L_{x_0}$ may be regarded as a perturbation of $L$, and there is a spectral shift formula
$${\hbox{trace}}(e^{-tL}-e^{-tL_{x_0}})=\int_0^\infty te^{-t\lambda }\xi (x_0, \lambda )\, d\lambda \qquad (t>0).\eqno(5.29)$$
\noindent  Gesztesy and Simon show that the xi function and diagonal Greens function mutally determine one another through the formula
$$\xi (x_0, \lambda )=\lim_{\varepsilon\rightarrow 0+} {{1}\over{\pi}}{\hbox{arg}} G(x_0, x_0; \lambda +i\varepsilon )\qquad (\lambda, x_0\in {\bf R}).\eqno(5.30)$$ 
\indent (ii) Suppose that $u$ is  real Schwartz class function $G(x,y;\zeta )$ be the Greens kernel of $(-\zeta I+L)^{-1}$. G)\asymp {{i}\over{2\sqrt{\zeta}}}\Bigl( 1+\sum_{j=1}^\infty {{\hat f_j}\over{\zeta^j}}\Bigr),\eqno(5.31)$$
\noindent uniformly for $x$ is bounded subsets of ${\bf R}$ as $\zeta\rightarrow i\infty$. Their proof uses the recursion formula (5.1) for the $\hat f_j$.\par
 \indent Suppose that $u\in C^2({\bf R}; {\bf R})$ is $\pi$-periodic with $\int_0^\pi u(x)dx=0$. Suppose that $\lambda_0=0$ and $\vert \lambda_{2j-1}-j^2\vert<1/(4j)$ and
 $\vert \lambda_{2j}-j^2\vert<1/(4j)$ . Let ${\bf X}=\{ 0\}\times\prod_{j=1}^\infty [\lambda_{2j-1},\lambda_{2j}]$.\par

\indent Let $\mu_j$ be the tied spectrum of Hill's equation, namely the spectrum of Hill's equation for the boundary conditions $f(0)=0=f(\pi )$. The tied eigenvalues interlace the periodic eigenvalues, so $\lambda_{2j-1}\leq \mu_j\leq \lambda_{2j}$.
Given ${\bf X}$, let ${\cal M}$ be the space of $\pi$ periodic $q\in C^2({\bf R}; {\bf R})$ that have (anti)-periodic spectrum $\{ \lambda_j ;j=0, 1, \dots \}$ and hence share a common discriminant $\Delta$. We regard ${\cal M}$ as a real torus based upon ${\bf X}$, or equivalently as the product 
$${\cal M}= \prod_{j=1}^\infty \{ (\mu_j, \sqrt{\Delta (\mu_j)^2-4}): \lambda_{2j-1}\leq \mu_j\leq \lambda_{2j}\}$$
\noindent of real tori based upon the intervals of instability of Hill's equation. We let $y_2(x, \lambda )$ be the solution of Hill's equation with initial conditions $y_2(0)=0, y_2'(0)=1$, so $y_2(\pi ; \lambda )=0$ if and only of $\lambda$ belongs to the tied spectum. Indeed one can write $y_2(\pi, \lambda )=\prod_{j=1}^\infty (\mu_j-\lambda )/(j^2\pi^2)$. We replace $u(x)$ by $u(x+t)$, then form the corresponding tied spectrum $\mu_j(t)$. By a classical trace formula of Gelfand and Levitan
$$u(t)-\int_0^\pi u(x) {{dx}\over{\pi}}=\lambda_0+\sum_{j=1}^\infty (\lambda_{2j-1}+\lambda_{2j}-2\mu_j(t)\bigr).$$
\noindent We let $p_j=p_j(t)=(\mu_j(t), \sqrt{\Delta (\mu_j(t))^2-4})$, and introduce the dynamical system
$${{d\mu_j(t)}\over{dt}}={{\sqrt{\Delta (\mu_j(t))^2-4}}\over{{{\partial y_2(\pi;\mu_j(t))}\over{\partial \mu}}}}.$$
\noindent Then McKean and Trubowitz show that the resulting solution $(\mu_j(t))_{j=1}^\infty )$ indeed gives the tied spectrum of $u(x+t)$. We can define
$${{\vartheta'(t)}\over{\vartheta (t)}}+\sum_{j=1}^\infty \int_{(\lambda_{2j-1}, 0)}^{p_j(t)} {{\Delta'(\mu )\, d\mu}\over{ \sqrt{\Delta (\mu)^2-4}}}=0.$$
\noindent Its and Matveev show that $u_0(t)=-{{d^2}\over{dt^2}}2\log\vartheta (t)$ gives a potential that has periodic spectrum $(\lambda_j)$ and tied spectrum $(\mu_j)$. By uniqueness, we deduce that 
$\tau (t)=\vartheta (t)e^{at+b}$ for some constants $a, b\in {\bf C}$. We now introduce the phase constant
$$\zeta (p)= \int_ {(\lambda_{0}, 0)}^{p} {{\Delta'(\mu )\, d\mu}\over{ \sqrt{\Delta (\mu)^2-4}}}.$$
\noindent A complete interval of instability contributes 
$$ \int_ {(\lambda_{2j-1}, 0)}^{(\lambda_{2j}, \varepsilon)} {{\Delta'(\mu )\, d\mu}\over{\sqrt{\Delta (\mu)^2-4}}}=0$$
\noindent since $\Delta (\lambda_{2j-1})=\Delta (\lambda_{2j})=\pm 2$. Hence $\zeta (p)$ is purely imaginary when $p$ lies in an interval of stability. Now we introduce
$$\psi (x, \zeta )=\exp \bigl( x\zeta ){{\tau_\zeta (x)}\over{\tau (x)}},$$
\noindent where $\tau_\zeta$ is defined for $(-A, (\zeta I+A)(\zeta I-A)^{-1}B, C).$ This gives the Baker--Akhiezer function.\par
\vskip.05in

\noindent {\bf 7 The Korteweg -de Vries equations}\par
 \vskip.05in
\indent In this section, we obtain solutions to the $KdV$ equations.\par
\vskip.05in
\indent  Let $(-A,B,C)$ as in Proposition 2.2, and suppose further that $(e^{-tA^j})_{t>0}$ defines a strongly continuous semigroup on $H$ for $j=1,3, \dots$. For any finite supported sequence $(t_j)$ with $t_j\geq 0$ let $U(t)=\exp (-{\sum_{j=1}^\infty t_{2j-1}A^{2j-1}})$, and let 
$${\bf \Sigma}=\{ \Sigma (t)=(-A, U(t)B, CU(t)): t\in({\bf T}^\infty)_0\}\eqno(7.1)$$
 with tau function $\tau  (x;t)$.\par  
\vskip.05in

\noindent {\bf Lemma 7.1} {\sl (i) There is a deformation given by the semigroup action of ${\bf U}$ on ${\bf \Sigma}$ by  $U(t): \Sigma (s)\mapsto \Sigma (s+t)$, so that $\tau (x;s)\mapsto \tau (x;s+t)$, and the the corresponding
homomorphism of  differential rings $({\bf A}, \ast ,\partial )$ is determined by $R_x\mapsto U(t)R_xU(t)$.}\par
\indent {\sl (ii)  Then $\Sigma (t)$ satisfies the
hypotheses of Propsition 2.2, and gives rise to a differential ring 
${\bf A}(t)$ for the derivatives ${{\partial}\over{\partial x}}$ and
${{\partial }\over{\partial t_j}}$ such that} 
$${{\partial F}\over{\partial t_j}}=A^jF+FA^j-2FA^jF\qquad (j=1,2 \dots
).\eqno(7.2)$$
\noindent {\sl (iii) The corresponding scattering functions $\phi (x;t_1, t_2,
\dots )$
satisfy ${{\partial\phi}\over{\partial t_j}}=2(-1)^{j-1}
{{\partial^j\phi}\over{\partial x^j}}.$}\par
\vskip.05in
\noindent {\bf Proof.} (i) By Theorem 2.2,  the linear system $\Sigma (s)$  satisfies the hypotheses of Lemma 4.1. The remaining statements are clear from the definitions since $d/dx$ commutes with the action of $U(t)$.\par
\indent (ii)  We repeat the proof of Proposition 4.2 with $A^j$ in place of $A$, so that 
$$P\ast_j Q=P(A^jF+FA^j-2FA^jF)Q\eqno(7.3)$$
$$\partial_jP=A^j(I-2F)P+{{\partial P}\over
{\partial t_j}}+P(I-2F)A^j, \quad \lfloor P\rfloor 
=CU(t)e^{-xA}FPFU(t)e^{-xA}B\eqno(6.23)$$
\noindent satisfy $\partial_j(P\ast_jQ)=(\partial_jP)\ast_jQ +P\ast_j (\partial_jQ).$ We also have 
$${{\partial}\over{\partial t_j}}\lfloor P\rfloor =\lfloor \partial_jP\rfloor.\eqno(7.4)$$ 
\par
\indent (iii) This is a simple computation.\par
\rightline{${\square}$}\par

\vskip.1in

\noindent {\bf Proposition 7.2} {\sl   (i) Then $\Sigma (t)$ has corresponding potential
$$u(x,t)=-2{{\partial^2}\over{\partial x^2}}\log\det (I+U(t)R_xU(t))\eqno(7.5)$$
\noindent satisfies the KdV equation}
$$4{{\partial u}\over{\partial t_3}}={{\partial^3u}\over{\partial
x^3}}-6u{{\partial u}\over{\partial x}};\eqno(7.6)$$
\indent {\sl (ii) generally, $u(x,t)$ satisfies the higher order $KdV$ equation}
$$(-1)^n{{\partial u}\over{\partial t_{2n+1}}}
={{\partial }\over{\partial x}}2f_{n+1}(u).\eqno(7.7)$$  
\vskip.05in
\noindent {\bf Proof.} (i) Let $v(x,t_3)$ be 
$$v(x,t)=-2Ce^{-xA-t_3A^3}(I+R)^{-1}e^{-xA-t_3A^3}B;\eqno(7.8)$$
\noindent and let $u(x,t_3)=- {{\partial v}\over{\partial x}}.$ By Lemma 6.1, we have the
following table of derivatives
$$\eqalignno{{{\partial v}\over{\partial t_3}}&=4\lfloor A^3\rfloor;\cr
{{\partial v}\over{\partial x}}&=4\lfloor A\rfloor;\cr
{{\partial^2 v}\over{\partial x^2}}&=8\lfloor A(I-2F)A\rfloor;\cr
{{\partial^3 v}\over{\partial
x^3}}&
=16\lfloor A^3\rfloor-48\lfloor A\rfloor\lfloor A\rfloor;&(7.9)\cr}$$
Hence
$$4{{\partial u}\over{\partial t_3}}=-4{{\partial v}\over{\partial x}}=-16{{\partial
}\over{\partial x}}\lfloor A^3\rfloor\eqno(7.10)$$
\noindent so 
$$4{{\partial u}\over{\partial t_3}}={{\partial^3u}\over{\partial x^3}}-3{{\partial
}\over{\partial x}} u^2.\eqno(7.11)$$ 
\indent (iii) We have
$${{\partial u}\over{\partial t_{2n+1}}}=-4{{\partial }\over{\partial t_{2n+1}}}\lfloor
A\rfloor =-4\lfloor A^{2n+1}(I-2F)A+A(I-2F)A^{2n+1}\rfloor,\eqno(7.12)$$
\noindent and 
$$\eqalignno{{{\partial }\over{\partial x}}2f_{n+1}(u)&
={{\partial }\over{\partial x}}4(-1)^n \lfloor
A^{2n+1}\rfloor \cr
&=4(-1)^n \lfloor A(I-2F)A^{2n+1}+A^{2n+1}(I-2F)A\rfloor,&(7.13)\cr}$$
\noindent hence the result.\par
\rightline{$\square$}\par
\vskip.05in

 Meanwhile,  let $u(x,t)$ be a solution of $KdV$ and consider the pair of simultaineous equations
$$\eqalignno{ -{{\partial^2\psi}\over{\partial x^2}}+u\psi &=\lambda\psi\cr
{{\partial \psi}\over{\partial t}}&={{\partial^3\psi}\over{\partial x^3}}+{{3u}\over{2}}{{\partial \psi}\over{\partial x}}+{{3\psi}\over{4}}.&(7.14)\cr}$$
\vskip.05in
\noindent {\bf Corollary 7.3} {\sl For $\lambda=\lambda_k$, let $\psi_k$ be a solution of (), as given by Theorem 6.1, and let
$$u_N=u-2{{\partial^2}\over{\partial x^2}}{\hbox{Wr}}(\psi_1, \dots , \psi_N).\eqno(7.15)$$
\noindent Then   
$$\Psi_N= {{{\hbox{Wr}}(\psi_1, \dots , \psi_N, \psi )}\over{{\hbox{Wr}}(\psi_1, \dots , \psi_N)}}\eqno(7.16)$$
\noindent gives a solution of () with $u_N$ instead of $u$.}\par
\vskip.05in
\noindent {\bf Proof.} This result is stated in section 3 of [31.]. The contribution of the current paper is matching together all the components together in terms of a common family of linear systems.\par
\rightline{$\square$}\par
\vskip.05in

This extends the classical notion of a theta function, according to the following table.
\vskip.05in
$$\matrix{ {\hbox{hyperelliptic curve}}&{\cal E}& {\hbox{spectral curve}}&\{ (z, \lambda ): z^2=4-\Delta (\lambda )^2\}\cr
{\hbox{Picard variety}}& {\bf X}={\bf C}^g/\Lambda& {\hbox{infinite torus}}& {\bf X}\cr
{\hbox{Riemann theta function}}& \vartheta& {\hbox{tau function}}& \tau\cr
{\hbox{Weierstrass function}}&\wp=D_{\bf X}^2 \log\vartheta & {\hbox{potential}}& u=-2(d^2/dx^2) \log\tau\cr
{\hbox{Fay identities}}&{}&\tau\quad {\hbox{ identities}}\cr
{\hbox{projective embedding}}& \{ (\wp, \dots ,\wp^{(2g-1)})\} & {\hbox{Jacobian }}& \{ (A, \partial A, \partial^2A, \dots )\}\cr}$$ 
\vskip.05in
\noindent {\bf Proof.}  We have already seen that $\psi_2(x \mid {\bf Z}2\omega_1 +{\bf Z}2\omega_2)$
arises from a periodic linear system. In this case, the identity () follows by
comparing poles and residues of each side.\par
\indent The rational, trigonometric and hyperbolic cases in this table satisfy the additional identity $u(x)=\psi(x)^2+\gamma (x)+c$ with $c$ constant, and are derived from the quantum Lax equation in [61, p.39]. However, we can derive all of these as limiting cases from the elliptic potential in the final line of the table. We write $\Lambda ={\bf Z}2\omega_1+{\bf Z}2\omega_2$ where $\omega_1, \omega_2/i>0$. Then we have the thermodynamic limit 
$$2\wp  (x\mid {\bf Z}2\omega_1 +{\bf Z}2\omega_2)\rightarrow 2(\pi/2\omega_2)^2
{\hbox{cosech}}^2(\pi x/2\omega_2) -\pi^2/6\omega_2^2\qquad 
(\omega_1\rightarrow\infty)\eqno(10.34)$$
\noindent and in contrast the high density limit 
$$2\wp  (x\mid {\bf Z}2\omega_1 +{\bf Z}2\omega_2)\rightarrow 2(\pi/2\omega_1)^2
{\hbox{cosec}}^2(\pi x/2\omega_1) -\pi^2/6\omega_1^2\qquad 
(\omega_2/i\rightarrow\infty );\eqno(10.35)$$
\noindent  in the latter case, ${\hbox{q}}\rightarrow 0$ and hence
 $\vartheta_1(x-\alpha \mid \Lambda )/\vartheta_1(x\mid \Lambda )\rightarrow (\sin \pi (x-\alpha ))/\sin \pi x,$ and one can check that the addition rule for $\psi_2$ reduces to the addition rule for the trigonometric function $\psi (x)=2\cot\pi x$. In each case, functions on ${\cal T}$ reduce to functions on the cylinder. When one limit is applied after the other, we have the limiting potential $u(x)=2/x^2$,  and the addition rule gives the formula for Carleman's 
operator. See [54].\par
\indent (iv) One can reverse the thermodynamic limit by taking a potential
 $u$ to a periodic potential with period $2\omega_1$, as in
 $u(x)\mapsto \sum_{n=-\infty }^\infty u(x+2\omega_1n)$. 
This observation is the basis for the scattering theory in [19, 33].\par

\rightline{$\square$}\par

\noindent By differentiating and then integrating by parts, one can easily check that
$$\eqalignno{{{\partial\psi}\over{\partial t}}(x;y,t) +{{\partial^3\psi}\over{\partial x^3}}(x;y,t)&={{3}\over{2}}\zeta u(x;y,t) e^{\zeta x-\zeta^2y/\beta -\zeta^3t}-3\Bigl( {{d}\over{dx}}{{\partial K}\over{\partial x}}(x,x;y,t)\Bigr) e^{\zeta x-\zeta^2y/\beta -\zeta^3t}
\cr&\quad+\int_x^\infty e^{\zeta z-\zeta^2y/\beta -\zeta^3t}\Bigl(  {{\partial^3K}\over{\partial x^3}}(x,z;y,t) +{{\partial^3K}\over{\partial z^3}}(x,z;y,t) \cr
&\quad+{{\partial K}\over{\partial t}}(x,z;y,t)\Bigr)\, dz.&(9.14)\cr}$$
\noindent Our next step is to obtain an expression for the integrand which is nonlinear, but involves only lower order derivatives.\par
\indent  Returning to the setting of Proposition 8.1, by differentiating the Gelfand--Levitan equation, we obtain
$${{\partial^3\phi}\over{\partial z^3}}(x,z;y,t)+{{\partial^3K}\over{\partial z^3}}(x,z;y,t)+\int_x^\infty K(x,w;y,t){{\partial^3\phi}\over{\partial z^3}}(w,z;y,t)dw=0\eqno(9.15)$$
and  

$$\eqalignno{{{\partial\phi}\over{\partial t}}(x,z;y,t)+{{\partial K}\over{\partial t}}(x,z;y,t)&+\int_x^\infty K(x,w;y,t){{\partial \phi}\over{\partial t}}(w,z;y,t)dw\cr
&\quad +\int_x^\infty {{\partial K}\over{\partial t}}(x,w;y,t) \phi (w,z;y,t)dw=0;&(9.16)\cr}$$
\noindent then from the scattering equations for $\phi$, we deduce that 

$$\eqalignno{ -{{\partial^3\phi }\over{\partial x^3}}(x,z;y,t) &+{{\partial^3 K}\over{\partial z^3}}(x,z;y,t) +{{\partial K}\over{\partial t}}(x,z;y,t)+\int_x^\infty {{\partial K}\over{\partial t}}(x,w;y,t)\phi (w,z;y,t)\, dw\cr
&-\int_x^\infty K(x,w;y,t){{\partial^3\phi}\over{\partial w^3}}(w,z;y,t)\, dw=0;&(9.17)\cr}$$
\noindent then by integrating repeatedly by parts, we obtain
$$\eqalignno{   -{{\partial^3\phi }\over{\partial x^3}}(x,z;y,t) &+{{\partial^3 K}\over{\partial z^3}}(x,z;y,t) +{{\partial K}\over{\partial t}}(x,z;y,t)+\int_x^\infty {{\partial K}\over{\partial t}}(x,w;y,t)\phi (w,z;y,t)\, dw\cr
&+K(x,x;y,t){{\partial^2\phi}\over{\partial x^2}}(x,z;y,t)-{{\partial K}\over{\partial z}}(x,x;y,t) {{\partial \phi}\over{\partial x}}(x,z;y,t) \cr
&+{{\partial^2K}\over{\partial z^2}} (x,x;y,t) \phi (x,z;y,t) +\int_x^\infty {{\partial^3K}\over{\partial w^3}}(x,w;y,t) \phi (w,z;y,t)dw\cr
& =0.&(9.18)\cr}$$

\noindent We also differentiate the Gelfand--Levitan equation three times with respect to $x$, obtaining
 $$\eqalignno{ {}&{{\partial^3\phi}\over{\partial x^3}}(x,z;y,t) +{{\partial^3K}\over{\partial x^3}} (x,z;y,t) -\Bigl({{d^2}\over{dx^2}} K(x,x;y,t)\Bigr)\phi (x,z;y,t) -2{{d K}\over{dx}}(x,x;y,t) {{\partial\phi}\over{\partial x}}(x,z;y,t) \cr
& \quad -K(x,x;y,t){{\partial^2\phi }\over{\partial x^2}}(x,y;z,t) -\Bigl({{d}\over{dx}} {{\partial K}\over{\partial x}}(x,x;y,t)\Bigr) \phi (x,z;y,t) -{{\partial K}\over{\partial x}}(x,x;y,t) {{\partial\phi}\over{\partial x}}(x,z;y,t)\cr
& \quad -{{\partial^2K}\over{\partial x^2}}(x,x;y,t)\phi (x,z;y,t)+\int_x^\infty {{\partial^3K}\over{\partial x^3}}(x,w;y,t)\phi (w,z;y,t)\, dw\cr
&=0;&(9.19)\cr}$$
  \noindent so by adding this to (9.18), we obtain

$$\eqalignno{ \Bigl( {{\partial^3K}\over{\partial x^3}}&(x,z;y,t) +{{\partial^3K}\over{\partial z^3}}(x,z;y,t) +{{\partial K}\over{\partial t}}(x,z;y,t)\Bigr)\cr
&\quad +\Bigl(-3{{d}\over{dx}}{{\partial K}\over{\partial x}}(x,x;y,t)\phi (x,z;y,t)-3{{dK}\over{dx}}(x,x;y,t) {{\partial \phi}\over{\partial x}}(x,z;y,t)\Bigr)\cr
&\quad +\int_x^\infty \Bigl( {{\partial^3K}\over{\partial x^3}}(x,w;y,t) +{{\partial^3K}\over{\partial z^3}}(x,w;y,t) +{{\partial K}\over{\partial t}}(x,w;y,t)\Bigr)\phi (w,z;y,t)\, dw\cr
&=0.&(9.20)\cr}$$
\noindent We  differentiate the Gelfand--Levitan equation (8.4) with respect to $x$, and multiply by $(3/2)u=-3{{d}\over{dx}}K$, obtaining
$$\eqalignno{{{3}\over{2}}&u(x;y,t) {{\partial\phi}\over{\partial x}}(x,z;y,t) +{{3}\over{2}}u(x;y,t) {{\partial K}\over{\partial x}}(x,z;y,t)\cr
&\quad -{{3}\over{2}}u(x;y,t)K(x,x;y,t)\phi (x,z;y,t)+\int_x^\infty {{3}\over{2}}u(x;y,t){{\partial K}\over{\partial x}}(x,w;y,t)\phi (w,z;y,t)\, dw\cr
&=0;&(9.21)}$$
\noindent then we add multiples of the Gelfand-Levitan equation to eliminate the third summand, thus obtaining
$$\eqalignno{{{3}\over{2}}&u(x;y,t) {{\partial\phi}\over{\partial x}}(x,z;y,t) -3\Bigl( {{d}\over{dx}}{{\partial K}\over{\partial x}}(x,x;y,t)\Bigr) \phi (x,z;y,t)\cr
&\quad+\Bigl( {{3}\over{2}} u(x;y,t) {{\partial K}\over{\partial x}}(x,z;y,t) +{{3}\over{2}}u(x;y,t) K(x,x;y,t)  K(x,z;y,t)\cr
&\quad-3\Bigl( {{d}\over{dx}}{{\partial K}\over{\partial x}}(x,x;y,t)\Bigr) K(x,z;y,t)\Bigr)\cr
&\quad +\int_x^\infty  \Bigl( {{3}\over{2}} u(x;y,t) {{\partial K}\over{\partial x}}(x,w;y,t) +{{3}\over{2}}u(x;y,t) K(x,x;y,t)  K(x,w;y,t) -3\Bigl( {{d}\over{dx}}{{\partial K}\over{\partial x}}(x,x;y,t)\Bigr) \cr
&\qquad \times K(x,w;y,t)\Bigr)\phi (w,z;y,t)\, dz\cr
&=0;&(9.22)\cr}$$
  \noindent comparing this with (9.20), we deduce that

$$ \eqalignno{{{\partial^3K}\over{\partial x^3}}(x,z;y,t) &+{{\partial^3K}\over{\partial z^3}}(x,z;y,t) +{{\partial K}\over{\partial t}}(x,z;y,t)\cr
&= {{3}\over{2}} u(x;y,t) {{\partial K}\over{\partial x}}(x,z;y,t) +{{3}\over{2}}u(x;y,t) K(x,x;y,t)  K(x,z;y,t)\cr
&\quad -3\Bigl( {{d}\over{dx}}{{\partial K}\over{\partial x}}(x,x;y,t)\Bigr) K(x,z;y,t),&(9.23)\cr}$$
\noindent by uniqueness of solutions of the integral equation. \par
\indent By substituting this into the integral (9.14), we obtain
$$\eqalignno{ {}&{{\partial\psi_\zeta}\over{\partial t}}+ {{\partial^3\psi_\zeta}\over{\partial x^3}}-{{3}\over{2}}u (x;y,t){{\partial\psi_\zeta}\over{\partial x}}-{{3}\over{4}} {{\partial u}\over{\partial x}}\psi_\zeta  (x;y,t) -3w(x;y,t)\psi_\zeta (x;y,t)\cr
&= \Bigl( -3{{d}\over{dx}}{{\partial K}\over{\partial x}}(x,x;y,t) -{{3}\over{4}}{{\partial u}\over{\partial x}}(x;y,t) -3\alpha w(x;y,t) +{{3}\over{2}}u(x;y,t) K(x,x;y,t)\Bigr)e^{\zeta x-\zeta^2y/\beta -\zeta^3t}\cr
& \quad +\int_x^\infty  \Bigl( -3{{d}\over{dx}}{{\partial K}\over{\partial x}}(x,x;y,t) -{{3}\over{4}}{{\partial u}\over{\partial x}}(x;y,t) -3\alpha w(x;y,t) +{{3}\over{2}}u(x;y,t)K(x,x;y,t)\Bigr)\cr
&\quad \times  K(x,z;y,t)e^{\zeta z-\zeta^2y/\beta -\zeta^3t}\, dz&(9.24)\cr}$$
\noindent where the terms in ${\partial K}/\partial x$ have canceled. Now we consider the coefficient of $e^{\zeta x-\zeta^3t}$ and use (8.6) to obtain 
$$\eqalignno{-3{{d}\over{dx}}{{\partial K}\over{\partial x}}(x,x;y,t)& -{{3}\over{4}}{{\partial u}\over{\partial x}}(x;y,t) -3\alpha w(x;y,t) +{{3}\over{2}}u(x;y,t) K(x,x;y,t)\cr
&= -{{3}\over{2}}{{d^2K}\over{dx^2}}(x,x;y,t) -{{3}\over{2}}\Bigl( {{\partial^2 K}\over{\partial x^2}}(x,x;y,t)- {{\partial^2 K}\over{\partial z^2}}(x,x;y,t)\Bigr) \cr
&\quad +{{3}\over{2}}{{d^2K}\over{dx^2}}(x,x;y,t) -3\alpha w(x;y,t) +{{3}\over{2}}u(x;y,t) K(x,x;y,t) \cr
&= {{3}\over{2}}{{\partial K}\over{\partial y}} (x,x;y,t) -3\alpha w(x;y,t)\cr
&=0.&(9.25)\cr}$$ 
\noindent Likewise, the integrand vanishes and hence the result (9.5).\par
\indent Let ${\bf \Sigma}_{-A}$ be the space of $(2,2)$ admissible linear systems with
generator $(-A)$ and input and output space ${\bf C}$ where we identify $(-A, B,C)$ with
$(-A, -B,-C)$. The starting with $\Sigma =(-A,B,C)$, with tau function $\tau$ and
potential $u_\infty$, we build 
$$\Sigma_{\zeta_1}=(-A,(\zeta_1I+A)(\zeta_1I-A)^{-1}B, C)\eqno(6.8)$$ 
and proceed to form 
$$\Sigma_{\zeta_1,\dots , \zeta_n }=\bigl(-A,(\zeta_1I+A)(\zeta_1I-A)^{-1}\dots  
(\zeta_nI+A)(\zeta_nI-A)^{-1}B, C\bigr)\eqno(6.9)$$
\noindent with tau function $\tau_{\zeta_1, \dots, \zeta_n}$ and potential
$$u_{\zeta_1, \dots, \zeta_n}(x)=-2{{d^2}\over{dx^2}}\log\tau_{\zeta_1, \dots,
\zeta_n}(x).\eqno(6.10)$$
\vskip.05in
\noindent {\bf Proposition 6.2} {\sl (i) The multipliers $M^\zeta$ on ${\bf
\Sigma}_{-A}$ satisfy $M^\zeta M^\eta=M^\eta M^\zeta$, and 
$M^\zeta M^{-\zeta}=M^\infty$, where
$M^\infty =Id$.\par
\indent (ii) $M^0$ corresponds to the Darboux transformation $(-A, B,C)\mapsto
(-A, B,-C)$.\par
\indent (iii) $M^\zeta$ has the effect on potentials of taking $u_\infty$ to
$u_\zeta =u_\infty -2(\log \psi_\zeta )''$.\par
\indent (iv) $M^{\zeta_n}\dots  M^{\zeta_1}$ has the effect of taking $u_\infty$ to
$u_\infty -2(\log Wr(\psi_{\zeta_n} ,\dots , \psi_{\zeta_1}))''$.}\par 
\vskip.05in
\noindent {\bf Proof.} (ii) Note that $\Sigma_{\zeta_1,\dots , \zeta_n,\infty }
=\Sigma_{\zeta_1,\dots , \zeta_n }$, whereas 
$\Sigma_{\zeta_1,\dots , \zeta_n,0 }$ is the Darboux
transformation of $\Sigma_{\zeta_1,\dots , \zeta_n }.$ Evidently
$M^{\zeta_1}\dots M^{\zeta_n}\Sigma=\Sigma_{\zeta_1,\dots , \zeta_n }$. \par
 We can therefore apply Theorems 3.4 and 
Theorem 6.1 repeatedly. 

\indent (iii) The function $\psi_\zeta (x)=e^{\zeta x}\tau_\zeta (x)/\tau_\infty (x)$ satisfies 
$$\psi_\zeta''(x)=(u_\infty (x)+\zeta^2)\psi_\zeta (x);\eqno(6.11)$$
\noindent in particular, $\psi_0(x)$ satisfies $\psi_0''(x)=u_\infty \psi_0(x)$ by the Darboux transform
relation (3.7). Hence 
$$u_\infty (x)-2{{d^2}\over{dx^2}}\log\psi_0(x)=u_0(x)\eqno(6.12)$$
\noindent which is a form of Darboux's addition rule, as in [19, p
484], and 414 of [42].\par 

\indent (iv) First we consider the effect of $M^\zeta M^\eta$ on $u_\infty$. Note that $\psi_{\zeta, \eta}
=e^{\eta x}\tau_{\zeta , \eta}/\tau_{\zeta , \infty}$ satisfies
with 
$$\psi_{\zeta, \eta}''=(u_{\zeta, \infty }+\eta^2)\psi_{\zeta, \eta}.\eqno(6.13)$$ 
The effect of $M^\zeta M^\eta$ on potentials is to update by adding $-2(\log\psi_\zeta )''$
and then $-2(\log \psi_{\zeta, \eta})''$, so 
$$M^{\eta}M^{\zeta}: u_\infty\mapsto 
u_\infty-2(\log\psi_\zeta )''-2(\log \psi_{\zeta, \eta})''=u_{\zeta, \eta,
\infty}.\eqno(6.14)$$  

For $\zeta_1\neq \zeta_2$, let $\Psi (x)=Wr( \psi_{\zeta_1}, \psi_{\zeta_2})/\psi_{\zeta_2}$, and observe that
$$\Psi''=\bigl(\zeta_2^2+u_\infty -2(\log \psi_{\zeta_1})''\bigr)\Psi .\eqno(6.15)$$
\noindent Hence we can write
$$\eqalignno{u_{\zeta_2, \zeta_1, \infty}&=u_\infty -2(\log \Psi )''-2(\log \psi_{\zeta_2})''\cr
&=u_\infty -2\bigl(\log Wr (\psi_{\zeta_1},
\psi_{\zeta_2})\bigr)''.&(6.16)\cr}$$
\indent Now we consider $n$-fold products. For brevity we write $W(\zeta_1)=\psi_{\zeta_1}$ and 
$$W(\zeta_1, \dots, \zeta_n)= 
Wr(\psi_{\zeta_1},\psi_{\zeta_2}, \dots , \psi_{\zeta_n}).\eqno(6.17)$$

\indent Adding the points $\zeta_n , \zeta_{n+1}$ to the list $\zeta_1, \dots, \zeta_{n-1}$ can be implemented
by forming the Wronskian of quotients formed by adding $\zeta_n$ and $\zeta_{n+1}$ separately.
Indeed, by Jacobi's identity
$${{W(\zeta_1, \dots, \zeta_n, \zeta_{n+1})W(\zeta_1, \dots, \zeta_{n-1})}\over{
W(\zeta_1, \dots, \zeta_n)^2}}={{d}\over{dx}}{{W(\zeta_1, \dots, \zeta_{n-1}, \zeta_{n+1})}\over{
W(\zeta_1, \dots,\zeta_{n-1}, \zeta_n)}};\eqno(6.18)$$
so 
$${{W(\zeta_1, \dots, \zeta_n, \zeta_{n+1})}\over{W(\zeta_1, \dots, \zeta_{n-1})}}=Wr\Bigl( 
{{W(\zeta_1, \dots, \zeta_{n-1}, \zeta_{n})}\over{W(\zeta_1, \dots, \zeta_{n-1})}} ,
{{W(\zeta_1, \dots, \zeta_{n-1},\zeta_{n+1})}\over
{W(\zeta_1, \dots, \zeta_{n-1})}}\Bigr).\eqno(6.19)$$
\noindent Thus we can build the product $M^{\zeta_1}\dots M^{\zeta_n}$ and obtain the
corresponding potential by adding one point at a time and forming second order 
Wronskians.\par
\rightline{$\square$}\par
We take the upper  half plane as a model of hyperbolic space, so ${\bf C}_+=\{z\in {\bf C}: \Im z>0\}$, and let $PSL(2, {\bf R})$ be the group of M\"obius transformations $z\mapsto (az+b)/(cz+d)$ with $a,b,c,d\in {\bf R}$ and $ad-bc=1$; thus $PSL(2, {\bf R})$ acts on ${\bf C}_+$. Let $\Gamma$ be a discrete subgroup of $PSL(2, {\bf R})$, known as a Fuchsian group, and suppose that $\Gamma$ is the fundamental group of ${\cal E}$, so that ${\cal E}$ is congruent to ${\bf C}_+/\Gamma$. Suppose further that ${\cal F}$ is a fundamental cell for the action of $\Gamma$ on ${\bf C}_+$, so ${\bf C}_+=\bigcup_{\gamma\in \Gamma} \gamma ({\cal F})$; the group $\Gamma$ is finitely generated, if and only if there exists ${\cal F}$ that is bounded by finitely many faces.\par
\indent The group $\Gamma$ is called of convergence type if the hyperbolic distance $\rho$ on ${\bf C}_+$ 
satisfies $\sum_{\gamma\in \Gamma} e^{-\rho (i, \gamma (i))}$ converges. This is equivalently to  
$\sum_{\gamma \in \Gamma} \Im \gamma (i)/\vert \gamma (i)+i\vert^2$ converging, or 
or equivalently ${\bf C}_+/\Gamma$ having a Greens function. Hence the product
$$\prod_{\gamma\in \Gamma} {{\gamma (z)-w}\over{\gamma (z)-\bar w}}$$
\noindent converges for all $z, w\in {\bf C}_+$, except for a countable set of $w$ that give zeros.\par
\indent Let $\sigma (L)=\cup_{n=0}^\infty [\lambda_{2n}, \lambda_{2n+1}]$ be the spectrum, then let $G=\cup_{j=1}^\infty (\lambda_{2j-1}, \lambda_{2j})$ be the set of bounded gaps in $\sigma_L$ relative to ${\bf R}$. 
There is a question as to which sets can arise as spectral gaps from Hill's potentials with $q\in L^2([0,1]; {\bf R})$. Marcentko and Ostrovskii characterized these in terms of certain mapping functions from the upper half plane into slit domains, and Garnett and Trubowitz refined their results. We are interested in the case where $1\leq g<\infty$ and the closure of $G$ is a compact set. The following condition has been used in works by Totik and Peherstorfer. \par 
\vskip.05in

\noindent {\bf Definition} Say that a compact subset $\Sigma$ of ${\bf R}$ is admissible if \par
\indent (i) there exist $a_1<a_2<\dots <a_{2g}$ such that $\Sigma=\cup_{j=1}^g [a_{2j-1}, a_{2j}]$;\par
\indent (ii) the equilibrium probability measure $\mu_\Sigma$ has $\mu_{\Sigma} ([a_{2j-1}, a_{2j}])$ rational for all $j$.\par
\vskip.05in
\noindent {\bf Proposition 6.5} {\sl (i) Given $G=\cup_{j=1}^g (\lambda_{2j-1}, \lambda_{2j})$ and $\varepsilon >0$, there exists an admissible $\Sigma$ within $\varepsilon$ of $\cup_{j=1}^g (\lambda_{2j-1}, \lambda_{2j})$ in the Hausdorff metric.\par
\indent (ii) The interior of $\Sigma$ is a gap set for Hill's equation with potential $q\in L^2([0, \pi ]; {\bf R})$.}\par
\vskip.05in
\noindent {\bf Proof.} (i) This follows from Theorem 2.1 and Lemma 2.2 of [66]. Given $\lambda_1<\lambda_2<\dots <\lambda_{2g}$ and $0<\varepsilon <\min_j(\lambda_{j+1}-\lambda_j)$, there exist $a_1<a_2<\dots <a_{2g}$ and a real polynomial $T_N$ of degree $N\geq 2$ with simple real zeros $x_k$ $(k=1, \dots, N$ such that
$$\Sigma=\bigcup_{j=1}^g [a_{2j-1}, a_{2j}]=\{ x\in {\bf R}: -1\leq T_N(x)\leq 1\}\eqno(6.23)$$
\noindent  and $\vert a_j-\lambda_j\vert<\varepsilon $ for all $j$;  we use the customary notation $T_N$ to emphasize the analogy with Chebyshev polynomials of the first kind. 
Furthermore, Totik has shown that there exist $t_j\in (a_{2j}, a_{2j+1})$ such that the polynomials
$$r(x)=\prod_{j=1}^{g-1}(x-t_j), \quad h(x)=\prod_{j=1}^{2g} (x-a_j)\eqno(6.24)$$
\noindent determine the equilibrium measure for $\Sigma$ by the formula
$$\mu_\Sigma (dx)={{\vert r(x)\vert}\over{\pi \sqrt{\vert h(x)\vert}}}{\bf I}_\Sigma (x)\, dx.\eqno(6.25)$$
\noindent  A crucial property is that we can choose $N\in {\bf N}$ such that $n_\ell =N\int_{a_{2\ell-1}}^{a_{2\ell}} \mu_\Sigma (dt)$ is a positive integer for all $\ell =1, \dots, g$, so $\Sigma$ is admissible. Note that the rationality condition in (ii) introduces a complicated coupling on the subintervals of $\Sigma$, and one cannot move them around at will. \par
\indent (ii) Hochstadt observed that, for a finite gap Hill's potential, the simple periodic spectrum $(\lambda_j^o)$ determines the double periodic spectrum and the nontrivial roots $\lambda_j'$ of $\Delta'(\lambda )=0,$ so we can write
$${{\Delta'(\lambda )}\over{\sqrt{4-\Delta (\lambda )^2}}}=c_3{{\prod_{j=1}^n (\lambda -\lambda_j')}\over{\sqrt{ \prod_{j=0}^{2n}(\lambda -\lambda_j^o)}}}\eqno(6.26)$$
\noindent for some $c_3$. We approximate this algebraic function, via the Totik's approximation theorem. We cannot in general choose the branch points to be the original $\lambda_j$; however, we can select branch points which are close to the $\lambda_j$. Let 
$$H(z)=-N\int_{\Sigma} \log (z-t)\mu_\Sigma (dt) +N {\hbox{logcap}}(\Sigma )+i\pi N,\eqno(6.27)$$
\noindent which has real part satisfying
$$\Re H(z)\quad\cases{= 0,& for $x
z\in \Sigma$;\cr 
<0, & for $z\in {\bf C}\setminus \Sigma$\cr}\eqno(6.28)$$
\noindent by Frostman's theorem, and imaginary part
$$\eqalignno{\Im H(x)&=-N\int_\Sigma {\hbox{Arg}}(x-t)\mu_\Sigma (dt) +\pi N\cr
& =\int_{-\infty}^x N\pi \mu_\Sigma (dt)\cr
&=\sum_{\ell=1}^j \pi n_\ell&(6.29)}$$
\noindent for $a_{2j}<x<a_{2j+1}$. Further, $T_N(z)=\cosh H(z)$, so introduce 
$$\varphi (z)=\int_0^z {{-T'_N(s)}\over{\sqrt {1-T_N(s)^2}}}ds,\eqno(6.30)$$
\noindent which gives a conformal mapping from the upper half plane $\{ z: \Im z>0\}$ to the slit domain 
$$\Omega =\{ w\in {\bf C}: \Im w >0, \Re w>0\} \setminus \bigcup_{\ell =1}^g [\pi (n_1+\dots + n_\ell ), \pi (n_1+\dots +n_\ell )+ih_\ell ]\eqno(6.31) $$
\noindent where $\Omega$ looks like the first quadrant with vertical spikes cut above the points $(n_1+\dots +n_\ell )\pi$ for $\ell=1, \dots g$; the last spike is at $N\pi$. 
\noindent We define $\nu_{2n-1}=\varphi^{-1}(n\pi-)$ and $\nu_{2n}=\varphi^{-1}(n\pi+),$ which may be equal. The simple eigenvalues are
$$\lambda_{2j-1}^o=a_{2j-1}=\varphi^{-1}((n_1+\dots +n_j)\pi-),\quad \lambda_{2j}^{o}=a_{2j}=\varphi^{-1}((n_1+\dots +n_j)\pi +).\eqno(6.32)$$
\noindent As in Garnett and Trubowitz, there exists $q\in L^2([0, \pi ]; {\bf R})$ such that the simple periodic eigenvalues of Hill's equation with potential $q$ are $(\lambda_j^o)$.\par

\rightline{$\square$}\par
\vskip.05in

\indent Now we introduce a real polynomial $U_{N-g}$ of degree $N-g$ by 
$$T'_N(x)=NU_{N-g}(x)r(x)$$
\noindent such that $T_N$ and $U_{N-g}$ satisfy Pell's equation in the form
$$ T_N(x)^2-H(x)U_{N-g}(x)^2=1.$$
 The curve with equation $Z^2=\prod_{j=1}^{2g}(X-\lambda_j)$ is approximated by the curve with equation  $Z^2=\prod_{j=1}^{2g}(X-a_j)$; then with $W=ZU_{N-g}(X)$, we introduce the hyperelliptic curve
$${\cal F}=\{ (X,W): W^2=-(T_N(X)-1)(T_N(X)+1)\}.$$
\indent The function  $\log g_0(x;\zeta )\sqrt {\lambda_0-\zeta }$ arises when we solve the Dirichlet problem on ${\bf C}_+$ with boundary values given by a particular square wave. With $\lambda_{2j-1}<\mu_j<\lambda_{2j}$, let $G_\ell =\cup_{j=1}^g (\lambda_{2j-1}, \mu_j]$ be the left parts of the bounded gaps in the spectrum, and let $G_r=\cup_{j=1}^g (\mu_j, \lambda_{2j})$ be the right parts of the bounded gaps. The Poisson kernel $P_z(dt)$ gives harmonic measure in the upper half plane, so by a simple calculation we have
$$\eqalignno{\omega (G_\ell , \zeta )-\omega (G_r, \zeta )&=\int_{G_\ell} P_z(dt)-\int_{G_r}P_z(dt)\cr
&={{1}\over{2\pi}}\sum_{j=1}^g \Im \log {{(\mu_j-\zeta )^2}\over{(\lambda_{2j-1}-\zeta )(\lambda_{2j}-\zeta )}},&(5.51)\cr}$$
\noindent and we recognize the logarithmic term from the expression for $g_0(x;\zeta )$. Let $z_0\in G_\ell$, and let $z_j\in {\bf C}_+$ have $z_j\rightarrow z_0$ nontangentially as $j\rightarrow\infty$; then $\varphi (z_0)\in \partial \Omega$ and $\varphi (z_j)\rightarrow \varphi (z_0)$ as $j\rightarrow\infty$; we say that $\varphi (z_0)\in S_\ell$. Likewise we identify a subset $S_r$ of $\partial \Omega$, arising from $G_r$. Under the mapping $\varphi : {\bf C}_+\rightarrow \Omega$, we have $S_\ell =\varphi (G_\ell )$ and 
$S_r=\varphi (G_r)$. By the subordination principle, the harmonic measure of $S_\ell$ and $S_r$ as boundary components of $\Omega$ satisfy
$$\omega (S_\ell , \varphi (\zeta ))-\omega (S_r, \varphi (\zeta) )=\omega (G_\ell , \zeta )-\omega (G_r, \zeta ).\eqno(5.52)$$